\documentclass[12pt]{article}
\usepackage{amsmath,amsfonts,amssymb,amsthm}
\usepackage{mathcommand}
\usepackage{stackrel}
\renewcommand{\baselinestretch}{1.1}
\usepackage{enumitem,kantlipsum}
\usepackage{color}   
\usepackage[toc, page]{appendix}
\let\originalpart=\part
\def\part{\cleardoublepage\originalpart}

\allowdisplaybreaks

\usepackage{xcolor,soul} %
\definecolor{mediumorchid}{rgb}{0.73, 0.33, 0.83}
\definecolor{lawngreen}{rgb}{0.49, 0.99, 0.0}
\definecolor{sacramentostategreen}{rgb}{0.0, 0.34, 0.25}
\definecolor{amaranth}{rgb}{0.9, 0.17, 0.31}

\newcommand{\BBLp}[1]{{{\color{amaranth}#1}}}

\usepackage[normalem]{ulem}
\usepackage[numbers]{natbib}

\newcommand{\bk}[3]{[\mspace{-2mu}#1\mspace{-1.5mu}]_{#2}^{#3}}

\newcommand{\MC}{\mathcal{A}}

\newcommand{\amu}{\mathfrak{m}}  % for atomic measure in the cumulant method
\newcommand{\frakX}{\mathfrak{X}} 
\newcommand{\frakf}{\mathfrak{f}} 
 
\newcommand{\gB}{\mathfrak{B}}
\newcommand{\gT}{\mathfrak{T}}

\newcommand{\BL}{\textnormal{BL}}
\newcommand{\TV}{\textnormal{TV}}
\newcommand{\GeoBM}{\textnormal{}}

\RequirePackage[colorlinks,citecolor=blue,urlcolor=blue,backref=page]{hyperref}
\usepackage[final]{showlabels}
\hypersetup{
	pdftitle={Limit theory for  Lipschitz-localized statistics in  random geometric models},
	pdfauthor={B.\ Blaszczyszyn\ -\ D.\ Yogeshwaran\ -\ J.\ E.\ Yukich}
}

\usepackage{datetime2}

\usepackage{bookmark}

\usepackage[bottom=0.4cm]{geometry}
\usepackage{xcolor}
\usepackage{amssymb}
\usepackage{graphicx}
\usepackage{soul}
\usepackage{footmisc}
\usepackage{verbatim,cprotect}

\usepackage{times}

\usepackage{accents}
\usepackage[utf8]{inputenc} 
\usepackage[T1]{fontenc}
\newcommand{\dbtilde}[1]{\tilde{\raisebox{0pt}[0.85\height]{$\tilde{#1}$}}}

\makeatletter
\newcommand{\leqnomode}{\tagsleft@true}
\newcommand{\reqnomode}{\tagsleft@false}
\makeatother

\newcommand{\nicesets}{\mathcal{B}_b(r_0,n_0)}

\newcommand{\Deg}{\mathrm{Deg}}

\newcommand{\tu}{\tilde{u}}
\newcommand{\tv}{\tilde{v}}
\newcommand{\tz}{\tilde{z}}

\newcommand{\cT}{\mathcal{T}}
\newcommand{\cM}{\mathcal{M}}

\newcommand{\K}{\mathbb{K}}
\newcommand{\cN}{\mathcal{N}}
\newcommand{\cNK}{\mathcal{N}_{\mathbb{R}^d \times \mathbb{K}}}
\newcommand{\be}{\begin{equation}}
\newcommand{\ee}{\end{equation}}
\newcommand{\bea}{\begin{eqnarray}}
\newcommand{\no}{\nonumber}
\newcommand{\non}{\nonumber}
\newcommand{\eea}{\end{eqnarray}}
\newcommand\bmodif{\begin{modif}}
	\newcommand\emodif{\end{modif}}

\newcommand{\sP}{\mathbb{P}}  % Probability
\newcommand{\sE}{\mathbb {E} } % Expectation

\newcommand{\E}{\mathbb{E} }
\newcommand{\B}{\mathbb{B} }
\newcommand{\Z}{\mathbb{Z} }

\newcommand{\EXP}[1]{\mathbb{E}\!\left(#1\right) }
\newcommand{\Var}{{\rm Var}}

\newcommand{\Vol}{{\rm Vol}}
\newcommand{\remove}[1]{}
\newcommand{\tod}{\stackrel{d}{\Rightarrow}}

\newcommand{\perc}[3]{#1 \, \cap \, #2_{*#3}}

\newtheorem{theorem}{Theorem}[section]
\newtheorem{corollary}[theorem]{Corollary}
\newtheorem{lemma}[theorem]{Lemma}
\newtheorem{proposition}[theorem]{Proposition}
\newtheorem{rem}[theorem]{Remark}

\newenvironment{remark}[1][]{\begin{rem}[#1]\rm}{\end{rem}}
\newtheorem{definition}[theorem]{Definition}
\newtheorem{defn}[theorem]{Definition}
\newtheorem{example}[theorem]{ Example}
\newenvironment{exe}[1][]{\begin{example}[#1]\rm}{\end{example}}

\numberwithin{equation}{section}
\newtheorem{innercustomass}{Assumption}
\newenvironment{customass}[1]
{\renewcommand\theinnercustomass{#1}\innercustomass}
{\endinnercustomass}

\numberwithin{equation}{section}

\def\0{{\bf 0}}

\def\Rar{\Rightarrow}
\def\rar{\rightarrow}

\newcommand{\ze}{\zeta}

\newcommand{\ka}{\kappa}

\newcommand{\bxi}{\boldsymbol{\xi}}
\newcommand{\bzeta}{\boldsymbol{\zeta}}
\newcommand{\ttau}{ \tilde{\tau}}

\newcommand{\F}{\Phi}
\newcommand{\G}{\Psi}
\def\T{\top}

\def\Ph{\Phi}
\def\N{\mathbb{N}}

\def\mR{\mathbb{R}}

\def\mP{\mathbb{P}}
\newcommand{\M}{\mathbb {M} }
\newcommand{\hx}{\hat x}

\def\mE{\mathbb {E} \,}
\def\mL{\mathbb{L}}

\def\X{{\cal X }}

\renewmathcommand{\v}{{v}} % old command  \def\v{{v}}
\def\z{{z}}

\newcommand{\ga}{\mathfrak{a}}
\newcommand{\gta}{\tilde{\mathfrak{a}}}

\newcommand{\gn}{\mathfrak{n}}

\newcommand{\oB}{\mathfrak{G}}

\newcommand{\md}{\mathrm{d}}

\newcommand{\cC}{{\cal C}}
\newcommand{\cO}{{\cal O}}
\newcommand{\cX}{{\cal X}}
\newcommand{\cW}{{\cal W}}
\newcommand{\cZ}{{\cal Z}}

\newcommand{\gC}{\mathfrak{C}}

\newcommand{\tA}{\tilde{A}}
\newcommand{\tB}{\tilde{B}}
\newcommand{\tXi}{\tilde{\Xi}}
\newcommand{\txi}{\tilde{\xi}}

\newcommand{\xitr}{\xi^{(r)}}

\def\hT{\hat{T}}

\def\hO{\hat{\omega}}

\def\R{\mathbb{R}}
\def\N{\mathbb{N}}
\def\K{\mathbb{K}}
\def\tP{\tilde{\P}}

\def\tC{\tilde{C}}
\def\tx{\tilde{x}}
\def\ty{\tilde{y}}
\def\tz{\tilde{z}}
\def\tmu{\tilde{\mu}}

\def\P{\mathcal{P}}
\def\Q{\mathcal{Q}}

\def\G{\mathcal{G}}

\newcommand{\1}[1]{\mathbf{1}(#1)}

\newcommand{\xil}{\Upsilon}
\newcommand{\bxil}{\boldsymbol{\xil}}
\newcommand{\bfL}{\mathbf{L}}
\newcommand{\bfM}{\mathbf{M}}

\def\Palm{\mathbb{P}}  %\mathfrak{P}

\gdef\GBkr{\G_{\bk x1p}^{\cup[k_r]}}
\gdef\VBkr{V_{\bk x1p}^{\cup[k_r]}} 
\gdef\GBr{\G_{\bk x1p}^{\cup(r)}}
\gdef\VBr{V_{\bk x1p}^{\cup(r)}} 
\gdef\VP{V_{\bk x1p}} 
\gdef\VBlkr{V_{\bk x1l}^{\cup[k_r]}} 
\gdef\VBpkr{V_{\bk x{l+1}p}^{\cup[k_r]}} 
\gdef\GBlr{\G_{\bk x1l}^{\cup(r)}}
\gdef\GBpr{\G_{\bk x{l+1}p}^{\cup(r)}}
\gdef\GBlkr{\G_{\bk x1l}^{\cup[k_r]}}
\gdef\GBpkr{\G_{\bk x{l+1}p}^{\cup[k_r]}}
\gdef\VBlr{V_{\bk x1l}^{\cup(r)}}
\gdef\VBpr{V_{\bk x{l+1}p}^{\cup(r)}}
\gdef\VPl{V_{\bk x1l}} 
\gdef\VPp{V_{\bk x{l+1}p}}

\newcommand{\gbeta}{b}

\newcommand{\sgn}{\mathrm{sgn}}

\def\qed{\hfill\hbox{${\vcenter{\vbox{
    \hrule height 0.4pt\hbox{\vrule width 0.4pt height 6pt
    \kern5pt\vrule width 0.4pt}\hrule height 0.4pt}}}$}}

\newcommand{\exclude}[1]{}

\newcommand{\nocontentsline}[3]{}
\let\origcontentsline\addcontentsline
\newcommand\stoptoc{\let\addcontentsline\nocontentsline}
\newcommand\resumetoc{\let\addcontentsline\origcontentsline}
\newcommand{\dy}[1]{\textcolor{magenta}{#1}}

\newcommand{\jy}[1]{\textcolor{blue}{#1}}

\title{Limit theory for  Lipschitz-localized statistics in  random geometric models}
\author{B.~B{\l}aszczyszyn\footnote{Inria/ENS, 48 rue Barrault
75013, Paris, France. Email~: blaszczy@ens.fr}  \ D.~Yogeshwaran\footnote{Theoretical Statistics and Mathematics Unit,  Indian Statistical Institute,  Bangalore.  Email~: d.yogesh@isibang.ac.in}\  \ and  \ J.\ E.~Yukich
\footnote{Department of Mathematics, Lehigh University, Bethlehem, PA 18015.  Email~: jey0@lehigh.edu}}
\date{\today}

\begin{document}

\maketitle

\begin{abstract}

We consider spatial models represented by marked point processes, which incorporate randomness through a random collection of points (referred to as {\em sites}) in \(\mathbb{R}^d\) and random marks (referred to as {\em states}) taking values in a general Polish space. Generally, the locations and states are assumed to be dependent, and the states may also evolve over time, introducing spatially dependent stochastic processes at these locations.

We focus on statistics \((H_n)_{n \in \N}\) of these models, which take the form of sums of locally dependent {\em score functions} of sites and states observed within expanding windows \(W_n\!\!:=\! \![-\frac{1}{2} n^{1/d}, \frac{1}{2} n^{1/d}]^d\allowbreak \subset \mathbb{R}^d\). As \(n \to \infty\), we establish Gaussian fluctuations for the centered and normalized statistics of \(H_n\),  which need not necessarily be linear functionals of the underlying marked point process.  Under the additional assumption of stationarity and translation invariant score functions, we  establish the asymptotics of the mean and variance of these statistics, leading to a multivariate central limit theorem for vectors of such statistics.

To achieve these general limit results, we introduce a foundational approach that captures asymptotic independence in random spatial models via a new {\em mixing  condition on the correlations of}  marked point processes. We also provide a more general stabilization paradigm, termed {\em \BL-localization}, for marks, which reinforces this new mixing condition. Both concepts rely on the bounded Lipschitz metric applied to the Palm distributions of marked point processes. Notably, our localization criterion is weaker than the classical stopping set stabilization---allowing for score functions whose classical stabilizing radii need not be bounded---but strong enough to yield the asymptotic normality of \((H_n)_{n \in \N}\) via factorial moment expansions and the cumulant method. Our approach does not require the scores to be defined on the infinite window. In the stationary setting, this also enables the derivation of expectation and variance asymptotics by showing existence of distributional limits of scores in the infinite windows.

Spatial models  falling 
within the scope of our general theorems include
spin systems on spatial random graphs, sparse networks of interacting diffusions, and particle systems with finite-time horizon.
Interacting particle systems considered here include continuum versions of exclusion models, 
generalized random sequential adsorption and ballistic deposition models,
finite time-horizon epidemic and voter models,  and majority models.
We also consider statistics of empirical random fields and geostatistical models defined on spatial random graphs.  
\end{abstract}



\phantomsection
% \addcontentsline{toc}{section}{}
\renewcommand{\baselinestretch}{0.8}\normalsize
	     {\hypersetup{linkcolor=black}\small
\pdfbookmark{\contentsname}{toc}
\tableofcontents
}	\renewcommand{\baselinestretch}{1.1}\normalsize	

%\newpage
%\cleardoublepage

%\part*{Part I \quad Introduction}
\part{Introduction}
\label{part:intro}
%\section{Introduction}

This introductory part presents key examples motivating our study of statistics of marked point processes 
and spatial random models.  
Besides these examples,  Section \ref{s:Intro} describes  the new foundational aspects underpinning the general limit theory developed here,  provides a concise review of the relevant literature, and outlines the overall structure of our work.
 Section \ref{s:notation}
provides the basic terminology and definitions.

\section{Overview and contributions}
\label{s:Intro}

This work considers the limit theory for  statistics of spatial random models which are asymptotically decorrelated over spatial domains, as  the spatial domain increases up to $\R^d$.
This encompasses Gaussian fluctuations  of sums of locally dependent score functions evaluated over expanding observation windows, as well as asymptotic results for the mean and variance of these sums in the case of stationary models.
The models involve multiple  sources of randomness, namely the random collection of sites in $\R^d$, their random initial states, and, in some cases,  additional randomness stemming from the evolution of these states. Our set-up is general and yields the limit theory for statistics of stochastic geometric structures such as continuum spin systems, interacting diffusions, interacting particle systems, geostatistical models and empirical random fields. To
illustrate, in the following section we present six 
models with differing interaction mechanisms,  
but nonetheless statistics of these models may be handled with a common framework.

\subsection{Motivating examples}
\label{ss.MotivExamples}
The common input in our representative examples is a point process $\P= \{x_i\}$, namely a countable set of random spatial locations in $\R^d$, also called {\em sites} in the case of particle systems. We may, for example,  take $\P$ to be either the Poisson point process on $\R^d$, a Neyman-Scott or 
Mat\'ern cluster process on $\R^d$,   the Ginibre point process on $\R^2$, or a stationary alpha-determinantal point process with an exponentially decreasing kernel. Given the sequence of windows $W_n:= [-\frac{1}{2} n^{1/d}, \frac{1}{2} n^{1/d}]^d, n \in \N,$ we denote the restriction of $\P$ to $W_n$  by $\P_n := \P \cap W_n$. In many of our models there is an underlying spatial random graph $\G(\P)$ or $\G(\P_n)$ defined on $\P$ or $\P_n$, respectively, and called the {\em interaction graph}. For example, one may take this to be the $k$-nearest neighbor graph, the Delaunay graph, or the Gilbert graph, i.e., the random geometric graph. More generally, the spatial random graphs $\G(\P)$ have 
finite average degree and the edges are defined by local rules but are not necessarily uniformly bounded.

\begin{exe}[Spin systems; Section \ref{s:gibbsmarking}]
\label{Ex:Spin}  Consider a Gibbs random field defined on an interaction graph  $\G(\P_n)$,  or equivalently a spin system $\{V_x\}_{x \in \P_n}$
defined with respect to an adjacency relation in  $\G(\P_n)$.  Assuming the spin system satisfies an averaged version of weak spatial mixing,  we show that the total spin
$\sum_{x \in \P_n}V_x$ is asymptotically normal as $n \to \infty$.  We use both
combinatorial and disagreement percolation methods to deduce that models such as the hard-core and Widom-Rowlinson spin models on certain proximity graphs satisfy this version of  weak spatial mixing, and consequently, establish that their total  spin
$H_n := \sum_{x \in \P_n}V_x$ is asymptotically normal.
\end{exe}

\begin{exe}[Interacting diffusions; Section \ref{NEW-s:id_sprg}]  \label{Ex:ID} Consider a system of interacting $\R^{d'}$-valued diffusions $M(x, t):= M(x,t, \P),  t > 0$, $x \in \P$, defined on the vertices  of a locally finite graph $\G(\P)$. The diffusion at $x \in \P$
interacts directly with only the diffusions at its neighbors in the graph.  Assuming that diffusions start inside a bounded ball in $\R^{d'}$ and the drift and diffusion coefficients satisfy Lipschitz conditions, we establish the limit theory for statistics of the diffusions, such as $H_n := \sum_{x \in \P_n} \|M(x,t)\| $, where $\|\cdot\|$ is the sup-norm for functions in the interval $[0,t_0]$. 
\end{exe}

\begin{exe}[Cooperative sequential adsorption on graphs; Section \ref{ss:csa}]
\label{Ex:CSA}  With $\P$ as above, we encode the interactions between the sites of $\P_n$ by a proximity graph, such as the $k$-nearest neighbor or  Delaunay graph.  
The dynamics of the model go as follows.  Let the initial state be $M(x,0) = M(x,0, \P_n) = 0, x \in \P_n$.  Each site is equipped with a unit rate Poisson clock. When the clock  at site $x$ rings, a particle arrives at site $x$. The particle is accepted at site $x$ with a probability depending on the configuration of occupied neighbors of $x$ in the interaction graph. Once accepted the particle stays forever.  We let $M(x,t):= M(x, t, \P_n)$  be the state of the site $x$ at time $t$. Thus  $M(x, t) = 1$ if a particle is accepted at site $x$ by time $t$, otherwise put $M(x, t) = 0$. Setting $H_n = \sum_{x \in \P_n}M(x, t_0)$ gives the total number of accepted particles at time $t_0$.  We establish the limit theory as $n \to \infty$  for 
 $(H_n)_{n \in \N}$ and related 
 statistics of the evolved system up to time $t_0 \in (0,\infty)$;
 the general framework shows there is a thermodynamic limit for the spatial average $n^{-1}H_n$ and that its fluctuations are Gaussian.
 \end{exe}

\begin{exe}[Continuum interacting particle systems; Section~\ref{s:applnsips}]
\label{Ex:EP}
As in Example~\ref{Ex:CSA}, let $\P_n$ be equipped with a local interaction
graph (e.g., $k$-nearest neighbor or Delaunay graph), and assume that each site
carries a unit rate Poisson clock.
Unlike cooperative sequential adsorption, the dynamics now allows particles to
\emph{move} between neighboring sites according to an exclusion-type rule that
may depend on the local configuration of occupied neighbors.
Initial states may be dependent.
The system is evolved on $W_n$ up to a fixed time horizon $t_0 \in (0,\infty)$.
In the limit $n \to \infty$, we establish limit theorems for additive spatial
statistics at time $t_0$, such as the fraction of occupied sites, or the
fraction of occupied sites having a neighbor within a fixed distance.

This framework covers a broad class of continuum interacting particle systems,
including exclusion processes, ballistic deposition, majority dynamics,
epidemic models (contact process, SIR, voter models), and toppling dynamics
(e.g., divisible sandpiles).
\end{exe}

Problems in spatial statistics falling within the scope of  our general theory include the following two examples, each involving two independent point processes.

\begin{exe}[Empirical random fields; Sections~\ref{ss.Empirical-field}--\ref{s:Covariogram}] 
 \label{Ex:ERF}
 Let $M=\{M(x)\}_{x\in\R^d}$ be a stationary real-valued random field on $\R^d$,
and let $\P$ be a stationary point process on $\R^d$, assumed to be independent
of $M$.
The \emph{empirical random field} of $M$ sampled at the points of $\P$ over a
window $W_n$ is defined by 
\(\{(x,M(x)): x \in \P_n\}
\).
Under suitable assumptions on the spatial correlation structure of  $\P$ and $M$,
we establish asymptotics for the
expectation and Gaussian fluctuations of this empirical random field as
$n\to\infty$. As an immediate consequence, we obtain limit theorems for empirical
distribution-type statistics of the real valued random field $M\in\R$, such as
\(n^{-1}\sum_{x \in \P_n} \1{M(x) < \tau}$, $\tau \in \R$.

More involved statistics can also be considered.
In particular, one may study a {\em covariogram} estimator of the random
field $M$ of the form
\[
\frac{1}{n}
\sum_{\substack{x,y \in \P_n \\ x-y \in B_\delta(\pm h)}}
\bigl| M(x) - M(y) \bigr|^2 ,
\]
where pairs of points are selected whose displacement lies in a ball of radius
$\delta>0$ around $\pm h \in \R^d$.
Under appropriate spatial correlation assumptions on $\P$ and on the random
field $M$, we derive the asymptotic behavior of the expectation and variance of
this estimator, as well as Gaussian confidence intervals for this estimator of  the covariogram of~$M$.
\end{exe}

\begin{exe}[Geostatistical Gilbert disc models; Section \ref{ss:geoBM}]
\label{Ex:GEO}  Consider a generalization of the Gilbert disc model, in which balls are positioned at the realization of a stationary point process $\P$ in $\R^d$ with bounded but dependent radii determined by an underlying stationary real valued random field
$M = \{M(x)\}_{x \in \R^d}$ which is independent of $\P$.  In this variant of the Boolean model, two points communicate with one another, i.e., are joined with an edge, if the balls at both points intersect, thereby giving rise to a random graph $\G(\P)$ on $\P$ with spatially correlated edge lengths. 
For example the field value $M(x)$, $x \in \R^d$, may be determined by the position of $x$ with respect to
an independent Poisson-Voronoi tesselation of $\R^d$.  The total edge length in the random graph 
 $\G(\P \cap W_n)$ exhibits Gaussian fluctuations and, after scaling, its average converges to an asymptotic limit as $n \to \infty$.
\end{exe}
This work studies the afore-mentioned systems and 
rigorously establishes the above statements. This is achieved by providing a general  framework for establishing the limit theory for statistics of a broad class of geometric structures, one encompassing the above examples among others.  

The question of limit theory common to the above examples can be framed in terms of establishing asymptotics for the sum of scores
\be \label{linearstatsum}
H_n = H_n^\xi:=  \sum_{x \in \P_n} \xi( \tx, \tP_n ), \quad   \ \ n \to \infty,
\ee
where \( \xi(\tx,\tP_n) \), assigned to the sites  \( x \in \P_n \), is a family of real-valued characteristics. These characteristics are referred to as {\em scores}, though they may also be called {\em marks} or {\em states}, depending on the context. Here, \( \tP_n = \{(x,U(x))\}_{x \in \P_n} \)  represents the sites \(x\), which may also be  associated with {\em initial states} or {\em pre-marks} \(U(x)\), collectively denoted as \( \tx = (x, U(x)) \). The configuration of sites and their pre-marks (if present) induces the real-valued characteristics \( \xi(\tx,\tP_n) \) in a potentially complex yet localized manner. For instance, by making suitable choices of the pre-marks $U$, we may define $\xi(\tx,\tP_n)$ to coincide
with $V_x$ in the case of spin systems,
with $\|M(x,t)\|$ in the setting of interacting diffusions, or with $M(x,t_0)$ in the context of cooperative sequential adsorption. Analogously, $\xi$ can be defined in other classes of models as well.

More generally,  we seek limit theory for linear statistics of the associated (possibly) signed measures
\be \label{linearstatB}
\mu_n^\xi: = \sum_{x \in \P_n} \xi( \tx, \tP_n ) \delta_{ n^{-1/d}x },  \quad
\ \ n \to \infty.
\ee 
Questions about spatial averages can be cast in terms of establishing expectation asymptotics as $n \to \infty$ for  \eqref{linearstatsum} and \eqref{linearstatB}.

The set-up taken here allows for the following sources of dependencies, each of which brings extra spatial correlations: 
\begin{itemize}
\item  the sites $\P$ may have spatial correlations and  the  initial states $U_x, x \in \P$, may be dependent;

\item the scores at sites  depend on pre-marks  (initial states) through an interaction graph whose edges may not be uniformly bounded; 

\item  the scores and pre-marks may also influence the scores of neighboring sites, introducing additional dependencies among them, as in spin systems.   In interacting diffusions or particle systems, for example, updates at a site  may also modify the states at all neighboring sites; this locally synchronous updating condition is sufficient  to cover many models. 

\item  in case of interacting diffusions or particle systems, the updates at a site $x$ need not be Markovian, but may be a function of the entire time-evolved history at $x$ and at nearby sites.   
\end{itemize}

Furthermore, as  mentioned in Section \ref{s:future}, there are potential applications to statistics of random graphs and networks with dynamically changing geometry,  shot-noise models on spatially correlated input,  dynamic Markov random fields in the continuum, as well as to population genetics models.  We do not develop these areas but leave them for future research.

\subsection{Contributions of this work}
\label{ss.Contribution}
The contributions of this work are three-fold. In addition to the already indicated diverse applications, we put forward two  new foundational aspects underlying the forthcoming general limit theorems for the statistics  \eqref{linearstatsum} and \eqref{linearstatB}:
\begin{itemize}
\item we introduce the new notion  of mixing involving {\em correlations for marked point processes}, which uses geometric criteria to systematically capture the asymptotic independence exhibited in various spatial random models
\item we provide a more general localization paradigm, {\em termed \BL-localization}  (with \BL~standing for bounded, Lipschitz), for marks, which broadly reinforces this new type of mixing and which is well suited for quantitatively describing dependence of scores on local data even in the presence of unbounded model dependencies. A key mechanism connecting the localization and mixing of correlations of marks is the systematic use of the bounded Lipschitz metric with respect to the Palm distributions of these marks. In examples such as spin systems, we use a variant termed  {\em \BL~cluster-localization}.
\end{itemize}
 In the following sections, we  shed light on this framework and the proof ideas behind some of the main results.

\vskip.3cm

\subsubsection{Mixing correlations of marked point processes}
\label{sss-intor.asymptotic-decorrelation} Dynamic geometric models have multiple sources of dependencies with each bringing in extra correlations. To capture these spatial dependencies, we introduce the concept of {\em mixing correlations of points (sites) and marks (states)}, which quantifies the asymptotic decorrelation of marked point processes $\{(x,\xi(\tx,\tP_n))\}_{x \in \tP_n}$ and, as such,  undergirds the entire work.  Using the shorthand $\xi_{i,n} := \xi(\tx_i,\tP_n)$, this new  mixing  property says that  the expected value of the product of any two bounded Lipschitz (BL) test functions  $f,g$, when respectively evaluated on  the marks of any two point sets $\{x_1,\ldots,x_p\}$ and $\{x_{p+1},\ldots,x_{p+q}\}$ in $W_n$, and denoted by 
$$\mE_{\bk x1{p+q}}\bigl[f(\xi_{1,n},\ldots,\xi_{p,n})g(\xi_{p+1,n}\ldots,\xi_{p+q,n})\bigr]\rho^{(p+q)}(x_1,\ldots,x_{p+q}),
$$
approximately factorizes into the product of two expectations
$\mE_{\bk x1{p+1}}\bigl[f(\cdots)\bigr]\rho^{(p)}(\cdots)$ and \\ $\mE_{\bk x1{p+1,p+q}}\bigl[g(\cdots)\bigr]\rho^{(q)}(\cdots)$ calculated separately on these two groups of points, where  $\mE_{[\cdots]}$ stands for 
the Palm formalism  involving  expectations $\mE_{\bk x1{p+q}}$ of $\P_n$ given the considered locations $x_1,\ldots,x_{p + q}$ and where $\rho^{(p+q)}$ denotes their correlation functions. We require that the  additive error of this approximate factorization decays to zero faster than any power of the separation distance 
$$s:=\inf_{i \in \{1,...,p\}, j \in \{p + 1,...,p + q\} } |x_i - x_j|$$
between these two group of  sites, uniformly 
over the windows $(W_n)_{n \in \N}$. Instead of the BL test functions one could use the class $\B$ of bounded test functions. This gives rise to two mixing conditions, one for each of the two classes BL and $\B$, formally described in Definitions~\ref{d.omegahmixing} and~\ref{B-omega-mixing}, respectively.

These two new types of mixing  for marked point processes  (called \BL-mixing correlations in case the test functions are bounded Lipschitz and $\B$-mixing correlations for all bounded functions) are inspired by the notion of `clustering of the correlation functions $\rho^{(p+q)}(x_1,\ldots,x_{p+q})$' of {\em unmarked} point processes introduced in statistical physics by \citet{Martin80} and reprised in \citet{Nazarov12}, \citet{BYY19},   and 
\citet{Fenzl2019asymptotic}. The former, involving \BL \,functions, can also be interpreted as a variant of the $(\BL,\theta)$-dependence for random fields indexed by $\Z^d$, as discussed in \citet[Section 10.1.6]{bulinski2013central}, albeit with less restrictive conditions on the decay constants but stronger requirements on the decay rate. Our mixing property presented here serves as an alternative to the classical mixing conditions considered by \citet{Heinrich_Schmidt_1985,Heinrich2013asymptotic,HLS,HP,Ivanoff82,JP}, among others. Unlike these classical approaches, our criterion is formulated through {\em geometric} conditions, which may in some cases be simpler to verify.

We prove in Theorem \ref{t:clt_linear_marks} that \BL-mixing correlations of marked point processes along with moment conditions suffice to prove asymptotic normality of the associated $\xi$-weighted measures $(\mu_n^\xi)_{n \in \N}$ at \eqref{linearstatB}. More precisely, for bounded measurable test functions $f$ on $W_1$, the integrals $\mu_{n}^{\xi}(f)%:= \langle \mu_n^\xi, f \rangle 
:= \int f d\mu_n^\xi$, centered and suitably normalized, converge to a standard normal distribution as $n \to \infty$. This is because \BL-mixing of correlations of the point process marked by~$\xi$  implies the Brillinger mixing of the random measures $(\mu_n^\xi)_{n \in \N}$ provided $\xi$ satisfies moment conditions.

This approach is applicable  for  general point processes without stationarity assumptions on $\P$.  The approach rests  on establishing Theorem \ref{l.CLT-Cumulants},  a central limit theorem for a sequence of purely atomic random signed measures.  This theorem, whose proof depends on the cumulant method, appears to have independent interest as it unifies and generalizes some of the earlier abstract central limit theorems for statistics of point processes.

Mixing correlations of marked point processes  arises whenever the points and marks separately possess a decorrelation structure; see Proposition \ref{p:pmix+mmix}.  Moreover, Theorem~\ref{t:mppmix_double} shows that this approach extends to the analysis of more complex marks constructed functionally from simpler ones (such as $\xi_{i,n} = \xi(\tx_i,\tP_n)$ mentioned above in \eqref{linearstatsum}), provided that the initial marked process $\tP$ satisfies stronger $\B$-mixing correlations, in particular with the additive error decaying exponentially to zero.  
For this, we use factorial moment expansion from ~\citet{Blaszczyszyn95,Bartek97}, together with localization or stabilization of scores (discussed in the next section) and moment conditions on $\xi$ .

\subsubsection{\BL-localization of marks and limit theory}
\label{sss.BL-contribution}
We introduce a geometric localization criterion, here called fast {\em \BL-localization}, for marks (scores) via the bounded Lipschitz metric on the space of probability distributions, here denoted $d_{\BL}$. 
Specifically, fast {\em \BL-localization} means that for any localization ball \( B_r(x) \) of radius \( r \),  under the Palm  expectation $\E_x$, the \( d_{\BL} \) distance between \( \xi(x, \tP_n \cap B_r(x)) \) and \( \xi(x, \tP_n) \) decreases faster than any inverse power of \( r \), uniformly in $n$.  
This criterion appears to be  new even in the setting of scores $\xi(x,\P)$ 
on unmarked point processes $\P$. 
Similar to the classical notion of stabilization via stopping sets, \BL-localization measures the dependence of marks on local (marked) data. However, this dependence is quantified using the distributional metric \( d_{\BL} \), rather than through stopping sets. Although weaker than classical geometric criteria such as the existence of stopping sets or even \( L^q \)-stabilization, the multi-site version of \BL-localization---where for each \( p \in \N \), the joint distribution of the scores \( \xi(x_1, \tP_n), \ldots, \xi(x_p, \tP_n) \) under \( \mP_{\bk{x}{1}{p}} \) depends only on the marked configuration inside the balls \( B_r(x_i) \)---is still strong enough to ensure that the family of marked point processes \( \{(x, \xi(\tx, \tP_n))\}_{x \in \P_n} \) exhibits fast \BL-mixing correlations uniformly in \( n \), provided the base marked process \( \tP \) satisfies summable exponential \( \B \)-mixing correlations, as shown in Theorem~\ref{t:mppmix_double}. 

Importantly, such \( \BL \)-mixing of the newly constructed marking is a key condition in Theorem~\ref{t:clt_linear_marks} for proving the asymptotic normality of the empirical measures \( (\mu_n^\xi)_{n \in \N} \). Therefore, \BL-localization or \BL~cluster-localization serve as a practical and verifiable criteria to guarantee asymptotic normality for both \( (H_n)_{n \in \N} \) and \( (\mu_n^\xi)_{n \in \N} \), as long as the score function \( \xi \) satisfies suitable moment assumptions; see Theorem~\ref{t:cltmarkedpp}.

The \BL-localization criterion does \emph{not} require any pathwise or
realization-wise relationship between the scores and their truncated versions.
Instead, it is formulated purely in terms of their probability distributions.
This allows us to establish closeness (in the \( d_{\BL} \) metric) between the laws of the random variables \( \xi(x, \tP_n \cap B_r(x)) \) and \( \xi(x, \tP_n) \), even when their realizations are pathwise non-comparable.   
This is why we use the term `localization', in contrast to `stabilization', which involves comparing specific realizations of these random variables---either exactly (via stopping sets) or approximately (in \( L^q \)).  
This flexibility in \BL-localization plays a key role in establishing asymptotic normality for statistics of spin models and interacting diffusions on spatial random graphs. 

Theorem \ref{t:cltmarkedpp} has wide applicability and in particular covers the case when the sites form a point process having fast mixing  correlations and the marks $\xi$ exhibit exponentially fast \BL-localization or \BL~cluster-localization. 
The approach requires neither that $\xi$ satisfy  power growth conditions
with respect to the cardinality of the underlying point set, nor does $\xi$ need to be defined on infinite point processes on all of $\R^d$; thus even in the setting of unmarked point processes and stopping set stabilization, we improve upon existing results; cf. Section \ref{s:remarksltthms}.

Subject to the stationarity of $\P$ and the fast \BL-localization of $\xi$, we also 
provide expectation and variance asymptotics
of $(H_n)_{n \in \N}$ and 
$(\mu_n^\xi)_{n \in \N}$; cf. Propositions \ref{expvar} and \ref{t:clt_linear_marks_new}. Here again, unlike  previous limit results for $(H_n)_{n \in \N}$ and 
$(\mu_n^\xi)_{n \in \N}$, we require neither existence of $\xi$ nor `some limiting version' of $\xi$ on the infinite point process but rather prove and use the existence of distributional limits via Palm theory, which resembles the approach via local weak convergence.

Theorem~\ref{t:cltmarkedpp} is a qualitative central limit theorem and we make no attempt here to assess the proximity between  \( (\frac{H_n - \E H_n} 
{\sqrt{\Var H_n }} )_{n \in \N} \)
and the standard normal.  In the special case that  $(H_n)_{n \in \N}$ are functionals  of a Poisson point process $\P$ on a metric measure space, then
 \citet{TY} show that the
Malliavin-Stein calculus, as developed by 
\citet{last2016normal, Trauthwein}, 
yields  presumably optimal rates of convergence of  \( (\frac{H_n - \E H_n} 
{\sqrt{\Var \,H_n }} )_{n \in \N} \) to the standard normal, assuming only fast \BL-localization of 
the score $\xi$, together with a fifth moment condition on  $\xi$.  The rates of convergence, given in terms of the Kolmogorov and Wasserstein distances,  improve upon 
\citet{Lachieze2019normal} which assumes 
that scores satisfy a stronger stopping set exponential stabilization criterion and which assumes the underlying metric measure space is either Euclidean space or a space satisfying a growth assumption on the measure of spheres. 

\subsubsection{\BL-localization and  stabilization in applications} To apply the general limit theorems for \eqref{linearstatsum} and \eqref{linearstatB} in applications as described in Examples \ref{Ex:Spin}-\ref{Ex:GEO}, one needs to verify suitable localization  and moment conditions for appropriate $\xi$ and thereby establish asymptotic decorrelation of $(\mu_n^\xi)_{n \in \N}$. This constitutes a significant part of the article as these applications are always not well explored in spatial settings and both the framework and verification of localization needs to be tailored to fit within our general theory described above.   

In case of spin systems as in Example \ref{Ex:Spin}, the \BL~cluster-localization of spins follows from an averaged weak spatial mixing assumption of the spin system and good stabilizing properties of the underlying interaction graphs. For bounded degree graphs, we rely upon bounds derived via the combinatorial approach 
in \citet{sinclair2017spatial} and disagreement percolation methods of \citet{van1994disagreement}. In case of unbounded degree graphs, we also use recent sharp phase transition results from continuum percolation theory in \citet{ziesche2018sharpness}, along with disagreement percolation bounds to show that weak spatial mixing holds with high probability and we show this suffices for our central limit theorem. Spin systems are an example where the
full reach of \BL~cluster-localization is exploited and where stronger notions of stabilization do not apply. 

For interacting diffusions as in Example \ref{Ex:ID}, we use Doob's inequality and the It\^{o} isometry to establish a decorrelation bound with respect to the graph distance as in \cite{lacker2019local}. Using again the stabilizing properties of the interaction graph, we obtain $L^2$-stabilization for the scores or summands comprising $H_n$, which is stronger than \BL-localization.

In the case of interacting particle systems as in Examples \ref{Ex:CSA} and \ref{Ex:EP}, 
 we establish classical stopping set stabilization for functionals $\xi(\tx,\tP_n) := h(M(x,t, \tP_n)_{t \in [0,t_0]})$
 that depend on the time evolution of particle states, where $h$ is a suitable function of the state $M(x,t, \tP_n)$ at time $t$. This involves a careful analysis of a graphical construction that captures the system’s full dependency structure. For each site $x$, the history of a particle’s state is determined by a {\em backwards in time cluster}, defined via chains in an oriented space-time graph. The diameters of these clusters exhibit exponentially decaying tails, which yields the required  asymptotic decorrelation. While the stabilization used here is classical, the generality of the framework makes the results non-trivial to prove---even when site locations form a Poisson process or the initial states are i.i.d.

Our general theory aligns  naturally with empirical random fields or geostatistical marking as in Examples \ref{Ex:ERF} or \ref{Ex:GEO}, respectively, as the two point processes in those examples are independent. The approach is advantageous as it only requires \BL-localization of the random field $\zeta$ with respect to the underlying point process $\tP$ for the limit theory to hold.

\vskip.3cm

\subsection{Existing literature}
\label{ss.Literature}

Many  spatial random models have been investigated either in the mean-field regime or for Erd\H{o}s-R\'enyi-type random graphs;  we do not detail these contributions here, as there is little overlap and still less relevance. We shall focus instead  on overlap with existing central limit theorems for models on lattice-based or other discrete sparse graphs.

\paragraph*{General limit theorems for marked and unmarked point processes:}
Up to now research  on limit theorems for sums of general score functions of point processes has focused primarily on {\em scores with strong (stopping set) stabilization assumptions and input  with independent initial states}, i.e. the study of ~\eqref{linearstatsum} and~\eqref{linearstatB} with $U(x)$ being  independent marks (states) of points~$x$ (sites) given the ground process $\P$ and $\xi$ satisfying certain strong stabilization criteria. In the case of Poisson input 
one typically obtains cleaner  conclusions, often under   weaker 
stabilization criteria.  We recall these before mentioning results on general point processes.
  
When the collection of  sites  $\P$ is a
{\em Poisson point process} with {\em independent initial states} (marks),  the limit theory for the random measures  at ~\eqref{linearstatB}  is established  in 
\citet{Baryshnikov05},
\citet{Penrose2007gaussian},  and \citet{Penrose2001central,Penrose2002limit}.
The seminal work of
\citet{last2016normal} employs the Malliavin-Stein calculus to establish quantitative central limit theorems for general functionals on the space of  Poisson input.  Their work is extended 
and refined by \citet{Trauthwein}. Additional  works establishing quantitative central limit theorems for the functionals $(H_n^\xi)_{n \in \N}$   
on both marked and unmarked Poisson input include  \citet{Lachieze2019normal,lachieze2019shotnoise,SchulteYukich2019, LRPY, SchulteYukich2023}; these papers assume that the
scores $\xi$ satisfy stabilization criteria.

The  case of Gibbsian point processes is considered in 
\citet{Schreiber13}, \citet{XiaYukich} and \citet{hirsch2025normal},
whereas the case of point processes satisfying classical mixing conditions is handled in
\citet{Heinrich2013asymptotic,HLS,HP,Ivanoff82,JP}. In \citet{BYY19}, the limit theory for the measures in~\eqref{linearstatB}
is  further extended to include a ground process $\P$ (of sites)
having fast mixing  correlations; see also~\citet{Fenzl2019asymptotic} and~\citet{cong2022convergence} for independent (static) marks
given such a ground process ~$\P$. Recent work of
\citet{CRX2020} uses Stein couplings to obtain 
Berry-Esseen types of bounds on normal approximation of
statistics of general random measures, but not necessarily 
those expressed as sums of marks.

The  fast mixing  correlations of ground processes is exhibited  by point processes which are 
either stationary $\alpha$-determinantal or permanental point processes  having a fast decreasing kernel, the point process given by the zero set of a Gaussian entire function, certain rarified Gibbsian point processes or some Cox point processes (see \citet{BYY19}).
%In all these papers, the {\em static} state $h(U(x,0))$ is allowed to depend on the spatial locations of points in $\P$ near to $x$. 
The general results of \cite{BYY19} extend those of Soshnikov \cite{Soshnikov00} and Nazarov and Sodin \cite{Nazarov12}, who restrict attention to  $\xi \equiv c$, $c$ a constant, that is to say restrict to the case when  the random measures  at~\eqref{linearstatB} are counting statistics.  The  work of \citet{benevs2020decorrelation}  gives sufficient conditions for fast mixing  correlations of a Gibbs particle process (a point process on the space of compact sets of $\R^d$) and uses them to prove a central limit theorem for $U$-statistics of such processes.

\paragraph*{Spin Systems:} The above-mentioned central limit theorems and ours have precursors in central limit theorems for discrete random fields (i.e., fields indexed by graphs or lattices) which involve some notion of mixing or quasi-association. Such theorems have a long history and we refer the reader to 
\citet{dobrushin1977central,holley1979central,kunsch1982decay,bolthausen1982central,bulinski2013central,Reddy2018central} among others. In particular, these theorems were also applied for spin systems on fixed graphs such as lattices and again work in high temperature or low activity regimes. 
In comparison, our notion of \BL-mixing accounts for the randomness of the underlying graphs.

The specific form of weak spatial mixing which  we require is taken  from \citet{martinelli1994approach},  itself a modification of the well-known `Dobrushin-Shlosman uniqueness' condition \cite{dobrushin1985constructive};
such a condition is crucial when proving  uniqueness, concentration inequalities, and fast mixing of Glauber dynamics for spin systems, among other properties; see for example references in \cite{dobrushin1985constructive,van1993uniqueness,van1994disagreement,chatterjee2005concentration,sinclair2017spatial}.  Many existing central limit theorems use mixing in one form or another (see \cite[pg. 469]{georgii2011gibbs}), but as far as we understand, weak spatial mixing appears weaker than many of these notions.  There is considerable literature devoted to existence and phase-transitions in the infinite-volume limit for spin systems on lattices \cite{friedli2017statistical,duminil2015geometric} but we bypass these aspects here.

To the best of our knowledge, there is relatively little literature on the central limit theorem for spin systems on sparse random graphs, i.e., graphs with finite average degree. \citet{GGHP15,giardina2016annealed} and \citet{can2019annealed} have proved asymptotic normality for Ising models on locally tree-like sparse random graphs. They prove a random quenched central limit theorem and annealed 
central limit theorem for rather general models of sparse random graphs but for an averaged quenched central limit theorem, which is what we study, they require restrictive assumptions on the underlying graphs. Their proofs proceed by  explicitly evaluating asymptotics for the second derivative of the cumulant-generating function and also yields explicit variance asymptotics at volume-order scale. This proof strategy is similar to that of \citet[Section~V7]{ellis2007entropy}, who uses it to prove a central limit theorem on lattices. In contrast, our approach for spin systems uses coarser information on the graph geometry and applies to many graphs and spin models, though is subject to lower bounds on the growth of the variance.

An  approach via a quantitative Marcinkiewicz type of theorem was recently used in \citet{dinh2021quantitative} to prove normal approximation  for continuous spin systems with a non-trivial external field on lattices. This approach requires showing zero-freeness of the characteristic function of total spins which relies upon zero-freeness of the partition function via the Lee-Yang theory. Zero-freeness of the partition function is known to imply (or is equivalent to) strong spatial mixing in lattice models \cite{dobrushin1987completely} and some other examples \cite{gamarnik2023correlation,liu2025correlation,regts2023absence} and hence it is in general stronger than the notion of weak spatial mixing employed in our proofs.

\paragraph*{Interacting particle systems:}  
 A significant part of the literature, see \citet{liggett2012interacting},  is devoted to existence,  long-time behavior, and phase-transitions in particle systems on lattices or deterministic graphs.  Our investigation is  focused on the evolution in finite-time horizons, as  inspired by \citet{Penrose2008existence}.  The latter also considers existence of particle systems on the infinite Poisson point process, whereas  we consider only particle systems on large windows and hence sidestep the question of existence and long-term behavior.   
By restricting attention to large windows, this allows us to consider more general particle system models than those in \cite{Penrose2008existence} or \cite{liggett2012interacting}.  Understanding the existence in the infinite volume limit as well as understanding the long-time behavior of our models on general point processes is an intriguing question lying beyond the scope of this work, though the techniques of \cite{Penrose2008existence,mazumder2024existence} combined with some of our derivations could be helpful in this direction.  

Central limit theorems for interacting particle systems on lattices and transitive graphs have a long history, going back to  \citet{holley1979central,holley1980generalized}, \citet{brox1984equilibrium}, \citet{rost1985central} and,  more recently \citet{doukhan2008functional}.  These works assume finite range interactions apart from other particular assumptions on either the state space, the initial states or Markovian dynamics. Specific models in the continuum, including random sequential adsorption, are treated  in \citet{PY2001},  but we are not aware of central limit theorems for general continuum interacting particle systems  beyond those of \cite{Penrose2008existence}. The papers \citet{qi2008functional} and \citet{onaran2022functional,onaran2023functional}  give functional central limit theorems for statistics of certain spatial birth death processes (and the latter incorporates diffusive dynamics as well), albeit under different assumptions,  whereas \cite{Penrose2008existence} establishes a functional central limit theorem for more general time-indexed interacting particle systems,  subject to finite range conditions and certain assumptions on the birth and death rates. For a more detailed comparison of proof techniques, see Remark \ref{rem:IPS}. The underlying point process in these works is the Poisson process and we are not aware of any results on general point processes.

Concerning  our other applications such as interacting diffusions or geostatistical marking, we are not aware of existing central limit theorems. Central limit theorems for empirical random fields under strong mixing conditions are given in \cite{Pa09}.

\vskip.3cm

\subsection{Organization of the work}
\label{ss.organization}

After introducing some notation and definitions in the next section (Section \ref{s:notation}), the rest of the work is divided into two parts. Part \ref{part:theory-foundations} lays out the general framework of mixing correlations of point processes and localization of scores, the main theoretical results on expectation, variance asymptotics and normal convergence as well as their proofs. We introduce the twin notions  of \BL-mixing correlations and \BL-localizing score functions and use these key concepts to formulate our main limit theorems. 
 The applications of our general theorems to the aforementioned spatial random models are not immediate and we reserve the entirety  of Part \ref{part:applications} to separately and thoroughly investigate each model. While these applications can be read independently of each other, 
the interactions in these models are commonly described via a stabilizing interaction graph. 

\section{Notation, terminology and definitions}
\label{s:notation}
We frame our approach in the setting of marked point processes. Relevant notation is introduced below, and additional background---including key definitions from Palm theory---is provided in the following subsection. For more on point processes and random measures, we refer the reader to \cite{BBK2020, Kallenberg2017random, last2017lectures}.
Section \ref{s:prelim} lists  the commonly used notation whereas  Section \ref{s:appppprelim} elaborates upon the definitions and terminology.

\subsection{List of notation}
\label{s:prelim}

\begin{itemize}

\item  $\K, \K',  \M$ and $\mathbb{L}$ are Polish  spaces, i.e. separable completely metrizable spaces.

\item  For integers $1\le a<b$, denote by $\bk zab$ the
 vector $(z_a,\ldots,z_b)$ with elements $z_i$ in $\R^d$, $\K$,
 $\N^0:=\{0,1,\ldots\}$, or $\R^d\times\K$, depending on the context.
 Denote the concatenation of such tuples by
 $\bk zab\sqcup\bk {z'}cd:=(z_a,\ldots,z_b,\allowbreak z'_c,\ldots,z'_d)$. 
 We also write  $[k]:=\{1,2,\ldots,k\}$ for $k\in\N$.

\item We consider  a compatible metric $d_\K$ on the Polish space $\K$ and we consider the {\em 1-product metric} $d_{\K^p}(\bk u1p,\bk v1p)=
\sum_{i=1}^p d_{\K}(u_i,v_i)$ for $\bk u1p,\bk v1p\in\K^p$. The same applies to other Polish spaces. 

\item For a vector $\bk x1p\in(\R^d)^p$, we abbreviate $B_r(\bk x1p) = \bigcup_{i=1}^pB_r(x_i)$, where $B_r(x)$ is a ball of radius $r>0$ centered at $x\in\R^d$. 

\item ${\cal B}_b := {\cal B}_b(\R^d)$---the class of bounded Borel subsets of~$\R^d$.

\item  ${\B}(\K^p)$---the class of real-valued, measurable functions $f$ on $\K^p$ with $|f|\le 1, p \in \N$.

\item $\BL(\K^p)$---the class of Lipschitz($1$) functions $f : \K^p \to \R$
with supremum norm bounded by $1$. This is the class of $f : \K^p \to \R$
such that 
$$
\sup_{\bk x1p,\bk y1p \in \K^p} \frac{ |f(\bk x1p) - f(\bk y1p)|}{d_{\K^p}(\bk x1p,\bk y1p)} \leq 1, \quad \sup_{\bk x1p\in\K^p} |f(\bk x1p) | \leq 1,
$$
where $d_{\K^p}$ is the $1$-product metric on $\K^p$ with respect to the given  metric $d_{\K}$ on  $\K$ (i.e., $d_{\K^p}(\bk x1p,\allowbreak\bk y1p)=\sum_{i=1}^p d_{\K}(x_i,y_i)$) 
These functions are sometimes referred to as bounded  (or double-bounded) Lipschitz functions, where both the Lipschitz constant and the function values are bounded.
The underlying domain $\K$ will be omitted whenever it is  implicit. 

\item Given $\K^p$-valued random vectors $\bk X1p$ and $\bk Y1p$, defined on possibly different probability spaces, the 
{\em bounded Lipschitz distance} between the laws of $\bk X1p$ and $\bk Y1p$  is
$$
d_{\BL}({\cal L}(\bk X1p), {\cal L}(\bk Y1p))= \sup_{f \in \BL(\K^p)} | \E f(\bk X1p) - \E f(\bk Y1p)|.
$$
Abusing notation we often write $d_{\BL}(\bk X1p,\bk Y1p)$ for $d_{\BL}({\cal L}(\bk X1p), {\cal L}(\bk Y1p))$.
Since \(\K\) (and hence \(\K^p\)) is a separable metric space, the bounded Lipschitz distance \(d_{\BL}\) metrizes convergence of probability measures on \(\K^p\) (see \citet[Theorem 8.3.2]{bogachev2007measure}). Moreover, if \(d_{\K}\) (and hence \(d_{\K^p}\)) is complete, then this metrization is also complete (see \cite[Theorem 8.10.43]{bogachev2007measure}). The bounded Lipschitz distance \(d_{\BL}\) is also equivalent to the Fortet–Mourier distance (see \citet{villani2009optimal}).

\item Given $\K$-valued random variables $X,Y$, defined on possibly different probability spaces, the {\em total variation distance} between the laws of $X$ and $Y$ is
$$
d_{\TV}({\cal L}(X), {\cal L}(Y)) := \frac{1}{2} \sup_{f \in \B(\K)} | \E f(X) - \E f(Y)|.
$$
When $\K$ is finite, we have 
$$d_{\TV}({\cal L}(X), {\cal L}(Y)) = \frac{1}{2}\sum_{k \in \K} |  \Palm(X = k) -  \Palm(Y=k)|.$$
\item The collection of random sites $\P$ is a simple point process on $\R^d$.

\item We let $\tilde \P = \sum_{x \in \P} \delta_{ \tx }$ be  a {\em
  simple, marked point process} on $\R^d \times \K$,  where $\tx :=
  (x, U(x))$  are random elements in  $\R^d \times \K$ and for any
  $x \in \R^d$, the set $\{x\} \times \K$ contains at most one
  element of $\tP$. Elements $x\in\R^d$ are  interchangeably referred to as {\em points} or {\em sites}, while $U(x)$ are  referred to as
  {\em pre-marks} or {\em initial states}. %\jy{Write $U(x):= U(x, \P, \Xi)$ to include 
  %geostatistical marks, random field examples?} 
  The projection of $\tilde \P$ on
  $\R^d$ is the collection of sites $\P$,  also called  the {\em
    ground process}.
\item  By $\P_n := \P \cap W_n, n \in \N,$ we denote the
    restriction of the ground process $\P$ to $W_n:=
    [-\frac{1}{2} n^{1/d}, \frac{1}{2} n^{1/d}]^d$, 
   and,  with a slight abuse of notation,  $\tP_n:=\tP\cap W_n$
  means the restriction  of the marked point process $\tP$ whose ground
  points belong to $W_n$. In certain contexts, we write
  $\P_\infty:=\P$ and $\tP_\infty=\tP$.

\item Let $\cN_{\R^d}$ and $\cNK$ denote {\em the spaces of locally finite subsets} of $\R^d$ and $\R^d \times \K$,  respectively, equipped with suitable topologies and $\sigma$-algebras.  Let $\hat{\cN}_{\R^d}$ and $\hat{\cN}_{\R^d \times \K}$ be {\em the space of finite subsets} of  $\R^d$ and $\R^d \times \K$ respectively.   $\hat{\cN}_{[0,t_0] \times \K  }$ denotes the collection of finite subsets
 of $[0,t_0] \times \K  $.
 
\item $\rho^{(p)}(\bk x1p)$---the  {\em $p$\,th order correlation function}  of the ground process $\P$;
for a simple point process \( \P \) it vanishes whenever \( x_i = x_j \) for any \( i \neq j \).

\item $\tilde \alpha^{(p)}(\md \bk {\tx}1p)$---the {\em $p$\,th order factorial moment measure} of $\tilde\P$,  where  $\bk {\tx}1p=\bk{x,u}1p \in(\R^d\times\K)^p$.   
 
 \item $\cM_{\bk x1p}(\md \bk u1p)$---the  {\em Palm probability
        distribution} on $\K^p$ of the random marks $\bk U1p=(U_1,\ldots,U_p) = (U(x_1),...,U(x_p))$
 respectively at the  fixed locations $\bk x1p = (x_1,...,x_p)$.
    
\item $\Palm_{\bk{\tx}1p},  \E_{\bk{\tx}1p}$---the {\em  Palm probability distribution} and  the corresponding Palm expectation of the marked point process $\tilde\P$  given points $\bk x1p$ and their marks $\bk u1p$ where  $\bk {\tx}1p=\bk{x,u}1p \in(\R^d\times\K)^p$. 
   
\item $\Palm^{!}_{\bk{\tx}1p},  \sE^{!}_{\bk{\tx}1l}$---the  {\em reduced Palm versions} 
  and the corresponding expectation are such that  the conditioning points $\bk x1p$ and marks  $\bk u1p$ are removed from  $\tP$.  
  
\item $\E_{\bk
x1p}[...]:=\int_{\K^p}\E_{\bk{\tx}1p}[...]\,
\cM_{\bk x1p}(\md \bk u1p)$---the {\em Palm expectation of
$\tilde\P$ given points $\bk x1p$} without fixing their marks
$\bk u1p$.  

\item  $\tilde{\alpha}^{(l)}_{\bk{\tx}1p}(\md \bk{\ty}1l )$---the {\em Palm  factorial moment measures}   of  $\tP$ of order $l$ on  $(\R\times\K)^l$ under  $\Palm^!_{\bk {\tx}1p}$, for fixed $\bk{\tx}1p$.

\item $\rho^{(p)}_{\bk x1p}$---the $p\,$th order  correlation functions of the ground process $\P$ under  its reduced Palm versions 
${\sP}^!_{\bk x1p}$ (without fixing
marks).
\item $G(\X)$ and ${\cal G}(\mu)$  denote graphs on the vertex set $\X$ and the  counting measure $\mu$ on~$\R^d$, respectively. The graph metric is denoted by $d_{G(\X)}$ or  $d_{\G(\mu)}$ and the graph boundary for $\X' \subset \X$ is denoted by $\partial \X'$. 

\item $\gB_{k}(\X')$ denotes the $k$-ball in the graph distance around $\X' \subset \X$.

\item $G^{[k]}_{\X'}(\X) := G(\gB_k(\X'))$ denotes the induced sub-graph of $G(\X)$ on $\gB_{k}(\X')$,
$\X' \subset \X$. Similar notation is used for ${\cal G}(\mu)$.

\item For $\X \subset \R^d$, $\G^{(r)}_{\X'}(\X) := \G(\X \cap B_r(\X'))$  denotes  the induced sub-graph of $\G(\X)$ on $\X \cap B_r(\X')$, the Euclidean $r$-ball around $\X' \subset \R^d$.

\item $Z$ denotes a standard normal random variable.

\item A function $\phi: [0, \infty) \to [0,\infty)$ is said to be {\em fast decreasing} if 
\be \label{varphibd}
\limsup_{r \to \infty} r^m \phi(r) \to 0,  \,  \,  \forall m \in \N.  
\ee
Further, $\phi$ is  {\em exponentially decreasing} if there exists $b \in (0, \infty)$ such that \\ $\limsup_{r \to \infty} r^{-b}\log \phi(r) < 0$.

\item For $k \in \N$ let  $\omega_k :  [0, \infty) \to [0, \infty)$. 
The collection $(\omega_k)_{k \in \N}$ of functions is said to be an {\em increasing-decreasing family} if for all $k \in \N$ we have $ \omega_k \leq \omega_{k + 1}$  and $\omega_k(s) \downarrow 0$ as $s \to \infty$. 
 \end{itemize}

\subsection{Formal definitions: \  marked point processes}
\label{s:appppprelim}
\begin{itemize}
\item The collection of random sites $\P$ is a simple point process on $\R^d$; it is a random element in the space $\cN_{\R^d}$ of {\em locally finite subsets} of $\R^d $ endowed with the evaluation $\sigma$-algebra, i.e., the $\sigma$-algebra generated by the maps $\mu \mapsto \mu(B)$ for all Borel subsets $B$.  Let $\hat{\cN}_{\R^d}$ denote the subspace of $\cN_{\R^d}$ consisting of all {\em finite subsets} of $\R^d$.  We note that $\cN_{\R^d}$ is a Polish space under the vague topology 
 (\cite[Lemma 4.6]{Kallenberg2017random}) and that   $\hat{\cN}_{\R^d}$ is likewise a Polish space under the vague topology, which is equivalent to the weak  topology \cite[Chapter 4]{Kallenberg2017random}.

\item We let $\tilde \P = \sum_{x \in \P} \delta_{ \tx }$ be  a {\em simple, marked point process}, where $\tx := (x, U(x))$, with points  $x \in\R^d$ and where the marks $U(x)$ 
are random elements in $\K$, a Polish space.  The $x$ are random points and we note that
$\tx$ are random elements in  $\R^d \times \K$ and that
      $\tilde \P$ is a point process on $\R^d \times \K$; i.e., it is  a
      random element in the space $\cNK$ of locally finite subsets of
      $\R^d \times \K$ (interpreted as counting measures) and  endowed with the evaluation $\sigma$-algebra, with any set $\{x\} \times \K$, $x \in \R^d$, containing at most one element of $\tP$.  Recall that the $U(x)$ are interchangeably referred to as marks or initial states. 
      The projection of $\tilde \P$ on $\R^d$ is the collection of
      sites $\P$, also called  the {\em ground process}. As above,  we
 take   $\hat{\cN}_{\R^d \times \K}$ to be the space of {\em finite} subsets
      of $\R^d \times \K$ and this again is a Polish space under the
weak topology, while  ${\cN}_{\R^d \times \K}$ is also  a Polish space under the vague  topology \cite[resp. Lemma 4.5  and Lemma 4.6]{Kallenberg2017random}.
 
 \item For $p\in\N$, let \( \P^{(p)} \) be the \( p \)-th factorial power of the simple point process \( \P \), i.e., the collection of tuples \( (x_1, \dots, x_p) \) of distinct points of \( \P \). It forms a point process on \( (\R^d)^p \) with the mean  measure \( \E[\P^{(p)}(\cdot)] \),   known as the {\em factorial moment measure} of order \( p \) of \(\P\). If the factorial moment measure of order $p$ is locally finite, then  the \( p \)-th order correlation function \( \rho^{(p)} \) is the  Radon-Nikodym density (if it exists) of this factorial moment measure with respect to the Lebesgue measure on \( (\R^d)^p \). For \( \P \) simple this ensures \( \rho^{(p)}(\bk x1p) = 0 \) when \( x_i = x_j \) for any \( i \neq j \).
 For any non-negative measurable function $f:(\R^d)^{p} \to [0, \infty)$,   we have
$$
\E[\sum_{(x_1,\ldots,x_p) \in \P^{(p)}}f(\bk x1p)] =
 \int_{(\R^d)^{p}}f(\bk x1p) \rho^{(p)}(\bk x1p) \md (\bk x1p).
 $$
We slightly abuse notation as the  term $ \bk x1p$ on the left-hand side is a random $p$-tuple whereas the term 
$\bk x1p$ on the right is deterministic.  Unless the context is specifically restricted, we assume throughout that the
correlation functions $\rho^{(p)}$, $p \in \N$, exist and are bounded.   We also write $\rho(x):=\rho^{(1)}(x), x \in \R^d$.  
  
\item Extending this approach, we define the factorial power \( \tP^{(p)} \) of order \( p \) for the simple marked point process \( \tP \) as the collection of \( p \)-tuples \( \bk{\tx}1p = \bk{x,u}1p \) of distinct marked points of \( \tP \). Its mean measure,  
\(
\tilde \alpha^{(p)}(\cdot) := \E[\tP^{(p)}(\cdot)],
\)
is the {\em \( p \)-th order factorial moment measure} of \( \tP \), considered on $(\R^d\times \K)^p$. The simplicity of the marked point process (i.e., $\tP(\{x\}\times \K)\le 1$)
implies that the factorial moment measures of $\tP$ satisfy the same finiteness
conditions as those of the ground process~$\P$.
Consequently, under the same moment assumptions as for $\P$,  for any non-negative measurable function $f :(\R^d \times \K)^{p} \to [0, \infty)$,  we have 
 $$
  \E[\sum_{(\tx_1,\ldots,\tx_p) \in \tP^{(p)}}f(\bk \tx1p)] =
     \E[\sum_{(\tx_1,\ldots,\tx_p) \in \tP^{(p)}}f(\bk{x,U}1p)]
     = \int_{(\R^d \times \K)^{p}}f(\bk{\tx} 1p) \tilde \alpha^{(p)}(\md \bk{\tx}1p),
$$
where in the middle expression we adopt the convention that $U_i = U(x_i)$ for 
all $i \in \N$, $\bk{x,U}1p=\bk{(x,U)}1p=((x_1,U_i),\ldots,\allowbreak (x_p,U_p))$
and where again on the left-hand side $ \bk {\tx}1p=\bk{x,U}1p$ is a random $p$-tuple of distinct marked
points, whereas the term $\bk {\tx}1p=\bk{x,u}1p$ on the right is deterministic.  
 
 \item For a fixed vector of distinct locations \( \bk x1p = (x_1, \dots, x_p) \in (\R^d)^p \), the {\em (marks-)Palm probability distribution} \( \cM_{\bk x1p}^{(p)}(\md \bk u1p) \) on \( \K^p \) describes the distribution of the random marks \( \bk U1p = (U_1, \dots, U_p) \) at these locations. It is defined by  
\[
\tilde \alpha^{(p)}(\md \bk{\tx}1p) = \cM_{\bk x1p}^{(p)}(\md \bk u1p) \rho^{(p)}(\bk x1p)\, \md \bk x1p,
\]
where \( \tilde \alpha^{(p)}\) is  the \( p\, \)th order factorial moment measure of \( \tilde\P \), and \( \bk{\tx}1p = \bk{x,u}1p = ((x_1, u_1), \dots, \allowbreak (x_p, u_p)) \).  
 This disintegration holds since the boundedness of the correlation function \( \rho^{(p)} \) ensures the \( \sigma \)-finiteness of \( \tilde \alpha^{(p)} \) on the space \( \K \) that is Polish; see \cite[Theorem 14.D.10]{BBK2020}.
Thus the Palm probability distribution \( \cM_{\bk x1p}^{(p)}(\md \bk u1p) \) exists for almost every \( \bk x1p \) and is only meaningful for distinct locations \( x_i \).
    
\item For fixed  distinct marked locations $\bk {\tx}1p=\bk{x,u}1p \in(\R^d\times\K)^p$, 
 we denote by $\Palm_{\bk{\tx}1p}$ the corresponding {\em  Palm probability distribution} of the (entire) process~$\tilde\P$  given points   $\bk x1p$ and their marks $\bk u1p$.  These are  again defined  via 
a disintegration of the higher order
Campbell's measure $C^{(p)}(\cdot\times\cdot)$ on $(\R^d\times\K)^p\times\cN_{\R^d \times \K}$ given by
$$
C^{(p)}(B\times L)
:= C^{(p)}_{\tilde\P} (B\times L)
:=\E\Bigl[\int_{(\R^d\times\K)^p}\1{\bk{\tx}1p\in
    B}\1{\tilde\P \in L}\,  \tilde\P^{(p)}(\md \bk{\tx}1p) \Bigr]\, 
$$
with respect to  the $p\,$th order factorial moment  measure $\tilde \alpha^{(p)}$ of $\tilde\P$ (from which $C^{(p)}$  inherits   $\sigma$-finiteness):
\begin{equation}
\label{e:palm_existence}
C^{(p)}(\md(\bk{\tx}1p, \tilde\mu))= \Palm_{\bk{\tx}1p} (\md \tilde\mu)
  \,\tilde \alpha^{(p)}(\md \bk{\tx}1p).
\end{equation}
The crucial existence of the  disintegration in~\eqref{e:palm_existence}
follows again from ~\cite[{ Theorem 14.D.10}]{BBK2020}  and  the
observation that $\cN_{\R^d \times \K}$ is Polish.  By $\E_{\bk{\tx}1p}$,   we mean the corresponding
expectation with respect to the Palm probability distribution ~$\Palm_{\bk{\tx}1p}$.  These again are defined only for $\tilde \alpha^{(p)}$ a.e. $\bk{\tx}1p$ and only meaningful for distinct marked locations. 

From the above definition,  it follows that the Palm expectation satisfies the following identity known as the Campbell-Little-Mecke formula.  For a measurable $f :  (\R^d\times\K)^p\times \cN_{\R^d \times \K} \to [0, \infty)$,  we have
\begin{equation} \label{CLM}
\E[\sum_{(\tx_1,\ldots,\tx_p) \in \tP^{(p)}}f(\bk \tx1p,\tP)] =  \int_{(\R^d\times\K)^p} \E_{\bk{\tx}1p}[f(\bk{\tx}1p,\tP)] \tilde \alpha^{(p)}(\md \bk{x,u}1p).
\end{equation}

\item The  {\em reduced Palm versions} $\Palm^{!}_{\bk{\tx}1p}$
  and their expectations $\sE^{!}_{\bk{\tx}1l}$  are such that  the conditioning points $\bk x1p$ and marks  $\bk u1p$ are removed from  $\P$.  These versions, which may be formally defined  via
  $$ \Palm^{!}_{\bk{\tx}1p} \left(\tP \in \cdot\right) = \Palm_{\bk{\tx}1p} \left(\tP - \sum_{i=1}^p \delta_{\tx_i} \in \cdot\right),
  $$
 satisfy the following modified  Campbell-Little-Mecke formula.  For a measurable $f :  (\R^d\times\K)^p\times\cN_{\R^d \times \K} \to [0, \infty)$,  we have 
\begin{equation}
\label{e:CLM}
\E[\sum_{(\tx_1,\ldots,\tx_p) \in \tP^{(p)}}f(\bk \tx1p,\tP-\sum_{i=1}^p \delta_{\tx_i})] =  \int_{(\R^d\times\K)^p} \E^{!}_{\bk{\tx}1p}[f(\bk{\tx}1p,\tP)] \tilde \alpha^{(p)}(\md \bk{\tx}1p).
\end{equation}

\item Let $\E_{\bk
x1p}[...]:=\int_{\K^p}\E_{\bk{\tx}1p}[...]\,
\cM_{\bk x1p}^{(p)}(\md \bk u1p)$ denote the Palm expectation of
$\tilde\P$ given points $\bk x1p$ without fixing their marks
$\bk u1p$.  Under this notation,  the Campbell-Little-Mecke formula can be rewritten as 
\begin{equation}
\label{e:CLM2}
\E[\sum_{(\tx_1,\ldots,\tx_p) \in \tP^{(p)}}f(\bk \tx1p,\tP)]  =  \int_{(\R^d)^p} \E_{\bk x1p}[f(\bk \tx1p,\tP)] \rho^{(p)}(\bk x1p)   \md \bk x1p. 
\end{equation}  
We recall and clarify that, when using $\Palm_{\bk{x}{1}{p}}$ or
$\E_{\bk{x}{1}{p}}$, we tacitly assume that the factorial moment measure
$\alpha^{(p)}$ of the ground process~$\P$ is locally finite. In particular,
this holds when its density $\rho^{(p)}$ exists (at least locally at $\bk{x}{1}{p}$).

\item The {\em Palm  factorial moment measures}   of order $l$ (on
  $(\R\times\K)^l$) of $\tP$ under  $\Palm^!_{\bk {\tx}1p}$ (for fixed $\bk{\tx}1p$)
  are  denoted by
  $\tilde{\alpha}^{(l)}_{\bk{\tx}1p}(\md \bk{\ty}1l )$. 
  The existence of {\em Palm correlation functions} is guaranteed by the following Palm algebra relation~\cite[Proposition 3.3.9]{BBK2020}
  $$\tilde{\alpha}^{(p+l)}(\md(\bk{\tx}1p\sqcup\bk{\ty}1l )) =
  \tilde{\alpha}^{(l)}_{\bk{\tx}1p}
  (\md\bk{\ty}1l)\tilde{\alpha}^{(p)}(\md
  \bk{\tx}1p), $$
for $\bk{\tx}1p\in (\R\times\K)^p$, $\bk{\ty}1l\in (\R\times\K)^l$.
 The above relation can also be rewritten in terms of correlation functions and mark distributions as 
 \begin{equation}
\label{e:palmalgebra}
\cM_{\bk x1p \sqcup \bk y1l}^{(p+l)}(\md(\bk u1p \sqcup \bk w1l)) \rho^{(p+l)}(\bk x1p \sqcup \bk y1l)  =  \tilde{\alpha}^{(l)}_{\bk {x,u}1p}(\md\bk {y,w}1l) \cM_{\bk x1p}^{(p)}(\md \bk u1p) \rho^{(p)}(\bk x1p),
 \end{equation}
where in $\cM^{(p+l)}_{\bk x1p\sqcup\bk y1l}(\md (\bk u1p\sqcup\bk
  w1l))$ we understand $\bk u1p, \bk w1l$ to denote the markings
  associated to $\bk x1p,\bk y1l$ respectively.
  Likewise, the  Palm distributions   $(\Palm^!_{\bk
      {\tx}1p})^!_{\bk{\ty}1l}$ of $\Palm^!_{\bk
      {\tx}1p}$ exist and satisfy  
      $(\Palm^!_{\bk {\tx}1p})^!_{\bk{\ty}1l}=\Palm^!_{\bk {\tx}1p
   \sqcup\bk {\ty}1l}$ for almost all $\bk {\tx}1p
   \sqcup\bk {\ty}1l$ with respect to
   $\tilde{\alpha}^{(p+l)}(\md(\bk{\tx}1p\sqcup\bk{\ty}1l))$. 

\item The {\em $l$\,th correlation functions $\rho^{(l)}_{\bk x1p}(\bk y1l)$ of the ground process $\P$ under  its reduced Palm versions ${\sP}^!_{\bk x1p}$ (without fixing marks)}.
Formally, these functions  are defined as the Radon-Nikodyn density of the measure
 $\sE^!_{\bk x1p}[\P^{(l)}(\cdot)]$ with respect to the Lebesgue measure. Dropping marks in the Palm algebra relation~\eqref{e:palmalgebra} we have 
\begin{equation}
\label{e:palmcorrprod}
\rho^{(p+l)}(\bk x1p \sqcup \bk y1l)  =  \rho^{(l)}_{\bk x1p}(\bk y1l) \rho^{(p)}(\bk x1p)
\end{equation}
\end{itemize}

\part{Theory and Foundations}
\label{part:theory-foundations}
This part systematically develops the  theoretical framework underpinning this work. In Section~\ref{s.decay-correlation}, we introduce the key concept of  asymptotic decorrelation of marked  point processes $\tP$ (not necessarily stationary) as well as the notion of \BL-mixing marked point processes.  We state a general central limit theorem for such point processes, one which implies most of our Gaussian fluctuation results.   Section~\ref{s.Stabilizing} provides a general way to describe and construct asymptotically decorrelated point processes having dependent marks via jointly \BL-localizing score functions.  Section~\ref{s:clta} establishes  the general limit theory for statistics of these point processes.
This includes central limit theorems, and,  under the additional assumption of stationarity and translation invariance of $\xi$, it provides 
 the first and second order limit theory for $( \mu_n^\xi(f))_{n \in \N}$. 
The proof of the key mixing correlation bound from Section \ref{s.Stabilizing} is in Section~\ref{s.proof-decay} whereas  the proofs of the general limit theorems are  in Section \ref{s.proof-limit}.

We make the following standing assumption. 
\begin{customass}{3.1} 
\label{Ass1} \, 

\begin{enumerate}[wide,label=(\roman*),labelindent=0pt]

\item $\P$ denotes a simple point process on $\R^d$ such that  whenever for some  $p \in \N$ the correlation function~$\rho^{(p)}$ exists, it is uniformly bounded, i.e., $\kappa_p:=\sup_{\bk x1p\in(\R^d)^p}  \rho^{(p)}(\bk x1p) < \infty$.   
We put $\kappa_0 :=  \max\{ \kappa_1, 1 \}$.

\item  Marked point processes $\tilde{\P}$ and
 $(\tilde\P_n)_{n \in \N}$ are assumed to be simple.
\end{enumerate}
\end{customass}

\section{Mixing correlations of marked point processes}
\label{s.decay-correlation}

We define the main notion of asymptotic decorrelation of marked point processes, called `mixing correlations', and use it to formulate our fundamental result (Theorem~\ref{t:clt_linear_marks}) concerning Gaussian fluctuations of some macroscopic characteristics of a large class of random spatial models. Next, in Section~\ref{ss.Mixing-separately} we elaborate upon the properties of mixing correlations,
and in particular show that they can be deduced from mixing correlations of points together with mixing correlations of marks. This is  useful when proving that the main mixing hypotheses follow from more elementary assumptions.

\subsection{Mixing correlations and asymptotic normality}
\label{ss:asnormmpp0}

We formulate the main result of this section, a fundamental tool for studying the Gaussian fluctuations of spatial models.
Following the first lines  of the  Introduction, consider a {\em ground point process} $\P=\{x_i\}$ on  $\R^d$,  and denote $\P_n := \P \cap W_n, n \in \N,$ where $W_n := [-\frac{1}{2} n^{1/d}, \frac{1}{2} n^{1/d}]^d$. Consider a \emph{family of real-valued marked point processes}
$(\tilde\P_n)_{n \in \N}
=\bigl(\{(x_i,\xi_{i,n})\}_{x_i \in \P_n}\bigr)_{n \in \N}$,
where the $\xi_{i,n}$ are real-valued random variables, commonly referred
to in the literature as \emph{scores}.
 We are interested in establishing  a central limit theorem for the random variables $\mu_{n}^{\xi}(f):=  \int f d \mu_{n}^{\xi}, n \in \N$,
where  
\begin{equation}\label{e:signed-measure}
    \mu_n^\xi: = \sum_{x_i\in \P_n} \xi_{i,n} \delta_{ n^{-1/d}x_i }
\end{equation}
is a (possibly signed) measure on $W_1$ 
%defined in~\eqref{e:signed-measure}
and where $f\in\B(W_1)$ is a bounded test function on $W_1$.

If $(\tP_n)_{n \in \N}$ are the respective restrictions to $(W_n)_{n \in \N}$
of a stationary ground point process~$\P$ equipped with i.i.d.\ scores,
that is, $\{\xi_{i,n}\}_{x_i \in \P_n}$ are i.i.d.\ given~$\P_n$ for each
$n \in \N$, then, under suitable moment conditions on $\{\xi_{i,n}\}$, as noted
in~\cite{Fenzl2019asymptotic}, straightforward modifications of~\cite{BYY19}
yield the asymptotic normality of $(\mu_{n}^{\xi}(f))_{n \in \N}$, for
$f \in \B(W_1)$.
Rather than imposing mixing conditions on $\P$ via its underlying
$\sigma$-algebra, this approach---motivated by the notion of
`clustering of correlation functions' in statistical physics---relies on
the factorization (or mixing) of correlation functions of~$\P$; see
Definition~\ref{d.omegamixing}.

However, this approach does not readily extend to the case of dependent scores,
which commonly arise in models of interest. Indeed, viewing a marked point
process $\tP_n$ with real-valued scores $\{\xi_{i,n}\}_{x_i \in \P_n}$ as a
(non-marked) point process in $\R^d \times \R$ for the purpose of analyzing its
correlation functions is often inconvenient, especially when the scores depend
on additional (possibly abstract) sources of randomness. If one wishes to avoid
delicate considerations involving the $\sigma$-algebra of $\tP_n$, an approach
adopted in this work is to define mixing correlations directly for abstractly
marked point processes, as follows.

By $(\tilde\P_l = \{(x,U_l(x))\}_{x \in \P})_{l \in \N}$ we mean a
\emph{family} of marked point processes sharing the same ground point process
$\P$ in $\R^d$, with marks $U_l(x) \in \K$, $x \in \P$, taking values in a
common mark space~$\K$. These marks may depend on the configuration~$\P$
and/or on auxiliary randomness, all defined on a common probability space on
which the processes $(\tilde\P_l)_{l \in \N}$ are constructed.
\begin{defn}[$\BL$-mixing correlations of marked point
    processes] 
\label{d.omegahmixing} \,
\begin{enumerate}[wide,label=(\roman*),labelindent=0pt]
\item \label{i.mixing} A marked point process $\tilde \P = \{(x, U(x))\}_{x \in \P}$ has  {\em $\BL$-mixing
  correlations} if there exists a family $( \omega_k)_{k \in \N}$ of  functions $\omega_k : [0, \infty) \to [0, \infty)$ with $\lim_{r \to \infty}\omega_k(r) = 0$ for all $k \in \N$ and 
   such that for all $p,q \in \mathbb{N}$, $x_1,\ldots,x_{p+q} \in \mR^d$, and all $f \in \BL(\K^p), g \in \BL(\K^q)$, 
    we have
\begin{align} \nonumber
& \left|\sE_{\bk x1{p+q}}[f(\bk U1{p})g(\bk U{p+1}{p+q})]\rho^{(p+q)}(\bk x1{p+q}) - \sE_{\bk x1p}[f(\bk U1p)] \rho^{(p)}(\bk x1{p})\sE_{\bk x{p+1}{p+q}}[g(\bk U{p+1}{p+q})]\rho^{(q)}(\bk
x{p+1}{p+q})\right|\\
&\hspace{0.1\linewidth}\le \omega_{p+q}(s),\label{e.mppmixing}
\end{align}
where
\be
\label{defs}
s := d(\bk x1p,\bk x{p+1}{p+q}):= \inf_{i \in \{1,...,p\}, j \in \{p + 1,...,p + q\} } |x_i - x_j|.
\ee
We refer to $(\omega_k)_{k \in \N}$  as the correlation decay functions for $\tilde \P$ and without loss of generality, we assume that the  collection $(\omega_k)_{k \in \N}$ is  an {\em increasing-decreasing family} of functions, by which we mean that for all $k \in \N$ we have $ \omega_k \leq \omega_{k + 1}$  and $\omega_k(s) \downarrow 0$ as $s \to \infty$.

\item  \label{i.fast-mixing} $\tilde\P$ has {\em 
 fast $\BL$-mixing correlations}  if in addition to  \eqref{e.mppmixing}, 
   for all $k\in \N$,  the correlation decay function $\omega_k$ is  fast decreasing, in accordance with property~\eqref{varphibd}. 
 \item \label{i.family-mixing} A {\em family of marked point processes}   $(\tilde\P_l)_{l \in \N},$ 
 has  $\BL$-mixing correlations,
 if there exists a family $( \omega_k)_{k\in\N}$ of increasing-decreasing functions $\omega_k :  [0, \infty) \to [0, \infty)$ 
  and such that  \eqref{e.mppmixing} holds uniformly in
 $l \in \N$ with respect to the Palm distributions  of marks $\sE_{\bk  x1{p+q}}
[f(\bk{U_l}1{p})g(\bk {U_l}{p+1}{p+q})]$ of $\tP_l$.
  The family  $(\tilde\P_l)_{l \in \N}$ has 
  {\em  fast} $\BL$-mixing correlations, if
  the increasing-decreasing family of correlation decay functions $(\omega_k)_{k \in \N}$
 is fast decreasing.
 \end{enumerate}
\end{defn}

Examples of marked point processes $\tP$ having \BL-mixing correlations include 
point processes $\P$ with mixing correlations and whose 
marks are generated either by spins,  states of interacting diffusions or  interacting particle systems on 
stabilizing interaction graphs on $\P$, as detailed in Part \ref{part:applications}. The upcoming 
Proposition \ref{p:pmix+mmix}, respectively Corollary \ref{c:mppmixcl_double},  provides general conditions under which a marked point process $\tP$, respectively family of marked point processes, has \BL-mixing correlations.

\medskip
We have  introduced mixing correlations for abstract marks, a concept which figures  prominently  throughout the work. We now return to scores (which are real-valued marks), for which we aim to establish Gaussian fluctuations of their sums. This is preceded by moment assumptions.
\begin{defn}[Moment conditions of scores]\label{d.moment}\ %
The real-valued marks of \\
$(\tilde\P_n)_{n \in \N}=\bigl(\{(x_i,\xi_{i,n})\}_{x_i \in \P_n}\bigr)_{n \in \N}$,
namely $\xi_{i,n}\in\R$, satisfy a $p$-moment condition, $p \in [1, \infty)$, if there are constants $M_p^{\xi} < \infty$
 such that 
\begin{equation}
\label{e:xinpmom}
\sup_{1 \leq n< \infty} \sup_{1 \leq q \leq p} \sup_{x_1,\ldots,
  x_{q} \in W_n} \sE_{\bk x1{q}}[ \max(1,|\xi_{1,n}|^p)]
\leq M_p^{\xi} < \infty,
\end{equation}
where $\xi_{1,n}$ is the score  of $x_1$ under the Palm distribution of  $\tP_n$ given the ground points $x_1,\ldots,x_q$. Here and elsewhere, we adopt the convention that the sup with respect to $\bk x1p$ is to be understood as $\text{ess}\sup$ with respect to the factorial moment  measures $\rho^{(p)}(\bk x1p) \md \bk x1p$ of the common ground process~$\P$, where  $\rho^{(p)}$ is tacitly assumed to be uniformly bounded by $\kappa_p$. Without loss of generality we assume that $M_p^{\xi}$
  is increasing in $p$ whenever \eqref{e:xinpmom} holds.
 \end{defn}
The formulation of the $p$-moment condition~\eqref{e:xinpmom} in terms of
higher-order Palm distributions, restricted to orders $q \in [1,  p]$, is motivated
by the cumulant method and justified by considering the moments
\(
\sE[(\sum_{x_i \in \P_n}\xi_{i,n})^k]\) with $k \in \N$,
of sums of scores over the window. Expanding these powers, one is led to
control expectations of sums of products of powers $\xi_i^{k_i}$, with
$k_i \le k$, taken over distinct $q \le k$ points of the ground process.
Such expressions are naturally related to $q$-fold Palm distributions.
Consequently, by applying H\"older's inequality, it suffices to impose
moment conditions of order $p = k + \epsilon$ under Palm measures of
order $q \le p$.

\BL-mixing correlations, defined via the class of \BL-test functions, 
are sufficient to establish the central limit theorem for 
real-valued marked point processes.

\begin{theorem}[CLT for sums of scores having fast \BL-mixing correlations]
\label{t:clt_linear_marks}
Consider the family  
$(\tilde\P_n)_{n \in \N}:=  \bigl(\{(x_i,\xi_{i,n})\}_{x_i \in \P_n}\bigr)_{n \in \N}$,
$\xi_{i,n}\in\R$,  of real-valued marked point processes sharing the same ground point process $\P$ in $\R^d$, all satisfying Assumption~\ref{Ass1}.
Assume that $(\tP_n)_{n \in \N}$ has fast $\BL$-mixing   correlations as in Definition \ref{d.omegahmixing}\ref{i.family-mixing} with scores $\xi_{i,n}$ 
satisfying  the $p$-moment condition 
    \eqref{e:xinpmom} for all $p \in \N$. 
    Furthermore,  if  
   $f \in \B(W_1)$ satisfies $\Var \,\mu_{n}^{\xi}(f) =
 \Omega(n^{\nu})$ for some $\nu > 0$,  then as $n \to \infty$, 
 \begin{equation}\label{e.CLT-Th2general}
(\Var \,\mu_{n}^{\xi}(f) )^{-1/2}\Big( \mu_{n}^{\xi}(f) - \sE \mu_{n}^{\xi}(f) \Big) \stackrel{d}{\Rightarrow} Z,
\end{equation}
where $Z$ is the standard normal random variable and $\stackrel{d}{\Rightarrow}$ denotes convergence in distribution.
\end{theorem}
Theorem \ref{t:clt_linear_marks}, {\em an umbrella central limit theorem}, is proved in Section~\ref{s:proofssclt}.  
It is  a central limit theorem for 
triangular arrays with the summands of $\mu_n^\xi$ forming the entries of the $n\,$th row of the array, $n \in \N$, and it does not require stationarity assumptions. The  result
encompasses the asymptotic normality of statistics of 
spin systems (Section \ref{s:gibbsmarking}), interacting diffusions (Section \ref{NEW-s:id_sprg}),  interacting particle systems (Section \ref{s:applnsips}),
as well as empirical random fields and geostatistical models (Section~\ref{s:erfgeostat}),
with all models considered on windows $W_n$, $n \to \infty.$

The main idea behind the proof is to show, together with moment conditions on the marks $\xi_{i,n}$, that
 the correlation functions  defined via Palm expectations $\E_{\bk x1k}=\E_{x_1,\ldots,x_k}$  by 
\begin{align}
m^{\bk k1{p+q}}(\bk x1{p+q};n)
& = m^{(k_1,...,k_{p + q})} (x_1,\ldots,x_{p + q};n) \nonumber \\
&:= \sE_{\bk x1{p+q}}\left[(\xi_{1,n})^{k_1} \ldots(\xi_{{p + q},n})^{k_{p + q} }\right] \rho^{(p+q)}(\bk x1{p+q})  \label{eqn:mixedmoment}
\end{align}
approximately factorize into $m^{\bk k1{p}}(\bk x1{p};n)m^{\bk
     k{p+1}{p+q}}(\bk x{p+1}{p+q};n)$
uniformly in $n \in \N$ up to an additive error decaying faster than any power of the separation distance $s$, defined at \eqref{defs}. Here $x_1,...,x_{p + q}$ are distinct points in $\R^d$ and $k_1,...,k_{p + q}\in\N$. The approximate factorization 
\eqref{eqn:mixedmoment} may be viewed as a {\em geometric mixing condition}  on correlation functions.  It  
  implies Brillinger mixing for the (possibly signed) measures  
$\sum_{x_i \in \P_n} \xi_{i,n}\delta_{x}$ (see  Lemma~\ref{l.moments-cumulants}) and hence it implies a central limit theorem for $(\mu_n^\xi)_{n \in \N}$; see Theorem~\ref{l.CLT-Cumulants}.

Our approach stands apart from others, where limit theorems for point processes with dependent scores  are typically established under mixing assumptions on both the ground process and the scores; see, for example, \cite{Heinrich_Schmidt_1985, Heinrich2013asymptotic,HP,HLS,Ivanoff82,JP,Pa09}. However, such assumptions can be difficult to verify in practice.

In contrast, we assume an asymptotic independence condition formulated in terms of the correlation structure of these models. While our Definition~\ref{d.omegahmixing} of mixing correlations for marked point processes may appear to be yet another form of mixing, it simplifies verification by allowing the correlations of the ground process and the marks to be treated separately (see Proposition~\ref{p:pmix+mmix}). Furthermore, it facilitates the analysis of more complex marks constructed functionally from simpler ones satisfying stronger mixing conditions (see Definition~\ref{B-omega-mixing}, Theorem~\ref{t:mppmix_double}, and Corollary~\ref{c:mppmixcl_double}).

\subsection{Mixing correlations of points, marks, and marked processes}
\label{ss.Mixing-separately}

In this section, we again consider point processes marked by abstract marks. We recall from \cite{BYY19} the definition of mixing correlation functions for the ground process $\P$, introduce a corresponding notion for the marks alone, and then show that, taken together, they imply the mixing correlations of marked point processes as defined in Definition~\ref{d.omegahmixing}.

\begin{defn}[Mixing correlations of the ground
    process $\P$]\label{d.omegamixing}\ 
 \begin{enumerate}[wide,label=(\roman*),labelindent=0pt]
 \item\label{i.mixing-ground} 
  The correlation functions
$\rho^{(p)}, p \in \N$, 
 of $\P$ are {\em mixing} if for all $k \in \N$
there exists an increasing-decreasing family of functions $\omega_k : [0, \infty)
\to [0, \infty)$ such that  for all $p,q \in \mathbb{N},
x_1,\ldots,x_{p+q} \in \mR^d$, we have
\be \label{e.gpmixing}
\left| \rho^{(p+q)}(\bk x1{p+q}) - \rho^{(p)}(\bk x1p)\rho^{(q)}(\bk x{p+1}{p+q}) \right| \leq \omega_{p+q}(s),
\ee
where $s := d(\bk x1p,\bk x{p+1}{p+q})$  is at \eqref{defs}.
Denote 
$C_k:= \sup_s \omega_k(s)$ and observe 
that $(C_k)_{k \in \N}$ is non-decreasing.

\item \label{i.mixing-ground-fast} $\P$ has  {\em fast mixing  correlations} if 
$\omega_k$ is fast decreasing  for all $k \in \N$. 
\end{enumerate}
\end{defn}

As seen in \cite[Section 2.2]{BYY19} (see also Appendix \ref{s.admppinteraction}),
 point processes having  fast mixing  correlations include determinantal and permanental point processes with a fast decreasing kernel, the zero set of a Gaussian entire function, and rarified Gibbs point processes. Recall from  Assumption \ref{Ass1} that
 $\kappa_p := \sup_{\bk x1p\in(\R^d)^p} \rho^{(p)}(\bk x1p)$ and $\kappa_0 :=  \max\{ \kappa_1, 1 \}$.
As noted in ~\cite[(1.12)]{BYY19}, fast mixing  correlations yields for all $p\in \N$
\begin{equation}\label{e.correlation-functions-bound}
\kappa_p  \leq p C_p \kappa_0^p.
\end{equation}

Recall that each point $x \in \P$ is associated with a mark $U(x) \in \K$, $\K$ a Polish space.
We formulate the notion of \BL-mixing correlations for the marks $U(x)$ of a point process $\tP$, without assuming any specific correlation properties of the ground process $\P$ itself.  

\begin{defn}[$\BL$-mixing correlations of marks]
\label{d.omegatmixing}\ 
 \begin{enumerate}[wide,label=(\roman*),labelindent=0pt]
\item \label{i.mixing-marks} The  marks $U$ of the point process $\tilde \P=\sum_{x \in \P}
  \delta_{\tx}$, $\tx=(x,U)$, $x\in\R^d$,  $U\in\K$,
satisfy  {\em $\BL$-mixing correlations} if
there exists an increasing-decreasing family of  functions $\omega_k :  [0, \infty) \to [0,2]$, $k \in \N$, called correlation decay functions,  such that for all $p,q \in \mathbb{N}, x_1,\ldots,x_{p+q} \in \mR^d$, and all $f \in \BL(\K^p)$, and
  $g \in \BL(\K^q)$ we have
\begin{align} \label{e.markmixing}
\left|\sE_{\bk x1{p+q}}[f(\bk U1{p})g(\bk U{p+1}{p+q})] - \sE_{\bk x1p}[f(\bk U1p)]\sE_{\bk x{p+1}{p+q}}[g(\bk U{p+1}{p+q})]\right|\le \omega_{p+q}(s),
\end{align}
where $s := d(\bk x1p,\bk x{p+1}{p+q})$ is at \eqref{defs}.

\item  The marks $U$ have  {\em fast $\BL$-mixing correlations} if $(\omega_k)_{k \in \N}$ are fast decreasing.

\end{enumerate}
\end{defn}

Before formulating the intuitively appealing result that mixing correlations of the points (sites), together with those of the marks (states) as developed in this section, imply the  mixing correlations of the marked point process (as in Definition~\ref{d.omegahmixing}), we introduce a stronger notion of  mixing correlations for the marks--- and consequently for the entire process. 
  Specifically, by replacing the class of bounded Lipschitz functions $\BL(\K)$ with the broader class of bounded functions $\B(\K)$, we obtain the corresponding definitions of $\B$-mixing correlations. This stronger form of mixing is instrumental in the applications of Theorem~\ref{t:clt_linear_marks}. 
    Specifically, if the sequence \((\tP_n)_{n \in \N}\) exhibits \(\B\)-mixing correlations that decay to zero exponentially fast (see Definition~\ref{def.A2} for more precise formulation of this exponential decay), and if the  marks \(\xi_{i,n}\)  are functionally constructed from this initial marking---i.e., \(\xi_{i,n} = \xi(x_i, \tP_n)\)---then the resulting marked process satisfies \BL-mixing correlations, provided that \(\xi\) fulfills certain localization or stabilization conditions; see Theorem \ref{t:mppmix_double}.  
\begin{defn}[$\B$-mixing correlations]
\label{B-omega-mixing} 

\begin{enumerate}[wide,label=(\roman*),labelindent=0pt]\

\item \label{B-mixing-marks} The  marks $U$ of the point process $\tilde \P=\sum_{x \in \P}
  \delta_{\tx}$, $\tx=(x,U)$, $x\in\R^d$,  $U\in\K$,
satisfy  {\em $\B$-mixing correlations} if
 \eqref{e.markmixing} holds when \BL\ is replaced by $\B$; i.e.,
  for all $p,q \in \mathbb{N}$, and all $f \in \B(\K^p), g \in \B(\K^q)$, 
 \begin{align*}
\left|\sE_{\bk x1{p+q}}[f(\bk U1{p})g(\bk U{p+1}{p+q})] - \sE_{\bk x1p}[f(\bk U1p)] \sE_{\bk x{p+1}{p+q}}[g(\bk U{p+1}{p+q})]\right|
\, & \le \, \omega_{p+q}(s), \qquad x_1,\ldots,x_{p+q} \in \mR^d,
\end{align*}
where $s = \inf_{i \in \{1,...,p\}, j \in \{p + 1,...,p + q\} } |x_i - x_j|$ and $\omega_k :  [0, \infty) \to [0, 2]$, $k\in\N$ 
is an increasing-decreasing family of functions.

\item \label{Bmixing} A marked point process $\tilde \P = \{(x, U(x))\}_{x \in \P}$ has  {\em $\B$-mixing
  correlations}  if \eqref{e.mppmixing} holds when \BL\ is replaced by $\B$ i.e.,
  for all $p,q \in \mathbb{N}$, and all $f \in \B(\K^p), g \in \B(\K^q)$, 
 \begin{align} \nonumber
& \left|\sE_{\bk x1{p+q}}[f(\bk U1{p})g(\bk U{p+1}{p+q})]\rho^{(p+q)}(\bk x1{p+q}) - \sE_{\bk x1p}[f(\bk U1p)] \rho^{(p)}(\bk x1{p})\sE_{\bk x{p+1}{p+q}}[g(\bk U{p+1}{p+q})]\rho^{(q)}(\bk
x{p+1}{p+q})\right|\\
\label{e.Bmppmixing} &\hspace{0.1\linewidth}\le \omega_{p+q}(s),  \qquad x_1,\ldots,x_{p+q} \in \mR^d,
\end{align}
with $s,\omega_k, k \in \N$ as in Item (i). Similarly, a {\em family of marked point processes}   $(\tilde\P_l)_{l \in \N},$ has  $\B$-mixing correlations if Definition~\ref{d.omegahmixing}\ref{i.family-mixing} holds when \BL\ is replaced by $\B$. 

\item \label{Bmixing-fast} If  $(\omega_k)_{k \in \N}$ are fast decreasing functions, then   we respectively speak of  the  {\em fast $\B$-mixing} of  marks, a marked point process, or a family of marked point processes.
\end{enumerate}
\end{defn}

To obtain   
 \BL-mixing, respectively $\B$-mixing, correlations of marked point processes, we show
 that it actually suffices to show \BL-mixing, respectively $\B$-mixing correlations, of the marks themselves
 provided that the ground process alone already satisfies these properties.
 This opens up a wide spectrum of models, as the marks exhibiting these mixing correlations are numerous. Beyond the obvious examples of independent marks on a ground process exhibiting mixing correlations, and marks depending  within a fixed deterministic distance of the points in a Poisson process, we can also handle marks with dependencies governed by graphs on $\P$ that involve local interactions, which are made precise through the concept of stabilizing interaction graphs, defined in Section \ref{s:modelgibbs}.

\begin{proposition}[Joint mixing correlations of points and marks]
  \label{p:pmix+mmix}
Let $\tP$ be a marked point process whose ground process $\P$ has mixing correlations as
in Definition \ref{d.omegamixing} with decay function $\omega$. If the marks
 have $\BL$-mixing  correlations  as in Definition
\ref{d.omegatmixing} with mixing function  $\tilde\omega$ (or $\B$-mixing  correlations  as in Definition~\ref{B-omega-mixing} 
\ref{B-mixing-marks} with mixing function $\tilde\omega$) then $\tilde\P$ has
$\BL$-mixing   correlations as   in Definition~\ref{d.omegahmixing}\ref{i.mixing}  
(respectively, $\B$-mixing  correlations as in Definition~\ref{B-omega-mixing}\ref{Bmixing})
with mixing function 
$$
\hO_k(s) :=
\omega_k(s)+\tilde\omega_k(s)\max(\kappa_p\kappa_q:p+q=k, p,q \in \N).
$$
\end{proposition}
\begin{proof}[Proof of Proposition~\ref{p:pmix+mmix}] 
We are given mixing functions $\omega_k$ and  $\tilde\omega_k$
 of the ground  point process and of the
marks, respectively.
For all $\bk x1p\in (\R^d)^p$, $\bk x{p + 1}{p + q} \in (\R^d)^q$ the triangle inequality gives
\begin{align*}
& \quad \Bigl|\E_{\bk x1{p + q}}[f(\bk U1p) g(\bk U{p + 1}{p + q})] \rho^{(p + q)} (\bk x1{p + q}) -
\E_{\bk x1p} [f(\bk U1p)] \rho^{(p)} (\bk x1p) \E_{\bk x{p + 1}{p + q}  }[ g(\bk U{p + 1}{p + q} )] \rho^{(q)} (\bk x{p + 1}{p + q})\Bigr|
\\
& \leq \Bigl| \E_{\bk x1{p + q}} [f(\bk U1p) g(\bk U{p + 1}{p + q})] \Big( \rho^{(p + q)} (\bk x1{p + q}) -  \rho^{(p)} (\bk x1p)\rho^{(q)} (\bk x{p + 1}{p + q}) \Big) \Bigr|
\\
 &\hspace{1em} + \Bigl|  \E_{\bk x{p + 1}{p + q}  } [ f(\bk U1p) g(\bk U{p + 1}{p + q})] \rho^{(p)} (\bk x1p)\rho^{(q)} (\bk x{p + 1}{p + q})\\
& \hspace{2em}- \E_{\bk x1p} [f(\bk U1p) ]\rho^{(p)} (\bk x1p) \E_{\bk
 x{p + 1}{p + q}  } [g(\bk U{p + 1}{p + q} )] \rho^{(q)} (\bk x{p + 1}{p + q})\Bigr|
\\
&
=: I + J.
\end{align*}

Now  $f \in \BL(\K^p)$ and
  $g \in \BL(\K^q)$) and the bound~\eqref{e.gpmixing} implies that term $I$ satisfies
\begin{align*}
  I \leq  | \rho^{(p + q)} (\bk x1{p + q}) -  \rho^{(p)} (\bk x1p)\rho^{(q)} (\bk x{p + 1}{p + q})| \leq \omega_{p + q}(s)
\end{align*}
where $s := d(\bk x1p,\bk x{p+1}{p+q})$ is at \eqref{defs}. For term~$J$, using 
~\eqref{e.markmixing}  one obtains
$J \leq \tilde \omega_{p + q}(s) \kappa_p\kappa_q.$
This proves Proposition~\ref{p:pmix+mmix}.  
\end{proof}

\section{Functionally constructed marks and their  localization}
\label{s.Stabilizing}

In this section, we will consider marks associated with points $x_i \in \P_n$, in particular real-valued scores denoted by $\xi_{i,n}$ in our umbrella Theorem~\ref{t:clt_linear_marks}, that are \emph{constructed as functions of $\P_n$}, and possibly also of pre-existing marks of this input process.
Theorem~\ref{t:mppmix_double} 
provides the key link between the umbrella central limit theorem and the forthcoming Theorem~\ref{t:cltmarkedpp}, which extends this result to functionally constructed scores and greatly simplifies the verification of their \BL-mixing correlations.
Crucially for this link, the  approach relies on the concepts of {\em localization} and {\em stabilization}, presented respectively in Sections~\ref{sec:marking-function}, \ref{s:blcluster}, and~\ref{s:strongstab}. Localization (in two versions) is based on quantitative distributional convergence of scores/marks when evaluated on the  ground point process restricted to balls increasing up to $\R^d$, whereas stabilization, a stronger property, relies on the stopping set structure of the marking functions.

We emphasize that our two notions of localization---{\em \BL-localization} and {\em \BL-cluster-localization}, which are central to our approach---contrast with stabilization in that they do not require the approximations (obtained by restricting the ground point process to growing balls) to eventually coincide with the original scores for large balls, nor even to converge to them pathwise. Instead, only convergence in distribution in the bounded Lipschitz metric is required.

More specifically, consider a (say `pre-marked') point process $\tP = \{(x,U(x))\}_{x \in \P}$, with marks in some general space~$\K$.  The models considered here
allow for a collection  of new $\K'$-valued marks on $\tP$, represented as
 $$
 \xi(\tx, \tP), \ \tx = (x,u) \in \tP
 $$
and formalized as follows. 
\begin{defn}[Marking function]
\label{d:mark_fn}
A {\em marking or score function} is a measurable function $\xi: (\mR^d \times\K)\times \hat{\cN}_{\R^d \times \K} \longrightarrow
\K'$ or $\xi:
(\mR^d \times\K)\times\cN_{\R^d \times \K} \longrightarrow \K'$  
 where $\K'$ is a Polish space. The value of the function $\xi((x,u),\tmu)$ is relevant only for $\tx=(x,u)\in\tmu$ and can be defined  arbitrarily for $(x,u)\not\in\tmu$.
\end{defn}
Although the term `score function' is more commonly used in the context of
limit theorems, in this section we adopt the terminology `marking function'
to emphasize a structural or modeling viewpoint, in which
$\xi(\tx,\tP)$, for $x \in \P$, is regarded as a mark attached to the point
process~$\P$. In contrast, when studying limit results—such as central limit
theorems—where $\xi$ appears as a real-valued weight in sums over the point
process, we revert to the term `score function,' reflecting its analytical
role. Accordingly, we will occasionally use the terms `marks' and `scores'
interchangeably, depending on the perspective adopted.

\subsection{Bounded Lipschitz localization} 
\label{sec:marking-function}

We  introduce a form of localization for the marks 
$\xi(\tx, \tP)$
taking values in a general Polish space $\K'$, 
which is weaker
than those used almost ubiquitously  in the models considered in the classical theory of stochastic geometry, but which is still strong enough to insure limit theorems for statistics of spin systems, 
interacting diffusions, and empirical random fields.
This localization, termed \BL-localization, is understood in terms of {\em rates of convergence in distribution with respect to the $d_{\BL}$ metric} of the {\em ($r$-)restricted version} 
\be \label{restricted}
\xitr(\tx_1, \tP):= \xi( \tx_1,\tP\cap B_r(x_1) )
\ee
to $\xi( \tx_1,\tP)$ as $r \to \infty$, provided $\xi(\tx_1, \tP)$ is defined. Here $B_r(x):= \{ y:|y-x|\le r\}$ denotes the closed ball of radius $r$ centered at $x$  and we shall slightly abuse notation by defining for $\tmu\in\cNK$ and $B\in {\cal{B}}(\R^d)$ the restriction $\tmu\cap B := \{ (x,u) \in \tmu : x \in B \}$. $B^c_r(x)$ will denote the complement of $B_r(x)$ in $\R^d$.

Unlike classical stabilization, discussed in Subsection \ref{s:strongstab}, this localization is well suited when one can only control distributional convergence of marks on balls of increasing radius,  and where the balls do not necessarily constitute stopping sets. We use this criterion, here called {\em bounded Lipschitz localization}, to establish asymptotic normality 
as $n \to \infty$
as well as expectation and variance asymptotics for  $(\mu_{n}^{\xi}(f))_{n \in \N}, f \in \B(W_1)$, under appropriate moment conditions on $\xi$ ; see Theorem \ref{t:cltmarkedpp} and Proposition~\ref{expvar}. 

Recall from Section~\ref{s:prelim} that  $d_{\BL}(X,Y)$  denotes the bounded Lipschitz distance between random  vectors $X$ and $Y$.  Furthermore,  the vector of the values of the  marking function $\xi$ evaluated at marked points $\bk {\tilde x}1p=\bk {(x,u)}1p \subset\tilde\mu$, is denoted by
$$
\bk{\xi}1p(\tx,\tilde\mu):=
(\xi(\tilde x_1,\tilde\mu),\ldots,\xi(\tilde x_p,\tilde\mu)).
$$
This additional notation is necessary because the following localization is in terms of convergence rates for the {\em joint distribution} of 
$$\bk{\xitr}1p(\tx, \tP)
= \Big( \xi(\tilde x_1,\tP\cap B_r(x_1)),\ldots ,\xi(\tilde x_p,\tP\cap B_r(x_p)) \Big)$$
 to $\bk\xi1p(\tx, \tP)$ as $r \to \infty$. 

\begin{defn}[\BL-localization of marking functions]
\label{def.Lp-stabilizing_marking}%
Let $\xi$ be a marking function and, depending on whether $\tP$ is a finite or infinite marked point process, $\xi$ is defined on either finite or infinite point sets, i.e., $\xi: \mR^d \times\K\times \hat{\cN}_{\R^d \times \K} \longrightarrow \K'$ or $\xi: \mR^d \times\K\times\cN_{\R^d \times \K} \longrightarrow \K'$.
\begin{enumerate}[wide,label=(\roman*),labelindent=0pt]
\item \label{i.BL-localizing}
We say that  $\xi$ is a {\em bounded Lipschitz-localizing marking function ({\em \BL-localizing} for short)} on the point process $\tP$ if 
there 
is an increasing-decreasing family of 
functions~$\varphi_p: [0, \infty) \to [0,1], p \in \N,$ such that 
\begin{align} \label{Lp-stab-infinite}
 \sup_{x_1,\ldots,x_p \in \R^d} 
 d_{\BL, \bk x1p}( \bk\xi1p(\tx,\tP),\bk{\xitr}1p( \tx,\tP))  \leq 2\varphi_p(r), \ r  \geq 0, 
\end{align}
where the extra notation $\bk x1p$ in $d_{\BL, \bk x1p}$ accounts for the distribution of the vector of scores  $\bk {\xi}1p(\tx,\tP)$ and $ \bk{\xitr}1p( \tx,\tP))$ under the Palm distribution $\mP_{\bk x1p}$ of $\tP$. 

\item \label{i.BL-localizing-windows} We say that   $\xi$  is {\em  \BL-localizing  on the finite  windows}
  of $\tP$ if there is an increasing-decreasing family of functions  $\varphi_p: [0, \infty) \to [0,1], p \in \N$,
  such that
\begin{align} \label{Lp-stab}
  \sup_{1 \leq n < \infty} \sup_{x_1,\ldots,x_p \in W_n} d_{\BL, \bk x1p}\big( \bk{\xi}1p(\tx,\tP_n),\bk{\xitr}1p( \tx,\tP_n)) \big)  \leq 2\varphi_p(r), \ r  \geq 0.
\end{align}

\item \label{i.BL-localizing-fast} We say that $\xi$ is {\em fast \BL-localizing} on $\tP$, or on the finite windows of $\tP$, if $(\varphi_p)_{p \in \N}$ are  fast decreasing. 

\item \label{i.BL-localizing-exponential} Exponential \BL-localization signifies that  $\varphi_p$ is exponentially decreasing to zero  for all $p \in \N$.   
\end{enumerate}
 As elsewhere, in~\eqref{Lp-stab-infinite} and~\eqref{Lp-stab} we adopt the
convention that the supremum over $\bk{x}{1}{p}$ is understood as an
$\mathrm{ess\,sup}$ with respect to the factorial moment measures of the
common ground process~$\P$, which are tacitly assumed to have  the form $\rho^{(p)}(\bk x1p)\md \bk{x}{1}{p}$ with $\rho^{(p)}$ uniformly bounded by $\kappa_p$.
\end{defn}

Note that \BL-localization requires that the marking function $\xi$  be neither bounded nor Lipschitz.  Furthermore, the condition $\varphi_p \in [0,1]$ is not restrictive since if there exists a non-negative $\varphi'_p$ satisfying \eqref{Lp-stab-infinite} or \eqref{Lp-stab}, then so does $\varphi_p = \min \{1,\varphi'_p \}$ and the decreasing or fast decreasing properties of $\varphi'_p$ are preserved by $\varphi_p$.

The \BL-localization property, as observed in
Remark~\ref{rem:comparison_stabilization} in
Section~\ref{s:strongstab}, follows readily from either classical
(stopping set) stabilization or $L^q$-stabilization, and therefore does
not restrict the class of models that can be treated. Moreover, the
\BL-localization criterion captures the closeness of {\em probability
distributions} of
\(\bk{\xitr}{1}{p}(\tx,\tP_n)\), 
\(\bk{\xi}{1}{p}(\tx,\tP_n)\), and  \(\bk{\xi}{1}{p}(\tx,\tP)\) if the latter is defined,
which need not be close in a pathwise sense. This weaker, purely distributional notion of localization provides additional
flexibility. A variant, introduced in the next section, is exploited in the
study of spin systems in Part~\ref{part:applications}, and both are  expected to be
useful in a broader class of models.

Fundamentally, \BL-localization of~$\xi$ on~$\P$ and~$\P_n$ yields
\[
\bk{\xi}{1}{p}(\tx,\tP_n) \tod \bk{\xi}{1}{p}(\tx,\tP),
\qquad n\to\infty,
\]
uniformly under all Palm measures~$\Palm_{\bk{x}{1}{p}}$, with a rate of
convergence in the \(d_{\BL,\bk{x}{1}{p}}\) distance governed by
\(4\varphi_p\).
Indeed, put $X_n :=  \bk{\xi}1p(\tx, \tP_n)$, $X_{n,r} :=\bk{\xitr}1p( \tx,\tP_n))$ , $X_r := \bk{\xitr}1p( \tx,\tP)$,
and $X := \bk{\xi}1p(\tx, \tP)$ with $r=d(\bk x1p,\partial W_n) := \min_{i=1,\ldots,p}d(x_i,\partial W_n)$ denoting the   distance between  $\bk x1p$ and  the boundary of $W_n$.
Then by \eqref{Lp-stab-infinite} and \eqref{Lp-stab}
we have
\begin{align}
d_{\BL, \bk x1p}(X_n,X) & \leq d_{\BL, \bk x1p}(X_n, X_{n,r})
+ d_{\BL, \bk x1p}(X_{n,r},X) \nonumber  \\
& = d_{\BL, \bk x1p}(X_n, X_{n,r})
+ d_{\BL, \bk x1p}(X_{r},X) \nonumber\\
& \leq 4 \varphi_p(d(\bk x1p,\partial W_n)), \label{weakconv}
\end{align}
where equality holds since the random variables  $X_{n,r}$ and $X_r$ coincide when $B_r(x_i) \subset W_n$ for all $i \in \{1,\ldots,p\}$.

This convergence can be extended to the expectations for real-valued $\xi$  under moment conditions
 of order $p = 1+\epsilon, \epsilon > 0$ (as in~\eqref{e:xinpmom}
 and with  $\xi_{i,n}$  representing $\xi( \tx_i, \tP_n )$ and  $\xi( \tx_i, \tP)$) with the rate 
$$| \E_{x} [ \xi( \tx, \tP_n)]- \E_{x} [\xi(\tx, \tP)] | \le \text{Const.} \times  \varphi_1(d(\bk x1p,\partial W_n))^{\epsilon / (2+\epsilon)}.
$$
See Lemma \ref{l:BL-limits-moments}
for details. 

\BL-localization on finite windows (namely condition ~\eqref{Lp-stab},
which does not require defining $\xi(\tx,\tP)$ on an infinite marked point
process on~$\R^d$) allows one to establish {\em distributional limits} for
\(\bk{\xi}{1}{p}(\tx,\tP_n)\) toward certain {\em probability kernels} as
\(n \to \infty\), uniformly under all Palm measures
\(\Palm_{\bk{x}{1}{p}}\), with the same rate of convergence, as well as for
their expectations. In this setting, these probability kernels
play exactly the role of  Palm distributions of  the (possibly
non-existent) infinite-volume object $\xi(\tx,\tP)$ in the analysis of
expectation and variance asymptotics of the {\em sums} of these \BL-localizing marks, which will be developed further in Sections~\ref{s:limitstatpp}, and \ref{s:stab_lemma_BL}--\ref{s:proofsexp}.
This relaxation makes it possible to
study statistics of spin systems on finite windows, as discussed in
Section~\ref{s:gibbsmarking}, without requiring the existence of an
infinite-volume spin system.

\subsection{Bounded Lipschitz cluster-localization}
\label{s:blcluster}

\BL-localization of marks, as formulated in
Definition~\ref{def.Lp-stabilizing_marking}, is sufficient for the theoretical
foundations and applications considered in
Part~\ref{part:applications}, except for the spin systems studied in
Section~\ref{s:gibbsmarking}. These systems involve yet another
relaxation of the classical stabilization framework, related to the
{\em irrelevance of the marking function} itself in the notion of
\BL-localization.

Indeed, \BL-localization of marks is a purely distributional notion: it does
not refer to a specific marking function~$\xi$, but rather to \emph{certain
laws} of the vectors
\(\bk{\xi^{(r)}}{1}{p}(\tx,\tP_n)\) for \(r>0\), which are required to be close
to the corresponding laws of
\(\bk{\xi}{1}{p}(\tx,\tP_n)\), uniformly in~$n$, and eventually to those of
\(\bk{\xi}{1}{p}(\tx,\tP)\), if the latter are well defined.

 However for Gibbs models there is typically no
canonical specification of the vectors
\(\bk{\xi^{(r)}}{1}{p}(\tx,\tP_n)\) for \(r>0\), i.e., the {\em joint law of spins} at \(\bk{x}{1}{p}\) given
{\em their separate local environments}
\(\P_n \cap B_r(x_i), 1 \leq i \leq p,\)
 which involve dependencies arising from overlap of balls \(B_r(x_i)\). While this does not hinder  the proof of \BL-localization
condition~\eqref{Lp-stab} for \(p=1\), it is a major obstacle in proving \BL-localization for \(p \ge 2\). 
The following version of \BL-localization, called  {\em cluster-localization}, is designed precisely to prove (fast) mixing
properties of correlations of~\(\xi\) despite the above complications.
It considers the vector
\begin{equation}
\label{e:xi-cl}    
 \bk{{\xi}^{\cup(r)}}{1}{p}(\tx,\tP)  =\bigl(\xi(\tilde x_1,\tP\cap B_r(\bk{x}{1}{p})),\ldots,\xi(\tilde x_p,\tP\cap B_r(\bk{x}{1}{p}))\bigr), 
 \end{equation}
whose components are the values of the  marking function~\(\xi\) computed with respect to  
the input process~\(\tP\) restricted to the \emph{union of the respective balls}
\(
B_r(\bk{x}{1}{p})=\bigcup_{i=1}^p B_r(x_i)
\).
\begin{defn}[\BL~cluster-localization of marks]
\label{def.Lp-stabilizing_marking-weak}%
In the setting of Definition~\ref{def.Lp-stabilizing_marking}, we say:
\begin{enumerate}[wide,label=(\roman*),labelindent=0pt]
\item  \label{i.BL-localizing-weak}
The marking function \(\xi\) is \emph{ \BL~cluster-localizing} on the point process~\(\tP\) if  
there 
is an increasing-decreasing family of  
functions $\varphi_p: [0, \infty) \to [0,1]$,  $p \in \N$,
and constants $\alpha_p\in[0,\infty)$ such that the following two conditions hold:
\begin{align}
\label{Lp-stab-infinite-weak1}
&\sup_{x_1,\ldots,x_p\in\R^d}
d_{\BL,\bk{x}{1}{p}}
\bigl(\bk{\xi}{1}{p}(\tx,\tP),\bk{\xi^{\cup(r)}}{1}{p}(\tx,\tP)\bigr)
\;\le\; 2\varphi_p(r),
\qquad r  \geq 0, %\\[1ex]
\end{align}
and for all $p \geq 2$, $l\in\{1,\ldots,p-1\}$ and $x_1,\ldots,x_p \in\R^d$ such that $d(\bk{x}{1}{l},\bk{x}{l+1}{p}) > r^{\alpha_p}$  
we have 
\begin{align}
\label{Lp-stab-infinite-weak2}
& 
d_{\BL,\bk{x}{1}{p}}
\Bigl(\bk{\xi^{\cup(r)}}{1}{p}(\tx,\tP),
\bigl(
\bk{\xi^{\cup(r)}}{1}{l}(\tx,\tP),
\bk{\xi^{\cup(r)}}{l+1}{p}(\tx,\tP)
\bigr)
\Bigr)
\;\le\; 2\varphi_p(r), \qquad r  \geq 0.
\end{align}
\end{enumerate}
(ii)~\emph{\BL~cluster-localization on finite windows},  
(iii)~\emph{fast \BL~cluster-localization},  
and (iv)~\emph{exponential \BL~cluster-localization}  
are defined \emph{mutatis mutandis} as in  
Definition~\ref{def.Lp-stabilizing_marking}.
\end{defn}

\begin{remark}[\BL-localization vs  \BL~cluster-localization]\label{r:two-BL-localization}
Replacing \(\bk{\xi^{(r)}}{1}{p}\) in Definition~\ref{def.Lp-stabilizing_marking}
by \(\bk{\xi^{\cup(r)}}{1}{p}\), which corresponds to retaining only the first
condition~\eqref{Lp-stab-infinite-weak1} of \BL~cluster-localization—a property
we sometimes refer to as \emph{\BL~union-localization}
—is often easier to verify.
Moreover, condition~\eqref{Lp-stab-infinite-weak1} alone already implies the
distributional limits of \(\bk{\xi^{\cup(r)}}{1}{p}(\tx,\tP_n)\) converge to
\(\bk{\xi}{1}{p}(\tx,\tP)\) and/or to the probability kernels discussed in the
previous section; see the detailed arguments in
Section~\ref{s:stab_lemma_BL}  and, in particular,  Remark~\ref{r:stab_lemma_BCL}.

However, \BL~union-localization does allow for cross-dependence between the
neighborhoods of distinct points \(x_i\), since these neighborhoods are taken
jointly. As a consequence, \BL~union-localization alone \emph{does not rule out
long-range interactions} between individual marks and therefore does not, by
itself, suffice to establish mixing correlations, nor even variance
asymptotics for sums of localized marks; see the counterexample
given in Example~\ref{ex:nonlocalizing} in
Section~\ref{s:stab_lemma_BL}.

From this perspective, condition~\eqref{Lp-stab-infinite-weak1} is complemented
in Definition~\ref{def.Lp-stabilizing_marking-weak} of
\BL~cluster-localization by an additional cluster-level localization condition,
namely~\eqref{Lp-stab-infinite-weak2}, which enforces a cut-off of interactions
between groups of points separated at scales larger than~\(r\).
Taken together, these two conditions imply all the distributional limits
developed in Section~\ref{s:clta} for sums of \BL-localizing marks, including
expectation and variance asymptotics in the stationary setting; see the details in Remark~\ref{r:stab_lemma_BCL-variance}, Section~\ref{s:proofsexp}.
\end{remark}

\subsection[Stabilization  via stopping sets]{Stabilization via stopping sets }
\label{s:strongstab}

We  consider marked point processes  with  marking functions satisfying  stopping set  stabilization, versions of which have appeared in  \cite{ Lachieze2019normal, SchulteYukich2023} in the setting of finite windows and also \cite{BYY19, Penrose2007laws, PY2003}  in the setting of infinite windows.

\begin{defn}[Radius of stabilization]\label{def.Stab.Radius}
Given  a marking  function $\xi :
W_n \times\K\times\hat{\cN}_{W_n \times \K }\longrightarrow \K'$, 
 $n\in\N \cup \{\infty\}$ (as usual for $n=\infty$, we take $W_\infty:=\R^d$) and $(x,u)\in\tmu \in\cN_{W_n \times \K}$  we 
define the radius of stabilization $R_{W_n}^\xi((x,u);\tmu)$ 
to be the smallest $r \in \N$    such that
\begin{equation}\label{eq:stopping set}
\xi((x,u), \tmu \cap B_r(x)) = \xi((x,u), (\tmu \cap B_r(x)) \cup ( \tilde \nu \cap B^c_r(x))) \,
\end{equation}
for  all $\tilde\nu\in \cN_{W_n \times \K}$. The definition is naturally extended to $\tmu \in\cN_{\R^d \times \K}$ by setting $R_{W_n}^\xi((x,u);\tmu) := R_{W_n}^\xi((x,u);\tmu \cap W_n)$.
\end{defn}
For $n\in\N$ the radius of stabilization satisfies the bound 
\begin{equation}\label{e.R_n<diam}
R_{W_n}^\xi((x,u);\tmu)\le \lceil \text{\rm diam}(W_n)\rceil.
\end{equation}
Indeed, \eqref{eq:stopping set} holds for $r=\text{\rm diam}(W_n)$. For  $n=\infty$ we abbreviate $R_{W_{\infty}}^\xi$ by 
$R^\xi$ and put $R^\xi:=\infty$ if no  finite $r$ satisfies~\eqref{eq:stopping set}.

We refer the reader to the discussion below \cite[Definition 2.1]{Penrose2007laws} for justification of measurability of the radius of stabilization. By definition, whether  the radius of stabilization 
satisfies the bound $R_{W_n}^\xi((x,u),\tmu) \leq r$ is determined by the realization of $\tmu \cap B_r(x)$ and not on the locations of
other points  in $W_n$ or their marks. In other words,  $B_{R_{W_n}^\xi((x,u),\tmu)}(x)$ is a {\em stopping set}, 
that is to say for all $r > 0$,
\begin{equation}
\label{e:stopsetstabrad}
\{ \tmu \in\hat\cN_{W_n \times \K} : x\in \mu,\, R_{W_n}^\xi((x,u),\tmu) \leq r \} = \{ \tmu \in\hat\cN_{W_n \times \K} : x\in \mu,\, R_{W_n}^\xi((x,u),\tmu \cap B_r(x)) \leq r \}
\end{equation}
and analogously for $W_\infty$ with $\tmu \in\cN_{\R^d \times \K}$.

\begin{defn}[Stabilizing marking
    functions via stopping sets]\label{def.stabilizing_marking}\ 
   \begin{enumerate}[wide,label=(\roman*),labelindent=0pt] 
\item \label{i.stabilizing}  We say that   $\xi$  {\em is  stabilizing on the marked
   point process $\tP$} if 
there 
is an increasing-decreasing family of   
functions $\varphi_p: (0, \infty) \to [0,1], p \in \N$   
   such that
\be \label{stab-infinite}
 \sup_{x_1,\ldots,x_p \in \R^d} \mP_{\bk x1p}
\bigl(R^\xi( \tx_1,\tP) > r \bigr) \leq \varphi_p(r), \ r > 0. 
\ee

\item \label{i.stabilizing-windows}  We say that   $\xi$  {\em is stabilizing on finite windows of $\tP$} if for all $p \in \N$ there 
is an increasing-decreasing family of   
functions $\varphi_p: (0, \infty) \to [0,1], p \in \N$,
such that
\be \label{stab}
 \sup_{1 \leq n < \infty} \sup_{x_1,\ldots,x_p \in W_n} \mP_{\bk x1p}
\bigl(R_{W_n}^\xi( \tx_1,\tP_n) > r \bigr) \leq \varphi_p(r), \ r > 0.
\ee
 \item \label{i.stabilizing-fast}  We say that $\xi$ is {\em fast stabilizing} on $\tP$, or on finite windows of $\tP$, if 
  $(\varphi_p)_{p \in \N}$ are  fast decreasing.
 Without loss of generality we always assume that $\varphi_p\le\varphi_{p+1}$ for all $p\in \N$.
\end{enumerate}
\end{defn}

\begin{remark}[Comparison between \BL-localization and   stabilization]
\label{rem:comparison_stabilization}
We justify our earlier remarks that \BL-localization is weaker than stabilization via stopping sets and $L^q$-stabilization. 

\begin{enumerate}[wide,label=(\roman*),labelindent=0pt]
\item (Stabilization via stopping sets is stronger than \BL-localization and \BL~cluster-localization).  The definition of the radius of stabilization gives  
$$ \sup_{f \in \BL((\K')^p)} \bigl|\mE_{\bk x1p}
\bigl[f\bigl(\bk\xi1p( \tx,\tP)\bigr) \bigr] - \mE_{\bk x1p}
\bigl[f\bigl(\bk{\xitr}1p( \tx,\tP)\bigr)
\bigr] \bigr| \leq 2\sum_{i=1}^p\mP_{\bk x1p}( R^{\xi}(\tx_i,\tP) > r) .$$
Thus, stabilization implies \BL-localization, fast stabilization implies fast
\BL-localization, and fast stabilization on finite windows implies fast
\BL-localization on finite windows.
The same line of arguments extends to establish
\eqref{Lp-stab-infinite-weak1} and~\eqref{Lp-stab-infinite-weak2}, namely
\BL~cluster-localization.

The stopping set stabilization property,  stronger than localization, confers benefits such as the explicit construction of scores on the infinite window $W_\infty :=\R^d$ as well as weaker moment conditions in the multivariate central limit theorem; this is described in Appendix \ref{s:fast_Stab_finite_windows}.
\item  \label{ii.L2->BL}($L^q$-stabilization is stronger than \BL-localization).
Given $q \geq 1$ and $\K'$-valued marks $\xi$, the marks satisfy $L^q$-stabilization with respect to $\tP$ if for all $p \in \N$
$$
\sup_{x_1,\ldots,x_p \in \R^d}\mE_{\bk x1p}|
\xi(\tx_{1},\tP) - \xi( \tx_{1},\tP\cap B_r(x_1))|^q \to 0
$$
as $r \to \infty$. It, too, is stronger than \BL-localization.  Indeed 
 \begin{align} \nonumber
 & \sup_{f \in \BL((\K')^p)}  \bigl|\mE_{\bk x1p}
\bigl[f\bigl(\bk\xi1p( \tx,\tP)\bigr) \bigr] - \mE_{\bk x1p}
\bigl[f\bigl(\bk{\xitr}1p( \tx,\tP))\bigr)
\bigr] \bigr| \\ \label{e.1-product-norm-used}
&
\leq \sum_{i=1}^p\mE_{\bk x1p}[|
\xi(\tx_i,\tP) - \xi( \tx_i,\tP\cap B_r(x_i))|]
 \\\nonumber 
&
\leq \sum_{i=1}^p \mE_{\bk x1p} \bigl[|
\xi(\tx_i,\tP) - \xi( \tx_i,\tP\cap B_r(x_i))|^q \bigr]^{1/q},
\end{align}
where in the first inequality we used 
the \BL ~property of $f$ with respect to the $1$-product metric on~$\R^p$.
$L^4$-stabilization has been used by \citet[(1.8)]{lachieze2019shotnoise} to obtain rates of normal convergence with respect to the Kolmogorov distance and variance asymptotics for geometric statistics of spatial shot-noise fields of Poisson processes. In our approach, \(L^1\)-stabilization, which is still weaker than stopping set stabilization, provides a more natural representation of the limits of the mean and variance established with the minimal assumption of \BL-localization. See the first remark in Section~\ref{s:remarksltthms} and further developments in Appendix~\ref{s:stab_lemma_Lp}.
 \end{enumerate}
\end{remark}

\subsection{Iterated functionally constructed marks}
\label{iterated}

We furnish  conditions showing how to verify the general assumptions of Theorem \ref{t:clt_linear_marks} when the marks $\xi_{i,n}$ are constructed as a \BL-localizing or \BL~cluster-localizing function of the input.
This result is established in Theorem \ref{t:mppmix_double} below, the main result of
Section \ref{s.Stabilizing}.

Given a marked point process $\tP$ having mixing correlations as in 
Definition~\ref{d.omegahmixing}, we equip $\P$ with a new set of functionally constructed marks via a marking function $\xi$ and establish that the new marked point process retains its   mixing correlations. 
To treat this case, we  require however the following strengthened version of $\B$-mixing correlations of $\tP$ (Definition~\ref{B-omega-mixing}), one controlling the growth rate in terms of the argument $s$ and the index $k$ of the decay functions $\omega_k(s)$. This strengthened version is a natural requirement when using FME expansions and the method of moments in general.  

\begin{defn}[Summable exponential  $\B$-mixing correlations of $\tP$ and $\P$]\label{def.A2}

\begin{enumerate}[wide,label=(\roman*),labelindent=0pt]\ 
\item \label{i.summ-exp-B-mixing} The marked point process  $\tP$, is  said to have {\em summable exponential $\B$-mixing correlations}  if it  has fast $\B$-mixing  correlations 
 as in Definition \ref{B-omega-mixing}\ref{Bmixing-fast} with functions
 $\omega_k(s) = C_k \phi(s)$, $k\in\N$,  
such that the $C_k$ are non-decreasing 
and where the $\omega_k(s)$ satisfy the twin growth conditions
\be \label{eqn:sum}
\sup_{k\in\N} \frac{C_k}{k^{ak}}  < \infty
\ee
for some constant $a\in[0,1)$ and  
\be \label{phibd}
\limsup_{r \to \infty} \frac{\log \phi(r)}{r^b } \in [-\infty,  0),
\ee
for some constant $b \in (0, \infty)$.

\item \label{i.summ-exp-mixingfamily} Similarly, the 
 family  
 $(\tilde\P_l = \{(x,U_l(x))\}_{x \in \P})_{l \in \N}$
 of marked point processes has 
{\em summable exponential $\B$-mixing correlations} 
if~\eqref{e.Bmppmixing} holds
for functions
 $\omega_k(s) = C_k \phi(s)$, $k\in\N$, with non-decreasing constants $C_k$,  satisfying~\eqref{eqn:sum} and~\eqref{phibd}
that is to say 
for all $l, p,q \in \mathbb{N}$, and all $f \in \B(\K^p), g \in \B(\K^q)$, we have
 \begin{align} \nonumber
& \left|\sE_{\bk x1{p+q}}[f(\bk {U_l}1{p})g(\bk {U_l}{p+1}{p+q})]\rho^{(p+q)}(\bk x1{p+q}) - \sE_{\bk x1p}[f(\bk {U_l}1p)] \rho^{(p)}(\bk x1{p})\sE_{\bk x{p+1}{p+q}}[g(\bk {U_l}{p+1}{p+q})]\rho^{(q)}(\bk
x{p+1}{p+q})\right|\\
\label{e.Bmppmixingnew} &\hspace{0.1\linewidth}\le C_{p + q} \phi(s),  \qquad x_1,\ldots,x_{p+q} \in \mR^d.
\end{align}

\item \label{i.summ-exp-mixing} The  point process $\P$ has
{\em summable exponential mixing correlations}
 if it has mixing correlations as in Definition \ref{d.omegamixing}\ref{i.mixing-ground} with 
$\omega_k(s) = C_k \phi(s)$, $k\in\N$, and where $\phi$ and $(C_k)_{k \in \N}$ satisfy the above conditions. 

\end{enumerate}
\end{defn}

The `exponential' in the terminology refers to the condition \eqref{phibd}, whereas `summability' refers to \eqref{eqn:sum} as it 
implies  finite bounds on the sums $\sum_k C_kw^k/k!$, $w > 0$. These are technical necessities for our proofs; see Remark \ref{Rmkconvsum}.

Appendix \ref{s.admppinteraction} and \cite[Section 2.2]{BYY19} provide examples of point processes $\P$ satisfying summable exponential mixing correlations. 
 Equipping such $\P$ with independent marks yields 
 a marked point process  $\tP$ having summable exponential $\B$-mixing correlations. 
 Such  point processes, which form the foundation of the spatial models in Part \ref{part:applications},  are enough to get us started in the
 iteration procedure and allows us to construct new point processes having valuable mixing properties. 
Indeed, if the ground point process \(\P\) exhibits summable exponential mixing correlations---or more generally, if a `pre-marked' process \(\tP = \sum_{x \in \P} \delta_{\tx} = \sum_{x \in \P} \delta_{(x, U(x))}\) exhibits summable exponential $\B$-mixing correlations, in the sense defined above---then the process
\(
\dbtilde\P = \sum_{x \in \P} \delta_{(x, \xi(\tx, \tP))},
\)
obtained by equipping points
$x \in \P$ with marks \(\xi(\tx, \tP) \in \K'\) that are fast stabilizing (resp. fast \BL-localizing or fast \BL~cluster-localizing), also exhibits fast \(\B \)-mixing (resp. fast \(\BL \)-mixing) correlations.

In other words, mixing properties of marked point processes are preserved under the construction of additional marks. The proof of this result and its corollary is deferred to Section~\ref{s:proof_exp_stab_marking}.

\begin{theorem}[Fast mixing properties of point processes with functionally constructed iterated marks]
  \label{t:mppmix_double}
  Let $\tilde\P$ be a marked point process on $\R^d\times\K$
  having summable exponential $\B$-mixing correlations  as in
Definition~\ref{def.A2}. Consider the functionally constructed mark $V(x) :=\xi(\tx, \tP)$, with 
  $\xi:
\mR^d \times\K\times\cN_{\R^d \times \K} \longrightarrow \K'$ and set $\dbtilde\P :=  \sum_{x \in \P}
\delta_{(x, V(x))}$, a marked point process on $\R^d\times \K'$. Fix $p,q \in \N$.
\begin{enumerate}[wide,label=(\roman*),labelindent=0pt]
\item \label{i.BL-iterated-marks} If $\xi$ satisfies 
   fast \BL-localization  
   \eqref{Lp-stab-infinite} on $\tP$ for $p,q$ and $p+q$, 
   then  $\dbtilde\P$ 
satisfies the fast \BL-mixing correlation condition for $p$ and $q$, that is to say \eqref{e.mppmixing} holds  for $f \in \BL((\K')^{p})$, $g \in \BL((\K')^{q}),$ and some fast decreasing function $\omega_{p+q}$. 
\item \label{i.BLc-iterated-marks}  Part~\ref{i.BL-iterated-marks}
 above remains valid if
$\xi$ satisfies fast \BL~cluster-localization---that is,
conditions~\eqref{Lp-stab-infinite-weak1} and~\eqref{Lp-stab-infinite-weak2}---on
$\tP$ for  $p,q$, and $p+q$.

 \item \label{i.iterated-marks} If 
  $\xi$ is fast stabilizing on $\tP$ as in \eqref{stab-infinite} (via stopping sets) for $p,q$ and $p+q$, 
then  $\dbtilde\P$ has fast $\B$-mixing correlations for $p$ and $q$, i.e., \eqref{e.mppmixing} holds for $f \in \B((\K')^{p})$, $g \in \B((\K')^{q}),$ and some fast decreasing function $\omega_{p+q}$.
\end{enumerate}
In  all cases~\ref{i.BL-iterated-marks}, \ref{i.BLc-iterated-marks}, and~\ref{i.iterated-marks}, the  fast decreasing function $\omega_{p+q}$ associated to $\dbtilde \P$ depends on the functions $\varphi_{p+q}$ and $\phi $ as well as on  the constants $(C_k)_{k = 1}^{\infty}$ and $\ka_{p+q}$ related to $\xi$ and $\tP$, respectively. In particular, there exist constants $C<\infty$, $c>0$ which  depend on $p,q$,  such that for all $s > 0$
\begin{equation}\label{e.iterated-marks-decay-function} \omega_{p+q}(s) \leq C \exp(-cs^b) + C\varphi_{p + q}\big( cs^{b(1-a)/2( p + q + d)} \big),
\end{equation}
where $a,b$ are as in~\eqref{eqn:sum} and \eqref{phibd}, respectively.
 \end{theorem}

\begin{corollary}[Mixing correlations of families of point processes with functionally constructed iterated marks]
  \label{c:mppmixcl_double} Assume that  the
  family $(\tP_{n})_{n \in \N}$ of marked point processes in $\R^d \times \K$ has summable exponential
$\B$-mixing  correlations as in Definition~\ref{def.A2}\ref{i.summ-exp-mixingfamily}. Put $V_{n}(x)=\xi(\tx, \tP_{n})$ and $\dbtilde\P_{n}:=\sum_{x \in \P_n}
\delta_{(x, V_{n}(x))}$, $n \in \N$ with $\xi:
\mR^d \times \K \times \hat{\cN}_{\R^d \times \K} \longrightarrow \K'$ a marking function.
  \begin{enumerate}[wide,label=(\roman*),labelindent=0pt]
\item If $\xi$ is fast \BL-localizing on finite windows of $\tP$, as in Definition~\ref{def.Lp-stabilizing_marking}\ref{i.BL-localizing-windows} and~\ref{i.BL-localizing-fast}, or fast \BL~cluster-localizing as in Definition~\ref{def.Lp-stabilizing_marking-weak}\ref{i.BL-localizing-weak} and (iii), then the family $(\dbtilde\P_{n})_{n \in \N}$ has fast \BL-mixing correlations as in Definition~\ref{d.omegahmixing}~\ref{i.mixing}, \ref{i.fast-mixing},  and~\ref{i.family-mixing}.

\item If $\xi$ is fast stabilizing on finite windows of $\tP$ as in Definition~\ref{def.stabilizing_marking}\ref{i.stabilizing-windows} and~\ref{i.stabilizing-fast}, then the family $(\dbtilde\P_{n})_{n \in \N}$ has fast $\B$-mixing correlations as in Definition~\ref{B-omega-mixing}~\ref{Bmixing} and~\ref{Bmixing-fast}.
 \end{enumerate}
 \end{corollary}

Corollary~\ref{c:mppmixcl_double}  will help establish the upcoming Gaussian fluctuations in Theorems \ref{t:cltmarkedpp} and \ref{t:multcltmarkedpp} for 
fast \BL-localizing or \BL~cluster-localizing statistics of point processes with summable exponential $\B$-mixing correlations. However, weaker assumptions suffice for  marking functions having a bounded interaction range, i.e. those which 
coincide with a restricted version $\xitr$ for some $r \in (0, \infty)$. 
Score/marking functions  having bounded interaction ranges include local  U-statistics.  Corollaries \ref{r.restricted-scores} and \ref{c.U-FME}
and the discussion therein further develops these remarks.

Finally, the procedure of constructing new marks $V(x)$ and $V_n(x)$, for
$x \in \P$, from existing pre-marks can in principle be iterated, yielding yet
another family of marks on $\P$. This procedure could be referred to as
\emph{iterated functionally constructed marks}.
However, a technical issue arises in this iteration. In the original setting, 
the pre-marks are required to satisfy \emph{summable exponential $\B$-mixing correlations},
encoded through
conditions~\eqref{eqn:sum} and~\eqref{phibd}, with decay functions of the form
$\omega_k(s)=C_k\phi(s)$.  
While the \emph{exponential decay} requirement in
condition~\eqref{phibd} can be readily ensured for the newly constructed marks
by imposing explicit assumptions on the localizing or stabilizing function
$\varphi_{p+q}$ appearing in~\eqref{e.iterated-marks-decay-function}, the verification of  the corresponding \emph{summability} condition~\eqref{eqn:sum} for
the associated constants $(C_k)_{k \in \N}$ is not always straightforward.

\section{Limit theory for sums of \BL-localizing scores}
\label{s:clta}
Having introduced the notions of marked point processes with mixing
correlations and \BL-localizing marking functions
(including the \BL-cluster-localization variant), we now establish the
limit theory for functionally constructed marks that are
localizing (in either sense) on point processes with summable
exponential $\B$-mixing correlations.
The general limit results established here form the foundation for the models
introduced in Section~\ref{s:Intro} and are further developed through
applications in Sections~\ref{s:gibbsmarking}--\ref{s:erfgeostat} of
Part~\ref{part:applications}.

In Section \ref{s:asnormmpp},  we  state a central limit theorem for a sequence of statistics of marked point processes constructed 
as sums of \BL-localizing or \BL~cluster-localizing marks on 
 expanding windows.  This result 
refines the central limit theorem of our general umbrella result,  Theorem \ref{t:clt_linear_marks}, 
using the localization property of functionally constructed marks to guarantee their \BL-mixing correlations. 
It neither requires  stationarity of the ground point process nor 
the existence of a marking function on an infinite marked  point
process, thus making the result  applicable in particular  to the study of spin systems and  particle systems on growing windows.  

When stationarity of the ground point process is assumed,  we  establish expectation and variance asymptotics for sums of \BL-localizing marks in Section \ref{s:limitstatpp}. These results can equally be obtained under the variant of
\BL~cluster-localization—specifically tailored to applications to spin
systems—but, for the sake of simplicity of exposition, we state them only  under
the default \BL-localization assumption.
Section \ref{s:ltthms_stat_marking} considers a variant of the above results when the marking function is defined on the entire point process, rather than on a finite window.

Finally, we comment on the assumptions underlying our results 
(such as moment conditions and variance bounds) 
and outline possible extensions in Section~\ref{s:remarksltthms}. 
That section shows that stronger assumptions—such as stopping set stabilization 
or $L^1$-stabilization---yield strengthened limit theorems. 
This is in contrast with Sections~\ref{s:asnormmpp}--\ref{s:ltthms_stat_marking}, 
where the limit theorems are derived under the weaker assumption of \BL-localization. 

\subsection{Asymptotic normality for sums of \BL-localizing scores}
\label{s:asnormmpp}
As in Section~\ref{s.Stabilizing}, consider a  `pre-marked' point process $\tP = \{\tx:=(x,U(x))\}_{x \in \P}$, with marks $U(\cdot)$ in a Polish space~$\K$. 
For a ground point process $\P$ on $\R^d$,  we have set $\P_n := \P \cap W_n, n \in \N,$ where $W_n := [-\frac{1}{2} n^{1/d}, \frac{1}{2} n^{1/d}]^d$. Consider an real-valued marking function $\xi(\tx, \tP_n)$ 
applied to marked points $\tx$ of $\tP_n$.
As in~\eqref{e:signed-measure} consider
(possibly signed) measures
\begin{equation}
\label{e:signed_measure_mark}
\mu_n^\xi: = \sum_{x \in \P_n} \xi( \tx, \tP_n ) \delta_{ n^{-1/d}x },
\end{equation}
where the marking function  $\xi(\tx, \tP_n)$ produces the (final) states $\xi_{i,n}$ in~\eqref{e:signed-measure}.

\begin{defn}[Moment conditions for  scores   on finite windows]\label{d.moment-function}\ %
We say that $\xi : \R^d \times \K \times \hat{\cN}_{\R^d \times \K} \to \R$ satisfies a $p$-moment condition with respect to $\tP$ 
 on {\em  finite windows}
 if $\xi( \tx_i, \tP_n )$
satisfy  the $p$-moment condition~\eqref{e:xinpmom}, that is to say
\begin{equation}
\label{e:xinpmomcopied}
\sup_{1 \leq n< \infty} \sup_{1 \leq q \leq p} \sup_{x_1,\ldots,
  x_{q} \in W_n} \sE_{\bk x1{q}}[ \max(1,|\xi(\tx_1, \tP_n )|^p)]
\leq M_p^{\xi} < \infty.
\end{equation}
\end{defn}
We refine the Gaussian fluctuations of Theorem~\ref{t:clt_linear_marks} for the sequence $(\mu_{n}^{\xi}(f))_{n \in \N}$, $f
\in \B(W_1)$.

\begin{theorem}[CLT for sums of fast \BL-localizing scores]
  \label{t:cltmarkedpp}
 Let $\tP$ be a marked point process on $\R^d\times\K$ having
 summable exponential $\B$-mixing correlations, as in Definition~\ref{def.A2}. Let $\xi: \R^d \times \K \times
  \hat{\cN}_{\R^d \times \K} \to \R$ be a marking function satisfying the  $p$-moment
  condition on finite windows~\eqref{e:xinpmomcopied}  for all
  \hbox{$p>1$}. 
Assume that $\xi$ is  fast \BL-localizing  on finite windows of $\tP$ as in Definition~\ref{def.Lp-stabilizing_marking}~\ref{i.BL-localizing-windows} and~\ref{i.BL-localizing-fast}  or  
that $\xi$ is fast \BL~cluster-localizing  on finite windows of $\tP$ as in Definition \ref{def.Lp-stabilizing_marking-weak}(ii) and (iii).  If $f \in \B(W_1)$ satisfies
 $\Var \,\mu_{n}^{\xi}(f)  = \Omega(n^{\nu})$ for some $\nu > 0$,
   then as $n \to \infty$,  the sequence $(\mu_{n}^{\xi}(f))_{n \in \N}$ satisfies  the central limit theorem ~\eqref{e.CLT-Th2general}.
 \end{theorem}

\begin{proof}  The proof is an immediate corollary of 
Theorem~\ref{t:clt_linear_marks} and  Theorem~\ref{t:mppmix_double}\ref{i.BL-iterated-marks} for $\xi$ fast \BL-localizing or, respectively~\ref{i.BLc-iterated-marks}, for $\xi$ fast \BL~cluster-localizing. Indeed, 
our assumptions on $\tP$ and $\xi$, together with part~\ref{i.BL-iterated-marks} or part~\ref{i.BLc-iterated-marks} of Theorem~\ref{t:mppmix_double} and part (i) of Corollary~\ref{c:mppmixcl_double}  yields that the family of marked point processes $\dbtilde \P_n =  \{(x,\xi(\tx,\tP_n))\}_{x \in \P_n},  n \in \N,$ has fast \BL-mixing  correlations as in Definition \ref{d.omegahmixing}\ref{i.mixing}, \ref{i.fast-mixing}, and~\ref{i.family-mixing}. 
By the assumed moment condition 
~\eqref{e:xinpmomcopied} on $\xi$ and assumed lower bound for $\Var \,\mu_{n}^{\xi}(f)$,  
it follows that $(\mu_{n}^{\xi}(f))_{n \in \N}$ satisfies all assumptions of Theorem~\ref{t:clt_linear_marks}, thus
yielding  the central limit theorem for $(\mu_{n}^{\xi}(f))_{n \in \N}$.
\end{proof}
\vskip.3cm
Theorem \ref{t:cltmarkedpp}  extends Theorems 1.13 and 1.14 of \cite{BYY19}, the main results of that paper, by 
allowing for a significantly broader class of input $
\P$ and scores $\xi$: 
(i) the point process $\P$ need not be stationary (ii) $\P$
may carry dependent marks,  (iii) the
score $\xi$ need not satisfy stopping set stabilization, and (iv) $\xi$ need not satisfy a power-growth condition. These comments also apply to the  upcoming theorems and propositions in this section. 
 Additionally,  if we take  $\xi \equiv 1$ in 
Theorem \ref{t:cltmarkedpp} then this gives the main result of \cite[Theorem 1.5]{Nazarov12}, which establishes a central limit theorem for the `number count' statistic of the zero set of a Gaussian analytic function. Such statistics are linear statistics of the ground point process.  We establish asymptotic normality for linear statistics of the marks themselves,  a more complex issue.  Our applications usually take $\P$ to be stationary,  but Theorems \ref{t:clt_linear_marks} and \ref{t:cltmarkedpp}   do not require this assumption. 

In the context of Theorem \ref{t:cltmarkedpp}, it is natural to inquire whether  one could deduce a central limit theorem  under  a $(2 + \varepsilon)$-moment condition on $\xi$, instead of
 assuming that all moments of $\xi$ exist.  As seen in the upcoming
Theorem \ref{t:multcltmarkedpp},  we may find general conditions (stationarity, volume-order variance growth and translation invariance of $\xi$) under which  a $(2 + \varepsilon)$-moment condition on $\xi$ would suffice.

 \subsection{Limit theory for sums of stationary \BL-localizing scores }
\label{s:limitstatpp}
Our theorems so far have not used any stationarity assumption on the
marked point process, i.e., invariance of the distribution of $\tP$ with respect to the translation of the sites.  However, doing
so yields expectation and variance asymptotics and thus also a
(multivariate) central limit theorem  with  explicit limiting covariances.  We establish the limit of the  expectation and variance of linear statistics of the measures $(\mu_n^\xi)_{ n \in \N}$ defined at \eqref{linearstatB}
under the stationarity of $\P$ and the
assumption that $\xi$ is  translation-invariant, i.e., for
all $y \in \R^d$ we have $\xi(\cdot+y,\cdot + y)=
\xi(\cdot,\cdot)$, where  ``$\cdot+y$'' acts
 only on  the locations of the ground process on $\R^d$.

The results of Section~\ref{s:asnormmpp}  assume  that the score function \(\xi\) is well-defined only on finite point processes, as is the case for example in spin systems and interacting particle systems. To develop explicit asymptotics of \(\xi(\tx, \tP_n)\) on \(\tP_n = \tP \cap W_n\) as \(n \to \infty\), one can use stopping set or \(L^1\)-stabilization as in Section~\ref{s:strongstab} ; see the remarks (i)-(iii) in Section~\ref{s:remarksltthms} and also Appendices \ref{s:fast_Stab_finite_windows} and~\ref{s:stab_lemma_Lp}. However,  the weaker framework of \BL-localization developed in Section~\ref{sec:marking-function} suffices to prove the limit of the 
expectation and variance of linear statistics of the measures $(\mu_n^\xi)_{ n \in \N}$.

Indeed, \BL-localization facilitates the derivation of certain {\em Palm distributional limits} of \(\xi(\tx, \tP_n)\), specifically for almost all \(x\in \R^d\) and for \(\rho^{(2)}\)-almost all \(\bk x12 \in (\R^d)^2\):
\begin{align}\label{e.xil-0}
\xi(\tx, \tP_n) &\xrightarrow[n \to \infty]{\text{BL}, x} \xil_{\0} \\
\bk{\xi}12(\tx, \tP_n) = \Bigl(\xi(\tx_1, \tP_n), \xi(\tx_2, \tP_n)\Bigr) &\xrightarrow[n \to \infty]{\text{BL}, \bk x12} \xil_{(\0, z)} = \Bigl(\xil_{(\0, z)}[1], \xil_{(\0, z)}[2]\Bigr),\label{e.xil-xy}
\end{align}
where $z=x_2-x_1$,  the \BL-convergence \(\xrightarrow[n \to \infty]{\text{BL}, \bk x1p}\) represents the distributional convergence of the vector of scores \(\bk {\xi}1p(\tx, \tP_n)\) under the Palm distribution \(\mP_{\bk x1p}\) of \(\tP\).
To demonstrate the existence of these limits, we prove the Cauchy property of the corresponding sequences in the bounded Lipschitz metric and exploit the completeness property of the metric; see Sections \ref{s:stab_lemma_BL}-\ref{s:proofsexp}.

In contrast to stopping set stabilization, which insures a.s. convergence to a limit random variable (cf. Appendix \ref{s:fast_Stab_finite_windows}), 
\BL-localization provides only distributional convergence to asymptotic probability distributions (more precisely {\em probability kernels}~\footnote{For a given order \(p \in \N\), \(\xil_{\bk{x}{1}{p}}\) represents a probability  kernel from  \((\R^d)^p\) to \(\R^p\). That is, for fixed \(\bk{x}{1}{p}\in (\R^d)^p\),  $ \xil_{\bk{x}{1}{p}}$ represents some probability distribution on \(\R^p\) (representing the  joint Palm \BL-limit of the scores at the locations~$\bk x1p$). Furthermore, given a measurable \(f: \R^p \to \R\), the mappings \((\R^d)^p \owns \bk{x}{1}{p} \mapsto \E[f(\xil_{\bk{x}{1}{p}})] \in \R\) are measurable. This is a consequence of $ \xil_{\bk{x}{1}{p}}$ being the  limit of the Palm kernels of \(\xi(\tx, \tP_n)\in\R\).}) represented by the random variable \(\xil_{\0}\) and the vector \(\xil_{(\0, x)}:=(\xil_{(\0, x)}[1], \xil_{(\0, x)}[2])\).
However, these are sufficient to formulate the asymptotics of the mean and variance as presented in the forthcoming Proposition~\ref{expvar}. More precisely, the two moments of these distributions correspond to the limits of two correlation functions of the scores  \(\xi\) (interpreted as marks of the ground process) evaluated on finite windows. As mentioned in Remark~\ref{r:two-BL-localization}, the same limits remain valid
when \BL~cluster-localization is assumed.

Equipped with these asymptotic probability kernel representations of the mean and  the correlation function of \(\xi\), we prove the expectation and variance asymptotics for the measures \((\mu_n^{\xi})_{n \in \N}\). We emphasize that these asymptotics heavily depend on the stationarity of the point processes and the translation invariance of the score functions. We do not seek minimal hypotheses to ensure these asymptotic results; this comment particularly applies to the assumed fast mixing of \(\xi\) and the summable exponential $\B$-mixing correlations of \(\tP\). Refinements of these assumptions are discussed in Section~\ref{s:remarksltthms}.

Denote by $\rho\equiv\rho^{(1)}(x)$
the intensity of the (stationary) point  process $\tP$.
 
\vskip.3cm
\begin{proposition}[Asymptotic mean and variance]
\label{expvar}
Let  $\tP$  be a {\em stationary} marked point process  
   on $\R^d\times\K$ with non-null, finite intensity
     $\rho$.  Let $\xi : \R^d \times \K \times \hat\cN_{\R^d \times \K} \to \R$ be a translation invariant marking function.
\begin{enumerate}[wide,label=(\roman*),labelindent=0pt]
\item \label{i.LLN-expvar} 
Let $\xi$ be fast \BL-localizing on all finite windows of $\tP$ as in 
\eqref{Lp-stab} for $p = 1$ and 
assume $\xi$
also
satisfies  a $(1 + \varepsilon)$-moment condition on {\em finite windows} 
as at~\eqref{e:xinpmomcopied} for {\em some} $\varepsilon > 0$. Then for almost all $x\in\R^d$ as $n\to\infty$ the Palm distributional limit of $\xi(\tx,\tP_n)$, 
given by  $\xil_{\0}$
in \eqref{e.xil-0}, exists and for all $f \in \B(W_1)$
\begin{equation}
\label{expasy}
\Big| n^{-1} \E \mu_n^\xi(f) -  \rho \E[\xil_{\0}]\int_{W_1}f(x)\,\md x\Big| = O(n^{-1/d}).
\end{equation}
If $\xi$ satisfies \eqref{Lp-stab}  for $p=1$ with
$\varphi_{1}$ only decreasing to zero, then the right-hand side of \eqref{expasy} is $o(1)$.

\item \label{i.Var-expvar}  
Assume $\xi$ is fast  \BL-localizing on finite windows of $\tP$  as in \eqref{Lp-stab} for $p \in \{1,2\}$  and assume 
$\xi$ satisfies a $(2 + \varepsilon)$-moment condition on finite windows as at~\eqref{e:xinpmomcopied} for some $\varepsilon > 0$. 
Then for \(\rho^{(2)}\)-almost all \(x_1,x_2\in\R^d\) as $n\to\infty$ the Palm distributional limit of $(\xi(\tx_1,\tP_n),\xi(\tx_2,\tP_n))$
given by $\xil_{(\0,z)}=(\xil_{(\0,z)}[1],\xil_{(\0,z)}[2])$, $z=x_2-x_1$,
in \eqref{e.xil-xy} exists.

Further, if $\tP$ has  summable exponential $\B$-mixing correlations as in Definition~\ref{def.A2} then for all $f \in  \B(W_1)$
\begin{equation}
\label{eqn:var}
\lim_{n \to \infty} n^{-1} \Var \,\mu_n^{\xi}(f) =  \sigma^2(\xil) \int_{W_1} f(x)^2 \,\md x \in [0,\infty),
\end{equation}
whereas for all $f, g \in \B(W_1)$
\begin{equation}
\label{eqn:Cov}
\lim_{n \to \infty} n^{-1}  {\rm{Cov}} ( \mu_n^{\xi}(f) , \mu_n^{\xi}(g) )  =  \sigma^2(\xil) \int_{W_1} f(x) g(x) \,\md x,
\end{equation}
where 
\begin{align}  \label{sigdef}
\sigma^2(\xil) =
\rho\sE[\xil_{\0}^2]  \ + \int_{\R^d} \left( \sE\Bigl[\xil_{
(\0,z)}[1]\xil_{(\0,z)}[2]\Bigr] \rho^{(2)}(\0,z) - \rho^2\sE[\xil_{\0}]^2 \right) \md z \in [0, \infty).
\end{align}
\end{enumerate}
\end{proposition}

The  proof of the above  result is  in Section~\ref{s:proofsexp}. Under stronger stabilization or localization assumptions, \(\E_{\0}[\xi(\tilde\0,\tP)]\), \(\E_{\0}[\xi^2(\tilde\0,\tP)]\), and \(\E_{\0,x}[\xi(\tilde{\0},\tP)\xi(\tx,\tP)]\), or variants thereof, can be used in place of expressions  \(\E[\xil_{\0}]\), \(\sE[\xil_{\0}^2]\), and \(\sE\bigl[\xil_{(\0,z)}[1]\xil_{(\0,z)}[2]\bigr]\), respectively. We elaborate in Remarks~\ref{ii.xi-extension-L1}, \ref{ii.xi-extension-stopping set}  and~\ref{ii.xi-infinite-window} of Section \ref{s:remarksltthms}.

\medskip 

To establish the multivariate central limit theorem  for the vector \((\mu_n^{\xi_i}(f))_{i=1,\ldots,{m}}\) formed by the statistics of different score functions \(\xi_i\), \(i=1,\ldots,m\), on the same input process \(\tP\), we need {\em joint \BL-localization of these score functions}. This means the vector \(\bxi = (\xi_1, \ldots, \xi_{m}) \in \R^{m}\) must satisfy the respective assumptions of Definition \ref{def.Lp-stabilizing_marking}~\ref{i.BL-localizing}, replacing \(\xi \in \K'\) therein by  $\bxi \in \R^m$.  Specifically, joint \BL-localization of \(\bxi\) of (Palm) order \(p\) on finite windows means that condition~\eqref{Lp-stab} is satisfied by \(\bk {\bxi}1p \in (\R^{m})^p\), where the bounded Lipschitz functions with respect to the \(\ell^1\)-metric on \((\R^{{m}})^p\) are considered under the \(d_{\BL, \bk x1p}\) metric on the probability measures on \(\R^{m}\). Recalling that  \BL-localization involves probability distributions and not specific realizations of random variables, it follows that  \BL-localization of the marginals \(\xi_i\) does not imply \BL-localization of \(\bxi\).

Similar to \eqref{e.xil-0} and \eqref{e.xil-xy},  under the two orders \(p \in \{1,2\}\) of this joint \BL-localization,  one can establish the convergence of two Palm kernels of \(\bk {\bxi}1p(\tx, \tP_n)\):
\begin{align}\label{e.xil-0-multi}
\bxi(\tx, \tP_n) &\xrightarrow[n \to \infty]{\text{BL}, x} \bxil_{\0} = \Bigl(\xil_{\0}[i]\Bigr)_{i=1,\ldots,m}, \\
\bk{\bxi}12(\tx, \tP_n) = \Bigl(\bxi(\tx_1, \tP_n), \bxi(\tx_2, \tP_n)\Bigr) &\xrightarrow[n \to \infty]{\text{BL}, \bk x12} \bxil_{(\0, z)} = 
\Bigl(\bigl(\xil_{(\0, z)}[1,i]\bigr)_{i=1,\ldots,m},
\bigl( \xil_{(\0, z)}[2,i]\bigr)_{i=1,\ldots,m}\Bigr)
,\label{e.xil-xy-multi}
\end{align}
for almost all \(x\) and for \(\rho^{(2)}\)-almost all \(\bk x12\),
where \(z = x_2 - x_1\).

The following result combines  the Gaussian fluctuations established in Theorem~\ref{t:cltmarkedpp} and the asymptotics of the mean and variance provided in Proposition~\ref{expvar} to establish a multivariate central limit theorem for \((\mu_n^{\xi_i}(f))_{i=1,\ldots,{m}}\). Importantly, this result relaxes the moment conditions of the former theorem to some moment \(p>2\), following an approach suggested to us by Matthias Schulte. However, this may result in a degenerate central limit theorem unless the variance is of volume order.
\begin{theorem}[Multivariate CLT  for statistics of point processes with fast \BL-localizing scores  under a $(2 + \varepsilon)$-moment condition]    \label{t:multcltmarkedpp}  
Consider a stationary marked point process  $\tP$ on $\R^d\times\K$ with non-null, finite intensity $\rho$, having summable exponential $\B$-mixing correlations as in Definition~\ref{def.A2}.  Let $\xi_i : \R^d \times \K \times \hat\cN_{\R^d \times \K} \to \R$,  $i =1,\ldots,{m}$ be translation invariant  score functions. Assume  that the  vector  \(\bxi = (\xi_1, \ldots, \xi_{m}) \in \R^{m}\) is jointly fast \BL-localizing on finite windows of $\tP$ in the sense that for all $p \in \N $, condition~\eqref{Lp-stab} is satisfied by \(\bk {\bxi}1p \in (\R^{m})^p\). 
If for {\em some} $\varepsilon > 0$, 
   $\xi_i$, $i=1,\ldots,{m}$ satisfy the $(2 + \varepsilon)$-moment condition
   \eqref{e:xinpmomcopied} on finite windows, then for  $f \in \B(W_1)$ 
\begin{equation} n^{-1/2}(\mu_n^{\xi_1}(f) - \sE\mu_n^{\xi_1}(f),\ldots, \mu_n^{\xi_{m}}(f) - \sE\mu_n^{\xi_{m}}(f)) \stackrel{d}{\Rightarrow} (Z_1,\ldots,Z_{m}),\label{e.M-CLT}
\end{equation}
where $(Z_1,\ldots,Z_{m})$ is a  multivariate normal random vector
(with possibly zero components) having zero mean and covariance 
$$  {\rm{Cov}}(Z_i,Z_j) = \sigma^2_{\xil}[i,j] \int_{W_1}f(x)^2 \md x \in (-\infty, \infty),$$
 where the matrix $(\sigma^2_{\xil}[i,j])_{1 \leq i \leq j \leq {m}}$ involves the vector $\bxil_{\0}$ and the matrix  $\bxil_{(\0,z)}$ of probability distributions  given in~\eqref{e.xil-0-multi} and~\eqref{e.xil-xy-multi}: 
\begin{align}  \label{sigcovdef}
& \sigma^2_{\xil}[i,j]:=
\rho\sE\bigl[\xil_{\0}[i] \xil_{\0}[j]\bigr]   \\
& \ \ \ + \int_{\R^d} \left( \sE\Bigl[\xil_{(\0,z)}[1,i] \xil_{(\0,z)}[2,j]\Bigr] \rho^{(2)}(\0,z) - \rho^2
\sE\bigl[\xil_{\0}[i]\bigr] \sE\bigl[\xil_{\0}[j]\bigr] 
\right) \md z \in (-\infty, \infty). \no
\end{align}
\end{theorem}
The proof is in Section \ref{s:proof_multclt}. 
Under the stronger assumption of stopping set stabilization, it suffices to assume 
stabilization under Palm conditioning only of order two,  i.e., we only require that 
$\xi_i, 1 \leq i \leq m,$ satisfy
\eqref{stab} 
for $p \in \{1,2\}$;  see Corollary \ref{c:multicltstab12}.

As seen in the proof in Section \ref{s:proof_multclt}, Theorem  \ref{t:multcltmarkedpp} 
allows for possibly  different exponents $\nu_i$ when
considering the variance lower bounds  $\Var \,\mu_{n}^{\xi_i}(f) = \Omega(n^{\nu_i})$ for  scores 
  $\xi_i,i=1,\ldots,{m}$. 
However such a variant would require that   
$\xi_i,i=1,\ldots,m,$ satisfy the $p$-moment condition
\eqref{e:xinpmomcopied} on finite windows for {\em all} $p > 2$.
In particular, if $\nu_i<1$ for some $i$ then this implies a vanishing limit in~\eqref{eqn:var}. As a
  consequence, the normal vector  $(Z_1,\ldots,Z_{m})$ might have
null components. If $\nu_i=1$ for some $i$,  then the above result
  establishes  the non-degenerate multivariate central limit theorem under only a $p$-moment condition for some $p > 2$.
  Extending the proof approach,  we may establish multivariate
   normal convergence for 
$$ 
n^{-\nu/2}\big(\mu_n^{\xi_1}(f) - \sE\mu_n^{\xi_1}(f),\ldots, \mu_n^{\xi_{m}}(f) - \sE\mu_n^{\xi_{m}}(f)\big)
$$
for some $\nu \in (0,\infty)$ provided  that for all $1 \leq i,j \leq {m}$,
$$ 
\lim_{n \to \infty} n^{-\nu}
{\rm{Cov}}(\mu_n^{\xi_i}(f),\mu_n^{\xi_j}(f)) =
\sigma_{\nu}(\xi_i,\xi_j ; f) \in [0,\infty).$$
One could likewise establish 
multivariate normal convergence for 
$$ 
n^{-\nu/2}\left(\mu_n^{\xi}(f_1) - \sE\mu_n^{\xi}(f_1),\ldots, \mu_n^{\xi}(f_{m}) - \sE\mu_n^{\xi}(f_{m})\right)
$$
for some $\nu \in (0,\infty)$ provided  that for all $1 \leq i,j \leq {m}$,
$$ 
\lim_{n \to \infty} n^{-\nu}
{\rm{Cov}}(\mu_n^{\xi}(f_i),\mu_n^{\xi}(f_j)) =
\sigma_{\nu}(\xi;f_i,f_j) \in [0,\infty).
$$

\subsection{Asymptotics for sums of scores on the infinite window}
\label{s:ltthms_stat_marking}
Unlike the previous section we now assume $\xi$ is defined on infinite input i.e., $\xi : \R^d \times \K \times \cN_{\R^d \times \K} \to \R$.
This assumption may not always hold, as when considering  statistics of spin systems and interacting particle systems, but it does apply to certain interacting diffusions and to empirical random fields considered in Part~\ref{part:applications}.
While the measures $(\mu^\xi_n)_{n\in \N}$ are defined in terms of the
scores $\xi(\tx, \tP_n)$, we now
  consider the measures induced by  {\em sums of scores  $\xi(\tx, \tP)$ with respect to input on  the infinite window}, namely
\be \label{linearstathat}
\hat{\mu}_{n}^{\xi} := \sum_{x \in \P_n} \xi( \tx, \tP) \delta_{ n^{-1/d}x }.
\ee
\begin{proposition}[Limit theory for marked point processes with scores on infinite windows]
\label{t:clt_linear_marks_new}
 Let $\tP$ be a marked point process on $\R^d\times\K$ having
 summable exponential $\B$-mixing correlations, as in Definition~\ref{def.A2} and let $\xi: \R^d \times \K \times
  {\cN}_{\R^d \times \K} \to \R$ be a marking function on $\tP$.
   \begin{enumerate}[wide,label=(\roman*),labelindent=0pt]
\item \label{i.CLT-infinite-xi}
Let $\xi$ be fast \BL-localizing on $\tP$
  as in Definition~\ref{def.Lp-stabilizing_marking}~\ref{i.BL-localizing},~\ref{i.BL-localizing-fast},
and assume it satisfies the $p$-moment
  condition 
  \begin{equation}
\label{e:xinpmomcopied-infty}
\sup_{1 \leq q \leq p} \sup_{x_1,\ldots,
  x_{q} \in \R^d} \sE_{\bk x1{q}}[ \max(1,|\xi( \tx_1, \tP )|^p)]
\leq M_p^{\xi} < \infty
\end{equation}
for all $p>1$. 
   If $f \in \B(W_1)$ satisfies
 $\Var \,\hat\mu_{n}^{\xi}(f)  = \Omega(n^{\nu})$ for some $\nu > 0$,
   then as $n \to \infty$, 
 \begin{equation}  (\Var \,\hat\mu_{n}^{\xi}(f))^{-1/2}\Big( \hat\mu_{n}^{\xi}(f) - \sE \hat\mu_{n}^{\xi}(f) \Big) \stackrel{d}{\Rightarrow} Z,
\end{equation}
where $Z$ denotes a standard normal random variable.
\item \label{i.EVar-infinite-xi}
If  moreover $\P$ is stationary and $\xi$ is translation invariant then the analogs of 
expectation,  variance asymptotics and multivariate central limit
theorems (as in Proposition~\ref{expvar} and Theorem \ref{t:multcltmarkedpp})
hold for $(\hat\mu_{n}^{\xi})_{n \in \N}$ where
\(\E_{\0}[\xi(\tilde\0,\tP)]\), \(\E_{\0}[\xi^2(\tilde\0,\tP)]\) and \(\E_{\0,x}[\xi(\tilde{\0},\tP)\xi(\tx,\tP)]\) can be used in place of \(\E[\xil_{\0}]\), \(\sE[\xil_{\0}^2]\) and \(\sE\bigl[\xil_{(\0,z)}[1]\xil_{(\0,z)}[2]\bigr]\), respectively, under a $p=(1+\epsilon)-$ or $p=(2+\epsilon)-$ moment condition~\eqref{e:xinpmomcopied-infty}.
\end{enumerate}
 \end{proposition}
The proof of Proposition \ref{t:clt_linear_marks_new} is in Section \ref{proofProp5.5}.

\subsection{Remarks}
\label{s:remarksltthms}\ 
In this section, we comment  on the results  in Sections~\ref{s:asnormmpp}--\ref{s:ltthms_stat_marking}. First, we  detail how the `minimal' representations of the limits of the two moments and the correlations of the score function \(\xi(\tx,\tP_n)\), as \(n \to \infty\) in Proposition~\ref{expvar}, available under \BL-localization on finite windows, can be replaced by the respective characteristics of some scores considered on the entire space $\R^d$. Next, we present more detailed extensions and remarks regarding Theorems~\ref{t:cltmarkedpp} and \ref{t:multcltmarkedpp} as well as    Propositions~\ref{expvar} and~\ref{t:clt_linear_marks_new}.

\begin{enumerate}[wide,label=(\roman*),labelindent=0pt]

\item (Asymptotic distribution of  scores via \(L^1\)-stabilization.) 
\label{ii.xi-extension-L1}
By assuming fast \(L^1\)-stabilization of $\xi$ on finite windows (see Appendix~\ref{s:stab_lemma_Lp}), which is stronger than fast \BL-localization, one may obtain a new {\em distribution} of scores \(\xil(x)\),  \(x \in \P,\) attached to points of the entire ground process \(\P\),  provided $\xi$ satisfies a $p=1$ moment condition on finite windows. These new scores  arise as  the distributional limit of the ground process marked in finite windows by the scores \(\xi(\tx, \tP_n)\):
$$
\sum_{x \in \P_n} \delta_{(x, \xi(\tx, \tP_n))}
\stackrel{d}{\Rightarrow}
\sum_{x \in \P} \delta_{(x, \xil(x))}.
$$
It is important to note that these scores (marks) \(\xil(x)\) are not obtained by applying some (existing or extended) score function \(\xi\) of the entire process. However, at least two Palm distributions of \(\xil(x)\) correspond to the two kernels \(\xil_x\) and \(\xil_{(x_1, x_2)}\), which are limits of the Palm kernels of \(\xi(\tilde{x}, \tilde{P}_n)\) as \(n \to \infty\), established using \BL-localization; see~\eqref{e.xil-0}, \eqref{e.xil-xy}, and Lemma~\ref{l:BL-limits}. In particular,  under respective moment conditions of $\xi$ on finite windows,  their moments can be used in Proposition~\ref{expvar}: 
$\E[\xil_{\0}] = \E_{\0}[\xil(0)]$,
$\E[\xil_{\0}^2] = \E_{\0}[\xil(0)^2]$, 
$\E\bigl[\xil_{(\0,x)}[1]\xil_{(\0,x)}[2]\bigr] = \E_{\0,x}[\xil(0)\xil(x)]$;
see the precise assumptions and formulation in Lemma~\ref{l:BL-limits-pp} in Appendix~\ref{s:stab_lemma_Lp}.

\item (Almost sure extension of $\xi$ via stopping set stabilization.)
\label{ii.xi-extension-stopping set}
In the case that \(\xi\) satisfies fast stabilization with respect to stopping sets, as in Definition~\ref{def.stabilizing_marking}\ref{i.stabilizing-windows}, the idea consists of extending the function \(\xi : \R^d \times \K \times \hat{\cN}_{\R^d \times \K} \to \R\) to the infinite domain \(\R^d \times \K \times \cN_{\R^d \times \K}\), denoted by \(\xi_\infty\), and satisfying
\begin{equation}\label{e:xi-to-infty}
\xi_{\infty}(\tx, \tmu) = \lim_{n \to \infty} \xi(\tx,\tmu_n),
\end{equation}
provided the limit exists for all \(\tmu\) (at least for those admitted by \(\tP\) with probability 1) and \(\tx \in \tmu\). This condition can be met under the additional condition that the stabilization radii are uniformly bounded, i.e.,
almost surely $\limsup_{n \to \infty} R^{\xi}_{W_n}(\tx,\tP_n) < \infty$.
We develop these arguments further in Appendix~\ref{s:fast_Stab_finite_windows}.
In particular, under this assumption, the extended score function \(\xi_\infty\) is fast stabilizing~\eqref{stab-infinite} admitting  $R^{\xi}_\infty(\tx,\tP):=\limsup_{n \to \infty} R^{\xi}_{W_n}(\tx,\tP_n)$ as the radius of stabilization and 
\(\xi_\infty\) inherits the moment conditions assumed for \(\xi\) in finite windows (see Lemma~\ref{l.xi-infnity}). 
Consequently, under the respective moment conditions of $\xi$ 
(on finite windows)
\(\E_{\0}[\xi_\infty(\tilde\0,\tP)]\), \(\E_{\0}[\xi_\infty^2(\tilde\0,\tP)]\), and \(\E_{\0,x}[\xi_\infty(\tilde{\0},\tP)\xi_\infty(\tx,\tP)]\) can be used in place of \(\E[\xil_{\0}]\), \(\sE[\xil_{\0}^2]\), and \(\sE\bigl[\xil_{(\0,z)}[1]\xil_{(\0,z)}[2]\bigr]\), respectively, in Proposition~\ref{expvar}.

\item (Score function $\xi$ is defined on the infinite window.) 
\label{ii.xi-infinite-window}
This is perhaps a  most natural  scenario where the initial score function is also well-defined on the entire input process \(\tP\).  If \(\xi\) is fast stabilizing on the infinite window as in \eqref{stab-infinite}, then \(\limsup_{n \to \infty} R^{\xi}_{W_n}(\tx,\tP_n) < \infty\) for all \(\tx \in \tP\) and  \(\xi_\infty = \xi\) a.s. (see Lemma~\ref{l.xi-infnity}); this assumption was  adopted in~\cite{BYY19}. If \(\xi\) is only fast \BL-localizing on the infinite window as in \eqref{Lp-stab-infinite} (this assumption is adopted  in Section~\ref{s:ltthms_stat_marking}) then at least the  Palm distributions  of $\xi$  in finite windows converge to it, i.e., $\bk{\xi}1p(\tx, \tP_n)
\xrightarrow[ n \to \infty]{\BL,\bk x1p} \bk {\xi}1p(\tx, \tP)\overset{d, \bk x1p}{=}\xil_{\bk x1p}$,
where the equality in distribution  holds under Palm distribution of $\bk {\xi}1p(\tx, \tP)$ given $\bk x1p$
(see Lemma~\ref{l:BL-limits}).
Consequently, in both cases, under the respective moment conditions of $\xi$ (at least on finite windows)
\(\E_{\0}[\xi(\tilde\0,\tP)]\), \(\E_{\0}[\xi^2(\tilde\0,\tP)]\), and \(\E_{\0,x}[\xi(\tilde{\0},\tP)\xi(\tx,\tP)]\) can be used in place of \(\E[\xil_{\0}]\), \(\sE[\xil_{\0}^2]\), and \(\sE\bigl[\xil_{(\0,z)}[1]\xil_{(\0,z)}[2]\bigr]\), respectively, in Proposition~\ref{expvar}.

\item   \label{i.rem-BLCuster} (\BL~cluster-localization.) Proposition~\ref{expvar}, Theorem~\ref{t:multcltmarkedpp}, and
Proposition~\ref{t:clt_linear_marks_new} remain valid when the marking functions
satisfy \BL~cluster-localization instead of \BL-localization, in complete
analogy with the two admissible choices in
Theorem~\ref{t:cltmarkedpp}.  
More details are in Remarks~\ref{r:stab_lemma_BCL} and~\ref{r:stab_lemma_BCL-variance}.

\item \label{i.rem-short-range-scores}
(\BL-localization via short-range scores.)
Another generalization of the concept of \BL-localization consists in
postulating, in Definition~\ref{def.Lp-stabilizing_marking}, the existence of
a family of short-range marking functions $(\xi^{[r]})_{r>0}$.
By a {\em short-range function}---i.e., a member $\xi^{[r]}$ of this family for some
$r>0$---we mean a marking function
\(
\xi^{[r]} : \R^d \times \K \times \cN_{\R^d \times \K} \longrightarrow \K'
\)
satisfying
\[
\xi^{[r]}(\tx,\tmu) = \xi^{[r]}(\tx,\tmu \cap B_r(x)),
\qquad \text{for all } \tx \in \tmu \in \cN_{\R^d \times \K}.
\]
{\em The family $(\xi^{[r]})_{r>0}$ of short-range functions is said to \emph{\BL-localize} a given marking
function $\xi$} if it satisfies \eqref{Lp-stab-infinite}, respectively
\eqref{Lp-stab}, with $\bk{\xitr}{1}{p}(\tx,\tP)$ replaced by
\[
\bk{\xi^{[r]}}{1}{p}(\tx,\tP)
=
\Big(
\xi^{[r]}(\tilde x_1,\tP),\ldots,
\xi^{[r]}(\tilde x_p,\tP)
\Big).
\]
Clearly, the $r$-restrictions $(\xi^{(r)})_{r>0}$ of the original marking $\xi$  defined in~\eqref{restricted}, form a natural example of such short-range functions.
Thus, \BL-localization via short-range scores extends
Definition~\ref{def.Lp-stabilizing_marking}. Moreover, under this alternative
notion of \BL-localization, based on real-valued short-range marks (scores),
Proposition~\ref{expvar}, Theorem~\ref{t:multcltmarkedpp}, and
Proposition~\ref{t:clt_linear_marks_new} remain valid, in complete analogy
with the admissible formulations in Theorem~\ref{t:cltmarkedpp}.
This follows directly from the proof of part~\ref{i.BL-iterated-marks} of
Theorem~\ref{t:mppmix_double}, which shows that the required mixing
correlation properties of the mark $\xi$ continue to hold when
\BL-localization is defined via such `external' short-range scores,
replacing the  short-range scores generated via natural restrictions of $\xi$; see
Remark~\ref{r.rem-short-range-scores}.
This relaxation is useful in situations where the natural  restrictions
$\xitr$ (which are themselves short-range scores) do not satisfy
\eqref{Lp-stab-infinite} or \eqref{Lp-stab}. A typical example arises when
the truncation must account for the geometry of the observation window
(e.g., balls $B_r(x)$ or windows $W_n$). For instance, Voronoi cells near
the boundary are unbounded, and one may wish to modify their associated
scores according to their distance from the boundary.
A specific example concerning extreme points in Laguerre
tessellations is developed in~\cite{TY}.

\item (Local statistics.) In the case of
scores which have a bounded interaction range, including local $U$-statistics, 
 one can prove limit theorems under weaker mixing correlation assumptions on $\tP$ as in~\cite{BYY19}. This point is further elaborated in Corollary~\ref{c.U-FME}.
 
\item (Quantitative CLTs---Rates of normal convergence.)

\begin{enumerate}[wide,label=(\alph*),labelindent=0pt]
\item %\noindent(a) 
{\em Poisson input}. 
When  $\tP$ is either a Poisson or binomial point process equipped with independent marks and when $\xi$ is exponentially stabilizing in the classical sense, i.e., when the decay functions  $\varphi_p$ in Definition~\ref{def.stabilizing_marking}~\ref{i.stabilizing}  and \eqref{stab-infinite}  are exponential functions, and when
$\Var{\mu_n^\xi(f)} = \Theta(n)$,
then
one obtains
in Theorem \ref{t:clt_linear_marks}, 
Theorem \ref{t:cltmarkedpp}, and Proposition~\ref{t:clt_linear_marks_new}
a rate of normal approximation in the Kolmogorov distance of $O(n^{-1/2})$ as shown in \citet{Lachieze2019normal}.  When 
When  $\tP$ denotes Poisson input and when the scores  are fast $\BL$-localizing and satisfy a fifth moment condition, then subject to the variance growth condition $\Var{\mu_n^\xi(f)} = \Theta(n)$, \citet{TY}
obtain rates of normal approximation for $(\mu_n^\xi(f))_{n \in \N}$ in the Kolmogorov and Wasserstein distances of order 
$O(n^{-1/2})$.  This extends the corresponding results of \citet{Lachieze2019normal} which is restricted to the setting of exponential stabilization of scores with i.i.d. marks. 

\item %\noindent(b) 
{\em Non-Poisson input}.
 In the case of general non-Poisson $\P$ with i.i.d. marks and when the variance $\Var{\mu_n^\xi(f)}$ is of volume order, i.e., when  $\Var{\mu_n^\xi(f)} = \Theta(n)$, then 
\citet{Fenzl2019asymptotic} refines the cumulant techniques of \cite{BYY19} and uses the general results from \citet{saulis1991limit} (see also \citet{doring2022method}) to deduce rates of normal convergence and also moderate deviations.  In their recent paper, \citet{cong2022convergence} quantify the rate of normal
  approximation for scores on an independently marked ground point process whenever the ground point process is 
  $\beta$-mixing, which is slightly stronger than 
  the fast mixing of Definition \ref{d.omegamixing}.
   The work of ~\citet{austern2018limit} can
  also be used to deduce normal approximation in terms of mixing
  coefficients.  \citet{dinh2021quantitative} have derived rates of convergence for linear statistics of stationary $\alpha$-determinantal point processes via a quantitative Marcinkiewicz theorem and 
  \citet{CRX2020} also deduce rates of convergence for
`locally defined' statistics using a coupling of the stationary point
process with its Palm version. In the case of Gibbs processes, this approach has been used to obtain more explicit bounds via `disagreement coupling'; see \cite{hirsch2025normal}. However it is unclear how to extend these latter results to  sums of localizing or stabilizing statistics of general marked point processes. 
\end{enumerate}

\item (Moment conditions.)   The moment conditions in Definitions  \ref{d.moment}  and   ~\ref{d.moment-function} hold for a range of input and scores.   When the  input 
is either determinantal, permanental, or the zero set of a Gaussian entire function we refer to Sections 2.1 and 2.2 of \cite{BYY19} for sufficient conditions insuring  the moment bound
\eqref{e:xinpmom} of Definition \ref{d.moment}.

\item(Variance asymptotics.) 
The proof of variance and covariance asymptotics in \eqref{eqn:var} and \eqref{eqn:Cov} actually only requires \BL-localization of \(\xi\) for \(p \in \{1, 2\}\) and the following fast \BL-mixing correlations of order 2 uniformly on \(\tP_n, \, n \in \N\). Specifically, the marks \(U_n(x_i) := \xi(\tx_i, \tP_n)\) must satisfy the bound:
\begin{equation}\label{xi-decor-2}
\left|\sE_{\bk x1{2}}[h_1(\xi(\tx_1, \tP_n)) h_2(\xi(\tx_2, \tP_n))] \rho^{(2)}(\bk x1{2}) - \sE_{x_1}[h_1(\xi(\tx_1, \tP_n))] \sE_{x_2}[h_2(\xi(\tx_2, \tP_n))] \rho^2 \right| \le \omega_2(|x_1 - x_2|),
\end{equation}
for all  \(h_1, h_2 \in \BL(\mathbb{R})\) and a fast-decreasing function \(\omega_2\); see  \eqref{e.variance-limit-BL} below. The summable exponential $\B$-mixing correlations of \(\tP\), assumed in Proposition~\ref{expvar} for presentational simplicity, ensures the fulfillment of this condition, as per Theorem~\ref{t:mppmix_double}\ref{i.BL-iterated-marks} and Corollary~\ref{c:mppmixcl_double}.

\item (Variance growth rates.) The variance lower bound assumption in Theorem \ref{t:clt_linear_marks}, 
Theorem \ref{t:cltmarkedpp}, and Proposition \ref{t:clt_linear_marks_new} is  crucial.  If $\sigma^2(\xi) > 0$,  then $\nu = 1$ in Theorem \ref{t:cltmarkedpp} and Proposition \ref{t:clt_linear_marks_new}. 
If $\sigma^2(\xi) = 0$, then $\Var{\mu_{n}^{\xi}(f)}$ could depend on the properties of $f$; see \citet[Section 1.5]{Nazarov12}.  Even in the simple case  $\xi \equiv 1$,  it is possible that for some $f$, $\Var \,\mu_{n}^{\xi}(f) \to 0$ arbitrarily fast; see for example \citet{gabrielli2008tilings,Nazarov12,mastrilli2024estimating,lachieze2024rigidity}.  However for sets $W \subset W_1$ of positive measure, one has $\Var \,\hat{\mu}_{n}^{\xi}(\mathbf{1}_W) = \Omega(n^{\nu})$ with $\nu \in \{1, \frac{d-1}{d}\}$ (see \cite[(3.8)]{KY20}) with explicit asymptotics in case of `nice' sets $W$; see \cite{sodin2023random,KY20,jalowy2025box,BYY19}. These provide the necessary variance lower bounds for Proposition \ref{t:clt_linear_marks_new} when $f = 1_W$, $W \subset W_1$ of positive measure. 
In the case of Poisson input $\P$ and exponentially stabilizing score functions $\xi$, Schulte and Trapp \cite{SchulteTrapp2024} provide conditions under which 
$\sigma^2(\xi) > 0.$  They survey previous work which  establishes  lower bounds for $\liminf_{n \to \infty} \Var{\mu_{n}^{\xi}(f)}$ and introduce a condition implying the strict positivity thereof,  and which amounts to checking that 
$D_z \sum_{x \in \P_n} \xi( \tx, \tP_n ), z \in W_n$, 
is positive on a sets with strictly positive probability, with $D_z$ being the difference operator as in \eqref{eqn:FME_Kernels1}.

\item (Covariance growth rates.) Observe that \eqref{eqn:Cov} provides covariance asymptotics only of  volume-order scale. If $\sigma^2(\xi) = 0$, it is possible to have different scaling for covariance depending on the test functions. Again for indicator functions, one can derive surface-order covariance asymptotics for $\hat{\mu}_{n}^{\xi}$; see for example \cite[Theorem 3.1]{KY20}, \cite[Theorem 1.3]{sodin2023random} and \cite{jalowy2025box} and the discussions therein. We refer to \cite{KY20} for more  details concerning the possible variance and covariance asymptotics based on the properties of $f$.

\item (Binomial input.)  We have assumed that $\P_n$ is the restriction of a point process $\P$ to the window $W_n$. 
This assumption facilitates establishing expectation and variance asymptotics, but it is not necessary when proving central limit theorems, provided that the family of marked point processes $(\tilde{\P}_n)_{n \in \N}$
satisfy $\hat{\omega}$-mixing correlations as in Definition \ref{d.omegahmixing}(iv) and provided the scores are fast \BL-localizing on $(\tilde{\P}_n)_{n \in \N}$ as in Definition \ref{def.Lp-stabilizing_marking}.  When $\P_n$ is the realization of $n$ i.i.d. uniform random variables on $W_n$ then the $k$-point correlations functions 
satisfy mixing, though with an extra error of order $O(n^{-1})$, bringing extra error terms in the subsequent analysis. A thorough 
study of the limit theory for dependently marked binomial input constitutes a separate project.

\end{enumerate}

\section{Proofs of results on  mixing  correlations---Theorem ~\ref{t:mppmix_double}}
\label{s.proof-decay}
In this section we prove Theorem~\ref{t:mppmix_double}, our main result on asymptotic decorrelation of marked point processes. Section \ref{ss.FME} introduces
the Factorial Moment Expansion (FME) for marked point processes and Section \ref{s.proofcorrstabmarks} applies the FME to prove  Theorem~\ref{t:mppmix_double}.
These results go a long way towards showing the factorization of mixed moments as at
\eqref{eqn:mixedmoment}.

\subsection{Factorial moment expansions of
  marked point processes}
\label{ss.FME}
In this subsection, we first provide conditions ensuring the factorial moment
expansion (FME), introduced in
\citet{Blaszczyszyn95} and~\citet{Bartek97}, for expectations of functions of point
processes with marks taking values in a Polish space.  We then apply this general result to obtain the FME  of  a {\em bounded} function $f$ of a $p$-tuple of score functions, each of which has a deterministically {\em bounded} radius of stabilization. 
        
Equip $\mR^d$ with a total order $\prec$ defined using the lexicographical ordering of the polar co-ordinates. A useful property of this order  is that $\{y : y \prec x\} \subset B_{|x|}(\0)$ and thus such sets are bounded. For
        $\tmu \in \cNK$ and $x \in \mR^d$, define the measure
	$\tmu_{|x}(\cdot) := \tmu(. \cap (\{y : y \prec x\}\times\K))$. Note that since $\tmu$
	is a locally finite measure and the ordering is defined via polar
	co-ordinates,  $\tmu_{|x}$ is a finite measure for all $x \in
	\mR^d$. Let $o$ denote the null-measure on $\mR^d\times\K$, i.e., $o(B) = 0$ for all
	Borel subsets $B$ of $\mR^d\times\K$.	For a measurable function $\psi : \cNK \to
	\mR$, $l \in \N$,
	and  $\bk {\ty}1l=
\bk{y, w}1l$, with distinct locations $\bk y1l \in
        (\R^d)^l$   and marks $\bk u1l \in \K^l$, we define the FME kernels
        \cite{Blaszczyszyn95,Bartek97} as follows. For $l \in \N$,
		\begin{align}
	\label{eqn:FME_Kernels1}
	D^{l}_{\bk {\ty}1l} \psi(\tmu) & = \sum_{i=0}^l(-1)^{l-i}\sum_{J \subset \binom{[l]}{i}}\psi(\tmu_{|y_*}
	+ \sum_{j \in J}\delta_{\ty_j}) \\
	& = \sum_{J \subset [l]}(-1)^{l-|J|}\psi(\tmu_{|y_*} + \sum_{j \in J}\delta_{\ty_j}), \nonumber 
	\end{align}	
where $\binom{[l]}{i}$ denotes the collection of all subsets of $[l]: =
\{1,\ldots,l\}$ with cardinality $i$, $y_* := \min
	\{y_1,\ldots,y_l\}$ with  the minimum  taken with respect to
	the order $\prec$, and $w_*=w_j$ with (unique) $j\in\{1,\ldots,l\}$ such that  $y_j=y_*$.	For $l = 0$, put $D^{0}\psi(\tmu) := \psi(o)$. Note that
	$D^{l}_{\bk {\ty}1l} \psi(\tmu)$ is invariant with
        respect to any joint permutation of marked points in $\bk{\ty}1l$. For $y_l \prec y_{l-1} \prec \ldots \prec y_1$  the functional $D^{l}_{\bk{\ty}1l}\psi(\tmu)$
		is equal to the {\em $l$th order iterated difference operator} (see \cite[Chapter 18]{last2017lectures}):
$$D^1_{\ty_1} \psi(\tmu)= 
\psi(\tmu_{|y_1} + \delta_{\ty_1})- \psi(\tmu_{|y_1}),  \quad D^{l}_{\bk{
\ty}1l}\psi(\tmu)= D^1_{\ty_l}(D^{l-1}_{\bk{\ty}1{l-1}}  \psi(\tmu) ).
$$
We say that $\psi$ is \textit{$\prec$-continuous at $\infty$} if for all $\tmu \in \cN$ we have
$$
\lim_{y \uparrow \infty }  \psi(\tmu_{|y} ) = \psi(\tmu).
$$

We shall require the following 
expansion for the expectation of  functionals of marked point processes. 
Recall from Section \ref{s:notation} that $\sE^{!}_{\bk {\ty}1l}$ denotes expectations with respect  to reduced Palm distributions of $\tP$ given  points   $\bk y1p$ and their marks $\bk w1p$.
\begin{theorem}[\cite{Blaszczyszyn95,Bartek97}; FME for functionals of marked point processes]
\label{FMEthm}
Let $\tP$ be a simple, marked point process on $\R^d\times\K$ and let $\psi : \cNK \to \mR$ be $\prec$-continuous at $\infty$.  Assume that for all $l \in \N$
\be \label{hyp1}
\int_{(\mR^{d}\times\K)^l}\E^!_{\bk {\ty}1l}[| D^{l}_{\bk {\ty}1l}\psi(\tP)|]\;\tilde{\alpha}^{(l)}(\md \bk {\ty}1l) < \infty
\ee
and
\be \label{hyp2}
		\frac{1}{l!}\int_{(\mR^{d}\times\K)^l}\E^!_{\bk {\ty}1l}[ D^{l}_{\bk {\ty}1l}\psi(\tP)]\;\tilde{\alpha}^{(l)}(\md \bk {\ty}1l) \to 0 \ {\rm{as}} \ l \to \infty.
\ee
Then $\E[\psi(\P)]$ has the following {\em factorial moment expansion} (FME)
\be
\label{FME}
\E[\psi(\tP)] = \psi(o) +
\sum_{l=1}^{\infty}\frac{1}{l!}\int_{
(\mR^{d}\times\K)^l }D^{l}_{\bk {\ty}1l}\psi(o)\;\tilde{\alpha}^{(l)}(\md \bk {\ty}1l).
\ee
\end{theorem}
The expansion \eqref{FME} is established in 
\cite[Theorem 3.1]{Bartek97} for $\K$ a locally compact, second countable
Hausdorff space and even earlier for unmarked point processes  in \cite[Theorem
3.2]{Blaszczyszyn95}.  The assumption on $\K$ in \cite{Bartek97} is restrictive whereas  Theorem
\ref{FMEthm} allows $\K$ to be Polish.  Indeed, the extension follows as the proof of \cite[Theorem 3.1]{Bartek97} depends only on the crucial existence and use of Palm distributions for marked point processes in this more general setting; see Section \ref{s:notation}.  Hence,  one can follow the proof of \cite[Theorem 3.1]{Bartek97} to obtain Theorem \ref{FMEthm} and so we shall not provide the proof here. The above FME also follows from \cite[Theorem A.1]{klatt2025invariant} but under (slightly) stronger conditions; see Remark \ref{Rmkconvsum}.

Consider now a marking function  $\xi$ and a $p$-tuple of distinct points $\bk
x1p\in (\R^d)^p$ with their marks $\bk u1p\in\K^p$.  Let
$\xi(\tx;\tmu)$ be a marking function with values in $\K'$ (a Polish space),  defined for
$\tx = (x,u)\in\tmu$, with $\tmu \in \cN_{\R^d \times \K}$.  
The proof of  Theorem~\ref{t:mppmix_double}  is based on 
the FME expansion for
\begin{equation}\label{eq:for-Palm-FME}
\E_{\bk x1p}[f( \xi((x_1,U_1),\tP),...,\xi((x_p,U_p), \allowbreak\tP))]\rho^{(p)}(\bk x1p),
\end{equation}
where $f \in \B((\K')^p).$
  Under $\mP_{\bk x1p}$ the point process $\tP$
     has  atoms at fixed  locations $x_1,\ldots,x_{p}$ with random marks
    $U_1,\ldots,U_p$, which complicates the form of its factorial
  moment measures.  We address this by considering these fixed
  locations and their marks as parameters of the following  modified functional
 
\begin{equation} \label{e.psip}
	\psi_{f}(\bk\tx1p;\tmu)
		 := f\Bigl( \xi(\tx_1, \tilde \mu + \sum_{i = 1}^p \delta_{\tx_i}),...,\xi(\tx_p, \tilde \mu+ \sum_{i = 1}^p \delta_{\tx_i})\Bigr)
	\end{equation}
and to not count these marked points $\tx_1,\ldots,\tx_{p}$ in $\tP$, i.e., to consider $\tP$ under the reduced Palm probabilities $\mP^!_{\bk \tx1p}$.  Obviously 
$$ \E_{\bk  \tx1p}[f( \xi(\tx_1, \tP),...,\xi(\tx_p, \tP))] = \E^!_{\bk \tx1p}[\psi_f(\bk {\tx}1p ;\tP)]. $$ 
For all 
$p \in \N$ recall that    $\kappa_p:=\sup_{\bk x1p\in(\R^d)^p}  \rho^{(p)}(\bk x1p)$ is the bound on the $p$th correlation function of $\P$.
To justify convergence of some upcoming series expansions, we shall need to assume that the ground process
$\P$ satisfies
\begin{equation}\label{e.kappa-bound-fine}
\lim_{l\to\infty}   \frac{w^l}{l!} \kappa_{k+l}=0\quad
\end{equation}
for certain values of $k \in \N, w > 0$. 

Now,  we state the following important consequence of Theorem~\ref{FMEthm}.  As usual, by $\mu$ we denote the natural projection of $\tmu$ on $\R^d$. 	Let $\theta_d:=\pi^{d/2}/\Gamma(d/2 + 1)$ denote the volume of the unit ball in $\R^d$. We recall for fixed $\bk x1p\in(\R^d)^p$,  the  Palm probability
        distribution $\cM_{\bk x1p}^{(p)}(\md \bk u1p)$
  on $\K^p$ of the random marks $\bk U1p=(U_1,\ldots,U_p)$
 respectively at the  fixed locations $\bk x1p$, is defined via the identity
  $$
  \tilde \alpha^{(p)}(\md \bk{\tx}1p)=\cM_{\bk x1p}^{(p)}(\md \bk u1p)\rho^{(p)}(\bk x1p)\, \md (\bk x1p),
  $$
 where $\tilde \alpha^{(p)}$ is the $p\,$th order factorial moment measure of $\tilde\P$. In the following expansion, the distribution 
 $\cM_{\bk x1p}^{(p)}(\md \bk u1p)$ incorporates  $\xi.$ 

Recall \(
B_r(\bk{x}{1}{p}) := \bigcup_{i=1}^p B_r(x_i) \subset \R^d,
\)
is the union of balls of radius $r>0$ centered at the points~$\bk{x}{1}{p}$.
\begin{lemma}[FME of a bounded function of a $p$-tuple of union-restricted marking  functions] \label{l.FMEf} 
Consider  a marking  function $\xi$ defined on  ~$\tP\in\cN_{\R^d \times \K}$. Assume $\tP$ satisfies~\eqref{e.kappa-bound-fine} for $k=p$ and $w = 2p\theta_dr^d$  for some  $r>0$.
Then for any $f \in \B((\K')^p)$,
the functional $\psi_f(\bk\tx1p;\tP\cap B_r(\bk{x}{1}{p}))$ at~\eqref{e.psip} admits  for $\alpha^{(p)}$ a.e. $\bk x1p \in (\R^{d})^p$ the following FME
				\begin{align}\nonumber
                  & \quad \E_{\bk x1p} \big[ f( \xi((x_1,U_1),
  \tP \cap B_r(\bk{x}{1}{p})),...,\xi((x_p,U_p), \allowbreak\tP \cap  B_r(\bk{x}{1}{p}))) \big] \, \rho^{(p)}(\bk x1p)
\\		&= \int_{\K^p}\psi_f(\bk {\tx}1p;o)\,\cM_{\bk
  x1p}^{(p)}(\md \bk u1p) \rho^{(p)}(\bk x1p)\nonumber\\
		&\;+\sum_{l=1}^{\infty}\frac{1}{l!}\int_{(B_r(\bk{x}{1}{p})\times \K)^l}
 D^{l}_{\bk {\ty}1l}\psi_f(\bk {\tx}1p;o)\,
		\cM^{(p+l)}_{\bk x1p\sqcup\bk y1l}(\md (\bk u1p\sqcup\bk w1l))\rho^{(p+l)}(\bk
                  x1p\sqcup \bk y1l)\,\md (\bk y1l)\,
		\label{FMEpsip}
		\end{align}
where  $\bk {\tx}1p=\bk {x,u}1p$, $\bk {\ty}1l=\bk {y,w}1l$.
\end{lemma}

\begin{proof} 
Throughout we fix  $f \in \B((\K')^p)$.
The bound $|f| \leq 1$  makes our proof simpler
than the corresponding result in the unmarked case, which does not use this assumption; see \cite[Lemma 3.2]{BYY19}.  We prove  Lemma~\ref{l.FMEf} by first establishing the validity of the FME at \eqref{FME} for $\E^!_{\bk \tx1p} \big[ \psi_f(\bk {\tx}1p;\tP \cap B_r(\bk{x}{1}{p})) \big] \, \rho^{(p)}(\bk x1p)$, with $\bk \tx1p=\bk{x,u}1p$ fixed. We then use the resulting FME to derive  \eqref{FMEpsip}.
  
Indeed, restricting $\xi$ to  $\tP\cap B_r(\bk{x}{1}{p})$ implies $\psi_f$
is $\prec$-continuous at~$\infty$.  
Regarding the integral assumptions \eqref{hyp1} and
\eqref{hyp2}
(with $\psi$ replaced by  $\psi_f(\bk {\tx}1p
    ;\tP \cap B_r(\bk{x}{1}{p})) \, \rho^{(p)}(\bk x1p)$, $\E^!_{\bk{\ty}1l}$ replaced  by
      $(\E^!_{\bk {\tx}1p})^!_{\bk{\ty}1l}=\E^!_{\bk {\tx}1p
   \sqcup\bk {\ty}1l}$ and  $\tilde{\alpha}^{(l)}(\md\bk {\ty}1l)$ by 
$\tilde\alpha^{(l)}_{\bk {\tx}1p}(\md \bk {\ty}1l)$), we  observe
    \begin{align}\nonumber
   & \int_{(\R\times\K)^l} \E^!_{\bk {\tx}1p
   \sqcup\bk {\ty}1l} \big[ |  D^l_{\bk{\ty}1l}\psi_f(\bk {\tx}1p;\tP \cap B_r(\bk{x}{1}{p})) | \big] \,\tilde{\alpha}^{(l)}_{\bk{\tx}1p}(\md\bk{\ty}1l)\,\rho^{(p)}(\bk
      x1p) \\
    & \le \int_{(B_r(\bk{x}{1}{p})\times\K)^l} 2^l\,\tilde{\alpha}^{(l)}_{\bk{\tx}1p}(\md\bk{\ty}1l)\,\rho^{(p)}(\bk
      x1p)\nonumber\\
      &= 2^l \int_{B_r(\bk{x}{1}{p})^l }\rho^{(p+l)}(\bk
                  x1p\sqcup \bk y1l)\,\md \bk y1l\nonumber\\
    &\le (2p \theta_d r^d)^l \kappa_{p+l}.\label{e:DL-bound}
   \end{align}
    The first inequality is justified since  {\em the difference operator $D^{l}_{\bk{\ty}1l}$ vanishes} as soon as $y_k \notin
B_r(\bk{x}{1}{p})=\bigcup_{i=1}^pB_r(x_i)$ for some $k \in \{1,\ldots,l\}$, that is to say $D^{l}_{\bk{\ty}1l}\psi_f(\bk{\tx}1p;\tmu \, \cap \, B_r(\bk{x}{1}{p}))= 0$. 
 On the other hand, when  $\{y_1,\ldots,y_l\} \subset
\bigcup_{i=1}^pB_r(x_i)$, we have
$D^{l}_{\bk{\ty}1l}\psi_f(\bk{\tx}1p;\tmu)\le 2^l$, using $|\psi_f(\bk{\tx}1p;\tmu)|\le 1$ for $|f|\le 1$
and~\eqref{eqn:FME_Kernels1}.
The equality follows by noting that integration with respect to  $\bk w1l$ in \eqref{e:palmalgebra} yields
$$ \int_{\K^l}  \tilde{\alpha}^{(l)}_{\bk{\tx}1p}(\md\bk{\ty}1l)\,\rho^{(p)}(\bk
      x1p)=\rho^{(p+l)}(\bk x1p\sqcup \bk y1l)\, $$
Finally, the last inequality in \eqref{e:DL-bound}
follows from the  bound  $\rho^{(p+l)}(\bk x1p\sqcup \bk y1l)\le  \kappa_{p+l}$ which holds uniformly in  
$\bk x1p$ and $ \bk y1l$, in view of  Assumption \ref{Ass1}.

The bound~\eqref{e:DL-bound},  together with the assumed finiteness of  $\kappa_{p+l}$ 
and the assumed validity of ~\eqref{e.kappa-bound-fine} for $k=p$, $w = 2p\theta_dr^d$,  justifies~\eqref{hyp1} and ~\eqref{hyp2}. Thus,  the FME at 
\eqref{FME} holds for $\E^!_{\bk \tx1p}[\psi_f(\bk {\tx}1p;\tP \cap B_r(\bk{x}{1}{p}))]\rho^{(p)}(\bk x1p)$.

Now we establish  \eqref{FMEpsip}.  First,  observe that 
 \begin{align}
           & \E_{\bk x1p}[f( \xi((x_1,U_1),
              \tP \cap B_r(\bk{x}{1}{p})),...,\xi((x_p,U_p),\tP \cap B_r(\bk{x}{1}{p})))]\rho^{(p)}(\bk x1{p}) \nonumber \\
            & \quad = \int_{\K^p}\E^!_{\bk \tx1p}[\psi_f(\bk {\tx}1p
     ;\tP \cap B_r(\bk{x}{1}{p}))]\,\cM_{\bk x1p}^{(p)}(\md \bk u1p)\rho^{(p)}(\bk
              x1{p})\,.    \label{e.inergral-psi}
   \end{align}
  The  terms in ~\eqref{FME} for $\E^!_{\bk \tx1p}[\psi_f(\bk {\tx}1p
      ;\tP \cap B_r(\bk{x}{1}{p}))]$ are integrals with respect to the $l\,$th order factorial moment measures $\tilde\alpha^{(l)}_{\bk {x,u}1p}(\md (\bk {y,w}1l))$ of $\tP$ under
          $\E^!_{\bk {\tx}1p}$. The subsequent          
          integration of these terms with respect to   $\cM_{\bk x1p}^{(p)}(\md \bk
          u1p)\rho^{(p)}(\bk x1p)$,  via~\eqref{e:palmalgebra}, yields the integrals in 
           \eqref{FMEpsip} with respect to the  measures
          $\cM_{\bk x1p \sqcup \bk y1l}^{(p+l)}$.  Indeed, the term-wise integration of  the FME series of $\E^!_{\bk \tx1p}[\psi_f(\bk {\tx}1p ;\tP \cap B_r(\bk{x}{1}{p}))]$ with respect to $\cM_{\bk x1p}^{(p)}(\md \bk u1p)\rho^{(p)}(\bk x1{p})$ is justified by the uniform bound~\eqref{e:DL-bound} with respect to the marks $\bk u1p$ and 
          hence the `remaining term'~\eqref{hyp2} of the FME expansion for~\eqref{e.inergral-psi} goes to zero by the assumption~\eqref{e.kappa-bound-fine}. 
           This  justifies  \eqref{FMEpsip} and completes the proof of Lemma~\ref{l.FMEf}. 
\end{proof}		

\subsection{Proof of mixing correlation properties of point processes with iterated marks} %--- Theorem~\ref{t:mppmix_double}
\label{s.proofcorrstabmarks}

Our goal in this subsection is to prove  Theorem~\ref{t:mppmix_double}. Recall that parts (i) and (ii), respectively part (iii), of this theorem states that summable exponential $\B$-mixing correlations of $\tP = \{(x,U(x))\}_{x \in \P}$ guarantees  fast $\BL$-mixing, respectively fast $\B$-mixing,  correlations of $\dbtilde \P = \{ (x,\xi(x,\tP)) \}_{x \in \P}$ provided the functionally constructed marking function $\xi$ satisfies fast $\BL$-localization or \BL~cluster-localization, respectively fast stopping set stabilization. The proof has three main steps:  
\begin{enumerate}
\item First, we establish in Lemma~\ref{p:mppmixcl_double} 
a mixing correlation bound \eqref{e.hatC-bound-fine-cl} for  union-restricted marking  functions  when $\tP$ has $\B$-mixing correlations. For this, we use the FME  at \eqref{FMEpsip} to prove the mixing correlation bound. As a consequence, assuming that $\xi$ is a short-range mark,
i.e.\
$\xi=\xi^{(r)}$ for some $r\in(0,\infty)$,
and that $\tP$ has $\B$-mixing correlations, we show that the marked point
process $\dbtilde{\P}$ inherits $\B$-mixing correlations of the same type.
This is established in Proposition~\ref{p:mppmix_double}; see
\eqref{e.hatC-bound-fine} below. In both  statements, the decay constants and
the decay function are given explicitly.

\item  When the point process $\tP$ exhibits summable 
exponential  $\B$-mixing correlations,
we demonstrate that $\dbtilde \P$ also has fast 
$\B$-mixing correlations. 
The assumption 
that \dy{$\xi$} is a short-range mark $\xi = \xitr$ for some $r \in (0, \infty)$ is
relaxed to a fast stopping set stabilization assumption. To achieve this, we leverage \eqref{e.hatC-bound-fine} in combination with a truncation approach to exploit these stabilization properties of~$\xi$. This constitutes the proof of part~\ref{i.iterated-marks} of Theorem~\ref{t:mppmix_double}.

\item If furthermore fast stabilization of  $\xi$ is relaxed to fast  $\BL$-localization or fast \BL~cluster-localization, 
then $\dbtilde \P$ exhibits only fast $\BL$-mixing correlations. This establishes the proof of parts~\ref{i.BL-iterated-marks} and~\ref{i.BLc-iterated-marks} of Theorem~\ref{t:mppmix_double}, where bound~\eqref{e.hatC-bound-fine-cl} is crucial for part~\ref{i.BLc-iterated-marks}.
\end{enumerate}

\subsubsection{Proof of Theorem~\ref{t:mppmix_double} for restricted marking  functions}
\label{ss.Proof-bounded-radius}
 
We recall the notation for the \emph{union-restriction operation}
introduced in~\eqref{e:xi-cl}. Given a marking function~$\xi$ and points
$x_1,\ldots,x_p \in \R^d$, we define
\[
\bk{{\xi}^{\cup(r)}}{1}{p}(\tx,\tP)
:=
\bigl(
\xi(\tilde x_1,\tP \cap B_r(\bk{x}{1}{p})),
\ldots,
\xi(\tilde x_p,\tP \cap B_r(\bk{x}{1}{p}))
\bigr),
\]
where \(
B_r(\bk{x}{1}{p}) := \bigcup_{i=1}^p B_r(x_i)
\)
denotes the union of balls of radius $r>0$ centered at the points $\bk{x}{1}{p}$.

\begin{lemma}[A mixing correlation bound for union-restricted marks]
\label{p:mppmixcl_double}
Fix $p,q \in \N$ and $r \in (0, \infty)$.
Assume that we are given
\begin{enumerate}[wide,label=(\roman*),labelindent=0pt,partopsep=0pt,topsep=0pt,parsep=0pt]
\item \label{iii:mppmix_double} a marking function $\xi: \R^d \times \K \times \cN_{\R^d \times \K} \to \K'$, and
\item \label{ii:mppmix_double} 
a marked point process $\tP$ on $\R^d\times\K$ having  $\B$-mixing  correlations
as in  Definition \ref{B-omega-mixing}\ref{Bmixing}
with correlation decay functions $\omega_{k}(s) = C_k  \phi(s)$, $k \in \N$, 
where $(C_k)_{k \geq 1}$ is a non-decreasing sequence such that for some fixed $\epsilon \in (0,1)$, we have
\be \label{weaksum} 
\sum_{l = 1}^{\infty}  \frac{ ({w}_0(1 + \epsilon)\ka_0^{1 + \epsilon})^l }{l!} C_{l + p + q}  < \infty,
\ee
with  ${w}_0= 4(p+q)\theta_d r^d$ and where $\kappa_0 = \max\{1, \sup_{x\in\R^d}\rho^{(1)}(x)\}$ is related to the ground process $\P$.
\end{enumerate}
\vskip.2cm
Then for  all $f \in \B((\K')^p),$   $g \in \B((\K')^q)$ and for almost all $\bk x1{p+q}\in\R^{d(p+q)}$ such that $s:=d(\bk x1p,\bk x{p+1}{p+q}) >4r$ we have
\begin{align} \nonumber
& \Big| \sE_{\bk x1{p+q}}[f(\bk{{\xi}^{\cup(r)}}{1}{p}(\tx,\tP))g(\bk{{\xi}^{\cup(r)}}{p+1}{p+q}(\tx,\tP))]\rho^{(p+q)}(\bk x1{p+q})  \\
&  \qquad - \,
\sE_{\bk x1p}[f(\bk{{\xi}^{\cup(r)}}{1}{p}(\tx,\tP))] \rho^{(p)}(\bk x1{p}) \; \sE_{\bk x{p+1}{p+q}}[g(\bk{{\xi}^{\cup(r)}}{p+1}{p+q}(\tx,\tP))]\rho^{(q)}(\bk x{p+1}{p+q}) \Big|  \nonumber \\
& \le \label{e.hatC-bound-fine-cl}
\phi( \frac{s}{2})\sum_{l=1}^\infty \frac{w_0^l}{l!}  C_{l+p+q}
=:\phi( \frac{s}{2})\hat{C}_{p+q}<\infty.
\end{align}
\end{lemma}
 Before proving the above proposition, we derive a simple but useful consequence
for marked point processes in which new (iterated) marks are functionally constructed via
a restricted marking function. 
\begin{proposition}[Mixing correlations of point processes with short-range marks]
\label{p:mppmix_double}
Fix $p,q \in \N$ and $r \in (0, \infty)$. Assume that $\tP$ satisfies condition~\ref{ii:mppmix_double}  of Lemma~\ref{p:mppmixcl_double}. Let $\xi: \R^d \times \K \times \cN_{\R^d \times \K} \to \K'$ be a marking function such that $\xi(\tx,\tmu) =  \xitr(\tx,\tmu) =\xi(\tx,\tmu\cap B_r(x))$ for all $(\tx,\tmu) \in \R^d \times \K \times \cN_{\R^d \times \K}$. Consider the marked point process $\dbtilde\P=\sum_{x \in \P} \delta_{(x,V(x))}$ on $\R^d\times \K'$ where $V(x) := \xi(\tx,\tP)$.  

Then, for  all $f \in \B((\K')^p),$   $g \in \B((\K')^q)$, for almost all $\bk x1{p+q}\in\R^{d(p+q)}$ such that $s:=d(\bk x1p,\bk x{p+1}{p+q}) >4r$, and with ${w}_0= 4(p+q)\theta_d r^d$
we have
\begin{align} \nonumber
& \Big| \sE_{\bk x1{p+q}}[f(\bk V1{p})g(\bk V{p+1}{p+q})]\rho^{(p+q)}(\bk x1{p+q}) -
\sE_{\bk x1p}[f(\bk V1p)] \rho^{(p)}(\bk x1{p}) \; \sE_{\bk x{p+1}{p+q}}[g(\bk V{p+1}{p+q})]\rho^{(q)}(\bk
x{p+1}{p+q}) \Big|
\\
&\le \label{e.hatC-bound-fine}
\phi( \frac{s}{2})\sum_{l=1}^\infty \frac{w_0^l}{l!}  C_{l+p+q}
=:\phi( \frac{s}{2})\hat{C}_{p+q}<\infty.
\end{align}
Thus the marked point process
 $\dbtilde\P$ satisfies the inequality~\eqref{e.Bmppmixing}  in Definition~\ref{B-omega-mixing}\ref{Bmixing} 
for $s>4r$, with  $\omega_{p+q}(s):=\hat C_{p+q}\phi(s/2)$ where  $p,q,\phi$ are as above.
\end{proposition}
\begin{proof}
 It is an easy observation that the union-restriction operation is neutral for a short-range mark $\xi=\xi^{(r)}$, for some $r>0$: for all $x_1,\ldots,x_p \in \R^d$,
\(
\bk{{\xi}^{\cup(r)}}{1}{p}(\tx,\tP)
=
\bk{\xi^{(r)}}{1}{p}(\tx,\tP)=\bk{\xi}{1}{p}(\tx,\tP).
\)
Therefore, the inequality~\eqref{e.hatC-bound-fine} follows directly from
\eqref{e.hatC-bound-fine-cl}, which completes the proof.
\end{proof}

We need another  immediate extension of the bound~\eqref{e.mppmixing} for the proof of Lemma~\ref{p:mppmixcl_double}.
\begin{corollary}[Extended version of $\B$-mixing correlations]
  \label{c.omegahmixing-mp}
If the  marked point process $\tilde \P$ 
satisfies $\B$-mixing
  correlations as in Definition \ref{B-omega-mixing}\ref{Bmixing} then the bound~\eqref{e.Bmppmixing}
  also holds with the same $\omega$ for all functions 
    $f=f(\bk {(x,u)}1p)$, $g=g(\bk {(x,u)}{p+1}{p+q})$
    in $\B((\R^d \times \K)^p)$ and $\B((\R^d \times \K)^q)$,
    respectively. %, with $|f|\leq 1$, $|g|\leq 1$
\end{corollary}

\begin{proof}[Proof of Lemma~\ref{p:mppmixcl_double}]
We will need to use the FME expansion of Lemma \ref{l.FMEf}.  We need to first verify that our assumptions, together with the bound ~\eqref{e.correlation-functions-bound} for $\ka_k$ and the non-decreasing sequence $(C_k)_{k \geq 1}$,  imply that 
$\P$ satisfies~\eqref{e.kappa-bound-fine} for all $k\in\{p,q,p+q\}$, $w = 2k\theta_dr^d$.  This is a straightforward exercise which is left to the reader. 

Now we are ready to apply Lemma \ref{l.FMEf}.
Again assume that $r >1$ without loss of generality. Consider  a $(p+q)$-tuple $\bk {\tx}1{p+k}=\bk {x,u}1{p+k}$
 of distinct points $\bk
x1{p+q}\in (\R^d)^{p+q}$ with marks $\bk u1{p+q}\in\K^{p+ q}$.
For given  functions $f,g$ on $(\K')^p$ and
$(\K')^q$, respectively, define  functionals of $\tmu\in\cNK$:
\begin{align}\label{e.psi-f}
\psi_{f}(\bk{\tx}1p;\tmu)&:=f(\bk v1p),\\
\psi_{g}(\bk{\tx}p{p+q};\tmu) &:=g(\bk v{p+1}{p+q}),\label{e.psi-g}\\
\psi_{fg}(\bk{\tx}{1} {p+q};\tmu)&:=f(\bk v1p)g(\bk v{p+1}{p+q}), \label{e.psi-fg}
\end{align}
where $v_i:=\xi(\tx_i; (\tmu+\sum_{\tx \in \bk{\tx}ab}\delta_{\tx}) \cap B_r(\bk xab))$,
with $a=1,b=p$ in~\eqref{e.psi-f}, $a=p+1,b=p+q$
in~\eqref{e.psi-g} and $a=1,b=p+q$ in~\eqref{e.psi-fg}.
 Note that the assumptions of Lemma~\ref{l.FMEf} are satisfied.
Hence we can represent each of 
the above three functionals as an
infinite sum of difference operators, as in the FME
representation ~\eqref{FMEpsip}.

Let $\bk x1{p+q}$ be such that $s:=d(\bk x1p,\bk x{p+1}{p+q})>4r$.   We apply the FME  \eqref{FMEpsip} to $\psi_{fg}(\bk{\tx}1{p+q};\tmu)$.  Using the definition of the FME kernels \eqref{eqn:FME_Kernels1} and that the FME kernels $D^l$ vanish outside $\bigcup_{i=1}^{p+q}B_r(x_i)$,  as noted in the proof of Lemma \ref{l.FMEf}, we may crucially replace the integration domain
$(\R^d)^l$ in \eqref{FMEpsip}, by a product of $B_r(\bk x1p)$ and $B_r(\bk x{p+1}{p+q})$. 

From the above vanishing property of difference operators, we obtain the following factorization of $\psi_{fg}$. Given $\bk{\ty}1j \subset B_r(\bk x1p)$ and $\bk{\ty}{j+1}l \subset B_r(\bk x{p+1}{p+q})$,
together with sets  $J_1\subset[j]$ and  $J_2\subset[l]\setminus[j]$, we have the factorizations
\[ \psi_{fg}(\bk{\tx}1{p+q};\textstyle{\sum_{i \in J_1 \cup J_2}\delta_{\ty_i}}) = \psi_f(\bk{\tx}1p;\textstyle{\sum_{i \in J_1}\delta_{\ty_i}}) \psi_g(\bk{\tx}{p+1}{p+q};\textstyle{\sum_{i \in J_2}\delta_{\ty_i}})\]
and 
$$ \sum_{J_1 \subset [j],J_2\subset[l]\setminus[j] }\!\!\!(-1)^{l-|J_1|-|J_2|}
\psi_f(\dots)\psi_g(\dots)= \sum_{J_1 \subset [j]}\!(-1)^{j-|J_1|}
\psi_f(\dots) \sum_{J_2\subset[l]\setminus[j] }\!\!\!(-1)^{l-j-|J_2|}\psi_g(\dots).
$$
Further, we use the definition of FME kernels $D^l$ as in \eqref{eqn:FME_Kernels1} to obtain 
\begin{align}
& \sE_{\bk x1{p+q}}[f(\bk{{\xi}^{\cup(r)}}{1}{p}(\tx,\tP))g(\bk{{\xi}^{\cup(r)}}{p+1}{p+q}(\tx,\tP))]\rho^{(p+q)}(\bk x1{p+q}) \no \\
& = \sum_{l=0}^{\infty}\sum_{j=0}^l\frac{1}{j!(l-j)!}\int_{(B_r(\bk x{1}{p}))^j \times(B_r(\bk x{p+1}{p+q}))^{l-j}}
		\int_{\K^{p+q+l}} \no \\
	& \no \hspace{2em} \times
	        \sum_{J_1 \subset [j]}(-1)^{j-|J_1|}\psi_f(\bk{\tx}1p;\textstyle{\sum_{i \in J_1}\delta_{\ty_i}})\\
	& \no \hspace{2em} \times
	 \sum_{J_2 \subset [l] \setminus[j]} (-1)^{l-j-|J_2|}
	 \psi_g(\bk{\tx}{p+1}{p+q};\textstyle{\sum_{i \in J_2}\delta_{\ty_i}}) \\
		& \hspace{4em} \times \cM_{\bk x1{p+q} \sqcup \bk y1l}^{(p+q+l)}(\md (\bk u1{p+q}\sqcup\bk {w}1l)) \,\rho^{(p+q+l)}(\bk x1{p+q}\sqcup\bk y1l)\, \md \bk y1l, \label{e.FME-fg}
\end{align}
where $\ty_i=(y_i,w_i)$ and where we emphasize that the marks depend also only on $\bk x1{p+q} \sqcup\bk y1l$. 
Similarly, 
\begin{align}
& \sE_{\bk x1{p}}[f(\bk{{\xi}^{\cup(r)}}{1}{p}(\tx,\tP))] \, \sE_{\bk x{p+1}{p+q}}[g(\bk{{\xi}^{\cup(r)}}{p+1}{p+q}(\tx,\tP))] \, \rho^{(p)}(\bk x1{p}) \, \rho^{(q)}(\bk x{p+1}{p+q})\, \no \\
	&\no =\sum_{l=0}^{\infty}\sum_{j=0}^l\frac{1}{j!(l-j)!}\int_{(B_r(\bk x{1}{p}))^j \times(B_r(\bk x{p+1}{p+q}))^{l-j}}
		\int_{\K^{p+j}} \int_{\K^{q+l-j}}\\
	&\no \hspace{2em}\times
	\sum_{J_1 \subset [j]}(-1)^{j-|J_1|}\psi_f(\bk{\tx}1p;\textstyle{\sum_{i \in J_1}\delta_{\ty_i}})\\
	&\no \hspace{2em}\times
	\sum_{J_2 \subset [l]\setminus[j]}(-1)^{l-j-|J_2|}
	\psi_g(\bk{\tx}{p+1}{p+q};\textstyle{\sum_{i \in J_2}\delta_{\ty_i}})\\
		&\hspace{4em}\times\cM_{\bk x1{p}\sqcup \bk y1j}^{(p+j)}(\md (\bk u1{p},\bk {w}1j))\rho^{(p+j)}(\bk x1{p},\bk y1j) \no
		\\[1ex]
		& \hspace{4em}\times	
		\cM_{\bk x{p+1}{p+q} \sqcup \bk y{j+1}l}^{(q+l-j)}(\md
                (\bk u{p+1}{p+q} \sqcup \bk {w}{j+1}l))\,
		\rho^{(q+l-j)}(\bk x{p+1}{p+q}\sqcup \bk y{j+1}l) \, \md \bk y1j\,\md \bk y{j+1}l. \label{e.FME-fxg}
\end{align}
These two derivations are the analogs of \cite[(3.20) and the display preceding (3.23)]{BYY19}, which were derived for  point processes which did not carry dependent marks. 

We bound the difference between the FME expressions \eqref{e.FME-fg}
and~\eqref{e.FME-fxg} as follows.  We match like terms in
the two expressions and bound their difference termwise
using the decorrelation~\eqref{e.mppmixing} of the underlying
  marked point process~$\tP$, observing that  
  $$
  d(\bk x1{p}\sqcup \bk y1j),\bk x{p+1}{p+q}\sqcup \bk  y{j+1}l) \geq d(\bk x1p,\bk x{p+1}{p+q}) -2r \geq
  d(\bk x1p,\bk x{p+1}{p+q}) - s/2 \geq s/2.
  $$

The details go as follows. For sets $J_1,J_2$ in~\eqref{e.FME-fxg} consider  functions
%\jy{(JY note to self:  these are not in BL)}\BB{BB: No, not %necessarily.}
\begin{align*}\hat f=\hat f_{J_1}=\hat f\bigl(\bk {(x,u)}1p\sqcup \bk {(y,w)}{J_1}{}\bigr) &:=\psi_f(\bk{\tx}1p;\textstyle{\sum_{i \in J_1}\delta_{\ty_i}})\\
\hat g=\hat g_{J_2}=\hat g \bigl( \bk {(x,u)}{p+1}{p+q}\sqcup \bk {(y,w)}{J_2}{} \bigr) &:=  \psi_g(\bk{\tx}{p+1}{p+q};\textstyle{\sum_{i \in     J_2}\delta_{\ty_i}}),
\end{align*}
with $\hat f$ defined  on $(\R^d \times \K)^{p+|J_1|}$ and evaluated for 
marks 
$\bk u1p$,  $\bk w{J_1}{}:=(w_i: i\in J_1)$
 and points $\bk x1p$,
$\bk y{J_1}{}:=(y_i:i\in J_1)$, respectively, and similarly for $\hat g$. 
In the following derivation,  the second inequality is due to~\eqref{e.mppmixing}
(in association with Corollary~\ref{c.omegahmixing-mp}) 
termwise  for
  $d(\bk x1{p} \sqcup \bk y1j,\bk x{p+1}{p+q}\sqcup\bk
  y{j+1}l) \geq s/2$. 
    We have 
\begin{align*}
&  \qquad \Big| \sE_{\bk x1{p+q}}[f(\bk{{\xi}^{\cup(r)}}{1}{p}(\tx,\tP))g(\bk{{\xi}^{\cup(r)}}{p+1}{p+q}(\tx,\tP))]\rho^{(p+q)}(\bk x1{p+q}) \nonumber \\
&  \quad \qquad - \,
\sE_{\bk x1p}[f(\bk{{\xi}^{\cup(r)}}{1}{p}(\tx,\tP))] \rho^{(p)}(\bk x1{p}) \; \sE_{\bk x{p+1}{p+q}}[g(\bk{{\xi}^{\cup(r)}}{p+1}{p+q}(\tx,\tP))]\rho^{(q)}(\bk x{p+1}{p+q}) \Big|  \nonumber \\
& \leq \sum_{l=0}^{\infty}\sum_{j=0}^l\frac{1}{j!(l-j)!}\int_{(B_r(\bk x{1}{p}))^j \times(B_r(\bk x{p+1}{p+q}))^{l-j}}
		 \no  \sum_{J_1 \subset [j]} \sum_{J_2 \subset                l]\setminus[j]}\\
& \quad \times \Bigg|\int_{\K^{p+q+l}}\hat f( \bk{(x,u)}1p \sqcup \bk {(y,w)}{J_1}{})\,\hat g(\bk {(x,u)}{p+1}{p+q}\sqcup \bk {(y,w)}{J_2}{})\\
&\hspace{0.2\linewidth}\times
        \cM_{\bk x1{p+q}\sqcup \bk y1l}^{(p+q+l)}(\md (\bk u1{p+q}\sqcup \bk {w}1l)) \,\rho^{(p+q+l)}(\bk x1{p+q}\sqcup \bk y1l) \no \\
&\quad-  \int_{\K^{p+j}} \hat f(\bk {(x,u)}1p \sqcup \bk
 {(y,w)}{J_1}{})\,
 \cM_{\bk x1{p}\sqcup\bk y1j}^{(p+j)}(\md (\bk
 u1{p}\sqcup\bk {w}1j))\rho^{(p+j)}(\bk x1{p}\sqcup\bk y1j)\no\\
 &\qquad\times \int_{\K^{q+l - j}}\hat g(\bk {(x,u)}{p+1}{p+q} \sqcup \bk {(y,w)}{J_2}{}) \cM_{\bk x{p+1}{p+q} \sqcup \bk y{j+1}l}^{(q+l-j)}(\md
                (\bk u{p+1}{p+q}\sqcup\bk {w}{j+1}l))\,
		\rho^{(q+l-j)}(\bk x{p+1}{p+q}\sqcup\bk y{j+1}l)\Bigg|\; \md \bk y1l \no \\
& = \sum_{l=0}^{\infty}\sum_{j=0}^l\frac{1}{j!(l-j)!}\int_{(B_r(\bk x{1}{p}))^j \times(B_r(\bk x{p+1}{p+q}))^{l-j}}
		 \no  \sum_{J_1 \subset [j]} \sum_{J_2 \subset
                   [l]\setminus[j]} \\
& \quad \times \Bigg| \sE_{\bk x1{p+q} \sqcup \bk y1l}\Big[\hat f(\bk {(x,u)}1{p} \sqcup \bk {(y,w)}{J_1}{}) \hat g(\bk {(x,u)}{p+1}{p+q} \sqcup \bk {(y,w)}{J_2}{})\Big]\rho^{(p+q+l)}(\bk x1{p+q} \sqcup \bk y1l) \\
& \qquad - \, \sE_{\bk x1p \sqcup \bk y1j}\Big[\hat f(\bk {(x,u)}1p \sqcup \bk {(y,w)}{J_1}{})\Big] \rho^{(p+j)}(\bk x1{p} \sqcup \bk y1j)  \\
&  \quad \qquad \times \, \sE_{\bk x{p+1}{p+q} \sqcup \bk y{j+1}l}\Big[\hat g(\bk {(x,u)}{p+1}{p+q} \sqcup \bk {(y,w)}{J_2}{})\Big]\rho^{(q+l-j)}(\bk
x{p+1}{p+q}\sqcup \bk y{j+1}l) \Bigg| \; \md \bk y1l \\
& \leq	 \phi( \frac{s}{2} ) \sum_{l=0}^{\infty}\sum_{j=0}^l\frac{
    C_{l+p+q}}{j!(l-j)!} \sum_{J_1 \subset [j]} \sum_{J_2 \subset [l]\setminus[j]}  \int_{(B_r(\bk x{1}{p}))^j \times(B_r(\bk x{p+1}{p+q}))^{l-j}}
		 \md \bk y1l \no \\
& \leq	 \phi(  \frac{s}{2} )\sum_{l=0}^{\infty}\sum_{j=0}^l\frac{
    C_{l+p+q}}{j!(l-j)!}(2(p+q)\theta_dr^d)^l  \no \\
& =   \phi( \frac{s}{2} ) \sum_{l=0}^{\infty} \frac{(4(p+q)
 \theta_dr^d)^l}{l!}  
  C_{l+p+q},
 \end{align*}
where the identity  $\sum_{j = 0}^l \frac{1}{j!(l-j)!} = \frac{2^l}{l!}$
justifies the last equality.
This establishes inequality~\eqref{e.hatC-bound-fine} and concludes the proof of Lemma~\ref{p:mppmixcl_double}. 
\end{proof}	

Before extending Proposition~\ref{p:mppmix_double} in the next section to
non-restricted marks (this result is  formulated in Theorem~\ref{t:mppmix_double}),
we comment on the assumptions and the applicability of the present proposition.

\begin{remark}[Summable exponential $\B$-mixing  ensures 
fast mixing of correlations]
\label{Rmkconvsum} 
The decay constants and decay function of the input process $\tP$ are specified in a minimal way 
in Lemma~\ref{p:mppmixcl_double} and  Propositions~\ref{p:mppmix_double}. However, in most applications of this result we shall assume that $\tP$ 
exhibits summable exponential $\B$-mixing correlations
as in 
Definition~\ref{def.A2}.
Indeed, it may be checked that the  {\em summable condition~\eqref{eqn:sum} of this definition} implies the following condition \begin{equation}\label{e.kappa-bound}
    \sum_{k}\frac{C_k w^k}{k!}<\infty,  
    \quad \forall \, w > 0
\end{equation}
and hence  \eqref{weaksum}, for the constants $C_k$ in Lemma~\ref{p:mppmixcl_double} and  Propositions~\ref{p:mppmix_double}. In addition, the fast nature of this mixing of correlations is guaranteed by the assumption~\eqref{phibd} on $\phi$ in Definition~\ref{def.A2}.   
The assumption of summable exponential $\B$-mixing correlations will be further used 
to relax the restriction on the marking function $\xitr$ in the next section. Also if  $\tP$ satisfies the summability condition \eqref{e.kappa-bound}, then one could alternatively  use the FME expansion from \cite[Theorem A.1]{klatt2025invariant};  we have instead required that $\tP$ satisfy the weaker assumption \eqref{e.kappa-bound-fine} when proving the the FME in Theorem \ref{FMEthm} via the intermediate Lemma \ref{l.FMEf}. 

\end{remark}

For real-valued short-range scores $\xi = \xi^{(r)}$ (viewed as short-range
marks in Proposition~\ref{p:mppmix_double}), the limit theory for $(\mu_n^\xi(f))_{n \in \N}$ stated in Section~\ref{s:asnormmpp} can already
be established using this proposition under assumption~\eqref{e.kappa-bound}
and merely fast $\B$-mixing correlations of the input process~$\tP$, without
requiring summable exponential-type mixing.
Moreover, if $\xitr$
 take the form of local $U$-statistics, 
the same conclusions hold under merely fast $\B$-mixing correlations of~$\tP$, 
as explained below.

\begin{corollary}[Limit theory for short-range scores]
\label{r.restricted-scores}
Consider a  marked point process $\tP$ on $\R^d\times\K$ having fast $\B$-mixing correlations as in Definition~\ref{B-omega-mixing}\ref{Bmixing} and~\ref{Bmixing-fast}, with correlation decay functions $\omega_{k}(s) = C_k  \phi(s)$, $k \in \N$, for $\phi$ fast decreasing,  and such that the constants $(C_k)_{k \in \N}$ satisfy  the constraints~\eqref{e.kappa-bound}. Assume that $\xi$ is a real-valued, {\em short-range score function}, i.e.,  
$\xi(\tx,\tP_n)= \xitr(\tx,\tP_n)$ for some $r\in (0,\infty)$ and all $n\in\N$.
Then the marked point process
 $\dbtilde\P=\sum \delta_{(x,\xi(\tx_i,\tP))}$  has 
 fast $\B$-mixing  correlations. 
 Consequently, with $
 \mu_n^\xi: = \sum_{x \in \P_n} \xi( \tx, \tP_n ) \delta_{ n^{-1/d}x }$,
 all statements in Theorems~\ref{t:cltmarkedpp} and~\ref{t:multcltmarkedpp},
as well as in Propositions~\ref{expvar} and~\ref{t:clt_linear_marks_new},
hold under their respective moment, variance, and stationarity conditions.
\end{corollary}
\begin{proof}
The result follows from Remark~\ref{Rmkconvsum} and, in particular, from the fact that 
condition~\eqref{e.kappa-bound} implies~\eqref{weaksum} for all $p,q \in \N$. 
Consequently, by Proposition~\ref{p:mppmix_double} and  fast decreasing $\phi$, the process~$\dbtilde{\P}$ 
exhibits fast $\B$-mixing correlations.

The assumption of summable exponential $\B$-mixing correlations of~$\tP$, 
required in all statements of the results mentioned in this corollary, 
serves solely to establish fast $\BL$-mixing correlations of the iterated scores (marks) 
$\xi(\tx,\tP_n)$, $\bxi(\tx,\tP_n)$, and $\xi(\tx,\tP)$ via Theorem~\ref{t:mppmix_double}. 
In the case where the scores are of short-range type, with radius
$r \in (0,\infty)$, Proposition~\ref{p:mppmix_double}, together with
condition~\eqref{e.kappa-bound} and fast $\B$-mixing correlations of the
input process~$\tP$, is sufficient, as explained above.
\end{proof}

\begin{corollary}[Limit theory for local $U$-statistics] 
\label{c.U-FME} 
Consider a marked point process $\tP$ on $\R^d \times \K$ having fast $\B$-mixing correlations 
as in Definition~\ref{B-omega-mixing}\ref{Bmixing} and~\ref{Bmixing-fast}.
Let the score function be a {\em local $U$-statistic} of the form
\begin{equation}\label{e:U-score}
\xi(\tx,\tP) = \xi_{h}(\tx,\tP) := \frac{1}{k!} \sum_{\bk \tx1{k-1} \in \tP^{(k-1)}}h(\tx,\tx_1,\ldots,\tx_{k-1}), 
\end{equation}
for some $k \in \N$ and a symmetric, bounded  function
$h : (\R^{d}\times\K)^{k} \to \R$,
satisfying $h(\tx_0,\tx_1 \ldots, \tx_{k-1}) = 0$ whenever either 
$\max_{1 \le i \le k-1} |x_i - x_0| > r$ for some $r>0$, 
or $x_i = x_j$ for some $i \neq j$, $i,j\in\{0,\ldots,k-1\}$.
Then the marked point process
\(
\dbtilde{\P} := \sum_{x_i \in \P} \delta_{(\tx_i, \xi_h(\tx_i, \tP))}
\)
has fast $\B$-mixing correlations. Consequently, 
with $
 \mu_n^\xi: = \sum_{x \in \P_n} \xi( \tx, \tP_n ) \delta_{ n^{-1/d}x }$, all statements in 
Theorems~\ref{t:cltmarkedpp} and~\ref{t:multcltmarkedpp}, 
as well as in Propositions~\ref{expvar} and~\ref{t:clt_linear_marks_new}, 
hold  under their respective moment, variance, and stationarity conditions.
\end{corollary}
\begin{proof}
This is a marked version of the local $U$-statistics model considered in \cite{BYY19}. 
As observed in that work, the FME expansion of $U$-statistics contains only a finite number of terms. 
Consequently, we do not need the convergence condition~\eqref{weaksum} for the FME, 
since the series in formulas~\eqref{e.FME-fg} and~\eqref{e.FME-fxg} involve only finitely many terms, 
and the same holds for the series in~\eqref{e.hatC-bound-fine}. 
Hence, Proposition~\ref{p:mppmix_double} applies without imposing any mixing condition 
on $\tP$ beyond fast $\B$-mixing.
\end{proof}

\subsubsection{Proof of Theorem~\ref{t:mppmix_double}
for general marking functions}
\label{s:proof_exp_stab_marking}

\begin{proof}[Proof of Theorem~\ref{t:mppmix_double}]
We first prove part~\ref{i.iterated-marks}  and then adjust the proof to obtain parts~\ref{i.BL-iterated-marks} and \ref{i.BLc-iterated-marks}. 

\vskip.2cm 

\noindent{\em Proof of part ~\ref{i.iterated-marks} of Theorem~\ref{t:mppmix_double}.} 
For $p,q \in \N$  consider again a  $(p+q)$-tuple of distinct points $\bk
x1{p+q}\in (\R^{d})^{(p+q)}$ with marks $\bk u1{p+q}\in\K^{p+q}$.
We recall  $s:= d(\bk x1p,\bk x{p+1}{p+q})$.  Without loss of
generality, we assume $s \in (4, \infty)$.	$\tP$ has summable exponential $\B$-mixing correlations 
  (Definition~\ref{def.A2}) with the exponents  $a \in [0,1)$
    in~\eqref{eqn:sum} and $b\in(0,d)$
    in~\eqref{phibd}.   Set
	\be \label{defsprime}
	r:=r(s):= (\frac{ s} {4})^{b(1-a)/(2(p+q+ d))}.
	\ee
 Since $s \in (4, \infty)$ and $p+q \geq 2$,  we easily have $r \in (1, s/4)$.	

Define a \emph{short-range mark} by truncating the domain of the score
$\xi$ to the corresponding balls
\begin{equation}\label{e.xi-t} \xitr(\tx,\tmu) :=  \xi(\tx,\tmu \cap B_r(x)),   \qquad  \tx=(x,u)\in\tmu.
\end{equation}

For  all $f \in \B((\K')^p)$ and
  $g \in  \B((\K')^q)$, we consider the variants
$\psi^r_f=f(\bk {v_r}1p)$, $\psi^r_g=g(\bk {v_r}{p+1}{p+q})$, $\psi^r_{gf}=f(\bk {v_r}1p) g(\bk {v^r}{p+1}{p+q})$
of the functionals $\psi_f,\psi_g,\psi_{gf}$ defined in~\eqref{e.psi-f}--\eqref{e.psi-fg}
with $\xi$ replaced by $\xitr$ and $v_{r,i} := \xitr(\tx_i,\tmu)$, $1 \leq i \leq p+q$.

We set 
\begin{align}\label{e.Axi}
  A&:=\sE_{\bk x1p}[f(\bk V1p)] \rho^{(p)}(\bk x1{p}),\\
B&:=\sE_{\bk x{p+1}{p+q}}[g(\bk V{p+1}{p+q})]\rho^{(q)}(\bk
x{p+1}{p+q}),\label{e.Bxi}\\
C&:=\sE_{\bk x1{p+q}}[f(\bk V1{p})g(\bk V{p+1}{p+q})]\rho^{(p+q)}(\bk x1{p+q}),\label{e.Cxi}
\end{align}
and similarly $\tA$, $\tB$, $\tC$  with $\bk V1{p+q}$ replaced by $\bk
{V_r}1{p+q}$, where $V_{r,i}:=\xitr(\tx_i,\tP)$, $i=1,\ldots,p+q$.

To prove Theorem~\ref{t:mppmix_double}~\ref{i.BL-iterated-marks}, we need to
show fast decay of $|AB-C|$ as $s \to \infty$.  Towards this, we observe first the following (general) inequalities :
$$|AB-C|\le |\tA \tB-\tC|+|AB- \tA \tB|+|C-\tC|$$
and
$$ |AB - \tA \tB| \leq |A(B - \tB)| + |(A -\tA)\tB|\le(|A|+|\tB|)(|A-\tA|+|B-\tB|).$$
Thus, we obtain 
\begin{equation}\label{ABC-bounds}
  |AB-C|\le  |\tA \tB-\tC|+
  (1+|A|+|\tB|)(|A-\tA|+|B-\tB|+ |C-\tC|).
  \end{equation}
In the following, we shall give bounds on the expressions in the right-hand side of~\eqref{ABC-bounds}. Note that since $|f|\le 1$, $|g|\le 1$, we have 
\begin{align}
\label{e:A+Bbds}
  (1+|A|+|\tB|)&\le 1+\kappa_p+\kappa_q<\infty,
\end{align}
where $\kappa_p:= \sup_{\bk x1p\in(\R^d)^p}  \rho^{(p)}(\bk x1p) < \infty$,
by Assumption \ref{Ass1}.

\textsc{Step 1:} (Bound for $|\tA \tB - \tC|$).
The marking function $\xitr$ and $\tP$ together satisfy 
the assumptions of Proposition~\ref{p:mppmix_double}.
Consequently, by~\eqref{e.hatC-bound-fine} with $r$ replaced by $r(s)$ 
we have  
\begin{align}\nonumber
|\tA \tB-\tC|&\le \phi\Bigl(\frac{s}{2}\Bigr)\sum_{l=1}^\infty \frac{(4(p+q) \theta_d r(s)^d)^l}{l!}  C_{l+p+q}\\
\label{e:decayfunction1} &\le   c_1 \exp\left(-\frac{s^b}{c_1}\right) c_2 \exp\left(c_3 ( \frac{s}{4})^{\frac{b}{2}}\right)
\end{align}
for some positive, finite constants $c_1,c_2,c_3$.
Indeed, \eqref{phibd} implies $\phi(s/2)\le
  c_1 \exp(-s^b/c_1) $ and the
  bound~\eqref{eqn:sum} on $ C_{l+p+q}$, by Stirling's formula,
gives the exponential  bound $c_2 \exp(c_3 (
s/4)^{b/2})$ on the  series in~\eqref{e:decayfunction1};  
  for details, we refer the reader to the arguments
in~\cite[display (3.26), page 872]{BYY19}.

\textsc{Step 2:} (Bounds for  $|A-\tA|$, $|B-\tB|$, and $|C-\tC|$).  We use that $\xi$ satisfies classical
stopping set stabilization as in Definition \ref{def.stabilizing_marking} with $r$ replaced by $r(s)$ and for $p,q,p+q$. Specifically,  observe that if  $R^\xi(\tx_i,\tP)\le r(s)$ for all $1,\ldots,p$  then $\xi(\tx_i,\tP) = \xitr(\tx_i,\tP)$ and $A=\tilde A$. Consequently, 
\begin{align}\nonumber
  |A-\tA|&\le 2\kappa_p \sP_{\bk x1p}
\Bigl(\max_{i=1,\ldots,p}R^\xi(\tx_i,\tP) > r(s) \Bigr) \\
\label{e:decayfunction2}   &\le 2p\kappa_p\varphi_p(r(s))\,,
\end{align}
which,  by~\eqref{defsprime} and the fast decreasing property of $\varphi_p$, is a fast
decreasing function of~$s$. The same arguments work for
$|B-\tB|$, and $|C-\tC|$.

Collecting all inequalities   \eqref{ABC-bounds},
 \eqref{e:A+Bbds},
 \eqref{e:decayfunction1}, and \eqref{e:decayfunction2},
we observe that the doubly marked point process
$\dbtilde\P$ has fast $\B$-mixing  correlations as in~\eqref{e.Bmppmixing} (but for  $\dbtilde \P$) with the decay function $\omega_{p+q}$ depending 
on $ \phi,  \{ C_k\}_{k = 1}^{\infty}$,  
$\varphi_{p+q}$,
and $\ka_{p+q}$. Here we also  use the assumption that $\max( \varphi_q,
\varphi_p) \le\varphi_{p+q}$ and $\max( \kappa_q,
\kappa_p) \le \kappa_{p+q}$ for all $p, q\in \N$.
This completes the proof of part~\ref{i.iterated-marks} of Theorem~\ref{t:mppmix_double}. 

\vskip.3cm

\noindent{\em Proof of part \ref{i.BL-iterated-marks} of Theorem~\ref{t:mppmix_double}.} We follow the arguments of the previous part,  now applied   to the \BL-localizing marking  function $\xi$, for which we again consider its short-range version  $\xitr$ in~\eqref{e.xi-t} with $r$ as at ~\eqref{defsprime}. Since we want to show 
that the marks $\xi$ have \BL-mixing correlations as in  Definition~\ref{d.omegahmixing}\ref{i.mixing} and~\ref{i.fast-mixing},
we consider only \BL \,functions $f, g$ respectively  on $(\K')^p$, $(\K')^q$. Given these functions we  consider the expressions $A,B,C$ as in~\eqref{e.Axi}--\eqref{e.Cxi} and their counterparts  $\tilde A,\tilde B,\tilde C$. The general bounds yield~\eqref{ABC-bounds}.
The bound~\eqref{e:decayfunction1}
holds  by  Proposition~\ref{p:mppmix_double}
since we have assumed $\tP$ has summable exponential $\B$-mixing correlations. Indeed~\eqref{e.hatC-bound-fine} holds for all bounded functions $f,g$ and in particular for our \BL \,functions $f,g$. Thus \textsc{step 1}  goes through verbatim. Only in \textsc{step 2} do we need to modify the proof, since~$\xi$  satisfies (the weaker)  \BL-localizing property.   
However, we may use this property to obtain   
\begin{align}
  |A-\tA|&
  \le  \kappa_p \Bigl| \sE_{\bk x1p}
[f(\bk V1p)] - \sE_{\bk x1p}[f(\bk {V_r}1p)] \Bigr| \le
  \kappa_p
  d_{\BL, \bk x1p}(\bk V1p,\bk {V_r}1p)
\le \kappa_p \varphi_p(r(s))\,,     \label{Lip-argumentA}
\end{align}
where the last inequality follows directly from \BL-localization~\eqref{Lp-stab-infinite} of the marking function $\xi$ for $p$.

Using \BL-localization~\eqref{Lp-stab-infinite} of the marking function $\xi$ for $q$ and $p+q$, the same argument extends to 
 $|B-\tB|$ and $|C-\tC|$ respectively, thereby completing the proof, exactly as in part \ref{i.iterated-marks}. For $|C-\tC|$, we have also used the $(\K')^{p+q} \owns v =(v_p,v_q) \mapsto f(v_p)g(v_q)$ is \BL~if $f$ and $g$ are. 

\noindent{\em Proof of part \ref{i.BLc-iterated-marks} of Theorem~\ref{t:mppmix_double}.} 
We follow the proof approach as for parts \ref{i.BL-iterated-marks} and \ref{i.iterated-marks} but now choose 
\begin{align*}
\tilde{A} &= \sE_{\bk x1p}[f(\bk{\xi^{\cup(r)}}{1}{p}(\tx,\tP))] \, \rho^{(p)}(\bk x1p), \\
\tilde{B} &=  \sE_{\bk x{p+1}{p+q}}[g(\bk{\xi^{\cup(r)}}{p+1}{p+q}(\tx,\tP))] \, \rho^{(q)}(\bk x{p+1}{p+q}),  \\
\tilde{C} &= \sE_{\bk x{1}{p+q}}[f(\bk{\xi^{\cup(r)}}{1}{p}(\tx,\tP))g(\bk{\xi^{\cup(r)}}{p+1}{p+q}(\tx,\tP))] \, \rho^{(p+q)}(\bk x{1}{p+q}),
\end{align*}
Also, modifying~\eqref{defsprime}, we define
\begin{equation}\label{defsprime-cl}
	r := r(s) := \Bigl(\frac{s}{4}\Bigr)^{\delta_{p+q}}, 
	\qquad 
	\delta_{p+q} := \min\!\left\{ \frac{b(1-a)}{2(p+q+d)},\, \alpha_{p+q}^{-1} \right\}.
\end{equation}
Here $\alpha_{p+q}\in[0,\infty)$ denotes the exponent governing the
separation of the two groups of points, as postulated  in
Definition~\ref{def.Lp-stabilizing_marking-weak}. By that definition, the bound
\eqref{Lp-stab-infinite-weak2} applies whenever
\[
s = d\bigl(\bk{x}{1}{p}, \bk{x}{p+1}{p+q}\bigr) > r^{\alpha_{p+q}}.
\]
In the present setting, this condition is satisfied by the above choice of $r$.
Indeed, since $s>4$ and $\delta_{p+q}>0$, we have $r=(s/4)^{\delta_{p+q}}>1$.
Moreover, by construction $\delta_{p+q}\le \alpha_{p+q}^{-1}$, and hence
\[
r^{\alpha_{p+q}} = \Bigl(\frac{s}{4}\Bigr)^{\alpha_{p+q}\delta_{p+q}}
\;\le\; \frac{s}{4}
\;<\; s,
\]
which yields the desired inequality $s>r^{\alpha_{p+q}}$.
Furthermore, since $\frac{b(1-a)}{2(p+q+d)}<1$, we have $\delta_{p+q}<1$ and thus 
\[
r = \Bigl(\frac{s}{4}\Bigr)^{\delta_{p+q}} < \frac{s}{4},
\]
as in the proof of the previous parts of Theorem~\ref{t:mppmix_double}.

Observe now that  \eqref{e.hatC-bound-fine-cl} implies  that $|\tilde{A}\tilde{B} - \tilde{C}|$ satisfies the bound \eqref{e:decayfunction1} with our new definition of $\tilde A$, 
$\tilde B$, and $\tilde C$ given above.  By \eqref{Lp-stab-infinite-weak1} and following verbatim the argument in \eqref{Lip-argumentA}, we have that $|A - \tilde{A}|$ and $|B - \tilde{B}|$ are bounded by $\kappa_p\varphi_p(r(s))$ and $\kappa_q\varphi_q(r(s))$, respectively.  Second, we have 
\begin{align*}
| C - \tilde{C}| & \leq \kappa_{p+q} d_{\BL,\bk{x}{1}{p+q}}
\Bigl(\bk{\xi}{1}{p+q}(\tx,\tP),
\bigl(
\bk{\xi^{\cup(r)}}{1}{p}(\tx,\tP),
\bk{\xi^{\cup(r)}}{p+1}{p+q}(\tx,\tP)
\bigr)
\Bigr) \\
& \leq \kappa_{p+q} d_{\BL,\bk{x}{1}{p+q}}
\Bigl(\bk{\xi}{1}{p+q}(\tx,\tP),\bk{\xi^{\cup(r)}}{1}{p+q}(\tx,\tP)\Bigr) \\
& \qquad + \kappa_{p+q} d_{\BL,\bk{x}{1}{p+q}}
\Bigl(\bk{\xi^{\cup(r)}}{1}{p + q}(\tx,\tP),
\bigl(
\bk{\xi^{\cup(r)}}{1}{p}(\tx,\tP),
\bk{\xi^{\cup(r)}}{p+1}{p+q}(\tx,\tP)
\bigr)
\Bigr) \\
& \leq 4\kappa_{p+q}\varphi_{p+q}(r(s)),
\end{align*}
where the second inequality follows from the triangle inequality and the third is due to \eqref{Lp-stab-infinite-weak1} and \eqref{Lp-stab-infinite-weak2}, in conjunction with the $r(s) =  (s/4)^{\delta_{p+q}}$. 

Thus, combining the above bounds for $|A-\tilde{A}|, |B-\tilde{B}|, |C - \tilde{C}|$ and $|\tilde{A}\tilde{B} - \tilde{C}|$ as in \eqref{ABC-bounds} and along with \eqref{e:A+Bbds}, we have the requisite fast decreasing bound on $|AB - C|$.
\end{proof}

\begin{remark}[Iterated \BL-localized marks via short-range scores]
\label{r.rem-short-range-scores}
Observe that the proof of part~\ref{i.BL-iterated-marks} of
Theorem~\ref{t:mppmix_double} presented above applies verbatim when replacing the  short-range scores $(\xi^{(r)})_{r \geq 1}$ generated via natural restrictions of $\xi$ by a family of `external' short-range scores $(\xi^{[r]})_{r \geq 1}$, satisfying 
\BL-localization as explained in Remark~\ref{i.rem-short-range-scores}
in Section~\ref{s:remarksltthms}.
\end{remark}

\noindent{\em Proof of Corollary~\ref{c:mppmixcl_double}.} 
This proof closely follows that of Theorem~\ref{t:mppmix_double} and makes use of the bound \eqref{e.Bmppmixing} describing 
summable exponential $\B$-mixing correlations for a family of marked point processes with the family being $(\tP_{n})_{n \in \N}$.  It suffices to note  that the
functions  $\phi$, $\{\varphi_p\}_{p = 1}^{\infty}$, and the constants 
 $\{ C_k\}_{k = 1}^{\infty}$,  $\{\ka_p\}_{p = 1}^{\infty}$ related to  $(\tP_{n})_{n \in \N}$ and $\xi$ and appearing implicitly in the proof  
 of Theorem~\ref{t:mppmix_double}
 uniformly satisfy the   definitions of the mixing correlations and localization, respectively, stabilization,  on the sequence $(\tP_{n})_{n \in \N}$. Indeed, these functions and constants 
can be tracked from the bounds in  \eqref{ABC-bounds}--\eqref{e:decayfunction2}.\qed
 
\section{Proofs of the  limit  theorems}
\label{s.proof-limit}

In this section we prove limit results for dependently marked
  point processes, as given in Sections \ref{ss:asnormmpp0} and \ref{s:clta}, namely central limit theorems and  expectation and variance asymptotics.
The proof of the main central limit theorem
(Theorem~\ref{t:clt_linear_marks}) uses the classical cumulant method
for integrals $\mu^\xi_n(f)= \int f d \mu_{n}^{\xi}$ of  random (possibly signed, atomic) measures   defined in~\eqref{e:signed-measure}.
 A key argument of this method
 is  the approximate factorization (or `mixing' in our terminology) of the correlation functions \eqref{eqn:mixedmoment} of these    measures.
This method was previously used in~\cite[Theorem~1.13]{BYY19}, where the
measures involve real-valued scores $\xi$ functionally constructed on a
non-marked point process. However, the cumulant method applies  to  general,
  purely atomic random signed
 measures satisfying moment conditions and exhibiting fast mixing  correlations. We believe this result is of independent interest and present it in Section~\ref{s:cumulant},  recollecting  some of the arguments  of the proof  of~\cite[Theorem 1.13]{BYY19}.  In Section \ref{s:proofssclt},  we use this general central limit theorem to prove our main central limit theorem ---Theorem \ref{t:clt_linear_marks}. Next, 
 in Section \ref{s:stab_lemma_BL}, we 
 show that if scores $\xi(x, \tP_n)$ 
 satisfy \BL-localization then they 
 have limits under Palm
 distributions as the domain $W_n$ increases up to $\R^d$.  These limits   determine the precise expectation and variance asymptotics given in 
 Theorem \ref{expvar}, whose proof is in Section \ref{s:proofsexp}. Section \ref{s:proof_multclt} proves the multivariate central limit theorem  (Theorem \ref{t:multcltmarkedpp}) whereas  Section \ref{proofProp5.5} establishes the  limit theory for scores on the infinite window, cf. Proposition \ref{t:clt_linear_marks_new}.

\subsection{Cumulant method for asymptotic normality of purely atomic random  signed measures}\label{s:cumulant}
We outline the cumulant method for random measures, but omit full details of some standard derivations, which may be found in e.g. \cite[Theorem 1.13]{BYY19} and \cite{Nazarov12}.

Let $\tP=\{(x_i,\xi_i)\}$ be a simple marked point process with ground process
(locations of atoms) $\P=\{x_i\}\subset \R^d$ and real-valued marks (scores)
$\xi_i\in\R$. Consider the random, purely atomic, possibly signed, locally
finite measure
\(
\amu=\sum_{(x_i,\xi_i)\in\tP}\xi_i\,\delta_{x_i}
\)
on $\R^d$. (Conversely, any such measure can be represented as a simple marked
point process.)

For all $p\in \N$, let  the ground process $\P$ admit bounded correlation functions
 $\rho^{(p)}$ (densities of the factorial moment measures). Then
define the   (generalized) correlation functions of $\amu$ for all $k_1,\ldots,k_p \in \N$
\begin{equation}
m^{\bk k1{p}}(\bk x1{p}) := \sE_{\bk x1{p}}\left[\xi_1^{k_1} \ldots\xi_p^{k_{p} }\right] \rho^{(p)}(\bk x1{p}),  \label{eqn:mixedmoment-general} 
\end{equation}
where   $\xi_i$ are  random marks (scores) of fixed, distinct locations of atoms  $x_1,\ldots,x_p
\in\R^d$  of $\tP_n$ under Palm probabilities  $\sE_{\bk x1{p}}$,
provided the expectations $\sE_{\bk x1{p}}[\xi_1^{k_1}
  \ldots\xi_p^{k_{p} }]$   are well defined.  Let $M^p[\langle g, \amu \rangle]:=\E[(\int_{\R^d}g(x)\,\amu(\md x))^p]$ be the moments of 
integrals with respect to $\amu$ of bounded real-valued functions $g$ with bounded support.  The correlation functions allow us to write
\begin{align}
	M^p[\langle g, \amu \rangle]&=\sum_{\gamma \in \Pi[p]}
        \int_{\R^{|\gamma|}}\prod_{j=1}^{|\gamma|}g^{|\gamma(j)|}(x_j)m^{(|\gamma(1)|,\ldots,|\gamma(|\gamma|)|)}(x_1,\ldots,x_{|\gamma|}) \md x_1 \ldots \md x_{|\gamma|} \no \\
\label{e:momintgmu}  &=\sum_{\gamma\in\Pi[p]}\langle
        \bigotimes_{i=1}^{|\gamma|}g^{|\gamma(i)|}m^{\gamma},\lambda^{|\gamma|}\rangle.
  \end{align}
 Here $\Pi[p]$ is the set of all unordered partitions of 
	$[p] := \{1,...,p\}$, and, given  a partition
$\gamma\in\Pi[p]$, we let $|\gamma|$ denote the number of partition classes and for convenience, we  order  its  elements
$\gamma=\{\gamma(1),\ldots,\allowbreak\gamma(|\gamma| )\}$
according to the number of   elements in $\gamma(i)$. Further,  $m^{\gamma}:=m^{(|\gamma(1)|,\ldots,|\gamma(|\gamma|)|)}$, $\lambda^l$ denotes the Lebesgue measure on $(\R^d)^l$, 
and $\bigotimes$ denotes the tensor product of functions.

The correlation functions $m^{\bk k1{p}}$ give rise to  the {\em Ursell functions} or {\em truncated correlation functions} $m_{\top}^{\bk k1{p}}$ of $\amu$ given by
\begin{equation}\label{e.Ursel-direct}
	m_\T^{\bk k1p}(\bk x1p)
	=\sum_{\gamma\in\Pi[p]}(-1)^{|\gamma|-1}(|\gamma|-1)!
	\prod_{i=1}^{|\gamma|} m^{(k_j:j\in \gamma(i))}(x_j:j\in\gamma(i))\,,
	\end{equation}
provided  $m^{(k_j:j \in J) }<\infty$ for all $J \subset [p]$;
see~\cite[Section 4.4.1]{BYY19}.  Ursell functions allow  one to
express the $p$-th order cumulants $S^p$  of the integrals $\langle g, \amu \rangle$ 
in terms of integrals of Ursell functions just as moments of $\langle g,\amu\rangle$ are expressed as integrals with respect to correlation functions $m^{\bk k1{p}}$ as in \eqref{e:momintgmu}. Indeed, for all $p \in \N$ we have the relation
 \begin{align}%\nonumber
 S^p[\langle g, \amu \rangle] %&= \sum_{\gamma \in \Pi[k]}
 %\int_{\R^{|\gamma|}}\prod_{j=1}^{|\gamma|}g^{|\gamma(j)|}(x_j)m_{\top}^{(|\gamma(1)|,\ldots,|\gamma(|\gamma|)|)}(x_1,\ldots,x_{|\gamma|}) \md x_1 \ldots \md x_{|\gamma|}\\
& =\sum_{\gamma\in\Pi[p]}\langle      \bigotimes_{i=1}^{|\gamma|}g^{|\gamma(i)|}m_{\top}^{\gamma},\lambda^{|\gamma|}\rangle\,,\label{e.cumulant-intergral}
 \end{align}
 provided $m^{\bk k1{p'}}<\infty$ for all  $k_i\in \N$, $k_1+\ldots+k_{p'}\le p$;
 see~\cite[Sections 4.3.2]{BYY19}
   for the $p$-th order cumulants $S^p[Y]$ of  a random variable $Y$
   and~\cite[(4.36)]{BYY19} for the expresion of  cumulants $S^p[\langle g, \amu \rangle]$ in terms of Ursell functions.

Next, we observe that the  property of  {\em mixing of correlations functions} $m^{\bk k1{p}}$ 
is equivalent to a {\em diameter bound on the  Ursell functions of $\amu$} of these measures.
\begin{lemma}[Mixing of generalized correlation functions of $\amu$ is equivalent to a diameter bound on Ursell functions of $\amu$]\label{l.moments-cumulants}

For a fixed  collection of integers $k_1,\ldots,k_p$, $p\ge2$,
 assume that  for all $J\subset[p]$,
 \begin{equation}
 \label{e.moments-uniformly-bounded}    
 \sup_{(x_j)_{j \in J} \in\R^{d|J|}} 
 m^{(k_j:j \in J) }( x_j : j \in J)<\infty.
 \end{equation}
 The following statements are equivalent:
 \begin{enumerate}[wide,label=(\roman*),  labelindent=0pt]
\item \label{i.mixing-moments}There exist fast decreasing functions $\tilde\phi=\tilde\phi_{k_1,\ldots,k_p}$
and constants
	$\tilde C=C(k_1,\ldots,k_p)$, such that
  for any subset  of integers  $J\subset [p]$, $|J|\ge2$, a nonempty, proper subset
  $I\subsetneq J$, and all configurations $\{x_j:  j \in J \}, x_j \in\R^d$, of distinct points
	we have
\begin{equation}\label{e.clustering-generalized-mixed-monents1}
	\Bigl|m^{(k_j : j\in J)}(x_j : j\in J)-
	m^{(k_j: j\in I)}(x_j:j\in I)\, m^{(k_j: j\in J\setminus I)}(x_j:j\in
	 J\setminus I)\Bigr|\le
	\tilde C\tilde\phi(s)\,,
	\end{equation}
	where  $s:=d\bigl(\{x_j:j\in I\}, \{x_j:j\in  J\setminus I\}\bigr)$.
      \item \label{i.mixing-cumulants} 	 There exist fast decreasing functions $\tilde\phi_\T=\tilde\phi_\T^{k_1,\ldots,k_p}$ and constants
	$\tilde C^\T=\tilde C^\T(k_1,\ldots,k_p)$, 
        such that for any subset  of integers $ J\subset [p]$, $|J|\ge2$ and all configurations  $x_j\in\R^d$, $j\in J$ of distinct points
	we have
	\begin{equation}\label{e.Usell-clustering-bound}
	|m_{\T}^{(k_j:j\in J)}(x_j : j\in J)|\le
	\tilde
	C^\T\tilde\phi_\T\bigl(\mathrm{diam}(x_j : j\in J)\bigr)
	\,,
	\end{equation}
where $\mathrm{diam}(x_j : j\in J):=\max_{i,j\in J}(|x_i-x_j|)$.
 \end{enumerate}
 Moreover, the bound~\eqref{e.Usell-clustering-bound}  for $J=[p]$ implies
 	\begin{equation}\label{e.Ursel-integral-bound}
\sup_{x_1\in \R^d}\int_{
(\R^d)^{(p-1)}}|m_\T^{\bk k1p}(\bk x1p)|\,\md x_2\cdots \md x_p< \infty.
	\end{equation}
Further, $\tilde\phi_\T$ and $\tilde C^\T$ in (ii) depend only on $\tilde \phi$, $\tilde C$, the dimension $d$ and the bound in \eqref{e.moments-uniformly-bounded} and similarly the $\tilde\phi$ and $\tilde C$ in (i) depend only on $\tilde \phi_\T$, $\tilde C^\T$ and the bound in \eqref{e.moments-uniformly-bounded}.
\end{lemma}
One may interpret \eqref{e.Ursel-integral-bound} as {\em Brillinger mixing} for the random measure $\amu$ and it can be shown to be equivalent to volume order growth of cumulants.  As we shall see shortly in Theorem \ref{l.CLT-Cumulants},  this suffices to derive a central limit theorem.  Brillinger mixing and its consequences have been studied for point processes (see \cite[Section 3.5]{Ivanoff82}, \cite{Heinrich_Schmidt_1985} and \cite[Theorem 3.2]{biscio2016brillinger}) and our results extend these to purely
 atomic random signed measures.
\begin{proof}[Proof of Lemma \ref{l.moments-cumulants}]
  Under  the (point-wise) assumption $m^{(k_j:j \in J) }<\infty$ for all $J\subset[p]$
  one establishes the following relation
  	\begin{align}
	m_\T^{(k_1,\ldots,k_p)}(x_1,\ldots,x_p)&=
	m^{(k_1,\ldots,k_p)}(x_1,\ldots,x_p)-
	m^{(k_j:j\in I)}(x_j:j\in I)\,
	m^{(k_j:j\in J\setminus I)}(x_j: j\in J\setminus I) \no \\
	&+\sum_{\genfrac{}{}{0pt}{2}{\gamma\in\Pi[p], |\gamma|>1}{		\gamma\,\text{mixes}\,\{I,I^c\}}}
	\prod_{i=1}^{|\gamma|} m_\T^{(k_j:j\in
		\gamma(i))}(x_j:j\in\gamma(i))\,\label{e.Ursel-cluster},
	\end{align}
 where a partition $\gamma$ of $\Pi[p]$ is said to mix $\{I,I^c:=[p]\setminus I\}$ when there exists  an element $\gamma(i)\in\gamma$ such that
 $\gamma(i)\cap I\not =\emptyset$ and $\gamma(i)\cap I^c\not =\emptyset$.             
 The proof of the equivalence property follows now from the fact that there exists constant
 $\tilde c^\T_{|J|}$ (depending on the dimension $d$) such that for each configuration $\{x_j:j\in J\}$, $x_j\in\R^d$ , there exists a partition
 $\{I,I^c\}$ of $J$ such that $d(\{x_j:j\in I\}, \{x_j:j\in I^c\})\ge \tilde c^\T_{|J|}\mathrm{diam}(x_1,\ldots,x_p)$ and, trivially,  for any partition 	$\{I,I^c\}$ of $J$,   $d(\{x_j:j\in I\}, \{x_j:j\in I^c\})\le\mathrm{diam}(x_1,\ldots,x_p)$.
Indeed, for  the first term  in the right-hand-side of  \eqref{e.Ursel-cluster} we have 
$$|m^{(k_1,\ldots,k_p)}-
	m^{(k_j:j\in I)}
	m^{(k_j:j\in I^c)}|\le 
\tilde C(k_1,\ldots,k_p)
\tilde\phi_{k_1,\ldots,k_p}(\tilde c^\T_{|J|}\mathrm{diam}(x_1,\ldots,x_p)).$$
Also, for each term of  the  sum $\sum_{\gamma}\dots$  
 in \eqref{e.Ursel-cluster} we 
consider   $\gamma(i)\in\gamma$ 
intersecting both $I$ and $I^c$ and  develop a similar bound by the induction, 
and using the uniform bounds~\eqref{e.moments-uniformly-bounded}
for other $\gamma(i')\in\gamma$.  
This way, $|\sum_{\gamma}\dots|$  
 can be bounded   by a {\em finite} sum of  expressions $C_\gamma \tilde\phi_{\gamma} (\tilde c_\gamma \mathrm{diam}(x_1,\ldots,x_p))$, with some constants $\tilde C_\gamma<\infty$, $\tilde c_\gamma>0$ and fast decreasing functions~$\tilde \phi_\gamma$. This justifies that ~\ref{i.mixing-moments} implies ~\ref{i.mixing-cumulants}. The proof of the converse again uses the expression \eqref{e.Ursel-cluster} in a similar way (but not requiring the induction).  
 
For the proof of \eqref{e.Ursel-integral-bound}, fix $x_1 \in \R^d$ and partition $(\R^d)^{p-1}$ into sets $(G_l)_{l \in \{0\} \cup \N}$ defined as follows:
 \begin{align*}G_0 &:= \{(x_2,\ldots,x_p)\in \R^{d(p-1)} : \text{diam}(x_1,\ldots,x_p)\le 1\}\\
G_l &:= \{(x_2,\ldots,x_p)\in \R^{d(p-1)} : 2^{l-1}<\text{diam}(x_1,\ldots,x_p)\le 2^l\},\quad l \in \N.
\end{align*}
Now use estimate~\eqref{e.Usell-clustering-bound}  to bound the integral on the left-hand side of~\eqref{e.Ursel-integral-bound} as
$$
\sum_{l=0}^{\infty} \int_{
G_l}|m_\T^{\bk k1p}(\bk x1p)|\,\md x_2\cdots \md x_p \leq \theta_d^{p-1}\tilde C^\T+\theta_d^{p-1}\tilde C^\T\sum_{l=1}^\infty 2^{dl(p-1)}\tilde \phi_\T( 2^{l-1})<\infty,
$$
since $\tilde \phi_\T$ is fast decreasing.
  \end{proof}

\begin{theorem}[Central limit theorem for purely atomic random signed measures]\label{l.CLT-Cumulants} 
 Let $(\amu_n)_{n \in \N}$  be a sequence of purely
 atomic random signed measures and assume they admit  correlation functions $m^{\bk k1{p}}(\bk x1{p})= m^{\bk k1{p}}(\bk x1{p};n)$ defined as at \eqref{eqn:mixedmoment-general}.   Assume  there exist fast decreasing functions $\tilde\phi_{k}$,
 and constants $M_{k}$,	$\tilde C_k$, 
 such that for all $p\ge 2$, $k_1,\ldots,k_p\in \N$, $k_1+\ldots+k_p=k$
 \begin{equation}\label{e.Moment-general}
   \sup_{n \in \N}\sup_{\bk x1p\in\R^{dp}}m^{\bk k1{p}}(\bk x1{p};n)\le M_{k}<\infty
 \end{equation}
 and statement (i) of Lemma~\ref{l.moments-cumulants} holds for
 the functions $\tilde\phi =\tilde\phi_k$ and constants  $\tilde C=\tilde C_k$ (uniformly in $n \in \N$).
 If the sequence  $(g_n)_{n \in \N}$ of real-valued, measurable functions on $ \R^d$, 
 is uniformly  bounded in $n$ (i.e., $\sup_{n \geq 1}  \|g_n\|_\infty\le M<\infty$) and  satisfies
$\|g_n\|_1 := \int_{\R^d}|g_n(x)|\,\md x=O(n^\kappa)$ for some $\kappa<\infty$, and
  $\Var[\langle g_n,\amu_n\rangle] = \Omega(n^{\nu})$ for some $\nu > 0$, then
$$
(\Var[\langle g_n,\amu_n\rangle])^{-1/2} ( \langle g_n,\amu_n\rangle - \sE[\langle g_n,\amu_n\rangle]) \stackrel{d}{\Rightarrow} Z
$$
where $Z$ denotes a standard normal random variable.
\end{theorem}

\begin{proof}
By the moment condition~\eqref{e.Moment-general}, the assumption of statement (i) in
Lemma~\ref{l.moments-cumulants} implies the  statement (ii) therein. Hence, by~\eqref{e.cumulant-intergral}  and~\eqref{e.Ursel-integral-bound}, we obtain the bounds on the cumulants of the integrals of order
$$
S^p[\langle g_n,\amu_n \rangle]\le
(1+\|g_n\|_{\infty})^{p-1} A_p\|g_n\|_1 =(1+A)^{p-1} A_p n^\kappa, 
$$
for some constants $A_p<\infty$. In other words,  cumulants have growth $S^p[\langle g_n, \amu_n\rangle]= O(n^{\BBLp{\kappa}})$ for all $p \geq 2$.
The  assumed variance lower bound 
 $\Var[\langle g_n,\amu_n\rangle] = \Omega(n^\nu)$   implies
$$
S^p\left[ \Var[\langle g_n,\amu_n\rangle]^{-1/2}\left(
\langle g_n,\amu_n\rangle -\sE[\langle g_n,\amu_n\rangle]\right)\right]= O(n^{\kappa-p\nu/2}),
$$
and hence the  cumulants or order $p \in (2\kappa/\nu, \infty)$ in the left-hand side above
vanish as  $n \to \infty$. This suffices to prove the central limit theorem by a classical result due to Marcinkiewicz (see e.g.  \cite[Lemma 3]{soshnikov2002gaussian}) which says that vanishing of cumulants of high-enough order implies normal convergence. 
$ $\par\nobreak\ignorespaces %This forces the Halmos symbol to appear at the end of the line (when there is no space left in the previous one).
\end{proof}
The above general CLT can be used to deduce CLTs such as \cite[Theorem 4.1]{klatt2022genuine} for certain linear statistics of point processes with fast mixing correlations and in fact, Brillinger mixing property (i.e., \eqref{e.Ursel-integral-bound}) as in their paper also suffices. However, in case the moment bound \eqref{e.Moment-general} fails as in \cite[Theorem 3.9]{mastrilli2024estimating}, the above approach needs to be suitably modified to derive appropriate cumulant bounds and hence a CLT.

\subsection{Proof of umbrella CLT in  Section~\ref{ss:asnormmpp0}%--- Theorem \ref{t:clt_linear_marks} 
}
\label{s:proofssclt}

We will use Theorem~\ref{l.CLT-Cumulants} to establish the central limit theorem
stated in Theorem~\ref{t:clt_linear_marks}. To this end, we first verify
condition~\ref{i.mixing-moments} of Lemma~\ref{l.moments-cumulants}, which
concerns the factorization of the correlation functions
$m^{\bk{k}{1}{p+q}}(\bk{x}{1}{p+q};n)$ defined in~\eqref{eqn:mixedmoment-general},
associated with the family of atomic measures
\[
\amu_n := \sum_{x \in \P_n} \xi_{i,n}\,\delta_x,
\]
where the scores $\xi_{i,n}$ act as weights attached to the atoms of the ground process.

Recall that in our umbrella Theorem~\ref{t:clt_linear_marks} these scores are
given \emph{per se}, whereas in applications they are typically constructed
functionally via score functions of the form $\xi_{i,n} = \xi(\tx_i,\tP_n)$.
In full generality, we seek to show that there exist a family  of fast decreasing functions $(\tilde \omega_k)_{k\in \N}$, which provide  the following bounds  for
 $p,q,k_1,\ldots,k_{p+q} \in \mathbb{N}$:
\be \label{dbtpptruncmixing}
 \sup_{n \in \N} \sup_{\bk x1{p+q} \in (W_n)^{p+q}  }| m^{\bk k1{p+q}}(\bk x1{p+q};n) -  m^{\bk k1{p}}(\bk x1p;n) m^{\bk k{p+1}{p+q}}(\bk x{p+1}{p+q};n)| \leq
 \tilde \omega_K(s), 
\ee
 where $K=k_1+\ldots+k_{p+q}$ and $s=d(\bk x1p,\bk x{p+1}{p+q})$ is the distance defined in~\eqref{defs}. This bound will be shown to be a consequence of the  fast \BL-mixing   correlations~\eqref{e.mppmixing} (with $\xi_{i,n}$ representing $U_n(x_i)$ and   $f,g \in \BL$) under a suitable moment condition. In other words, the existence of the function $\tilde \omega_K$ in~\eqref{dbtpptruncmixing} is guaranteed under
finiteness of the following moments 
\begin{equation}\label{e.moment-cond-Lemma}
\widetilde{M}_{K,\epsilon, p + q} := \sup_{n \in \N} \sup_{r \le p+q} \sup_{\bk x1r \in (W_n)^r}   \E_{\bk x1{r}}\left[ \max(1,|\xi_{1,n}|^{K(1+\epsilon)} ) \right]
\end{equation}
for  some $\epsilon>0$ where,
we recall,   $\xi_{1,n}$  is the random score (mark) of the fixed point~$x_1$  under the Palm distribution $\Palm_{\bk{x}1{r}}$. The $p$-moment condition~\eqref{e:xinpmom} for $p=K(1+\epsilon)$ implies the finiteness of $\widetilde{M}_{K,\epsilon, p + q}$ since  $\widetilde{M}_{K,\epsilon, p + q} \leq M^{\xi}_{K(1+\epsilon)}$. This last inequality  holds since, in view of  \( K \geq p+q \), the right-hand side of
~\eqref{e.moment-cond-Lemma} is upper bounded by the same expression with 
 \( \sup_{r \leq p+q} \) replaced by  \( \sup_{r \leq K(1+\epsilon)} \) in~\eqref{e.moment-cond-Lemma}.

\begin{lemma}[Factorization of correlation functions under moments and \BL-mixing conditions of the scores]
\label{l:non-bdd-mppmixing} 
 Let $(\tilde\P_n)_{n \in \N} = \bigl(\{(x_i,\xi_{i,n})\}\bigr)_{n \in \N}$
be a family of marked point processes on $\R^d$ sharing the same ground
point process $\P$, and equipped with real-valued marks (scores)
$\xi_{i,n} \in \R$. Fix $p,q,k_1,\ldots,k_{p+q} \in \N$ with $K = k_1 + \ldots + k_{p+q}$. Assume that $(\tilde\P_n)_{n \in \N}$ satisfies the fast \BL \,mixing correlation condition~\eqref{e.mppmixing} for $p,q$ with fast decay function $\hat\omega_{p+q}$ (and replacing $U(x_i)$ with $\xi_{i,n}$ in \eqref{e.mppmixing}) and for all functions $f \in \BL(\R^p)$, $g \in \BL(\R^q)$.

For some $\epsilon>0$ we assume finiteness of the moment $\widetilde{M}_{K,\epsilon, p + q}$ as in \eqref{e.moment-cond-Lemma} associated to the marks $\xi_{i,n}$.  
Then there is a fast decreasing function $\tilde \omega_K$ such that
the  correlation functions $m^{\bk k1{p+q}}(\bk x1{p+q};n)$ of $\tilde\P_n$ satisfy~\eqref{dbtpptruncmixing} uniformly in $n \in \N$ and $\bk k1{p+q}$ given $k_1+\ldots+k_{p+q}=K$. The function $\tilde \omega_K$  depends on the value of $K$, the function $\hat\omega_{K}$,  and the constants  $M^\xi_{K+\epsilon}$, 
$C_K$, $\kappa_0$   related to $\xi$ and $\tP$.
\end{lemma}

\begin{proof}
For $ i \in \{1,...,p + q\}$ define the functions 
$g_{i}(u_i):= u_i^{k_i}$ on $\R$.  Consider the functions  $h_1(\bk u1p) = \Pi_{i  = 1}^p g_i(u_i)$  and  $h_2(\bk u{p+1}{p+q}) = \Pi_{i  = p + 1}^{p + q} g_i(u_i)$, respectively, on $\R^p$ and $
\R^q$.  By the definition~\eqref{eqn:mixedmoment} of correlation functions~$m^{\bk k1{p+q}}$ 
the left-hand-side of~\eqref{dbtpptruncmixing}
corresponds to 
\begin{align}\label{e.local}
& \Bigl|\E_{\bk x1{p + q}} [h_1(\xi_{1,n},\ldots,\xi_{p,n}) h_2(\xi_{p+1,n},\ldots,\xi_{p+q,n})] \\
&-
\E_{\bk x1p} [h_1(\xi_{1,n},\ldots,\xi_{p,n})] \rho^{(p)} (\bk x1p) \E_{\bk x{p + 1}{p +
    q}  }[h_2(\xi_{p+1,n},\ldots,\xi_{p+q,n})] \rho^{(q)} (\bk x{p + 1}{p +
  q})\Bigr|. \nonumber
\end{align}
To use~\eqref{e.mppmixing} with $U(x_i)=\xi_{i,n}$ we have to truncate $h_1, h_2$ in a Lipschitz way to fulfill the role of $f,g$ in the \BL-mixing condition assumed to hold uniformly for all $\tP_n$, $n \in \N$. For this purpose
 we choose  $\epsilon > 0$ as in~\eqref{e.moment-cond-Lemma}, we define the function $a(s)= a_{p+q}(s):= \max \{\hat\omega_{p+q}(s)^{-(1+\epsilon)/(2+3\epsilon)},1\}$, and we truncate the functions $g_i, 1 \leq i \leq p$ as follows :
\begin{equation}
\bar g_{i} :=\begin{cases}  g_i  & \quad \text{when\ } |g_i| \leq a(s)^{1/p} \\
\sgn(g_i)a(s)^{1/p}  & \quad\text{when\ } |g_i| > a(s)^{1/p}.
\end{cases}
\end{equation}
Similarly, we 
truncate the functions $g_i, p + 1 \leq i \leq p + q$ as follows :
\begin{equation}
\bar g_{i} :=\begin{cases}  g_i  & \quad \text{when\ } |g_i| \leq a(s)^{1/q} \\
\sgn(g_i)a(s)^{1/q}  & \quad\text{when\ } |g_i| > a(s)^{1/q}.
\end{cases}
\end{equation}

We assert that each $\bar g_{i}, 1 \leq i \leq p,$ is Lipschitz  with Lipschitz constant $k_i \cdot a(s)^{(k_i-1)/(pk_i)}$.  Indeed notice for $|x|, |y| \leq a(s)^{1/(pk_i)}$ by the mean value theorem
$$
| \bar g_{i}(x) - \bar g_{i}(y) | \leq |x - y| \cdot 
g_i'(a(s)^{1/(pk_i)})
= k_i \cdot a(s)^{(k_i-1)/(pk_i)}|x - y|.
$$
On the other hand, if $|x| \leq a(s)^{1/(pk_i)} \leq |y|$ then
$$
| \bar g_{i}(x) - \bar g_{i}(y) | = |x^{k_i} - \sgn(y)a(s)^{1/p} | 
\leq k_i \cdot a(s)^{(k_i-1)/(pk_i)}   |x - y|,
$$
which shows the assertion. Likewise, each $\bar g_{i}, p + 1 \leq i \leq p + q,$ is Lipschitz  with Lipschitz constant $k_i \cdot a(s)^{(k_i-1)/(qk_i)}$.

Put
$$
\bar h_{1} = \Pi_{i  = 1}^p \bar{g}_i(u_i),  \quad \bar h_{2} = \Pi_{i  = p + 1}^{p + q} \bar{g}_i(u_i).
$$

Observe, $\bar h_{1} \leq a(s)$  is a bounded, Lipschitz function (with respect to the $\ell^1$ norm on $\R^p$) with Lipschitz constant 
bounded by $ a(s)^{(p-1)/p}\sum_{i=1}^p  k_i \cdot a(s)^{(k_i-1)/(pk_i)}= a(s)\sum_{i=1}^p k_i\le Ka(s)$.\footnote{\label{f.Lip-Rd} Indeed, a $p$-fold product of Lipschitz functions on $\R$ with respective Lip constants $L_1,...,L_p$ and common sup norm bound $S$ is a Lipschitz function on
$\R^p$ with Lipschitz constant equal to $\sum_{i = 1}^p L_i S^{p - 1}$ with respect to the $\ell^1$ norm on $\R^p$.}
Similarly $\bar h_{2} \leq a(s)$,  and is Lipschitz with Lipschitz constant bounded by 
$Ka(s)$. Consequently 
 $\bar h_1/(K a(s))$ and $\bar h_2/(K a(s))$ are  bounded Lipschitz functions useable in ~\eqref{e.mppmixing}. 

In the  evaluation of the truncation error we  simplify the notation by writing  $H_1:=h_1(\xi_{1,n},\ldots,\xi_{p,n})$,
$H_2:=h_2(\xi_{p+1,n},\ldots,\xi_{p+q,n})$ and similarly for  
$\bar H_1:=\bar h_1(...)$ and
$\bar H_2:=\bar h_2(...)$.
With this notation, we write the first expectation in~\eqref{e.local} as the sum of three terms
%\blue{with $V^p:=\bk V1p$ and $V^q:=\bk V{p+1}q$ to simplify the notation}
\begin{align*}
& T_1:=  \rho^{(p + q)} (\bk x1{p + q}) \E_{\bk x1{p + q}} [H_1(H_2 - \bar H_2)],\\
& T_2:=  \rho^{(p + q)} (\bk x1{p + q}) \E_{\bk x1{p + q}} [( H_1 - \bar H_1) \bar H_2], \\
& T_3:=  \rho^{(p + q)} (\bk x1{p + q}) \E_{\bk x1{p + q}} [\bar H_1 \bar H_2].
\end{align*}
Write the second expectation in~\eqref{e.local} as the sum of four terms:
\begin{align*}
& T_4:=  \rho^{(p)} (\bk x1p)  \rho^{(q)} (\bk x{p + 1}{p + q})  \E_{\bk x1p} [H_1] \E_{\bk x{p + 1}{p + q}} [H_2 - \bar H_2], \\
& T_5:=  \rho^{(p)} (\bk x1p)  \rho^{(q)} (\bk x{p + 1}{p + q})  \E_{\bk x1p} [ H_1 - \bar H_1] \E_{\bk x{p + 1}{p + q}} [ \bar H_2],\\
& T_6:=  \rho^{(p)} (\bk x1p)  \rho^{(q)} (\bk x{p + 1}{p + q})  \E_{\bk x1p} [ \bar H_1]  \E_{\bk x{p + 1}{p + q}} [ \bar H_2].
\end{align*}
Applying  \eqref{e.mppmixing} to the Lip(1)  functions $\bar H_1/(K a(s))$ and  $\bar H_2/(K a(s))$,   we find that the difference of terms $T_3$ and $T_6$ is bounded by
$$
|T_3 - T_6| \leq K^2 a(s)^{2} \hat\omega_{p+q}(s) \leq K^2
\hat\omega_{p+q}(s)^{\epsilon/(2+3\epsilon)},
$$
by  choice of $a(s)$.    %
We now control the terms $T_1, T_2, T_4, T_5,$ as follows. Put
$\widetilde{M}_{K,\epsilon}:= \widetilde{M}_{K,\epsilon, p + q}$.
By the  H\"older and the Markov inequalities and the bound on $\kappa_p$  in~\eqref{e.correlation-functions-bound} we have 
\begin{align*}
|T_1| & \leq \kappa_{p + q}  \E_{\bk x1{p + q}}[| H_1 H_2 |\1{\exists i\in\{p+1,\ldots,p+q\}: |g_i(\xi_{i,n})|>a(s)^{1/q}}] \\
& \leq \sum_{i=p+1}^{p + q} \kappa_{p + q}  \E_{\bk x1{p + q}}[|H_1H_2|^{1+\epsilon}]^{1/(1+\epsilon)} a(s)^{- \epsilon/(q(1+\epsilon))}
\E_{\bk x1{p + q}}[ |g_i(\xi_{i,n})| ]^{\epsilon/(1+\epsilon)}\\
&\le q\kappa_{p+q} \widetilde{M}_{K,\epsilon}^{2}\,a(s)^{- \epsilon/(q (1+\epsilon))}, 
\end{align*}
where we use  
the moment condition~\eqref{e.moment-cond-Lemma}  twice (for 
$\E_{\bk x1{p + q}}[|H_1H_2|^{1+\epsilon}]$
and $\E_{\bk x1{p + q}}[ |g_i(\xi_{i,n})|]$) in the third
inequality as follows. Specifically, we  apply the  H\"older inequality to obtain 
\begin{align*}
\E_{\bk x1{p + q}}[|H_1H_2|^{1+\epsilon}]
&\le \E_{\bk x1{p + q}}\Bigl[\prod_{i=1}^{p+q}|\xi_{i,n}|^{(1+\epsilon)k_i}\Bigr] \le \prod_{i=1}^{p+q}\E_{\bk x1{p + q}}\Bigl[\max(1,|\xi_{i,n}|)^{(1+\epsilon)K}\Bigr]^{k_i/K} \\
&\le \prod_{i=1}^{p+q}(\widetilde{M}_{K,\epsilon})^{k_i/K} = \widetilde{M}_{K,\epsilon}
\end{align*}
and 
\begin{align*}
\E_{\bk x1{p + q}}[ |g_i(\xi_{i,n})|] & \le \E_{\bk x1{p + q}}[\max(1,|\xi_{i,n}^{k_i}|)] \le  \E_{\bk x1{p + q}}[\max(1,|\xi_{i,n}|)^{(1+\epsilon)K}]^{k_i/(K(1+\epsilon))} \\
&= \big( \widetilde{M}_{K,\epsilon}\big)^{k_i/(K(1+\epsilon))} \le \widetilde{M}_{K,\epsilon}. 
\end{align*}
By symmetry (interchanging the functions $h_1$ and
$h_2$ as well as $p$ and $q$ and using the bound $|\bar H_2|\le |H_2|$)
we also have  
$$
|T_2| \le p\kappa_{p+q} \widetilde{M}_{K,\epsilon}^{2} a(s)^{-\epsilon/(q(1+\epsilon))}.
$$
Applying the H\"older and Markov inequalities as when  bounding $T_1$,  we obtain 
\begin{align*}
|T_4| & \leq \kappa_{p} \kappa_{q}   \E_{\bk x1p} [|H_1|]
\E_{\bk x{p + 1}{p + q}}  [|H_2| \1{\exists i\in\{p+1,\ldots,p+q\}: |g_i(u_i)|>a(s)^{1/q}} ]\\
&\leq q\kappa_{p} \kappa_{q} \widetilde{M}_{K,\epsilon}^{2} a(s)^{- \epsilon/(q(1+\epsilon))}, 
\end{align*}
where we used
\begin{align*}
\E_{\bk x1{p}}[|H_1|] =\E_{\bk x1{p }}\Bigl[\prod_{i=1}^{p}|\xi_{i,n}|^{k_i}\Bigr]
&\le \E_{\bk x1{p }}\Bigl[\prod_{i=1}^{p}|\max(1,|\xi_{i,n}|)|^{(1+\epsilon)k_i}\Bigr]\\
&\le\prod_{i=1}^{p}\E_{\bk x1{p }}\Bigl[|\max(1,|\xi_{i,n}|)|^{(1+\epsilon)(k_1+\ldots+k_p)}\Bigr]^{\frac{k_i}{k_1+\ldots+k_p}}\\
&\le\prod_{i=1}^{p}\E_{\bk x1{p + q}}\Bigl[|\max(1,|\xi_{i,n}|)|^{(1+\epsilon)K}\Bigr]^{\frac{k_i}{k_1+\ldots+k_p}}\\
& \le \prod_{i=1}^{p}(\widetilde{M}_{K,\epsilon})^{\frac{k_i}{k_1+\ldots+k_p}} = \widetilde{M}_{K,\epsilon}, \\
\end{align*}
and mutatis mutandis for $\E_{\bk x{p + 1}{p + q}}  [|H_2|^{(1+\epsilon)}]\le \widetilde{M}_{K,\epsilon}$.
By symmetry,  
$$|T_5| \leq p \kappa_{p} \kappa_{q} \widetilde{M}_{K,\epsilon}^{2}  a(s)^{- \epsilon/(p(1+\epsilon))}
.$$
By~\eqref{e.correlation-functions-bound} we have 
$\kappa_p\kappa_q\le  (p + q)^2 C_{p + q}^2  \kappa_0^{p + q}$. Collecting all estimates above we the bound of  expression in~\eqref{e.local} and using that $a(s) \geq 1$ and $a(s) \geq \hat\omega_{p+q}(s)^{-(1+\epsilon)/(2+3\epsilon)}$, we derive 
\begin{align*}
 & |T_3 - T_6| + |T_1| + |T_2| + |T_4| + |T_5| \\
 &\le (K^2 + 4(p + q)^{3} C_{p + q}^2  \kappa_0^{p + q} \widetilde{M}_{K,\epsilon}^{2}) a(s)^{- \epsilon/(K(1+\epsilon))}\\
&\le (K^2 + 4 K^{3} C_{K}^2  \kappa_0^{K} \widetilde{M}_{K,\epsilon}^{2}) \hat\omega_{p+q}(s)^{\epsilon/(K(2+3\epsilon))}
  \end{align*}
where the first inequality uses  the fact that $\widetilde{M}_{K,\epsilon} \geq 1$
and the second uses  $p+q\le K$, and the increasing
property of $\hat \omega_k(s)$ in the argument $k$. The fast decay of $(K^2 + 4 K^3 C_{K}^2  \kappa_0^{K}\widetilde{M}_{K,\epsilon}^{2}) \hat\omega_{p+q}(s)^{\epsilon/(K(2+3\epsilon))}$ follows from the fast decay of $\hat \omega$, completing the proof of Lemma~\ref{l:non-bdd-mppmixing}.
\end{proof}

%\fi

\begin{proof}[Proof of Theorem \ref{t:clt_linear_marks}]
We deduce  this result 
from Theorem ~\ref{l.CLT-Cumulants} as follows.
 Recall that the  correlation functions  $m^{\bk k1{p+q}}(\bk x1{p+q};n)$  of the atomic measures
$\amu_n:= \sum_{x \in \P_n} \xi_{i,n} \delta_{x }$ 
 are defined  at~\eqref{eqn:mixedmoment-general}.  (In contrast to $\hat\mu_n^\xi$,  the measure
 $\amu_n$ is not scaled to $W_1$.) 
 Under the moment assumption \eqref{e:xinpmom},  the H\"{o}lder
inequality implies that the correlation functions $m^{\bk k1{p}}(\bk x1{p};n)$ exist.
Moreover,  the scores   $\xi_{i,n}$ have finite moments as in
  \eqref{e.moment-cond-Lemma} for any $\epsilon > 0$ (uniformly for $\tP_n$) as $\widetilde{M}_{K,\epsilon,p+q} \leq M^{\xi}_{K(1+\epsilon)}$ and that the moment assumption \eqref{e:xinpmom} holds for all $p \in \N$. 
By assumption, we have that $(\dbtilde \P_n)_{n \in \N}$ are fast \BL-mixing as in Definition \ref{d.omegahmixing}(iii).
Hence, applying 
  Lemma \ref{l:non-bdd-mppmixing} with
  $h_1(\bk t1p) := \prod_{i=1}^pt_i^{k_i}$ and $h_2(\bk t1q) :=
  \prod_{i=1}^qt_i^{k_i}$ for $p,q,k_1,\ldots,k_{p+q} \in
  \mathbb{N},  t_i \in \R$,
 we obtain  the fast decay of the correlation
  functions $m^{\bk k1p}(\cdot;n)$  uniformly in $n \in \N$, i.e., for any $p,q,k_1,\ldots,k_{p+q} \in \mathbb{N}$ and $x_1,\ldots,x_{p+q} \in \mR^d$, we have 
    \be \label{e.dbtpptruncmixing}
 | m^{\bk k1{p+q}}(\bk x1{p+q};n) -  m^{\bk k1{p}}(\bk x1p;n) m^{\bk k{p+1}{p+q}}(\bk x{p+1}{p+q};n)| \leq \hat \omega_K(s)
\ee
where $s := d(\bk x1p,\bk x{p+1}{p+q})$, $K := \sum_{i=1}^{p+q}k_i$, and where $\hat \omega_{K}$ is a  fast decreasing function for all $K \in \N$.

Thus the correlation
  functions $m^{\bk k1p}(\cdot;n)$ are fast  mixing, uniformly in $n \in \N$, in accordance with statement (i) of  Lemma \ref{l.moments-cumulants} 
by taking $\tilde C = 1$ and $\phi_{k_1,...,k_p}(s)= \hat \omega_{K}(s)$ where $K = \sum_{i = 1}^p k_i$.  
Hence, we obtain the central limit theorem \eqref{e.CLT-Th2general} 
from Theorem ~\ref{l.CLT-Cumulants}  with $g_n(x)=f(n^{-1/d}x)$ and $f\in\B(W_1)$,
using the $p$-moment condition~\eqref{e:xinpmom} to verify \eqref{e.Moment-general}
  and the assumed lower bound for $\Var{\mu_n^\xi(f)}$.
 \end{proof}

\subsection{Distributional limits for \BL-localizing score functions}
\label{s:stab_lemma_BL}

In this section, we justify the existence of distributional limits for $\xi(\tx_1,\tP_n)$ as $n \to \infty$ under Palm distributions as in \eqref{e.xil-0-multi} and \eqref{e.xil-xy-multi}, of which  \eqref{e.xil-0} and \eqref{e.xil-xy} are special cases. These are essential  when proving  expectation and variance asymptotics (Proposition~\ref{expvar} and Theorem \ref{t:multcltmarkedpp}) in the next subsections.

We operate under the premise that the score function $\xi$ exhibits fast \BL-localization across all finite windows of $\tP$, as  in Definition \ref{def.Lp-stabilizing_marking}~\ref{i.BL-localizing-windows}. We start  by presenting key findings concerning the distributional limits resulting from \BL-localization. While these limits hold for general marking functions $\xi$ taking values in a Polish space (the setting in which \BL-localization is originally defined), we shall in the present context---where moments are considered too--- restrict throughout this section to real-valued scores $\xi$.

Recall the notation
$$
\xitr(\tx_1, \tP)= \xi( \tx_1,\tP\cap B_r(x_1) ) \, \text{and} \ \,  \bk{\xitr}1p(\tx, \tP)=\Bigl(\xi(\tx_1,\tP\cap B_r(x_1)),\ldots,\xi(\tx_p,\tP\cap B_r(x_p)\Bigr).
$$
\begin{lemma}[Distributional limits for vectors of \BL-localizing scores]
\label{l:BL-limits} 
Let  $\tP$  be a marked point process on $\R^d\times\K$.  Let $\xi : \R^d \times \K \times \hat\cN_{\R^d \times \K} \to \R$ be an
real-valued score function, and assume that $\xi$ is \BL-localizing on all
finite windows of $\tP$ as in~\eqref{Lp-stab} for a given $p \in \N$
(tacitly assuming the existence of the uniformly bounded correlation function $\rho^{(p)}(\bk{x}{1}{p})$).
\begin{enumerate}[wide,label=(\roman*),labelindent=0pt]
\item \label{i.BL-r-limit}
For a given vector  $\bk x1p \in(\R^d)^p$, 
the  random vectors  $\bk{\xitr}1p(\tx, \tP)$ 
(under the Palm distributions ${\sP}_{\bk x1p}$ of $\tP$) 
converge in the $d_{\BL, \bk{x}{1}{p}}$ metric (and hence in law) as $r\to\infty$ to a probability distribution on $\R^p$, represented here by the random vector  $\xil_{\bk x1p}=(\xil_{\bk x1p}[1],\ldots,\xil_{\bk x1p}[p])\in \R^p$, i.e.,
$$\bk{\xitr}1p(\tx, \tP)
\xrightarrow[r \to \infty]{\BL,\bk x1p} \xil_{\bk x1p}
$$
and the speed of convergence is bounded by the same decreasing function $\varphi_p$ as in~\eqref{Lp-stab}, namely
\begin{align} \label{e.BL-r-rate}
 \sup_{x_1,\ldots,x_p \in \R^d} 
 d_{\BL, \bk x1p}(\bk{\xitr}1p(\tx, \tP), \xil_{\bk x1p})  \leq 4\varphi_p(r), \ r > 0.
\end{align}
\item \label{i.BL-r-limit-kernel}
The probability distributions represented by $\xil_{\bk x1p}$ are probability kernels from $(\R^d)^p$ to $\R^p$; formally, for any measurable function $f$ on $\R^p$, the function $(\R^d)^p\owns\bk x1p\mapsto\E[f(\xil_{\bk x1p})]\in\R$ is measurable.
\item  \label{i.BL-n-limit} The convergence holds also for $\bk{\xi}1p(\tx,\tP_n)$:
\begin{equation}\label{e:BL-n-limit}
\bk{\xi}1p(\tx,\tP_n) \xrightarrow[n \to \infty]{\BL,\bk x1p} \xil_{\bk x1p}
\end{equation}
with the speed of convergence bounded by the same decreasing function $\varphi_p$ of the  distance $ d(\bk x1p,\partial W_n)$ from  $\{x_1,\ldots, x_p\}\subset W_n$ to the boundary of $W_n$:
\begin{equation}\label{e.BL-n-rate} d_{\BL, \bk x1p}(\bk{\xi}1p(\tx, \tP_n),\xil_{\bk x1p})  \leq
 6\varphi_{p}(d(\bk x1p,\partial W_n)).
 \end{equation}
\item \label{i.BL-stat} If $\tP$ is stationary and $\xi$ is translation invariant, then the  distribution of $\xil_{\cdot}$ is translation invariant:
i.e., for all $y\in\R^d$
 $$\xil_{\bk x1p}
\overset{\text{law}}{=}\xil_{\bk x1p-y}=\xil_{(x_1-y,\ldots,x_p-y)}.
$$

\end{enumerate}

\begin{enumerate}[wide,label=(\roman*),  labelindent=0pt]
\setcounter{enumi}{4}
\item \label{i.BL-2limits} When  $\xi$ is defined on  $\R^d \times \K \times
  \cN_{\R^d \times \K}$ and satisfies \BL-localization~\eqref{Lp-stab-infinite}  for some $p\ge 1$, then for $\{x_1,\ldots, x_p\}\subset W_n$ we have
\begin{equation}\label{e:BL-2limits}
d_{\BL, \bk x1p}( \bk{\xi}1p(\tx, \tP_n), \bk{\xi}1p(\tx, \tP))\le 12\varphi_{p}(d(\bk x1p,\partial W_n)),
\end{equation}
and hence, under $\mathbb{P}_{\bk x1p}$, $\bk{\xi}1p(\tx, \tP)$ has the law of 
$\xil_{\bk x1p}$.
\end{enumerate}
\end{lemma}

\begin{proof}
\noindent \ref{i.BL-r-limit}
Taking $0<r\le r'$ and $n$ large enough such that  $B_{r'}(\bk x1p) :=\bigcup_{i=1}^p B_{r'}(x_i) \subset W_n$,
by the triangle inequality,  \BL-localization \eqref{Lp-stab-infinite} and the decreasing property of $\varphi_p$, we obtain
\begin{align}\label{e.BL-r-limit}
&\hspace{-3em}d_{\BL, \bk x1p}\Bigl(\bk{\xitr}1p(\tx, \tP), \bk{\xi^{(r')}}1p(\tx, \tP)\Bigr)\\
&\le
d_{\BL, \bk x1p}\Bigl(\bk{\xitr}1p(\tx, \tP_n), \bk{\xi}1p(\tx,\tP_n)\Bigr)+
d_{\BL, \bk x1p}\Bigl(\bk{\xi^{(r')}}1p(\tx, \tP_n), \bk{\xi}1p(\tx,\tP_n)\Bigr) \nonumber\\
&\le 4\varphi_p(r),\nonumber
\end{align}
where $\varphi_p(r)\to0$ as $r\to\infty$.
Under the Palm distributions $\mathbb{P}_{\bk x1p}$, this shows that any subsequence of random vectors $\bk{\xi^{(r_k)}}1p(\tx, \tP)$, where $r_k \to \infty$, is a Cauchy sequence in the $d_{\BL, \bk{x}{1}{p}}$ distance.  This justifies the first statement in~\ref{i.BL-r-limit}. Indeed, $\mathbb{R}^p$ is a complete separable metric space, and hence the weak convergence of probability measures on it is completely metrizable by the  $d_{\BL, \bk{x}{1}{p}}$ metric (see \cite[Theorem 8.10.43]{bogachev2007measure}). To  bound the speed of convergence, we  take $r,r'$ as in~\eqref{e.BL-r-limit} and 
by the triangle inequality we obtain
\begin{align*}
&\hspace{-3em}d_{\BL, \bk x1p}\Bigl(\bk{\xitr}1p(\tx, \tP),\xil_{\bk x1p}\Bigr)\\
&\le
d_{\BL, \bk x1p}\Bigl(\bk{\xitr}1p(\tx, \tP), \bk{\xi^{(r')}}1p(\tx,\tP)\Bigr)+
d_{\BL, \bk x1p}\Bigl(\bk{\xi^{(r')}}1p(\tx, \tP),\xil_{\bk x1p}\Bigr)
\\
&\le 4\varphi_p(r)+ d_{\BL, \bk x1p}\Bigl(\bk{\xi^{(r')}}1p(\tx, \tP),\xil_{\bk x1p}\Bigr).
\end{align*}
When   $r'\to\infty$, the last term goes to~0  by the previous  statement and  we obtain~\eqref{e.BL-r-rate}.
\medskip

\noindent\ref{i.BL-r-limit-kernel} This follows from the fact that Palm distributions are probability kernels. Indeed, for a bounded measurable function \(f:\R^p\to\R\), the mapping \((\R^d)^p \owns \bk{x}{1}{p} \mapsto \E_{\bk{x}{1}{p}} [f(\bk{\xitr}{1}{p}(\tx, \tP))] \in \R\) is measurable. Since the law of \(\bk{\xitr}{1}{p}(\tx, \tP)\) under $\mP_{\bk x1p}$ converges in the $d_{\BL, \bk x1p}$ distance, and hence weakly, to \(\xil_{\bk{x}{1}{p}}\) for all \(\bk{x}{1}{p}\), the aforementioned mapping converges point-wise to \((\R^d)^p \owns \bk{x}{1}{p} \mapsto \E[f(\xil_{\bk{x}{1}{p}})] \in \R\). Consequently, the limit  is measurable as a pointwise limit of measurable functions. Extensions to unbounded measurable functions \(f\) are straightforward.

\medskip

\noindent  \ref{i.BL-n-limit}
Taking $r=d(\bk x1p,\partial W_n)$, by the triangle inequality, \BL-localization \eqref{Lp-stab} and~\eqref{e.BL-r-rate} we obtain
\begin{align*}
&\hspace{-3em}d_{\BL, \bk x1p}( \bk{\xi}1p(\tx, \tP_n),\xil_{\bk x1p})  \\
& \leq d_{\BL, \bk x1p}\Bigl(\bk{\xi}1p(\tx, \tP_n),\bk{\xitr}1p(\tx,\tP_n)\Bigr)+
d_{\BL, \bk x1p}\Bigl(\bk{\xitr}1p(\tx, \tP_n),\xil_{\bk x1p}\Bigr)  \\
&= d_{\BL, \bk x1p}\Bigl(\bk{\xi}1p(\tx, \tP_n),\bk{\xitr}1p(\tx, \tP)\Bigr)+
d_{\BL, \bk x1p}\Bigl(\bk{\xitr}1p(\tx, \tP),\xil_{\bk x1p}\Bigr)
\\
&\le 6\varphi_p(r)
\end{align*}
with the equality due to  $\{x_1,\ldots,x_p\}\subset W_n$, which gives~\eqref{e.BL-n-rate}  and consequently also proves the statement of the convergence \eqref{e:BL-n-limit}.
\medskip 

\noindent \ref{i.BL-stat} By the stationarity of $\tP$ and translation invariance of $\xi$, the law of $\bk{\xitr}1p(\tx, \tP)$  (under $\P_{\bk x1p}$) is invariant by the translation of the location $\bk x1p$ with respect to any vector $y \in\R^d$. Hence, the statement of~\ref{i.BL-stat} follows from the weak convergence in Part~\ref{i.BL-r-limit}.
\medskip

\noindent \ref{i.BL-2limits} Similarly as in the proof of  part ~\ref{i.BL-n-limit},
take $r=d(\bk x1p,\partial W_n)$ and use the triangle inequality
\begin{align*}
&d_{\BL, \bk x1p}( \bk{\xi}1p(\tx, \tP_n), \bk{\xi}1p(\tx, \tP))  \\
&\leq d_{\BL, \bk x1p}\Bigl(\bk{\xi}1p(\tx, \tP_n),\xil_{\bk x1p} \Bigr)
+ d_{\BL, \bk x1p}\Bigl(\xil_{\bk x1p},\bk{\xitr}1p(\tx, \tP)\Bigr)
+d_{\BL, \bk x1p}\Bigl(\bk{\xitr}1p(\tx, \tP),\bk{\xi}1p(\tx, \tP)\Bigr).
\\
&\le 12\varphi_p(r),
\end{align*}
where the last inequality follows from~\eqref{e.BL-n-rate}, \eqref{e.BL-r-rate}  and \BL-localization~\eqref{Lp-stab}.
\end{proof}

Now, we formulate and prove some results regarding the moment conditions. The key point is that fast \BL-localization of $\xi$ together with $(p + \epsilon)$-moments implies $p$-moments of $\xil$. It's important to note this higher moment condition in the \BL \,setting in comparison to the stopping set setting (cf. Lemma~\ref{l.xi-infnity}~\ref{i.SS-limit-moment}).
\begin{lemma}[Moment bounds on the distributional limits of \BL-localizing scores]
\label{l:BL-limits-moments}
Let  $\tP$  be a marked point process on $\R^d\times\K$ with non-null, finite intensity $\rho$.  Let $\xi : \R^d \times \K \times \hat\cN_{\R^d \times \K} \to \R$ be an real-valued score function which is \BL-localizing on all finite windows  for Palm distributions of some order $q\in \N$ (that is~\eqref{Lp-stab} holds  with  $p$ replaced by  $q$,  with  uniformly bounded correlation function $\rho^{(q)}$).
\begin{enumerate}[wide,label=(\roman*),labelindent=0pt]
\item \label{i.BL-limit-p-moment} 
If $\xi$ satisfies the $(p + \epsilon)$-moment condition on finite windows~\eqref{e:xinpmomcopied} for some $p\geq q$ and some $\epsilon > 0$, then the random variables $\xil_{\bk{x}{1}{q}}[i],1  \leq i \leq q$ have uniformly bounded $p$-moments, namely 
\begin{equation}\label{e.BL-limit-p-moment}
\sup_{x_1,\ldots,
  x_{q} \in \R^d}\max_{i=1,\ldots,q} \sE[\max(1,|\xil_{ \bk x1q}[i]|^p)]\le M^\xi_p.
  \end{equation}
\item \label{i.BL-limit-p-bound}
Furthermore, under the assumptions of part ~\ref{i.BL-limit-p-moment} above, we also obtain the speed of convergence of (natural) moments $p\in\N$; namely for $x_1,\ldots,x_q \in W_n$ we have
\begin{equation}\label{e.BL-limit-p-bound}
\max_{i=1,\ldots,q} 
      \Big | \sE_{\bk x1{q}}[ (\xi(\tx_i,\tP_n))^p]-\sE[(\xil_{\bk x1q}[i])^p]\Bigr| \le  (2M_{p+\epsilon/2}^\xi+6p)\varphi_{q} (d(\bk x1q,\partial W_n))^{\epsilon/(2p+\epsilon)},
  \end{equation}
  where $d(\bk x1q,\partial W_n)$ is the distance  between $x_1,\ldots,x_q \in W_n$ and the boundary of $W_n$.
  \item \label{i.BL-correlation-bound}
  If $\xi$ is \BL-localizing on all finite windows    of Palm order $q=2$ as in~\eqref{Lp-stab} and  satisfies the $(2+\epsilon)$-moment condition on finite windows~\eqref{e:xinpmomcopied} for some $\epsilon > 0$, then the  vector $\xil_{\bk{x}{1}{2}}$ has uniformly bounded correlations, i.e.,  
\begin{equation}\label{e.BL-correlation-bound}
\sup_{x_1,x_2 \in \R^d} 
 \sE[\max(1,|\xil_{\bk x12}[1]\xil_{\bk x12}[2]|)] 
  \le M^\xi_2.
  \end{equation}
\item \label{i.BL-correlation-rate}
Furthermore, under the assumptions of part ~\ref{i.BL-correlation-bound}, the correlations have the speed of convergence:
\begin{equation}\label{e.BL-correlation-rate}
  \Big | \sE_{\bk x1{2}}[\xi(\tx_1,\tP_n)\xi(\tx_2,\tP_n)]-\sE\bigl[\xil_{\bk x12}[1]\xil_{\bk x12}[2]\bigr]\Bigr| \le  (4M^\xi_{2+\epsilon/2} + 12)\varphi_{2}(d(\bk x12,\partial W_n))^{\epsilon/(4+\epsilon)}\,,
  \end{equation}
  where $ d(\bk x12,\partial W_n)$ is the distance between $\{x_1, x_2\}\subset W_n$ and the boundary of $W_n$.
  \end{enumerate}
\end{lemma}

\begin{proof}
\noindent \ref{i.BL-limit-p-moment}
It is well known that convergence in distribution of random variables, along with a uniform bound on $(p+\epsilon)$- moments, implies the convergence of $p$-moments. In our specific case, as stated in Lemma~\ref{l:BL-limits}~\ref{i.BL-n-limit}, this implies 
$$\E[\max (1,|\xil_{\bk x1q}[i]|^p)]=\lim_{n\to\infty}\sE_{\bk x1q}[\max (1,|\xi(\tx_i,\tP_n)|^p)]\le M^\xi_p.$$ 
\noindent \ref{i.BL-limit-p-bound}
We use a truncation as in the proof of Lemma~\ref{l:non-bdd-mppmixing}. For a function $ f_a: \mathbb{R} \rightarrow \mathbb{R} $ where $a = a(r) \geq 1$ (the specific value of $a$ will be determined later), consider the following definition 
$$
f_a(t):=\begin{cases}  
t^p& \quad \text{when\ } |t| \leq  a\\
\sgn(t^p)a^p  & \quad\text{when\ } |s|> a.
\end{cases}
$$
The function $f_a/(pa^{p})$  is a bounded,  Lipschitz(1) function, and thus by the triangle inequality we have 
\begin{align}
&\hspace{-3em}\Big | \sE_{\bk x1{q}}[ (\xi(\tx_i,\tP_n))^p]-\sE[(\xil_{\bk x1q}[i])^p]\Bigr| \no \\
\le & \sE_{\bk x1{q}}[ |\xi(\tx_i,\tP_n)|^p\1{|\xi(\tx_i,\tP_n)|>a}]+
\sE[ |\xil_{\bk x1 q}[i]|^p\1{|\xil_{\bk x1q}[i]|>a}] \no\\
\label{e.xipmom_trunc} &+
pa^{p}d_{\BL,\bk x1q}\Bigl(\xi(\tx_i,\tP_n),\xil_{\bk x1q}[i]\Bigr).
\end{align}
Using \eqref{e.BL-n-rate}, the last term is bounded by \( 6pa^{p}\varphi_{q}(r) \)
where $r=d(\bk x1q,\partial W_n)$.
Next, let's choose \( p' = p + \epsilon/2 \). Both random variables \( \xi(\tx_i,\tP_n) \) and \( \xil_{\bk{x}{1}{q}}[i] \) have \( p' \)-moments bounded by \( M_{p'}^\xi \) (for the former by assumption and for the latter by the statement in point \ref{i.BL-limit-p-moment} with \( p \) replaced by \( p' \)). 
Therefore, using Hölder's inequality, both truncated moments are bounded by \( M_{p'}^\xi a^{-\epsilon/2} \) and consequently
\[ \Big| \mathbb{E}_{\bk x1q}[ (\xi(\tx_i,\tP_n))^p] - \mathbb{E}[(\xil_{\bk x1q}[i])^p] \Big| \le  2M_{p'}^\xi a(r)^{-\epsilon/2} + 6pa(r)^{p}\varphi_{q}(r) .\]
Setting \( a(r) = \varphi(r)^{-2/(2p+\epsilon)} \), the right-hand side simplifies to:
\[ (2M_{p+\epsilon/2}^\xi + 6p)\varphi_{q}(d(\bk x1q,\partial W_n))^{\epsilon/(2p+\epsilon)}. \]
This completes the statement of part \ref{i.BL-limit-p-bound}.
\medskip 

\noindent \ref{i.BL-correlation-bound}
Using the same reasoning as in the proof of part~\ref{i.BL-limit-p-moment}, the convergence in distribution of the vector  $\bk{\xi} 12(\tx,\tP_n)$ to $\xil_{\bk x12}$, as per the statement of Lemma~\ref{l:BL-limits}~\ref{i.BL-n-limit} and the continuous mapping theorem, implies convergence in distribution of the products $\xi(\tx_1,\tP_n)\xi(\tx_2,\tP_n)$ to $\xil_{\bk x12}[1]\xil_{\bk x12}[2]$. Along with uniform $(2+\epsilon)$-moments~\eqref{e:xinpmomcopied}, by H\"older's inequality, these products possess bounded $(1+\epsilon)$-moments: 
$$\sup_{n \in \N}\sup_{\bk x12 \in W_n^2}\sE_{\bk x12}[|\xi(\tx_1,\tP_n)\xi(\tx_2,\tP_n)|^{1+\epsilon/2}]<\infty.$$
This implies the convergence of expectations $\sE_{\bk x12}[|\xi(\tx_1,\tP_n)\xi(\tx_2,\tP_n)|]$ to $\sE[|\xil_{\bk x12}[1]\xil_{\bk x12}[2]|)]$ as $n\to\infty$. The explicit bound follows from~\eqref{e:xinpmomcopied} with $p=1$, using H\"older's inequality.
\medskip

\noindent \ref{i.BL-correlation-rate}
We use the same reasoning as in the proof of part~\ref{i.BL-limit-p-bound}, with the truncation of the product of functions as in the proof of Lemma~\ref{l:non-bdd-mppmixing} with $p=2$ and $k_i=1$, $i=1,2$. Specifically, consider a function defined as:
$$
f_a(t):=\begin{cases}  
t & \text{when\ } |t| \leq  a^{1/2}\\
\sgn(t)a^{1/2} & \text{when\ } |t|> a^{1/2}.
\end{cases}
$$
where the specific value of $a=a(r)\ge 1$ will be determined later.
The product $f_a(t_1)f(t_2): \R^2\rightarrow [0, \infty)$ is a function bounded by $a$ and Lipschitz with constant $2a^{1/2}$ (cf Footnote~\ref{f.Lip-Rd}). Consequently, $2af_a(t_1)f(t_2)$ is a bounded Lipschitz function applicable to~\eqref{e.BL-n-rate}.
This truncation, along with H\"older's inequality used twice and $(2+\epsilon)$- moment for $\xi$, leads to the following bound:
\begin{align*}
\Big | \sE_{\bk x1{2}}[ \xi(\tx_1,\tP_n)\xi(\tx_2,\tP_n)]-\sE_{\bk x1{2}}[ f_a(\xi(\tx_1,\tP_n))f_a(\xi(\tx_2,\tP_n))]\Big| 
&\le 2M_{2+\epsilon/2}^\xi a^{-\epsilon/4}.
\end{align*}
The same approximation by the truncation $f_a$ holds true for $|\xil_{\bk x1p}[1]\xil_{\bk x12}[2]|$ because $\xil_{\bk x12}[1]$ and $\xil_{\bk x12}[2]$ have $(1+\epsilon/4)$-moments with the same bound $M^\xi_{1+\epsilon/4}$. Here, we use  the same argument as in the proof of Item~\ref{i.BL-limit-p-bound}  with $p=1+\epsilon/2$ and $p'=1+\epsilon/4$.
Using the triangle inequality systematically as in \eqref{e.xipmom_trunc}, also using \eqref{e.BL-n-rate} and collecting the developed bounds, we obtain:
$$
\Big | \sE_{\bk x1{2}}[\xi(\tx_1,\tP_n)\xi(\tx_2,\tP_n)]-\sE\bigl[\xil_{\bk x12}[1]\xil_{\bk x12}[2]\bigr]\Bigr|\le 4M^\xi_{2+\epsilon/2}a^{-\epsilon/4}+12 a\varphi_{2}(r),$$
 where $r=d(\bk x12,\partial W_n)$.
Setting $a(r) = \varphi_{2}(r)^{-4/(4+\epsilon)}$, the right-hand side simplifies to:
$$
(4M^\xi_{2+\epsilon/2} + 12)\varphi_{2}(d(\bk x12,\partial W_n))^{\epsilon/(4+\epsilon)}.
$$
This completes the statement of part \ref{i.BL-correlation-rate}.
\end{proof}

\begin{remark}[Distributional limits for \BL~cluster-localizing score functions]\label{r:stab_lemma_BCL}
The statements of Lemmas~\ref{l:BL-limits} and~\ref{l:BL-limits-moments}
remain valid for \(\bk{\xi^{\cup(r)}}{1}{p}\) (in place of
\(\bk{\xi^{(r)}}{1}{p}\)) under the sole assumption of
\emph{\BL~union-localization}~\eqref{Lp-stab-infinite-weak1}, together with the
corresponding moment conditions and/or stationarity assumptions.
The limiting probability kernels induced by \BL~cluster-localization may exist
even in cases where \BL-localization fails; see
Example~\ref{ex:nonlocalizing} below.
In situations where the score function~\(\xi\) satisfies both
\BL-localization and \BL~union-localization, the two types of limits coincide.
The proofs of these results carry over verbatim, starting from the Cauchy
property derived via the triangle inequality in~\eqref{e.BL-r-limit}.
\end{remark}

\begin{exe}[\BL~union-localization without cluster-localization]
 \label{ex:nonlocalizing}
We  give an example showing that it is possible for \eqref{Lp-stab-infinite-weak1} to hold but not \eqref{Lp-stab-infinite-weak2} in the definition of \BL-cluster localization. This shows that marking (score) functions may be neither  \BL~cluster-localizing nor \BL~localizing.

Let $\mu$ be a point process on $\R^d$ and consider $\tmu = \{(x,u_x)\}_{x \in \mu}$, a marked point process in $\R^d \times [0,1]$. Define the marking function 
$$\xi(\tx,\tmu) =
\begin{cases}
\textrm{sgn}(\frac{1}{\mu(\R^d)}\sum_{y \in \mu}u_y - 1/2) \, \quad \textrm{if $\mu$ is finite  }\\
 \textrm{sgn}(u_x - 1/2) \,  \hskip2cm   \quad \textrm{if $\mu$ is infinite,}
\end{cases}
%\label{e.T-KNN-1}
$$
where $\textrm{sgn}(x) := x/|x|$ is the sign function with $\textrm{sgn}(0)$ chosen arbitarily. Let  $\tP$ be a homogeneous Poisson process equipped with  random variables $U_x, x \in \P,$ which are i.i.d. uniform on $[0,1]$. The  $U_x$'s are continuous random variables and hence for all Borel sets $B \subseteq \R^d$, $\xi(\tx,\tP \cap B)$, $x \in \tP \cap B,$ are symmetric Rademacher random variables and  are independent iff $|B| = \infty$.  If $|B| < \infty$, then $\xi(\tx,\tP \cap B) = \xi(\ty,\tP \cap B)$ for all $x, y \in \tP \cap B$.

For all $n \in \N \cup \{\infty\}$, we have for all $r > 0$ that
$$ d_{\BL,\bk{x}{1}{p}}
\bigl(\bk{\xi}{1}{p}(\tx,\tP_n),\bk{\xi^{\cup(r)}}{1}{p}(\tx,\tP_n)\bigr) = 0,$$
that is \eqref{Lp-stab-infinite-weak1} holds,
whereas for all $r < d(\bk{x}{1}{q},\bk{x}{q+1}{p})/2$, there are constants $a_p$ such that 
\begin{align}
& 
d_{\BL,\bk{x}{1}{p}}
\Bigl(\bk{\xi^{\cup(r)}}{1}{p}(\tx,\tP_n),
\bigl(
\bk{\xi^{\cup(r)}}{1}{q}(\tx,\tP_n),
\bk{\xi^{\cup(r)}}{q+1}{p}(\tx,\tP_n)
\bigr)
\Bigr)
\;\ge\; a_p > 0,
\end{align}
and hence $\xi$ fails to satisfy \BL~cluster-localization.
This is because the co-ordinates of $\bk{\xi^{\cup(r)}}{1}{p}(\tx,\tP_n)$ are all $+1$ or $-1$ but $\bk{\xi^{\cup(r)}}{1}{q}(\tx,\tP_n)$ and $\bk{\xi^{\cup(r)}}{q+1}{p}(\tx,\tP_n)$ are independent. It can be verified similarly that $\xi$ satisfies neither \eqref{Lp-stab-infinite} nor \eqref{Lp-stab}, i.e., $\xi$ is 
neither \BL-localizing on $\R^d$ nor on the finite windows of $\R^d$.

Further,  $H_n^{\xi} = \sum_{x \in \P_n} \xi(\tx,\tP_n) = |\P_n|X$ where $X$ is an independent symmetric Rademacher random variable. Thus $\E \, H_n^{\xi} = 0$ and  $\Var \, H_n^{\xi} = \E |\P_n|^2 = n^2 + n$. The strong law of large numbers implies that as $n \to \infty$
$$ \frac{H_n^{\xi}}{\sqrt{\Var \,H_n^{\xi}}} = X\frac{|\P_n|}{\sqrt{n^2+n}}\overset{a.s.}{\to} X,$$
and so $\frac{H_n^{\xi}}{\sqrt{\Var \,H_n^{\xi}}} \Rar X$ as $n \to \infty$. Thus the CLT fails in this case.
\end{exe}

\subsection{Proofs of expectation and variance asymptotics in  Section~\ref{s:limitstatpp}}
\label{s:proofsexp}
\begin{proof}[Proof of Proposition~\ref{expvar}]~
The claims of existence of $\xil_{\0}$ and $\xil_{\0,z}$ in Parts (i) and (ii) respectively follows from Lemma~\ref{l:BL-limits}~\ref{i.BL-r-limit}. We now prove the convergence claims in the two parts. \\

\noindent \ref{i.LLN-expvar}  We establish  the mean  asymptotics \eqref{expasy}. 
Recall that $\tx = (x, U(x))$.
The Campbell-Little-Mecke formula yields
\begin{align*}
&\left| n^{-(d -1)/d}\E \mu_n^\xi(f) -
  n^{1/d}\rho \E[\xil_{\0}]\int_{W_1}f(x)\,\md x\right|   \\
&= n^{-(d -1)/d}\left|\E\Bigl[\sum_{ x \in
    \P_n} f(n^{-1/d} x)\xi( \tx, \tP_n)\Bigr] - n\rho \E[\xil_{\0}]\int_{W_1}f(x)\,\md x \right| \\
& \leq  n^{-(d -1)/d} \rho \int_{W_n} | f(n^{-1/d} x)| \,
    \big| \E_{x} [ \xi( \tx, \tP_n)]- \E[\xil_{x}]\ \big| \, \md x \\
&\le  (2M^{1+\epsilon/2}_\xi+6)
 \|f\|_{\infty}\rho \,n^{-(d -1)/d} \int_{W_n}  \varphi_1(d(x,\partial W_n))^{\epsilon/(2+\epsilon)} \,\md  {x}. 
\end{align*}
Here, we used the distributional equality $\xil_{x}\overset{{\rm{law}}}{=}\xil_{\0}$
(see Lemma~\ref{l:BL-limits}~\ref{i.BL-stat}) and in  the last inequality, we used the estimate~\eqref{e.BL-limit-p-bound} to bound  $\big| \E_{x} [ \xi( \tx, \tP_n)]- \E[\xil_{x}]\ \big|$.
The fast decreasing property of $\varphi_1$
gives 
$$
n^{-(d -1)/d} \int_{W_n} \varphi_1(d(x,W_n))^{\epsilon/(2+\epsilon)} \,\md x = O(1),
$$
and thus \eqref{expasy} follows. If $\varphi_1$ is a function decreasing to $0$, but not necessarily fast decreasing (as in~\eqref{varphibd}), then
$$
n^{-1} \int_{W_n}\varphi_1(d(x,W_n))^{\epsilon/(2+\epsilon)} \,\md  {x}=o(1).
$$

\noindent \ref{i.Var-expvar}
 Now we establish  variance asymptotics \eqref{eqn:var}.
Using the Campbell-Little-Mecke formula we obtain
\begin{align}
\Var \,\mu_n^{\xi}(f) & = \E \sum_{x \in \P_n} f(n^{-1/d} x)^2 \xi^2(\tx, \tP_n) \nonumber \\
& \qquad + \E \sum_{x,y \in \P_n, x \neq y} f(n^{-1/d}x)  f(n^{-1/d}y) \xi(\tx, \tP_n) \xi(\ty, \tP_n) \nonumber \\
& \qquad \qquad   - \Bigl( \E \sum_{x \in \P_n}f(n^{-1/d}x) \xi(\tx, \tP_n) \Bigr)^2 \nonumber	\\
& =  \int_{u \in W_n} f(n^{-1/d}u)^2 \sE_x(\xi^2(\tu, \tP_n ))  \rho  \, \md u \nonumber  \\
\label{eqn:variance_exp-BL} 
& \qquad + \int_{W_n \times W_n} f(n^{-1/d} u) f(n^{-1/d}v) \bigl( m_{(2)}( u,v ;n) - m_{(1)}( u;n)m_{(1)}( v;n) \bigr) \, \md u \md v, 
\end{align}
where 
\begin{align*}
m_{(1)}( u;n) &:= \E_{u}  \xi(\tu, \tP_n ) \, \rho\\
m_{(2)}( u,v;n)& := \E_{u,v} [\xi(\tu, \tP_n) \xi( \tv, \tP_n)] \, \rho^{(2)}(u,v).
\end{align*}

The first term in ~\eqref{eqn:variance_exp-BL}, multiplied by $n^{-1}$, converges
to~
\begin{equation}\label{e:first-term-variance-BL}
    \rho \E[\xil_{\0}^2]  \int_{W_1}f(x)^2\,\md x\;<\infty.
\end{equation}
Indeed, we employ similar arguments as in the proof of expectation asymptotics in point \ref{i.LLN-expvar}, now using the estimate \eqref{e.BL-limit-p-bound} with $q=1$ and $p=2$ for   $\left| \E_{x} [ \xi^2( \tx, \tP_n)]- \E[\xil_{x}^2]\ \right|$  under the $(2+\epsilon)$-moment condition.

For the second term in \eqref{eqn:variance_exp-BL} multiplied by $n^{-1}$, we set $x = n^{-1/d}u$, $z = v-u= v-n^{1/d}x$ and we rewrite it as:
\begin{equation}\label{e:second-term-variance-BL}
  \int_{W_1} \int_{W_n-n^{1/d}x}\hspace{-2em}  f(x)f(x + n^{-1/d}z)
\Bigl[ m_{(2)}(n^{1/d}x, n^{1/d}x+z;\tP_n)-m_{(1)}(n^{1/d}x;\tP_n)m_{(1)}(n^{1/d}x+z;\tP_n)\Bigr]\,\md z\md x.
\end{equation}
We claim now that \eqref{e:second-term-variance-BL} converges to the desired value:
\begin{equation}\label{e.variance-limit-BL}
    \int_{W_1}  f^2(x)\,\md x\times \int_{\R^d}  \left( \sE\bigl[\xil_{
(\0,z)}[1]\xil_{(\0,z)}[2]\bigr] \rho^{(2)}(\0,z) - \rho^2\sE[\xil_{\0}]^2 \right)\,\md z\;<\infty.  
 \end{equation}

To establish this convergence as $n\to\infty$, we use the Lebesgue dominated convergence theorem. Indeed, \(\xi(\tx, \tP_n)\) exhibits fast BL-mixing correlations of order two as in \eqref{e.mppmixing} with $p=q=1$; specifically, it satisfies \eqref{xi-decor-2}. The $\BL$-mixing of correlations is guaranteed by Theorem~\ref{t:mppmix_double}~\ref{i.BL-iterated-marks} and Corollary~\ref{c:mppmixcl_double} under the conditions of fast \BL-localization of \(\xi\) for \(p = \{1, 2\}\) and summable exponential $\B$-mixing correlations of \(\tP\) as in Definition~\ref{def.A2}.  Then the $(2 + \varepsilon)$-moment condition and Lemma~\ref{l:non-bdd-mppmixing} ensure that the integrand  in~\eqref{e:second-term-variance-BL} is bounded uniformly for \(x \in W_1\), \(z \in W_n - n^{1/d}x\), and \(n \in \N \), by a fast-decreasing function \(\tilde\omega_{2}(|z|)\), which is integrable over \(\R^d\). 

Under the  integral, we can exploit the point-wise convergence of the expectations and correlations
developed in Lemma~\ref{l:BL-limits-moments} to obtain
\begin{align}
\label{e:ptwiseexp} \lim_{n\to\infty} \E_{n^{1/d}x} [\xi(\widetilde{n^{1/d}x}, \tP_n )]=\lim_{n\to\infty} \E_{n^{1/d}x+z} [\xi(\widetilde{n^{1/d}x+z}, \tP_n )]&=\E[\xil_{\0}],\\
\label{e:ptwisecorr} \lim_{n\to\infty}\E_{n^{1/d}x,n^{1/d}x+z}[\xi(\widetilde{n^{1/d}x}, \tP_n) \xi( \widetilde{n^{1/d}x+z}, \tP_n)]&=\sE\bigl[\xil_{
(\0,z)}[1]\xil_{(\0,z)}[2]\bigr].
\end{align}
 Specifically,  the stationarity (refer to Lemma~\ref{l:BL-limits}~\ref{i.BL-stat}), for any $x,z\in\R^d$, we observe the following equality in distributions: $\xil_{\0}\overset{{\rm{law}}}{=}\xil_{n^{1/d}x}\overset{{\rm{law}}}{=}\xil_{n^{1/d}x+z}$, which consequently holds  for expectations as well. Additionally, based on estimate~\eqref{e.BL-limit-p-bound}, the differences 
$$\Bigl|\E_{n^{1/d}x} [\xi(\widetilde{n^{1/d}x}, \tP_n )]-\E[\xil_{n^{1/d}x}]\Bigr| \quad \text{and}  \quad \Bigl|\E_{n^{1/d}x+z} [\xi(\widetilde{n^{1/d}x+z}, \tP_n )]-\E[\xil_{n^{1/d}x+z}]\Bigr|$$
are bounded by a decreasing function of the  distance of $n^{1/d}x$ and $n^{1/d}x+z$, respectively, to the boundary of $W_n$. As $n$ approaches infinity, this distance scales as $n^{1/d}$ for fixed $x$ and $z$ and thus we have shown \eqref{e:ptwiseexp}.

Similar arguments can be applied to deduce the convergence of correlations in \eqref{e:ptwisecorr}, utilizing the stationarity for $\xil_{(n^{1/d}x,n^{1/d}x+z)}\overset{{\rm{law}}}{=}\xil_{(\0,z)}$ and the bound \eqref{e.BL-correlation-rate}. Again, this involves a decreasing function of the minimum of the distances $d(n^{1/d}x, \partial W_n)$ and $d(n^{1/d}x+z, \partial W_n)$.

Finally, the continuity of $f$ would suffice for the convergence of $f(x+z/n^{1/d})$ to $f(x)$. In the absence of this assumption, we employ the argument of Lebesgue points of $f$  which is borrowed from the proof of \cite[Theorem~2.1]{Penrose2007gaussian} and was also used in a similar context in \cite{BYY19}. This completes the proof of the variance asymptotics.

 The asymptotic of the covariance  follows from  the polarization identity
\[ {\rm{Cov}}(\mu_n^{\xi}(f),\mu_n^{\xi}(g)) = \frac{1}{2} \left( \Var{\mu_n^{\xi}(f+g)} - \Var{\mu_n^{\xi}(f)}  - \Var{\mu_n^{\xi}(g)} \right).\]
\end{proof}

\begin{remark}[Expectation and variance asymptotics under \BL~cluster-localization]
\label{r:stab_lemma_BCL-variance}
For \(p = 1\), \BL-localization coincides with \BL~cluster-localization—
more precisely, with its first condition~\eqref{Lp-stab-infinite-weak1},
referred to as \BL~union-localization.
Consequently, Part ~\ref{i.LLN-expvar} of Proposition~\ref{expvar}
(the asymptotic behavior of the mean)
remains valid under \BL~cluster-localization.

The variance asymptotic, stated in Part ~\ref{i.Var-expvar} of the same
proposition, also remains valid under \BL~cluster-localization,
provided the localization holds for \(p \in \{1,2\}\).
For \(p = 2\), the Palm limits and the required moment bounds already
exist under the first condition~\eqref{Lp-stab-infinite-weak1},
as noted in Remark~\ref{r:stab_lemma_BCL}.
However, the second condition of \BL~cluster-localization,
namely~\eqref{Lp-stab-infinite-weak2}, is genuinely required.
Indeed, the convergence in~\eqref{e:ptwisecorr} relies on fast
\BL-mixing correlations of order two.
Such mixing is ensured only when the two conditions defining
\BL~cluster-localization are satisfied jointly, as formalized in
Theorem~\ref{t:mppmix_double}\ref{i.BLc-iterated-marks}.
Example~\ref{ex:nonlocalizing} in Section~\ref{s:stab_lemma_BL}
provides a counterexample illustrating the necessity of the second
condition.
\end{remark}

\subsection{Proof of multivariate central limit theorem in  Section~\ref{s:limitstatpp}%--- Theorem \ref{t:multcltmarkedpp}
}
\label{s:proof_multclt}

\begin{proof}[Proof of Theorem \ref{t:multcltmarkedpp}]

The asymptotics of the mean and variance of the random variables $(\mu_n^{\xi^i}(f))_{n \in \N}$,  $i=1,\ldots,{m}$ follow from Proposition~\ref{expvar}. Since $(x,y) \mapsto f(x+y)$ is a bounded Lipschitz function in the $1$-product metric when $f \in \BL(\R)$, the sum $\xi_i + \xi_j$ forms a new $\BL-$localizing score function. Further, $\BL$-localization of order $p =2$ suffices for us. Now the limiting covariance can be derived from the variance-covariance identity, 
\begin{equation} \label{cova}
  \mathrm{Cov}(\mu_n^{\xi_i}(f), \mu_n^{\xi_j}(f)) = \frac{1}{2} \left( \mathrm{Var}(\mu_n^{\xi_i + \xi_j}(f)) - \mathrm{Var}(\mu_n^{\xi_j}(f)) - \mathrm{Var}(\mu_n^{\xi_i}(f)) \right), \quad 1 \leq i \leq j \leq k,
\end{equation}
together with variance asymptotics \eqref{eqn:var} for $\xi_i + \xi_j, \xi_i, \xi_j$ and the assumed $(2 + \varepsilon)$-moment condition.

The proof of the multivariate central limit theorem for the vector \((\mu_n^{\xi_1}(f), \ldots, \mu_n^{\xi_{m}}(f))\) uses the Cramér-Wold device and involves relaxing the moment conditions in Theorem~\ref{t:cltmarkedpp} to some moment \(p>2\). For this, we employ an approach suggested by Matthias Schulte, which focuses on reducing the problem of asymptotic normality to establishing three variance asymptotics, specifically \eqref{e.var-xi}, \eqref{V1} and \eqref{V2} below. 
\smallskip 

{\em Cramér-Wold linear combination:}
Consider a linear combination \(\xi := \sum_{i=1}^{m} t_i \xi_i\) of the score functions \(\xi_i\) with given fixed coefficients \(t_i \in \mathbb{R}\), \(i=1,\ldots,{m}\). We aim to prove the central limit theorem for 
\(\mu_n^{\xi}(f) = \sum_{i=1}^{m} t_i \mu_n^{\xi_i}(f)\) for some given function \(f \in \B(W_1)\).
Observe, by the joint fast \BL-localization of \(\bxi = (\xi_1, \ldots, \xi_{m}) \in \mathbb{R}^{m}\) assumed in Theorem \ref{t:multcltmarkedpp}, the linear combination  
\(\xi\) also satisfies fast \BL-localization on finite windows on \(\tP\) as in \eqref{Lp-stab} with 
fast decreasing functions \(\varphi_p\) (depending on  coefficients \(t_i\), \(i=1,\ldots,{m}\)) for 
all \(p \in \N\). Consequently, by the assumption of \(\tP\) having summable exponential $\B$-mixing correlations, Theorem~\ref{t:mppmix_double}~\ref{i.BL-iterated-marks}, and Corollary~\ref{c:mppmixcl_double}, \(\xi\) has fast \BL-mixing correlations as in \eqref{e.mppmixing} on this input process (with \(U\) replaced by \(\xi\)), with some fast decreasing functions \(\hat\omega_k\) (again involving \(t_i\), \(i=1,\ldots,{m}\)). Clearly, \(\xi\) inherits the stationarity of \(\bxi\) and the \(p > 2\) moment condition for \(\xi_i\), \(i=1,\ldots, {m}\) (however, the values of the constants \(M^\xi_p\) involve the linear coefficients \(t_i\), \(i=1,\ldots,{m}\)).
\smallskip

{\em BL-truncation of $\xi$:} To bypass the requirement of all moment conditions for the central limit theorem in  Theorem~\ref{t:cltmarkedpp}, for a given \(M > 0\) we define a function \(g_M : \mathbb{R} \to \mathbb{R}\) by 
\( g_M(x) := x\mathbf{1}_{|x| \leq M} + \text{sgn}(x)M\mathbf{1}_{|x| > M} \)
and consider two auxiliary score functions which are transformations of  \(\xi\):
\begin{align*}
\xi_M(x, \tP_n) &:= g_M(\xi(x, \tP_n)), \\
\xi_{\Delta M}(x, \tP_n) &:= \xi(x, \tP_n) - \xi_M(x, \tP_n).
\end{align*}
Observe the equivalent representation
\[\xi_{\Delta M}(x, \tP_n) = g_{\Delta M}(\xi(x, \tP_n)),\]
with \( g_{\Delta M}(x) :=(x - \text{sgn}(x)M)  \mathbf{1}_{|x| > M}\).
Since \(g_M\) and \(g_{\Delta M}\) are Lipschitz(1), we have \(f \circ \bigotimes_{i=1}^p g_M, f \circ \bigotimes_{i=1}^p g_{\Delta M} \in \BL(\mathbb{R}^p)\) for any \(f \in \BL(\mathbb{R}^p)\). Therefore, \(\xi_M\) and \(\xi_{\Delta M}\) satisfy the same fast \BL-localization properties \eqref{Lp-stab} as does \(\xi\), with the same fast-decreasing function \(\varphi_p\), for all \(p \in \N\), uniformly in \(M > 0\). Similarly, \(\xi_M\) and \(\xi_{\Delta M}\) have fast \BL-mixing correlations as in \eqref{e.mppmixing} on this input process (with \(U\) replaced by \(\xi_M\) and \(\xi_{\Delta M}\), respectively), with the same fast decreasing functions \(\hat{\omega}_k\) as for \(\xi\), uniformly in \(M\).
This is again a consequence of the fact that the functions \(f \circ \bigotimes_{i=1}^p g_M\), \(g \circ \bigotimes_{i=1}^q g_M\), \(f \circ \bigotimes_{i=1}^p g_{\Delta M}\), and \(g \circ \bigotimes_{i=1}^q g_{\Delta M}\) are Lipschitz(1), when \(f \in \BL(\mathbb{R}^p)\) and \(g \in \BL(\mathbb{R}^q)\). Clearly, \(\xi_M\) and \(\xi_{\Delta M}\) are stationary just like $\xi$.
Furthermore, \(\xi_M\) satisfies \(p\)-moment conditions~\eqref{e:xinpmom} for all \(p > 1\), with constants \(M^{\xi_M}_p \le \max(1, M^p)\), while \(\xi_{\Delta M}\) has the same \(p > 2\) moment condition existing for \(\xi\), with a bound \(M^{\xi_{\Delta M}}_p \le M^\xi_p\) because \(|g_{\Delta M}(x)| \le |x|\).
\smallskip

{\em  Variance asymptotics of $\xi$, $\xi_M$, and  $\xi_{\Delta M}$:}
By the above observations, the three score functions \(\xi\), \(\xi_M\), and \(\xi_{\Delta M}\) satisfy the assumptions for the variance asymptotics formulated in Proposition~\ref{expvar}~\ref{i.Var-expvar}. We denote these limits by:
\begin{align}
\lim_{n \to \infty} n^{-1} \Var \,\mu_n^{\xi}(f) &=: \sigma^2, \label{e.var-xi}\\
\lim_{n \to \infty} n^{-1} \Var \,\mu_n^{\xi_M}(f) &=: \sigma^2_M, \label{e.var-xi_M}\\
\lim_{n \to \infty} n^{-1} \Var \,\mu_n^{\xi_{\Delta M}}(f) &=: v^2_M. \label{e.var-xi-DM}
\end{align}
Moreover, by Theorem~\ref{t:cltmarkedpp}, \(\mu_n^{\xi_M}(f)\) satisfies the (possibly degenerate) CLT as \(n \to \infty\). Thus, to establish \eqref{e.M-CLT} for \(\mu_n^{\xi}(f)\), we assert that it is enough to show the following limits 
\begin{align}
\label{V1} 
\limsup_{M\to\infty} \;& v^2_M\to 0, \\
\label{V2} 
\lim_{M\to\infty} & \sigma^2_M\to \sigma^2.
\end{align}
\smallskip

{\em Central limit theorem for  $(\mu_n^{\xi}(f))_{n \in \N}$  via three variance asymptotics:}
Before establishing~\eqref{V1} and~\eqref{V2}, we demonstrate how these limits and the asymptotic normality of 
$(\mu_n^{\xi_M}(f))_{n \in \N}$ implies that of $(\mu_n^{\xi}(f))_{n \in \N}$. Denote
$$
\overline{\mu}_n^\xi(f) := n^{-1/2}(\mu_n^{\xi}(f) - \mathbb{E}[\mu_n^{\xi}(f)]),
$$
and similarly \(\overline{\mu}_n^{\xi_M}(f)\) for \(\xi_M\).
The triangle inequality gives, for all \(m \in \N\) and $\sigma^2\ge 0$,
\begin{align}\label{e.3variances1} 
& d_{\BL}(\overline{\mu}_n^\xi(f), N(0,\sigma^2))\\
 & \leq d_{\BL}(\overline{\mu}_n^\xi(f), \overline{\mu}_n^{\xi_M}(f)) +
 d_{\BL}(\overline{\mu}_n^{\xi_M}(f), N(0,\sigma^2_M)) +
 d_{\BL}(N(0,\sigma^2_M), N(0,\sigma^2)).\nonumber
\end{align}
In view of
\begin{align*}
& d_{\BL}(\overline{\mu}_n^\xi(f), \overline{\mu}_n^{\xi_M}(f)) \nonumber \\ & \leq 
\mathbb{E}|\overline{\mu}_n^\xi(f) - \overline{\mu}_n^{\xi_M}(f)|
\leq \sqrt{n^{-1} \Var (\mu_n^{\xi}(f) - \mu_n^{\xi_M}(f))} = \sqrt{n^{-1} \Var \,\mu_n^{\xi_{\Delta M}}(f)},
\end{align*}
we can rewrite~\eqref{e.3variances1} as
\begin{align} \label{triangleDBL}
& \quad d_{\BL}(\overline{\mu}_n^\xi(f), N(0,\sigma^2)) \\
& \leq \sqrt{n^{-1}\Var \,\mu_n^{\xi_{\Delta M}}(f)} +
d_{\BL}(\overline{\mu}_n^{\xi_M}(f), N(0,\sigma^2_M)) +
d_{\BL}(N(0,\sigma^2_M), N(0,\sigma^2)).\nonumber
\end{align}
With \(M\) fixed, we take \(\limsup_{n \to \infty}\) of both sides of \eqref{triangleDBL}.
On the right-hand side: the first term goes to \(v_M\) by~\eqref{e.var-xi-DM}, while the second term goes to 0 either by  $\sigma^2_M=0$ or the asymptotic normality of \((\overline{\mu}_n^{\xi_M}(f))_{n \in \N}\). 
Consequently,
$$ 
\limsup_{n \to \infty} d_{\BL}(\overline{\mu}_n^\xi(f), N(0,\sigma^2)) \leq v_M +
d_{\BL}(N(0,\sigma^2_M), N(0,\sigma^2)).
$$
Now letting $M \to \infty$ on both sides and using \eqref{V1} and~\eqref{V2} establishes  \eqref{e.M-CLT} for \((\mu_n^{\xi}(f))_{n \in \N}\). It is therefore enough to show~\eqref{V1} and~\eqref{V2}, as asserted.
\smallskip 

\noindent{\em Proof of variance limits \eqref{V1} and~\eqref{V2}:}
By~\eqref{eqn:var} and~\eqref{sigdef}, the asymptotic variance of \(\xi\) can be represented as
\begin{align}\label{eqn:var-M}
\sigma^2 &= \sigma^2(\xil) \int_{W_1} f(x)^2 \,\md x, \\
\noalign{with}
\sigma^2(\xil) &= \rho \mathbb{E}[(\xil_{\0})^2] + \int_{\R^d} \left( \mathbb{E}\bigl[\xil_{(\0,z)}[1]\xil_{(\0,z)}[2]\bigr] \rho^{(2)}(\0,z) - \rho^2 \mathbb{E}[\xil_{\0}]^2 \right) \md z, \label{sigdef-M}
\end{align}
where \(\xil_{\0}\) and \(\xil_{(\0,z)} = (\xil_{(\0,z)}[1], \xil_{(\0,z)}[2])\) represent the corresponding Palm distributional limits of \(\xi\) defined in~\eqref{e.xil-0} and~\eqref{e.xil-xy}.

Analogously, \(\sigma^2_M\) and \(v^2_M\) can be represented by~\eqref{eqn:var-M}-\eqref{sigdef-M} with \(\xil_{\0}\) and \(\xil_{(\0,z)}\) replaced, respectively, by \(\xil^{M}_{\0}\), \(\xil^{M}_{(\0,z)}\) and \(\xil^{\Delta M}_{\0}\), \(\xil^{\Delta M}_{(\0,z)}\), representing the corresponding Palm distributional limits of \(\xi_M\) and \(\xi_{\Delta M}\).
Observe in~\eqref{e.xil-0} and~\eqref{e.xil-xy}, these limits are for the auxiliary score functions \(\xi_M = g_M(\xi)\) and \(\xi_{\Delta M} = g_{\Delta M}(\xi)\), which are produced by composing \(\xi\) with Lipschitz(1) functions \(g_M\)
and \(g_{\Delta M}\).
Using the continuity of Lipschitz functions with respect to \BL-convergence,
it follows that with $M$ fixed, the Palm distributional \BL-limits as $n \to \infty$ of
$\xi_M$ and $\xi_{\Delta M}$
can be computed by applying the corresponding Lipschitz(1) functions to the limit of $\xi$:  
\begin{align*}
\xi_M(\tilde\0, \tP_n) &\xrightarrow[n \to \infty]{\text{BL}, \0} \xil^{M}_{\0} \stackrel{\rm{law}}{=} g_M(\xil_{\0}), \\
\Bigl(\xi_M(\tilde\0, \tP_n), \xi_M(\tz, \tP_n)\Bigr) &\xrightarrow[n \to \infty]{\text{BL}, (\0,z)} \xil^{M}_{(\0,z)} \stackrel{\rm{law}}{=} \Bigl(g_M(\xil_{(\0, z)}[1]), g_M(\xil_{(\0, z)}[2])\Bigr),
\end{align*}
and similarly,
\begin{align*}
\xil^{\Delta M}_{\0} &\stackrel{\rm{law}}{=} g_{\Delta M}(\xil_{\0}), \\
\xil^{\Delta M}_{(\0,z)} &\stackrel{\rm{law}}{=} \Bigl(g_{\Delta M}(\xil_{(\0, z)}[1]), g_{\Delta M}(\xil_{(\0, z)}[2])\Bigr).
\end{align*}
Now, by the dominated convergence theorem and the  \(p>2\) moment condition (see Lemma~\ref{l:BL-limits-moments}),
\begin{align*}
\lim_{M \to \infty} \mathbb{E}[\xil^{\Delta M}_{\0}] &= \mathbb{E}[\lim_{M \to \infty} g_{\Delta M}(\xil_{\0})] = 0, \\
\lim_{M \to \infty} \mathbb{E}[\xil^{\Delta M}_{(\0,z)}[i]] &= \mathbb{E}\bigl[\lim_{M \to \infty} g_{\Delta M}(\xil_{(\0, z)}[i])\bigr] = 0, \quad i=1,2, \\
\noalign{and also}
\lim_{M \to \infty} \mathbb{E}\bigl[\xil^{\Delta M}_{(\0,z)}[1] \xil^{\Delta M}_{(\0,z)}[2]\bigr] &= 0.
\end{align*}

Analogously, we have for \(\xi_{M}\):
\begin{align*}
\lim_{M \to \infty} \mathbb{E}[\xil^{M}_{\0}] &= \mathbb{E}[\lim_{M \to \infty} g_{M}(\xil_{\0})] = \mathbb{E}[\xil_{\0}], \\
\lim_{M \to \infty} \mathbb{E}[\xil^{M}_{(\0,z)}[i]] &= \mathbb{E}[\lim_{M \to \infty} g_{M}(\xil_{(\0, z)}[i])] = \mathbb{E}\bigl[\xil_{(\0, z)}[i]\bigr], \quad i=1,2, \\
\lim_{M \to \infty} \mathbb{E}\bigl[\xil^{M}_{(\0,z)}[1] \xil^{M}_{(\0,z)}[2]\bigr] &= \mathbb{E}\bigl[\xil_{(\0,z)}[1] \xil_{(\0,z)}[2]\bigr].
\end{align*}

To establish limits~\eqref{V1} and~\eqref{V2}, it remains  to justify the passage of \(\lim_{M \to \infty}\) under the integral \(\int_{\R^d} [\ldots] \, \md z\) in~\eqref{sigdef-M}. This is justified by the dominated convergence theorem, evoking, for both \(\xi_M\) and \(\xi_{\Delta M}\), \BL-mixing correlations uniformly in \(M>0\); that is, these properties involve the common functions \(\omega_k\) in~\eqref{e.mppmixing}. Indeed, by Lemma~\ref{l:non-bdd-mppmixing} with \(p=q=1\), \(K=2\),
we have
$$
\Bigl|\mathbb{E}\bigl[\xil^{M}_{(\0,z)}[1] \xil^{M}_{(\0,z)}[2]\bigr] \rho^{(2)}(\0,z) - \rho^2 \mathbb{E}[\xil^{M}_{\0}]^2 \Bigr| \le \hat\omega_{2}(|z|) 
$$
and similarly for $\xil^{\Delta M}$, 
with a fast-decreasing function \(\hat\omega_{2}\) that depends on \(\omega_{2}\) (and other constants) uniformly in \(M\). 

This completes the proof of limits~\eqref{V1} and~\eqref{V2}, and thus completes the proof  of Theorem~\ref{t:multcltmarkedpp}.
$ $\par\nobreak\ignorespaces %This forces the Halmos symbol to appear at the end of the line (when there is no space left in the previous one).
\end{proof}

\begin{remark}[Bypassing the summable exponential mixing condition in the multivariate CLT]
\label{r.M-CLT-bypassing-SExp}
In Theorem~\ref{t:multcltmarkedpp}, which establishes a multivariate CLT for statistics of point processes with fast \BL-localizing marks under a $(2+\epsilon)$-moment condition, 
the assumption of summable exponential mixing correlations of the input process~$\tP$ 
can be relaxed to merely fast mixing of~$\P$, 
provided that one can directly verify the fast \BL-mixing correlations of the (iterated) marks 
$\bxi(\tx, \tP_n) = (\xi_1, \ldots, \xi_m) \in \R^m$ associated with~$\P$. 
Indeed, the stronger assumption on~$\tP$ is only required in the application of 
Theorem~\ref{t:mppmix_double} 
to ensure the corresponding mixing property of these new marks. 
All the crucial steps in the proof of Theorem~\ref{t:multcltmarkedpp}---namely, the use of  the Cramér–Wold linear combination $\xi$ and the justification of the CLT for its 
\BL-truncation $\xi_M$---rely solely on the fast \BL-mixing of~$\bxi$ (of all orders), 
while the derivation of the variance asymptotics depends only on the \BL-mixing of~$\bxi$ 
of order~2. \end{remark}

\subsection{Proof of limit theory for scores on infinite window in Section~\ref{s:ltthms_stat_marking}}
 \label{proofProp5.5} 

\begin{proof}[Proof of Proposition~\ref{t:clt_linear_marks_new}]
\noindent (i) Taking  $\xi_{i,n}=\xi_{i}:=\xi( \tx_i, \tP)$,  
the central limit theorem for $(\hat\mu_{n}^{\xi}(f))_{n \in \N}$ is an easier version of the  triangular array central limit theorem given by Theorem ~\ref{t:clt_linear_marks},
where the summands of $\hat\mu_n^\xi$, namely the terms $\xi_{i,n}:=\xi( \tx_i, \tP)$, are equal for  $n \in \N$. 
Indeed, the arguments in the proof of Theorem~\ref{t:cltmarkedpp} apply  straightforwardly via Theorem~\ref{t:mppmix_double}.
The assumptions on $\tP$ and $\xi$ imply that the family of marked point process $\dbtilde \P_n:=  \{(x,\xi(\tx,\tP))\}_{x \in \P_n}$  has fast \BL-mixing  correlations as in Definition 
\ref{d.omegahmixing}(iv).
%\ref{i.family-mixing}.
Hence   Theorem~\ref{t:clt_linear_marks}
yields   the central limit theorem for $(\hat\mu_{n}^{\xi}(f))_{n \in \N}$.  
  
\medskip 

\noindent (ii) For the limit of the expectation, observe by Campbell-Little-Mecke's formula and substituting \(x\) for \(n^{-1/d} x\) that
\[  \E\left[\sum_{x \in \P_n} f(n^{-1/d} x)\xi(\tx, \tP)\right] = n\rho \E_0[\xi(\tilde {\0}, \tP)]\int_{W_1} f(x) \, dx.
\]
Hence, the right-hand side of \eqref{expasy} is equal to zero when \(\E_0[\xi(\tilde{\0}, \tilde{\P})]\) is substituted for \(\E[\xil_{\0}]\).

For the variance, observe in the proof of Proposition ~\ref{expvar} in Section~\ref{s:proofsexp}
 that for all $n \in \N$
 $$
  m_{(1)}(u;n) = \E_0[\xi(\tilde{\0}, \tilde{\P})] \, \rho \text{ and } m_{(2)}(u,v;n) = \E_{\0,z} [\xi(\tilde{\0}, \tP) \xi(\tz, \tP)] \, \rho^{(2)}(0,z) \text{ with } z = v - u. 
 $$
Consequently, the expression in the bracket \([ \ldots ]\) in \eqref{e:second-term-variance-BL} is equal to 
 \[\E_{\0,z}[ \xi(\tilde{\0}, \tP) \xi(\tz, \tP) \rho^{(2)}(0,z)\allowbreak - \E_0[\xi(\tilde{\0}, \tilde{\P})]^2 \rho^2], \] 
 which is integrable with respect to \(z\) on \(\R^d\) thanks to the fast mixing correlations of \(\xi(\tx, \tP)\) guaranteed by Theorem~\ref{t:mppmix_double}~\ref{i.BL-iterated-marks}. 
Then, the limit of \(f(x + n^{-1/d}z)\) to \(f(x)\) under the integral relies on the same arguments  employed in  proof of Propositions~\ref{expvar}.
\end{proof}

\part{Applications}
\label{part:applications}

In this part, the theoretical framework developed in the previous parts is applied to various spatial random models. The focus is on spin systems (Section \ref{s:gibbsmarking}), interacting diffusions (Section \ref{NEW-s:id_sprg}), and particle systems on spatial random graphs (Section \ref{s:applnsips}).  Each application illustrates the practical relevance of the theory, offering detailed examples and proofs of how these systems can be studied within the framework established in Part \ref{part:theory-foundations}. This part also explores applications to empirical random fields and geostatistical models (Section \ref{s:erfgeostat}) as well as indicates applications to other spatial random models (Section \ref{s:future}) highlighting the versatility of the methods introduced.

Carrying out these various applications requires establishing localization: either directly, as when showing \BL~cluster-localization for spin systems, or \BL-localization for certain empirical random fields models or indirectly, as when showing \BL-localization for statistics of other models by leveraging more classical stabilization techniques (such as $L^2$-stabilization for interacting diffusions and stopping set stabilization for interacting particle systems).
Moreover, models originally defined on a fixed graph, when considered in a random geometric setting, typically also require the use of stabilization techniques for the associated graph structures.

For convenience, we recall our standing assumption from Part \ref{part:theory-foundations}. 
\begin{customass}{3.1} \,
\begin{enumerate}[wide,label=(\roman*),labelindent=0pt]

\item $\P$ denotes a simple point process on $\R^d$ such that whenever for some $p \in \N$ the correlation function~$\rho^{(p)}$ exists, it is uniformly bounded i.e., $\kappa_p:=\sup_{\bk x1p\in(\R^d)^p}  \rho^{(p)}(\bk x1p) < \infty$.  
We put $\kappa_0 :=  \max\{ \kappa_1, 1 \}$.

\item  All  marked point processes in this article are assumed to be simple.
\end{enumerate}
\end{customass}
Our applications hold for many point processes and spatial random graphs with explicit assumptions made precise in the respective results. However, for simplicity and to get a quick understanding of the applications, the reader may take  $\P$ to be  a stationary Poisson process or a determinantal process with fast-decaying kernel, such as the Ginibre process, and one may take the spatial random graph to be the Gilbert graph, $k$-nearest neighbor graph or Delaunay graph. Some of these applications---most notably spin systems involving
\BL~cluster-localization---are new even in the context of Poisson processes.
We also alert the reader that one may consider a more extended framework for random graphs than the one investigated here. In particular, this could presumably include a random connection model or a Boolean model with random grains in our framework. We remark on this in Section~\ref{s:future}.

\section{Spin systems on spatial random graphs}
\label{s:gibbsmarking}

Spin systems (or Gibbs random fields) defined on spatial random graphs fall within the scope of our main general result, namely Theorem~\ref{t:cltmarkedpp}. Specifically, we establish the asymptotic normality of the total sum of spins. Our Gibbsian models differ from those of Gibbs point processes (e.g., see \citet{dereudre2019introduction}), where probability measures are defined on locally finite configurations of points, $\cN_{\R^d}$. In contrast, we consider a spatial random \emph{interaction graph} on a point process. Given its realization, we define a spin system (or Gibbs random field) on the points based on the adjacency relations within the graph.  
This approach is inspired by lattice spin systems, except that here, the underlying graph is random. Our models serve as spatial analogs to statistical mechanics models on (non-spatial) random graphs, such as those studied by \citet{georgii2001random}, \citet{van2017stochastic}, \citet{duminil2015geometric}, and \citet{friedli2017statistical}.

In Section~\ref{s:modelgibbs}, we introduce the spin model on spatial random graphs and two key assumptions: \emph{stabilization} of the interaction graph and \emph{averaged weak spatial mixing} of spins on this graph.  
Stabilization refers to the property that the construction of a random graph, based on an input process of points (interpreted as graph vertices), is governed by the stopping sets of these points. In this context, long-range interactions in the graph structure statistically decay at a fast rate. This property, which also ensures sparsity, is a common characteristic of many proximity graphs.  
The second assumption, averaged weak spatial mixing, requires that the dependence of spins on the boundary condition in a spatial random graph (accounting for both the Gibbs measure and the randomness of the graph) decays rapidly with respect to the Euclidean distance to the boundary, as the latter tends to infinity. This is an adaptation of the classical weak spatial mixing property for deterministic graphs and typically holds in the high-temperature or low-activity regime of Gibbs models.

In Section~\ref{mainspinresults}, we combine these two key assumptions 
to establish a central limit theorem for the total sum of spins. Our central limit theorem accounts for the randomness of  both the spin variables and underlying spatial random graph. While the spin variables themselves do not exhibit a stopping set property, the combination of averaged weak spatial mixing and graph stabilization guarantees the  \BL~cluster-localization of the spins.  Though establishing sufficient conditions for averaged weak spatial mixing  is a distinct challenge, there have been notable contributions of \citet{dobrushin1987completely}, \citet{weitz2005combinatorial}, \citet{sinclair2017spatial}, and \citet{regts2023absence}, among others. In Sections~\ref{s:Ex_spin1} and \ref{s:hard-core}, we employ two complementary approaches   providing examples of spin models satisfying the  averaged weak spatial mixing condition. 

The first is the disagreement percolation method introduced by \citet{van1994disagreement}, which is applicable to a wide range of spin models and is detailed in Section~\ref{s:gen_disagreement}. The required disagreement percolation bounds may be obtained using path-counting arguments involving the  (averaged) `connective constant' as described in Section~\ref{s:gen_disagreement} whereas for random graphs exhibiting sharp phase transitions, we may use more direct arguments from continuum percolation theory; see Section~\ref{s:spin_Poisson}. The use of sharp phase transition offers two key advantages. First, it applies to graphs with unbounded degrees and second, since the sharp phase transition assumption is, in principle, weaker than a bound on the (quenched and even  averaged) connective constant, it has the potential to yield central limit theorems for a wider range of temperatures or activity parameters. We apply this to the Poisson Gilbert graph in Section \ref{s:spinPoissonGilbert}.

The second approach to showing averaged weak spatial mixing utilizes a combinatorial technique developed by \citet{weitz2005combinatorial} and \citet{sinclair2017spatial}, specifically designed to show a stonger (and quenched) spatial mixing for Gibbs  models and is illustrated with an application to the hard-core model in Section~\ref{s:hard-core} . For these analyses, the range of the temperature or activity parameter is linked to the (quenched) connective constant, a concept adapted from \citet{sinclair2017spatial} to measure how the degree distribution of a graph influences spatial mixing. 

We refer the reader to Section \ref{ss.Literature} for more on the existing literature of limit theory for spin systems, which is mainly restricted to lattices, save for \cite{GGHP15,giardina2016annealed,can2019annealed}, which treats locally-tree like random graphs with local weak limits.

\subsection{The spin model}
\label{s:modelgibbs}
Let $\K$ be a finite or countably infinite subset of $\R$, which we call the spin space. Given a point set of sites (nodes) denoted by $\X : = \{x_1,\ldots,x_m\}$ and  a graph  $\G(\X)$,   we associate a $\K$-valued random variable to each
site, called the {\em spin variable}. The spin variables $\bk v1m:=(v_1,\ldots,v_m)$ of points $\bk x1m:=(x_1,\ldots,x_m)$ are chosen randomly, with joint probability proportional to
$$
\pi_{\bk x1m}(\bk v1m):=\exp\left[ \beta \sum\limits_{(x_i,x_j) \in \G(\X) } \Psi(v_i,v_j) +\, \gamma\sum\limits_{x_i \in \X} \Phi(v_i)\right],$$
where $\beta \in (0, \infty),  \gamma \in \mR$, $\Psi: \K \times \K \to [-\infty,\infty)$ is  the symmetric {\em pairwise potential} 
and  where $\Phi: \K \to \R$ is the  {\em external field}.
More formally, given $\G(\X)$, we consider 
the collection of $\K$-valued random variables $\{V_x\}_{x \in \X}:= \{V_x(\G(\X))\}_{x \in \X}$ distributed according to 
the following probability measure $\mP_{\G(\X)}$ on the configuration space $\K^{\X}$
 \be
\mP_{\G(\X)}(V_{x_1} = v_1,\ldots, V_{x_m} = v_m) =
\frac{\pi_{\bk x1m}(\bk v1m)}{Z(\X)} 
\label{eqn:gibbs}
\ee
 where $Z(\X) := Z(\G(\X))$ is the {\em normalizing constant} 
$$ Z(\X) := \sum_{\bk v1m:=(v_1,\ldots,v_m) \in \K^m}  \pi_{\bk x1m}(\bk v1m). $$
The probability measure $\mP_{\G(\X)}$ at ~\eqref{eqn:gibbs} is the {\em Gibbs measure} (or the {\em specification of the spin model}) on $\G(\X)$,  and the collection of $\K$-valued random variables $ \{V_x\}_{x \in \X}$ distributed as $\mP_{\G(\X)}$ is  the {\em spin configuration}. The probability measure $\mP_{\G(\X)}$ exists whenever $Z(\X) > 0$.
 In the language of statistical physics, $\beta$ represents the inverse temperature. We refer to  \citet{friedli2017statistical} and \citet{duminil2015geometric} for more details on such models. 
 Suitable choices of $\K,\Psi$ and $\Phi$ yield
  the hard-core and soft-core models, the Widom-Rowlinson model, the random cluster model, the Ising model and the  Potts model (see \cite[Section 1.2]{duminil2015geometric}). Examples are in Sections \ref{s:Ex_spin1} and \ref{s:hard-core}.

Spatial mixing controls the correlations between spin values. Additionally, we need to control the interaction graph, as it impacts the spatial mixing bounds. This is achieved by  requiring that the interaction graph has a stabilizing property, specifically by assuming that the radii of stabilization for vertex neighborhoods are not too large. We now formally introduce both the  notions of stabilizing interaction graphs and spatial mixing on such graphs, and then,  in Section~\ref{mainspinresults}, we  state our main central limit theorem for the random measures  $\mu_n^{\G} :=  \sum_{x \in \P_n} V_x(\G(\P_n)) \delta_{n^{-1/d}x}$ induced by the spin systems.

\subsubsection{Stabilizing  interaction graphs}\label{sss.graph-stabilization}
In this section we explain how we construct our random graphs  on the input (ground) process $\P$ within finite windows.   We refer to Section \ref{s:notation} for all notational formalism regarding marked point processes.  In particular, recall that $\mP_{\bk x1p}$ denotes the Palm probability distribution of the point process $\P$ conditional on the existence of  points $x_1,...,x_p$. 

For any finite point set \(\mu \in \hat{\cN}_{W_n}\), with \(n \in \N\), let  \(\G(\mu)\) be a graph constructed on \(\mu\), where the edges are determined by a deterministic rule \(\sim\) 
that takes the entire point set \(\mu\) into account.
(Thus $\G(\mu)$ encodes the interactions between sites in the spin system~\eqref{eqn:gibbs}, putting an edge between sites if they may potentially interact.)
A deterministic construction of the graph $\G(\mu)$ on $\mu$ 
 is a modelling restriction. It excludes, for example, graphs built on a Boolean model with random grains or those from a random connection model. In Section~\ref{s:future}, we discuss how this framework can be extended to include such graphs. 
 
In order to introduce the notion of stabilization for these graphs (as well as to address the measurability issues of these objects), it is customary to represent the neighbors of a point $x \in \mu$ within $\G(\mu)$ as a subset $N_x(\mu) \subset \mu$, where $N_x(\mu) := \{y \in \mu : y \stackrel{\G(\mu)}{\sim} x\}$. For the sake of translation invariance (if assumed), we center these finite neighborhood-sets  as finite counting measures $\bar N_x(\mu) := \sum_{y \in N_x(\mu)} \delta_{y - x} \in \hat\cN_{\R^d}$ and treat them as auxiliary marks of the ground points $x \in \mu$. In essence, 
when assuming the existence of a graph $\G(\mu)$ together with its neighborhood structure, we view this as an auxiliary marking $\{\delta_{(x, \bar N_x(\mu))}\}_{x \in \mu}$ of the ground measure $\mu$, where the marks satisfy $\bar N_x(\mu) + x \subset \mu$ for all $x \in \mu$ (neighborhood relation), and $y \in \bar N_x(\mu)$ if and only if $-y \in \bar N_{x + y}(\mu)$ (undirected).

A consequence of assuming a deterministic rule  for constructing the graph $\G(\mu)$ on $\mu$ is that we may view $\bar N_\cdot(\cdot) : W_n \times \hat{\cN}_{\R^d} \longrightarrow \hat{\cN}_{\R^d}$ as a marking function in the sense of Definition~\ref{d:mark_fn} (though here the ground process has no marks).  More formally, for us an undirected graph corresponds to a {\em neighborhood marking function} $\bar N$ as above satisfying the two conditions of neighborhood relation and undirectedness and thus measurability of the map $\mu \mapsto \G(\mu)$ is equivalent to measurability of the marking function $\bar N$.

Consequently, we define the {\em interaction range} $S_n : W_n \times \hat{\cN}_{W_n} \longrightarrow \N$ of $x\in\mu\in\hat\cN_{W_n}$ in the graph $\G(\mu)$ as the radius of stabilization of this neighborhood marking $\bar N$, as per Definition~\ref{def.Stab.Radius}, with
\begin{equation}
S_n(x,\mu):=R^{\bar N}_{W_n}(x;\mu).\label{e:S_n=R_Wn}
\end{equation}
By this definition, all graphs $\G(\mu)$ on $\mu \in \hat{\cN}_{W_n}$  have bounded  interaction range  \( S_n \le 
\lceil \text{\rm diam} (W_n)+1\rceil\). However, we are specifically interested in graphs whose interaction ranges  \( S_n(x, \P_n) \) decay rapidly and uniformly in $x \in W_n, n \in \N$. This concept, which by the way, implies asymptotic sparsity of  the  family of graphs $(\G(\P_n))_{n \in \N}$, is formalized as follows (in keeping with the general Definition~\ref{def.stabilizing_marking}\ref{i.stabilizing-windows}
of stabilizing marking functions on finite windows):
\begin{defn}[Stabilizing graph \(\G\) on the input process \(\P\)]
\label{d:stabilizing-graph}
    We say that a family of graphs 
    $(\G(\P_n))_{n \in \N}$  (where each graph is generated by the same rule \(\mu \mapsto \G(\mu)\) represented by the marking neighborhood $\bar N$) forms a {\em stabilizing interaction graph} \(\G\) on finite windows with respect to a point process \(\P\) if   for every \(p \in \mathbb{N}\), there exist a fast decreasing function $\varphi_p'$ (see~\eqref{varphibd}), such that the  interaction range satisfies:
\begin{equation} \label{sdecayspin}
\sup_{1 \leq n < \infty} \sup_{x_1, \ldots, x_p \in W_n} \mP_{\bk x1p}
\bigl(S_n(x_1, \P_n) > s\bigr) \leq \varphi_p'(s), \quad s > 0.
\end{equation}
Without loss of generality, we assume that $\varphi_p'(\cdot)\le \varphi_{p+1}'(\cdot)$ for $p \in \N$.
\end{defn}
Examples of stabilizing graphs include the Gilbert disc graph (geometric graph), whose range  $S_n$ is  bounded by  the 
greatest integer in the sum of $1$ plus the disc radius, for all $n \in \N$.
 Proximity graphs from computational geometry, including those defined by $k$ nearest neighbors, Voronoi tessellations, or spheres of influence,  are stabilizing
with respect to homogeneous Poisson point processes and stationary $\alpha$-determinantal point processes.  This is spelled out in Appendix~\ref{s.admppinteraction}.  

We conclude this section by listing several properties of the interaction ranges of graphs (regardless of their probabilistic stabilization~\eqref{sdecayspin}). These properties will be used in various places when analyzing stabilizing graphs. %\BBLp{NOT YET COMPLETELY EDITED Appendix~A.}

Assuming  $x\in\mu \in \hat{\cN}_{W_n}$, the neighborhood of
$x$ is determined by the points in the ball  $B_{S_n(x,\mu)}(x)$ of radius $S_n$ around it, i.e., 
\begin{equation}
\label{e:cons_adm_gr}   
y \sim x  \, \, \mbox{in} \, \,  \mathcal{G}(\mu)  \, \, \, \mbox{iff}  \, \, \, y \sim x \, \, \mbox{in} \,  \, \mathcal{G}(\mu \cap B_{S_n(x,\mu)}(x)).
\end{equation}
Moreover, the ball $B_{S_n(x,\mu)}(x)$ is a stopping set, i.e.,  for all $s \geq 0$, we have 
\be \label{stopspin}
\{ \mu \in \hat{\cN}_{W_n} : S_n(x,\mu) \leq s \} = \{ \mu \in \hat\cN_{W_n} : S_n(x,\mu \cap B_s(x)) \leq s \},
\ee
or, equivalently,  
\be \label{stopspin-new}
S_n(x,\mu)=S_n\Bigl(x,(\mu\cap B_{S_n(x,\mu)}(x))\cup (\nu\cap B^c_{S_n(x,\mu)}(x))\Bigr)
\ee
for all $\mu,\nu\in \hat{\cN}_{W_n}$.
Finally, note that \eqref{e:cons_adm_gr} implies
\be \label{neighborulespin}
\mbox{$x \sim y$ in $\G(\mu)$ } \quad \Rightarrow  \quad 
|x- y| \leq \min\{S_n(x,\mu),S_n(y,\mu )\}.
\ee
\begin{remark}[Consistency of finite graphs]
In general, the stabilizing property~\eqref{sdecayspin} (defined only for finite windows) does not imply the consistency of the family of graphs $\G(\P_n)$, $n \in \N$, in the sense of the existence of some graph $\G$ on the entire input process $\P$ (e.g. the Gilbert disc graph or Voronoi tessellation) such that $\G(\P_n)$ suitably converges to it. Our rule $\mu \mapsto \G(\mu)$ applies only to finite input $\mu$, and consistently extending this rule to infinite input $\mu \in \cN_{\R^d}$ is possible under the additional condition of uniform boundedness of the interaction ranges, i.e., $\limsup_{n \to \infty} S_n(x, \P_n) < \infty$ almost surely for all $x \in \P$. 
Such an extension is possible, and
follows from arguments similar to those used  when extending the domain of  marking functions  on finite windows and possessing stopping set stabilization on finite windows to infinite input on infinite domains; see remark~\ref{ii.xi-extension-stopping set} in Section \ref{s:remarksltthms}.
\end{remark}

\subsubsection{Averaged weak spatial mixing of spins} \label{sss.spatial-mixing}
Dating back to the work of \citet{dobrushin1985constructive}, spatial mixing properties have been used to characterize the decay of correlations between spins induced by the Gibbs measure~\eqref{eqn:gibbs} on a graph $\G(\X)$ or a family of graphs. Among the several formulations of this concept, we adopt a relatively weak version and formulate it for  spin models on a random  stabilizing  graph \( \G \) on \( \P \). 

For subsets $\X''\subset\X'\subset \X$, by  
${\cal L}(V_{\X''} \mid {V}_{\partial \X'} = \v_{\partial \X'})$  
we mean the law of the spins on $\X''$ conditioned on the spin configuration $\v_{\partial \X'} \in \K^{\partial \X'}$ on $\partial \X'$, assuming that $\sP_{\G(\X)}({V}_{\partial \X'} = \v_{\partial \X'}) > 0$. Here, $\partial \X' := \partial_\X \X'$ is the set of all elements in $\X \setminus \X'$ that are at unit graph distance from $\X'$ and if   $\partial \X' = \emptyset$, we adopt the convention that ${\cal L}(V_{\X''} \mid {V}_{\partial \X'} = \v_{\partial \X'}) = {\cal L}(V_{\X''})$. When considering conditional laws, we implicitly assume that the boundary condition is non-degenerate as above i.e., $\sP_{\G(\X)}({V}_{\partial \X'} = \v_{\partial \X'}) > 0$. In general, when referring to boundary conditions we mean specifying spin values on a set $\cZ \subset \X$.

Recall the {\em total variation distance} $d_{\TV}$ between probability measures, as defined in Section \ref{s:prelim}.  {\em Weak spatial mixing}  of the Gibbs specification~\eqref{eqn:gibbs} on $\G(\X)$ says  that two different boundary conditions $\v_{\partial \X'}, \z_{\partial \X'}\in \K^{\partial \X'}$ on $\partial \X'$ result in two distributions of spins on $\X''$ which differ in the  total variation distance by a fast decaying function of the 
graph distance $d_{\G(\X)}(\X'', \partial \X')$ between $\X''$ and $\partial \X'$. In what follows, we will primarily take \(\X'\) to be the union of graph balls centered at the points in \(\X''\). Specifically, for \(x \in \X\) and \(k \in \N\), the \(\G(\X)\)-graph ball of radius $k$ is  
\[
\gB_k(x) := \{y \in \X : d_{\G(\X)}(x, y) \leq k\},
\]  
where \(d_{\G(\X)}(x, y)\) is the graph distance between \(x\) and \(y\). For a subset \(\X'' \subset \X\),  define the union of these graph balls by  
\[
\gB_k(\X'') := \bigcup_{x \in \X''} \gB_k(x).
\]
Now we give our formal definition of the spatial mixing property for spin models on random graphs.
\begin{defn}[Averaged weak spatial mixing of the Gibbs spin model on a random spatial  graph]
\label{d:ssm-pp}
The Gibbs spin model~\eqref{eqn:gibbs} satisfies {\em averaged weak spatial mixing} on the graph $\G = \G(\cdot)$ with respect to the input process $\P$ if
the following conditions hold: 
\begin{enumerate}[wide,label=(\roman*),labelindent=0pt]
%\vskip.2cm
\item \label{i.spatial-mixing} for any finite subset $\X \subset \R^d$, there exists a spatial mixing function 
$\MC_{\G(\X)}: 2^\X \times 2^\X \to [0, \infty)$ 
such that   
\begin{equation} \label{spatmix}
d_{\TV} \big( \mathcal{L}(V_{\X''} \mid V_{\partial \X'} = \v_{\partial \X'}),  
\mathcal{L}(V_{\X''} \mid V_{\partial \X'} = \z_{\partial \X'}) \big)  
\leq \MC_{\G(\X)}(\X'', \partial \X'),
\end{equation}  
for all subsets $\X'' \subset \X' \subset \X$ and boundary conditions $\v_{\partial \X'}, \z_{\partial \X'}\in \K^{\partial \X'}$, and 
\vskip.2cm
\item  \label{i.averaged-spatial-mixing} 
for all $p \in \N$ the functions 
\begin{equation} \label{e.spin-mixing-constant-bound}
\ga_p(k):=\sup_{B \in {\cal B}_b} \sup_{x_1,\ldots, x_p\in B}  
\sE_{\bk x1p} \big[\MC_{\G(\P \cap B)} (\bk x1p, \partial \gB_k(\bk x1p))\big],  \ \  k \in \N,
\end{equation}  
are fast decreasing as \( k \to \infty \).  
Here, \(\gB_k(\bk x1p) = \bigcup_{i=1}^p \gB_k(x_i)\) and the boundary $\partial$ of \(\gB_k(\bk x1p)\) is again with respect to~\(\G(\P \cap B)\);
${\cal B}_b$ are the bounded Borel subsets of $\R^d$. 
We  assume that the mapping $((\R^d)^p,\hat\cN_{\R^d\times\cN_{\R^d}})\owns(\bk x1p,\tilde\mu)=(\bk x1p,\G(\mu)) \mapsto \MC_{\G}(\bk x1p,\partial^\G\gB_k(\bk x1p))\in\R$ is measurable for all $k\in\N$. Without loss of generality, we assume that $\ga_p$ are increasing in $p$.
\end{enumerate}
\end{defn}

 Condition~\ref{i.spatial-mixing}  is a way of controlling the decay of spin correlations and is mainly dependent on the spin system whereas  condition~\ref{i.averaged-spatial-mixing} ensures fast decay with respect to graph distance in expectations and so depends on the graph model as well as the underlying point process.
More precisely, the function \(\MC_{\G(\X)}\allowbreak(\X'', \partial \X')\) in ~\eqref{spatmix} measures the sensitivity of the spin distribution at \(\X''\) to the boundary conditions of the Gibbs measure at \(\partial \X'\) and the property of spatial mixing requires the sensitivity to diminish rapidly with respect to the graph distance between $\X''$ and $\partial \X'$. Typically, in the setting of deterministic graphs such as lattices, this decay is assumed to be exponential i.e., 
$$ \MC_{\G(\X)}(\X'', \partial \X') \leq C |\X''| |\partial\X'|  \exp(-a \cdot d_{\G(\X)}(\X'', \partial \X')),$$
for some $C, a \in (0, \infty)$; see \cite[(2.5)]{martinelli1994approach}. Often this is also taken to be the definition of weak spatial mixing in such settings. While it would have been convenient to assume such an  exponential bound a.s. for spatial random graphs,
we have opted in  \eqref{e.spin-mixing-constant-bound} for a weaker version, one averaging the 
weak spatial mixing constant over the point process.  Hence we call this {\em  averaged weak spatial mixing} and this suffices to establish our general central limit theorem for Gibbs models on spatial random graphs. 

This weaker condition brings additional flexibility and enlarges the scope of applications.  Indeed, as seen in examples in Section~\ref{s:Ex_spin1},  weak spatial mixing bounds are derived via connective constant or site percolation probability bounds; see Sections~\ref{s:gen_disagreement} and~\ref{s:spin_Poisson}. Requiring such graph parameters to be almost surely 
well-controlled is sometimes very restrictive, possibly requiring almost sure degree bounds on  the graph. However, our requirement of  averaged weak spatial mixing imposes control  on the connective constant (also in its  averaged version) or site percolation probability only in expectation and hence it is easier to establish this condition for  a wider class of spin models and spatial random graph models as well as for a larger range of temperatures. 

In Section~\ref{s:hard-core}, using an example of  hard-core spin model, we will show how a  stronger version of spatial mixing for  Gibbs models, accounting for the exact location of the differences between the boundary conditions,
may allow one to  extend the range of model parameters under which the central limit theorem holds for specific Gibbs models. 

As highlighted in the introduction, our definition of spatial mixing is inspired by the notion of mixing introduced in \citet{dobrushin1985constructive}, a cornerstone in the study of spin systems. For infinite configurations (\(\mathrm{card}(\X) = \infty\)) and finite spin spaces \(\K\), weak spatial mixing guarantees the existence of a unique Gibbs measure on the graph configuration \(\G(\X)\); see \cite[Proposition 2.2]{weitz2005combinatorial}. For more on different notions of spatial mixing and their implications, see \cite[Section 2]{martinelli1994approach}. 

\subsection{Main result for spin models}
\label{mainspinresults}
We state a general central limit theorem for spatial spin models satisfying 
 averaged weak spatial mixing. Recall that $\rho^{(1)}_{\bk x1p}$ are the correlation functions of order $1$ under reduced Palm distributions $\sP^!_{\bk x1p}$; see Section~\ref{s:prelim}.

\begin{theorem}[CLT for the sum of spins]
\label{spinmain}   
Let $\P$ be a simple  point process on $\R^d$ with summable exponential mixing correlations as in Definition \ref{def.A2} and bounded reduced  Palm intensity function of the ground process $\P$, i.e., for all $p \in \N$
\be \label{bdedreduced}
\sup_{x_1,\ldots,x_p \in \mR^d} \sup_{y \in \R^d} \rho^{(1)}_{\bk x1p}(y) \leq \hat{\kappa}_p<\infty.
\ee
Let $\G = \G(\cdot)$ be a stabilizing interaction graph on finite windows with respect to $\P$ as in Definition~\ref{d:stabilizing-graph}. 
Assume the Gibbs specification~\eqref{eqn:gibbs} satisfies  averaged weak spatial mixing on~$\G(\P)$,  as in  Definition~\ref{d:ssm-pp}.
Denote by $\{V_x\}_{x \in \P_n}= \{V_x(\G(\P_n))\}_{x \in \P_n}$  the spin configuration on $\P_n$
assuming for all $n \in \N$ that  $Z(\P_n) > 0$.
Assume moreover  for all $p>1$ the following moment condition 
\begin{equation}\label{e:moment-conditions-spin}
\sup_{1 \leq n< \infty} \sup_{q\le p } \sup_{x_1,\ldots,x_q \in W_n} \sE_{\bk x1{q}}[ \max(1,|V_{x_1}(\G(\P_n))|^{p})]
< \infty\,.
\end{equation}

Put
 $$
 \mu_n^{\G} :=  \sum_{x \in \P_n} V_x \delta_{n^{-1/d}x}.
 $$
If  $f \in \B(W_1)$ satisfies
$\Var{\mu_n^{\G}(f)} = \Omega(n^{\nu})$ for some $\nu > 0$, then as $n \to \infty$
$$ (\Var{\mu_n^{\G}(f)})^{-1/2}(\mu_n^{\G}(f) - \sE[\mu_n^{\G}(f)]) \stackrel{d}{\Rightarrow} Z.
$$
\end{theorem}
We shall deduce this result from Theorem~\ref{t:cltmarkedpp} in the next section.  Before doing so,  we comment on extensions and refinements.  
\begin{remark}[Mean and variance asymptotics and further comments]\
\label{r:expvarspin}\
\begin{enumerate}[wide,label=(\roman*),labelindent=0pt]\item If the input process $\P$ is stationary and the interaction rule $\mu \mapsto \G(\mu)$  is translation invariant, then under the assumptions of Theorem~\ref{spinmain}, with the moment condition~\eqref{e:moment-conditions-spin} satisfied only for $p=2+\epsilon$ for some $\epsilon>0$, the mean asymptotics for  $n^{-1}\sE[\mu_n^{\G}(f)]$ and variance asymptotics for $n^{-1}\Var[\mu_n^{\G}(f)]$ can be expressed as in~\eqref{expasy} and~\eqref{eqn:var}, as stated in Proposition~\ref{expvar}. These results involve Palm distributional limits, given in~\eqref{e.xil-0} and~\eqref{e.xil-xy}, for $V_{x}(\G(\P_n))$ and $(V_{x_1}(\G(\P_n)),V_{x_2}(\G(\P_n)))$, without requiring the existence of an infinite spin system consistent with $(\G(\P_n))_{n \in \N}$. The proof of the central limit theorem in Theorem~\ref{spinmain} 
relies on 
\BL~cluster-localization for spins as defined in Definition \ref{def.Lp-stabilizing_marking-weak}, see Lemma~\ref{stablifting}.
\BL~cluster-localization
 is also sufficient to derive the mean and variance asymptotics
of Proposition~\ref{expvar}; see Remark~\ref{i.rem-BLCuster} in
Section~\ref{s:remarksltthms} and
Remarks~\ref{r:stab_lemma_BCL}, and~\ref{r:stab_lemma_BCL-variance}.

\item Still assuming stationarity of the model (both the input process and the graph) as above, the central limit theorem for $(\mu_n^{\G})_{n \in \N}$, as formulated in Theorem~\ref{spinmain}, holds under only a $(2 + \epsilon)$-moment assumption on the spins in \eqref{e:moment-conditions-spin}, provided 
that $\lim_{n \to \infty} n^{-1}\Var[\mu_n^{\G}(f)]$ 
exists and is non-zero. 
 This follows from  Theorem~\ref{t:multcltmarkedpp}, which moreover extends the result to multivariate spin systems.  Many classical spin systems are finite-valued and in these cases, the moment condition is trivially satisfied.

\item A general strategy is outlined in 
\cite[Section 1.6]{giardina2016annealed} for deriving more explicit expectation and variance asymptotics in terms of susceptibility and pressure but, as noted there, this strategy is more challenging in the averaged quenched setting which corresponds to our framework. 

\item It is sometimes easier  to derive volume-order variance lower bounds in `certain spin systems' than to establish limiting variance asymptotics. For example, using the approach as in \cite[Section 3.4.3 of v2]{dinh2021quantitative}, one can derive variance lower bounds for ferromagnetic or positively associated spin systems. 

\item For more general spin systems, one derive variance bounds as follows. Suppose ${\cal I}_n \subset \P_n$ is an independent set of vertices in the graph $\G(\P_n)$, i.e., $x,z \in {\cal I}_n$ implies that $x \not \sim z$. Using the conditional variance formula and the spatial Markov property twice, we derive that
\begin{align*}
    \Var{(\sum_{x \in \P_n} V_x)} & \geq \EXP{ \Var{(\sum_{x \in {\cal I}_n} V_x \mid V_y, y \in I^c_n)}} = \EXP{\sum_{x \in {\cal I}_n} \Var{(V_x \mid V_y, y \in I^c_n)}} \\
    & = \EXP{\sum_{x \in {\cal I}_n} \Var{(V_x \mid V_y, y \neq x)}}.
    \end{align*}
Next, we find an independent set ${\cal I}_n$ such that for a positive constant $c$, a.s. for all $x \in {\cal I}_n$ we have $\Var{(V_x \mid V_y, y \neq x)} \geq c$. This can be established under conditions on the spin system and possibly requiring uniform bound on degrees of ${\cal I}_n$; for example, see \cite[Lemma 3.6 of v3]{dinh2021quantitative}. Then we obtain 
$$  \Var{(\sum_{x \in \P_n} V_x)} \geq c \, \EXP{|{\cal I}_n|}.$$
One may show that $\EXP{|{\cal I}_n|} \geq c'n$ with $c' > 0$ for many stabilizing interaction graphs, including, for example, bounded degree graphs such as $k$-nearest neighbor graphs or graphs with bounded interaction range such as random geometric graph.  Combining the above observations,  we obtain volume order variance lower bounds for certain spatial spin systems. 

\item Averaged weak spatial mixing is here formulated with respect to the graph distance in order to exploit spatial mixing results on deterministic graphs; we may then transfer this to 'spatial mixing with respect to the Euclidean distance' via the stabilization of the interactions graphs. However in certain applications it may be easier to directly establish spatial mixing in the Euclidean distance; see Section~\ref{s:spin_Poisson}.

\item \label{i.particular-sets}
The conclusion of Theorem~\ref{spinmain} remains valid under a slightly weaker
form of the averaged weak spatial mixing assumption. Namely, the functions
$\ga_p$ in~\eqref{e.spin-mixing-constant-bound} may be defined by restricting
the supremum over bounded Borel sets $B \in \mathcal{B}_b$ to sets in
\begin{equation}
\nicesets
:=
\{\, W_n \cap B_r(\bk{x}{1}{l})
: n \ge n_0,\ r \ge r_0,\ l \in \N,\ x_1,\ldots,x_l \in W_n \,\},
\label{e.particular-sets}
\end{equation}
for some $n_0 \in \N$ and $r_0 \in (0,\infty)$. In other words, it suffices to
consider domains given by large windows $W_n$, possibly intersected with finite
unions of balls of sufficiently large radius; cf.\
\eqref{e.BL-spin-stab-1*} and~\eqref{e.BL-spin-stab-2*} in
Lemma~\ref{stablifting}. This restriction simplifies the verification of the
condition, especially when the graph $\G$ exhibits irregular behavior on
bounded sets that do not arise in the subsequent analysis. 

\item In Section~\ref{s:future}, we will discuss how this framework can be extended to include more general randomized graph models.

\item While the measure \( \mu_n^{\G}(\cdot) \) may be interpreted as the {\em empirical magnetization}---even if this is not immediately evident---we expect that our strategy of the proof of Theorem~\ref{spinmain} could also be applied to derive limit theorems for {\em internal energy} and {\em susceptibility}; see \cite[Chapter 5]{van2017stochastic}.
\end{enumerate}
\end{remark}

\subsubsection{Proof of CLT for spin models %(Theorem \ref{spinmain}) 
}
\label{s:gibbs_proof}

 The specification~\eqref{eqn:gibbs} of the laws of the spins
\( V_x = V_x(\G(\P_n)),\; x \in \P_n, \)
given $\P_n$ does not involve any particular marking function~$\xi$.
Without such a representation, a direct application of our umbrella
Theorem~\ref{t:clt_linear_marks} to prove a CLT would require verifying
fast \BL-mixing correlations of these spins on the input process~$\P$
in the sense of Definition~\ref{d.omegahmixing}, which appears to be a
delicate task.

Instead, we introduce auxiliary i.i.d.\ pre-marks on~$\tP$ and define a
\emph{spin marking function} on this extended input process, which is
used to construct spins whose joint laws are faithful to the given
Gibbs specification.
It is clear that, in general, any such marking function cannot satisfy stopping set
stabilization, nor even standard \(L^{q}\)-stabilization.
Consequently, our approach relies crucially on a purely distributional
notion of localization—namely bounded Lipschitz localization—of the
constructed marking function, and more precisely on
\BL~cluster-localization, developed in
Section~\ref{s:blcluster}.
Establishing this property for our marking function leverages the
averaged weak spatial mixing assumption of the spin model on the given
input process.
As a result, the proof of the CLT for the sum of spins stated in
Theorem~\ref{spinmain} follows as a direct application of
Theorem~\ref{t:cltmarkedpp} for \BL~cluster-localizing marks.

A more precise working plan for this section is as follows.
First, we present Lemma~\ref{mixingimpliesdbl}, which establishes a link
between averaged weak spatial mixing (holding almost surely with respect
to the input process~$\P$) and \emph{\BL~union-localization} of spins, an analogue of condition~\eqref{Lp-stab-infinite-weak1},
with localization measured in terms of the graph distance.
Next, Lemma~\ref{stablifting} establishes \BL~cluster-localization,
under the assumption of independence between well-separated clusters
 in~\eqref{Lp-stab-infinite-weak2}.
We then introduce a spin marking function and show in
Lemma~\ref{l:markfnspin} that it is consistent with the Gibbs
specification and, moreover, satisfies the independence
property between well-separated clusters assumed above.

This result allows us to conclude the proof of
Theorem~\ref{spinmain} by an application of
Theorem~\ref{t:cltmarkedpp}.
Finally, auxiliary lemmas—some of which are used in this programme are
 stated and proved in Section~\ref{sss.graph-lemmas}.

For $x \in \X$ and $m \in \N$ recall 
the notation  $\gB_m(x):= \{y \in \X: \ d_{\G(\X)}(x,y) \leq m\}$ for the $\G(\X)$-graph ball of radius $m$, 
and their union centered at points in $\X'' \subset \X$, namely 
$\gB_m(\X'') = \bigcup_{x \in \X''}\gB_m(x)$. We denote by
$\G^{\cup[m]}_{\X''}(\X):=\G(\gB_m(\X''))$ the {\em induced subgraph} of $\G(\X)$ on  $\gB_m({\X''})$ and we let
$V_{\X''}^{\cup[m]}:= V_{\X''}(\G^{\cup[m]}_{\X''}(\X))$ be the spins with respect to the Gibbs measure on $\G(\gB_m(\X''))$. 

To emphasize that the $d_{\BL}$ distance is a property of the laws of the random variables, we shall use here the $d_{\BL}$ distance between laws of random variables rather than random variables themselves.
\begin{lemma}[Averaged weak spatial mixing implies fast 
\BL~union-localization of spins] \label{mixingimpliesdbl}
The spatial mixing function $\MC_{\G(\X)}$ at ~\eqref{spatmix} controlling the  decay of spin correlations induced by the Gibbs measure~\eqref{eqn:gibbs}  on the graph $\G(\X)$ implies  for all $\X'' \subset \X$ and $m \in \N$
\begin{align}\label{spamixa}
d_{\BL} \big( {\cal L}(V_{\X''}),
{\cal L}(V_{\X''}^{\cup[m]})  \big) 
\leq 
d_{\TV} \big( {\cal L}(V_{\X''}),
{\cal L}(V_{\X''}^{\cup[m]})  \big) \leq \MC_{\G(\X)}(\X'',\partial \gB_{m-1}(\X'')).
\end{align}
\end{lemma}
\begin{proof}
The first inequality follows since $d_{\BL}$ is bounded by $d_{\TV}$. 
For the second inequality, put $\X' = \gB_{m-1}(\X'')$,  the $\G(\X)$-graph ball of radius $m-1$ around $\X''$ and so $\partial \X' = \partial \gB_{m-1}(\X'')$ . If $\partial\X'=\emptyset$ (when  $\gB_{m-1}(\X'')=\gB_{m}(\X'')$) then $d_{\TV} \big( {\cal L}(V_{\X''}),{\cal L}(V_{\X''}^{\cup[m]})  \big)=0$. Otherwise, considering  probabilities $\mP=\mP_{\G(\X)}$  and $\mP'=\mP_{\G^{\cup[m]}_{\X''}(\X)}$ of the Gibbs model at ~\eqref{eqn:gibbs},
respectively, on $\G(\X)$  and $\G^{\cup[m]}_{\X''}(\X)$,
we have
\begin{align*}
&d_{\TV} ( 
{\cal L}(V_{\X''} ),{\cal L}(V^{\cup[m]}_{\X''}\big))  \\
& =
\sum_{v_{\X''} \in \K^{\X''}} \Big|\sum_{ \z_{\partial \X'},w_{\partial \X'}
\in \K^{\partial \X'}  }
\Bigl(  \mP(V_{\X''} = v_{\X''} | {V}_{\partial \X'} = \z_{\partial \X'} ) 
   - \mP'(V^{\cup[m]}_{\X''} = v_{\X''} | {V}^{\cup[m]}_{\partial \X'} = w_{\partial \X'})\Bigr)\\[-2ex]
&\hspace{20em}\times \mP({V}_{\partial \X'} = z_{\partial \X'}) \mP'(V^{\cup[m]}_{\partial \X'} = w_{\partial \X'}) \Big|\\
& =
\sum_{v_{\X''} \in \K^{\X''}} \Big|\sum_{ \z_{\partial \X'},w_{\partial \X'}
\in \K^{\partial \X'}  }
\Bigl(  \mP(V_{\X''} = v_{\X''} | {V}_{\partial \X'} = \z_{\partial \X'} ) 
   - \mP(V_{\X''} = v_{\X''} | {V}_{\partial \X'} = w_{\partial \X'})\Bigr)\\[-2ex]
&\hspace{20em}\times \mP({V}_{\partial \X'} = z_{\partial \X'})
\mP'(V^{\cup[m]}_{\partial \X'} = w_{\partial \X'})\Big|\\[2ex]
&\le    
  \sum_{ z_{\partial \X'},w_{\partial \X'}
\in \K^{\partial \X'}  }
\sum_{v_{\X''} \in \K^{\X''}} \Big|  \mP(V_{\X''} = v_{\X''} | {V}_{\partial \X'} = \z_{\partial \X'} ) 
   - \mP(V_{\X''} = v_{\X''} | {V}_{\partial \X'} = w_{\partial \X'})\Big| \\[-2ex]
&\hspace{20em}\times \mP({V}_{\partial \X'} = z_{\partial \X'})
\mP'(V^{\cup[m]}_{\partial \X'} = w_{\partial \X'})\\[2ex]
&\le \MC_{\G(\X)}(\X'',\partial\X'),
  \end{align*}
  where  we  used the  consistency
  of conditional spin distributions  (Lemma \ref{consistency}) in the second equality and the averaged weak spatial mixing in the last inequality.
\end{proof}

We introduce some notation for the forthcoming crucial lemma, which
establishes \BL~cluster-localization of our spin marking function on
stabilizing random graphs. 
 The vector $(V_{x_1}(\G(\X)), \ldots, \allowbreak V_{x_p}(\G(\X)))$ is  abbreviated by $V_{\bk x1p}(\G(\X))$, where $x_i \in \X$, $i=1,\ldots,p$.
In the following lemma, we will take $\X$ to be  either the entire input process $\P_n$ in the window $W_n$,  or the process restricted to the union of Euclidean balls $B_r(\bk xab) = \bigcup_{i=a}^b B_r(x_i)$, $a,b\in\{1,p\}$ ; i.e., $\X = \P_n \cap B_r(\bk xab)$. 

\begin{lemma}[Fast  \BL~cluster-localization of spins  on stabilizing graphs]
\label{stablifting}
Let \(\P\) be a simple, stationary point process on \(\R^{d}\) with bounded
(reduced) Palm intensity function, as in~\eqref{bdedreduced}.
Let \(\G=\G(\cdot)\) be a stabilizing interaction graph on finite windows with
respect to~\(\P\), in the sense of Definition~\ref{d:stabilizing-graph}.
Assume that the Gibbs specification~\eqref{eqn:gibbs} satisfies averaged weak
spatial mixing on~\(\G\) with respect to~\(\P\), as in
Definition~\ref{d:ssm-pp}.
Then for all $p \in \N$ there are fast decreasing functions  $\varphi_p$ and constants such that 
\begin{align} \label{Lp-spin-stab}
\sup_{n \in \N} \sup_{x_1,\ldots,x_p \in W_n} d_{\BL, \bk x1p}(V_{\bk x1p}(\G(\P_n)), V_{\bk x1p}(\G(\P_n \cap B_r(\bk x1p))))  \leq 2\varphi_p(r), \ r > 0
\end{align}
and for all $l \in \{1,\ldots,p-1\}$ and $x_1,\ldots,x_p\in\R^d$ such that  $0< r < d(\bk{x}{1}{l},\bk{x}{l+1}{p})$, we have 
\begin{align}
\label{Lp-spin-stab2}
& 
d_{\BL,\bk{x}{1}{p}}
\left(V_{\bk x1p}(\G(\P_n \cap B_r(\bk x1p))),
\left(V_{\bk x1l}(\G(\P_n \cap B_r(\bk x1l))), V_{\bk x{l+1}p}(\G(\P_n \cap B_r(\bk x{l +1}p)))
\right)\right)
\;\le\; 2\varphi_p(r),
\end{align}
where $V_{\bk x1l}(\G(\P_n \cap B_r(\bk x1l))), V_{\bk x{l+1}p}(\G(\P_n \cap B_r(\bk x{l+1}p)))$ are conditionally independent (given $\P_n$) spin vectors with laws $\pi_{\bk x1l}$ and $\pi_{\bk x{l+1}p}$ induced respectively by graphs $\G(\P_n \cap B_r(\bk x1l))$ and $\G(\P_n \cap B_r(\bk x{l+1}p))$.
The above   bound also holds under the stronger total variation distance.
\end{lemma}

\begin{proof}
Fix $x_1,\ldots,x_p \in W_n$ and consider $x_i \in \P_n$, $i=1,\ldots,p$, with the aim of studying $\P_n$ under the Palm distribution $\Palm_{\bk{x}1p}$. Set $k_r = \lceil r^{\gbeta} \rceil$ for a fixed $\gbeta \in (0,1)$ and let $r \geq 1$ be large such that $k_r < r/3$.

Denote by $\GBkr=\GBkr(\P_n)$
the induced subgraph of $\G(\P_n)$ on $\gB_{k_r}(\bk x1p)$---the union of the $k_r$-hop balls around $x_1,\ldots,x_p$ in the $\G(\P_n)$-graph distance. Denote  by $\GBr=\G(\P_n \cap B_r(\bk x1p))$ the graph constructed on $\P_n$ restricted to the union of Euclidean balls $B_r(\bk x1p)$. We adopt also the shorthand notation  $\G=\G(\P_n)$.

We will assume the following probability bounds about stabilizing interaction graphs and complete the proofs of ~\eqref{Lp-spin-stab} and ~\eqref{Lp-spin-stab2}. After that we will derive these two bounds. For all $p \in \N$, there exist fast-decreasing functions $\varphi^{'''}_p$ such that
\begin{align}
\sup_{n \in \N} \sup_{x_1,\ldots,x_p \in W_n} \sP_{\bk x1p}\{\GBr \not \supset \GBkr\} & \leq \varphi^{'''}_p(r), \, \label{e:graphEucstabilization} 
\end{align}
and for all $l \in \{1,\ldots,p-1\}$ and $x_1,\ldots,x_p\in\R^d$ such that  $0< r < d(\bk{x}{1}{q},\bk{x}{q+1}{p})$, we have 
\begin{align}
\sup_{n \in \N} \sup_{x_1,\ldots,x_p \in W_n}  \sP_{\bk x1p}\{\GBkr
\neq \GBlkr \, \dot{\cup} \, \GBpkr\} & \leq \varphi^{'''}_p(r),  \, \label{e:graphdecoupling}
\end{align}
where $\dot{\cup}$ denotes disjoint union of graphs. We assume without loss of generality that $\varphi'_p$ are increasing in $p$. The bounds show that with high probability, graph balls are contained in large Euclidean balls and graph balls on well-separated vertices are disjoint.

We consider the vector of spins $\VP=\VP(\G)$, $\VBr=\VP(\GBr)$ and $\VBkr=\VP(\GBkr)$ at $\bk x1p\in \P_n$ in the respective graphs.
\bigskip

\noindent \uline{{\em Proof of~\eqref{Lp-spin-stab}.}}
Given the inclusion  \(\GBr \supset \GBkr\), we can apply the fact that averaged weak spatial mixing implies fast \BL-localization of spins from Lemma~\ref{mixingimpliesdbl} and we have assumed that the inclusion holds with high Palm probability \(\Palm_{\bk x1p}\) in~\eqref{e:graphEucstabilization}.  

Indeed, we apply the triangle inequality 
\begin{equation}
\label{e:dbltrineq}    
\quad d_{\BL, \bk x1p}(\VP, \VBr) 
\leq d_{\BL, \bk x1p}(\VP, \VBkr) + d_{\BL, \bk x1p}(\VBkr, \VBr).
\end{equation}
For the first  term, we condition on \(\mathcal{P}_n\) and apply averaged weak spatial mixing via~\eqref{spamixa} with $\X = \P_n$ and obtain: 
\begin{align}\nonumber%\label{e.BL-spin-stab-1}
 d_{\BL, \bk x1p}(\VP, \VBkr)&=
  \sup_{f \in BL} \left| \E_{\bk x1p} \left[\E[f(\VBkr) \mid \P_n]- 
\E[ f(\VP) \mid \P_n]\right]\right|\\
%&\le \sup_{f \in BL}  \E_{\bk x1p} \left[\left|\E[f(\VBkr) -  f(\VP) \mid \P_n]\right|\right]\no\\
&\le   \E_{\bk x1p} \left[\sup_{f \in BL}\left|\E[f(\VBkr) - 
 f(\VP) \mid \P_n]\right|\right]\no\\
 &\le   \E_{\bk x1p} \left[d_{\BL}(\mathcal{L}(\VBkr|\P_n),  
\mathcal{L}(\VP \mid \P_n))\right] \no\\
 & \leq  \sE_{\bk x1p}[\MC_{\G(\P_n)}(\bk x1p,\partial \gB_{k_r}(\bk x1p))] 
\label{e.BL-spin-stab-1*}\\
&\le \ga_p(k_r)\,.\no
\end{align} 
Recalling $k_r:=\lceil r^\gbeta \rceil$, the last expression is fast decreasing in $r$ by condition (ii) of Definition \ref{d:ssm-pp}. 

Now, for the second term in~\eqref{e:dbltrineq}, on the event
\(\GBr \supset \GBkr\), we again apply the averaged weak spatial mixing
property via~\eqref{spamixa}, with respect to
\(\GBr = \G(\P_n \cap B_r(\bk{x}{1}{p}))\);
on the complementary event, we use the trivial bound \(d_{\BL}\le 2\):
\begin{align}\nonumber 
 & \quad d_{\BL, \bk x1p}(\VBkr, \VBr) \no\\
& \leq E_{\bk x1p} \left[d_{\BL}(\mathcal{L}(\VBkr|\P_n),  
\mathcal{L}(\VBr \mid \P_n))  \right] \no \\
& \leq \E_{\bk x1p} \left[d_{\BL}(\mathcal{L}(\VBkr|\P_n),  
\mathcal{L}(\VBr \mid \P_n)) \1{\GBr \supset \GBkr} \right] + 2 \sP_{\bk x1p}\{\GBr \not \supset \GBkr\}\no\\
& \leq
   \E_{\bk x1p} \left[\MC_{\GBr}(\bk x1p,\partial \gB_{k_r}(\bk x1p))\1{\GBr \supset \GBkr}\right] + 2 \sP_{\bk x1p}\{\GBr \not \supset \GBkr\}\no\\
&  \le   \sE_{\bk x1p}[\MC_{\GBr}(\bk x1p,\partial \gB_{k_r}(\bk x1p))]+ 2 \sP_{\bk x1p}\{ \GBr \not \supset \GBkr \} \label{e.BL-spin-stab-2*}\\
&\le \ga_p(k_r)+ 2 \sP_{\bk x1p}\{\GBr \not \supset \GBkr\}\,, \label{e.BL-spin-stab-2**}
\end{align}
where in the third inequality we use~\eqref{spamixa} with $\X = \P_n \cap B_r(\bk{x}{1}{p})$
knowing that 
$$\G(\X) = \GBr \supset \GBkr = \G^{\cup[k_r]}_{\bk x1p}(\P_n) = \G^{\cup[k_r]}_{\bk x1p}(\X).$$
The first term, again by condition (ii) of Definition \ref{d:ssm-pp}, is fast decreasing in $r$ and the second term  is fast decreasing in $r$ by ~\eqref{e:graphEucstabilization}. Thus, we have proven \eqref{Lp-spin-stab}.\\

\noindent \uline{{\em Proof of~\eqref{Lp-spin-stab2}.}}
We keep the notation of the previous part regarding the points 
$\bk x1p \in W_n$ and the localization radius $r>0$. 
Moreover, for $l \in \{1,\ldots,p-1\}$ we assume that
\(
0 < r  <d(\bk{x}{1}{l},\bk{x}{l+1}{p})\).
Given $\P_n$ under the Palm distribution $\sP_{\bk x1p}$, in addition to $\VBr=\VP(\GBr)$, we consider the  
random vectors of spins
\[
\VBlr=\VPl(\GBlr), 
\qquad 
\VBpr=\VPp(\GBpr)
\]
defined on the respective graphs
\[
\GBlr:=\G(\P_n \cap B_r(\bk x1l)),\qquad
\GBpr:=\G(\P_n \cap B_r(\bk x{l+1}p)).
\]
Recall that we have assumed that $\VBlr$ and $\VBpr$ are conditionally independent given $\P_n$.

To evaluate the distance between the law of $\VBr$ and the product law of
$(\VBlr,\VBpr)$, we observe that the approximation $\VBkr$ of $\VBr$,
constructed in the preceding part, also admits a product structure.
More precisely, on the event
\(
\GBkr = \GBlkr \,\dot{\cup}\, \GBpkr
\),
which is assumed to hold with high Palm probability
$\Palm_{\bk x1p}$ (see~\eqref{e:graphdecoupling}),
the two groups of spins located at $\bk{x}{1}{l}$ and $\bk{x}{l+1}{p}$
do not interact.
Consequently, $\VBkr$ decomposes into two independent components.
Moreover, these components provide marginal approximations of
$\VBlr$ and $\VBpr$, respectively.
Here are the details.

Substituting \eqref{e:graphEucstabilization} in \eqref{e.BL-spin-stab-2**}, we conclude that
 $\VBkr$ approximates \(\VBr\) and,
similarly,  $\VBlkr$ approximates $\VBlr$
 and likewise $\VBpkr$  approximates $\VBpr$,
 namely we have the bound 
\begin{equation}
\label{e.BL-spin-stab-3}
      \max \{ d_\BL(\VBlr, \VBlkr), d_\BL(\VBpr, \VBpkr), d_{\BL}(\VBr,\VBkr) \} \le \ga_p(k_r) + 2 \varphi'''_p(r)
     \end{equation}
Above, we have used the monotonicity of $\ga_p$ in $p$, imposed without loss
of generality in Definition~\ref{d:ssm-pp}. We make the same assumption for
$\varphi'''_p$, namely that it is increasing in $p$.

Further, on the event $\GBkr = \GBlkr \, \dot{\cup} \, \GBpkr$
the two groups of spins
\[
(\VBlkr,\VBpkr) := \VBkr
\]
located at respective sites \(\bk x1l\) and \(\bk x{l+1}p\) and generated on the graph
\(\GBkr\) are independent. Indeed, here, we construct the spin system on the union of the two components.
The resulting independence is an intrinsic
property of the spin system itself and follows from the definition \eqref{eqn:gibbs} since the two
components of the underlying graph are disconnected on the event \(\GBkr = \GBlkr \, \dot{\cup} \, \GBpkr\).

Collecting the above arguments, we have
\begin{align}\nonumber%\label{e.BL-spin-stab-1}
&d_{\BL,\bk{x}{1}{p}}
\left(\VBr,
\left(\VBlr,\VBpr\right)\right) \no\\
&\le \E_{\bk x1p} \left[d_{\BL}\left(\mathcal{L}(\VBr\mid \P_n)\;,\;
\mathcal{L}(\VBlr \mid \P_n)\times \mathcal{L}(\VBpr \mid \P_n)\right)\right]\no \\
 &\le   
 \E_{\bk x1p} \left[d_{\BL}\left(\mathcal{L}(\VBr|\P_n), \mathcal{L}(\VBkr|\P_n)\right)\right.\no\\
&\quad\qquad\left.+d_{\BL}\left(\mathcal{L}(\VBkr|\P_n)\;,\; \mathcal{L}(\VBlr \mid \P_n)\times \mathcal{L}(\VBpr \mid \P_n)\right)\right]\no\\
 &\le
 \E_{\bk x1p} \left[\left\{d_{\BL}\left(\mathcal{L}(\VBr|\P_n), \mathcal{L}(\VBkr|\P_n)\right)\right.\right.\no\\
&\quad\qquad \left.\left.+ d_{\BL}\left(\mathcal{L}(\VBlkr \mid \P_n)\times \mathcal{L}(\VBpkr \mid \P_n)\;,\; \mathcal{L}(\VBlr \mid \P_n)\times \mathcal{L}(\VBpr \mid \P_n)\right)\right \} \right. \no \\
& \qquad \qquad \left. \; \1{\GBkr
= \GBlkr \, \dot{\cup} \, \GBpkr}\right]\no\\
&\hskip1.5cm + 2\sP_{\bk x1p} \bigl( \{\GBkr
\neq \GBlkr \, \dot{\cup} \, \GBpkr\} \bigr)\no\\
 &\le
 \E_{\bk x1p} \left[\left\{d_{\BL}\left(\mathcal{L}(\VBr|\P_n), \mathcal{L}(\VBkr|\P_n)\right) + d_{\BL}\left(\mathcal{L}(\VBlkr \mid \P_n),\mathcal{L}(\VBlr \mid \P_n)\right) \right.\right. \no \\
&\quad\qquad\left.\left.+d_{\BL}\left(\mathcal{L}(\VBpkr \mid \P_n), \mathcal{L}(\VBpr \mid \P_n)\right)
\right\} \right] + 2\sP_{\bk x1p} \bigl( \{\GBkr
\neq \GBlkr \, \dot{\cup} \, \GBpkr\} \bigr)\no\\
&\le 3\ga_p(k_r)+8\varphi_p^{'''}(r)\,,\no
\end{align}
where in the last equality we apply \eqref{e.BL-spin-stab-3} and \eqref{e:graphdecoupling}. This justifies the validity  of~\eqref{Lp-spin-stab2}
and completes the proof of Lemma~\ref{stablifting} assuming \eqref{e:graphEucstabilization} and \eqref{e:graphdecoupling}. \\

\noindent \underline{{\em Proof of~\eqref{e:graphEucstabilization} and \eqref{e:graphdecoupling}.}}
We  control the event in \eqref{e:graphEucstabilization} by a larger event and then observe that this larger event also controls the event in \eqref{e:graphdecoupling}. Then, we derive probability bounds for the  larger event to derive both bounds.

Towards this, first observe that
\begin{equation}
\{\GBr \not \supset \GBkr\}  \subset \bigcup_{i=1}^p 
\{ \GBr\not\supset \G_{x_i}^{+[k_r]} \},
\label{e.BL-spin-stab-3*}    
\end{equation}
where $\G_{x_i}^{+[k_r]}$ represents the induced subgraph of $\G(\P_n)$ on $\gB_{k_r}(x_i)$, including additionally the edges between points in $\gB_{k_r}(x_i)$ and $\gB_{k_r}(x_j)$, $j \neq i$, that exist in $\G(\P_n)$ and therefore in $\GBkr$. (These edges would otherwise not be included if we simply took the union of edges within each $\gB_{k_r}(x_i)$, for $i=1, \ldots, p$.)

We  control the last event in \eqref{e.BL-spin-stab-3*} for $x_1$ and note the same approach works for all $x_1,\ldots,x_p$. We  control the last term for $x_1$ by exploiting  the fact that if a path originating from $x_1$ of length $k_r$ in the graph distance is not contained within $\GBr$, then at least one point in $\P_n \cap B_r(x_1)$ must have an interaction range which is `very large' 
in the Euclidean metric. The details go as follows. 

Set $s_r = r^{1-\gbeta}/4$ with $\gbeta$ as above.  Recall the family of stabilization ranges \((S_n)_{n \in \N}\) associated to the graph $\G(\P_n)$  as required  by our assumption; see~\eqref{e:S_n=R_Wn}. Denote by $A_r = A_{r,n}(x_1)$
the event (configurations of points of $\P_n$)   such  that for all $x \in \P_n \cap B_{r/2}(x_1)$,  $S_n(x,\P_n) \leq s_r$. 

By the stopping set property~\eqref{stopspin-new} of $S_n$, and  $s_r\le r/2\le r$ (recalling $r \geq 1$), we have for all $x \in \P_n \cap B_{r/2}(x_1)$
\begin{equation}\label{e.stoping-Sn}
S_n(x,\P_n) = S_n(x,\P_n \cap B_{S_n(x,\P_n)}(x_1))=S_n(x,\P_n\cap B_{s_r}(x_1) )=S_n(x,\P_n\cap B_r(\bk x1p)).
\end{equation} 
Thus, for all $x \in \P_n \cap B_{r/2}(x_1)$,  the graph neighborhood consistency condition \eqref{e:cons_adm_gr}, implies  that the graph neighborhoods of $x$ in $\mathcal{G}(\P_n\cap B_r(\bk x1p))$ and $\mathcal{G}(\P_n)$ coincide.  Hence  
$$
\GBr=\G(\P_n \cap B_r(\bk x1p)) \supset \G(\P_n) \cap B_{r/2}(x_1)\quad  \text{on } A_r.
$$
Furthermore, on $A_r$, any path of (graph) length $k_r$ originating  from $x_1$ in $\mathcal{G}(\P_n)$  has Euclidean length of  at most $(k_r+1) s_r \leq r/2$.
Extending the length $k_r$ by $1$ accounts for possible  additional edges going from $\gB_{k_r}(x_1)$ to $\gB_{k_r}(x_j)$, $j\not=1$, included in  $\G_{x_i}^{+\cup[k_r]}$. In other words if $x_1,y_1,\ldots,y_{k_r}$ is a path in $\mathcal{G}(\P_n)$, then $|y_i-x_1| \le r/2$, $i = 1,\ldots,k_r$ and so $y_i \in B_{r/2}(x_1)$. Hence, by the previous argument,  on $A_r$ all paths in $\mathcal{G}(\P_n)$ starting from $x_1$, of graph length $k_r$,  are in $\GBr$, i.e., $ \GBr  \supset \G_{x_1}^{+[k_r]}$ on $A_r$. So, we have shown that
\begin{equation}
 \label{e:graphstab_Arn_bound}
 \{\GBr \not \supset \GBkr\} \subset \bigcup_{i=1}^p A_{r,n}(x_i)^c. 
\end{equation}

Furthermore, on the event
\(A_{r,n}(x_i)\) the subgraph \(\G_{x_i}^{[k_r]}\) of \(\G(\P_n)\)
induced by \(\gB_{k_r}(x_i)\) is contained in $B_{r/2}(x_i)$.  
Furthermore, if there exist edges in \(\G(\P_n)\) connecting
\(\gB_{k_r}(x_i)\) and \(\gB_{k_r}(x_j)\) for \(i \neq j\), then their
endpoints must lie in
\(B_{r/2}(x_i) \cap B_{r/2}(x_j)\).
Under the assumption that
\(r < d\bk x1l, \bk x{l+1}p)\), this intersection is empty for $i\in\{1,\ldots,q\}$ and $j\in\{q+1,p\}$. Hence,
on the event \(\bigcap_{i=1}^p A_{r,n}(x_i)\), the graph \(\GBkr\) is a
disjoint union $\GBlkr \, \dot{\cup} \, \GBpkr$. Thus, we also have

\begin{equation}
 \label{e:graphdecoupling_Arn_bound}
 \{\GBkr \neq \GBlkr \, \dot{\cup} \, \GBpkr \} \subset \bigcup_{i=1}^p A_{r,n}(x_i)^c.
\end{equation}

Hence, to prove both \eqref{e:graphEucstabilization} and \eqref{e:graphdecoupling}, it suffices to bound $\sP_{\bk x1p}(\bigcup_{i=1}^p A_{r,n}(x_i)^c)$,  which, by the union bound, reduces to bounding $\sP_{\bk x1p}(A_{r,n}(x_1)^c)$. Recalling assumption \eqref{bdedreduced}, which uniformly bounds the Palm intensity functions of~$\P$, we obtain
\begin{align}
\label{e:adm_graph_euc_dist}
\sP_{\bk x1p}(A_{r,n}(x_1)^c) & \leq \sP_{\bk x1p} \left( \bigcup_{x \in \P_n \cap B_{r/2}(x_1)} \{S_n(x,\P_n) > s_r \} \right) \nonumber \\
& \leq \int_{B_{r/2}(x_1)} \sP^{}_{\bk x1p,x}(  S_n(x,\P_n) > s_r)) 
\rho_{\bk x1p}^{(1)}(x) \md x + \sum _{i=1}^p\sP_{\bk x1p}(  S_n(x_{i},\P_n) > s_r)  \nonumber \\
&\leq  \theta_d \hat{\kappa}_p (\frac{r}{2})^d \varphi_{p+1}'(s_r)+ p\varphi_{p}'(s_r) \nonumber\\ 
& \leq (\theta_d \hat{\kappa}_p (\frac{r}{2})^d+p)  \varphi_{p+1}'(s_r)\nonumber  \\
& \leq \varphi_p^{''}(r),
\end{align}
for a fast decreasing function $\varphi_p^{''}$.  Thus \eqref{e:adm_graph_euc_dist} together with \eqref{e:graphstab_Arn_bound} and \eqref{e:graphdecoupling_Arn_bound} yields that \eqref{e:graphEucstabilization} and \eqref{e:graphdecoupling} hold with the fast-decreasing function $\varphi^{'''}_p = \min \{ p \varphi_p^{''},1\}$.
\end{proof}

Having established in Lemma~\ref{stablifting} the \BL~cluster-localization of spins, to complete the proof of Theorem~\ref{spinmain}, we need to represent spins as a marking function of a suitably marked point process. 
In this regard, we augment the ground  process
\(\P=\sum_{x\in\P}\delta_x\) to
\(\tP=\sum_{x\in\P}\delta_{(x,U(x))}\) by auxiliary pre-marks \(U(x)\),
where $U(x)$ are i.i.d. uniform random variables on $[0,1]$,
which allows the spins to be generated in accordance with the Gibbs
specification. We define a marking function~\(\xi\) that generates the
spin system~\eqref{eqn:gibbs} on the graph \(\G(\P\cap B)\) at all sites
\(x\in\P\cap B\), by leveraging the randomness of the auxiliary marks
\(U(x)\) attached to the points of the process.

To this end, we first introduce a deterministic set-function of the
pre-marks, obtained by summing them modulo~\(1\):
\begin{equation}\label{e.U-spins}
U(\tilde{\P}\cap B)
:= \sum_{x\in\P\cap B} U(x) \;\; \text{mod } 1,
\end{equation}
where \(B\subset\R^d\) is a bounded Borel set.

Next, to encode the law of the spins on \(\G(\P\cap B)\), recalling that $\K$ is countable, we construct a
\emph{partition $\oB$ of the unit interval}~\([0,1]\),
\[
\oB=\bigl\{ \oB_{(v_x : x\in\P\cap B)} \bigr\}_{(v_x : x\in\P\cap B)\in\K^{\P(B)}},
\]
indexed by all possible joint spin configurations on \(\P\cap B\).
Each element of the partition has Lebesgue measure equal to the Gibbs
probability of the corresponding configuration:
\begin{equation}\label{e.Gibbs-partition}
\bigl| \oB_{(v_x : x\in\P\cap B)} \bigr|
=
\mP_{\G(\P\cap B)}\bigl(
V_x(\G(\P\cap B)) = v_x \; : \; x\in\P\cap B
\bigr).
\end{equation}
This partition of \([0,1]\) is chosen to be a deterministic and
measurable function of the graph \(\G(\P\cap B)\),
that is,
\(\oB = \oB(\G(\P\cap B))\).
Its construction depends on the interaction potential and external
field functions \(\Psi\) and \(\Phi\) appearing in the Gibbs
specification~\eqref{eqn:gibbs}.
Moreover, since the graph construction is translation invariant, the
partition may be chosen to inherit this invariance.

Using these ingredients, we assign spins to all points
\(x\in\P\cap B\) according to the index of the partition element
\(\oB_{(v_x : x\in\P\cap B)}\) containing the value
\(U(\tilde{\P}\cap B)\).
That is, for a given point \(x\in\P\cap B\), we define the spin-marking
function 
$\xi:\R^d\times[0,1]\times\hat\cN_{\R^d\times[0,1]}\longrightarrow \K$ 
by
\begin{equation}\label{e.spin-marking}
\xi(\tx,\tilde{\P}\cap B)
:= v_x
\quad \text{whenever} \quad
U(\tilde{\P}\cap B) \in \oB_{(v_x : x\in\P\cap B)}.
\end{equation}

The following lemma justifies the relevance of this spin marking
function.

\begin{lemma}[Marking function for spin systems]
\label{l:markfnspin}
Consider the marked point  process  $\tilde{\P} = \sum_{x\in\P}\delta_{(x,U(x))}$
where $U(x)$ are i.i.d. uniform random variables on $[0,1]$. Consider the  spin-marking function $\xi$ given in~\eqref{e.spin-marking}, related to the Gibbs specification~\eqref{eqn:gibbs} via~\eqref{e.Gibbs-partition}.
For any bounded Borel set $B\subset\R^d$, jointly the values of this function 
$\{\xi(\tx, \tP_n \cap B)\}_{x \in \P\cap B}$, 
given a realization of $\P\cap B$,  are distributed as the spin configuration 
$\{V_x(\G(\P \cap B))\}_{x \in \P\cap B}$ under this Gibbs model. Further, if $B_1 \cap B_2 = \emptyset$ then the two collections of marks $\{\xi(x,\P \cap B_1)\}_{x \in \P\cap B_1}$ and $\{\xi(x,\P \cap B_2)\}_{x \in \P\cap B_2}$ are conditionally independent given $\P$. 
Moreover, if the mapping $\mu \mapsto \G(\mu)$ is translation invariant, 
then the marking function can also be constructed to be translation invariant.
\end{lemma}

\begin{proof}
Through the  construction, we have
\begin{align*}
\sP_{\P\cap B}(\,\xi(\tx,\tilde{\P}\cap B)=v_x:x\in\P\cap B) &= \sP(\,U(\tilde{\P}\cap B)\in \oB_{(v_x: x\in\P\cap B)}\,)\\
&= \mP_{\P\cap B}(V_x(\G(\P\cap B)) = v_x: x\in\P\cap B)
\end{align*}
for all spin configurations $(v_x: x\in\P\cap B)\in\K^{\P(B)}$ and hence the marking function has the distribution~$\mP_{\P \cap B}$ of the spin system~\eqref{eqn:gibbs}.

Let $B_1,B_2$ be disjoint sets. The distributions of $\{\xi(\tx,\tilde{\P}\cap B_1)\}_{x \in \P \cap B_1}$ and $\{\xi(\tx,\tilde{\P}\cap B_2)\}_{x \in \P \cap B_2}$ are functions of $U(\tP \cap B_1)$ and $U(\tP \cap B_2)$ respectively. $U(\tP \cap B_1)$ and $U(\tP \cap B_2)$ are conditionally independent given $\P$ for disjoint sets $B_1,B_2$ as they depend on independent collections of uniform random variables $\{U(x)\}_{x \in \P \cap B_1}$ and $\{U(x)\}_{x \in \P \cap B_2}$ respectively. Thus, the assertion on conditional independence for marking functions follows.
Finally, as explained below \eqref{e.Gibbs-partition}, the asserted  translation invariance also holds by construction.
\end{proof}

The marking function introduced above and considered  in Lemma~\ref{l:markfnspin} does not satisfy
\BL-locali\-zation~\eqref{Lp-stab}.  
The obstruction stems from overlaps of the balls \(B_r(x_i)\) around points
\(x_i \in \bk{x}{1}{p}\), which induce arbitrary dependencies among the components
of \(\bk{\xi^{(r)}}{1}{p}(\tx,\tP_n)\).  
These dependencies need not match those present in
\(\bk{\xi}{1}{p}(\tx,\tP_n)\) and typically become more pronounced as \(r\) increases.

Although the construction of spins via a marking function \(\xi\) on the same
input process~\(\tP\) is not unique---for instance, one may modify the function
\(U(\tP \cap B)\) in~\eqref{e.U-spins} while preserving uniform marginals on
\([0,1]\)---it appears difficult to devise a representation that avoids this
issue at the level of \BL-localization. Nevertheless, Lemma \ref{stablifting} shows that the proposed marking function does satisfy \BL~cluster-localization, positioning us to complete the proof of Theorem \ref{spinmain}.

\begin{proof}[Proof of Theorem~\ref{spinmain}] 
We deduce the result from Theorem~\ref{t:cltmarkedpp}.  
Let $\tP = \{(x,U(x))\}_{x \in \P}$, where $U(x)$ are i.i.d.\ uniform $[0,1]$-valued random variables.  
The assumed summable exponential mixing correlations of $\P$, together with Proposition~\ref{p:pmix+mmix}, imply that $\tP$ has summable exponential $\B$-mixing correlations as in Definition~\ref{def.A2}.  
On this extended input process, we consider the spin marking function $\xi$ of Lemma~\ref{l:markfnspin}, which reproduces the law of the spins:  
\[
\bk \xi1p(\tx,\tP_n ) \;\overset{d}{=}\; V_{\bk x1p}(\G(\P_n )).
\]
Furthermore by Lemma~\ref{l:markfnspin}, if $B_r(\bk x1l) \cap B_r(\bk x{l+1}p)=\emptyset$ for some $l \in \{1,\ldots,p-1\}$, then the random vectors $\bk \xi1l(\tx,\tP_n \cap B_r(\bk x1l))$ and $\bk \xi{l+1}p(\tx,\tP_n \cap B_r(\bk x{l+1}p))$ are conditionally independent give $\P_n$ and are distributed as $V_{\bk x1l}(\G(\tP_n \cap B_r(\bk x1l)))$ and $V_{\bk x{l+1}p}(\G(\tP_n \cap B_r(\bk x{l+1}p)))$, respectively.

Thus, Lemma~\ref{stablifting} guarantees fast \BL~cluster-localization of this spin marking function on finite windows of~$\tP$, as per Definition~\ref{def.Lp-stabilizing_marking-weak}\ref{i.BL-localizing-weak} (comprising the two conditions~\eqref{Lp-stab-infinite-weak1} and~\eqref{Lp-stab-infinite-weak2}), together with the further qualifications (ii) and (iii).  
Hence, all the assumptions of Theorem~\ref{t:cltmarkedpp} are satisfied, and the conclusion of Theorem~\ref{spinmain} follows.
\end{proof}

\subsubsection{Auxiliary statements  for spin systems}\label{sss.graph-lemmas}

Recall  \( \X'' \subset \X' \subset \X \) as well as the notation of graph-ball unions $\gB_m(\X'')$ centered at points in $\X''$ and the subgraph $\G^{[m]}$ induced by them, 
introduced prior to Lemma \ref{mixingimpliesdbl}. The factorization property of Gibbs distributions over cliques implies {\em the spatial Markov property}, a fundamental tool in the analysis of Gibbs processes. This implies that the law of the spin vector $V_{\X''}= V_{\X''}(\G(\X))$ at ${\X''}$, when conditioned on spins outside a region is influenced only by the values on the boundary of the region. A consequence is that the spin vector $V_{\X''}= V_{\X''}(\G(\X))$ on ${\X''}$ conditioned on the boundary specification on $\partial \gB_{m-1}(\X'')$, coincides with the law of spin $V_{\X''}^{[m]}$ conditional on the same boundary specification.  Formally, this is stated as follows. 
\begin{lemma}[Spatial Markov property; \protect{\cite[Proposition 1]{lauritzen1990independence}}]
\label{l:spmarkspin}
Let $\cZ \subset (\X' \cup \partial \X')^c$. Then $V_{\X''}$ is conditionally independent of $V_{\cZ}$ given $V_{\partial \X'}$.
Thus, if $\sP_{\G(\X)}(V_{\cZ \cup \partial \X'} = v) > 0$ for $v = v_{\cZ \cup \partial \X'}$ then 
$$ {\cal L}(V_{\X''} \mid V_{\cZ \cup \partial \X'} = v) = {\cal L}(V_{\X''} \mid V_{\partial \X'} = v_{\partial \X'}).$$
\end{lemma}
\begin{lemma}[Consistency of conditional spin distributions]
\label{consistency}
When $\sP_{\G(\X)}(V_{\partial\gB_{m-1}(\X'')} =v) > 0$ for $v=v_{\partial\gB_{m-1}(\X'')}$ we have 
$$
{\cal L}(V_{\X''} | V_{\partial\gB_{m-1}(\X'')} = v) = {\cal L}(V_{\X''}^{[m]} | V_{\partial\gB_{m-1}(\X'')} = v).
$$ 
\end{lemma}
The spatial mixing function \(\MC_{\G(\X)}(\X'',\partial \X')\) in~\eqref{spatmix}, controlling the decay of spin correlations on the graph \(\G(\X)\), may be refined by introducing a function
 \(\MC'_{\G(\X)}: 2^\X \times 2^\X \to [0, \infty)\)
 which  explicitly accounts for disagreements in the boundary conditions on \( \cZ \subset \X\)  with $\cZ \cap \X'' = \emptyset$ and which satisfies an analog of ~\eqref{spatmix}, namely
\begin{equation}\label{ssmsingle}
d_{\TV} \left( \mathcal{L}(V_{\X''} \mid V_{\cZ} = \v_{\cZ}), 
\mathcal{L}(V_{\X''} \mid V_{\cZ} = \z_{\cZ}) \right) \leq 
\MC'_{\G(\X)}(\X'', \cZ_{\neq})\,,
\end{equation}
where \( \cZ_{\neq} := \cZ_{v \neq z} := \{y \in \cZ : \v_{\cZ}(y) \neq \z_{\cZ}(y)\) \}. Though \(\cZ_{\neq}\) depends on the boundary conditions $v$ and $z$, we shall suppress this dependence for notational convenience. Compared to~\eqref{spatmix}, the function 
\( \MC'_{\G(\X)} \) in~\eqref{ssmsingle} considers a potentially strict subset of \( \cZ \), which may be farther from \( \X'' \) than \( \cZ \) itself.  The function \( \MC'_{\G(\X)} \) is typically decreasing and often exhibits exponential decay
\begin{equation}
\label{e:ssmliterature} 
\MC'_{\G(\X)}(\X'', \cZ_{\neq}) \leq C|\X''| \,  \,\exp(-ad_{\G(\X)}(\X'', \cZ_{\neq})),
\end{equation}
for some constants $a, C >0$. In this case, the spin model is said to exhibit {\em strong spatial mixing} and we call \( \MC'_{\G(\X)} \)  the {\em strong spatial mixing bound}. We demonstrate that the strong spatial mixing bound can be expressed as a sum of strong spatial mixing bound over individual sites.
\begin{lemma}[Additive representation for the strong spatial mixing bound] \label{p:ssmsingle-wsm}
If  the function \(\MC'_{\G(\X)}\) satisfies
~\eqref{ssmsingle}  for all singletons  \( \X'' = \{x\}  \subset \X \)
then the following inequality  holds
for any configuration of sets \( \X'', \cZ \subset \X \)  such that $\cZ \cap \X'' = \emptyset$,
\be \label{ssmsingle-mutiple}
d_{\TV} \left( \mathcal{L}(V_{\X''} \mid V_{\cZ} = \v_{\cZ}),
\mathcal{L}(V_{\X''} \mid V_{\cZ} = \z_{\cZ}) \right) \leq 
\sum_{y\in\X''} \MC'(\{y\},\cZ_{\neq}).
\ee
\end{lemma}
\begin{proof}
We will show that if~\eqref{ssmsingle-mutiple} 
hold for  all subsets \( \X'' \subset \X \) with \( |\X''| \le k-1 \), for $k\ge 2$, then it also holds for \( |\X''|= k \). Thus the  lemma follows by induction, since by the  assumption for a single site, ~\eqref{ssmsingle-mutiple} holds for \( k = 1 \).

Let $\X'' \subset \X$. Let $w_{\X''}=(w_x:x\in\X'')\in \K^{\partial \X''}$. We abbreviate the events $\{{ V}_{\cZ} = v_{\cZ}\}$ and $\{{ V}_{\cZ} = z_{\cZ}\}$ by $A_v$ and $A_z$ respectively. 

Considering the probability $\mP:=\mP_{\G(\X)}$ of the Gibbs model~\eqref{eqn:gibbs},  conditioning on $V_{x}$ and using the triangle inequality, we obtain %
\begin{align*}
& |\sP(V_{\X''} = w_{\X''} | A_v) - \sP(V_{\X''} =w_{\X''} | A_z) |  \\
& \leq  | \sP(V_{\X''\setminus\{x\}} =  w_{\X''\setminus\{x\}}  | V_{x} = w_x, A_v) - \sP(V_{\X''\setminus\{x\}} =  w_{\X''\setminus\{x\}} | V_{x} = w_x, A_z) | \cdot \sP(V_{x} = w_x | A_v) \\
& \quad  +  \sP(V_{\X''\setminus\{x\}}  = w_{\X''\setminus\{x\}}| V_{x} = w_x, A_z) \cdot| \sP(V_{x} = w_x | A_v) -  \sP(V_{x} = w_x | A_z) |. 
\end{align*}
Summing over $w_{\X''}$'s and using that ~\eqref{ssmsingle-mutiple} holds for $|\X''| \in \{1,...,k-1\}$ we obtain 
\begin{align*}
& d_{\TV} \big( {\cal L}(V_{\X''} | { V}_{\cZ} = v_{\cZ}),
{\cal L}(V_{\X''} | {V}_{\cZ} = z_{\cZ}) \big)   \\
&  \leq d_{\TV} \big( {\cal L}(V_{x} | { V}_{\cZ} = v_{\cZ}),
{\cal L}(V_{x} | {V}_{\cZ} = z_{\cZ}) \big) \\
& \quad + \sum_{w_x\in\K}d_{\TV} \big( {\cal L}(V_{\X''\setminus\{x\}} | V_{x} = w_x, { V}_{\cZ} = v_{\cZ}) ,
{\cal L}(V_{\X''\setminus\{x\}} | V_{x} = w_x, {V}_{\cZ} = z_{\cZ}) \big)
\\[-2ex]&\hspace{25em}\times \sP(V_{x} = w_x | {V}_{\cZ} = v_{\cZ})   \\
& \leq \MC'(\{x\}, \cZ_{\neq})  + \sum_{w_x\in\K} \MC'(\X''\setminus\{x\}, (\cZ \cup \{x\})_{\neq})
\sP(V_{x} = w_x | { V}_{\cZ} = v_{\cZ}) \\[2ex]
& \leq \MC'(\{x\}, \cZ_{\neq})+  \MC'(\X''\setminus\{x\}, \cZ_{\neq})  \\
&= \sum_{x\in\X''} \MC'(\{x\},\cZ_{\neq}),
\end{align*}
where in the penultimate inequality, we have used that $\cZ_{\neq} = (\cZ \cup \{x\})_{\neq}$. This proves \eqref{ssmsingle-mutiple}. 
\end{proof}

If strong spatial mixing~\eqref{ssmsingle} holds with bounds as in \eqref{e:ssmliterature}, then it implies averaged weak spatial mixing~\eqref{spatmix} with \( \MC(\X'', \partial\X') = \MC'(\X'', \partial\X') \).  Furthermore, if one establishes the stronger version only for single sites and uses~\eqref{ssmsingle-mutiple}, relying  on a straightforward extension of the proof of \BL~union-localization~\eqref{spamixa}, we obtain the following :
\begin{align}\label{spamixa-strong}
&d_{\BL} \big( {\cal L}(V_{\X''}),
{\cal L}(V_{\X''}^{[m]})  \big) \\
&\le    
  \sum_{ z_{\partial \X'},w_{\partial \X'}
\in \K^{\partial \X'}  }
 \MC'(\X'',\partial_{\neq}\X') \, \mP_{\G(\X)}({V}_{\partial \X'} = z_{\partial \X'}) \, \mP_{\G(\X)}(V^{[m]}_{\partial \X'} = w_{\partial \X'}), \no 
\end{align}
where $\partial_{\neq}\X'$ denotes $(\partial \X')_{\neq}$ and notation has been slightly modified for convenience.
However, the right-hand side of~\eqref{spamixa-strong} depends not only on the graph structure but also on the probabilities associated with the Gibbs specification. In Section~\ref{s:hard-core}, we will apply this approach to the hard-core model.

\subsection{Averaged weak spatial mixing models via disagreement percolation}
\label{s:Ex_spin1}

We take up the question of determining which spin models satisfy the averaged weak spatial mixing condition  in Definition \ref{d:ssm-pp}.  In this section we exploit the  {\em disagreement percolation} approach introduced by \citet{van1994disagreement} and then use it to establish  averaged weak spatial mixing for spin models on random graphs.
This method involves analyzing an auxiliary independent site percolation model on the given graph, where open paths represent disagreements in the spin configurations caused by differing boundary conditions. More precisely, this auxiliary model allows one to stochastically dominate any path of disagreement from the boundary to a given site by an independent Bernoulli site percolation on the graph. The probability of this dominating percolation, defined in what follows, is the essential component in establishing averaged weak spatial mixing on random graphs.

To prevent percolation in this auxiliary model (and to ensure the fast decay of the functions \( (\ga_p)_{p\geq 1} \) in ~\eqref{e.spin-mixing-constant-bound}), we consider two alternative  conditions
in this general disagreement percolation approach:
 In Section~\ref{s:gen_disagreement}, we use bounds on the {\em  averaged connective constant}, cf. \citet{sinclair2017spatial}, and a standard self-avoiding walk counting argument from percolation theory. In Section~\ref{s:spin_Poisson}, we refine this approach by incorporating the concept of a {\em sharp phase transition} (cf. \citet{duminil2017new}) into the disagreement percolation framework.

Henceforth we assume that the spin space $\K$ is finite. The probability of an open site in the auxiliary disagreement percolation model is given by the {\em maximum influence} of the neighbors on sites, formalized in
\cite{van1994disagreement} as
\be \label{defq}
 q(\G(\X)) :=  \max_{x \in \X} \sup_{\v, \z  \in \K^{|\X| -1}} d_{\TV} \big( {\cal L}(V_{x} | {V}_{N_x} = \v) ,
{\cal L}(V_{x} | {V}_{N_x} = \z) \big),
\ee
where the distribution of the spins $V$ follows  the Gibbs specification~\eqref{eqn:gibbs} on the graph $\G(\X)$, and where $N_x$ denotes the set of neighbors of $x$. The maximum influence depends on both the graph structure and  also on the specification of the spin model.
We state the key result from disagreement percolation \cite[Corollary~1]{van1994disagreement}
{\em which controls the decay of spin correlations at 
 ~\eqref{spatmix} by the decay of disagreement percolation probabilities}. 
  
Let $\X'' \subset \X' \subset \X$ and assume that $\partial \X' \neq \emptyset$.
Then 
\begin{align}
\label{eq:gendisagbd}
d_{\TV} \Big(
\mathcal{L}( V_{\X''} \mid V_{\partial \X'} = \v_{\partial \X'} ),
\mathcal{L}( V_{\X''} \mid V_{\partial \X'} = \z_{\partial \X'} )
\Big)
\le
\mathbb{P}_q\!\left(
\X'' \stackrel{\X'_{*q}}{\longleftrightarrow} \partial \X'
\right),
\end{align}
where, $\X'_{*q}$ denotes the independent thinning of $\X'$, in which vertices
are retained (i.e., declared open) with probability
$q = q(\G(\X))$ as defined in~\eqref{defq}.
The event $\{ \X'' \stackrel{\X'_{*q}}{\longleftrightarrow} \partial \X' \} = \{ \X'' \stackrel{\perc{\G(\X)}{\X'}{q}}{\longleftrightarrow} \partial \X' \}$  in~\eqref{eq:gendisagbd} thus corresponds to the existence of a path
in the graph $\G(\X)$ that connects $\X''$ to $\partial \X'$ and that is entirely
contained in the vertex set $\X'_{*q} \cup (\X \setminus \X')$.
Equivalently, this defines a site percolation model on $\G(\X)$ in which
vertices in $\X'$ are independently open with probability~$q$, while all
vertices in $\X \setminus \X'$ are deterministically open. Moreover, when the underlying graph $\G(\X)$ is clear from the context, we omit it in the above connectivity event. The probability measure $\sP_q$ is taken with respect to this thinning (percolation) randomness.
If $\partial \X' = \emptyset$, then the left-hand side of
\eqref{eq:gendisagbd} vanishes and the inequality holds trivially.

We now bound maximum influence when $\X$ is replaced by $\P$ and when $\G(\P)$ is a stabilizing graph with not too many self-avoiding paths.

\subsubsection{Controlling disagreement percolation probabilities via averaged connective constant}
\label{s:gen_disagreement}
Connective constants allow one to bound the (usually exponential) growth in the number of self-avoiding paths of length \(m\) starting from an arbitary vertex in a graph. We use this notion to bound disagreement percolation probabilities and hence obtain the decay of spin correlations through \eqref{eq:gendisagbd}.
We adapt these concepts to our context of a spatial random graph \(\G\) defined over \(\P\), defining an  averaged connective constant as follows.

\begin{defn}[Averaged connective constant for spatial random graph]\label{d:connective-constant-pp}
For a graph $\G = \G(\cdot)$ on the input process $\P$, the  
 {\em averaged connective constant} $\bar{\Delta} \in [0,\infty]$
is defined as 
\begin{equation}\label{e:bardelta}
\bar{\Delta} := \inf\Bigl\{\delta\ge 0: \forall\, p \in \N \;  \exists \, m_p \in \N \text{\ such that for all $m\ge m_p$}, \,
\gn_{p}(m)\le \delta^m\Bigr\},
\end{equation}
where \begin{equation}\label{e.connective-constant-pp}
\gn_p(m):=\sup_{B\in\nicesets}\sup_{x_1,\ldots, x_p\in B}
  \sE_{\bk x1p}[N(x_1, m; \G(\P\cap B))],
\end{equation}
 with \(N(x_1, m;\G(\P\cap B))\) being the number of  self-avoiding paths of length \(m\) starting from \(x_1 \in \P \cap B\) in  $\G(\P\cap B)$, where $\nicesets$, defined in~\eqref{e.particular-sets}, are the sets of the form large windows
$W_n$, possibly intersected with finite unions of balls of sufficiently large
radius.
\end{defn}
The {\em quenched version} of this constant,
given in the sequel by $\Delta$ as in \eqref{connconstant}, is defined by removing the expectation in~\eqref{e.connective-constant-pp} and hence yielding the almost sure inequality \(N(x_1, m;\G(\P\cap B)) \leq \Delta^m\) for large $m \in \N$; this shall be used in applications to hard-core models in Section \ref{s:hard-core}. Bounds for this quenched version, which upper-bounds the  averaged version, exist for certain random graphs. Notably, this holds for the \(k\)-nearest neighbor graph, as it has bounded degree (see \citet[Lemma 8.2]{yukich1998probability}).  \citet{lacoin2014non} provides an example illustrating the sharp inequality  between  averaged and quenched connective constants.

The next result, derived using Theorem \ref{spinmain}, proves normal convergence for spin systems on spatial random graphs whose maximal influence is strictly smaller than the inverse of the averaged connective constant.
\begin{proposition}[CLT for sum of spins for  graphs with finite averaged connective constant] 
\label{c:cltspindisag}  
Let \(\P\) be a simple, stationary point process on \(\R^d\) with summable exponential mixing correlations as in Definition \ref{def.A2}, and assume that the (reduced) Palm intensity function of the ground process \(\P\) is bounded at all orders, as specified in \eqref{bdedreduced}. Let \(\G = \G(\cdot)\) be a stabilizing interaction graph on finite windows with respect to \(\P\) as in Definition~\ref{d:stabilizing-graph}, with finite averaged connective constant \(\bar\Delta<\infty\)  defined in Definition~\ref{d:connective-constant-pp}. If, for some $\epsilon>0$,  all sets \( B \in\nicesets \)
defined in~\eqref{e.particular-sets}, and almost every realization \(\P \), the maximum influence \( q(\G(\P \cap B)) \) as in \eqref{defq} of the Gibbs specification~\eqref{eqn:gibbs} on $\G(\P\cap B)$ with finite spin space $\K$ satisfies the inequality 
\be \label{influencebound}
q(\G(\P \cap B)) \le (\epsilon+\bar{\Delta})^{-1}
\ee
 then the spin configuration \( \{V_x(\G(\P_n))\}_{x \in \P_n} \) satisfies the CLT as stated in Theorem~\ref{spinmain}, provided the variance condition there is also met.
\end{proposition}
Before giving the proof, we provide several examples where one may verify the assumption \eqref{influencebound} on the influence parameter $q(\G(\P\cap B))$. We choose classical examples from \cite[Section 3]{van1994disagreement} and \cite[Section 4.3]{chatterjee2005concentration}, as upper bounds on the maximum influence $q(\G(\X))$ have been computed in these models for general graphs. More examples are in \cite[Section 3]{georgii2001random} and \cite[Section 8]{georgii2011gibbs}. In the examples below, we take $\X = \{x_1,\ldots,x_m\}$ and $v_i$'s are taken to be $\K$-valued.
\begin{enumerate}[wide,label=(\roman*),  labelindent=0pt]
\item {\em Hard-core model.} We associate \(\{0,1\}\)-valued random variables (spins) to each vertex in \(\G(\X)\), ensuring that no neighboring vertices both take the value \(1\),   with joint probability proportional to
\begin{equation}
 \label{eqn:hardcore}   
\pi_{\bk x1m}(\bk v1m) = \prod\limits_{(x_i,x_j) \in \G(\X) } \1{v_i v_j = 0} \prod\limits_{x_i \in \X} \lambda^{v_i},  \, v_i \in \{0,1\} \, ,
\end{equation}
where \(\lambda > 0\) is the {\em activity parameter}.
This is  the Gibbs model~\eqref{eqn:gibbs}, with  \(\K = \{0,1\}\),
\(\beta = 1\),
$\Psi(v_i,v_j) = \log \1{v_iv_j=0}$, $\lambda = e^{\gamma}$, and $\Phi(v) = v$. Loosely speaking, small $\lambda$-values are expected to result in weak correlations between spins at distant sites, a chief requirement for our general central limit theorems. 
Indeed, it can be shown that \( q(\G(\X)) \leq \frac{\lambda}{1 + \lambda} \) for any graph \( \G(\X) \), so Proposition~\ref{c:cltspindisag} already applies for 
\( \lambda < (\max(0,\bar{\Delta}-1))^{-1} \) for graphs with finite averaged connective constant
$\bar{\Delta} <\infty$. In Section~\ref{s:hard-core}, we will revisit  the hard-core model 
and prove the central limit theorem using a strong version of the spatial mixing approach, though under assumptions on the quenched, rather than  averaged, connective constant.

\item {\em Widom-Rowlinson model.} Here $\K = \{-1,0,+1\}$ and the spin variables are chosen randomly with joint probability proportional to
\begin{equation}
 \label{eqn:WR}   
\pi_{\bk x1m}(\bk v1m)=\prod\limits_{(x_i,x_j) \in \G(\X) }\1{v_iv_j \neq -1} \prod\limits_{x_i \in \X}\lambda^{|v_i|}, 
\end{equation}
where $\lambda$ is the {\em activity parameter}. This is  the Gibbs model in \eqref{eqn:gibbs} with $\beta = 1$, $\Psi(v_i,v_j) = \log \1{v_iv_j \neq -1}$, $\lambda = e^{\gamma}$, $\Phi(v) = |v|$. Here if $\G(\X)$ has a vertex of degree at least $2$ then $q(\G(\X)) \leq \frac{2\lambda}{1 + 2 \lambda}$.  Thus, Proposition
\ref{c:cltspindisag}
holds in this model for $\lambda < (2\max(0,\bar{\Delta}-1))^{-1}$. 

\item {\em Dobrushin's interdependence matrix.} Consider a  model that satisfies
$$ d_{\TV} \big( {\cal L}(V_{x} | {V}_{\X \setminus {x}} = \v) ,
{\cal L}(V_{x} | {V}_{\X \setminus {x}} = \z) \big) \leq \sum_{y \in \X \setminus {x}}A_{x,y}\1{\v_y \neq \z_y},$$
for some matrix $A$ on $\X \times \X$ with non-negative entries and zero entries on the diagonal.  Then  
$$  q(\G(\X))\leq \|A\|_{\infty} := \max_{x \in \X} \sum_{y \in \X}A_{x,y}.$$
The matrix $A$ is called `Dobrushin's interdependence matrix' and there are many instances  where an upper bound for $\|A\|_{\infty}$ can be derived for bounded degree graphs; for example see \cite[Section 8]{georgii2011gibbs}. The work \cite[Theorem 4.1]{kunsch1982decay} provides conditions on $\|A\|_{\infty}$ guaranteeing correlation decay bounds as well as a central limit theorem for spin systems (with values in a compact space) on lattices with variance being of volume order. In the case of a pairwise potential as in our work, the {\em Dobrushin uniqueness condition} $\|A\|_{\infty} < 1$ suffices.

We now present two specific examples next where this has been computed in \cite[Section 4.3]{chatterjee2005concentration} and in particular can be applied to the $k$-nearest neighbor graph.

\item {\em Ising model.} Here $\K = \{-1,+1\}$, $\Psi(u,v) = uv$ and $\Phi(v) = v$. The Gibbs measure is given by

$$ \pi_{\bk x1m}(\bk v1m)=\prod\limits_{(x_i,x_j) \in \G(\X) } \exp\{ \beta \sum_{(x_i,x_j) \in \G(\X)} v_iv_j  + \gamma \sum_{x_i \in \X} v_i \},
 $$
In this case $\|A\|_{\infty} \leq 4 \, \beta \, \textrm{deg}_{\max}$, where $\textrm{deg}_{\max}$ stands for the maximal degree. Thus in case of interaction graphs with degrees bounded almost surely, say by $\textrm{deg}_{\max}$, then Proposition~\ref{c:cltspindisag} applies for the Ising model with $\beta < (4\textrm{deg}_{\max}(\textrm{deg}_{\max} - 1))^{-1}.$  We have used that $\bar{\Delta} \leq \textrm{deg}_{\max} -1$.

\item {\em Proper $k$-colorings.} Here $\K = \{1\ldots,k\}$, $\Psi(v_i,v_j) = \log \1{v_i \neq v_j}$, $\beta = 1$ and $\gamma = 0$. The Gibbs measure is given by
$$ \pi_{\bk x1m}(\bk v1m)= \prod\limits_{(x_i,x_j) \in \G(\X) } \1{v_i \not= v_j}.$$
For graphs whose degrees are bounded by $\textrm{deg}_{\max}$, we again have $\|A\|_{\infty} \leq %\frac{\xout{\textrm{deg}_{\max}}}{\xout{k - \textrm{deg}_{\max}}}$ 
\textrm{deg}_{\max}/(k - \textrm{deg}_{\max})$. Thus in case of interaction graphs with degrees bounded almost surely (say by $\textrm{deg}_{\max}$) then we can verify that Proposition \ref{c:cltspindisag} applies for the proper $k$-colorings with $k > \textrm{deg}_{\max}^2$.
\end{enumerate}

\begin{proof}[Proof of Proposition~\ref{c:cltspindisag}] 
To deduce the result from Theorem~\ref{spinmain} and \eqref{eq:gendisagbd}, it suffices to show that the Gibbs model satisfies averaged weak spatial mixing as in Definition~\ref{d:ssm-pp} with some spatial mixing function \( \MC \), since all other assumptions are explicitly assumed and the moment condition follows trivially due to the finiteness of the spin space $\K$.
This involves verifying the inequality ~\eqref{spatmix} and showing the functions 
$(\ga_p(k))_{p \in \N}$ in~\eqref{e.spin-mixing-constant-bound}
are fast decreasing in $k$, both involving \( \MC \).
In what follows we will show that  this double condition is met in  our model for  the spatial mixing function 
\begin{equation}
\label{e:spamixfn_Nxd}
\MC_{\G(\X)}(\X'', \partial \X') := \sum_{x \in \X''} N(x, d_{\G(\X)}(x, \partial \X'); \G(\X))
a^{d_{\G(\X)}(x, \partial \X')}\,,
\end{equation}
where $a := 1/(\bar\Delta + \epsilon)$, with $\epsilon > 0$ fixed, provided $\partial\X'\not=\emptyset$, otherwise  we put $\MC_{\G(\X)}(\X'', \emptyset)=0$.
 We start by studying the functions $(\ga_p)_{p \in \N}$ at ~\eqref{e.spin-mixing-constant-bound} and, following Remark~\ref{r:expvarspin}\ref{i.particular-sets}, we assume $B\in\nicesets$. Observe that
\[
 \sE_{\bk x1p}[\MC_{\G(\P \cap B)}(\bk x1p, \partial \gB_k(\bk x1p))] \le  \sum_{i=1}^p \sE_{\bk x1p}[N(x_i, k+1, \G(\P \cap B))]a^{k+1} \leq p \, \gn_p(k+1)a^{k+1}.
\]
Given $\epsilon > 0$, by ~\eqref{e:bardelta}  we have  $\gn_p(k+1) \leq (\bar{\Delta}+ \epsilon/2)^{k+1}$ for $k\ge m_p$ and thus
$$\sE_{\bk x1p}[\MC_{\G(\P \cap B)}(\bk x1p, \partial \gB_k(\bk x1p))]
\le
p((\bar{\Delta}+\epsilon/2) a)^{k+1}=p\Bigl(\frac{\bar{\Delta}+\epsilon/2}{\bar{\Delta}+\epsilon}\Bigr)^{k+1}\,,$$
which shows that $(\ga_p(k))_{p \in \N}$ are fast decreasing in $k$,
as required. 

We show that ~\eqref{spatmix} holds with $\MC$ as in \eqref{e:spamixfn_Nxd} above by bounding the right-hand side of \eqref{eq:gendisagbd}  using standard percolation arguments. Consider $x \in \X''$ and denote by $\Gamma$ a path  from $x$ to $\partial \X'$ on $\G(\X)$, and let $|\Gamma|$  be the number  of sites in $\Gamma$. Observe that $|\Gamma| \geq d_{\G(\X)}(\X'',\partial \X') =: l$. 

For a given $q\in[0,1]$, in what follows, consider open paths in  site percolation on $\G(\X)$ in which vertices in~\(\X'\) are independently open with probability~\(q\), while all other vertices are always open. With probability $\sP_q = \sP_{\perc{\G(\X)}{\X'}{q}}$ of this model,  Markov's inequality and the above observations yield 
\begin{align}
& \sP_{q}(\text{there exists a self-avoiding (s-a) open path from $x$ to $\partial \X'$ })  \no \\
&  \leq \quad \sP_{q}(\text{there exists a  s-a open path of length $l$ from $x$}) \no \\
& \leq \sum_{\mbox{$\Gamma$:  s-a  path of length $l$ in $\X'$ from $x$}} \sP_{q}(\mbox{all sites  in $\Gamma$ are open}) \no \\
& \leq \, q^{l}N(x,l;\G(\X)), \label{e:percbdsaw}
\end{align}
where we recall that $N(x,l;\G(\X))$ is the number of self-avoiding walks in
$\G(\X)$ of length $l$ starting at $x$.
Note that a path of length~\(l\) consists of \(l+1\) vertices. Since  only the vertices in
\(\X'\) may be removed with probability $1-q$, while those in~\(\partial\X'\) are always open, it follows that the probability that the entire path is open is generously bounded by ~\(q^{\,l}\).

Assume now $q = q(\G(\X))$. The disagreement percolation bound \eqref{eq:gendisagbd},  the union bound and \eqref{e:percbdsaw} yield
\begin{align*}
& \quad d_{\TV} \big( {\cal L}({V}_{\X''}| {V}_{\partial \X'} = v_{\partial \X'}),
{\cal L}({V}_{\X''}| {V}_{\partial \X'} = \z_{\partial \X'}) \big) \\ & \leq \sum_{x \in \X''} \sP_{q}(\text{there exists a s-a open path from $x$ to $\partial \X'$}) \\
& \leq \sum_{x\in\X''}N(x,l;\G(\X))q^{l}\\
&\le  \sum_{x\in\X''} N(x,d_{\G(\X)}(\X'',\partial \X');\G(\X)) a^{d_{\G(\X)}(\X'',\partial \X')}
\\
&= \MC_{\G(\X)}(\X'',\partial \X')\,,
\end{align*}
where, the last inequality holds provided the inequality  $q(\G(\X))\le 1/(\epsilon+\bar\Delta)=:a$. This inequality   is assumed  almost surely 
for $\G(\X)=\G(\P\cap B)$ .
This shows that \eqref{spatmix} holds with $\MC$ as in \eqref{e:spamixfn_Nxd} and  completes the proof of Proposition~\ref{c:cltspindisag}.
\end{proof}

\subsubsection{Controlling spin correlations via disagreement percolation and sharp phase transition}
\label{s:spin_Poisson}

In this section, we refine the self-avoiding walk counting argument used in Section~\ref{s:gen_disagreement} along with disagreement percolation to establish averaged weak spatial mixing for the spin systems. Specifically, we  present an alternative argument for the exponential decay of disagreement probabilities, one exploiting the concept of a {\em sharp phase transition} in an auxiliary independent Bernoulli site percolation model. Essentially, if a sharp phase transition exists in this model, then by definition, for subcritical values \(q\) of the Bernoulli probability, {\em the $1$-arm probabilities} decay exponentially i.e., the probability of paths connecting a point \(x \in \X\) to points outside the ball \(B_r(x)\) of radius \(r\) decays exponentially with \(r\). (This property, combined with a mean-field lower bound for the percolation function at supercritical \(q\), formally defines the sharp phase transition in the model.) The sharp phase transition has been proved for many graph models, including the graph of the Boolean model; see \cite{duminil2019exponential,duminil2020subcritical,ziesche2018sharpness,meester1996continuum}.
In this section, we detail how this concept can be applied to establish the averaged weak spatial mixing property for these graphs. 

Since we are dealing with spatial random graphs, sharp phase transition results are often with respect to the Euclidean distance and hence we use a version of averaged weak spatial mixing with respect to the Euclidean distance to prove the same. As will be illustrated by the theorem below, this variant of averaged spatial mixing can be used in Theorem~\ref{spinmain} instead of condition~\ref{i.averaged-spatial-mixing} of Definition~\ref{d:ssm-pp}. For generality, we begin by defining `sharp' critical probability \(q_c\) in the context of disagreement percolation models on random graphs. Analogously to Definition~\ref{d:connective-constant-pp}, we define it as follows:  
\begin{defn}[Critical probability in disagreement percolation]\label{d:sharp-dperc}
For a graph $\G = {\G(\cdot)}$ defined on the input process \( \P \), the {\em critical probability $q_c$ for  percolation} is defined as  
\begin{equation}\label{e:qc-sharp}
q_c := \sup\Bigl\{q \in [0,1] :  \forall p \in \N , \, 
\text{the function $\gta_{p,q}(r)$ is fast decreasing as $r\to\infty$}\Bigr\},
\end{equation}
where for $r > 0$
\begin{equation}\label{e.gta-pqr} \gta_{p,q}(r):=\sup_{n \geq 1} \sup_{x_1,\ldots,x_p\in W_n}
\mathbb{P}_{\bk x1p,q}\bigl( x_1 \stackrel{(\P_n\cap B_{r}(x))_{*q}}{\longleftrightarrow} \P_n \cap B_r(x_1)^c \bigr), \,
\end{equation}
where the underlying graph is $\G(\P_n)$ and recall that
\(
\{x \stackrel{{\X'}_{*q}}{\longleftrightarrow} \mathcal{Y}\} = \{x \stackrel{\perc{\G(\X)}{\X'}{q}}{\longleftrightarrow} \mathcal{Y} \}
\)
means that there exists an open path in~\(\G(\X)\) from~\(x\in\X'\subset\X\) to
\(\mathcal{Y}\subset \X\)  through vertices in $\X'_{*q}$, the thinned subset of $\X'$ where vertices in~\(\X'\) are
independently open/retained with probability~\(q\), while all other vertices in $\X$ are
always open, and where \(\mathbb{P}_{\bk x1p,q}\) incorporates this additional percolation randomness on top of Palm distribution. \end{defn}

The next result, proved similarly to  Theorem \ref{spinmain}, proves normal convergence for spin systems on spatial random graphs whenever the maximal influence is smaller than the critical probability $q_c$.
\begin{proposition}[CLT for sum of spins for graphs with  sharp phase transition] 
\label{c:cltsharp}  
Let \(\P\) be a simple, stationary point process on \(\R^d\) with summable exponential mixing correlations as in Definition \ref{def.A2}, and assume that the (reduced) Palm intensity function of the ground process \(\P\) is bounded at all orders, as specified at \eqref{bdedreduced}. Let \(\G = \G(\cdot)\) be a stabilizing interaction graph on finite windows with respect to \(\P\) as per Definition~\ref{d:stabilizing-graph},
with the  critical probability $q_c$ as in Definition~\ref{d:sharp-dperc}.
Assume that for some $\epsilon>0$, the maximum influence \( q(\G(\P_n)) \) defined in~\eqref{defq} of the Gibbs specification~\eqref{eqn:gibbs} with finite spin space $\K$ satisfies the inequality
\be \label{maxinfbound}
q(\G(\P_n)) <q_c-\epsilon,
\ee
almost surely for all  $\P_n$, for all large enough $n\in\N$.  Then the system of spins \( \{V_x(\G(\P_n))\}_{x \in \P_n} \) satisfies the CLT as stated in Theorem~\ref{spinmain}, provided the variance condition there is also met.
\end{proposition}
Using bounds as in \eqref{e:percbdsaw} 
and the quenched version $\Delta$ of the connective constant 
(see the comment after Definition~\ref{d:connective-constant-pp})
one obtains \(\Delta^{-1} \leq q'_c\), where $q'_c$ is a (quenched) version of the sharp critical percolation in the graph $\G(\P)$ given a realization $\P$.
The strict inequality has been demonstrated in certain lattice cases where it is known that $q'_c > \Delta^{-1}$; see~\cite{wierman2022new,Duminil-Copin2012,sinclair2017spatial}. 
Although averaged connective constants or critical probabilities are not easy to compute for spatial random graphs models, we expect the strict inequality $q_c>\bar\Delta^{-1}$ to hold for some models and hence Proposition~\ref{c:cltsharp} would give a central limit theorem  covering a wider range of activity parameters than that given by Proposition~\ref{c:cltspindisag}. Moreover, as demonstrated in Section~\ref{s:spinPoissonGilbert}, we can still apply Proposition~\ref{c:cltsharp} for unbounded degree graphs without needing to compute bounds on the averaged connective constant.

\begin{proof}[Proof of Proposition~\ref{c:cltsharp}]
We exploit the critical probability threshold in disagreement percolation,
together with the maximum influence bounds of the Gibbs specification, to
establish Lemma~\ref{stablifting} directly, thereby bypassing the verification
of averaged weak spatial mixing as defined in
Definition~\ref{d:ssm-pp}. With this key lemma in place, the proof of
Theorem~\ref{spinmain} then follows by reusing the same core arguments.
The main difference between the proofs of~\eqref{Lp-spin-stab} and
\eqref{Lp-spin-stab2} in Lemma~\ref{stablifting} lies in the approximation
strategy for the spin vector
\(
\VBr
= V_{\bk x1p}\bigl(\G(\P_n \cap B_r(\bk x1p))\bigr).
\)
Instead of passing through the intermediate approximation~$\VBkr$, which
relies on graph-based neighborhoods of the points~$\bk x1p$ in
$\G(\P_n)$, we approximate~$\VBr$ directly by
\(
V_{\bk x1p}\bigl(\G(\P_n \cap B_{r/2}(\bk x1p))\bigr)
\)
using a Euclidean truncation of the input configuration. This direct
approximation is made possible by the disagreement-percolation bounds
described above and avoids the need for graph-neighborhood expansions.
The details are given below.

Fix $n \in \N$ and $\{x_1,\ldots, x_p\} \subset W_n$. We set $\X'' = \bk x1p$ and under $\mathbb{P}_{\bk x1p}$ consider the realization $\P_n$. Set $\X' = \P_n \cap B_{r/2}(\bk x1p)$ 
where we have suppressed the dependence on $r, \bk x1p$ in $\X'$ for convenience. 
Define $A_{r,n}(\bk x1p)$ as
the event (configurations of points of $\P_n$)  such  that for all $x \in \P_n \cap B_{r/2}(\bk x1p)$,  $S_n(x,\P_n) \leq r^b/4$ for a fixed $b \in (0,1)$. Recall that a similar event  was defined in the proof of  Lemma~\ref{stablifting} and we have used the same notation for convenience. We know that under the event $A_{r,n}(\bk x1p)$, the $r/2$ neighborhood of $\X'$ remains the same in $\P_n$ and $\P_n \cap B_r(\bk x1p)$ when $r \geq 1$.  Thus, we have that 
\begin{equation}
\label{eq:bdequiv}
\partial^{\G(\P_n)} \X' = \partial^{\G(\P_n \cap B_r(\bk x1p))} \X' = : \partial \X' \subset \P_n \cap B_{r/2}(\bk x1p)^c.
\end{equation}
Assume that $A_{r,n}(\bk x1p)$ holds. Following the proof of Lemma \ref{mixingimpliesdbl} verbatim with the above notation, for   $V_{\X''}=V_{\X''}(\G(\P_n))$ and  $V^{\cup(r)}_{\X''}=V_{\X''}(\G(\P_n \cap B_{r}(\bk x1p)))$, we have 
\begin{align*}
&d_{\TV} ( 
{\cal L}(V_{\X''} ),{\cal L}(V^{\cup(r)}_{\X''}\big))  \\
&\le    
  \sum_{ z_{\partial \X'},w_{\partial \X'}
\in \K^{\partial \X'}  }
\sum_{v_{\X''} \in \K^{\X''}} \Big|  \mP(V_{\X''} = v_{\X''} | {V}_{\partial \X'} = \z_{\partial \X'} ) 
   - \mP'(V_{\X''} = v_{\X''} | {V}^\cup(r)_{\partial \X'} = w_{\partial \X'})\Big| \\[-2ex]
&\hspace{20em}\times \mP({V}_{\partial \X'} = z_{\partial \X'})
\mP'(V^\cup(r)_{\partial \X'} = w_{\partial \X'}),
\end{align*}
with probabilities $\mP=\mP_{\G(\P_n)}$  and $\mP'=\mP_{\G(\P_n \cap B_{r}(\bk x1p)))}$ corresponding to  the Gibbs model  on the respective graphs. Here we have used that the boundary equivalence as in \eqref{eq:bdequiv} holds under $A_{r,n}(\bk x1p)$. Now the spatial Markov property from Lemma~\ref{l:spmarkspin}  yields the conditional equality under the  two probabilities $\mP$ and $\mP'$
$$
{\cal L}(V_{\X''} | {V}_{\partial \X'} = w_{\partial \X'}) = {\cal L}(V_{\X''}^{\cup{(r)}} | 
{V}_{\partial \X'} = w_{\partial \X'})
$$
and thus 
\begin{equation}\label{e.dtv-V-Vcup}
 d_{\TV} ( 
{\cal L}(V_{\bk x1p} ),{\cal L}(V^{\cup(r)}_{\bk x1p}\big)) \le 
\sup_{ z_{\partial \X'},w_{\partial \X'}}d_{\TV} \big( \mathcal{L}({V}_{\X''} \mid {V}_{\partial \X'} = z_{\partial \X'}),  
\mathcal{L}({V}_{\X''} \mid {V}_{\partial \X'} = w_{\partial \X'}) \big).
\end{equation}
To further bound the last expression, we consider site percolation on $\G(\P_n)$ where vertices in~\(\X'\) are independently open with probability $q\in[0,1]$ and where all other vertices are always open. Denote the probability measure of this model by $\mathbb{P}_q$ and omit $\G(\P_n)$ from the path notation $\longleftrightarrow$. Applying~\eqref{eq:gendisagbd} with retention probability $q=q(\P_n)$ and then increasing it to $q' := q_c - \varepsilon$, we obtain
\begin{align}
d_{\TV} \big( \mathcal{L}({V}_{\X''} \mid {V}_{\partial \X'} = z_{\partial \X'}),  
\mathcal{L}({V}_{\X''} \mid {V}_{\partial \X'} = w_{\partial \X'}) \big)  
& \le \mathbb{P}_{q(\G(\P_n))}\big( \X'' \stackrel{\X'_{*q(\G(\P_n))}}{\longleftrightarrow} \partial \X'\big)
\nonumber  \\
& \le \mathbb{P}_{q'}\big( \X'' \stackrel{\X'_{*{q'}}}{\longleftrightarrow} \partial \X'\big)
\nonumber  \\
& \leq \sum_{x \in \X''} \mathbb{P}_{{q'}}\big( x \stackrel{\X'_{*{q'}}}{\longleftrightarrow} \partial \X'\big)
\nonumber  \\
& \le \sum_{x \in \X''} \mathbb{P}_{ {q'}}\big( x \stackrel{(\P_n\cap B_{r/2}(x))_{*{q'}}}{\longleftrightarrow} \P_n\cap B_{r/2}(x)^c \big), \nonumber
\end{align} 
where the second inequality follows from the assumption~\eqref{maxinfbound} (almost surely for all $\P_n$ and $B=W_n$),
the third one  by the union bound, and the last 
by~\eqref{eq:bdequiv}, considering shorter paths and  considering  $q'$-thinning only on points in $\P_n\cap B_{r/2}(x)$.

Thus, using~\eqref{e.dtv-V-Vcup}, we have derived that under $A_{r,n}(\bk x1p)$
\begin{equation}
d_{\TV} ( 
{\cal L}(V_{\X''} ),{\cal L}(V^{\cup(r)}_{\X''}\big)) \leq  \sum_{x\in\X''}\mathbb{P}_{ {q'}}\big( x\stackrel{(\P_n\cap B_{r/2}(x))_{*{q'}}}{\longleftrightarrow} \P_n\cap B_{r/2}(x)^c \big).
\end{equation}
Using the above bound, averaging  with respect to the  input 
$\P_n$, and following the lines of the proof of Lemma \ref{stablifting}, we find
\begin{align*}
d_{BL,\bk x1p} \big({\cal L}(V_{\bk x1p} ),{\cal L}(V^{\cup(r)}_{\bk x1p})\big)  
& \leq  \E_{\bk x1p} \Big[ d_{BL} \big({\cal L}(V_{\bk x1p} \mid \P_n),{\cal L}(V^{\cup(r)}_{\bk x1p} \mid \P_n) \big) \Big] \\
& \leq  \E_{\bk x1p} \Big[ d_{\TV} \big({\cal L}(V_{\bk x1p} \mid \P_n),{\cal L}(V^{\cup(r)}_{\bk x1p} \mid \P_n) \big) \1{A_{r,n}(\bk x1p)} \Big]+  2 \sP_{\bk x1p}(A_{r,n}(\bk x1p)^c) \\
& = \sum_{i=1}^p \sP_{\bk x1p,{q'}}\big( x_i \stackrel{(\P_n\cap B_{r/2}(x_i))_{*{q'}}}{\longleftrightarrow} \P_n \cap B_{r/2}(x_i)^c \big) +  2 \sum_{i=1}^p \sP_{\bk x1p}(A_{r,n}(x_i)^c)\\
& \leq p \, \gta_{p,{q'}}(r/2) +  2 \sum_{i=1}^p \sP_{\bk x1p}(A_{r,n}(x_i)^c),
\end{align*}
where in the last inequality we used $\gta_{p,{q'}}(r)$ function defined in~\eqref{e.gta-pqr} and  the fact that 
%$A_{r,n}(\bk x1p)^c \subset \bigcup_{i=1}^pA_{r,n}(x_i)\BBLp{{}^c}$. 
$A_{r,n}(x_i)\subset A_{r,n}(\bk x1p)$ for $i=1,\ldots,p$.
By Definition~\ref{d:sharp-dperc} of $q_c$, the fast decay of the first term follows by the assumption of ${q'} = q_c - \epsilon$.

The fast decay of the second term can be shown by observing for all $i \in \{1,\ldots,p\}$, similarly as in~\eqref{e:adm_graph_euc_dist},
$$\sP_{\bk x1p}(A_{r,n}(x_i)^c) \leq (\theta_d \hat{\kappa}_p (\frac{r}{2})^d+p)  \varphi_{p+1}'(r^b/4).$$
Thus, we have established the fast \BL~union-localization bound~\eqref{Lp-spin-stab}.

We now  prove~\eqref{Lp-spin-stab2}. Following closely the proof of
Lemma~\ref{stablifting}, we fix $p \ge 2$, choose
$q \in \{1,\ldots,p-1\}$, and consider points
$\bk x1p \in W_n$.
Let $r>0$ satisfy
\(
0 < r < \BBLp{d\xout{s}}\bk x1l, \bk x{l+1}p),
\)
and, under the Palm distribution $\sP_{\bk x1p}$, consider the spin vectors
$\VBr$, $\VBlr$, and $\VBpr$ defined as before, with $\VBlr$, and $\VBpr$ assumed to be conditionally independent given $\P_n$. 

On the event $A_{r,n}(\bk x1p)$, the vectors $\VBlr$ and $\VBpr$
are marginally approximated by their counterparts obtained on the graph 
$\G(\P_n \cap B_{r/2}(\bk x1p))$ i.e., by $V_{\bk x1l}(\G(\P_n \cap B_{r/2}(\bk x1l)))$ and $V_{\bk x{l+1}p}(\G(\P_n \cap B_{r/2}(\bk x{l+1}p)))$, respectively.
The same approximation holds for $\VBr$.
Moreover, the approximation of $\VBr$ coincides with the product of the
corresponding approximations of $\VBlr$ and $\VBpr$ i.e.,
$$
\mathcal L\left(V_{\bk x1p}(\G(\P_n \cap B_{r/2}(\bk x1p))) \right) = 
\mathcal L\left(V_{\bk x1l}(\G(\P_n \cap B_{r/2}(\bk x1l))) \right) \times\mathcal L\left(V_{\bk x{l+1}p}(\G(\P_n \cap B_{r/2}(\bk x{l+1}p))) \right).
$$
Using the fact that the probability of the complement of
$A_{r,n}(\bk x1p)$ decays fast in $r$, together with the fast decreasing of the error in the above approximations
of $\VBr$, $\VBlr$, and $\VBpr$, we conclude that
\[
d_{\BL,\bk x1p}
\bigl(
\mathcal L(\VBr),
\mathcal L(\VBlr)\times\mathcal L(\VBpr)
\bigr)
\]
is fast decreasing in~$r$.
This establishes~\eqref{Lp-spin-stab2}.

Finally, applying Theorem~\ref{t:cltmarkedpp} as in the proof of
Theorem~\ref{spinmain}, we conclude the proof.
\end{proof}

\subsubsection{Spin systems on the Poisson-Boolean model}
\label{s:spinPoissonGilbert}

 We consider  an example of a general spin model on the graph generated by a Boolean model, where one can identify the critical probability for disagreement percolation $q_c$  and also furnish an explicit bound 
 for the maximal influence $q(\G(\P \cap B))$ of the spin system 
 in \eqref{maxinfbound}. 

Let \(\P\) be a Poisson point process of intensity \(\rho\) on \(\R^d\), $K \subset \R^d$ a compact set with a non-empty interior and without loss of generality, we assume that $\0 \in K$. 
The Poisson-Boolean model is $\cO = \cO(\P) := \bigcup_{x \in \P} (x \oplus K)$.
For $B\subset \R^d$, the Poisson-Boolean graph \(\G(\P \cap B)\) has vertex set \(\P\cap B\), with vertices $x, y$  connected by an edge if \((x \oplus K) \cap (y \oplus K) \neq \emptyset\). 

The critical intensity for percolation in the Boolean model, which determines whether an infinite connected component exists or not in the graph defined above, is denoted by \(\rho_c\). In other words, \(\rho_c> 0\) denotes the critical intensity below which an infinite connected component does not exist in the Boolean model. The non-triviality of  $\rho_c$ (i.e.,  \(0 < \rho_c < \infty\)) in dimensions \(d \geq 2\) as well as a sharp phase-transition is known; see \citet{hall1985continuum,meester1996continuum}.

In this section, we exploit the existence of a sub-critical regime where the graph \(\G(\P)\) `sharply' fails to percolate when \(\rho \in (0, \rho_c)\). 
Even if the actual intensity \(\rho\) is arbitrarily large, we thin the Poisson process with parameter \(q < \rho_c / \rho\), 
yielding $\P_{*q}$ (a Poisson point process  of intensity $q \rho$),
placing our auxiliary model of disagreement percolation in the subcritical regime. By leveraging the exponential decay of the probability of long paths in this regime, which arises from the sharp phase transition in the Boolean model proved in~\cite{ziesche2018sharpness}, we identify the value \(\max\{1, \rho_c / \rho\}\) as an upper bound for the critical probability $q_c$ i.e., \(q_c \geq \max\{1, \rho_c / \rho\}\).

Consequently, although the graph has unbounded degree and  we cannot establish bounds for the averaged connective constant in~\eqref{e.connective-constant-pp}, but we are able to apply Proposition~\ref{c:cltsharp}. This allows us to prove a central limit theorem  for the Gibbs specification~\eqref{eqn:gibbs}, provided the maximum influence remains below \(q_c\).
Importantly, as observed in models like the hard-core model or the Widom-Rowlinson model, this result holds across the entire parameter range $\rho$ of the Boolean model including the percolating regime (\(\rho > \rho_c\)), i.e., $q_c > 0$ for all $\rho > 0$. Furthermore, in the subcritical regime of the Boolean model (\(\rho < \rho_c\)), we find that \(q_c = 1\), meaning that any activity level for the Gibbs parameters is admissible, thereby extending the applicability of the central limit theorem without imposing additional constraints. 
\begin{corollary}
    \label{t:cltspindisag1} 
(CLT for sum of spins on Poisson-Boolean model)
Let \(\P\) be a stationary Poisson process on $\R^d$ of intensity \(\rho\) and $K$ be a compact set containing the origin and with a non-empty interior. Let  \(\G = \G(\cdot)\) be the Poisson-Boolean graph induced by the Poisson-Boolean model on \(\P\), as described above. Consider a spin configuration model with the Gibbs distribution given by \eqref{eqn:gibbs}, defined on the graph~\(\G\) and with finite spin space $\K$.
Assume that for some \(\epsilon > 0\), the maximum influence \eqref{defq} satisfies
\[
q(\G(\P_n)) < \min \{1, \rho_c / \rho\} - \epsilon
\]
almost surely for all $\P_n$, for all large enough $n\in\N$, and where $\rho_c$ is the critical intensity for percolation of the Poisson-Boolean model.
Under these conditions, the spin model \(\{V_x(\G(\P_n))\}_{x \in \P_n}\) satisfies the central limit theorem as stated in Theorem~\ref{spinmain}, provided that the moment and variance conditions specified in the theorem are also met. 
\end{corollary}%
One can consider a more general Boolean model with random grains, i.e., the fixed compact set $K$ is replaced by random `shapes' $Z_x, x \in \P$. In case $Z_x$'s are bounded, our approach can be adapted by extending the framework of stabilizing graphs from Section~\ref{sss.graph-stabilization} to allow for graphs constructed on marked point processes; see Remarks in Section \ref{s:future}. In particular, the percolation-theoretic results that we shall use are already available in this case; see \cite{ziesche2018sharpness}. The case of unbounded grains requires more careful investigation and in particular suitable moment assumptions on the grains may need to be imposed; see \cite{duminil2020subcritical}. It may be possible to consider Poisson-Delaunay graphs, leveraging the results of \cite{duminil2019exponential} on site percolation in Poisson-Delaunay graphs.
\begin{proof}[Proof of Corollary \ref{t:cltspindisag1}]
We will apply Proposition~\ref{c:cltsharp} showing that \(q_c\ge \min \{1, \rho_c /\rho\}\), where $q_c$  is defined in Definition~\ref{d:sharp-dperc}.~\footnote{In fact, there is equality: \(q_c = \min\{1, \rho_c / \rho\}\). Indeed, in the subcritical regime (\(\rho < \rho_c\)), the arguments developed below apply for \(q = 1\), resulting in \(q_c = 1\). Conversely, in the supercritical regime (\(\rho > \rho_c\)), any \(q\)-thinning of the Boolean model with \(q > \rho_c / \rho\) keeps it in the supercritical regime, causing the functions \(\gta_{p,q}(m)\) in~\eqref{e:qc-sharp} to converge to a positive value.}   
A key result (see, for example, \cite[Theorem~1.1; see also Theorem~3.1]{ziesche2018sharpness})
establishes an exponential bound on the probability that, in the subcritical
Boolean model (\(\rho < \rho_c\)), there exists a \emph{continuum path}  in  $\cO$ from the origin to the exterior of the Euclidean ball \(B_r(\0)\): 
\begin{equation}\label{e.Ziesche}  
\mathbb{P}_{\0} \left( \0\stackrel{\cO}{\longleftrightarrow}  B_r(\0)^{c} \right) \leq e^{-cr} ,
\end{equation}  
for some constant \( c = c(\rho) \).

This result yields the exponential decay of the functions \(\gta_{p,q}(r)\) in ~\eqref{e.gta-pqr}
for $q<q_c$, seen as follows.
Specifically, we recall the notation from Definition~\ref{d:sharp-dperc}: for
\(\X' \subset \X\),
\(
\{ x \stackrel{\X'_{*q}}{\longleftrightarrow} \mathcal{Y} \} = \{ x \stackrel{\perc{\G(\X)}{\X'}{q}}{\longleftrightarrow} \mathcal{Y} \}
\)
means that there exists an open path in the graph \(\G(\X)\) (here induced by the Boolean model) connecting
\(x \in \X'\) to \(\mathcal{Y} \subset \X\) under site percolation, where
vertices in \(\X'\) are independently declared open with probability
\(q\), while all other vertices in \(\X\) are always open.
The probability measure \(\mathbb{P}_{\bk{x}{1}{p},q}\) incorporates
this additional percolation randomness on top of Palm distribution. 

With this notation and assuming also without loss of generality that $K \subset B_1(\0)$, for \(n \in \N\), \(\bk{x}{1}{p} \in W_n\), and \(r > 0\), we have
\begin{equation}\label{e.expm0}
\mathbb{P}_{\bk{x}{1}{p},q}\!\left(
x_1 \stackrel{(\P_n \cap B_r(x_1))_{*q}}{\longleftrightarrow}
\P_n \cap B_r(x_1)^c
\right)
\le
\mathbb{P}_{\bk{x}{1}{p},q}\!\left(
x_1 \stackrel{\cO(\P_{*q})}{\longleftrightarrow}
B_{\max(0,r-2)}(x_1)^c
\right).
\end{equation}
Indeed,  applying independent
\(q\)-thinning to the original Boolean model of intensity~\(\rho\) is same as considering a Boolean model of intensity $q\rho$. The inequality follows from the fact that we allow continuum paths to reach the
entire complement of \(B_{\max(0,r-2)}(x_1)\), rather than restricting them to grains of points in
\(\P_n \cap B_r(x_1)^c\); this relaxation can only increase the probability.
Moreover, thinning of \(\P_n\) outside \(B_r(x_1)\) is immaterial after this
relaxation.
Finally, considering a Boolean model outside  $W_n$ can only further increase the
upper bound for a given $n$, in case  \(B_r(x_1)\cap W_n^c\not=\emptyset\).

Note that $\P_{*q}$ is distributed as a Poisson process $\P_{\circ}$ of intensity $q \rho$. Thus, we have that
$$ \mathbb{P}_{\bk{x}{1}{p},q}\!\left(
x_1 \stackrel{\cO(\P_{*q})}{\longleftrightarrow}
B_{r}(x_1)^c
\right) = 
\mathbb{P}_{\bk{x}{1}{p}}\!\left(
x_1 \stackrel{\cO(\P_{\circ})}{\longleftrightarrow}
B_{r}(x_1)^c
\right).$$
We  assert that, uniformly over
\(\bk{x}{1}{p}\subset\R^d\)
\begin{equation}\label{e.expm}
%\mathbb{P}_{\bk{x}{1}{p},q}\!\left(x_1 \stackrel{\G(\P_n)\cap (\P_n \cap B_r(x_1))_q}{\longleftrightarrow}\P \cap B_r(x_1)^c\right)
\mathbb{P}_{\bk{x}{1}{p}}\!\left(
x_1 \stackrel{\cO(\P_{\circ})}{\longleftrightarrow}
B_{r}(x_1)^c
\right)
\le e^{-c r},
\end{equation}
for some constant \(c>0\) and all sufficiently large \(r>0\), provided the Boolean model is thinned with probability \(q\) such that \(q < \rho_c/\rho\), where $\rho_c$ is the critical intensity for percolation of this model.

To derive the bound \eqref{e.expm} from~\eqref{e.Ziesche}, we need to a control the impact of several fixed atoms of the input process induced by  higher-order Palm probabilities, which may result in looser bounds compared to the analysis in~\cite{ziesche2018sharpness}, where Palm conditioning is applied at a single point placed at the origin. 

Fix $r > 4p$. By the Slivnyak-Mecke theorem, (see e.g.~\cite[Theorem 3.2.4]{BBK2020}) we have that
\[
\mathbb{P}_{\bk{x}{1}{p}}\!\left(
x_1 \stackrel{\cO(\P_{\circ})}{\longleftrightarrow}
B_r(x_1)^c
\right)=
\mathbb{P}\!\left(
x_1 \stackrel{\cO(\P_{\circ}+ \sum_{i=1}^p \delta_{x_i})}{\longleftrightarrow}
B_r(x_1)^c
\right)
\]
where \( \stackrel{\cO(\P_{\circ}+ \sum_{i=1}^p \delta_{x_i})}{\longleftrightarrow} \) denotes the existence of a path in the Boolean model augmented with the fixed points \(x_i\) and grains $K$ at these sites i.e., $\cO(\P_{\circ}+ \sum_{i=1}^p \delta_{x_i}) := \bigcup_{x \in \P_{\circ} \cup \{x_1,\ldots,x_p\}} (x \oplus K)$.

Under the assumption \( q < \min \{1, \rho_c / \rho\} \), the model \(\P_{\circ}\) is in its subcritical regime.
To leverage equation~\eqref{e.Ziesche}, we must address the `undesired' extra points arising from \(\sum_{i=1}^p \delta_{x_i}\). To do this, we subdivide \(B_{r}(x_1)\) into annuli \(A_1, \ldots, A_{10p}\) of width \(\frac{r}{10p}\), where:  
\[
A_1 := B_{r / (10p)}(x_1), \quad A_i := B_{i r / (10p)}(x_1) \setminus B_{(i-1) r / (10p)}(x_1), \quad i = 2, \ldots, 10p.
\]  
By the pigeonhole principle, the choice of \(r\), and the set inclusion $K \subset B_1(\0)$, for sufficiently large \(r\) (precisely when \(2 < \frac{5 r}{10p}\)), there exist five consecutive annuli, say \(A_j, \ldots, A_{j+4}\) for some \(j \in \{1, \ldots, 10p-4\}\), such that  
\[
\bigcup_{i=1}^p B_{1}(x_i) \cap (A_j \cup \ldots \cup A_{j+4}) = \emptyset.
\]
As a result, these five annuli remain unaffected by the extra points and their grains in \(\sum_{i=1}^p \delta_{x_i}\).
If we further increase \(r\) such that \(2 < \frac{r}{10p}\), then any path from \(x_1\) to \(B_{r}(x_1)^c\) must pass through a Poisson point \(x \in \P_{\circ} \cap A_{j+2}\) since the grain $K$ is within a unit ball. Consequently, we have:
\begin{align*}
&\mathbb{P} \left( x_1 \stackrel{\cO(\P_{\circ}+ \sum_{i=1}^p \delta_{x_i})}{\longleftrightarrow} B_{r}(x_1)^c \right) \\
&\leq \sum_{j=3}^{10p-2} \mathbb{P} \left(\text{there exists } x \in \P_* \cap A_j \text{ such that } x\stackrel{\cO(\P_{\circ})}{\longleftrightarrow} B_{r / (10p)}(x)^c \,\right).
\end{align*}
For a fixed \(j = 3, \ldots, 10p-2\), we now bound the probability using the Campbell-Little-Mecke formula~\eqref{CLM} and the bound~\eqref{e.Ziesche}:
\begin{align*}
 &\hspace{-5em}\sP\left(\text{there exists $x\in\P_*\cap A_{j+2}$ such that  }
x\stackrel{\cO(\P_{\circ})}{\longleftrightarrow} B_{r/(10p)}(x)^c\,\right)\\
 & \leq \sE\Bigl[\sum_{x \in \P_* \cap A_{j+2}}\1{x\stackrel{\cO(\P_{\circ})}{\longleftrightarrow}B_{r/(10p)}(x)^c}\Bigr]\\
 & = q\rho \int_{A_j} \sP_x\left( x\stackrel{\cO(\P_{\circ})}{\longleftrightarrow}B_{r/(10p)}(x)^c\,\right)\, \md x \\
 & \leq  q \rho \, \Vol(A_j) \, e^{-cr/(10p)},
\end{align*}
where in the last inequality, we used stationarity and equation~\eqref{e.Ziesche} with \(\rho\) replaced by \(q\rho < \rho_c\). Here, \( \Vol(A_j)\) denotes the volume of the annulus \(A_j\).

Summing over \(j = 3, \ldots, 10p-2\) and combining with the above derivations, we obtain the necessary exponential bound for the left hand side in~\eqref{e.expm}:
\[
q\rho \, \Vol(B_{r}) \, e^{-c r / (10p)} = q\rho \theta_d \, r^d e^{-c r / (10p)} = Cq\rho e^{-c r / (10p)},
\]
for some constant $C$ depending on $p$.
Accounting for the maximal grain diameter by shifting the argument \(r\) by
\(2 \) in the exponential function above yields the
exponential bound for~\eqref{e.expm0}.

This establishes  the  exponential decay of the functions \(\gta_{p,q}(r)\) in~\eqref{e.gta-pqr} and thus the corollary as argued before \eqref{e.expm}.
\end{proof}

\subsection{Strong spatial mixing for the hard-core model via combinatorial techniques}
\label{s:hard-core}
Recall that the hard-core model is specified by the probability distribution
\begin{equation}
 \label{eqn:re-hardcore}   
\pi_{\bk x1m}(\bk v1m) = \prod\limits_{(x_i,x_j) \in \G(\X) } \1{v_i v_j = 0} \prod\limits_{x_i \in \X} \lambda^{v_i},  \, v_i \in \{0,1\} \, ,
\end{equation}
where \(\lambda > 0\) is the {\em activity parameter}, $\G(\X)$ is the interaction graph on $\X = \{x_1,\ldots,x_m\}$. In this section, we revisit this model using the combinatorial techniques of \citet{weitz2005combinatorial} and \citet{sinclair2017spatial}. Their analysis establishes a stronger form of spatial mixing, extending the range of low-activity parameters for the Gibbs hard-core model beyond the   maximum influence conditions typically employed to demonstrate weak spatial mixing via disagreement percolation. This is achieved by making use of a stronger spatial mixing property established for this Gibbs model in the aforementioned works. Specifically, this approach demonstrates the strong mixing bound~\eqref{ssmsingle},
 at the level of individual sites \(\X'' = \{x\}\). Using the additive representation provided by Lemma~\ref{p:ssmsingle-wsm} we revisit  the proof of  Lemma~\ref{stablifting}, a key argument in  the main result of Theorem~\ref{spinmain}, and which is based on  averaged weak spatial mixing.
Consequently, we  consider the range of the activity parameter \( \lambda \)  above the critical value 
\begin{equation}\label{e.hc-critical-lambda}
 \lambda_c(\Delta) := \frac{\Delta^\Delta}{(\Delta - 1)^{\Delta + 1}} > \frac{1}{\Delta - 1},\quad \text{for $\Delta\ge1$ and} \, \lambda_c(0) = \infty, 
\end{equation}
where the {\em (quenched) connective constant} \(\Delta \in [0,\infty]\) is defined as the smallest constant such that for all $m \geq m_0 \in \N$ and for all $x \in \X$, it holds that
\be \label{connconstant}
N(x, m; \G(\X)) \leq \Delta^m,
\ee
with \(N(x, m; \G(\X))\) denoting the number of self-avoiding paths of length \(m\) starting from \(x \in \X\) in the graph \(\G(X)\); see \cite[Section~2.6]{sinclair2017spatial}. However working with this quenched connective constant comes at the cost of requiring almost sure (quenched) control over the graph's connectivity.  Observe that the quenched connective constant cannot take fractional values in the range $\Delta \in (0,1)$. Additionally, a value of $\Delta = 0$ indicates that the graph consists of uniformly bounded components and one can take large activity parameters as well, which is consistent with $\lambda_c(0) = \infty$. The combinatorial techniques can also be applied to Monomer-Dimer model and Ising model with zero external field and hence our results also can be extended to these models; see the discussion at the beginning of Subsection~\ref{r.Hard-core-Sinclair}. The critical value $1/{\bar{\Delta}}$ established in Proposition~\ref{c:cltspindisag} may not be comparable to $\lambda_c(\Delta)$ as $\bar{\Delta} \leq \Delta$.
 
For deterministic graphs, the critical value \(\lambda_c(\Delta)\) was first identified for hard-core models on infinite \(d\)-ary trees with \(\Delta = d\), as it characterizes the regime where correlations decay exponentially and consequently enabled efficient algorithms for sampling and computation of the partition function.
The work ~\cite{sinclair2017spatial} further generalized it to all finite or infinite graphs with a connective constant bounded by \(\Delta\). 
These works, however, do not address limit theorems for
statistics of random.

\begin{proposition}[CLT for the sum of spins in the hard-core model]
\label{c:clthardcore}
 Let \(\P\) be a  stationary process on \(\R^d\) with summable exponential mixing correlations as in Definition \ref{def.A2} and having bounded Palm intensity functions as at ~\eqref{bdedreduced}. Let \(\G = \G(\sim, \cdot)\) be a stabilizing interaction graph on \(\P\)
 and  satisfying the following a.s. connective constant bound: there exists $\Delta  \in [1,\infty)$ and $m_0 \in \N$ such that for all sets $B\in\nicesets$ with defined in~\eqref{e.particular-sets}, $x\in\P\cap B$, and $m\ge m_0$ 
\begin{equation}\label{e.Delta-as-bound}
N(x, m; \G(\P\cap B)) \le\Delta^m.
\end{equation}
Consider the hard-core Gibbs model~\eqref{eqn:re-hardcore} with activity parameter \(\lambda\) such that 
\[
\lambda < \lambda_c(\Delta), 
\]
where the critical value \(\lambda_c(\Delta)\) is at ~\eqref{e.hc-critical-lambda}. Under these conditions, the spin configuration \\
\( \{V_x(\G(\P_n))\}_{x \in \P_n} \) satisfies the CLT as stated in Theorem~\ref{spinmain}, provided that the %\dy{\xout{moment and} 
variance condition there is also satisfied.
\end{proposition}
\begin{proof}
We cannot directly use Theorem \ref{spinmain} but rather follow the approach therein by establishing Lemma \ref{stablifting} by suitably adapting the proof by borrowing the bounds from \cite{sinclair2017spatial}.

The case of $\Delta = 0$ is straightforward, as the central limit theorem follows directly from the presence of uniformly bounded components in the graph. For $\Delta \geq 1$, the proof hinges on the single-site strong spatial mixing property of the hard-core model on graphs, as demonstrated in~\cite{sinclair2017spatial}. This property relies on the `self-avoiding walk representation' introduced by~\cite{Weitz2006}, which simplifies the problem by reducing it to trees. The (quenched) spatial mixing property is established and analyzed through the study of occupation ratios. Specifically,  
in this model, the distribution of the spin at a given site can be described by the occupation probability:  
\[
p(x, \X; \v_{\cZ}) := \mathbb{P}(V_x = 1 \mid V_{\cZ} = v_{\cZ}),
\]  
which represents the probability that the spin at node \( x \in \X \setminus \cZ \) in the graph \( \G(\X) \) equals 1, conditioned on the spins at the subset \( \cZ \not\owns x \) being fixed to the values \( \v_{\cZ} \).
We analyze the decay of the dependence of the spin \( V_x \) on the boundary condition \( V_{\cZ} = v_{\cZ} \) through the  {\em occupation ratio} at node \( x \), namely  
\[
R(x, \X; \v_{\cZ}) := \frac{p(x, \X; \v_{\cZ})}{1 - p(x, \X; \v_{\cZ})}.
\]  
The total variation distance between the conditional distributions of \( V_x \), given two different boundary conditions, is bounded as follows:  
\begin{align}
d_{\TV} \big( \mathcal{L}(V_x \mid V_{\cZ} = \v_{\cZ}), \mathcal{L}(V_x \mid V_{\cZ} = \z_{\cZ}) \big) 
&= |p(x, \X; \v_{\cZ}) - p(x, \X; \z_{\cZ})| \nonumber \\
&\leq |R(x, \X; \v_{\cZ}) - R(x, \X; \z_{\cZ})|. \label{e.tv-occupation-ratio}
\end{align}  
Thus, a `spatial mixing property of the occupation ratios'  implies spatial mixing for the model.

In the hard-core model~\eqref{eqn:re-hardcore} on a general graph \(\G(\X)\), the strong spatial mixing of occupation ratios'  as established in~\cite{sinclair2017spatial}, may be stated as follows:   For any vertex  \(x \in \X\) and any subset $\cZ \subset \X$ not containing $x$, it is the case that there exists $\alpha < 1/\Delta$ such that
\begin{equation}\label{e.smor}
| R(x, \X; v_{\cZ}) - R(x, \X; z_{\cZ}) | 
\leq C N(x,l-1;\G(\X))^{\frac{1}{q}} \alpha^{\frac{l-1}{q}},
\end{equation}
where $l =  d_{\G(\X)}(x, \cZ_{\neq})$ is the distance between \(x\) and points in \(\cZ\) where the boundary conditions \(v_{\cZ}\) and \(z_{\cZ}\) differ,
and where it is required that \(l\)  is non-zero. 
The constants $C,q$ depend on $\lambda$ and are specified explicitly in the sequel. The decay constant
$\alpha$ is defined to be
\(\alpha := 1 / \Delta_c(\lambda)\), where \(\Delta_c(\lambda) \) is determined as the unique solution to the equation:  
\begin{equation}\label{e.lambda-solution-d}
\lambda = \frac{t^t}{(t - 1)^{t + 1}}.
\end{equation}  
Observe that the function \(t \mapsto t^t/(t - 1)^{t + 1}\) in~\eqref{e.lambda-solution-d} is decreasing. Recall that \(\Delta \in [1, \infty)\) is a bound on the connective constant, as in~\eqref{e.Delta-as-bound} and \(\lambda_c(\Delta)\) is the critical value as in~\eqref{e.hc-critical-lambda}. Thus, our assumption \(\lambda < \lambda_c(\Delta)\) implies that \(\Delta_c(\lambda)\), the solution in \(t\) of the equation~\eqref{e.lambda-solution-d}, satisfies \(\Delta_c(\lambda) > \Delta \geq 1\). Consequently, \(\alpha := 1 / \Delta_c(\lambda) < 1 / \Delta\). 

This formulation~\eqref{e.smor} provides a single-site strong spatial mixing bound based on occupation ratios and as demonstrated in~\eqref{e.tv-occupation-ratio}, it can be used in further considerations. The key ideas and steps behind the derivation of ~\eqref{e.smor} are  in Section \ref{r.Hard-core-Sinclair}.

We recognize that \(\alpha\) serves as a refinement of the requirement that the maximum influence \(q = q(\G(\X)) \) \footnote{The maximum influence $q = q(\mathcal{G}(\mathcal{X}))$, defined in~\eqref{defq} and not used any more in this proof, should not be confused with the constant $q$ in~\eqref{e.smor}. Both depend on the activity $\lambda$ of the hard-core model, but the former also depends on the underlying graph.}
 in Proposition~\ref{c:cltspindisag} is smaller than the (albeit averaged) connective constant~\(\bar\Delta\).
 Although multi-site weak spatial mixing can be derived (see \eqref{spamixa-strong}) from the strong spatial mixing single-site bound~\eqref{e.smor}, the resulting estimate would be too crude to establish \BL-stabilization \eqref{Lp-spin-stab}, crucial to proving  Theorem~\ref{spinmain}.  Instead, we shall establish  \eqref{Lp-spin-stab} by combining the multi-site bound \eqref{spamixa-strong}  with the additive representation for strong spatial mixing provided in Lemma~\ref{p:ssmsingle-wsm}, using a stronger, quenched (almost sure) bound on the connective constant \(\Delta\).

To prove \BL-stabilization~\eqref{Lp-spin-stab} under the new assumptions, we only need to refine the inequalities~\eqref{e.BL-spin-stab-1*} and~\eqref{e.BL-spin-stab-2*} in the proof of Lemma~\ref{stablifting}. Using this strong spatial mixing bound~\eqref{spamixa-strong} and letting $\X' = \gB_{k_r}(\bk x1p)$, we have: 
\begin{align*} 
&\E_{\bk x1p} \left[d_{\BL}(\mathcal{L}(\VBkr\mid \P_n), \mathcal{L}(\VP \mid \P_n))\right] \\
&\le \E_{\bk x1p} \Bigl[ \sum_{ z_{\partial \X'},w_{\partial \X'} \in \K^{\partial \X'}  } 
\MC'(\bk x1p,\partial_{z\not=w}\X')
\mP_{\G(\P_n)}({V}_{\partial \X'} = z_{\partial \X'})
\mP_{\G(\P_n)}(V^{[k_r]}_{\partial\X'} = w_{\partial \X'})\Bigr]\,.
\end{align*}
Next, applying the additive representation~\eqref{ssmsingle-mutiple} for strong spatial mixing along with~\eqref{e.smor} (linked to ~\eqref{e.tv-occupation-ratio}) and taking the function 
$$
\MC'_{\G(\P_n)}(\{x\}, \cZ_{\neq}) = C N(x, d_{\G(\P_n)}(x, \cZ_{\neq})-1;\G(\P_n))^{1/q} 
\alpha^{(d_{\G(\P_n)}(x, \cZ_{\neq})-1)/q},
$$ 
with constant $C$ as in \eqref{e.smor}, this  gives:  
\begin{align}\label{e.dBL-strong}
&\E_{\bk x1p} \left[d_{\BL}(\mathcal{L}(\VBkr\mid \P_n), \mathcal{L}(\VP \mid \P_n))\right] \\
&\le C \, \E_{\bk x1p} \Bigl[ \sum_{ z_{\partial \X'},w_{\partial \X'} \in \K^{\partial \X'}  } \sum_{i=1}^p
N(x_i, d_{\G(\P_n)}(x_i, \partial_{z\not=w} \X')-1;\G(\P_n))^{1/q} 
\alpha^{(d_{\G(\P_n)}(x_i, \partial_{z\not=w}\X')-1)/q} \no \\[-2ex]
&\hspace{20em}\times \mP_{\G(\P_n)}({V}_{\partial \X'} = z_{\partial \X'})
\mP_{\G(\P_n)}(V^{[k_r]}_{\partial\X'} = w_{\partial \X'})\Bigr]\,. \no 
\end{align}
Since \(\X' = \gB_{k_r}(\bk x1p)\) is the subgraph induced in~\(\P_n\) by the graph balls centered at \(x_i\), \(i=1,\ldots,p\), of radius \(k_r\), we have \(d_{\G(\P_n)}(x_i, \partial_{z\not=w} \X') > k_r\) (interpreted as $=\infty$ if $\partial_{z\not=w} \X'=\emptyset$), for all \(z_{\partial \X'}, w_{\partial \X'} \in \K^{\partial \X'}\) and \(i=1,\ldots,p\). By assumption~\eqref{e.Delta-as-bound},  for \(k_r > m_0\),
\[
N(x_i, d_{\G(\P_n)}(x_i, \partial_{\not=} \X')-1; \G(\P_n)) \leq \Delta^{d_{\G(\P_n)}(x_i, \partial_{\not=} \X')-1},
\]
 Substituting this into~\eqref{e.dBL-strong}, we obtain:  
\begin{align*}
&\E_{\bk x1p} \left[d_{\BL}(\mathcal{L}(\VBkr\mid \P_n), \mathcal{L}(\VP \mid \P_n))\right] \\
&\le C \, \E_{\bk x1p} \Bigl[ \sum_{ z_{\partial \X'},w_{\partial \X'} \in \K^{\partial \X'}  } \sum_{i=1}^p
(\Delta\alpha)^{(d_{\G(\P_n)}(x_i, \partial_{\not=}\X')-1)/q} \, 
\mP_{\G(\P_n)}({V}_{\partial \X'} = z_{\partial \X'})
\mP_{\G(\P_n)}(V^{[k_r]}_{\partial\X'} = w_{\partial \X'})\Bigr]\\
&\le C \, p(\Delta\alpha)^{(k_r)/q},
\end{align*}
where, in the last inequality, we used  \(d_{\G(\P_n)}(x_i, \partial_{\not=} \X') \ge d_{\G(\P_n)}(x_i, \partial \X') > k_r\) and the fact the $\Delta\alpha<1$
(as a consequence of our assumption  \(\lambda < \Delta^\Delta / (\Delta - 1)^{\Delta + 1}\)
and the representation of  \(\alpha = 1 / \Delta_c(\lambda) < 1 / \Delta\)).
Choosing \(k_r = \lceil r^\gbeta \rceil \), with \(\gbeta \in (0, 1)\), the term \((\Delta\alpha)^{k_r/q} \, \leq \, (\Delta\alpha)^{(r^\gbeta)/q}\) decreases rapidly as \(r\) increases. Thus, we achieve the same conclusion as in~\eqref{e.BL-spin-stab-1*}.  

The same reasoning applies when modifying ~\eqref{e.BL-spin-stab-2*}, thereby confirming the result of Lemma~\ref{stablifting} whenever the strong spatial mixing bound in~\eqref{spamixa-strong} decays rapidly. The latter is guaranteed by our assumption on $\lambda$ and \eqref{e.smor}. This allows us to use  Lemma~\ref{stablifting} in the proof of Theorem~\ref{spinmain} giving 
the asymptotic normality for the sum of 
spins in the hard-core model.
\end{proof}

\stoptoc
\subsubsection{Derivation of \eqref{e.smor}---Strong spatial mixing for the hard-core model; \protect\cite{sinclair2017spatial}, \protect\cite{Weitz2006}}
\label{r.Hard-core-Sinclair}
\resumetoc

In the following, we provide additional details on how the strong spatial mixing~\eqref{e.smor} is established in~\cite{sinclair2017spatial} for the hard-core model. Since the specific form of the bound we use may not be transparent in the proof of Theorem 1 in \cite{sinclair2017spatial}, we describe the key steps of the derivation here for self-containment. Broadly the proof of \eqref{e.smor} involves two steps, firstly reducing the right-hand side in \eqref{e.smor} to that of a tree via the `Weitz reduction' and then a very involved analysis of the same for a tree using message passing recursions. In particular, for the second step, we refer many times to specific claims in \cite{sinclair2017spatial} for proofs. We shall explain these steps and then describe how the proof follows from these two steps. 

 A similar approach was used in~\cite{sinclair2017spatial} to establish strong spatial mixing for the monomer-dimer model, where instead of independent vertices in the hard-core model, random matchings on the graph are considered, also controlled by the activity parameter; see also~\cite{quitmann2024decay}. One can also apply this approach to the Ising model and
 random cluster model; see \cite[Theorem 28]{shao2025zerofreenessneedweitztypefptas}. This suggests that the 
 central limit theorem should hold for these models as well under an analogous range of the activity parameters.

\vskip.3cm

\noindent \textsc{First step:} \, \, The `Weitz reduction' or `self-avoiding walk (SAW) representation' as stated in~\cite[Theorem 4]{sinclair2017spatial} is based on the original result of~\cite{Weitz2006}. The SAW tree (where SAW stands for Self-Avoiding Walk) \( \cT_{SAW} = \cT_{SAW}(x, \G(\X)) \) represents all possible self-avoiding paths in \( \G(\X) \) that start from \( x \). Its root corresponds to the trivial self-avoiding walk that begins and ends at \( x \). The tree is built recursively: for each node \( w \) in \( \cT_{SAW} \), its children represent all possible extensions of the self-avoiding path associated with \( w \), each extended by exactly one additional step. The theorem asserts that, for any vertex \( x \) in a (finite or locally finite, infinite) graph \( \G(\X) \), the tree \( \cT_{SAW} \) rooted at $x$ may be equipped with a boundary condition \( w_{\cW} \) on a subset \( \cW \subset \cT_{SAW} \) (see \cite[Appendix~A]{sinclair2017spatial} and ~\cite[Section~3]{Weitz2006})  such that, for the hard-core Gibbs model on \( \G(\X) \) specified by~\eqref{eqn:re-hardcore} with activity parameter \( \lambda > 0 \), and any boundary condition \( v_{\cZ} \) on \( \cZ \subset \G(\X) \), we have:
\begin{equation}\label{e:Weitzred1}    
R(x, \X; \v_{\cZ}) = R(x, \cT_{SAW}; w_{\cW} \cup \v_{\cZ}),
\end{equation}
where \( v_{\cZ} \) on the right-hand side denotes the translation of the boundary condition from \( \G(\X) \) to \( \cT_{SAW} \). The combined boundary condition \( w_{\cW} \cup v_{\cZ} \) is defined by first applying \( w_{\cW} \) and then \( v_{\cZ} \), with \( v_{\cZ} \) overriding \( w_{\cW} \) wherever they both apply and differ. Thanks to the above reduction, to prove \eqref{e.smor} it now suffices to bound
\begin{equation}\label{e:Weitzred2} 
|R(x,\cT_{SAW}; w_{\cW} \cup v_{\cZ}) - R(x,\cT_{SAW}; w_{\cW} \cup z_{\cZ})|.
\end{equation}
Thus the difference on the boundary due to the conditions $w_{\cW} \cup v_{\cZ}$, $w_{\cW} \cup z_{\cZ}$ arises from difference between $v_{\cZ}$ and $v_{\cZ}$. \\

\noindent \textsc{Second step:} \, \,We describe how to derive bounds for \eqref{e:Weitzred2}  for any tree. Consider a tree~$\cT$ rooted at $x$ and denote the graph metric on $\cT$ by $d_{\cT}$. For a finite tree $\cT$, a {\em cutset} $\gC$ is a set of vertices in~$\cT$ such that (i) any path from the root $x$ to a leaf $y$ where the graph distance \( d_{\cT}(x,y)\ge \max_{y' \in \gC} d_{\cT}(x,y') \) must pass through $\gC$, and (ii) the vertices in $\gC$ form an `antichain', meaning that for any two vertices in $\gC$, neither is an ancestor of the other in $\cT$. For rooted, locally finite, infinite trees, condition (i) is replaced by the condition that any infinite path starting from the root $x$ must pass through $\gC$. A key example of a cutset we will use later is the set of all vertices at a fixed distance \( l \) from $x$ in $\cT$. For a cutset $\gC$, let $\cT_{\le \gC}$ be the subtree of $\cT$ obtained by removing the descendants of vertices in $\gC$, and let $\cT_{<\gC}$ denote the subtree rooted at $x$ obtained by  removing the vertices in $\gC$. Also, let $\gC^+$ denote the set of the children of the vertices in  $\gC$ (here $\gC$ may not necessarily form a cutset). For $y \in \cT$, we denote by $\cT_y$ the sub-tree of $\cT$ rooted at $y$.

We assert that the following adaptation of~\cite[Lemma~3]{sinclair2017spatial} applies:
For the hard-core Gibbs model specified by~\eqref{eqn:re-hardcore} with activity parameter $\lambda>0$, there exist constants \( C=C(\lambda)<\infty \), \( q=q(\lambda)>1 \), and \( \alpha=\alpha(\lambda) > 0 \) such that this specification, when applied to any finite or infinite locally finite tree $\cT$ rooted at $x$, satisfies the following decay of the occupation ratio at the root:
\begin{equation}\label{eq:hard-core-tree-mixing}
|R(x,\cT;v_{\cZ}) - R(x,\cT;z_{\cZ})|^q \le C^q \sum_{y \in \gC} \alpha^{d_{\cT}(x,y)},
\end{equation}
where $\gC$ is a cutset in $\cT$ located at a distance of at least 1 from the root \( x \).
Here, $\cZ$ is a subset of $\cT_{\le \gC^+}$ not containing \( x \), and \( v_{\cZ} \) and \( z_{\cZ} \) are two boundary configurations on $\cZ$ that differ only on $\gC^+$ (agreeing on the spins fixed in \( \cZ \cap \cT_{\le\gC} \)).

Note that the above claim constitutes a slight, two-fold modification of \cite[Lemma~3]{sinclair2017spatial}: we have dropped the requirement that $\gC$ has no leaves, and we formulate the conditional probability of the spin (and hence the occupation ratio) at $x$ using our boundary condition terminology---represented by $v_{\cZ}, z_{\cZ}$---rather than the `initial conditions' used in the original work, which refer to fixed values of the occupation ratios $R$ at the vertices in $\cZ$. Note that the occupation ratio takes real values in $[0, \lambda]$ at all vertices, except it is infinite at points where the boundary condition sets the spin to $1$. To avoid this issue in the proof, $\gC$ was  assumed to have no leaves in \cite[Lemma~3]{sinclair2017spatial}.

The rest of this step is to explain choices of $C,q,\alpha$ along with  details showing how \eqref{eq:hard-core-tree-mixing} can be derived from \cite[Lemma~3]{sinclair2017spatial} by substituting specific choices from \cite[Section 4]{sinclair2017spatial} pertaining to the hard-core model. Our notation coincides with that in \cite{sinclair2017spatial}.  As in  \cite{sinclair2017spatial}, 
we choose $\phi(x)= \sinh^{-1}(\sqrt{x})=\text{arcsinh}(\sqrt{x})$ with derivative $\phi'(x) = 1/(2\sqrt{x(1+x)})$ and so we have that $M := \sup_{x \in [0,\lambda]}\phi(x) = \sinh^{-1}(\sqrt{\lambda})$ and $L := \inf_{x \in [0,\lambda]}\phi'(x) = 1/(2\sqrt{\lambda(1+\lambda)})$. Since the occupation ratio at any internal node of the tree is upper bounded by $\lambda$, we have  
via the mean-value theorem  that
$$|R(x,\cT;v_{\cZ}) - R(x,\cT;z_{\cZ})| \leq L^{-1} |\phi(R(x,\cT;v_{\cZ})) - \phi(R(x,\cT;z_{\cZ}))|.$$
Putting the  constant $C$ in \eqref{eq:hard-core-tree-mixing} to be 
$$C:=\frac{M}{L}=2\sinh^{-1}(\sqrt{\lambda}) \, \sqrt{\lambda(1+\lambda)}$$
it suffices to show
\begin{equation}
\label{eq:tree-mixing-phi}
|\phi(R(x,\cT;v_{\cZ})) - \phi(R(x,\cT;z_{\cZ}))|^q \leq M^q \sum_{y \in \gC} \alpha^{d_{\cT}(x,y)}.
\end{equation}
This claim is the analogue of \cite[(20)]{sinclair2017spatial} but with 
boundary conditions on  $\cZ\subset\cT_{\le\gC^+}$  (instead of $R$-related initial condition on $\cT_{\le\gC^+}$, as  discussed above) and is similarly  proved inductively on $d_{\cT}(x,\gC) = \inf_{y \in \gC}d_{\cT}(x,y)$. Before sketching the proof, we need to mention the choices of $q,\alpha$ for which the above claim holds.

Let $\lambda_c(t) : (1,\infty) \to (0,\infty)$ be the strictly decreasing function $t \mapsto \frac{t^t}{(t-1)^{t+1}}$ and let $\Delta_c(\lambda)$ be the unique solution in $t$ to the equation $\lambda_c(t) = \lambda$ for $\lambda \in (0,\infty)$ i.e., $\lambda_c(\Delta_c(\lambda)) = \lambda$. As explained below \eqref{e.smor}, $\Delta_c(\lambda) > \Delta \geq 1$ if $\lambda < \lambda_c(\Delta)$ where $\Delta$ is the upper bound on the connective constant. Assuming this, the decay constant $\alpha = \alpha(\lambda)$ (defined in \cite[(9)]{sinclair2017spatial}) and the exponent $q = q(\lambda)$ are defined as follows 
$$ \alpha(\lambda) := \frac{1}{\Delta_c(\lambda)} \, \, \textrm{and} \, \, q := 1 - \frac{\Delta_c(\lambda) - 1}{2}\log(1 + \frac{\Delta_c(\lambda) - 1}{2}),$$
where the well-definedness follows because $\Delta_c(\lambda) > 1$; see Lemma 6 of \cite{sinclair2017spatial}. Additionally, as in \cite{sinclair2017spatial}, take $1/a = 1 - 1/q$, though this exponent will not appear explicitly in our derivations. \\

Now we sketch why \eqref{eq:tree-mixing-phi} holds with the above choices of $\alpha,q$. Consider the base case $d_{\cT}(x,\gC) = 1$. Let $y_1,\ldots,y_k$ be the immediate descendants of $x$ in the tree $\cT$. Suppose $R(y_i,\cT_{y_i};v_{\cZ}) = \infty$ for some $y_i \in \cZ \cap \gC$. Then by the assumption that boundary conditions are differing only on $\gC^+$, $R(y_i,\cT_{y_i};z_{\cZ}) = \infty$ as well and in this case $R(x,\cT;v_{\cZ}) = R(x,\cT;z_{\cZ}) = 0$. Thus  \eqref{eq:tree-mixing-phi} holds trivially here. Otherwise, we have $R(y,\cT_{y};v_{\cZ}) \in [0,\lambda]$  for all $y \in \gC$ and \cite[Lemma 2]{sinclair2017spatial} along with definition of $M$ gives that
$$|\phi(R(x,\cT;v_{\cZ})) - \phi(R(x,\cT;z_{\cZ}))|^q \leq \alpha \sum_{y_i \in \gC} |\phi(R(y_i,\cT_{y_i};v_{\cZ})) - \phi(R(y_i,\cT_{y_i};z_{\cZ}))|^q \leq M^q \sum_{y \in \gC} \alpha,$$
thereby immediately proving \eqref{eq:tree-mixing-phi} in this case as well. Thus the base case holds in \eqref{eq:tree-mixing-phi}. 

Now, for some $k\in \N$, assume~\eqref{eq:tree-mixing-phi} for $d_{\cT}(x,\gC)\le k$.
As above, consider  immediate descendants $y_i$ of $x$ in the tree $\cT$.  Either accounting for possibility that $R(y_i,\cT_{y_i};v_{\cZ}) = R(y_i,\cT_{y_i};z_{\cZ}) = \infty$ for $y_i \in \cZ \cap \cT_{\leq C}$ or using \cite[Lemma 2]{sinclair2017spatial} as in the first inequality above and induction hypothesis, one can obtain \eqref{eq:tree-mixing-phi}. As argued before \eqref{eq:tree-mixing-phi}, this suffices to prove \eqref{eq:hard-core-tree-mixing}. \\

\noindent \textsc{Proof of \eqref{e.smor} from the first two steps}:  
Choose $\gC$ to be set of vertices at distance $l-1$ from $x$ in $\cT_{SAW}$; 
by the construction of $\cT_{SAW}$ these are all self-avoiding paths of length~$l-1$ from $x$ in $\G(\X)$ and therefore
$|\gC| = N(x,l-1;\G(\X))$.  So $\gC^+$  
corresponds to the set of all self-avoiding paths from $x$ to
$\partial\gB_{l-1}(x)$ where $\gB_{l-1}(x)=\{y \in \G(\X): \ d_{\G(\X)}(x,y) \leq l-1\}$. Since the boundary conditions $v,z$ differ on $\cZ$ in $\G(\X)$  in~\eqref{e.smor} 
only when the graph distance to $x$ is greater than or equal to $l$, we have that $v_{\cZ},z_{\cZ}$ may differ on paths being elements in $\gC^+$ but not in $\gC$. Moreover, due to the cutset nature of $\gC$ and the spatial Markov property, the spin values assigned in $\cZ$ on paths beyond $\gC^+$ do not affect the distribution of the spin at $x$.
Hence  \eqref{eq:hard-core-tree-mixing} holds in this case and using that $d_{\cT_{SAW}}(x,y) = l-1$ for $y \in \gC$, we obtain
\begin{align*}
 |R(x,\cT_{SAW};w_{\cW} \cup v_{\cZ}) - R(x,\cT_{SAW};w_{\cW} \cup z_{\cZ})|^q & \leq C^q \sum_{y \in \gC}\alpha^{d_{\cT_{SAW}}(x,y)} \\
& = C^q |\gC|\alpha^{l-1} = C^q N(x,l-1;\G(\X))\alpha^{l-1}.
\end{align*}
Thus \eqref{e.smor} now follows from \eqref{e:Weitzred1} and \eqref{e:Weitzred2}. The choices of $C,q$ and $\alpha$ have been explained in the Second Step above.

\section{Interacting diffusions on spatial random graphs}
\label{NEW-s:id_sprg}

In this section we define a model of interacting diffusions on spatial random graphs and establish limit theorems for statistics of these diffusions. In contrast to mean field models, the diffusions considered here interact only with finitely many neighbors. Similar to the spin variables encountered in the previous section, diffusions do not have  a stopping set property, but  they nonetheless satisfy fast \BL-localization, as well as the even  stronger $L^2$-stabilization, when defined on stabilizing interaction graphs. 
The diffusion model on graphs considered here was studied in \cite{lacker2019local}, which provides conditions ensuring the existence of the model in infinite space.  
Building on this prior work---and in contrast to the finite-window approach commonly used for spin systems---we study in this section the limiting behavior of sums of scores of the interacting diffusions over the infinite window, using the results from Section~\ref{s:ltthms_stat_marking}.
We refer to \cite[Section 1]{lacker2019local} and \cite{ramanan2023interacting} regarding the motivation and applications of this model. An interacting diffusion model on deterministic point sets has been considered in \cite{chargaziya2025stochastic} and one may  consider such models on random graphs using the methods here.

\subsection{The diffusion model}
\label{NEW-ss:model_thm}

Using the same graph notation as for spin systems in  Section~\ref{sss.graph-stabilization}, we now consider a finite or countable  set of sites (nodes), denoted by $\X$, and a locally finite graph $\G(\X)$ defined on $\X$. For each $x \in \X$, we denote by 
$$
N_x = N_x(\X):= \{x' \in \X : x' \sim x\}
$$
the (finite) set of neighbors of $x$ in $\G(\X)$.
Fix a time horizon $t_0 \in [0, \infty)$ and a diffusion dimension $d' \in \N$. We consider a system of interacting $\R^{d'}$-valued diffusions
$$
M(x,t) := M^{\G(\X)}(x,t), \quad x \in \X, \; t \in [0, t_0],
$$
evolving on vertices of the graph  $\G(\X)$, interacting only via its edges, and defined by the stochastic differential equations:
\begin{equation}
\label{NEW-e:sys_ID_V}
\md M(x,t) = b(t, M[x,t], M[N_x,t])\, \md t + \sigma(t, M[x,t], M[N_x,t])\, \md Z_x(t),
\end{equation}
where 
 $b$ and $\sigma$ are the {\em drift} and {\em diffusion} coefficients, respectively, understood as $\R^{d'}$-valued and $\R^{d'} \times \R^{d'}$-valued functions, and 
$\{Z_x(\cdot)\}_{x \in \X}$ are i.i.d. standard Brownian motions in $\R^{d'}$; in this notation $M[x,t]$ denotes the path of $M(x,\cdot)$ in the interval $[0,t]$ and $M[N_x,t]:=\sum_{x'\in N_x}\delta_{M[x',t]}$;  $\sigma$ here is matrix-valued but treated equivalently as a vector in $\R^{d'} \times \R^{d'}$. Indeed,  following~\cite{lacker2019local}, we allow path-dependence in the coefficients $b$ and $\sigma$ in~\eqref{NEW-e:sys_ID_V} `both because such interactions arise in applications and because this does not complicate the arguments'. The SDE \eqref{NEW-e:sys_ID_V} says that  particles interact directly only with  (finitely many) neighbors in the graph. As  with  interacting particle system models considered in the next section,  this constitutes the sparse/thermodynamic regime, in contrast to the mean field regime, where particles interact with all other particles.

 When the {\em initial conditions $\{M(x,0)\}_{x \in \X}$ 
 are uniformly  bounded} and the drift and diffusion {\em coefficients $b$, $\sigma$ are Lipschitz functions}, as detailed in the next section, the system~\eqref{NEW-e:sys_ID_V} admits a {\em unique solution with the  family of processes $\{M(x,t)\}_{x \in \mathcal{X}}$},  taking values in the space of continuous paths with $t \in [0, t_0]$, as established in \cite{lacker2019local}. In what follows we will call 
this collection of processes {\em the family of interacting diffusions} on the graph  $\G(\X)$.

Large ensembles of interacting diffusions are used to model complex dynamical behavior arising in statistical physics, biology, and neural networks; see \cite{ramanan2023interacting} and also \citet{luconstannat}.  Mean field models have received considerable attention and their dynamics are relatively well understood when $\G(\X)$ is the complete graph.  The `$P$-nearest-neighbor
model', $P \in (0,1)$, has also been studied \cite{luconstannat}, where here the neighbors of a site consist of all particles within a window centered at the site and whose volume equals $Pn$. In this model and others, one often takes the diffusion $\sigma$ to be a fixed constant.  Here we consider diffusion models which almost surely satisfy  $P = o(1)$, a regime which is not well studied. 

\begin{exe}[Markovian drift]
\label{ex:markovdrift}
We  give an example of a simple Markovian model of drift and diffusion. This model has been investigated in mean-field settings; see for example \cite[Example 2.1]{redig2010short} or \cite[Section 1.2.2]{luconstannat}. In this model, we put $\sigma \equiv 1$ (i.e., the constant function) and take $b$ as follows:
\[b(t, M[x,t], M[N_x,t]) = b(M(x,t),M(N_x,t)) = F_1(M(x,t)) + \frac{1}{|N_x|} \sum_{y \sim x}F_2(M(x,t),M(y,t)),\]
where $F_1,F_2$ are measurable functions on $\R^{d'}$ and $\R^{d'} \times \R^{d'}$ respectively.
\end{exe}

\subsubsection{Lipschitz assumption on coefficients}  
\label{NEW-sss.ID-G}
For classical diffusion processes defined by stochastic differential equations, Lipschitz conditions on drift and diffusion coefficients are essential to ensure  existence and uniqueness of strong solutions via Gr\"onwall's inequality.
A similar principle applies to interacting diffusions on graphs, where such conditions not only guarantee well-posedness but also imply (spatial) $L^2$-stabilization: 
given a fixed time horizon,
the state at a given node becomes increasingly independent of distant nodes as their separation grows, capturing how influence disperses across the graph.

Before introducing this crucial Lipschitz condition, 
we establish some notation. 
Let \\ $\M := \cC([0,t_0],\R^{d'})$ denote the path space of continuous functions, equipped with the topology of uniform convergence on compact sets. For truncated paths up to time $t$, where $0\le t\le t_0$, we use the notation $\M(t) := \cC([0,t],\R^{d'})$.
The paths of the processes $M(x,\cdot)$,
$x\in\X$, defined in equation~\eqref{NEW-e:sys_ID_V} and subsequently referred to as $M[x]$, are elements of $\M$ (i.e., $M[x]\in\M$). Their truncations are denoted by $M[x,t]\in\M(t)$ in~\eqref{NEW-e:sys_ID_V}. We denote the zero function by $0[x]$ i.e., $0(x,t) = \0$ for all $t \leq t_0$.

More generally, for a given function $m[x]\in M(t)$ representing a possible entire
trajectory $m(x,s)$ over $s \in [0,t_0]$, at $x\in\X$, 
its truncation to $s\in[0,t]$ is denoted by $m[x,t]$ and 
its norm is $|m[x,t]| := \sup_{s \in [0,t]} |m(x,s)|$, where $|\cdot|$ denotes the Euclidean norm in~$\R^{d'}$.

With the above notation, following \cite[Lipschitz Assumption A']{lacker2019local}, we impose the following Lipschitz conditions on the coefficient functions $b$ and $\sigma$. These conditions hold trivially for the Markovian drift model in Example \ref{ex:markovdrift} if $F_1,F_2$ therein satisfy the Lipschitz assumptions.
\begin{definition}[Lipschitz diffusion coefficients]
\label{NEW-d:lip_assum}
The jointly measurable functions $(b,\sigma) :[0,\infty) \times \M \times \hat\cN_{\M}\to \R^{d'}\times (\R^{d'} \times \R^{d'})$ are said to be \emph{Lipschitz diffusion coefficients} if they satisfy the non-anticipative property, meaning that for all $t \in [0,t_0]$, $x \in \X$, any finite subset $\X' \subset \X$, and for $m_1[x]$, $m_2[x] \in \M$, with $m_1[\X']:=\sum_{x'\in\X'}\delta_{m_1[x']}\in\hat \cN_{\M}$, and similarly for $m_2[\X']\in\hat\cN_{\M}$, we have: 
\begin{align*}b(t,m_1[x], m_1[\X']) &= b(t,m_2[x],m_2[\X']),\\
\sigma(t,m_1[x], m_1[\X']) &= \sigma(t,m_2[x], m_2[\X'])\\ \noalign{whenever}
m_1[x,t] &= m_2[x,t],\\
m_1[\X',t] &= m_2[\X',t].
\end{align*}
Furthermore, there exists a {\em Lipschitz diffusion constant} $K_t < \infty$ such that:
\begin{align*}
& \sup_{s\in[0,t]}|b(s,m_1[x], m_1[\X']) - b(s,m_2[x], m_2[\X])|
+ |\sigma(s,m_1[x],m_1[\X']) - \sigma(s,m_2[x], m_2[\X'])| \\
& \qquad \leq K_t \Big( \|m_1[x,t] - m_2[x,t]\| + \frac{1}{|\X'|}\sum_{x' \in \X'} \|m_1[x',t] - m_2[x',t]\| \Big),
\end{align*}
with the convention that $\frac{1}{|\X'|}\sum_{x'\in\X'}(\dots)=0$ if $\X'=\emptyset$. We assume, moreover, that  $\sup_{t\in[0,t_0]}K_t<\infty$ and that $K_t$ are increasing in $t$.
Finally, using the zero-function notation $0[x]\in\M$ and $0[\X']):=\sum_{x'\in\X'}\delta_{0[x']}$, the coefficients $(b,\sigma)$ satisfy (uniformly in $x$ and $\X'$)
$$\int_0^{t_0}(|b(s,0[x],0[\X'])|^2 + |\sigma(s,0[x],0[\X'])|^2) \md s < K_{t_0}.$$
\end{definition}

The Lipschitz assumptions imposed on the diffusion coefficients lead to a crucial property, which can be thought of as a {\em pre-Grönwall inequality for interacting diffusions}. This  property is formulated in Lemma~\ref{NEW-l:pre-Gron} in Section~\ref{NEW-ss:proofs_sid}.
This pre-Grönwall property is a cornerstone of the
Grönwall's inequality, used in the 
standard Picard iteration argument  to  establish the existence and uniqueness of the diffusion described in~\eqref{NEW-e:sys_ID_V} on any arbitrary, locally finite graph $\G(\X)$, provided that the initial conditions are uniformly bounded. This is demonstrated in \cite[Theorem 3.1, see also Appendix~C]{lacker2019local}, which we shall recall in Lemma~\ref{NEW-l:strong-solution-diffusion}.
Furthermore, this property will be instrumental in demonstrating {\em $L^2$-stabilization of the diffusions on graphs}, see Lemma~\ref{NEW-l:spat-mixing-diffusion} in Section~\ref{NEW-ss:proofs_sid}.

\subsubsection{Stabilizing  random graphs}
\label{NEW-sss.ID-graphs}

In this section, we shift our focus to the graph $\G(\P)$ whose vertices are the points $x$ of a simple point process $\mathcal{P}$ in $\mathbb{R}^d$, each hosting diffusion paths in dimension $d'$, which may differ from the ambient space dimension $d$.
We refer to Section \ref{s:notation} for all notational formalism regarding marked point processes, including  $\mP_{\bk x1p}$, which denotes the Palm probability distribution of the point process $\P$ conditional on the existence of  points $x_1,...,x_p$.

Specifically, as in Section~\ref{sss.graph-stabilization},
we represent $\G(\P)$ as an auxiliary  measurable marking $\{\delta_{(x,\bar N_x)}\}_{x \in \P}$, where the marks $\bar N_x$ represent neighbors of $x$ in $\G(\P)$
whose locations are expressed with respect to $x$; i.e., they satisfy
$\bar N_{x} + x \subset \P$ for all $x \in \P$, and $y \in \bar N_x$ if and only if $-y \in \bar N_{x+y}$. Assuming $\bar N_x \in \hat\cN_{\R^d}$ ensures local finiteness of the graph $\G(\P)$.

We make also the simplifying assumption that the marks $\bar N_x$ (and hence the graph $\G(\P)$) are (measurable) {\em deterministic functions of solely the input process} $\P$,  $\bar N_x=\bar N_x(\P)$
considering
$\bar N_\cdot(\cdot):\R^d \times {\cN}_{\R^d}\longrightarrow \hat{\cN}_{\R^d}$  
as a marking function as per Definition~\ref{d:mark_fn},
admitting this time all input configurations of points $\mu \in \cN_{\R^d}$, including those with infinite input.

Again as in Section~\ref{sss.graph-stabilization}, but for the infinite window $W_\infty=\R^d$, this marking function $\bar N$
 admits an interaction range of stabilization $S: \R^d \times \cN_{\R^d} \longrightarrow \N\cup\{\infty\}$ of the graph $\G(\mu)$ as per Definition~\ref{def.Stab.Radius}, with
$S(x,\mu):=R^{\bar N}(x;\mu)$.
Clearly, $S$  satisfies~\eqref{e:cons_adm_gr}, \eqref{stopspin}, and~\eqref{stopspin-new}.

In this graph construction scenario, a key stabilization assumption for the graph $\G(\P)$ is formulated analogously to Definition~\ref{d:stabilizing-graph}, but applied directly in the infinite setting via~\eqref{stab-infinite}. Specifically, for every $p \in \N$, we assume that there exists a fast decreasing function $\varphi_p'$, as defined in~\eqref{varphibd}, such that
\begin{equation} \label{sdecayspin-infty}
\sup_{x_1, \ldots, x_p \in \R^d} \mP_{\bk x1p}
\bigl(S(x_1, \P) > s\bigr) \leq \varphi_p'(s), \quad s > 0.
\end{equation}
Without loss of generality, we further assume that the functions form an
increasing--decreasing family; cf.\ the terminology in
Section~\ref{s:prelim}.

We remark that  the above assumption of stabilization of the graph $\G(\P)$ allows it to be approximated by local graphs. This makes  a subtle difference from the setting of spin systems in Section~\ref{sss.graph-stabilization}, where the local graphs $\G(\P\cap W_n), n \in \N$, are, in fact, the only available representations of the interaction, and their stabilization only on finite windows ($W_n$, $n\in\N$) ensures consistency of these representations across increasing windows.

In Section~\ref{s:future}, we will discuss how this framework can be extended to include more general randomized graph models.

\subsection{Main results for the diffusion model}
\label{NEW-sss.ID-results}
Our goal is to study  scores associated with the individual paths 
of the diffusion on the infinite graph $\G(\P)$.
 Specifically, the sums of these scores over points $x\in\P \cap W_n$ will be the subject of the Gaussian fluctuations presented in this section.

To fully define the diffusion on this graph, we introduce the following 
marked point process as an input  (upon which we will construct  diffusion trajectories $M[x]$, for $x\in\P$): For ambient and diffusion dimensions $d,d'\in\N$, consider a marked point process
$\tP:=\{\delta_{(x, M(x),Z_x)}\}_{x\in\P}$,
where $x\in\R^d$ are points (sites)  of the diffusions, 
$M(x)=:M(x,0)\in\R^{d'}$ are  initial condition of the diffusions at these sites,  $Z_x\in\M=\cC([0,t_0],\R^{d'})$ refer to standard Brownian motions on the interval $[0,t_0]$ in $\R^{d'}$, assumed to be independent given the object  $\{\delta_{(x, M(x))}\}_{x\in\P}$. In other words, this latter entire object is measurable at time $t=0$ with respect to a filtration that accommodates the Brownian motions.

 Recall that the graph $\G(\P)$ is a function of the ground process $\P$ (represented by the neighborhood marks $\bar N_x$ introduced  in Section~\ref{NEW-sss.ID-graphs}, which are not part of the constitutive marks of the input process $\tP$ directly related to the diffusion).

Assuming the Lipschitz coefficients $(b,\sigma)$ as in Definition \ref{NEW-d:lip_assum}, and given that the initial conditions $M(x)$ are almost surely deterministically, uniformly bounded, then based on  \cite[Theorem 3.1]{lacker2019local} (see also Lemma \ref{NEW-l:strong-solution-diffusion}), almost surely for $\tP$, a path-wise unique strong solution $M = M^{\G(\P)}$ to the spatial system \eqref{NEW-e:sys_ID_V} on $\G(\P)$ exists. By its very construction, the individual paths $M[x]=M^{\G(\P)}[x]\in\M$ become new, measurable marks for points  $x\in\P$ 
(deterministically constructed given the input marked process  $\tP$ and the graph $\G(\P)$). 

Finally, for $\tx=(x,M(x),Z_x)\in\tP$, we consider a real-valued score  of the trajectory of the solution $M$ at $x$:
\begin{equation} \label{NEW-xisysID}
\xi (\tx, \tP) := \xi^{(h,\G(\P))} (\tx, \tP) := h(M^{\G(\P)}[x]),
\end{equation}
where $h : \M \to \R$ is a Lipschitz function with respect to the sup-norm on $\M = \cC([0, t_0], \R^{d'})$.

Given $p \in [1, \infty)$, we say that $\xi = \xi^{(h,\G(\P))}$ satisfies the  $p$-moment condition if
\begin{equation}
\label{NEW-momh1}
\sup_{1 \leq q \leq p} \sup_{x_1,\ldots,
x_{q} \in \R^d} \sE_{\bk x1{q}}\left[ \max\left(1, |\xi(\tx,\tP)|^p\right) \right] \leq M_p^{\xi} < \infty
\end{equation}
for some constant $M_p^{\xi}$, which is assumed to be no smaller than $M_{p'}$ for all $p' \in [1,p]$.

Uniform boundedness of the initial conditions $M(x)$ and the Lipschitz form of the diffusion coefficient functions ensure that the diffusion trajectories $M[x]$ are uniformly $L^2$-bounded; see Lemma~\ref{NEW-l:strong-solution-diffusion}. Consequently, given a Lipschitz function $h$, this implies $M_p^\xi < \infty$ for $p \in [1,2]$ and hence guarantees that $\xi$
satisfies the requisite $(1 + \varepsilon)$-moment condition needed to prove expectation asymptotics given 
in the next limit theorem.

Now, we present a  main result of this section---the  central limit theorem for interacting diffusions on spatial random graphs.

\begin{theorem}[Limit theory for statistics of interacting diffusions]
\label{NEW-t:sysID} 

For $d,d'\in\N$, $t_0\in(0,\infty)$, let $\tP:=\{\delta_{(x, M(x),Z_x)}\}_{x\in\P}$ be a 
marked point process 
on $\R^d\times\R^{d'}\times\M$
such that   $Z_x\in\M=\cC([0,t_0],\R^{d'})$ are i.i.d. standard Brownian motions on the interval $[0,t_0]$ in $\R^{d'}$, independent given the  marked point process $\{\delta_{(x, M(x))}\}_{x\in\P}$ of sites with initial conditions. We assume that  the initial conditions $M(x)$, $x\in\P$, are almost surely deterministically, uniformly bounded, i.e.,  $\sup_{x\in\P}|M(x)| \le L<\infty$.  
Moreover, we assume that the marked point process $\tP$ has summable exponential $\B$-mixing correlations as in Definition \ref{def.A2} and has bounded reduced  Palm intensity function as in~\eqref{bdedreduced}.
Consider a graph $\G(\P)$ constructed deterministically on $\P$ 
which is fast stabilizing  with the interaction range 
satisfying ~\eqref{sdecayspin-infty}.
Given Lipschitz diffusion coefficients $(b,\sigma)$ as in  Definition \ref{NEW-d:lip_assum} consider the strong solution to the interacting diffusions system \eqref{NEW-e:sys_ID_V} on $\G(\P)$, with initial conditions $M(x)$, realized by  individual paths $M^{\G(\P)}[x]\in\M$ for points  $x\in\P$. Finally,
define the measures
\begin{equation}
\label{NEW-SysID}
\hat{\mu}_n^{\xi} := \sum_{x \in \P \cap W_n} \xi (\tx,\tP)  \delta_{n^{-1/d}x},
\end{equation}
where $\xi =\xi^{(h, \G(\P))}$ are real-valued scores of the  individual paths $M^{\G(\P)}[x]\in\M$  given by~\eqref{NEW-xisysID} with  a Lipschitz function $h : \M \to \R$.
\begin{enumerate}[wide,label=(\roman*),  labelindent=0pt]
\item \label{NEW-i:CLT-difussion} If 
the  score $\xi$ in \eqref{NEW-xisysID} satisfies the $p$-moment condition~\eqref{NEW-momh1} for all  $p \in (1, \infty)$, then the random measures $(\hat{\mu}_n^{ \xi})_{ n \in \N}$  at \eqref{NEW-SysID} satisfy the central limit theorem,  i.e., for  $f \in \B(W_1)$ such that $\Var{\hat{\mu}_{n}^{\xi}(f)} = \Omega(n^{\nu})$ for some $\nu > 0$, we have
\[ (\Var{\hat{\mu}_{n}^{\xi}(f)})^{-1/2}(\hat{\mu}_{n}^{\xi}(f) - \sE\hat{\mu}_{n}^{\xi}(f)) \stackrel{d}{\Rightarrow} Z. \]
\item \label{NEW-ii:mean-varianc-diffusion} If the input process $\{\delta_{(x,M(x))}\}_{x \in \P}$ is stationary and the graph construction $\G(\P)$ is translation invariant, then the mean asymptotics can be expressed for all $f \in \B(W_1)$ as:
\begin{align}
 \Big| n^{-1} \E \hat{\mu}_n^{\xi}(f) -  \rho \E_{\0} \xi(\tilde\0, \tP) \int_{W_1}f(x)\,\md x\Big| = O(n^{-1/d}), 
\end{align}
whereas under an extra $p=(2+\epsilon)$-moment condition, for some $\epsilon>0$,  the variance asymptotics  can be expressed as
 \begin{align}
\label{NEW-sysyIDlln} & \lim_{n \to \infty} n^{-1} \Var \,\hat{\mu}_n^{\xi}(f) =  \sigma^2(\xi) \int_{W_1} f(x)^2 \,\md x \in [0,\infty),
\end{align}
where 
\begin{align}
\label{NEW-sigdefsysyID} \sigma^2(\xi) & :=
\rho\sE_{\0}\xi^2(\tilde{\0},\tP)   \\
 &+ \int_{\R^d} \left( \sE_{\0,x}(\xi(\tilde{\0}, \tP) \xi(\tilde{x}, \tP)) \rho^{(2)}(\0,x) - \rho^2\sE_{\0}(\xi(\tilde{\0},\tP)) \sE_{x}(\xi(\tilde{x},\tP))  \right) \md x \in [0, \infty). \nonumber 
\end{align}
 \end{enumerate}
 \end{theorem}

\begin{proof}
We deduce all statements from Proposition~\ref{t:clt_linear_marks_new}, relying upon the upcoming Lemmas \ref{NEW-l:strong-solution-diffusion} and \ref{NEW-l:spat-mixing-diffusion} for stabilization properties of interacting diffusions as well as the graph stabilization argument used in Lemma \ref{stablifting}. Given the mixing and moment conditions on $\tilde\P$ and $\xi$, it remains to
localize the scores $\xi(\tx,\tP)=\xi^{(h,\G(\P))}(\tx,\tP)$, defined via
diffusion trajectories in the infinite model, by means of short-range marking
functions of the form $\xi(\tx,\tP\cap B)$, where $B\subset \R^d$ is bounded.
These localized versions should approximate $\xi(\tx,\tP)$ in the \BL-sense as
$B \uparrow \R^d$.  This corresponds exactly to the fast \BL-localization of $\xi$ as required in Definition~\ref{def.Lp-stabilizing_marking}\ref{i.BL-localizing} and~\ref{i.BL-localizing-fast}, and invoked in Proposition~\ref{t:clt_linear_marks_new}.
To achieve this, we consider the following score function 
\begin{equation}\label{NEW-e:local-diffusion}
\xi(\tx, \tP\cap B):=\xi^{(h,\G(\P\cap B))}(\tx,\tP\cap B):=h(M^{\G(\P\cap B)}[x]), \qquad x\in\P\cap B.
\end{equation}
Note that the function $h$ is applied to the trajectory at $x$ generated by
interacting diffusions evolving on the local graph $\G(\P \cap B)$, which
approximates $\G(\P) \cap B$ due to the stabilization property of the graph
construction. Indeed, for sufficiently large $B$, the multi-hop neighborhood of
$x$ in $\G(\P \cap B)$ coincides with that in $\G(\P) \cap B$. This allows us to
invoke Lemma~\ref{NEW-l:spat-mixing-diffusion} to establish $L^2$-stabilization
of the interacting diffusions. The details of this strategy are developed below.

We aim to prove the following \BL-localizing property 
\begin{align}
 \label{MEW-e:BLstab_sysID_pp}   
\sup_{x_1,\ldots,x_p \in \R^d } \E_{\bk x1p}\left[\Bigl|f([\xi]_1^p({\tx},\tP)) - 
f\Big( \xi(\tilde x_1,\tP\cap B_r(x_1)),\ldots ,\xi(\tilde x_p,\tP\cap B_r(x_p)) \Big)
\Bigr|\right] & \leq 2 \varphi_p(r), \quad r \geq 1,
\end{align}
satisfied for all test functions  $f \in \BL(\R^p)$, all $p \in \N$, and some  fast decreasing functions $\varphi_p$.
Here, $[\xi]_1^p({\tx},\tP)$ stands for the $p$-vector of scores at $x_1,\ldots,x_p$ in the system defined on the full graph $\G(\P)$:
\[
(\xi^{(h,\G(\P))}(\tilde{x}_1,\tP),\ldots,\xi^{(h,\G(\P))}(\tilde{x}_p,\tP))
\]
while $\xi(\tilde x_i,\tP \cap B_r(x_i))$, in the form of~\eqref{NEW-e:local-diffusion}, denotes the localized version of this score (evaluated in the  system defined on localized graphs).

By Remark \ref{rem:comparison_stabilization}, it will be enough to show the stabilization in $L^2$; i.e.,  
that 
\begin{align}
 \label{e:L2stab_sysID_pp}   
\sup_{x_1,\ldots,x_p \in \R^d} \E_{\bk x1p}[|\xi^{(h,\G(\P))}(\tilde{x}_1,\tP) - \xi^{(h,\G(\P \cap  B_r(x_1)))}(\tilde{x}_1, \tP \cap  B_r(x_1)) )|^2] & \leq 2\varphi_p(r), \quad r \geq 1,
\end{align}
for all $p \in \N$, and some fast decreasing functions $\varphi_p$.
By the Lipschitz(1) property of $h$, it is enough to show
$$
\sup_{x_1,\ldots,x_p \in \R^d} \E_{\bk x1p}
\| M^{\G(\P)}[x_1] - M^{\G(\P \cap  B_r(x_1))}[x_1]\|^2
\leq 2 \varphi_p(r), \quad r \geq 1.
$$
We may bound the required difference in the left-hand side 
by using  the spatial $L^2$-stabilization of the diffusion on graphs (see Lemma~\ref{NEW-l:spat-mixing-diffusion}) and the stabilization of the graph.

Indeed, set $k_r := r^{\gbeta}$ for a fixed $\gbeta \in (0,1)$
and let $\G_{x_1}^{(r)} := \G(\P\cap  B_r(x_1))$ and  
$\G_{x_1}^{[k_r]}$ be a subgraph of $\G(\P)$
induced by its graph ball
$\gB_{k_r}(x_1)$ of radius $k_r$ centered at $x_1$
(both graphs are  considered  for the ground process $\P$ under its  Palm conditioning $x_1,\ldots, x_p\in\P$). Note that on the event $\{ \G_{x_1}^{(r)} \supset \G_{x_1}^{[k_r]} \}$, we have that $\G_{x_1}^{[k_r]} = \G_{x_1}^{[k_r]}(\P) = \G^{[k_r]}_{x_1}(\P\cap  B_r(x_1))$. Thus via triangle inequality and intersecting with the event $\{ \G_{x_1}^{(r)} \supset \G_{x_1}^{[k_r]} \}$ as in the derivation of \eqref{e.BL-spin-stab-1*} and \eqref{e.BL-spin-stab-2*},  we have
for all $x_1,\ldots,x_p \in \R^d$
\begin{align}
& \quad  \E_{\bk x1p}
\| M^{\G(\P)}[x_1] - M^{\G(\P \cap  B_r(x_1))}[x_1]\|^2 \no \\
 & \quad \leq   \E_{\bk x1p}\big[\E [|| M^{\G(\P)}[x_1] - M^{\G(\P \cap  B_r(x_1))}[x_1]||^2 \mid \P]\1{\G_{x_1}^{(r)}\supset \G_{x_1}^{[k_r]}} \big] +2C'\Palm_{\bk x1p}\{\G_{x_1}^{(r)}\not\supset \G_{x_1}^{[k_r]}\}  \no\\
 &\quad  \le  2 \E_{\bk x1p}\big[\E [|| M^{\G(\P)}[x_1] - M^{\G_{x_1}^{[k_r]}}[x_1]||^2 \mid \P]\1{\G_{x_1}^{(r)}\supset \G_{x_1}^{[k_r]}} \big]   \no \\ 
  & \qquad + \,  2 \E_{\bk x1p}\big[\E [|| M^{\G_{x_1}^{[k_r]}}[x_1] - M^{\G(\P \cap  B_r(x_1))}[x_1]||^2 \mid \P]\1{\G_{x_1}^{(r)}\supset \G_{x_1}^{[k_r]}} \big]   +2C'\Palm_{\bk x1p}\{\G_{x_1}^{(r)}\not\supset \G_{x_1}^{[k_r]}\}  \no\\
&\quad \le  \frac{4C^{k_r}}{k_r!}+  2C'\Palm_{\bk x1p}\{\G_{x_1}^{(r)}\not\supset \G_{x_1}^{[k_r]}\}.  \label{e:L2_stab_sysID_bound}
\end{align}
Here the constant $C'$ uniformly bounds in $L^2$ all diffusion trajectories,
i.e., $\sup_{y\in\X}\E [\| M[y]\|^2]<C'$, given the bound of the initial conditions 
and Lipschitz coefficients but regardless of the graph (see  Lemma \ref{NEW-l:strong-solution-diffusion})
whereas the first term in~\eqref{e:L2_stab_sysID_bound} results from $L^2$ stabilization of the diffusion with respect to the graph distance, specifically from~\eqref{NEW-e:spat-mixing-diffusion} in   Lemma~\ref{NEW-l:spat-mixing-diffusion}.

To show that the probability bound in~\eqref{e:L2_stab_sysID_bound} is a fast-decreasing function of $r$, and thus conclude the proof, we follow the arguments used in Lemma~\ref{stablifting} (concerning interaction graphs for spin systems). Despite the slightly different context---here we localize a single score and not cluster localization of multiple scores, together with the setting of an infinite graph---the same line of reasoning remains applicable.

Specifically, we bound the probability in this expression by exploiting the following observation: if a path of length $k_r$ from $x_1$ in $\G(\P)$ is not contained within $\G_{x_1}^{(r)}$, then there must exist a point $x \in \P \cap B_{r/2}(x_1)$ with a `very large' interaction range $S(x,\P)$. This is unlikely, since the graph is stabilizing.

More precisely, following the arguments in the proof of Lemma~\ref{stablifting}, consider the events $A_r(x_1)$ (configurations of points of $\P$) such that, for all $x \in \P \cap B_{r/2}(x_1)$, the stabilization interaction range $S(x,\P)$ in $\G(\P)$ is bounded by $s_r = r^{1-\gbeta}/4$, where $\gbeta$ is as above (related to the definition of $k_r$); that is, $S(x,\P)\le s_r$. Recall that the stabilization radii $S(x,\P)$ are assumed to satisfy~\eqref{sdecayspin-infty}.
Following the same arguments as in the proof of Lemma~\ref{stablifting}, the event $\G_{x_1}^{(r)} \not\supset \G_{x_1}^{[k_r]}$ can occur only on the complement of the event $A_r(x_1)$ (with $\P$ in place of $\P_n$, corresponding to the events $A_{r,n}(x_1)$ considered there). The probability $\sP_{\bk x1p}(A_r(x_1)^c)$ admits the same bounds as those derived for $A_{r,n}(x_1)^c$ in~\eqref{e:adm_graph_euc_dist}, and thus decays rapidly in~$r$.

This completes the argument for the bound in~\eqref{e:L2_stab_sysID_bound}, and hence the proof of Theorem~\ref{NEW-t:sysID}.
\end{proof}
We conclude  this section with a few  remarks.
\begin{remark}
\label{NEW-r.diffusions}
\begin{enumerate}[wide,label=(\roman*),labelindent=0pt]
\item\label{NEW-i.B-mixing-G}
{\em Summable exponential $\B$-mixing} of points and initial marks is sufficient.
In Theorem~\ref{NEW-t:sysID}, we assume that the Brownian motions $Z_x$ are i.i.d. given $\{\delta_{(x, M(x))}\}_{x \in \P}$. Consequently, in the summable exponential $\B$-mixing assumption for $\tP$, one can omit these Brownian marks and verify the condition only for  $\{\delta_{(x, M(x))}\}_{x \in \P}$. 

\item\label{NEW-i.diffusion-vis-windows} {\em CLT for diffusions on finite windows.}
One can establish the central limit theorem  for diffusions considered only on finite windows, similarly to the spin systems discussed in Section~\ref{s:gibbsmarking}, whether or not a limiting infinite graph exists. 

\item \label{NEW-i.diffusion-graph-convergence}
{\em Convergence of interacting diffusions.} Our work on interacting diffusions was inspired by  \cite{lacker2019local}, which focused on diffusion limits over sequences of converging graphs. The typical setting there involves finite, uniformly rooted graphs converging in the local weak sense to infinite unimodular graphs, viewed as equivalence classes under graph isomorphism (e.g., sparse Erdős–Rényi graphs converging to a Bienaymé–Galton–Watson tree). A key assumption in their analysis is the uniform boundedness of the initial conditions and the Lipschitz continuity of the diffusion coefficients.
\indent Our work---of independent interest---extends this framework to interacting diffusions on finite spatial random graphs, which converge to their limiting infinite graphs in a stabilization sense (see Definition~\ref{d:stabilizing-graph}). While our main convergence notion relies on graph stabilization, weaker forms of convergence (e.g., in $L^2$) should also be considered as natural in this context.
\end{enumerate}
\end{remark}

\subsection{Auxiliary statements for diffusions on graphs}
\label{NEW-ss:proofs_sid}
In this section, we collect results concerning interacting diffusions, including their  existence, uniqueness, and the boundedness of their trajectories. These results were established in \cite{lacker2019local}, but specifically in the context of graphs considered as elements of the space $\G^*$ of equivalence classes under the isomorphism relation on locally finite, rooted graphs. This setting is too restrictive for many random graphs built on point processes. However, the key arguments developed in that work remain valid in our more general setting. For completeness and self-containment, we briefly recall these arguments.

The main result of this section is Lemma~\ref{NEW-l:spat-mixing-diffusion}, which establishes spatial $L^2$-stabilization of the diffusion on graphs. As a preparation, we require a pre-Gr\"onwall inequality and pathwise existence and uniquess for interacting diffusions.
\begin{lemma}[(Pre-)Gr\"onwall inequality for interacting diffusions]
\label{NEW-l:pre-Gron}
Consider Lipschitz diffusion  coefficients $(b,\sigma)$ as defined in Definition~\ref{NEW-d:lip_assum}. Let $\{X[x]\}_{x\in\X}$ and $\{X'[x]\}_{x\in\X}$ be two families of  continuous  processes in $\R^{d'}$ adapted to the filtration generated by the family of i.i.d. standard Brownian motions $\{Z_x(\cdot)\}_{x \in \X}$ in $\R^{d'}$ on $[0,t_0]$ 
and $\E[\|X[x]\|^2]<\infty$, $\E[\|X'[x]\|^2]<\infty$ for all  $x\in\X$.
Consider  the family of stochastic processes $\{Y[x]\}_{x\in\X}$ on $\M$, defined on $[0,t_0]$ by
\begin{equation}
\label{e:int_rep_ID}   
Y(x,t):= Y(x,0)+\int_0^t b(t, X[x], X[N_x])\, \md t + \sigma(t, X[x], X[N_x])\, \md Z_x(t),
\end{equation}
along with the family of stochastic processes $Y'(x,t)$ on $\M$ defined as above but with $X'$ instead of $X$.
These processes are indexed by the vertices $x\in\X$, they utilize the graph neighborhood $N_x=\{y\in\X:y\sim x\}$ on a given graph $\G(\X)$
and they share the same initial conditions $Y(x,0)=Y'(x,0)$.
Then for any $x\in\X$ and $t\in[0,t_0]$, we have
\begin{equation}\label{NEW-e:pre-Gron}
%\Delta[x,t] := 
\E [\| Y[x,t] - Y'[x,t] \|^2]
\le 32\max(t,1)K_t^2\int_{0}^t \max_{y\in\gB_1(x)} \E [\| X[y,s] - X'[y,s] \|^2]\, \md s,
\end{equation}
where $\gB_1(x)$  is graph ball on  $\G(\X)$ of radius $1$ centered at~$x$
and where $K_t$ is the Lipschitz diffusion constant given in Definition \ref
{NEW-d:lip_assum}.
\end{lemma}

\begin{proof}
For $x\in\X$, 
using the  triangle inequality, we have for all $t \in [0, t_0]$
\begin{align*}
& \quad \|Y[x,t] - Y'[x,t]\|^2\\ & \leq 2 \left[ \left( \int_0^t \left|b(s,X[x],X[N_x]) - b(s,X'[x],X'[N_x])\right| \, \md s \right)^2 \right. \\    
& \quad + \sup_{s \leq t} \left. \left( \int_0^s[\sigma(u,X[x],X[N_v]) - \sigma(u,X'[x],X'[N_v]) \, \md Z(x,u) \right)^2 \right]\,.
\end{align*}
By  the progressive measurability of $\sigma$ in Definition~\ref{NEW-d:lip_assum}, the finite second moments of $\|X[x]\|$ and $\|X'[x]\|$, and  the Lipschitz assumption on $\sigma$ (making $\E[\int_0^{t_0}|\sigma(u,X[x],X[N_v])|^2\,\md u<\infty$ and similarly for $\sigma(u,X'[x],X'[N_v])$)
 the second integral in the right-hand side  is a martingale; see \cite[Section 3.2(A,B)]{karatzas2014brownian}. Thus the square of the integral is a submartingale. 
Using the Cauchy-Schwarz inequality for the first integral and Doob's maximal inequality for submartingales (see \cite[Section 1.3, Theorem 3.8(iv)]{karatzas2014brownian} for example) followed by the It\^{o} isometry for the second integral, we obtain for all $x\in\X$
and $t \in [0,t_0]$
\begin{align*}
\Delta_{Y}[x,t] &:= \E [\| Y[x,t] - Y'[x,t] \|^2]
= \E [ \sup_{s \leq t}|Y(x,s) - Y'(x,s) |^2]\\
 & \leq 2 \left[ t\int_0^t \E \Bigl[\Big|b(s,X[x],X[N_x]) - b(s,X'[x],X'N_x])\Big|^2 \Bigr] \md s \right. \\
& \quad +  \left.  4\int_0^t \E \Bigl[ \Big|\sigma(s,X[x],X[N_x]) - \sigma(s,X'[x],X'[N_x])\Big|^2 \Bigr] \md s \right].
\end{align*}
Using the  Lipschitz property of $(b,\sigma)$ in Definition~\ref{NEW-d:lip_assum} (including their anticipative property) and then applying Cauchy-Schwarz to the second term, we obtain
\begin{align*}
\Delta_Y[x,t]& \leq 16\max(t,1) K_t^2 \int_0^t \Bigl( \E[\|X[x,s]-X'[x,s]\|^2] + \frac{1}{|N_x|} \sum_{y \in N_x} \E[\|X[y,s]-X'[w,s]\|^2] \Bigr) \md s. 
\end{align*}
This yields 
\textbf{}%
\begin{align*}
\Delta_Y[x,t] & \leq 16\max(t,1)K_t^2 \int_0^t \Bigl( \Delta_X[x,s] +  \frac{1}{|N_x|} \sum_{y \in N_x} \Delta_X[y,s] \Bigr) \md s \\
& \leq 32 \max(t,1) K_t^2 \int_0^t \max_{y \in\gB_1(x)} \Delta_X[y,s] \md s,
\end{align*}
where $\Delta_{X}[x,t]:= \E [\| X[x,t] - X'[x,t] \|^2]$.
This completes the proof.
\end{proof}

\begin{lemma}[Unique, path-wise  solution of the interacting diffusions and uniform $L^2$-norm bound; \protect{\cite[Theorem 3.1, Lemma 4.1]{lacker2019local}}]
\label{NEW-l:strong-solution-diffusion}
Consider the finite or locally finite, infinite graph $\G(\X)$ on $\X$ 
equipped with  initial conditions $\{M(x)\}_{x \in \X}$ in $\R^{d'}$, $d'\in\N$,  and
uniformly  bounded over  $\X$, i.e.,
 $\sup_{x\in\X}|M(x)|\le L$. 
Consider i.i.d. standard Brownian motions $\{Z_x(\cdot)\}_{x \in \X}$ in $\R^{d'}$ on $[0,t_0]$, for $t_0\in(0,\infty)$ and  
 the Lipschitz drift and diffusion  coefficients $b$, $\sigma$
 as per Definition~\ref{NEW-d:lip_assum}. Then the system~\eqref{NEW-e:sys_ID_V}, with initial conditions as above, admits a {\em unique solution with the family of processes $\{M[x]\}_{x \in \mathcal{X}}$},  taking values in the space of continuous paths $\M\owns M[x]=(M(x,t); t \in [0, t_0])$. Moreover, the trajectories $M[x]$ are uniformly bounded in $L^2$:
 \begin{equation}\label{e.M-trajectories-ubounded}
 \sup_{x\in\X}\E [\| M[x]\|^2] <C',
 \end{equation}
 where  the constant $C'$ depends only on $t_0$, the Lipschitz constant $K_{t_0}$, as well as on the bound $L$ of  the initial conditions, and {\em not} on the graph structure.
\end{lemma}
\begin{proof}
To construct a solution $M[x]$, $x\in\X$, we  consider the Picard iterations $M_m[x]$ for $m \in \N \cup \{0\}$. For $t\in[0,t_0]$ we put
\[
M_m(x,t) = M(x,0) + \int_0^t b\bigl(s, M_{m-1}[x], M_{m-1}[N_x]\bigr)\,\md s
+ \sigma\bigl(s, M_{m-1}[x], M_{m-1}[N_x]\bigr)\,\md Z_x(s),
\]
with $M_{-1}[x] := 0[x]$, that is, $M_{-1}(x,t) \equiv 0$ on $[0,t_0]$.

Starting from the integrability assumptions on the coefficients $b$ and
$\sigma$ (see the last condition in Definition~\ref{NEW-d:lip_assum}),
together with the boundedness of the initial condition $M(x)$, and using
arguments similar to those in the proof of Lemma~\ref{NEW-l:pre-Gron}, we
obtain that all successive iterations  $m \in \N$ satisfy
\(
\E\bigl[\|M_m[x]\|^2\bigr] < \infty.
\)
Applying Lemma~\ref{NEW-l:pre-Gron} to two successive iterations $M_{m}[x]$ and $M_{m-1}[x]$ 
and unfolding the integral in~\eqref{NEW-e:pre-Gron} 
we obtain
\begin{align}\no
&\E\bigl[\|M_m[x,t] - M_{m-1}[x,t]\|^2\bigr]\\
&\le (32\max(t,1)K_t^2)^m  \no
\int_0^t \int_0^{s_1} \cdots \int_0^{s_{m-1}}
\max_{y \in \gB_m(x)}
\E\bigl[\|M_0[y,s_{m}] - M_{-1}[y,s_{m}]\|^2\bigr]\,\md s_{m} \cdots \md s_2\md s_1 \\
&\le \frac{(32\max(t,1)^2K_t^2)^m}{m!} (L + K_{t_0})^2, \label{e.bound-cauchy-diffusion}
\end{align}
where the last inequality follows from the boundedness of $M(x)$ and the
integrability of $b$ and $\sigma$. 

Hence, by the Cauchy property, $M_m[x]$ converges in $L^2$ (uniformly in $x\in\X$)
to a limit process $M[x]$, i.e.,
\[
\sup_{x \in \X} \E\bigl[\|M_m[x] - M[x]\|^2\bigr] \to 0
\quad \text{as } m \to \infty.
\]
In view of the bound in~\eqref{e.bound-cauchy-diffusion}, which is independent
of the underlying graph structure, together with the assumed
uniform boundedness of $M(x)$ and the integrability of $b$ and $\sigma$,
the limit process satisfies the uniform bound~\eqref{e.M-trajectories-ubounded}.
Uniqueness follows from another application of
Lemma~\ref{NEW-l:pre-Gron}.
\end{proof}

\begin{lemma}[$L^2$-stabilization of the diffusion on graphs]
\label{NEW-l:spat-mixing-diffusion}
Under the assumptions of Lemma~\ref{NEW-l:strong-solution-diffusion}, consider
the interacting diffusions $\{M[x]\}_{x \in \X}$ on the graph $\G(\X)$ over
the time interval $t \in [0,t_0]$.
For $x\in\X$, consider a graph 
$\G^{[m]}_{x}(\X):=\G(\gB_m(x))$  induced by  $\G(\X)$ on  the ball $\gB_m(x)$ of radius~$m$ centered at~$x$ in~$\G(\X)$, and consider  
the interacting diffusions $\{M^{\G^{[m]}_{x}(\X)}[y]\}_{y\in\gB_m(x)}$ over the interval $t\in[0,t_0]$ with the same coefficients $(b,\sigma)$, generated by the same  Brownian motions  and sharing the same initial conditions $M^{\G^{[m]}_{x}}(y,0)=M(y)$ as the original diffusion $\{M[x]\}_{x \in \X}$.
Then there exists a constant $C=C(L,t_0,K_{t_0})<\infty$  (not depending on the graph $\G(\X)$) such that  
for all $t\in[0,t_0]$
\be \label{NEW-e:spat-mixing-diffusion}
\E [\| M[x,t] - M^{\G^{[m]}_{x}(\X)}[x,t] \|^2]\leq 
\frac{C^m}{{m!}}.
\ee
\end{lemma}
%The proof is deferred  to Section~\ref{NEW-ss:proofs_sid}.
The proof below in fact shows that \eqref{NEW-e:spat-mixing-diffusion} holds  if $\G^{[m]}_{x}$ is replaced by any graph $G$ such that $\G^{[m]}_{x} \subset G \subset \G(\X)$, with the inclusions denoting subgraph containment. 

\begin{proof}
We again apply the pre-Gr\"onwall inequality from
Lemma~\ref{NEW-l:pre-Gron}, this time to the difference between the two
diffusions
\[
M := \{M[x]\}_{x \in \X}
\quad \text{and} \quad
M^{\G^{[m]}} := \{M^{\G^{[m]}_{x}(\X)}[y]\}_{y \in \gB_m(x)},
\]
whose existence and uniqueness are guaranteed by
Lemma~\ref{NEW-l:strong-solution-diffusion}.
 
More precisely, we set
$X(y,t) := M(y,t)$, 
$X'(y,t) := M^{\G^{[m]}_{x}(\X)}(y,t)$ for
$t \in [0,t_0]$, $y \in \gB_m(x)$, and denote
$\Delta[y,t]:=\E [\| X[y,t] -  X'[y,t]\|^2]$.
Since $M$ and $M^{\G^{[m]}_{x}(\X)}$ are the respective solutions of the same
diffusion equation on the graphs $\G(\X)\supset \G^{[m]}_{x}(\X)$, they both
satisfy the integral representation~\eqref{e:int_rep_ID} for $y \in \gB_{m-1}(x)$
with  $Y(y,t)=X(y,t)$ and $Y'(y,t)=X'(y,t)$. Consequently, by~\eqref{NEW-e:pre-Gron}, we obtain
$$
\Delta[x,t]
\le 32\max(t,1)K_t^2\int_{0}^t \max_{y\in\gB_1(x)} \Delta[y,s] \, \md s,
$$
and iterating, 
\begin{align*}
\Delta[x,t] & \leq  (32\max(t,1)K_t^2)^m\int_0^t \int_0^{s_1} \ldots \int_0^{s_{m-1}} \max_{y \in\gB_{m}(x)} \Delta[y,s_{m}]\, \md s_{m} \ldots \md s_2 \md s_1   \\
& \leq \frac{(32\max(t,1)^{2}K_t^2)^m}{m!}
\max_{y \in\gB_{m}(x)} \Delta[y,s_{m}].
\end{align*}
By Lemma~\ref{NEW-l:strong-solution-diffusion} we have 
$\sup_{y\in\X}\E [\| M[y]\|^2]<C'$, where  the constant $C'$ depends only on $t_0$, the Lipschitz constant $K_{t_0}$, as well on the bound $L$ of  the initial conditions, and {\em not} on the graph structure. 
This implies
$\max_{y \in\gB_{m}(x)} \Delta[y,s_{m}]\le 2C'$ completing the proof
of Lemma~\ref{NEW-l:spat-mixing-diffusion}.
\end{proof}

\section{Interacting particle systems on spatial random graphs}
\label{s:applnsips}

This section uses the main theoretical results of Section \ref{s:clta} to establish the limit theory for summary statistics of 
discrete and continuous time interacting particle systems in the continuum.  
Interacting particle systems are usually studied on fixed geometries, including lattices and trees. Here we study statistics of such systems on random geometries described by  a graph on a random point process.  There are thus two sources of randomness, namely the random set of sites and the interacting particle system evolution.  Statistics of interacting particle systems are expressed as sums of real-valued marking  functions with stopping sets satisfying stabilization criteria, putting us in the set-up of Section \ref{s:clta}. This contrasts with the statistics encountered in Sections \ref{s:gibbsmarking} and \ref{NEW-s:id_sprg}, where the underlying  point processes were equipped with scores/marking functions which did not have stopping sets.

We work in a particle system framework, which is more general than that considered in most classical models. As mentioned in Section \ref{s:Intro},  the framework allows for these conditions: (i) initial states  may be dependent,  (ii) interaction neighborhoods may be unbounded,   (iii) states of particles may take values in Polish spaces and they may be a function of the entire time-evolved history of states in the interaction neighborhood and (iv) particle locations may be correlated. 

Roughly speaking, the dynamics of these general interacting particle system go as follows.  The particles of a point process having dependent initial states are equipped with independent Poisson clocks, whose `rings'
at a given site trigger an update at that site (through a specified update function) as well as at neighboring sites. The update may depend on the prior history at the site as well as on  the histories at neighboring sites, formalized below at \eqref{e:updaterulef}. 
Interactions  depend on neighbors which are not `too far away',
meaning that the neighbors should be within some random distance having a decaying tail; this is formalized via the notion of stabilizing interaction graphs given in the previous two sections.  The particle locations are required to belong to the realization of a point process
having exponentially decaying correlations and satisfying an additional growth condition on their Palm correlation functions, a generalization of
 condition~\eqref{bdedreduced} also needed when considering spin systems or interacting diffusions.

 This approach enlarges the scope of previous studies, e.g.  in \citet{Penrose2008existence},   which  assumes bounded interaction ranges, independent initial states,  and independent particle locations, namely those given by a Poisson point process.  Our set-up does not require discretization or lattice-based particle systems and  allows for models defined by
a graphical structure based on geometry,  arguably more realistic than
the Erd\"{o}s-R\'enyi graphical structure. We shall work in the sparse/thermodynamic regime, where the particles interact with only neighboring particles
and not all the particles as in mean-field models.

 In Section \ref{s:genmodelassum},  we introduce a general continuous-time interacting particle system model and state and prove limit theorems for statistics of such models in Section \ref{sss.ContPartRes}.  Markovian particle systems and  examples are discussed in Section \ref{s:exIPS}.   Discrete-time interacting particle systems, which  are of independent interest, are taken up in Section \ref{s:discIPS}.  The discrete-time model features  globally synchronous updates whereas the continuous-time model will have only locally-synchronous updates. Our  list of explicit models is 
 illustrative and not exhaustive.

\subsection{General model assumptions}
\label{s:genmodelassum}

We introduce the key ingredients in our framework for continuous-time spatial interacting particle systems. We refer to Section \ref{s:notation} for all notational formalism regarding marked point processes, including those related to Palm theory.

\subsubsection{Sites and states of particles }
\label{sss.Clocks}
We consider a simple marked point process $\tP=\{(x,U(x))\}_{x \in \P}$ on $\R^d \times \K$, where the ground process $\P$   defines the {\em sites} ({\em locations}) of the particles on $\R^d$ and their marks $U(x)=(M(x),\ttau_x)\in\K$  take values in a Polish space   having the following form
\begin{equation}\label{e:K}
\K :=  \M  \times \hat{\cN}_{[0,t_0] \times \mL  },
\end{equation}
where
\begin{itemize}
\item $ M(x) \in\M$ is  an {\em initial state} at site  $x \in \P$ at time $0$; $\M$ is a Polish space and represents for example, the occupation of the site by a particle of  some  color, infection status, etc;  $\{M(x)\}_{x \in \P}$ are possibly dependent given the ground process $\P$.
\item 
$\ttau_x\in\hat{\cN}_{[0,t_0]\times \mL}$   
is an independently {\em marked Poisson clock} at site $x\in \P$; $\ttau_x$ are 
 i.i.d. given $\{(x,M(x))\}_{x\in\P}$.
More precisely, 
$\ttau_x := \{(T_i=T_i(x), L_i=L_i(x))\}$, where $\tau_x := \{T_i(x)\}$ is a  unit rate Poisson point process on $[0,t_0]$ for some finite time horizon $t_0<\infty$---it denotes the clock whose rings trigger updates of the states of site $x$ and its neighbors in $\M$ during the $[0, t_0]$ window (cf. admissible update rules below)---and $L_i(x)$ are i.i.d. (given $\tau_x$) random elements with values in a Polish space $\mL$. These {\em time-marks} $L_i(x)$ allow for additional independent randomness (given $\tau_x$)  in the admissible update rule (for example, to define Markovian updates).  In many examples $\mL=[0,1]$ but for some
particle systems such as ballistic deposition,  it is useful to allow $\mL$ to be a space of shapes as in Section~\ref{s:exIPS}. 
To lighten the 
exposition we often  suppress the dependence on $x$ when writing  $T_i$ and $L_i$. 
\end{itemize}

In Section \ref{s:clta} we considered simple marked point processes on $\R^d \times \K$, $\K$ a Polish space. 
Taking $\K$ of the form~\eqref{e:K}, with both Polish $\M$ and $\hat{\cN}_{[0,t_0] \times \mL  }$ (the latter being  the space of finite subsets of $ [0,t_0] \times \mL  $, equipped with the weak  topology and corresponding Borel $\sigma$-algebra), 
provides an input  framework for studying the limit theory for continuous-time interacting particle systems. 

We emphasize that we neither assume that   $\P$ is a Poisson point process nor do we assume independence of the collection of initial states $\{M(x)\}_{x \in \P}$.

\subsubsection{Update rules and the evolution  of states} 
\label{sss.Apdates}
 The ringing of clocks $\tau_x$ are the times that trigger  updates of the state in $\M$ 
(for example  signalling  arrivals or departures of particles at the considered locations or updating the particles' color, infection status, etc.)
at  site  $x$, and possibly at its neighbors $y\sim x$.
Here the  relation $y\sim x$ is in terms of graph neighbors, where the graph belongs to the set of
 stabilizing interaction graphs recalled in the next  Section \ref{sss.Graph}. Also, we shall write $y \simeq x$ if $y \sim x$ or $ y = x$.
The times at which a site $x \in \P$ updates (or re-updates) its state in $\M$ are denoted by $t_j \in [0, t_0]$. The corresponding states at these times are denoted $M(x, t_j)$, with $M(x, 0) = M(x)$ representing the initial state.
The {\em history (of the evolution)} of these states visited by site $x$ is formally represented as a time-marked point process
$M(x,\cdot) :=\{(t_j, M(x, t_j))\} \in \hat{\cN}_{[0, t_0] \times \M}$,
with analogous notation used for marked clocks $\hat{\cN}_{[0, t_0] \times \mL}$.
However, it is often more intuitive to view this history  as a càdlàg function
$M(x, t) := M(x, \max\{t_j : t_j \le t\})$ for $ t \in [0, t_0]$,
i.e., a function that is right-continuous with left limits, defined over continuous time.

The process of the evolution on the states  $M(x,\cdot)=M_n(x,\cdot)$ 
will be  constructed on   {\em  finite windows $W_n$}, $n \in \N$; i.e., for locations $x\in\P_n=\P\cap W_n$, with their neighborhood  $\sim$  possibly specific for  $W_n$; see   Section \ref{sss.Graph} below. This will be constructed recursively, 
 with an  {\em updating  function} $\F$ that will be applied recurrently 
at all times $t\in\tau_x$ of all $x\in\P_n$.

We define the inductive rules of the updates of states. 
The {\em update rules are admissible} if the following assumptions are met: 
\begin{enumerate}[label=(\alph*)]
 \item   For $t=0$ we set  $M_n(x,\cdot):=M(x,\cdot)=\{(0,M(x))\}$ for all $x \in \P_n$; that is, the history of the evolution of the states  in $\M$ of all sites $x\in\P_n$ consists  only of the input initial states $M(x,0)=M(x)$.
    \item For any  marked clock event $(T,L)\in\ttau_x$ of   some site  $x\in\P_n$,
    the prior histories $\{(t,M(y,t))\}_{t<T}$ 
    at sites $y\simeq x$, i.e. at $x$ and its neighbors in $\P_n$, are expanded to
   $\{(t,M(y,t))\}_{t<T}\cup  \{(T,M(y,T))\}$, which include events $(T,M(y,T))\in \hat{\cN}_{[0,t_0] \times \M  }$  given by an updating function~$\F$ as specified below in \eqref{e:updaterulef}.
\item 
The updating  function $\F:   
\R^d\times \mL  \times \hat{\cN}_{\R^d \times \hat{\cN}_{[0,t_0] \times \M  }}\to \hat{\cN}_{\R^d\times \hat{\cN}_{[0,t_0] \times \M}}$ is measurable.
Here the argument in   $\R^d$ identifies the site $x\in\P$, $\mL$ contains  time-marks corresponding to the clock ringing  at $x$ and 
$\hat{\cN}_{\R^d \times \hat{\cN}_{[0,t_0] \times \M  }}$ denotes a marked point process on $\mathbb{R}^d$, representing the relative locations $y-x$ of the neighbors $y$ of $x$ (including $x$),
each marked by its own process of state history, which lies in $\hat{\cN}_{[0,t_0] \times \M  }$.  
The function $\F$ is applied to these arguments and returns a new point process in the same space, representing the {\em new states}: 
\begin{align} \label{e:updaterulef}
 \F\left(x,L,\Bigl\{\Bigl(y-x,\{(t,M(y,t))\}_{t <T}\Bigr)\Bigr\}_{\{y \in \P_n : y \simeq  x \}}\right) 
& =: 
\Bigl\{\Bigl(y -x,\bigl(T,M(y,T)\bigr)\Bigr)\Bigr\}_{\{y\in\P_n:y\simeq x\}}
\end{align}
{\em to be incorporated} into the updated state histories of the neighboring sites.
Sometimes, we also explicitly  assume that $\F$ is  translation-invariant, meaning 
$\F(x+z,\cdot,\{(y,\cdot)\})\equiv\F(x,\cdot,\{(y,\cdot)\})$
 for all $x,z\in\R^d$.
\end{enumerate}
Observe, for 
a marked clock event $(T,L)\in\ttau_x$ at  the site  $x\in\P_n$, besides the location of the site $x$ and the time mark $L$, the function $\F$  takes as its third  argument 
the previous histories  $M(y,t)$, $t<T$, 
at all neighboring locations $y\simeq x$ and 
its  value consists in producing the new events $(T,M(y,T))$,  which  represent  entering (possibly re-entering) these states at time $T$ to the neighboring sites.
In other words, at time $T$, the  histories are updated only for  $y\simeq x$ to $\{(t,M(y,t))\}_{t<T}\cup  \{(T,M(y,T))\}$, where $M(y,T)$ are defined in the right-hand side of~\eqref{e:updaterulef}.

The {\em complete history} of a given site $x \in \P_n$---resulting from the execution of the function $\Phi$ at all clock rings $\ttau_y$ for $y \in \P_n$---is denoted by $M_n[x] \in \hat{\cN}_{[0, t_0] \times \M}$.  These  histories
are treated as an additional marking of  the point $x \in \P_n$. Note, this marking is  measurable as a measurable mapping
\begin{equation}\label{e.g-measurability}
g : (\R^d \times \K) \times \hat{\cN}_{\R^d \times \K} \ni (\tx, \tmu) \mapsto M[x]=(M(x, t))_{t \in [0, t_0]} \in \hat{\cN}_{[0, t_0] \times \M},
\end{equation}
constructed recursively using the measurable function $\Phi$ described above. This mapping describes the history of marks $M_n[x], x \in \P_n,$ when applied simultaneously on all $x \in \P_n$. 
For more discussion  see Remark~\ref{rem:IPS}~\ref{i.update} and~\ref{i.xi-measurability} below.

%\vskip.3cm
\subsubsection{Stabilizing  interaction graphs}
\label{sss.Graph}  
Recall from Section \ref{s:notation} that $\mP_{\bk x1p}$ denotes the Palm probability distribution of the point process $\P$ conditional on the existence of  points $x_1,...,x_p$. The admissible update rules
with the update function $\F$ depend on the interaction relation specified via graphs on $\P_n$. 
As in Section~\ref{sss.graph-stabilization}, we restrict attention to stabilizing graphs---that is, graphs whose interaction ranges $S_n$ on finite windows, as defined in~\eqref{e:S_n=R_Wn} and satisfying properties~\eqref{e:cons_adm_gr}--\eqref{neighborulespin}, fulfill the fundamental assumption stated in Definition~\ref{d:stabilizing-graph},  recalled here for convenience:
\begin{equation} 
\sup_{1 \leq n < \infty}\sup_{x_1,\ldots,x_p\in W_n} \mP_{\bk x1p}
\bigl(S_n( x_1,\P_n) > s\bigr) \leq \varphi_p'(s), \quad s>0, 
\label{e:ddecay}
\end{equation}
for all $p\in\N$, where the functions $(\varphi'_p)_{p \in \N}$ are fast-decreasing
as at ~\eqref{varphibd}.

Recall these conditions imply that the expected degree of a vertex in $G(\P_n)$ is finite, i.e., they imply sparsity. As noted in Section \ref{sss.graph-stabilization} and Appendix \ref{s.admppinteraction}, the class of  stabilizing  interaction graphs includes certain proximity graphs from computational geometry.
 
\subsection{Asymptotic normality for statistics of particle systems}\label{sss.ContPartRes}
Our goal is to study certain scores associated with the history of the  marks $M_n[x]$ for all $x \in \P_n$ (describing the evolution of the updates of the states in $\M$ of $x$ during the time interval $[0,t_0]$, $t_0<\infty$). 
We restrict  to the real-valued score function
\be \label{xiU}
\xi (\tx, \tP_n) =\xi^{(h,\F,\G)} (\tx, \tP_n):=  h(M_n[x]),
\ee 
where  $h: \hat{\cN}_{[0,t_0] \times \M  } \to \R$ is measurable. 
We aim to establish Gaussian fluctuations for the summary statistics
\begin{equation}
\label{IPSfunctional}
H_n^{\xi} := \sum_{x \in \P \cap W_n} \xi (\tx,\tP_n)
\end{equation}
and the (possibly signed) random measures
\begin{equation}
\label{IPS}
\mu_n^{\xi} := \sum_{x \in \P \cap W_n} \xi (\tx,\tP_n)  \delta_{n^{-1/d}x}.
\end{equation}
For example, for a suitable choice of $h$, the statistic $H_n^{\xi}$ could either count the
total number of sites in a
 given state $M_0 \subset \M$ at time $t_0$,  the total number of sites spending
at least time $t_{\text{min}}$ in the state $M_0$, or the total time sites are in a given state
$M_0$. Section \ref{s:exIPS} provides  detailed examples. 

Finally, given $p \in [1, \infty)$, we say that the score function $\xi$ in~\eqref{xiU} satisfies the {\em $p$-moment condition on finite windows} with respect to $\tP$ if
\begin{equation}
\label{momh}
\sup_{1 \leq n < \infty} \sup_{1 \leq q \leq p} \sup_{x_1,\ldots, x_{q} \in W_n} \sE_{\bk x1{q}}[ \max(1,|\xi(\tx,\tP_n)|^p)] \leq M_p^{\xi} < \infty,
\end{equation}
where $M_p^{\xi}$ is assumed to be non-decreasing in $p$.

All applications in Section \ref{s:exIPS} showing asymptotic normality  depend on  the  following main result or its variants. Recall the notion of Palm distributions and correlations from Sections~\ref{s:prelim} and~\ref{s:appppprelim}.
\begin{theorem}[CLT for continuous time interacting particle systems] \label{thmIPS} 
Let  the input process to the interacting particle system be 
\begin{equation}\label{e:IPS-input}
\tP = \left\{(x, M(x), \ttau_x)\in\R^d \times \M \times \hat{\cN}_{[0,t_0] \times \mL}\right\}, 
\end{equation}
where $\ttau_x$ are i.i.d. and independently marked Poisson clocks, which are also independent of sites and initial states, i.e., independent of $\sum_x \delta_{(x,M(x))}$.
Assume that:
\begin{enumerate}[wide,label=(\roman*),labelindent=0pt]
\item \label{i.thmIPS}
$\tP$ has {\em summable exponential $\B$-mixing correlations}, as in Definition~\ref{def.A2}, and has  uniformly bounded correlation functions under reduced Palm distributions, i.e., for all $p \in \N$, there exists a constant $\hat{\kappa}_p$ such that for all $m \in \N$,
\begin{equation}
\label{e:bd_palmcorr}
\sup_{x_1,\ldots,x_p \in \R^d} \sup_{y_1,\ldots,y_m \in \R^d} \rho^{(m)}_{\bk x1p}(\bk y1m) \leq (\hat{\kappa}_p)^m,
\end{equation}
where $\rho^{(m)}_{\bk x1p}$ denotes the $m$-th order correlation function of $\P$ under its reduced Palm distribution $\sP^!_{\bk x1p}$.

\item \label{ii.thmIPS} {\em An admissible update rule} given by a measurable function $\F$ as in~\eqref{e:updaterulef} acts on a {\em stabilizing interaction graph} $\G$ (on finite windows) with respect to $\P$, such that the interaction range satisfies the decay assumption~\eqref{e:ddecay}.
\item \label{iii.thmIPS} A {\em score function} $\xi = \xi^{(h, \F, \G)}(\tx,\tP_n) := h(M_n[x])$ defined in~\eqref{xiU}, satisfies the {\em $p$-moment condition}~\eqref{momh} for all $p \in [1, \infty)$.
\end{enumerate}
Then the sequence of random measures $(\mu_n^{\xi})_{n \in \N}$ defined in~\eqref{IPS} satisfies the central limit theorem. That is, for every  $f \in \B(W_1)$ such that
$
\Var[\mu_n^{\xi}(f)] = \Omega(n^{\nu})$ for some $\nu > 0$,
we have
$$
(\Var[\mu_n^{\xi}(f)])^{-1/2} \big(\mu_n^{\xi}(f) - \sE[\mu_n^{\xi}(f)]\big)  \stackrel{d}{\Rightarrow} Z,$$
where $Z$ is a standard Gaussian random variable.
 \end{theorem}
In the sequel, when we say that Theorem \ref{thmIPS} holds, we implicitly mean that it is valid subject to a variance lower bound $\Var{\mu_{n}^{\xi}(f)} = \Omega(n^{\nu})$ for some $\nu > 0$, which, as always, we regard as a separate problem.

The Palm correlation bound condition \eqref{e:bd_palmcorr} implies that all factorial (and ordinary) moment measures under Palm distributions are dominated by those of a Poisson point process with intensity $\hat{\kappa}_p$. This motivates the definition of {\em sub-Poisson processes} in the sense of moment measures, as introduced in~\cite{Blaszczy14}, with determinantal processes being a prominent example. See Appendix~\ref{s.admppinteraction} for further examples.

\begin{remark}[Further comments and comparison with previous work]
\label{rem:IPS}\
\begin{enumerate}[wide,label=(\roman*),  labelindent=0pt]
\item\label{NEW-i.B-mixing-P-Z}
{\em (Summable exponential $\B$-mixing} of points and initial marks is sufficient.)
In Theorem~\ref{thmIPS}, we assumed that the Poisson marked clocks $\ttau_x\in\hat\cN_{[0,t_0]\times\mL}$  are i.i.d. given sites and initial conditions. Consequently, in the summable exponential $\B$-mixing assumption for $\tP$, one can omit these clocks  and verify the condition only for  $\{\delta_{(x, M(x))}\}_{x \in \P}$. 
\item (LLN and variance asymptotics.) \label{i.LLN-Var-IPS}
Under the assumptions of Theorem~\ref{thmIPS}, but with the $p$-moment condition~\eqref{momh} weakened to hold for some $p \in (2, \infty)$, and without any assumption on the variance rate of $\mu_n^{\xi}(f)$, suppose further that the input process $\tP$ is stationary and that both the graph $\G = \G(\P)$ and the update function $\F$ in~\eqref{e:updaterulef} are translation invariant.
Then, by Proposition~\ref{expvar}, the mean and variance asymptotics of $\mu_n^{\xi}(f)$ are given by~\eqref{expasy} and~\eqref{eqn:var}, involving Palm-distributional limits given in~\eqref{e.xil-0} and~\eqref{e.xil-xy}.
If, in addition, the interaction ranges $S_n(x, \P_n)$ are uniformly bounded in $n$ almost surely for all $x \in \P$,
i.e. satisfy $\sup_n S_n(x, \P_n)< \infty$ a.s., then the score function $\xi = \xi^{(h, \F, \G)}$ has a limit 
$$
\lim_{n \to \infty} \xi(\tx, \tP_n) = \xi_\infty(\tx, \tP)
\quad \text{a.s. for all } x \in \P.
$$
Interaction ranges for the spatial random graph examples in Appendix~\ref{s.admppinteraction} can be shown to satisfy the uniform bound required above.

Moreover, the mean and variance asymptotics admit representations involving the limit $\xi_{\infty}$:
\begin{align}
&  \Big| n^{-1} \E \mu_n^{\xi}(f) -  \rho \E_{\0} \xi_{\infty}(\tilde{\0}, \tP) \int_{W_1}f(x)\,\md x\Big| = O(n^{-1/d}),  \no \\
\label{IPSlln} & \lim_{n \to \infty} n^{-1} \Var \,\mu_n^{\xi}(f) =  \sigma^2(\xi_{\infty}) \int_{W_1} f(x)^2 \,\md x \in [0,\infty),
\end{align}
where 
\begin{align}
\label{e:xi-infty0again}
\xi_{\infty}(\tx,\tmu)&=\xi^{(h,\F,\G)}_{\infty}(\tx,\tmu)  := \lim_{n \to \infty}\xi(\tx,\tmu \cap W_n), \\
\label{sigdefagain} \sigma^2(\xi_{\infty}) & :=
\rho\sE_{\0}\xi_{\infty}^2(\tilde{\0},\tP)   \\
 &+ \int_{\R^d} \left[ \sE_{\0,x}(\xi_{\infty}(\tilde{\0}, \tP) \xi_{\infty}(\tilde{x}, \tP)) \rho^{(2)}(\0,x) - \rho^2 \, \sE_{\0}\xi_{\infty}(\tilde{\0},\tP) \, \sE_{x}\xi_{\infty}(\tilde{x},\tP)  \right] \md x \in [0, \infty). \nonumber 
\end{align}
This representation,  described in Remark~\ref{ii.xi-extension-stopping set} (Section~\ref{s:remarksltthms}), can be justified via Lemma~\ref{l.xi-infnity}. This requires deriving a uniform in $n \in \N$ bound on the stabilization radius $R^\xi_{W_n}$ of $\xi$ (see~\eqref{e:stabradips} in the proof of Theorem~\ref{thmIPS}). The bound on $R^\xi_{W_n}$ is derived using a graphical construction for interacting particle systems and controlling the diameter of the corresponding backward space-time cluster, which tracks all the nodes influencing updates at a given node. Though not straightforward, one can derive uniform bounds for the diameter by exploiting the uniform boundedness of the interaction range $S_n$ for $n \in \N$.  

\item \label{i.update}
(Local synchronous,  global asynchronous  update rules.)
The function  $\F$ in~\eqref{e:updaterulef}
 is applied successively with each ringing of some clock $\tau_x$, $x\in\P_n$.
At time $T\in\tau_x$, the application of $\F$ creates  a {\em local, synchronous update} of the history at $x$ and its neighbors $y\sim x$:   by the stopping set property~\eqref{stopspin} of $S_n$, this update is {\em local} because the  neighborhood of $x$ has a  radius bounded by  $S_n(x,\P_n)$, also, it is  {\em synchronous} because the new states at $x$ and its neighbors are evaluated 
jointly based on all histories at $x$ and its neighbors (with, if necessary, priorities assigned according to the location of $x$, its time-mark $L$ triggering this update, and the positions of the neighbors).
At the same  time $T$, all locations $y'$ which are not neighbors of $x$  are not updated and their histories do not contribute to the evaluation of the local update around $x$. This is hence  {\em asynchronous  updating occurring globally}. It is in contrast to a discrete, globally synchronous updates considered in~Section~\ref{s:discIPS}.
 
 \item (Non-Markovian updates.)
 Locations $x\in\P_n$ are typically updated multiple  times in the period $[0,t_0]$. Hence, the process of the history $M[x]=M(x,\cdot)$ may be viewed  as  {\em a pure jump type process}, but it  need not be  Markovian as in \cite[Section 2.1]{Penrose2008existence}. This is because the update rule depends upon the entire history of the process at a site and its neighboring sites. (This was also the case with   the diffusion model in Section~\ref{NEW-ss:model_thm}). Allowing the update  to depend only upon the current configuration makes the process Markovian and hence encompasses well-studied examples in the literature.  Such examples are described in detail in Section \ref{s:exIPS}.

\item (Comparison with  \citet{Penrose2008existence}.)
We compare the theorems of this section with the closely related results of \cite{Penrose2008existence}.
Theorem \ref{thmIPS} requires neither Poisson input,   independence of initial states, nor bounded interaction range.   In this way we extend upon the thermodynamic (LLN) and Gaussian limit theory of   \cite{Penrose2008existence}
 and address  some  questions raised in this article  (see especially Remark 6 in section 3 and also Section 5.2 of \cite{Penrose2008existence}). 
 Penrose's  proof of  the thermodynamic and Gaussian limits rely upon an add-one cost stabilization condition for the functionals $H_n^{\xi}$ (which is distinct from the stabilization and localization of scores described in Section \ref{s:strongstab}), discretization, and a martingale central limit theorem, an approach relying upon the existence of the infinite particle system.  By contrast,  our approach, which uses the method of cumulants,  only requires  stabilization  of the scores $\xi$ on particle systems confined to the window $W_n$, bypassing the question of existence.  We show stabilization of the scores $\xi$ via a graphical construction.  Showing existence of the limit of the finite systems to a limiting dynamic system on the infinite input $\tP$ lies outside the scope of this work though it may be possible to use the approaches from \cite{Penrose2008existence,mazumder2024existence} in combination with our derivations. Though graphical constructions have been previously used to study interacting particle systems (see \cite{Penrose2008existence,penrose2005exact} for example),  the analysis of our graphical construction is more involved due to the spatial dependencies of the underlying point processes and the unbounded interaction range. 

\item (Related literature.) With the exception of \citet{qi2008functional}, the limit theory for interacting particle systems has assumed finite range interactions.  The paper \cite{qi2008functional} proves a functional central limit theorem (FCLT) for 
spatial birth-death processes with possibly unbounded interaction range, provided the processes are obtained as a solution to stochastic equations.  The study of asymptotic Gaussianity was initiated by \citet{holley1979central},  who used classical tightness arguments to establish the FCLT for latttice based systems.  Later,  
\citet{doukhan2008functional} worked on general transitive graphs, sometimes requiring bounded degree, but still assuming finite range interactions.   More recently, \citet{onaran2022functional} establish a FCLT for spatial birth-death processes on sites given by  Poisson input. They deduce convergence of finite-dimensional distributions via  Malliavin-Stein methods from \citet{Lachieze2019normal} and establish tightness using a classical criteria from \citet{ethier2009markov}. In case of birth-death processes with diffusive dynamics, finite-dimensional convergence of local statistics has been shown in \citet{onaran2023functional}.

\item (Extensions.) \label{i.rem-extensions} Setting $\xi_i =\xi_i^{(h,\F,\G)} = h(M(x,t_i)), 1 \leq i \leq k,$ for $h : \K \to \R$
bounded and measurable and given $t_1,\ldots,t_k \in [0,\infty)$, we may obtain a multivariate central limit theorem for $(\mu_n^{\xi_i})_{i=1,\ldots,k}$ provided 
the covariances of the entries converge at some scale $n^{-\nu}$ for some positive $\nu$; see 
the discussion below Theorem ~\ref{t:multcltmarkedpp}.  As in 
Item~\ref{i.LLN-Var-IPS} above, when $\nu = 1$, one can deduce
 convergence of the covariances provided 
the interaction ranges $S_n$ of the graph $\G$ are  uniformly bounded. 
One expects that  central limit theorems could  be extended to functional central limit theorems yielding convergence of scaled versions of the process $(\sum_{x \in \P_n} h((M_n(x,t)_{t \in [0, t_0]}))$, 
thus extending Theorem 3.3 of \cite{Penrose2008existence}, which assumes independent initial marks and
bounded interaction range,  and also  extending the functional 
central limit theorem of \cite{qi2008functional} given for certain birth death processes having polynomial decay of interactions.  Furthermore, one may expect a  functional central limit theorem for the function-indexed process $(\mu_n(f))_{f \in {\cal F}}$, with $\cal F$ being a suitable function class and $t = t_0$ fixed, but to our knowledge
such a general functional limit theorem is not available even
in the simpler  setting consisting of Poisson input,  independent marks,
and unbounded interaction range.  We expect that such a  general  functional central limit theorem would 
extend those given in more specialized settings considered by~\citet{BCHS,BS}.

\item \label{i.xi-measurability}(Measurability of the score function $\xi^{(h,\F,\G)}$.) 
We justify the measurability of $g$ defined at~\eqref{e.g-measurability}.  We need to show that the map $g$ from $\R^d \times \K \times \hat{\cN}_{\R^d \times \K}$ to 
 $\hat{\cN}_{[0,t_0] \times \M  }$ given by
$(\tx,\tmu) \mapsto (M(x,t))_{t \in [0,t_0]}$ is measurable.  
Note $M[x]=(M(x,t))_{t \in [0,t_0]}$
is the final mark of $x$ in 
$\hat{\cN}_{[0,t_0] \times \M}$, measurable as a consequence of the  iteration of the measurable  marking function $\F$ given in~\eqref{e:updaterulef}, each time  being executed  based on the measurable input consisting of:
the location $x\in\P_n$, its  clock event $(T,L)\in\ttau_x$, and the set of (prior) histories $\{(y,M(y,\cdot))\}_{y\in\P_n:y\simeq x}\}\in \hat{\cN}_{\R^d \times\hat{\cN}_{[0,t_0] \times \M}}$, measurably by the induction starting from the measurable initial states $M(y,0)$. Here, the   measurability of the neighborhood $\sim$ (in the sense of the neighborhood-sets described in Section~\ref{sss.graph-stabilization})  is required in the specification of the arguments of~$\F$. 
\end{enumerate}
\end{remark}

\subsubsection{Proofs}
\label{sss.ProofsIPS}

We next prove Theorem \ref{thmIPS}.  The proof relies upon deviation inequalities for number counts of the point process $\P$ as well as for degrees and stabilization radii of spatial random graphs. These estimates will be proved at the end of section in Lemmas \ref{l:palm_conc}, \ref{l.degrees} and \ref{l.expdecay}.

\begin{proof}[Proof of Theorem \ref{thmIPS}.]
We show that $(\mu^{\xi}_n)_{n \in \N}$ satisfy
 the assumptions of Theorem \ref{t:cltmarkedpp}, with
$ \xi (\tx, \tP_n) := h(M[x])= h((M_n(x,t))_{t \in [0,t_0]})$ as at \eqref{xiU}. 
The  $p$-moment condition~\eqref{momh} of $\xi=\xi^{(h,\F,\G)}$ for all $p \in (1, \infty)$  with respect to $\P$ is a verbatim repetition of the moment condition \eqref{e:xinpmomcopied} required in this theorem. The variance lower bound holds by assumption and thus  to verify the assumptions in Theorem \ref{t:cltmarkedpp}, we need only show that  $\xi$ is fast $\BL$-localizing on finite windows as in Definition ~\ref{def.Lp-stabilizing_marking}
\ref{i.BL-localizing-windows} and~\ref{i.BL-localizing-fast}. In this regard, recall from Remark~\ref{rem:comparison_stabilization} that this is implied by fast stabilization on finite windows, as specified in Definition~\ref{def.stabilizing_marking}~\ref{i.stabilizing-windows} and \ref{i.stabilizing-fast}, which we will aim to establish.

 Our proof is based on the {\em study of  backward clusters of `information' propagation.  For a site $x$  this cluster} 
  consists of the set of sites in~$\P_n$ from which its history $M_n(x,t)$ can be influenced at any   time $t\in [0,t_0]$. This is introduced in {\em Step 1}. Next in {\em Step 2}, we {\em control  the diameter of the backward clusters}, i.e., we upper-bound the (appropriate Palm) probability of the event that the state at $x$
  is influenced by sites outside  the ball $B_m(x)$ of radius~$m\ge1$. We show that this probability is fast-decreasing to zero as a function  of $m\to\infty$, in the sense of~\eqref{varphibd}. This is 
essentially because there are  not enough  `opportunities', i.e. not enough  ordered events with respect to clock rings
 on paths from $B^c_m(x)$ to $x$, especially taking into account  
probabilistic constraints on edge lengths and on degrees of sites 
 in our stabilizing  interaction graph $\G$ on  $\P$ satisfying decay
 bounds \eqref{e:ddecay}. Finally in {\em Step 3}, we add {\em a buffer stabilizing region} around $B_m(x)$ to account for stabilization of graph edges making $\xi$ fast stabilizing at $x$. The proof details are as follows.

\vskip.3cm

\noindent{\em Step 1. Definition of backward clusters.}
Given the graph $G(\P_n)$, $n \in \N$, consider the `super-clocks' at locations $x\in\P_n$, 
$$ \hat{\tau}_{x}:= \hat{\tau}_{n,x}= \tau_{x} \cup \bigcup_{y \in \P_n:y \sim x} = \bigcup_{y \in \P_n:y\simeq x} \tau_y \subset [0,t_0], $$  
i.e., the union of the clock at $x$ and all neighboring clocks,
with the neighbor relation respecting
\eqref{e:cons_adm_gr}.
Since the (original) Poisson clocks are i.i.d.,  we may assume  that a.s.  the clocks do not `ring' together,  i.e.,  $\tau_x \cap \tau_y = \emptyset$ for all  $x, y \in \P_n$, $x\not=y$.   Observe that  $M({x},\cdot)$ is updated only at ringing times of the {\em `super-clock'} $\hat\tau_{x}$,  i.e., the history  $M(x,\cdot)$ has events $(t,M(x,t))$ only for $t \in \hat{\tau}_{{x}}$ created according to the update function \eqref{e:updaterulef}, apart from its own history $\{M({x},s)\}_{s<t}$  and influenced by the prior  histories $\{M(y',s)\}_{s <t}$ of the sites $y'$ whose `super-clocks' ring at the same time $t\in\hat\tau_{y'}$. 

In order to track sites (the points of $x\in\P_n$) and times of their `super-clocks' $\hat\tau_x$ which mutually influence 
their updates we define
an {\em  oriented graph} $\hat{\G}(\P_n)$
as follows: The vertices of the oriented graph $\hat{\G}(\P_n)$ are $V=V_n :=\{(x,t) : x \in \P_n, t \in \hat{\tau}_x \cup \{0\} \}$ and  the oriented edges are $(x,t) \rar (x',t')$ if $t < t'$ and $(x,t') \in V$ (i.e., $t'\in \hat{\tau}_x$). Since the clocks $\tau_x,  x \in \P_n$ have no common `rings',   the condition that $(x',t'),  (x,t') \in V$ is equivalent to either (i) $t' \in \tau_x \cup \tau_{x'}$ and $x \sim x'$ or (ii) $t' \in \tau_{x''}$ for some $x''$ where $x'' \sim x'$ and $x'' \sim x$.  In the first case,  the clock at $x$ or $x'$ rings at time $t'$ and because of update rule \eqref{e:updaterulef},   the update at $x$ or $x'$ depends on previous updates at neighboring points including the other point.  In the second case,  the clock at a mutual neighbor $x''$ rings at time $t'$ and the update at $x$ and $x'$ is influenced by previous updates at neighboring points including the point $x''$
 due to the update rule \eqref{e:updaterulef}.  Thus the oriented graph captures space-time dependency between updates. We shall now make this more precise. 

Define $C'(x,t)=C_n'(x,t)$, $(x,t)\in V=V_n$,
 to be the {\em backward (in time) cluster} of $(x,t)$,  
i.e., $(x',t') \in C'(x,t)$ if there exists a path in $\hat{\G}(\P_n)$  from $(x',t')$ to $ (x,t)$.
The events at times $s$ prior to   $t'$ of the histories of  sites $x'\in\P_n$ given by
$$
\{((s,M(x',s)))_{s<t'} : \ (x',t') \in C'(x,t)\}
$$
are precisely those which  contribute to the evaluation of the   value of the updated $(t,M(x,t))\in M(x,\cdot)$ by the function~\eqref{e:updaterulef}, given the initial states $M(x',0)$.

By cylinder of radius $r$ with axis through $x \in \R^d$ we mean the set $B_r(x) \times [0, t_0]$.
Let $\hat{t}_{{x}}$ be the last time of the `super-clock' $\hat{\tau}_{{x}}$.
The radius of stabilization $R_{W_n}^\xi(\tx,\tP_n)$ of  $\xi=\xi^{(h,\F,\G)}$ as in Definition~\ref{def.Stab.Radius} may be upper bounded 
by the radius of the 
smallest cylinder
containing the backwards cluster $C'({x}, \hat{t}_{{x}})=C_n'({x}, \hat{t}_{{x}})$, centered at $x\in\P_n$, plus the maximum of the interaction ranges $S_n(x',\P_n)$ (defined in \eqref{e:S_n=R_Wn})
among 
points $x'\in\P_n$ in the basis of this cylinder, i.e.,  
\begin{equation}
\label{e:stabradips}
R_{W_n}^{\xi}( \tx, \tP_n) \leq \max_{(x',t')  \in C_n'({x}, \hat{t}_{{x}})}|x' - {x}| + 
\max_{(x',t')  \in C_n'({x}, \hat{t}_{{x}})} S_n(x', \P_n).
\end{equation}
Indeed, the histories $(s,M(x',s))_{s\le t'}$
for  $(x',t') \in C_n'({x}, \hat{t}_{{x}})$
are sufficient to construct the complete history $(s,M(x,s))_{s\le t_0}$ of $x\in\P_n$,
and the second term in~\eqref{e:stabradips} corresponds to a buffer region of $C_n'({x}, \hat{t}_{{x}})$, such that the addition of points outside this 
buffer region will not affect 
the  backwards cluster $C_n'({x}, \hat{t}_{{x}})$, and hence will not affect the score at ${x}$.

The  two terms in the right-hand-side of~\eqref{e:stabradips}  may be controlled using the  Palm correlation functions of the ground process~$\P$ as in~\eqref{e:bd_palmcorr}, the fast  decay of the graph stabilization radius  $S_n$ as in~\eqref{e:ddecay}, and the `asymptotic independence' of `super-clocks'.  This goes as follows: We shall bound the first term on the right-hand side of \eqref{e:stabradips} in Step~2, which forms the bulk of the remaining proof.
In Step 3, using  Lemma~\ref{l.expdecay} at the end of this section, we 
bound the second term on the right-hand side of \eqref{e:stabradips}. The requisite deviation bounds and degree bounds are derived at the end of the section in Lemmas \ref{l:palm_conc} and \ref{l.degrees}.

\vskip.3cm

\noindent{\em Step 2. Controlling the diameter of the backward clusters.} In this step we control the growth of 
\begin{equation}
 \label{e:diambackwardcluster}   
D(\tx, \tP_n):= \max_{(x',t') \in C_n'({x}, \hat{t}_{{x}})} |x' - {x}|.
\end{equation}
Fix $\alpha\in(0,1)$ and
$\beta \in (0,(1-\alpha)/(d+1))$ (where $d$ is the dimension of $\R^d$) and for $m \in \N$ set $r_m := m^{\beta}$, $k_m:=m^\alpha$. For $x\in \P_n$,  consider the 
event $A_m(x) :=A_{n,m}(x)$ such that there is a site $y\in \P_n\cap B_{m}({x})$  such that either 
(i) its degree  in the graph $\G(\P_n)$, denoted $\Deg(y, \P_n)$, is at least $k_m$, i.e., $\Deg(y,\P_n)\ge k_m$,  {\em or} (ii) one of the edges incident to it  has  length at least~$r_m$.
Conditioning on $A_m(x_1)$ and its complement $A_m^c(x_1)$ we get
\begin{align}
\label{e:radstab_IPS_bd1}
& \quad \sP_{{x_1},\ldots,x_p}(D(\tx_1,\tP_n) \geq m) 
 \leq \sP_{\bk x1p}(A_m(x_1)) + \sP_{\bk x1p}(\{D(\tx_1,\tP_n) \geq m \} \cap A_m^c(x_1)).
\end{align}
For the event $A_m(x_1)$ to occur under $\sP_{{x_1},\ldots,x_p}$, in view of the  rule for the existence of edges \eqref{stopspin},  there must be a site $x \in \P_n\cap B_{m}({x_1})$ such that 
$S_n(x, \P_n) \geq r_m$ or $\Deg(x, \P_n) \geq k_m$.  Applying the union bound,  the Campbell-Little-Mecke formula \eqref{e:CLM2} and the bound \eqref{e:bd_palmcorr} with $m = 1$, we have  
\begin{align}
\sP_{\bk x1p}(A_m(x_1))  & \leq \int_{B_{m}({x})} \sP_{x,\bk x1p}\bigl(\,S_n(x, \P) \geq r_m \text{\;or\;} \Deg(x, \P_n) \geq k_m\,\bigr) \rho^{(1)}_{\bk x1p}(x) \md x\no \\
&\qquad 
\no \\
&\qquad +\sum_{i=1}^p \sP_{\bk x1p}\Bigl((  S_n(x_i,\P_n) > r_m) +
 p\sP_{\bk x1p}( \Deg(x_i, \P_n) \geq k_m)\Bigr)
\no \\[1.5ex]
& \leq \theta_d \hat{\kappa}_p m ^d \bigl( \varphi'_{p + 1}(r_m)+
\varphi''_{p + 1}(k_m)
\bigr)+p\bigl(\varphi'_p(r_m)+ \varphi''_p(k_m)\bigr)\,,
 \label{bdAm}
 \end{align}
where $\varphi_p'$  is at \eqref{e:ddecay} and where
$\varphi_p''$, defined at \eqref{ddecay}, controls the degree bound
given in
Lemma~\ref{l.degrees} at the end of this section. Recalling that $r_m = m^{\beta}$ and $k_m = m^{\alpha}$ we find that for any $\alpha, \beta > 0$ the expression in~\eqref{bdAm} is a fast decreasing function of $m\in\N$ and so is the first term on the right-hand side of \eqref{e:radstab_IPS_bd1}.

We now bound the second term in \eqref{e:radstab_IPS_bd1}, as a fast decreasing function of $m$.  We define the  sequence  $l_m :=\lfloor m^{1-\beta}/(2(p+1)) \rfloor$, $m\in\N$. 
Observe, 
on $A_m^c(x_1)$, if  
$D(\tx_1,\tP_n) \geq m$ then    there exist {\em distinct} sites $y_1,\ldots,y_{l_m}\in\P_n\cap B_m(x_1)$,
{\em different from} $x_1,\ldots,x_p$,
such that $|y_{i}-y_{i-1}|<2(p+1)r_m$, 
$i=1,\ldots,l_m$, where we set $y_0:=x_1$, whose super-clocks  (when they update their states) ring  at some  strictly positive times in $\hat\tau_{y_i}$, with the exception of the earliest time being possibly the initial time zero. 

Specifically, looking backwards in time
from $t_0$, we  enumerate these nodes $y_i$ and their corresponding updating times in 
reverse chronological order 
$ t_0\ge \hat T_1>\hat T_2>\ldots>\hat T_{l_m}$ 
with  $\hat T_i\in \hat\tau_{y_i}$ for $i=1,\ldots,l_{m-1}$
and $\hat T_{l_m}\in\hat\tau_{y_i}\cup\{0\}$.

Note, the degrees at all sites $y_i$ are bounded  by $\Deg(y_i,\P_n)< k_m$.
Indeed,  when all edges of  $\G(\P_n)$ in $B_m(x_1)$ have length smaller than  $r_m$ (complement  of the condition (ii) of the event $A_m(x_1)$), then the corresponding  oriented edges in $\hat\G(\P_m)$
(with chronologically oriented times of clocks), when projected on $\R^d$, are shorter than  $2r_m$, and hence,  one  needs at least $l'_m:=\lfloor m/(2r_m)\rfloor $ distinct sites 
$\{y'_i:i=1,\ldots,l'_m\}$ in $\P_n$ to cover the distance from $x_0$ to $B^c_m(x_1)$ on $\hat\G(\P_n)$. 
Among these sites $y'_i$ there might be fixed atoms $x_1,\ldots,x_p$ under $\sP_{\bk x1p}$;  removing them one  obtains the considered path $\{y_i:i=1,\ldots,l_m\}$. 
Clearly $\Deg(y_i,\P_n)< k_m$ results from the complement of  the condition (i) of $A_m(x_1)$.

Considering the path $\{y_i:i=1,\ldots,l_m,\}$, using the Markov inequality and the Campbell-Little-Mecke formula one  bounds the considered probability by
\begin{align}
& \quad \sP_{\bk x1p}\bigl(\,\{D(\tx_1,\tP_n) \ge m \} \cap A_m^c(x_1)\bigr) \no \\
& \leq \sE_{\bk x1p}\left[\sum_{
\genfrac{}{}{0pt}{2} {(y_1,\ldots,y_{l_m})\in(\P_n\setminus\{x_1,\ldots,x_p\})^{(l_m)}}
{\hat T_i\in\hat{\tau}_{{y_i}}: \,i=1,\ldots,l_m}}
\prod_{i=1}^{l_m}\1{|y_{i}-y_{i-1}|<2(p+1)r_m
}\1{\hT_{i}>\hT_{i+1}}\1{\Deg(y_i,\P_n)<k_m}\right] \no \\
& \leq  \int_{\bk y1{l_m} \in (\mR^d)^{l_m}} \rho^{(l_m)}_{\bk x1p}(\bk y1{l_m}) \prod_{i=1}^{l_m}\1{|y_{i}-y_{i-1}<2(p+1)r_m}
    \no \\
& \hskip3cm \times \sE_{\bk x1p \sqcup\bk y1{l_m}}
\Bigl[\sum_{ \hat T_i\in\hat{\tau}_{{y_i}}: \,i=1,\ldots,l_m}
\prod_{i=1}^{l_m}\1{\hT_{i}> \hT_{i+1}}\1{\Deg(y_i,\P_n)<k_m}\Bigr]\,\md \bk y1{l_m} \,\label{e:bd_stab_IPS}
\end{align}
where, for the convenience, we set $y_0:=x_1$ and $\hat T_{l_{m}+1}:=-1$. Note the super-clocks  $\hat\tau_{y_i}$ are, in general, not independent since the sets of neighbors of sites $y_i$ may not be disjoint. However, under $\sP_{\bk x1p \sqcup\bk y1{l_m}}$ with  $\Deg(y_i,\P_n)<k_m$, these aggregated clocks are still homogeneous Poisson point processes on $[0,t_0]$ of intensities  bounded by $k_m$,
 with the last (index $l_m$)  appended by the initial time event $t=0$. 
First  assuming  $\hat T_{l_m}>0$, using the strong Markov property for the  aggregated clocks  one gets 
\begin{align*}
&\sE_{\bk x1p \sqcup\bk y1{l_m}}
\Bigl[\sum_{ \hat T_i\in\hat{\tau}_{{y_i}}: \,i=1,\ldots,l_m}
\prod_{i=1}^{l_m}\1{\hT_{i}> \hT_{i+1}}\1{D(y_i,\P_n)<k_m}\Bigr]\\
&\le k_m^{l_m}\int_0^{t_0}\md s_{l_m}\int_{s_{l_m}}^{t_0}\md s_{l_m-1}\dots
\int_{s_3}^{t_0}\md s_{2}\int_{s_2}^{t_0}\md s_{1}
\\
&=\frac{k_m^{l_m}t_0^{l_m}}{l_m!}.
\end{align*}
Considering separately $\hat T_{l_m}=0$
one obtains similarly that the expected number of the considered paths is equal to  $k_m^{l_m-1}t_0^{l_m-1}/(l_m-1)!$.
Combining and substituting into~\eqref{e:bd_stab_IPS}
and using~\eqref{e:bd_palmcorr} one obtains 
\begin{align}
 \sP_{\bk x1p}\bigl(\,\{D(\tx_1,\tP_n) \ge m \} \cap A_m^c(x_1)\bigr) 
&\le\hat\kappa_p^{l_m}\bigl(\theta_d(2(p+1)r_m)^d\bigr)^{l_m}t_0^{l_m}\bigl(1+\frac{l_m}{k_mt_0}\bigr)\frac{k_m^{l_m}}{l_m!}.\label{e.but-last-Ac}
\end{align}
Recall 
$r_m=m^\beta$, 
$l_m =\lfloor m^{1-\beta}/(2(p+1)) \rfloor$,
and $k_m=m^\alpha$, with $\alpha\in(0,1)$ and
$\beta \in (0,(1-\alpha)/(d+1))$.
Using  a crude Stirling bound $l_m!\ge Ce^{\frac{1}{2}l_m\log l_m}$,
with $C > 0$ some small constant, and applying it to ~\eqref{e.but-last-Ac}, 
one obtains  the following bounds, with some further  positive, finite  constants $C_1,C_2$ (not dependent on $n,m$ and $\alpha,\beta$)
\begin{align}\no
\sP_{\bk x1p}\bigl(\,\{D(\tx_1,\tP_n) \ge m \} \cap A_m^c(x_1)\bigr) 
&\le C_1\bigl(1+\frac{l_m}{k_m}\bigr)
\exp\Bigl(C_2l_m+ \frac{1}{2} l_m\log\frac{r_m^dk_m}{l_m}\Bigr)\\
&\le C_1\bigl(1+m^{1-\beta-\alpha}\bigr)
\exp\Bigl(m^{1-\beta}(C_2+\frac{1}{2} \log m^{\beta(d+1)+\alpha-1})\Bigr).
\label{e.final-D>m}
\end{align}
Observe  $\beta(d+1)+\alpha-1<0$, whence  the last bound is a fast-decreasing function of $m$. This shows that the second term in \eqref{e:radstab_IPS_bd1} is fast-decreasing in $m$ and together with the same statement for  the bound at ~\eqref{bdAm}, this yields that the random variable
$D(\tx_0,\tP_n)$ has a fast decreasing tail.

\vskip.3cm
\noindent{\em Step 3. Controlling the width of the buffer region.} Here we wish to control the term  
$$
{\cal W}(\tx, \tP_n):= \max_{(x',t')  \in C'_n(x, \hat{t}_{x})} S_n(x', \P_n)
$$
in~\eqref{e:stabradips}. Notice that by \eqref{e:diambackwardcluster},
$$
{\cal W}(\tx, \tP_n) \leq \sup_{|x' - x| \leq D(\tx,\tP_n)} S_n(x', \P_n).
$$
Since  $D(\tx,\tP_n)$ and $S_n(x, \P_n)$ are fast decreasing, it follows by  
Lemma~\ref{l.expdecay} with 
$Q_1$ set to  $D$
and $Q_2$ set to $S_n$, 
that ${\cal W}(\tx, \tP_n)$ is fast decreasing. 

Since $R_{W_n}^\xi(\tx, \tP_n) \leq D(\tx,\tP_n) +
{\cal W}(\tx, \tP_n)$, where both $D(\tx,\tP_n)$ and ${\cal W}(\tx, \tP_n)$ are fast decreasing, it 
follows that $R_{W_n}^\xi(\tx, \tP_n)$ is fast decreasing on finite windows as in Definition \ref{def.stabilizing_marking}~(ii),~(iii). 
$ $\par\nobreak\ignorespaces %This forces the Halmos symbol to appear at the end of the line (when there is no space left in the previous one).
\end{proof}
We prove three lemmas used in the proof of Theorem \ref{thmIPS}. 

\begin{lemma}[Deviation inequality for number counts]
\label{l:palm_conc} 
Let $\P$ be a simple point process satisfying the Palm correlation bound \eqref{e:bd_palmcorr}. Then for a bounded Borel set $B \subset \R^d$,
\begin{align}\no \sup_{x_1,\ldots,x_p \in \R^d} \sP_{\bk x1p} \bigl(\P(B) \geq \hat{\kappa}_p|B| + p + t \bigr) \le
& \exp\Bigl(-t\Bigl(\log(1+\frac{t}{\hat{\kappa}_p|B|})-1\Bigr)\Bigr), \ \ \ p \in \N, \ \ \ t > 0.\label{e:decay-number-points}
\end{align}

\end{lemma}
\begin{proof}
By the relation between Palm and reduced Palm, we have 
\[ \sP_{\bk x1p} \bigl(\P(B) \geq \hat{\kappa}_p|B| + p + t \bigr) \leq \sP^{!}_{\bk x1p} \bigl(\P(B) \geq \hat{\kappa}_p|B| + t \bigr), \]
and so we shall bound the latter probability. To do so, we bound $\P(B)$ by  a
Poisson random variable $N(B)$ with mean $\hat{\kappa}_p|B|$ 
whose (factorial, hence non-factorial) exponential moments are larger; cf.~\cite{Blaszczy14,Blaszczyszyn15}. Indeed, using definition of Palm correlation functions and \eqref{e:bd_palmcorr}, we have that
$$ \sE^{!}_{\bk x1p}\P(B)^{(k)} = \int_B^k \rho_{\bk x1p}(\bk y1k) \, \md y_1 \ldots \md y_k \leq \bigl( \hat{\kappa}_p|B| \bigr)^k = \sE N(B)^{(k)}$$
and hence for $u \geq 0$ we have $\E^{!}_{\bk x1p}[e^{u\P(B)}] \leq \E[e^{uN(B)}]$. The rest of the proof now follows from the classical Chernoff bounds for a Poisson random variable:  
\begin{align*}
  \sP^{!}_{\bk x1p} \bigl(\P(B) \geq \hat{\kappa}_p|B| + t \bigr)&\le \exp\Bigl(-\sup_{u\ge0}
  \Bigl\{u(\hat{\kappa}_p|B|+t)-\log\E^{!}_{\bk x1p}[e^{u\P(B)}]\Bigr\}\Bigr)\\
&\le \exp\Bigl(-\sup_{u\ge0}
  \Bigl\{u(\hat{\kappa}_p|B|+t)-\log\E[e^{uN\dy{(B)}}]\Bigr\}\Bigr) \\
%&= \exp\Bigl[-\kappa_p|B|h\Bigl(\frac{\hat{\kappa}_p|B|+t}{\hat\kappa_p|B|}\Bigr)\Bigr]\\
&=\exp\Bigl(-(\hat\kappa_p|B|+t)\log\Bigl(\frac{\hat{\kappa}_p|B|+t}{\hat\kappa_p|B|}\Bigl)+t\Bigr)\\
&=\exp\Bigl(-t(\log\Bigl(\frac{\hat{\kappa}_p|B|+t}{\hat\kappa_p|B|}\Bigl)-1)-
\hat\kappa_p|B|\log\Bigl(\frac{\hat{\kappa}_p|B|+t}{\hat\kappa_p|B|}\Bigl)\Bigr)\\
&\le\exp\Bigl(-t(\log\Bigl(\frac{\hat{\kappa}_p|B|+t}{\hat\kappa_p|B|}\Bigl)-1)\Bigr),
\end{align*}
where we used in the first  equality $\sup_{u\ge0}(ua - \log\E[ e^{uN(B)}])= a \log(a/\E[N(B)]) - a + \E[N(B)]$ (see e.g. \cite[Example 1.10]{draief2010epidemics}). 
\end{proof}

\begin{lemma}[Degree bounds in stabilizing interaction graphs]\label{l.degrees} 
Let $\P$  be a simple point process on $\R^d$
satisfying Palm correlation bound~\eqref{e:bd_palmcorr} and let $\G$ be a stabilizing interaction  graph defined on the finite  windows of $\R^d$, with
$\G(\P_n)$ satisfying~\eqref{e:cons_adm_gr}-\eqref{sdecayspin}.
Denote by $\Deg(x,\P_n)$ the degree of the site $x$ in $\G(\P_n)$.
Then for all  
$p \in \N$ there are fast decreasing functions 	$\varphi_p''$  such that 
\begin{equation}
\sup_{1 \leq n < \infty}\sup_{x_1,\ldots,x_p\in W_n} \mP_{\bk x1p} \bigl(\Deg(x_1,\P_n) > k\bigr) \leq \varphi_p''(k), \quad k \in \N. \label{ddecay}
\end{equation}
\end{lemma}
\begin{proof}
Fix $p \in \N$ and assume without loss of generality that $k \geq 2(p+1).$
Given such a $k \geq 1$, choose such that $\hat{\kappa}_p\theta_dr^d + p = k/2$ and set $t = k/2$. By stabilization of the graph, we have that if $S_n(x_1,\P_n) \leq r$ then $\Deg(x_1,\P_n) \leq \P(B_r(x_1))$. By the choice of $t,r$,  we have 
\[ \mP_{\bk x1p} \bigl(\Deg(x_1,\P_n) > k \bigr) \leq \sP_{\bk x1p} \bigl(\P(B_r(x_1)) \geq \hat{\kappa}_p\theta_dr^d + p + t \bigr) + \sP_{\bk x1p}(S_n(x_1,\P_n) \geq r) .\]
Using again the choice of $t,r$, 
the proof follows from the fast stabilization of $S_n$ as at \eqref{e:ddecay} and the deviation inequality Lemma \ref{l:palm_conc} applied to $B = B_r(x_1)$. 
\end{proof}

 \begin{lemma}[Tail bounds for the maximum of fast decreasing marks]\label{l.expdecay} 
 Let $Q_1, Q_2$ be two marking functions on finite marked point processes. Let $\tP$ be a marked point process such that the ground process $\P$ satisfies the Palm correlation bound~\eqref{e:bd_palmcorr} for $m =1$.
 Suppose for $p\in\N$, $i=1,2$ that there are   fast decreasing functions $\varphi_{i,p}$ as in \eqref{varphibd} such that 
\begin{equation}\label{e:Q12tails}
\sup_{1 \leq n < \infty} \sup_{x_1,\ldots,x_p \in W_n} \mP_{\bk x1p}
\bigl(Q_i(\tx_1, \tP_n)>r\bigr) \le \varphi_{i,p}(r),\quad r > 0.
\end{equation}
 Then 
for all $p \in \N$ there are fast decreasing functions ~$\tilde\varphi_p$ such that
$$
\sup_{1 \leq n < \infty} \sup_{x_1,\ldots,x_p\in W_n}\sP_{\bk x1p}\left( \max_{x \in \P_n :  |x - x_1| \leq Q_1(\tx_1, \tP_n)} Q_2(\tx, \tP_n) > r \right) \leq   \tilde{\varphi}_p(r), \quad r > 0.
$$
If the assumption \eqref{e:Q12tails} holds with the supremum taken up to $n = \infty$, i.e., $\sup_{1 \le n \le \infty}$, then the conclusion  remains valid with  $\sup_{1 \le n \le \infty}$.
\end{lemma}

\begin{proof}
Assume that $Q_i(x, \tP_n),  i  \in \{1,2\}$ satisfy \eqref{e:Q12tails}. Recall that $\theta_d$ is  the volume of the $d$-dimensional unit ball. Thus using the union bound and the Campbell-Little-Mecke formula, we obtain 
\begin{align}
& \quad \sP_{\bk x1p}\bigl(\max_{x \in \P_n :|x - x_1| \leq Q_1(\tx_1, \tP_n)} Q_2(\tx, \tP_n) > r\bigr) \\
& \leq  \E_{\bk x1p} \left[ \sum_{x \in \P \cap  B_{Q_1(\tx_1, \tP_n)}(x_1) } { \bf{1}}(Q_2(\tx, \tP_n) > r) \right] \nonumber\\
& \le \sum_{x\in\bk x1p}\E_{\bk x1p}[ {\bf{1}}(Q_2(\tx, \tP_n) > r)] \no\\
& \quad +
\E_{\bk x1p} \left[ \sum_{x \in (\P\setminus \bk x1p) \cap B_{Q_1(\tx_1, \tP_n)}(x_1) } { \bf{1}}(Q_2(\tx, \tP_n) > r)   { \bf{1}}(Q_1(\tx_1, \tP_n) \in [0,1) )\right]  \nonumber \\
& \quad + \E_{\bk x1p} \left[ \sum_{j = 1}^{\infty} \sum_{x \in (\P\setminus \bk x1p)\cap B_{Q_1(\tx_1, \tP_n)}(x_1) } { \bf{1}}(Q_2(\tx, \tP_n) > r)   { \bf{1}}(Q_1(\tx_1, \tP_n) \in [2^{j - 1}, 2^j) )\right] \nonumber\\
&  \le \sum_{x\in\bk x1p}\E_{\bk x1p}[ {\bf{1}}(Q_2(\tx, \tP_n) > r)] \no \\
&\quad + \E_{\bk x1p} \left[ \sum_{x \in \P \cap  B_{Q_1(\tx_1, \tP_n)}(x_1) } { \bf{1}}(Q_2(\tx, \tP_n) > r) \right]  \nonumber \\
& \quad + \sum_{j = 1}^{\infty}\E_{\bk x1p} \left[ \sum_{x \in \P_n \cap B_{2^{j}}(x_1) } { \bf{1}}(Q_2(\tx, \tP_n) > r)   { \bf{1}}(Q_1(\tx_1, \tP_n) \geq  2^{j - 1}) \right] \nonumber \\
& \le  p\varphi_{2,p}(r) + \int_{B_{1}(x_1) } 
 \sP_{x,\bk x1p}(Q_2(\tx, \tP_n) > r) \rho^{(1)}_{\bk x1p}(x) \,\md x  \nonumber \\
&  \quad +  \sum_{j = 1}^\infty \int_{B_{2^{j}}(x_1) }
 \sP_{x,\bk x1p}(Q_2(\tx, \tP_n) > r, Q_1(\tx_1, \tP_n) \geq  2^{j - 1}) \rho^{(1)}_{\bk x1p}(x) \,\md x  \label{e:CLM+sum}
\end{align}
with  $\sP_{x,\bk x1p}=\sP_{x,x_1,\ldots,x_p}$, and where the first  term in the last bound  is a bound for  the
 sum of the probabilities 
 $$
 \sP_{\bk x1p}\bigl( Q_2(\tx_j, \tP_n) > r\bigr)\le \varphi_{2,p}(r)
 $$ of the tails of the $Q_2$ marks of the atoms $x_j$, $j=1,\ldots p$, 
 some of which  might be in $B_{2^{j}}(x_1)$ but  not counted in the 
Campbell-Little-Mecke integral formula involving the intensity  $\rho^{(1)}_{\bk x1p}(x)$ of the reduced Palm $\sP^!_{\bk x1p}$.
This in turn is a fast decreasing function,  as $\varphi_{2,p}$ is fast decreasing as well.
 
The sum of the  integral formulas $ \sum_{j = 1}^\infty \int_{B_{2^{j}}(x_1) }\sP_{x,\bk x1p}(\dots)\rho^{(1)}_{\bk x1p}(x) \,\md x$, using  the Cauchy-Schwartz inequality and  the condition \eqref{e:bd_palmcorr}, can be further  bounded  by 
\begin{align*} 
& \sum_{j = 0}^{\infty} \int_{B_{2^j}(x_1) }
\sP_{x,\bk x1p}(Q_2(\tx, \tP_n) > r)^{1/2} \sP_{x,\bk x1p} (Q_1(\tx_1, \tP_n) \geq  2^{j - 1})^{1/2} \rho^{(1)}_{\bk x1p}(x) \,\md x\\
 & \leq  \hat\kappa_p( \varphi_{2, p + 1}(r))^{1/2}
\sum_{j = 0}^{\infty} \theta_d 2^{jd} (\varphi_{1, p + 1}( 2^{j - 1}))^{1/2}.
\end{align*}
Since $\varphi_{1, p + 1}$ is fast decreasing, the last sum is finite and so the right-hand side in the last line is also fast decreasing as $\varphi_{2, p + 1}$ is fast decreasing as well.  

The proof of the version of the result  with  $\sup_{1\le n \le \infty}$  replacing $\sup_{1\le n < \infty}$ follows the same lines. 
$ $\par\nobreak\ignorespaces %This forces the Halmos symbol to appear at the end of the line (when there is no space left in the previous one). 
\end{proof}

\subsection{Markovian interacting particle systems}
\label{s:exIPS}

\vskip.3cm

We give examples of well-known continuum interacting particle systems, as well as some of their variants,  all falling within the framework presented in Section \ref{s:genmodelassum}.  We  describe update rules in detail and refer the reader to  Appendix \ref{s.admppinteraction} for examples of admissible point processes and stabilizing interaction graphs.   We will  highlight some examples of simple admissible score functions that are of interest.  Concerning admissible update rules, we shall focus only on Markovian updates,  though the framework allows for analysis of non-Markovian updates.  

We again remind the reader that we are studying interacting particle systems in the finite time horizon regime. However, in some  models where update rules satisfy certain monotonicity conditions (e.g., the RSA model in Section \ref{ss:csa}) it may be possible to adapt our proof techniques to cover the infinite time horizon regime. As noted before, there is considerable literature on interacting systems focussed on long-term behaviour and phase-transitions. Such studies are
often model-specific and involve challenging questions outside  the scope of this work.

\vskip.3cm
\noindent{\bf Examples of admissible score functions.}
We list some relevant choices of the function $h$ used to define the scores $\xi(\tx, \tP_n) = h(M(x,t)_{t \in [0,t_0]})$ and the  functional $H_n^{\xi}= \sum_{x \in \P \cap W_n}  \xi(\tx, \tP_n)$ at \eqref{IPSfunctional}. These choices of $h$ yield the following summary statistics, which  figure prominently in the upcoming applications. 

\begin{itemize}[wide,label=(\roman*),  labelindent=0pt]
\item[S1.] Let $\gT \subset [0,t_0]$ be  a finite set.
Put $ h(M(x,t)_{t \in [0,t_0]}):=  {\bf 1}(M(x,t) \in M_0, \, \forall \, t \in \gT)$ for some
measurable subset $M_0 \subset \M$, giving that  $H_n^{\xi}$ is the total number of sites in $\P \cap W_n$ which are in state $M_0$  at (discrete) times in  $\gT$. 
This can be extended to Borel subsets $\gT \subset [0,t_0]$ by taking 
 $ h(M(x,t)_{t \in [0,t_0]}):=
 \exp\{\int_{\gT}\log\1{M(x,t) \in M_0}\,\md t\}$.

\item[S2] Put $h(M(x,t)_{t \in [0,t_0]}) :=  \int_0^{t_0} {\bf 1}(M(x,t) \in M_0)\,\md t$ for some measurable subset $M_0 \subset \M$, giving that  $H_n^{\xi}$ is the total time the constituent sites in $\P \cap W_n$ are in state $M_0$.

\item[S3.] Given  $t_{\min} \in  [0,t_0]$ and a measurable subset $M_0 \subset \M$, 
put $ h(M(x,t)_{t \in [0,t_0]}) := {\bf 1}( \int_0^{t_0} {\bf 1}(M(x,t) \in M_0)\,\md t \geq t_{\min})$  giving that  $H_n^{\xi}$ is the total number of  sites  in $\P \cap W_n$ which are in state $M_0$ for at least time $t_{\min}$.

\item[S4.] When $\M = \R$, put $ h(M(x,t)_{t \in [0,t_0]}) :=  M(x,t_0)$,  giving that  $H_n^{\xi}$ is the total sum of all states at time $t_0$.

\item[S5.] Put  $h(M(x,t)_{t \in [0,t_0]}) :=  h'(M(x,t_1),...,M(x,t_m))$, where $t_1,...,t_m \in [0,t_0]$ are fixed and $h':\M^m \to \R$ is a measurable function.  This yields the study of finite-dimensional distributions of the process $(M(x,t))_{t \in [0,t_0]}$.

\end{itemize}

Note that  $\xi$ in S1--S3 trivially satisfy the moment condition \eqref{momh} and, in case of S4 and S5, we refer to  Sections 2.1 and 2.2 of \cite{BYY19} for conditions on the point process $\P$ insuring that $\xi$ satisfy the requisite moment conditions. 
%Remark \ref{r:stationarized_version}(3) for conditions insuring that $\xi$ n satisfy the requisite moment %conditions in case of S4 and S5.

We illustrate how measurability may be verified in the above examples.  Consider  the example of score function S1.  Remember, 
the history   of state modifications at site $x$  is formalized as a counting measure $M(x,t)_{t \in [0,t_0]}= \sum_j \delta_{(t_j,M(x,t_j))} \in \hat{\cN}_{[0,t_0] \times \M}$ where $t_j$'s are the update times of $M(x,\cdot)$ and are ordered in increasing order. Observe that this increasing ordering is a measurable function of $M$ from $\hat{\cN}_{[0,t_0] \times \M}$ to $\hat{\cN}_{[0,t_0]}$. With this representation, for any $t \in [0,t_0]$ and $M_0 \subset \M$, we may write 
$$ \1{M(x,t) \in M_0} =  \sum_{j} \1{t_j \leq t < t_{j+1}}\1{M(x,t_j) \in M_0},$$
and note that the sums and indicators are also compositions of measurable real-valued functions from $\hat{\cN}_{[0,t_0] \times \M}$. Thus $\1{M(x,t) \in M_0}$ is a measurable real-valued function from $\hat{\cN}_{[0,t_0] \times \M}$ and hence the finite product as in score function S1 for finite $\gT$ is also a measurable function. Similarly the other score functions can be checked to be measurable.

\vskip.3cm

\vskip.3cm
\noindent{\bf Markovian  update rules.} The evolution is said to be {\em Markovian} if the update rule $\F$  at \eqref{e:updaterulef}
depends only on the configuration at the time of update,  i.e., it does  not depend on the entire history but only the current state at the site and its neighbors. To be more precise, for any  marked clock event $(T,L)\in\ttau_x$ of some site  $x\in\P_n$, we set $M(x,T-) \in \M$ to be the value of the mark 
of the last  event in $M(x,\cdot)$ strictly before time $T$, and hence $M(x,T-)$
is a measurable function of $(M(x,t))_{t < T}$.
If the update rule $\F$ in~\eqref{e:updaterulef} has the form
\begin{align}
 \F\left(x,L,\Bigl\{\Bigl(y-x,\{(t,M(y,t))\}_{t <T}\Bigr)\Bigr\}_{\{y\in\P_n:y\simeq x\}}\right) 
&= \tilde{\F}\left(x,L,\Bigl\{\bigl(y-x,M(y,T-)\bigr)\Bigr\}_{\{y\in\P_n:y\simeq x\}}\right)   \label{e:updaterulem}\\
&= \Bigl\{\Bigl(y-x,\bigl(T,M(y,T)\bigr)\Bigr)\Bigr\}\in \hat{\cN}_{\R^d\times \hat{\cN}_{[0,t_0] \times \M}}  \no
\end{align}
for some  function 
$\tilde{\F}: \R^d\times \mL  \times \hat{\cN}_{\R^d \times \M  }\to \hat{\cN}_{\R^d\times \hat{\cN}_{[0,t_0] \times \M}}$
%\footnote{\dy{Range should be $\hat{\cN}_{\R^d\times %[0,t_0] \times \M}$}}  
which is jointly measurable  with respect to the product $\sigma$-algebra, then we say that both the rule $\F$ and  the 
interacting particle system are {\em Markovian}.
The specification of the update rule corresponds to the random mapping representation of a Markov process.  We do not specify the generator as it is not necessary for our analysis.  

In the coming subsections,  we shall explicitly specify update rules satisfying condition \eqref{e:updaterulem}.   These extend some well-known Markovian models of interacting particle systems by allowing them to have dependent initial states, unbounded admissible interaction neighborhoods, and continuum particle locations given by admissible input $\P$, and hence the locations  may be spatially correlated.   We will also show that Theorem~\ref{thmIPS} 
%\BBLp{\sout{and \ref{IPSLLN}}} 
immediately establishes  Gaussian fluctuations  for statistics of these models described by admissible score functions of the type S1--S5,  including  univariate and multivariate central limit theorems
given by Theorem \ref{t:cltmarkedpp} and Theorem \ref{t:multcltmarkedpp}. 
To show that Theorem~\ref{thmIPS} is  indeed applicable,  {\em we only need to show that the update rules  $\tilde{\F}$ for these models are admissible, the only  missing ingredient.} 

Our first three models do not involve locally synchronous updating and thus do not utilize the full generality of our update rule  \eqref{e:updaterulem}. On the other hand,  the last three models in Section \ref{ss:balld} and Section \ref{ss:exclusion} involve locally synchronous updating in an essential way.

\subsubsection{Cooperative and random sequential adsorption (CSA and RSA)}
\label{ss:csa}
Cooperative sequential adsorption (CSA) is a model of adsorption of particles in the continuum~\cite{shcherbakov2023probabilistic}.  We may establish Gaussian fluctuations for this model for unbounded
interaction ranges, addressing an open problem in \cite[Section 5.2]{Penrose2008existence}.
Informally, particles arrive at random times at sites given by an admissible point process
$\P$ and the probability that an arriving particle centered at $x$ is accepted depends on the 
configuration of previously arrived particles in the neighborhood of $x$ in a stabilizing interaction graph and the acceptance rule may have a stochastic component as well. 
The choice of unbounded interaction range is in keeping with many classical models; see e.g. \cite{Evans}. We now present more formally a CSA model introduced in \cite[Section 7]{penrose2001random}  and which includes RSA as a special case. RSA models on the  random geometric graph  are of interest in statistical physics (see \cite{sanders2015sub}) and  in wireless networks. These models on fixed bounded degree graphs were studied in \cite{penrose2005exact}, but we do not require bounded degree. In case of Poisson input $\P$, a space-time variant of this model has been investigated in \cite{gloria2013random}, \cite{Penrose2002limit}, \cite{PY2001}, \cite{SPY}.

Let $p : \hat{\cN}_{\R^d} \to [0,1]$ be a measurable function.  Informally,  particles arrive at random times at sites  on the  substrate $\R^d$.  The sites are 
 given by an admissible point process $\P$  and if the configurations of  occupied neighbors of a site $x$ is $\cX$  then the particle is placed at the site (i.e., it is adsorbed) with probability $p(\cX-x)$.  We consider states $M(x,t) \in \{0,1\}$, with $1$ denoting the presence of a particle at site $x$ and time $t$  and $0$ denoting its absence.  

Denote by  $M(\cdot,0) = M(\cdot)$ the initial states, not necessarily independent. When a marked clock rings at $(T,L)\in\ttau_x$ at some site  $x\in\P_n$, where $L$ is a uniform $[0,1]$ random variable,
we define  the update rule as follows: If $M(x,T-) = 0$ and $\{y : M(y,T-) = 1,  y \sim x\} =: N_{1,x}$
then 
we set $M(x,T) = 1$ with probability 
$p(N_{1,x}-x)$, i.e., when 
$L \le p(N_{1,x}-x)$. In all other cases we have $M(x,T):=M(x,T-)$, i.e., the state does  not change at $x$. The neighbors of $x$ are not updated. 
More formally, the Markovian  update rule $\tilde{\F}$ in \eqref{e:updaterulem} is considered with
\begin{equation} M(x,T)  :=  \1{  M(x,T-) = 1 }+ \1{M(x,T-) = 0}\1{L \leq p(N_{1,x}-x)}.\label{e.CSA2}
\end{equation}
The measurability of $\tilde{\F}$ follows since the indicator 
restricts a point process to a measurable subset of $\R^d \times \{1\}$.  Thus  $\tilde{\F}$ is admissible.  

As a corollary of Theorem~\ref{thmIPS},  
for a given $t_0 \in (0, \infty)$,  we obtain the limit theory for $(\mu^\xi_n)_{n \in \N}$ with
$\xi(\tx,\tP_n) = M(x,t_0)$.  In this case $H_n^{\xi}$ is simply the total number of accepted particles in the window $W_n$ under CSA dynamics up to time $t_0$.

The following result is an immediate consequence of Theorem~\ref{thmIPS}.
\begin{corollary}[Limit theory for the number of accepted particles
 in CSA] \label{c:CSAgen} 
Consider the CSA model with input  point process $\tP$ given in~\eqref{e:IPS-input}, with state space $\M=\{0,1\}$ and clocks having uniformly distributed marks in $\mL=[0,1]$. As in Theorem~\ref{thmIPS},we assume $\tP$ has {\em summable exponential $\B$-mixing correlations},  and  Palm moment bounds~\eqref{e:bd_palmcorr} for all $p, m \in \N$. We assume  Markovian update rules~ \eqref{e:updaterulem} with \eqref{e.CSA2} on a stabilizing interaction graph $\G$ on finite windows with respect to $\P$ as in~\eqref{e:ddecay}.
Put $\xi(\tx,\tP_n) = M_n(x,t_0)$  for
$t_0 \in (0, \infty)$.
Then the measures $(\mu_n^{\xi})_{n \in \N}$ at \eqref{IPS} satisfy the conclusions of Theorem~\ref{thmIPS}.
When $f \equiv 1$ this yields the asymptotic normality for the total number of accepted particles in $W_n$ up to time~$t_0$. 
 \end{corollary}

Choosing $p(\cX) = \1{\cX = \emptyset}$
we obtain {\em the random sequential adsorption (RSA) model on a graph}.  In other words,  a particle is adsorbed at an empty site if there are no occupied neighboring sites.  
 Note that the first clock ring of the Poisson process definitively fixes the state (0 or 1) of the site $x$, provided the clock rings before the finite time $t_0 < \infty$. Otherwise, $x$ retains its initial state $M(x) = M(x,0)$.  
\begin{remark}[Further comments and comparison with previous work]
\begin{enumerate}[wide,label=(\roman*),  labelindent=0pt] \
\item (Comparison with literature.) 
In the case of Poisson point  with bounded range of interactions, Corollary~\ref{c:CSAgen} %\ref{i.CSA-1}
can be deduced from \cite{Penrose2008existence}. However the case of unbounded interactions, even in the case of Poisson input, was left open in \cite{Penrose2008existence}.  %\sout{\dy{With additional effort one could generalize this model to allow for locally synchronous updates.}}
Corollary~\ref{c:CSAgen} %\ref{i.CSA-1}
adds to the central limit theory of previously studied RSA models having finite time horizon (\cite{Penrose2002limit}, \cite{PY2001} and references therein),  which is confined to Poisson input,  finite range interactions and independent initial states.

\item  (Multivariate CLT.) More generally,  given $t_0 \in (0, \infty)$ ,  
subsets $\gT_i \subset [0, t_0],$ scores $\xi_i(\tx, \tP_n) = \int_{\BBLp{\gT}_i} {\bf 1}(M(x,t) = 1)\, \md t$ and assuming convergence of covariances,  we find that Theorem \ref{t:multcltmarkedpp} (along with Remarks~\ref{rem:IPS}\ref{i.LLN-Var-IPS} and~\ref{i.rem-extensions}) establish a multivariate central limit theorem for the $k$-vector whose $i$th entry is the number of
accepted particles in the time period $T_i$.

 \item  (Graph based RSA models  with random particle sizes.)  
Let marked clocks have rings $(T,L)\in\ttau_x$ with $L$ being a random particle shape, taking value in the space~$\mathcal{K}$  of  nonempty, compact subsets of $\R^d$, with diameter bounded by a fixed constant $r$.   
 The update Markovian  rule~\eqref{e:updaterulem} for $M(x,t)\in\mathcal{K}$  is as follows:
 If $M(x,T-) = \emptyset$ and for all $y\in\P_n$, $y\sim x$,
$M(y,T-)\cap (x \oplus L)=\emptyset$, 
then
we set $M(x,T) := x \oplus L$.
Otherwise, the state of $x$ does not change and we put  $M(x,T) = M(x,T-)$. Consider the score function $\xi(\tx,\tP_n) = \1{M(x,t_0)\neq \emptyset}$ for $t_0 \in (0, \infty)$.
 Assuming admissible $\tP$ and letting the stabilizing interaction graph be the  Gilbert (geometric)  graph with radius $2r$,
we may check that all assumptions 
of Theorem \ref{thmIPS}  are satisfied,
thus giving the limit theory in particular for the total number $H_n^{\xi}$ of accepted particles in this model up to time $t_0<\infty$ in $W_n$.

\item  (Further extensions.) Variants of the CSA model including dimer RSA, the annihilation process, and bootstrap percolation may be analysed using the  framework of this subsection; see also \cite[Section 2]{penrose2005exact}. 

\end{enumerate}
\end{remark}

\subsubsection{Epidemic spread and voter models}
\label{ss:epidemic} 

In the SIR epidemic model,  particles are of three types: susceptible (S), infected (I) and recovered (R). 
Particle locations and initial states are given by the realization of $\tP$. 
Whether a susceptible particle becomes infected is a function of the number of infected neighbors.  However, once infected, a particle
may reover on its own, and similarly,  a recovered particle may again become susceptible.  We shall consider more general dynamics where the updates are based on the relative locations and states of the neighboring sites.  We provide here
conditions under which the number of infected individuals in the window $W_n$ at time $t_0$ has Gaussian fluctuations.
 We adopt a Markovian set-up, but we could also consider a non-Markovian example, whereby updates are based on the entire time-evolved history of neighbors.  

Formally,  set $\M = \{\mathrm{S},\mathrm{I}, \mathrm{R}\}$.  Let $p : (\M \times \hat{\cN}_{\R^d \times \M}) \times \M \to [0,1]$ be a jointly measurable probability kernel governing the transition from $\mathrm{S},\mathrm{I},\mathrm{R}$ to $\mathrm{S},\mathrm{I},\mathrm{R}$ and depending on the configuration of neighbors,  i.e., for $\tmu \in \hat{\cN}_{\R^d\times {\M}}$,  $u \in \M$, $p(u,\tmu,\cdot)$ is a probability distribution function   on $\M$.  As before,  $L_i, i \geq 0,$ are i.i.d.  uniform random variables with values in $\mL = [0,1]$.  Let the initial states be given by $M(x,0), x \in \P$, which are not assumed to be independent.  

Set
\begin{align*}
p_\mathrm{S} &:= p(M(x,t-), \bigl\{(y-x,M(y, t-))\bigr\}_{y \sim x},\mathrm{S}),  \\
p_\mathrm{I} &:= p(M(y,t-), \bigl\{(y-x,M(y, t-))\bigr\}_{y \sim x},\mathrm{I})
\end{align*}
and put $p_\mathrm{R} = 1 - p_\mathrm{S} - p_\mathrm{I}$. 
Define the function $\tilde{\F}$ in  \eqref{e:updaterulem} as follows: 
\begin{align}\label{e:M-RSI}
& M(x,t) := \mathrm{S} \1{L \leq p_\mathrm{S}} + \mathrm{I} \1{ p_\mathrm{S} < L \leq p_\mathrm{S} + p_\mathrm{I}} + \mathrm{R} \1{L> p_\mathrm{S} + p_\mathrm{I}},
 \end{align}
where the sum is to taken as a formal sum.  Only one of the indicators is non-zero and hence $\tilde{\F}$ is well-defined.  

Put $\xi^{(1)}(\tx, \tP_n) = {\bf{1}}(M(x,t_0) = \mathrm{I})$ so that $H_n^{\xi^{(1)}}$ is the total number of infected individuals at time $t_0$.   Letting $\xi^{(2)}(\tx, \tP_n) = \int_0^{t_0}{\bf{1}}(M(x,t_0) = \mathrm{I})\,\md t$
gives that $H_n^{\xi^{(2)}}$ is the total time all individuals are infected. To deduce the limit theory
for these statistics we apply  directly  Theorem~\ref{thmIPS}.  
 
\begin{corollary}[Gaussian fluctuations for the number of infected individuals and for the total infection time] \label{thmepidemic} 
Consider the  SIR epidemic model with input  point process $\tP$ given in~\eqref{e:IPS-input}, with state space $\M = \{\mathrm{S},\mathrm{I}, \mathrm{R}\}$ and  clocks having uniformly distributed marks in $\mL=[0,1]$. 
As in Theorem~\ref{thmIPS}, we assume $\tP$  has  {\em summable exponential $\B$-mixing correlations} and  Palm moment bounds~\eqref{e:bd_palmcorr} for all $p, m \in \N$.  We assume the Markovian update rules~\eqref{e:updaterulem} along with \eqref{e:M-RSI} on  a stabilizing interaction graph $\G$ on finite windows with respect to $\P$ as in~\eqref{e:ddecay}.
 Then the measures $(\mu_n^{\xi^{(i)}})_{n \in \N}, i = 1,2,$   
at \eqref{IPS} induced by the score functions $\xi^{(i)},i=1,2$ above,   respectively satisfy the conclusion of Theorem~\ref{thmIPS}.  When $f \equiv 1$ this yields the asymptotic normality for the total number of infected individuals at time $t_0$ as well as the total time the individuals are infected during the period $[0,t_0]$. 
 \end{corollary}

This corollary extends upon Section 4.5 of  \cite{Penrose2008existence}, which is constrained to models having finite range and independent initial states.
The models described here arise in the study of epidemics and also in the propagation of viruses in  networks. Suitably  choosing the probability kernel  yields classical models including  the Richardson model,  the contact process, the voter model, and the chase-escape model  (see  \cite[Chapter 8]{Jahnel2020} and \cite{draief2010epidemics}).   The more general kinetically constrained models which include bootstrap percolation and its stochastic counterpart given by the Fredrickson-Andersen model (see \cite{hartarsky2022fredrickson,hartarsky2024kinetically}), can also be handled by the above framework. 

\subsubsection{Majority dynamics in the continuum}
\label{ss:maj}

This model is a specific version of the SIR epidemic model. We again denote the sites by $\P$, and now
the state at each site $x \in \P$ is $\{+1,-1\}$-valued and at each clock ring, the state is updated to that of the majority of its neighbors. In case of a tie, the state is unchanged. More formally, we take $\M=\{-1,+1\}$
(and in general our setting, the initial states $M(x)$ , $x \in \P$, and they are not necessarily independent). 
 When the  clock at site $x$ rings at time $t$,  it undergoes a Markovian update as follows : If $\sum_{y \sim x} M(y,t-) \neq 0$, then $M(x,t) = \text{sgn}(\sum_{y \sim x} M(y,t-))$ where $\text{sgn}$ is the sign function.  Otherwise we put $M(x,t) = M(x,t-)$.  Formally, the update rule  $\tilde{\F}$  in \eqref{e:updaterulem} is defined as follows:
 \begin{align} M(x,t) :=  \text{sgn}\big(\sum_{y \sim x}M(y,t-)\big)\1{\sum_{y \sim x}M(y,t-) \neq 0} + M(x,t-)\1{\sum_{y \sim x}M(x,t-) = 0}.\label{e:M-Majority}
\end{align}
The measurability of $\tilde{\F}$ here can  be verified by writing it as composition of projections,  sums, and the sign function.   Thus $\tilde{\F}$ is admissible. 

Put $\xi(\tx, \tP_n) = {\bf{1}}(M(x,t_0) = 1)$ so that $H_n^{\xi}$ is the total number of sites in state $1$  at time $t_0$.
 Theorem~\ref{thmIPS} immediately gives the following corollary.

\begin{corollary}[Gaussian fluctuations for the size of the majority/minority] \label{thmmajority} Consider the  majority/minority model with input  point process $\tP$ given in~\eqref{e:IPS-input}, with state space $\M=\{-1,+1\}$  and clocks
 with no marks. 
As in Theorem~\ref{thmIPS}, we assume $\tP$  has  {\em summable exponential $\B$-mixing correlations},  and  Palm moment bounds~\eqref{e:bd_palmcorr} for all $p, m \in \N$.  We assume  Markovian update rules~\eqref{e:updaterulem} with \eqref{e:M-Majority} on a stabilizing interaction graph $\G$ on finite windows with respect to $\P$ as in~\eqref{e:ddecay}.
Then the measures $(\mu_n^{\xi})_{n \in \N}$   at \eqref{IPS} induced by the score functions $\xi$ above  satisfies the conclusions of Theorem~\ref{thmIPS}.  When $f \equiv 1$ this yields the asymptotic normality for the total number of sites in state  $1$ at time $t_0$. 
 \end{corollary}

One could also  consider this model of majority dynamics in discrete time with synchronous updates; see Section~\ref{s:discIPS}.  These results add to 
the  central limit theorem  for majority  dynamics in a dense Erd\"{o}s-R\'enyi random graph (cf. \cite[Theorem 1]{berkowitz2020central}). More general opinion dynamics such as spatial graph versions of {\em edge-averaging dynamics} or {$\ell^p$-minimization dynamics} can also be considered; see for example \cite{aldous2012averaging,amir2025convergence}.

\subsubsection{Ballistic deposition} \label{ss:balld}
We  sketch  the model and refer the reader to \cite{Penrose2002limit} and  \cite{penrose2008growth} for a discussion of similar models. Particles rain down at sites given by the realization of a  point process~$\P$ and they stick to the  first point of contact. 
We assume that the  particles are random shapes
with random heights. We thus consider clock marks $L = (L^{(1)}, L^{(2)})\in [0,\infty)\times \mathcal{K}$, where $L^{(1)}$ denotes the height of the particle (arriving at some clock event time $T_i$ of some site $x$), and $L^{(2)}$ is a shape variable---a random element of the space of compact sets $\mathcal{K}$ of $\R^d$---assumed to be  nonempty, contain the origin~$\0$.
As a technical assumption, we suppose that the height $L^{(1)}$ is possibly unbounded but has all moments finite, while the diameter of $L^{(2)}$ (the particle shapes in $\mathcal{K}$) is bounded by a fixed finite constant $r>0$, and $L^{(1)}$, $L^{(2)}$ may be mutually dependent.

We now formally define a specific model. As the state space, we consider $\M := \R \times \mathcal{K}$, where each state $M = (M_1, M_2) \in \M$ represents a particle $M_2$ positioned at level $M_1$, with $M_1$ indicating the height at which $M_2$ is aligned at the top.
The assumption that the particle diameter is bounded by a given $r>0$ allows one 
to consider the Gilbert graph on $\P_n$, $\G(\P_n)$, where $x \sim y$ if and only if $|x-y| < 2r$. With this specification, when a clock rings at time $t$ at a site $x$ with a mark $L = (L^{(1)},L^{(2)})$, we update $M(x,t)$ and $M(y,t)$ for $y \sim x$ by the function $\tilde{\F}$ in~\eqref{e:updaterulem}
with 
\begin{align}\label{e:deposition-new1}
 M_1(x,t) &:=   L^{(1)} + \max_{y: y\simeq x,\;   (y\oplus M_2(y,t-))\cap (x \oplus L^{(2)})\not=\emptyset}  M_1(y,t-), \qquad 
M_2(x,t):=L^{(2)},\\
M_1(y,t)&:= M_1(y,t-),\qquad   M_2(y,t) := \Bigl((y \oplus M_2(y,t-)) \setminus (x \oplus L^{(2)})\Bigr) \oplus (-y).\label{e:deposition-new2}
\end{align}
Note, the neighbors $y\sim x$ keep their  top level $M_1(y,t-)$ unchanged while 
their particles $M_2(y,t-)$ are eroded  by the shadow of the particle $L^{(2)}$ deposited at $x$. 
The {\em height function}
$$
H(y,t) := \max_{x \in \P_n :\, x \oplus M_2(x,t) \ni y} M_1(x,t), \quad y \in \R^d,
$$
represents the  growth level of the  hyper-surface
 in this ballistic deposition model. Putting $\xi^{(1)}(\tx, \P_n) = M_1(x, t_0)|M_2(x,t_0)|$ gives that $H_n^{\xi^{(1)}}= \sum_{x \in \P_n}\xi^{(1)}(\tx, \P_n)
$ is the total  volume under the height function $H$ at time $t_0$.  Alternatively, putting 
$\xi^{(2)}(\tx, \P_n) = |M_2(x,t_0)| \allowbreak\1{M_1(x, t_0) \allowbreak\in [h_{\text{min}},  \infty))}$,  we obtain that 
$H_n^{\xi^{(2)}}$  is the total  surface area visible at height at least $h_{\text{min}}$ at time 
$t_0$.

\begin{corollary}[Gaussian
fluctuations for the volume under the height function and the
surface area at a given height] \label{thminterface} 
 Consider the ballistic deposition model with input  point process $\tP$ given in~\eqref{e:IPS-input}, with initial states  $M(x)=(0,\emptyset)$ and  the clocks as described above, with the particles' diameter satisfying $\text{diam}(L^{(2)})<r$ for some  $r<\infty$, and with
the height having moments $\E[{(L^{(1)}})^p]<\infty$ for all $p>1$. 
We assume  $\P$  has  {\em summable exponential mixing correlations} and  Palm moment bounds~\eqref{e:bd_palmcorr} for all $p, m \in \N$.  We assume  Markovian update rules~\eqref{e:updaterulem} along with~\eqref{e:deposition-new1}--\eqref{e:deposition-new2} (assuming $x\sim y$ iff $|x-y|\le 2r$).
Then the measures $(\mu_n^{\xi})_{n \in \N}$ at \eqref{IPS}   
induced by the score functions $\xi=\xi^{(i)}$,  $i=1,2$  above  satisfy the conclusions of Theorem~\ref{thmIPS}. When $f \equiv 1$,  this  yields the asymptotic normality for the total volume $H_n^{\xi^{(1)}}$ under the  height function
and the total  surface  $H_n^{\xi^{(2)}}$ visible above some given  height $h_{{\rm{min}}}$.
\end{corollary}
\begin{proof}
This result follows from Theorem~\ref{thmIPS}, provided the considered score functions $\xi^{(i)}$ satisfy moment conditions~\eqref{momh} for all $p>1$. Indeed, the other assumptions are met: the Gilbert graph stabilizes---$x \sim y$ iff $|x-y|\le 2r$---for any input, and the initial conditions $M(x)$ together with the marked clocks are i.i.d., thus ensuring summable exponential $\B$-mixing correlations for $\tP$.

We now prove the moment conditions for $M_1(x,t_0)=:M_1(x,\tP_n)$, which suffices for those of $\xi^{(i)}(\tx,\tP_n), i = 1,2$ since the particle surface area $|M_2(x,t)|$ is bounded for all $t \in [0,t_0]$.
Thus,  the latter assumption is equivalent to 
showing that the probability 
$\sP_{\bk x1p}(M_1(x_1,\tP_n) > m)$ is fast decreasing to~0 as a function  of $m\to\infty$.

For proving this, observe that
\begin{equation}
M_1(x_1,\tP_n)\le \sum_{y:\; |y-x_1|\le D(\tx_1,\tP_n)}\sum_{T_i(y)\in\hat\tau_y}L^{(1)}_i(y)\,,
\label{e:ballistic-height-bound}
\end{equation}
that is, the height function of $x_1$ at time $t_0$ is bounded by the sum of all particle heights $L^{(1)}_i(y)$ corresponding to arrivals at the clock times $T_i(y)$ during the interval $(0,t_0]$, over all sites $y$ in the backward cluster of $x_1$, which are within distance $D(\tx_1,\tP_n)$ as defined in~\eqref{e:diambackwardcluster}.

Also note that $L_i^1(y),\hat{\tau}_y$'s are 
i.i.d. for $y \in \P_n$ given $\P$. In particular for all $y \in \P_n$, $\sum_{T_i(y)\in\hat\tau_y}L^{(1)}_i(y)$ has distribution of $\sum_{j = 1}^{N_0}L^{(1)}_i$, where $L_i^{(1)}$ are i.i.d. random variables distributed as $L^{(1)}$ and $N_0$ is an independent Poisson random variable with mean $t_0$. Thus, from the assumed finiteness of moments of $L^{(1)}$ and tail property of $N_0$, we have the fast decay of the tail of the above sum i.e., 
$$  \textrm{$\sP(\sum_{j = 1}^{N_0}L^{(1)}_i \geq m)$ decays fast in $m$ as $m \to \infty$.}$$
Now using the above fast decay, Lemma \ref{l.expdecay}, and the fast decay of $D(\tx_1,\tP_n)$ (see Step 2 of the proof of the Theorem \ref{thmIPS}; in particular \eqref{e:radstab_IPS_bd1}, \eqref{bdAm} and \eqref{e.final-D>m}), we deduce that for any $\alpha > 0$, there exists a fast decreasing function $\tilde \varphi_p$ such that 
\begin{equation}
 \sP_{\bk x1p} \Big(\max_{y \in \P_n : |y-x_1| \le D(\tx_1,\tP_n)} \sum_{T_i(y) \in \hat\tau_y} L^{(1)}_i(y) > m^{\alpha} \Big) \leq \tilde \varphi_p(m), \ \quad m \in \N.   
\label{e:fastdecaym_sumN_oL2}  
\end{equation}
Denote by $N_1=|\{ |y-x_1| \le D(\tx_1,\tP_n)\}|$ the number of sites under the maximum. We fix some $\alpha \in (0,1)$ and $\beta = (1-\alpha)/(d+1)$. If the above fast decay event does not happen and $M_1(x_1,\tP_n) > m$, then we must have $N_1 > m^{1-\alpha}$. Combining these observations we have
\begin{align*}
& \sP_{\bk x1p}(M_1(x_1,\tP_n) > m) \\ 
& \leq \sP_{\bk x1p} \Big( \sum_{T_i(y)\in\hat\tau_y} L^{(1)}_i(y) > m^{\alpha} 
\textrm{ for some $y \in \P_n$ with $|y-x_1| \le D(\tx_1,\tP_n)$} \Big) \\
& \quad + \sP_{\bk x1p} \big( N_1 > m^{1-\alpha} \big)  \\
& \leq \tilde \varphi_p(m)  + \sP_{\bk x1p}( N_1 > m^{1-\alpha}, D(\tx_1,\tP_n) > m^{\beta}) + \sP_{\bk x1p}( N_1 > m^{1-\alpha}, D(\tx_1,\tP_n) \leq m^{\beta}) \\ 
& \leq \tilde \varphi_p(m)  + \sP_{\bk x1p}(D(\tx_1,\tP_n) > m^{\beta}) + \sP_{\bk x1p}( \P(B_{m^{\beta}}(x_1)) \geq m^{1-\alpha}) ),
\end{align*}
where in the penultimate inequality we have divided into cases whether $D(\tx_1,\P_n) > m^{\beta}$ or not and in the last inequality we have used the simple fact that if $D(\tx_1,\tP_n) \leq m^{\beta}$ and $N_1 > m^{1-\alpha}$, 
then $B_{m^{\beta}}(x_1)$ must contain  at least $m^{1-\alpha}$ points from $\P$.

It remains to justify the fast decay of second and third terms in the right-hand side of the last inequality. The fast decay of the second term again follows from Step 2 of the proof of the Theorem \ref{thmIPS}; see \eqref{e:radstab_IPS_bd1}, \eqref{bdAm} and \eqref{e.final-D>m}. The fast decay of the third term follows from the deviation inequality in Lemma \ref{l:palm_conc}.
\end{proof}

This result extends the central limit theorem which could be deduced 
from the methods of \cite[Section 4.4]{Penrose2008existence} which requires the stronger assumptions that the interaction ranges are bounded,  the initial states are i.i.d., and the input $\P$ is Poisson. One could likewise modify the score $\xi$ to obtain Gaussian fluctuations for other statistics, including the  total number of contacts, further extending results
of \cite[Section 4.4]{Penrose2008existence} to non-Poisson input and unbounded interaction ranges.

The above model is but a simple prototype which  admits extensions and variants.  A further generalization, amenable to our analysis, involves including  queueing dynamics whereby  the particles start slowly diminishing in height and vanish \cite{baccelli2011poisson}.

\subsubsection{Further models}
\label{ss:exclusion}

We describe two additional models which may be treated by Theorem \ref{thmIPS} and Remark \ref{rem:IPS}\ref{i.LLN-Var-IPS}. For examples of further update rules, and thus additional examples of
spatial interacting particle systems satisfying our framework, 
 we refer  to \cite[Section 4]{Penrose2008existence}, \cite[Chapter 8]{Jahnel2020}, \cite[Chapters 8 and 9]{draief2010epidemics}, \cite[Chapter 1]{swart2017course}, \cite{mairesse2014around}, \cite{evans2005nonequilibrium} and \cite{hartarsky2024kinetically}.  
Section \ref{s:future} describes additional   models falling within our framework but which require separate investigation.

\begin{enumerate}[wide,label=(\roman*),labelindent=0pt]
\item {\bf Exclusion processes.} 
In this model, within the framework of the update rules at  ~\eqref{e:updaterulem}, some sites are 
initially occupied by particles which attempt to hop at independent times 
according to an exclusion rule. In the simplest case, when the clock at an 
occupied site $x$ rings, the particle at $x$ chooses a neighboring site at 
random and jumps there if it is unoccupied, otherwise it remains at $x$. Hence, 
each site contains at most one particle, and the total number of particles is 
conserved. Neighboring sites need not lie within a bounded interaction range, and hence the model goes beyond  of~\cite[example (c) in Section~4.5]{Penrose2008existence}. 
Moreover, the sites are not assumed to form a homogeneous Poisson process, 
but only a realization of an admissible point process under the assumptions of 
Theorem~\ref{thmIPS}, and the initial states need not be independent. This 
prototype can be extended to allow multiple particles per site and more general 
jump rules, covering, for instance, some of the birth–death–migration models 
in~\cite{dawson2014spatial} and~\cite[Section~4.1]{Penrose2008existence}.

 To handle this simple
exclusion process, we proceed as follows. 
We set $\M = \{0,1\}$ and let the $L$'s be  i.i.d. uniform random variables on $[0,1]$.  Consider a jointly
measurable probability kernel $p : \hat{\cN}_{\R^d} \setminus \{\emptyset\} \times \R^d \to [0,1]$ which selects a site in $\mu$ based on its location and assume that $p(\mu,\cdot)$ is supported on $\mu$ for all $\mu \in \hat{\cN}_{\R^d} \setminus \{\emptyset\}$,  i.e.,  $p(\mu, \cdot)$ gives a probability distribution on $\mu$.  For example,  we may choose $p(\mu,y)$ to be proportional to $|y|$ for $y \in \mu$.  

Given $\mu$, we define $X'(\mu)$ to be the random element in $\mu$  with
 distribution $p(\mu, \cdot)$. For convenience, set $X' := X'(\{y -x: y \simeq x\} =X'(\{y -x: y \sim x\} \cup \{\0\})$.  
  Define $\tilde{\F}$ in \eqref{e:updaterulem} according to the  rule  for $x$ and $y\sim x$ : 
 \begin{align*}
 M(x,t) &:=M(x,t-)M(x+ X',t-),   \\
  M(y,t) &:=  \Bigl(M(x,t-)(1-M(x+ X',t-))+(1-M(x,t-))M(x+ X',t-)\Bigr)\1{x + X' = y}\\
 &\ \hspace{2em}+ M(y,t-)\1{x + X' \neq y},
\end{align*}
where $L$ is used to choose the random element $X'$  (see a similar random mapping representation as in Section \ref{ss:epidemic}).  

Choosing a score function related to the occupations $M(x,t)$ at $t\in[0,t_0]$ of the site $x$ (see  examples given at the beginning of Section \ref{s:exIPS})
we obtain the asymptotic normality for the  measures $(\mu_n^{\xi})_{n \in \N}$ at~\eqref{IPS} leveraging Theorem \ref{thmIPS}.

\item {\bf Divisible sandpile dynamics.}
\label{ex:sandpile}
In the classical Abelian sandpile model (see, e.g., \cite{jarai2018sandpile}), 
each site contains an integer number of grains. When the number of grains 
at a site exceeds a fixed threshold, the site becomes unstable and topples, 
sending one grain to each of its finitely many neighbors. Such topplings 
may trigger further instabilities, leading to a cascade of updates known as 
an avalanche, which continues until the system reaches a stable configuration.

Within the framework of~\eqref{e:updaterulem}, we assume that at time $t=0$ 
each site carries an initial amount of mass, and at independent exponential 
times the site topples part of this mass to its neighbors. Importantly, 
the topplings are not followed by the complete avalanche mechanism; 
that is, we do not enforce further updates until stability is restored, 
but rather there is only a single, local redistribution at the ringing site.  
However, unlike the classical sandpile, the set of neighbors updated in a toppling need not lie within a fixed deterministic distance.

We put $\M = \R$ and  we interpret $M(x,t) \in \mR$  as the amount of (signed) mass at site $x$, 
with $x \in \P$.   The initial masses are given by $M(x,0)$, not necessarily independent.
When the clock at $x$ rings at time $t \in (0,t_0]$, we update  $\tilde{\F}$ in \eqref{e:updaterulem} according to the  rule  for $x$ and $y\sim x$:  
 \begin{align*}
 M(x,t) &:=  M(x,t-){\bf{1}}(M(x,t-) < 1) \, + \,  {\bf{1}}(M(x,t-) \geq 1),  \\
 M(y,t) &:= M(y,t-) + \frac{M(x,t-)-1}{|N_{x}|} {\bf{1}}(M(x,t-) \geq 1) ,
 \end{align*}
where $N_{x} := \{x : y \sim x\}$. In other words,  $x$ retains unit mass for itself and distributes the remaining mass equally among its neighbors. 
If $M(x,t-) < 1$, then there is no update at $x$ or its neighbors. 

For score functions~$\xi$ related to the mass evolution $M(x,t)$ at $t\in[0,t_0]$ of the site $x$ (see  examples given at the beginning of Section \ref{s:exIPS}) we obtain the asymptotic normality for the  measures $(\mu_n^{\xi})_{n \in \N}$ at~\eqref{IPS} leveraging Theorem \ref{thmIPS}.
In particular, if we set   
$$\xi(\tx, \tP_n)  = \sum_{t \in \tau_x}[M(x,t-) - 1]\cdot {\bf 1}[M(x,t-) \geq  1],
$$
then  $H_n^{\xi} := \sum_{x \in \P \cap W_n} \xi (\tx,\tP_n)$ yields the `odometer' up to time $t_0$,  i.e., the total amount of mass distributed to neighbors until time $t_0$.  

The above model, where excess mass is proportionally dispatched to neighbors, is not Abelian 
since the order of topplings affects the mass distribution.  
Our framework, however, also accommodates more classical Abelian dynamics and variants such as the chip-firing game, the dollar game as well as more general Abelian networks; see \cite{bond2016abelian,fey2009,jarai2018sandpile,klivans2018mathematics,corry2018divisors}.

\end{enumerate}

\subsection{Discrete-time interacting particle systems}
\label{s:discIPS}
We now consider discrete-time particle systems. In many models such as probabilistic cellular automata or kinetically constrained models (see for example \cite{mairesse2014around,hartarsky2024kinetically}), it is natural to consider globally synchronous dynamics and these are modelled better by discrete-time particle systems. We shall now  modify our framework  and then state and prove our limit theorems for discrete-time particle systems.

Let $\M$ and $\mL$ be Polish spaces as before, hosting respectively 
the states of the sites and the marks of the clocks (see Section~\ref{sss.Clocks}). 
In contrast to the previous setting with individual (independent Poisson) clocks $\tau_x$, $x \in \P$, 
we now replace them by a universal clock ringing at times $1,2,\ldots,J$ for some finite $J<\infty$, 
simultaneously for all $x \in \P$. 
That is, the input process 
$
\tP=\{(x,U(x))\}_{x \in \P} \in \cN_{\R^d \times \K}
$ 
is defined by the ground process $\P$ (specifying the locations of the sites), 
equipped with marks $U(x)=(M(x),\bfL_x)\in\K$ taking values in the Polish space 
$\K := \M \times \mL^J$.
Here, as before, $M(x)=M(x,0)\in\M$ denotes the initial state of site $x$, while
\[
\bfL_x = (L_1(x), \ldots, L_J(x))\in\mL^J
\]
is an individual vector of marks where $L_i(x), i \leq J, x \in \P$ are i.i.d. (conditionally on $\{(x,M(x))\}_{x\in\P}$) $\mL$-valued random elements.   
For simplicity of notation, we often suppress the dependence on $x$ when writing $L_i$.  
We emphasize again that we do not assume $\P$ to be a Poisson point process, 
nor do we assume independence of the initial states $\{M(x)\}_{x \in \P}$.

We assume  the framework of stabilizing interaction graphs from Section~\ref{sss.Graph}. For the update rules and state evolution, we refer to Section~\ref{sss.Apdates}, modifying it as follows:  
At each universal clock tick $j \leq J$, the state at all sites $x \in \P_n$ is updated by adding event $M(x,j)$ to the previous history  
\[
\bfM_{j-1}(x)=\bfM_{j-1}(x;n) := \bigl(M(x,0), M(x,1), \dots, M(x,j-1)\bigr) \in \M^{j}.
\]  
The update is performed via the function $\F_j$, which depends on the clock mark  $L_j(x)$ and the histories of $x$ and its neighbors $y \sim x$, possibly accounting for their relative positions.
More formally, for each $j \leq J$, we assume the existence of a measurable  updating  function 
 $\F_j:   
\R^d\times \mL  \times \hat{\cN}_{\R^d \times  \M^{j} }\to \hat{\cN}_{\R^d\times  \M^{j+1}}$
and put 
\begin{align} \label{e:updaterulef-discrete}
M(x,j):=\F_j\left(x,L,\Bigl\{\Bigl(y-x,\bfM_{j-1}(y)\Bigr)\Bigr\}_{\{y\in\P_n:y\simeq x\}}\right).
\end{align}
The  complete history of a given site $x \in \P_n$---resulting from the execution of the functions $\Phi_j$ at all clock rings $j\in[J]$---is denoted by $M_n[x]=
(M(x,0),\ldots M(x,J))\in \M^{J+1}$.  These  histories
are treated as an additional marking of  the point $x \in \P_n$, which is  measurable as a measurable mapping
\begin{equation}\label{e.g-measurability-discrete}
g : (\R^d \times \K) \times \hat{\cN}_{\R^d \times \K} \ni (\tx, \tmu) \mapsto M[x] \in \M^{j+1},
\end{equation}
constructed iteratively  using the measurable functions $\Ph_j$ described above. 
This mapping describes the history of marks $M_n[x], x \in \P_n,$ when applied simultaneously on all $x \in \P_n$. 

Note that, unlike the continuous-time model of Section~\ref{s:genmodelassum}, the states at sites $x \in \P$ are updated at globally synchronous times $t=1,\ldots,J$.  
However, each site $x$ updates only its own state, using the previous histories of its neighbors.

Put $\xi(\tx, \tP_n):=  h(M_n[x])$, where $h : \M^{J+1} \to \R$ is a  measurable function.  
  We now state our central limit theorem for linear statistics of  $(\mu_n^{\xi})_{n \in \N}$ in~\eqref{IPS}, similar to Theorem~\ref{thmIPS}, and then  we shall show how certain Markovian interacting particle systems (including all of those in Section \ref{s:exIPS}) fall within this framework. 
\begin{theorem}[CLT for discrete-time interacting particle systems]\label{discreteIPS} 
Let the input process for the discrete-time interacting particle system be 
$\tP = \left\{(x, M(x),\bfL_x)\in\R^d \times \M \times \mL^J\right\}$,
where $J \in \N$ and $\bfL_x$ are vectors of i.i.d.
elements of $\mL$, independent of $\sum_x \delta_{(x,M(x))}$. Put
\begin{equation}
\label{IPSagain}
\mu_n^{\xi} := \sum_{x \in \P \cap W_n} \xi (\tx,\tP_n)  \delta_{n^{-1/d}x}.
\end{equation}
If 
$\tP$ and $\xi$ respectively satisfy assumptions~\ref{i.thmIPS} and~\ref{iii.thmIPS} in Theorem ~\ref{thmIPS}, 
$\F_j$ are update rules as at ~\eqref{e:updaterulef-discrete} with respect to a stabilizing interaction graph $\G$ on finite windows with respect to $\P$ as in \eqref{e:ddecay}, then
 the random measures $(\mu_n^{\xi})_{n \in \N}$ satisfy the central limit theorem. That is, for every test function $f \in \B(W_1)$ such that
$
\Var[\mu_n^{\xi}(f)] = \Omega(n^{\nu})$ for some $\nu > 0$,
we have
$$
(\Var[\mu_n^{\xi}(f)])^{-1/2} \big(\mu_n^{\xi}(f) - \sE[\mu_n^{\xi}(f)]\big)  \stackrel{d}{\Rightarrow} Z,$$
where $Z$ is a standard Gaussian random variable.
\end{theorem}

\begin{remark}[Further comments and comparison with previous work]
\label{r:discIPS}\ 
\begin{enumerate}[wide,label=(\roman*),  labelindent=0pt]
\item By imposing additional assumptions of translation-invariance on score-functions, stationarity on $\tP$, and a uniform bound for the interaction ranges across all windows for every $x\in\P$ one may obtain expectation and variance asymptotics for $(\mu_n^\xi)_{n \in \N}$ as in Remark \ref{rem:IPS}\ref{i.LLN-Var-IPS} .

\item When $\P$ is Poisson and when the states $M(x), x \in \P$, are independent,  then this result follows from the central limit theorems of \cite{Penrose2007gaussian} and \cite{Lachieze2019normal}.  We  require neither Poisson input nor independent initial states. 

\item(Markovian interacting particle systems)
Let $M(x) \in \M$.  Let the updates at discrete times $j\in[J]$ be as follows
$$ M(x,j):=\F\left(x,L,\Bigl\{\Bigl(y-x,M(y,j-1)\Bigr)\Bigr\}_{\{y\in\P_n:y\simeq x\}}\right).
$$
where $\F: \R^d \times \mL \times \hat{\cN}_{\R^d \times   \M}   \to \M$ is a measurable function.  Informally,  at each step,  a vertex is updated according to its own state,  the  states of its neighbors and also the location of  its neighbors. 

\item Update rules defined in Sections \ref{ss:epidemic}--\ref{ss:maj} can be trivially adapted in the above Markovian framework but with the caveat that the evolution of the particle systems can be considerably different.  For example,  in the RSA model in discrete time,  it is possible that neighboring sites can become occupied together.  

\item Examples in Sections \ref{ss:balld} and \ref{ss:exclusion} need to be modified suitably to fit into the above framework as the update rules therein are more general than those in Sections \ref{ss:epidemic}--\ref{ss:maj}.   For example in the divisible sandpile dynamics,  there could be an exchange of mass as every site will re-distribute its excess mass in parallel.  More specifically,  we may define the updates as follows.  For $j = 1,\ldots,J$ 
\begin{align*}
M(x,j) & = M(x,{j-1})\mathbf{1}(M(x,{j-1}) < 1) \, + \, \1{M(x,t_{j-1}) \geq 1} \\
& \quad +  \sum_{y \, : \, y \sim x} \frac{M(y,j-1)-1}{|N_x|}\1{M(y,t_{j-1}) \geq 1}.
\end{align*}
\end{enumerate}
\end{remark}

\begin{proof}[Proof of Theorem \ref{discreteIPS}.] 
We show that $(\mu_n^{\xi})_{n \in \N}$
satisfy the assumptions of Theorem \ref{t:cltmarkedpp}. We only need to establish the fast stabilization of $\xi$ on finite windows.
Due to the updates occurring in parallel at discrete times, the proof is simpler than in the continuous-time case.  
Precisely, for every $x \in \P_n$ we consider a backward cluster $C_n'(x,j)$ at discrete time instances $j \leq J$, consisting of all space-time points $(y,j')$, $y \in \P_n$, $j' \in \{0,\ldots,j\}$, that (possibly indirectly) contribute to the update of $x$ up to  time $j$.  
By the nature of discrete time and the update rule~\eqref{e:updaterulef-discrete}, this backward cluster grows backward from $j'$ to $j'-1$, $j' \leq J$, by adding neighbors of points $y$ in the cluster at time $j'$.  

Consequently, the projection of the space-time cluster $C_n'(x,J)$ onto $\P_n$ gives, besides $x$, the union of all $j$-nearest neighbors of $x$ in $\P_n$, for $j \leq J$. Call this set $\cN_J(x,\P_n)$.  
As in \eqref{e:stabradips}, we enlarge this set by the graph stabilization radius of its elements, $S_n(y,\P_n)$, $y\in\cN_J(x,\P_n)$ defined in~\eqref{e:S_n=R_Wn}. We obtain for $x \in \P_n$ the bound
\begin{equation} \label{e:stabradips-discrete} 
R_{W_n}^{\xi}( \tilde{x}, \tilde{\P}_n) \leq \max_{y \in \cN_J(x,\P_n)} \bigl(|y-x| + S_n(y, \P_n) \bigr). 
\end{equation}  
Using~\eqref{neighborulespin} and the  stabilization of the graph on finite windows,   iteratively,  together with Lemma~\ref{l.expdecay}, this yields the required fast stabilization of $\xi$ on finite windows for general $J$ and completes the proof of Theorem~\ref{discreteIPS}.  \end{proof}

\section{Empirical random fields and geostatistical models} \label{s:erfgeostat}

In the three subsequent subsections, we consider models involving two random objects defined on \(\R^d\): a random field \(M\) with values in a general Polish space, and a point process \(\P\). A key feature of these objects is that they possess fast mixing correlations.  We consider the following two scenarios:

\begin{itemize}
    \item The point process \(\P\) is assumed independent of the random field \(M\). 
Viewed as a \emph{sampling process} in an increasing window \(W_n\), 
it allows us to estimate statistics of the random field \(M\), 
such as its distribution function or covariogram. 
We study the properties of the corresponding estimators, 
including their asymptotic mean, variance, and Gaussian fluctuations. 
These results appear  in Sections~\ref{ss.Empirical-field} and~\ref{s:Covariogram}.
    
    \item In Section~\ref{ss:geoBM}, still assuming independence between \(M\) and \(\P\), we consider models based on \(\P\) where the random field \(M\) is seen as  providing  `external' pre-marks in the construction of the model. These are referred to as \emph{geostatistical models} (see~\cite{SRD04,ahlberg2018gilbert}). We illustrate this with a Boolean model of this type and present results related to the total length of the edges induced by the model.
\end{itemize}

In both of the above scenarios the theoretical framework of  Part~\ref{part:theory-foundations} allows us to achieve these goals, providing either new results or reaffirming existing results under assumptions that are easier to verify in practice. Also the independence between the field $M$ and point process $\P$ eases the formulation of  mixing or localization criteria without reference to Palm formalism; our general framework can possibly work for the more general case of dependent random fields $M$; see Remark \ref{r.TGGM}.

Specifically, let $M := \{M(x)\}_{x \in \R^d}$ be a random field taking values in a Polish space~$\M$. 
As mentioned above, the assumptions and results regarding this random field $M$ 
concern its mixing correlations, which we formulate below, 
naturally extending the definitions given in Definitions~\ref{d.omegatmixing} and~\ref{B-omega-mixing}.

\begin{defn}[Mixing correlations of the random field]
\label{d.omegatmixing-field}\ 
 \begin{enumerate}[wide,label=(\roman*),labelindent=0pt]
\item \label{i.BL-mixing-field} The random field $M=\{M(x)\}_{x\in\R^d}$ with values in $\M$  
satisfies  {\em $\BL$-mixing correlations}  if
there exists an increasing-decreasing family of 
 functions $\omega_k :  [0, \infty) \to [0,2]$, $k \in \N$, 
(called correlation decay functions), such that $\lim_{s \to \infty}\omega_k(s) = 0$ for all $k \in\mathbb{N}$ and for all $p,q \in \mathbb{N}$, $x_1,\ldots,x_{p+q} \in \mR^d$, and all $f \in \BL(\M^p)$, and
  $g \in \BL(\M^q)$ we have
\begin{align} \label{e.field-mixing}
\left| \sE[f(\bk {M(x_i)}1{p}) \, g(\bk  {M(x_i)}{p+1}{p+q})] - \sE[f(\bk  {M(x_i)}1p)] \, \sE[g(\bk  {M(x_i)}{p+1}{p+q})] \right| \le \omega_{p+q}(s),
\end{align}
where $s := d(\bk x1p,\bk x{p+1}{p+q})$ is at \eqref{defs}.
\item \label{i.B-mixing-field} The random field $M$
satisfies {\em $\B$-mixing correlations}  if~\eqref{e.field-mixing} holds when \BL\ is replaced by $\B$.
 \item  Under Item~(i)  (resp. Item (ii)), the random field  $M$  has  {\em fast $\BL$-mixing correlations}  (resp.  {\em fast $\B$-mixing correlations}) 
 if $(\omega_k)_{k \in \N}$ are fast decreasing. \label{i.mixing-field-fast}
\end{enumerate}
\end{defn}
In  Section~\ref{ss:Geostat}, we present explicit constructions of random fields 
$M$ satisfying the above mixing correlation properties, which are required for the results in Sections~\ref{ss.Empirical-field}--\ref{ss:geoBM}. This is done in Lemma \ref{l.decorrelation-field-new} by constructing the random field as a marking function of an independent marked point process $\tXi$ which serves as a \emph{covariate} for the random field \(M\). This is analogous to Theorem \ref{t:mppmix_double}.

\subsection{Limit theory for empirical  random fields} 
\label{ss.Empirical-field}
%%then fast mixing of the point process alone is sufficient.
For a real valued measurable function $\zeta : \R^d \times \M \to \R$, we consider the following  random field $\zeta(M):=\{\zeta(x,M)\}_{x\in\R^d}$
sampled at $\P$:
\be \label{linearstatBsigma}
\hat\sigma_{n}^{\ze(M)} := \sum_{x \in \P_n} \ze(x,M(x)) \delta_{n^{-1/d}x},
\ee
where $\P_n := \P \cap W_n$. We call $\hat\sigma_{n}^{\ze(M)}$ {\em empirical random field.}

In what follows, we establish Gaussian fluctuations for~ $(\hat\sigma_{n}^{\ze(M)})_{n \in \N}$
under fast mixing correlations of $\P$ and $M$, and suitable moment assumptions on~$\zeta(M)$. We say that $\zeta$ satisfies the $p$-moment condition with respect to $M$ if
\begin{equation}
\label{e:M-p-moment}
\sup_{x \in \R^d} 
\sE\big[ \max\big(1, |\zeta(x, M(x))|^p \big) \big]
\leq M_p^{\zeta} 
< \infty.
\end{equation}

\begin{proposition}[limit theory for the empirical distribution function of a random field]
  \label{t.randomfieldCLT-new}
Let $M=\{M(x)\}_{x\in\R^d}$ be a random field in $\R^d$ with values in a Polish space
$\M$ and $\zeta: \R^d \times \M \to \R$ a  measurable function. 
Let $\P$ be a simple point process on $\R^d $,
independent of the random field $M$. 
Consider empirical random measures $(\hat{\sigma}_{n}^{\zeta(M)})_{n \in \N}$ as in \eqref{linearstatBsigma}.
\begin{enumerate}[wide,label=(\roman*),labelindent=0pt]
\item \label{i.CLT-infinite-zeta-field-mixing}
Assume that $\P$ has fast mixing correlations 
as in Definition~\ref{d.omegamixing}\ref{i.mixing-ground-fast}. Assume that either (a) the random field $M$ satisfies the fast $\B$-mixing condition of Definition~\ref{d.omegatmixing-field}\ref{i.B-mixing-field}-\ref{i.mixing-field-fast} or (b) 
$M$  satisfies the fast \BL-mixing condition of Definition~\ref{d.omegatmixing-field}\ref{i.BL-mixing-field}-\ref{i.mixing-field-fast}
and $\zeta(x,\cdot)$ is Lipschitz. If $\zeta$ satisfies the $p$-moment condition~\eqref{e:M-p-moment} with respect to $M$ for all $p>1$ then the central limit theorem holds for
    $(\hat\sigma_{n}^{\zeta(M)}(f))_{n \in \N}$, i.e., as $n \to \infty$
\[   (\Var \, \hat\sigma_{n}^{\zeta(M)}(f) )^{-1/2}\Big( \hat\sigma_{n}^{\zeta(M)}(f) - \sE \hat\sigma_{n}^{\zeta(M)}(f) \Big) \stackrel{d}{\Rightarrow} Z,
  \]
for $f \in \B(W_1)$,  provided
$\Var \, {\hat{\sigma}_{n}^{\zeta(M)}(f)} = \Omega(n^{\nu})$ for some $\nu > 0$.

\item \label{i.EVar-infinite-zeta-new} 
Suppose the random field~$M$ and the point process~$\P$ are stationary 
(necessarily jointly so, due to their independence) with~$\P$ having intensity~$\rho$, 
and suppose that~$\zeta$ is translation invariant, 
i.e., $\zeta(x,\cdot) = \zeta(\cdot)$. Also assume $\zeta$ satisfies the $(1+\epsilon)$-moment condition~\eqref{e:M-p-moment} with respect to $M$ for some $\epsilon > 0$.
Then the mean empirical field satisfies for all $f \in \B(W_1)$
\begin{equation}\label{eq:mean-limit}
n^{-1} \, \E \, \hat\sigma_n^{\zeta(M)}(f)
= \rho\, \sE \, \zeta(\0,M(\0)) \int_{W_1} f(x)\, \md x, \quad n\in\N,
\end{equation}
which holds  without assuming mixing conditions on~$\P$ and~$M$ nor  Lipschitz conditions on $\zeta$. Moreover, under either of the mixing assumptions (a) or (b) of part~\ref{i.CLT-infinite-zeta-field-mixing} we obtain for all $f \in \B(W_1)$
\[
\lim_{n\to\infty} n^{-1} \Var \, \hat\sigma_{n}^{\zeta(M)}(f)
= \sigma^2(\zeta) \int_{W_1} f(x)^2\,\md x,
\]
where
\begin{align}  \label{sigdef-emp}
\sigma^2(\zeta) 
&:= \rho\, \sE \, \zeta^2(\0,M(\0)) \\
&\quad + \int_{\R^d} 
\Big( 
\sE[\zeta(\0,M(\0))\, \zeta(x,M(x))]\, \rho^{(2)}(\0,x)
- \rho^2\, \sE[\zeta(\0,M(\0))]\, \sE[\zeta(x,M(x))]  
\Big)\, \md x \in [0, \infty) \no
\end{align}
provided that~$\zeta$ satisfies the $(2+\epsilon)$-moment condition~\eqref{e:M-p-moment} for some $\epsilon > 0$.
\item \label{i.zeta-MulitCLT} Let
$\bzeta := (\zeta_1, \ldots, \zeta_m)$ be a vector of translation-invariant, real-valued measurable functions 
on $\R^d \times \M$, with each $\zeta_i$ satisfying the $(2+\epsilon)$-moment condition~\eqref{e:M-p-moment} with respect to $M$ for some $\epsilon > 0$. 
Under either of the mixing assumptions (a) or (b) of part~\ref{i.CLT-infinite-zeta-field-mixing}, and under the stationarity assumption of part~\ref{i.EVar-infinite-zeta-new}
 the multivariate central limit theorem holds as in Theorem~\ref{t:multcltmarkedpp}.
\end{enumerate}
\end{proposition}

\begin{remark}
Central limit theorems for empirical measures of random fields were established by~\citet{Pa09},
who assumes strong mixing conditions (in terms of $\sigma$-fields) on both the point process~$\P$ 
and the random field~$M$. 
He also points out that, in practice, these conditions are often difficult to verify. 
In contrast, assuming fast mixing correlations of~$\P$ 
together with fast $\B$- or $\BL$-mixing correlations of~$M$ 
constitutes a weaker and often more tractable alternative in applications; see Section \ref{ss:Geostat} for examples of such $M$.
\end{remark}

\begin{proof}[Proof of Proposition~\ref{t.randomfieldCLT-new}]
We first make some common observations on  the proofs of the three statements. Consider the marked point process $\tilde{\P}:=\{(x,\zeta(x,M(x)))\}_{x\in\P}$.
Due to the independence between  $\P$ and $M$, the corresponding Palm distributions of $\zeta(x,M(x))$
(with respect to multiple locations of $x$ in $\P$) 
coincide with the (unconditional) joint finite-dimensional distributions 
of $\zeta$,  i.e., for all bounded functions $f$ and $x_1,\ldots,x_p \in \R^d$, we have
\begin{equation}
\label{e:ZetaPalmindep}
\sE_{\bk x1p}\big[ \, f(\bk{\xi(x_i,\zeta(x_i,M(x_i)))}{1}{p}) \, \big] = \sE\big[ \, f(\bk{\xi(x_i,\zeta(x_i,M(x_i)))}{1}{p}) \, \big],
\end{equation}
where in the left-hand side the Palm expectation is with respect to $\P$.

 We prove now  statement~\ref{i.CLT-infinite-zeta-field-mixing} where we  directly apply the  umbrella result, 
Theorem~\ref{t:clt_linear_marks}. It suffices to show that $\tilde{\P}$ has fast 
\BL-mixing correlations as in 
Definition~\ref{d.omegahmixing}\ref{i.mixing} and~\ref{i.family-mixing} 
(since here $\P_n=\P\cap W_n$). 
For this, we rely on the result on joint mixing correlations of points and marks 
given in Proposition~\ref{p:pmix+mmix}. 
Due to the independence of  the random field $M$ and the point process $\P$,  
the assumed fast $\B$-mixing correlations of $M$ 
translate directly into the fast  \BL-mixing of the marks 
$\zeta(x,M(x))$, $x\in\P$, as functions of $M$ inherit the fast $\B$-mixing correlations of $M$, 
i.e. fast $\B$-mixing correlations of $\{M(x)\}_{x\in\R^d}$ implies
fast $\B$- (hence \BL-)~mixing   correlations of the  field~$\{\zeta(x,M(x))\}_{x\in\R^d}$. 
Also, if the function $\zeta(x,\cdot)$ is Lipschitz, then fast \BL-mixing of $M$ is sufficient for this last observation. Thus, in view of \eqref{e:ZetaPalmindep} and under either of the assumptions (a) or (b) of \ref{i.CLT-infinite-zeta-field-mixing}, we have fast \BL-mixing  correlations of the  marked point process~$\tP = \{\zeta(x,M(x))\}_{x\in\P}$. This justifies statement~\ref{i.CLT-infinite-zeta-field-mixing} via Theorem~\ref{t:clt_linear_marks}.

We now establish statements~\ref{i.EVar-infinite-zeta-new}
and~\ref{i.zeta-MulitCLT}. 
The first of these, concerning expectations in \eqref{eq:mean-limit}, follows directly from the Campbell formula 
under the stationarity and independence assumptions on~$\P$ and~$M$. 

All other statements (i.e., variance asymptotics \eqref{sigdef-emp} and the multivariate CLT in part~\ref{i.zeta-MulitCLT})  share a common feature: 
the marked point process $\tP = \{(x, \zeta(x, M(x)))\}_{x \in \P}$ 
exhibits fast $\BL$-mixing correlations. 
This property follows directly from the mixing properties of~$\P$ and~$M$,
as explained in the proof of part~\ref{i.CLT-infinite-zeta-field-mixing}.
The asymptotics of the variance and the multivariate central limit theorem
under the $(2+\epsilon)$-moment condition 
then rely on the fast $\BL$-mixing of~$\tP$ in the framework of Theorem~\ref{t:clt_linear_marks_new}. 
Note that we do not require summable fast \BL-mixing correlations of $\P$ as is to be expected from Theorem \ref{t:clt_linear_marks_new} or Proposition \ref{expvar}. In those results, such an assumption on $\tP$ (which corresponds to $\P$ here) is to establish fast $\BL$-mixing of $\dbtilde \P$ (which corresponds to $\tP$ here)  but here our assumptions on $M,\P$ already guarantee fast \BL-mixing correlations of $\tP$.
\end{proof}

\subsection{Limit theory for covariograms}
\label{s:Covariogram}
In this section, we are interested in the \emph{covariogram}, 
a tool capturing the spatial dependence structure of a real-valued, stationary random field.  
Let \( M = \{M(x)\}_{x \in \R^d} \) denote such a real-valued random field.  
The \emph{theoretical covariogram} of \(M\) is the deterministic, real-valued function  defined by 
\[
\gamma(h) = \gamma_M(h) \, := \frac{1}{2} \sE\big[(M(\0) - M(h))^2\big], \quad h \in \R^d,
\]
see, e.g., \cite[Section~2]{SRD04}.

To estimate the value of \(\gamma(h)\), we consider sampling the random field \(M\) 
at locations given by a stationary, simple point process \(\P\) on \(\R^d\), independent of $M$.  
A straightforward approach would be to average the quantities  
\(\frac{1}{2} |M(x) - M(y)|^2\)  
over all pairs of points \(x, y\) in some observation window \(W_n\) such that \(y - x = \pm h\), where we admit $\pm h$ since the function $\gamma$ is even.
However, unless the point process \(\P\) possesses a lattice structure aligned with 
the specific vector \(h\), there will typically be no such pairs satisfying this exact relation.

To overcome this issue, we introduce a \emph{tolerance parameter} $\delta>0$, such that \(\delta < |h|\),  and instead consider pairs \(x, y \in W_n\) such that \(y-x \in B_\delta(\0\pm h)\), with the notation $B_\delta(a\pm h):=B_\delta(a+h)\cup B_\delta(a-h)\) for $a\in\R^d$. Define the following score function for the points of~\(\P\), involving also $M$:
\begin{equation}\label{e.mark-covar}
\xi(x, \P_n, M) := \frac{1}{4\theta_d \delta^d} 
\sum_{y \in \P_n: y-x\in B_\delta(\0\pm h)} (M(x) - M(y))^2,
\end{equation}
where \(\theta_d\) denotes the volume of the unit ball in \(\R^d\);  
the factor \(4\) (instead of \(2\)) compensates for the fact that the term $(M(x)-M(y))^2$  
is counted twice in the sum, and $\delta<|h|$ prevents taking $x=y$.
 It is natural to postulate that this score satisfies the limit 
\begin{align}\no
\lim_{n \to \infty} \frac{1}{n}\sE\Bigl[ \sum_{x \in \P_n} \xi(x, \P_n, M)\Bigr]
&= \frac{\rho}{4\theta_d \delta^d} \sE_0 \Bigl[\int_{B_\delta(\0\pm h)} (M(\0) - M(y))^2 \P(\md y)\Bigr]\\
&= \frac{\rho}{4\theta_d \delta^d} \int_{B_\delta(\0\pm h)} \sE\big[(M(\0) - M(y))^2\big] \rho^{(2)}(\0, y)\, \md y\no\\
&=\frac{\rho}{\theta_d \delta^d} \int_{B_\delta(h)} \gamma(y)\, \rho^{(2)}(\0, y)\, \md y=:\rho \, \bar\gamma(h),
\label{e:covario-mean}
\end{align}
where \(\rho = \rho^{(1)}\) and \(\rho^{(2)}\) denote, respectively, 
the first- and second-order correlation functions of \(\P\) and, 
in the last equality, we used the translation invariance and symmetry of \(\rho^{(2)}\), 
and the evenness of~\(\gamma\). 
For small \(\delta > 0\), this limit is expected to be close to 
\(\rho^{(2)}(\0, h)\, \gamma(h)\), 
thus providing a practical means to estimate the covariogram using the sampling process~$\P$.

The following result provides a more complete limit theory for the \emph{empirical covariogram} of the random field $M$ given by
\begin{equation}\label{e.sigma-covariogram-new}
\sigma^{\gamma}_n =\sigma^{\gamma(h)}_n := \sum_{x \in \P_n} \xi(x, \P_n, M)\, \delta_{n^{-1/d} x},
\end{equation}
for $h\in\R^d$, based on the independent sampling point process $\P$, where we suppress the dependence on the parameter \(\delta\), $0<\delta<|h|$, related to the tolerance.

\begin{proposition}[Limit theory for empirical covariogram] \label{p:covario-new} 
Let $M=\{M(x)\}_{x\in\R^d}$ be a real-valued random field in $\R^d$ having fast $\B$-mixing correlations as in  Definition~\ref{d.omegatmixing-field}\ref{i.B-mixing-field}\ref{i.mixing-field-fast}.
Let $\P$ be a simple point process on $\R^d$ having fast mixing correlations as in Definition~\ref{d.omegamixing}\ref{i.mixing-ground-fast}, bounded reduced  Palm intensity function~\eqref{bdedreduced}, and independent of the random field $M$.
Consider empirical covariogram $\sigma_{n}^{\gamma}$ of $M$ based on $\P$, as in \eqref{e.sigma-covariogram-new}, for some value $h$ and tolerance $\delta>0$ .
\begin{enumerate}[wide,label=(\roman*),labelindent=0pt]
\item \label{i.CLT-infinite-gamma-new}
If  $M$ has $p$-moment conditions for all $p>1$ (i.e., ~\eqref{e:M-p-moment} holds with  $\zeta(x,M(x))=M(x)$) then the central limit theorem holds  for
    $(\sigma_{n}^{\gamma}(f))_{n \in \N}$, i.e., as $n \to \infty$
\[   (\Var \, \sigma_{n}^{\gamma}(f) )^{-1/2}\Big( \sigma_{n}^{\gamma}(f) - \sE \, \sigma_{n}^{\gamma}(f) \Big) \stackrel{d}{\Rightarrow} Z
  \]
for $f\in\B(W_1)$, provided the random    $\Var \, \sigma_{n}^{\gamma}(f) = \Omega(n^{\nu})$ for some $\nu > 0$.

\item \label{i.EVar-infinite-gamma-new}
Further, if the random field $M$ and the point process $\P$ are  stationary (necessarily jointly so, due to their independence), with $\rho$ the intensity  of $\P$,
then the  expectation and variance asymptotics and multivariate central limit
theorem hold {\em mutatis mutandis} for 
$(\sigma_{n}^{\gamma}(f))_{n \in \N}$
as in Proposition~\ref{expvar} and Theorem \ref{t:multcltmarkedpp}
under $p=(2+\epsilon)-$ or $p=(4+\epsilon)$--~moment condition for $M$ i.e.,~~\eqref{e:M-p-moment} with $\zeta(x,M(x))=M(x)$, respectively.
Precisely, the limiting expectation takes the form $\rho \, \bar\gamma(h)$
of~\eqref{e:covario-mean}
and the covariance of the empirical covariogram at two values  
$\sigma_{n}^{\gamma(h)}(f)$ and $\sigma_{n}^{\gamma(h')}(f)$, with $h,h'\in\R^d$, (and for the simplicity, we assume  the same tolerance $\delta>0$ such that $\delta<\min\{|h|, |h'|, |h-h|', |h+h'|\}$) takes the form
$$\lim_{n\to\infty}\frac 1n {\rm Cov}(\sigma_{n}^{\gamma(h)}(f)),\sigma_{n}^{\gamma(h')}(f)))=
\sigma^2(h,h')\int_{W_1}f(x)^2\,\md x\in (-\infty,\infty),$$
 where 
\begin{align}  \label{sigcovdef-S11}
&\sigma^2(h,h') :=\\
&\quad \frac{\rho}{16\theta^2_d \delta^{2d}} \int_{B_\delta(\0\pm h)}\int_{B_\delta(\0\pm h')}\sE\big[(M(\0) - M(y))^2(M(\0) - M(y'))^2\big]\rho^{(3)}(\0, y,y')\, \md y \, \md y' \no\\
&\quad +
\int_{\R^d}\biggl(\frac1{16\theta^2_d \delta^{2d}}\int_{B_\delta(\0\pm h)\cap B_\delta(z\pm h')}\sE\big[(M(\0) - M(y))^2(M(z) - M(y))^2\big] \rho^{(3)}(\0,z,y)\md y \no\\
&\qquad +\frac1{16\theta^2_d \delta^{2d}}\int_{B_\delta(\0\pm h)}\int_{B_\delta(z\pm h')}\sE\big[(M(\0) - M(y))^2(M(z) - M(y'))^2\big] \rho^{(4)}(\0,z,y,y') \, \md y \, \md y'  \no\\
&\hspace{10cm} - \rho^2 \, \bar\gamma(h) \, \bar\gamma(h') 
\biggr) \md z \in (-\infty, \infty). \no
\end{align}
\end{enumerate}
 \end{proposition}

\begin{proof} We deduce this from Theorem~\ref{t:cltmarkedpp}.
First, we  verify the moment conditions on the score $\xi$.  
Using the independence of~$\P$ and~$M$ we can apply 
the Minkowski inequality and then  the  Campbell-Little-Mecke formula \eqref{e:CLM2} to derive the bound:
\begin{align*}
&\sE_{\bk x1{q}}[ \max(1,|\xi(\tx_1,\tP_n|^p)]\\
&\le \biggl(1+\frac{1}{4\theta_d\delta^d}
\Bigl(\sum_{i=1}^{q}\sE[(M(x_1)-M(x_i))^{2p}]^{1/p}+\int_{B_\delta(x_1\pm h)}\sE[(M(x_1)-M(y))^{2p}]^{1/p}\,\rho^{(1)}_{\bk x1l}(y) \, \md y \Bigr)\biggr)^p.
\end{align*}
The finiteness of the $p$th moment of the score function~\eqref{e.mark-covar} follows from the existence of $2p$-moments of the random field~$M$ and the moment assumption~\eqref{bdedreduced} on the  reduced  Palm intensity of~$\P$.  

We now prove points~\ref{i.CLT-infinite-gamma-new} and~\ref{i.EVar-infinite-gamma-new}.  The score function~$\xi$ in~\eqref{e.mark-covar} has the specific form of a local $U$-statistic, whose general expression is given in~\eqref{e:U-score}.  
The limit theory for such scores can be established under fast $\B$-mixing correlations of the marked input process~$\tP$, as in Corollary~\ref{c.U-FME}.  
More precisely, under these fast mixing correlation the central limit theorem in part~\ref{i.CLT-infinite-gamma-new} follows from Theorem~\ref{t:cltmarkedpp}.

In the stationary setting, the limits of the mean and variance follow from Proposition~\ref{expvar}, under the $(1+\epsilon)$- and $(2+\epsilon)$-moment conditions for the score~$\xi$, respectively.  
Since this score involves the square of the values of the random field~$M$, these requirements correspond to the $p = (2+\epsilon)$ and $p = (4+\epsilon)$ moment conditions~\eqref{e:M-p-moment}, with $\zeta(x, M(x)) = M(x)$, respectively.  
Similarly, Theorem~\ref{t:multcltmarkedpp}, under these reduced mixing assumptions, establishes the multivariate central limit theorem under the $p = 2 + \epsilon$ moment condition for~$\xi$, which again translates to the $p = (4+\epsilon)$ moment condition~\eqref{e:M-p-moment} with $\zeta(x, M(x)) = M(x)$.  

Finally, the explicit expressions of the limiting quantities follow from a  suitable application of the Campbell-Little-Mecke formula \eqref{e:CLM2} involving~$\P$ and~$M$.
\end{proof}
\subsection{Geostatistical Boolean models} 
\label{ss:geoBM}
In the two previous sections, our main object of interest was a random field $M = \{M(x)\}_{x \in \R^d}$.
%taking values in a general Polish space or simply in~$\R$. 
We considered a point process~$\P$ in~$\R^d$, independent of this field, 
used to estimate distributional properties of~$M$.
In this context, the process~$\P$ was regarded  as a `sampling point process'.
In the present section, our focus shifts to more complex models that are based on the process~$\P$ as a ground process 
and involve pre-marks~$U(x)$ attached to points~$x \in \P$. 
This terminology was introduced in Section~\ref{s.Stabilizing} and has been used throughout the applicative part of our work. 
However, here, these pre-marks take the form~$U(x) = M(x)$, where the values are driven by an `external' random field~$M$, 
independent of~$\P$. 
In this setting, the marks are referred to as \emph{geostatistical marks}, 
and the corresponding models are called \emph{geostatistical models}, 
as suggested in the seminal paper~\cite{SRD04}. 
Despite this paradigmatic change in the roles of~$M$ and~$\P$, 
one can easily recognize in the covariogram considered in Section~\ref{s:Covariogram} 
a particular instance of a geostatistical model. 
In this section, we revisit another model of this type, namely the \emph{geostatistical Boolean model} studied in~\citet{ahlberg2018gilbert}, 
thus illustrating the relevance of the  limit theory developed in Section~\ref{s:clta} in this  context.

%\subsubsection{\dy{\sout{Geostatistical %Boolean model}}}

As in  previous subsections,  consider a real, non-negative valued random field~$M$ on~$\R^d$, 
and let~$\P$ be a point process on~$\R^d$, independent of~$M$. 
Following~\cite{ahlberg2018gilbert}, we define the 
\emph{geostatistical random geometric graph}~$\G(\P, M)$, 
which places an edge between two points~$x, y \in \P$ whenever 
$|x - y| \leq M(x) + M(y)$,
i.e., whenever the balls at both the points intersect each other.
Equivalently, this is the graph induced by the \emph{spherical geostatistical Boolean model (GeoBM)} 
on~$\P$, where to each point~$x \in \P$ we attach a ball of radius~$M(x)$, 
given by the value of the (independent) random field~$M$. See the upcoming Example \ref{ex:PV_GBM} for a Poisson-Voronoi driven geostatistical Boolean model.

\vskip.3cm

We establish below the limit theory for the total edge length of~$\G(\P_n, M)$, 
thus taking up a statistic  suggested by~\citet{ahlberg2018gilbert}, but whose limit theory was not analyzed.

Specifically, consider the  score function
\begin{equation}\label{e.mark-Geo-BM}
\xi(x, \P_n, M) := \frac{1}{2} 
\sum_{y \in \P_n} |x - y| \, \1{|x - y| \le M(x) + M(y)},
\end{equation}
which corresponds to one half the sum of the lengths of the edges in~$\G(\P_n, M)$ 
that are incident to~$x \in \P_n$. Define {\em the edge length measure}
\begin{equation}\label{e:Geo-edges}
\sigma_n^\GeoBM := \sum_{x \in \P_n} \xi(x, \P_n, M) \, \delta_{n^{-1/d}x}.
\end{equation}
Note that $\sigma_n^\GeoBM(1)$ equals the total edge length of the graph~$\G(\P_n, M)$.

\begin{proposition}[Limit theory for total edge lengths in GeoBM] 
\label{thm:GeoBM}
Let $M = \{M(x)\}_{x \in \R^d}$ be a non-negative, real-valued random field on $\R^d$ having fast $\B$-mixing correlations as in Definition~\ref{d.omegatmixing-field}\ref{i.B-mixing-field}-\ref{i.mixing-field-fast}.
Assume that $M$ is bounded by some constant $R$, that is,
$\sup_{x \in \R^d} M(x) \le R < \infty$.
Let $\P$ be a simple point process on $\R^d$ having 
fast mixing correlations  as in Definition~\ref{d.omegamixing}\ref{i.mixing-ground-fast}, bounded reduced  Palm intensity function~\eqref{bdedreduced},
 and independent of the random field $M$.
Consider the  edge length measure $\sigma_{n}^{\GeoBM}$ of the geostatistical graph $\G(\P_n, M)$ as in \eqref{e:Geo-edges}.
\begin{enumerate}[wide, label=(\roman*), labelindent=0pt]
\item \label{i.CLT-GeoBM}
A central limit theorem holds for the family $\big(\sigma_{n}^{\GeoBM}(f)\big)_{n \in \N}$; that is, as $n \to \infty$,
\[
  \big(\Var \, \sigma_{n}^{\GeoBM}(f)\big)^{-1/2}
  \Big( \sigma_{n}^{\GeoBM}(f) - \sE \, \sigma_{n}^{\GeoBM}(f) \Big)
  \stackrel{d}{\Rightarrow} Z,
\]
for all $f \in \B(W_1)$, provided $\Var \, \sigma_{n}^{\GeoBM}(f) = \Omega(n^{\nu})$ for some $\nu > 0$.

\item \label{i.stat-GeoBM}
Further, if the random field $M$ and the point process $\P$ are stationary, with $\rho$ denoting the intensity of $\P$,
then the expectation and variance asymptotics
hold for 
$\big(\sigma_{n}^{\GeoBM}(f)\big)_{n \in \N}$,
as stated in Proposition~\ref{expvar}. % and Theorem~\ref{t:multcltmarkedpp}.
\end{enumerate}
\end{proposition} 
\begin{proof} 
The boundedness of the random field~$M$ implies that the score function~$\xi$ 
in~\eqref{e.mark-Geo-BM} is a local $U$-statistic, 
whose general expression is given in~\eqref{e:U-score}.
The limit theory for such scores can be established under fast $\B$-mixing correlations
of the marked input process~$\tP=\{(x,M(x))\}_{x\in\P}$, as stated in Corollary~\ref{c.U-FME},
in the same way as for the covariogram considered in Proposition~\ref{p:covario-new}
(the bound for the  reduced Palm intensity of $\P$ is only required and not the bound for higher-order Palm moments
due to the boundedness of~$M$).
 \end{proof}

 \begin{exe}[Poisson–Voronoi–driven GeoBM]
 \label{ex:PV_GBM}
 Consider the following Boolean model whose grains are defined according to a geostatistical marking. 
Given a Poisson point process  $\X$, let  ${\rm Vor}(\X)$ be the  Poisson-Voronoi tessellation of $\R^d$ induced by $\X$, with cells 
$C(y,\X), y \in \X$. For $x \in \R^d$, $C(x,\X)$ denotes the cell containing $x$. 
In other words, if $y \in \X$ is the nearest point of $\X$ to $x \in \R^d$ (which is a.s. defined), then we have that $x\in C(y,\X)$.
We assign i.i.d. real-valued random variables $M(y)$ to each  $C(y,\X)$, $y \in \X$. We assume the random variables are bounded by $R \in (0, \infty).$

Next, given the ground point process $\P$ on $\R^d$, 
  assumed independent of $\X,$ we construct a Boolean model by putting a ball around each  point $x \in \P$ and whose radius  equals the random variable $M(x)$ assigned to $C(x,\X)$ (i.e., $C(y,\X)$ such that $x\in C(y,\X)$).
  The grains $B(x, M(x)), x \in \P,$ are  spatially dependent, as their radii are a function of the underlying Poisson-Voronoi tessellation ${\rm Vor}(\X)$. As in the second example  of the introduction to \cite{ahlberg2018gilbert}, this is a geostatistical Boolean model on $\P$, with cell radii determined by the external Poisson point process $\X$.  This model generates a geostatistical random geometric graph
  $\G(\P_n, M)$.  Proposition \ref{thm:GeoBM} gives conditions on $\P$ such that  the total edge length of $\G(\P_n, M)$ satisfies the CLT. 
  Under the same conditions on $\P$, one could similarly show that the edge count of $\G(\P_n, M)$ is also asymptotically normal.  Questions regarding such statistics of $\G(\P_n, M)$ are indicated at the conclusion  of the introduction to \cite{ahlberg2018gilbert}. 
 \end{exe}

\begin{remark}[Towards more general geostatistical models]\ 
\label{r.TGGM}
\begin{enumerate}[wide,label=(\roman*),labelindent=0pt]
\item {\em Multivariate statistics of GeoBM.} 
Under the assumptions of Proposition~\ref{thm:GeoBM}, 
the result can be straightforwardly extended to a larger class of score functions~$\xi$
related to GeoBM, as given in Corollary~\ref{c.U-FME}.
Moreover, in the stationary setting, 
one obtains a multivariate central limit theorem for several such scores 
by leveraging Theorem~\ref{t:multcltmarkedpp}.

\item {\em Unbounded random field in GeoBM.} 
\label{i.URF}
Extending the result of Proposition~\ref{thm:GeoBM} to the case 
where $\sup_{x\in\R^d} M(x)$ has fast-decreasing tail probabilities is possible, 
but requires summable exponential $\B$-mixing of 
$\tilde \P = \{(x,M(x))\}$ as in Definition~\ref{def.A2}.

\item {\em Beyond GeoBM.} 
The assumption of summable exponential $\B$-mixing of 
$\tilde \P = \{(x,M(x))\}$ also allows one to obtain the same limit theorems 
for \BL-localizing score functions in fairly general geostatistical models.

\item {\em Geostatistical models  with a dependent random field.} 
The paper~\cite{SRD04}, in the first lines of page~80,
discusses the desirability of relaxing the fundamental
assumption of independence in geostatistical marking, 
namely between the sampling locations (ground process $\P$) 
and the values of the spatial random field $M$.
Our framework, in its full generality, allows one  
to address such generalized models, capturing the dependence between $M$ and $\P$ 
through the Palm distributions of the pre-marks $M(x)$ for $x\in\P$.
(Recall that in the classical geostatistical models these Palm distributions 
coincide with the original finite-dimensional distributions of $M$.)
In this setting, again, summable exponential $\B$-mixing of 
$\tilde \P = \{(x,M(x))\}$ and \BL-localizing score functions 
constitute the right tools.
It then becomes primarily a modeling question 
what kind of dependence between $M$ and $\P$ 
still allows one to obtain the required mixing properties of~$\tilde \P$.
\end{enumerate}
\end{remark}

\subsection{Mixing correlations of random fields with a covariate point process} 
\label{ss:Geostat}
In the previous section, the mixing of the random field as in Definition~\ref{d.omegatmixing-field} was taken as a constitutive assumption,  allowing one to obtain various limiting results in geostatistical models.  
In this section, we provide a tool (Lemma \ref{l.decorrelation-field-new}) to establish the mixing of the random field 
when it is itself constructed from a {\em covariate point process} and we illustrate it with the example of the shot-noise random field in Corollary \ref{cor:shot-nois-field}.

Consider a marked point process~\(\tXi = \{(x, U(x)) \}_{x \in \Xi} \) on~\(\cN_{\R^d \times \K}\) with marks in a Polish space~\(\K\).
Let 
\[
\digamma : \R^d \times \cN_{\R^d \times \K} \longrightarrow \M
\]
be a measurable function, where~\(\M\) is a Polish space.
Define the random field on~\(\R^d\) with values in~\(\M\) by
\[
M(x) = \digamma(x, \tXi), \qquad x \in \R^d.
\]
The process~\(\tXi\) plays the role of a \emph{covariate point process} for this field.

If the function~\(\digamma\) is either \BL-localizing or stabilizing on the marked point process~\(\tXi\),
then the random field~\(M\) exhibits, respectively, \BL-mixing or \(\B\)-mixing correlations, 
in the sense of Definition~\ref{d.omegatmixing-field}.

\begin{lemma}[Fast mixing correlations of random fields with a covariate point process] \label{l.decorrelation-field-new}
Let $\tXi$ be a simple, marked point process on $\R^d\times\K$ having
summable exponential $\B$-mixing correlations  as in Definition~\ref{def.A2},
and with decay function $\hat \omega_k(s) = \hat C_k \hat \phi(s)$.
 Consider the random field $M(x):=\digamma(x, \tXi)$, $x \in \R^d$, where
$\digamma$ takes values in the Polish space $\M$ and is   \BL-localizing as in Definition~\ref{def.Lp-stabilizing_marking}~\ref{i.BL-localizing-windows}-\ref{i.BL-localizing-fast}
(resp. fast stabilizing as in Definition~\ref{def.stabilizing_marking}\ref{i.stabilizing} -\ref{i.stabilizing-fast}) but under the unconditional probability $\mP$ of $\tXi$ in both cases. 
Then $M$  exhibits fast \BL-mixing correlations as in Definition~\ref{d.omegatmixing-field}\ref{i.BL-mixing-field}-\ref{i.mixing-field-fast})  (respectively, fast $\B$-mixing correlations as in  
Definition~\ref{d.omegatmixing-field}\ref{i.B-mixing-field}-\ref{i.mixing-field-fast}).
  \end{lemma}

  \begin{proof}
Consider the restricted score $\digamma^{(r)}$ (i.e., $\digamma^{(r)}(x,\tXi) = \digamma(x,\tXi \cap B_r(x)), x \in \R^d$), so that the  radius of stabilization $R^{\digamma^{(r)}}$ satisfies $R^{\digamma^{(r)}} \le r$.
 We use the FME expansion 
(Theorem~\ref{FMEthm}) and follow the proof of
Proposition~\ref{p:mppmix_double}
to  obtain
$$
\left| \sE[ f(\bk {\digamma^{(r)}}{1}{p}) g(\bk {\digamma^{(r)}}{p+1}{p+q})] - \sE[f(\bk {\digamma^{(r)}}{1}{p})]\;\sE[g(\bk
  {\digamma^{(r)}}{p+1}{p+q}) \right| \leq
\hat \phi(\frac s2)\sum_{l=1}^\infty \frac{(4(p+q)\theta_dr^d)^l}{l!}\hat C_{l}.
$$
The difference between this and ~\eqref{e.hatC-bound-fine} 
is that the  ground point process $\tXi$ is not conditioned to have points at $x_1,\ldots,x_{p+q}$.
Consequently
the constants $\hat C_{\cdot}$ in ~\eqref{e.hatC-bound-fine} have no dependence on 
$(p+q)$ and only depend on $l$.
 We may extend the arguments to unrestricted score functions $\digamma$ satisfying fast \BL-localization, 
by following the proof of part~\ref{i.BL-iterated-marks} of Theorem~\ref{t:mppmix_double} 
(resp. for fast stabilization of $\digamma$, part~\ref{i.iterated-marks} of the same theorem); 
see Section~\ref{s.proofcorrstabmarks} and in particular, the bound~\eqref{Lip-argumentA} (resp.~\eqref{e:decayfunction2}).
\end{proof}
Random fields constructed via a point process include the shot-noise random field and we apply the above lemma to this case, thus furnishing an explicit example of a random field $M$ as in Definition~\ref{d.omegatmixing-field}.
\begin{corollary}[Shot-noise random field with fast mixing correlations]\label{cor:shot-nois-field}
Consider a simple marked point process 
\(\tXi = \{ \tilde{x} = (x, U(x)) \}_{x \in \Xi}\) on \(\R^d\) 
with i.i.d.\ marks \(U(x)\) taking values in a Polish space \(\K\).  
Let the shot-noise field \(M\) on \(\R^d\) with values in \( \R\) be defined by
\[
M(x) := \sum_{y \in \Xi} \zeta(x, \tilde{y}), \quad x \in \R^d,
\]
for some measurable function \(\zeta : \R^d \times (\R^d \times \K) \longrightarrow  \R\) 
such that \(\zeta(\tilde{x}, \tilde{y}) \le \phi(|x-y|)\) 
for a fast-decreasing function \(\phi\) as in~\eqref{varphibd}.  
If the ground point process \(\Xi\) has fast summable exponential mixing correlations as in Definition~\ref{def.A2}\ref{i.summ-exp-mixing}, 
then the random field \(M\) exhibits fast \BL-mixing correlations 
as in Definition~\ref{d.omegatmixing-field}\ref{i.BL-mixing-field} and~\ref{i.mixing-field-fast}.   
\end{corollary}

\begin{proof}
The result can be derived from Lemma~\ref{l.decorrelation-field-new} by setting $\digamma(x,\tXi) := M(x), x  \in \R^d$. Indeed, the i.i.d. marking ensures that \(\tXi\) has  fast summable exponential \(\B\)-mixing correlations via Proposition~\ref{p:mppmix_double}. Further, we assert that the score function \(\digamma\) is fast \(L^1\)-stabilizing (hence \BL-localizing, 
see Remark~\ref{rem:comparison_stabilization}\ref{ii.L2->BL}) on \(\tXi\), leveraging the boundedness of the first intensity function of the ground process. To see this, let $p \in \N$, $f$  a bounded Lipschitz function on $\R^p$ and $x_1,\ldots,x_p \in \R^d$. Using the Lipschitz property of $f$, the definition of $\digamma$, and the Campbell-Little-Mecke formula, we obtain
\begin{align*}
\big| \E f(\bk {\digamma(x_i,\tXi)}{1}{p}) - \E f(\bk {\digamma(x_i,\tXi \cap B_r(x_i)}{1}{p}) \big| & \leq \sum_{i=1}^p \E \, \big| \digamma(x_i,\tXi) - \digamma(x_i,\tXi \cap B_r(x_i) \big| \\
& = \sum_{i=1}^p \E |\sum_{y \in \Xi \cap B_r(x_i)^c} \zeta(x_i,\ty)| \\
& \leq \kappa_0 \sum_{i=1}^p \int_{y \in B_r(x_i)^c} \phi(|y-x_i|) \, \md y \\
& \leq p \kappa_0 \int_{y \in B_r(\0)^c} \phi(|y|) \, \md y \\
& \leq p \kappa_0 \theta_{d-1} \int_r^{\infty} \phi(s) s^{d-1}  \md s,
\end{align*}
which is fast-decreasing in $r$ as $\phi$ is fast-decreasing. Thus $\digamma$ is fast $L^1$-stabilizing and hence is also fast \BL-localizing. 
\end{proof}

\section{Further applications and directions}
\label{s:future}

In this section  we indicate some random geometric structures that deserve separate investigation and where we expect  that the general limit theory of Part~\ref{part:theory-foundations} will be applicable. 
The general theorems 
of Section~\ref{s:clta}
may be used to establish the limit theory for statistics of other geometric structures, including those having an underlying graph which is randomized or which has  
a dynamic component.  We list three such structures and include five directions for future research.

\begin{enumerate}[wide,label=(\roman*),  labelindent=0pt]
\item {\em Statistics of weighted Voronoi tessellations.} Recall, given a weight  function $\rho: \R^d \times \R^d \to \R$, 
 each $x \in \P$ generates weighted cells
$$
C(x, \P) := \left\{ y \in \R^d: \rho(y, x) \leq \rho(y, z) \  \text{for all}  \ z \in \P \right\},
$$
the collection of which generates the weighted Voronoi tessellation.  It is not difficult to establish conditions on the weights and the cell centers $\P$ under which geometric characteristics 
 (including total number of edges, total length of edges) may be expressed as a sum of 
fast \BL-localizing score functions satisfying the conditions of Proposition  \ref{expvar} and Theorem \ref{t:multcltmarkedpp}.  The work of \citet{FPY} shows that if the weight $\rho$ is bounded, then this gives rise to score functions satisfying stopping set stabilization;  if an  unbounded $\rho$ satisfies  a suitable integrability condition, then  it is likely that this would give rise to score functions satisfying fast \BL-localization. 

\item  {\em Euclidean minimal spanning trees.}   One could apply the  law of large numbers given by Proposition~\ref{expvar}
 to deduce the first order limit theory for statistics (total edge length, number of vertices of a given degree)  of the minimal spanning tree on the ergodic  point sets $\P$ in Subsection \ref{s:limitstatpp} having fast mixing  correlations. 
 This extends 
\citet[Theorem 2.3]{PY2003} and is proved by making  modifications to the proof of that theorem, which is confined to Poisson input. In particular we expect to 
 require  that (i) there is a constant $\rho \in (0, \infty)$ such that the complement of the union of closed balls of radius $\rho$ centered at points of $\P$ a.s. has no unbounded component and (ii) the infinite component in continuum percolation on $\P$ is a.s. unique. 
 
\item {\em Spatial birth-growth models.}   In the classical spatial birth-growth model,  seeds (grains) arrive independently and at random on a substrate.  Upon contact with the substrate the seeds  grow radially in all directions with independent growth rates.  One is typically interested in geometric statistics of this model, including the total volume covered at time $t_0$.  Our general results provide the limit theory for statistics of these models under  conditions which are less restrictive than the classical ones.   This includes possibly dependent arrival locations on the substrate, including arrival locations specified by the realization of a point process having fast mixing  correlations.   It also allows for  the initial grain sizes to be a function of  
 previously arrived grains falling within some fixed neighborhood and it also allows the growth of the grains to be dependent on growth rates of grains within some fixed distance.
 This gives rise to a  Boolean model having  dependencies in grain locations, initial grain sizes, as well  as grain growth rates.  Theorem \ref{t:cltmarkedpp}
  may be used to establish the central limit theory for statistics of  such  models, including, for example, statistics involving sums of scores $\xi(x,\P):= h((\Vol(x, \P, t))_{t \in [0, t_0]}), x \in \P$, where $\Vol(x, \P, t)$ is the volume contribution at time $t$ to the  birth-growth model associated to the seed at $x$ and $h$ is an appropriate measurable function.  
\end{enumerate}

We  list additional  models where we expect our general results would lead to limit theorems for relevant statistics of these structures.  These models merit an independent investigation.
\begin{enumerate}[wide,label=(\roman*),  labelindent=0pt]

\item 
\label{NEW-i.randomized_graph}
{\em Randomized graphs $\G$ on $\P$.}  The graphs underlying the spin models, diffusion models and interacting particle systems, as 
developed in Sections~\ref{sss.graph-stabilization},~\ref{NEW-sss.ID-graphs} and \ref{s:genmodelassum} respectively, 
can be naturally extended to accommodate extra randomization ---for instance, Laguerre tessellations with random weights or, more generally, weighted Delaunay/Voronoi graphs.
 In these models, each site is assigned a random weight that modifies the effective distance in the classical Euclidean construction; see \citet[Section 9.2]{chiu2013stochastic}.
 In these cases, the graph $\G$ is $\G(\hat\P)$ where $\hat \P = \{(x,W(x))\}_{x \in \P}$ is a marked point process including auxiliary marks $W(x)$, for example, the weights of the aforementioned tessellations.  
Using our framework of neighbor-marking functions $\bar{N}$ of $\hat \P$, one can extend the notion of stabilizing interaction graphs to include these randomized graphs $\G(\hat\P)$.
Graphs arising from the random connection model and from Boolean models with random grains may not exhibit this type of stabilization if the weights or grains are unbounded.
In such cases, a promising alternative approach goes as follows.

We sketch this approach in the case of interacting diffusions and it can be adapted suitably to the case of interacting particle systems as well. We extend the input process to  $\tP = \{(x,M(x),\allowbreak Z_x,\bar{N}_x)\}_{x \in \P}$, with $M(x)$  representing the initial conditions of the model (states of particles or diffusions) and $Z_x$ representing Brownian motions in the case of diffusions, and $\bar N_x$ representing the neighborhood of $x$ in the spatial random graph $\G = \G(\hat \P)$ where $\hat \P$
is marked point process. Assuming summable exponential $\B$-mixing for this  extended input process $\tP$, it would be enough to assume an `$L^1$-type stabilization' condition for the neighborhood distance
instead of the graph stabilization condition at \eqref{sdecayspin-infty}; namely it would be enough to assume a condition of the form
\begin{equation} \label{NEW-e:L1-stabl-G(P)}
\sup_{x_1, \ldots, x_p \in \R^d} \mP_{\bk x1p} \Big( \max_{y \in \bar{N}_{x_1}} |y| > s \Big) \leq \varphi'_p(s), \quad s > 0,
\end{equation}
instead of our usual stabilization condition on the interaction range
$$
\sup_{x_1, \ldots, x_p \in \R^d} \mP_{\bk x1p}
\bigl(S(x_1, \P) > s\bigr) \leq \varphi_p'(s), \quad s > 0,
$$
 for a rapidly decaying function $\varphi'_p$. However, we caution that checking the summable exponential $\B$-mixing of $\tP$ and the above `$L^1$-type stabilization' condition  may be non-trivial beyond the Poisson-Boolean or random connection model.

\item {\em Dynamic random geometric graphs (RGG).} It would be worthwhile to use our general results to establish the limit theory for statistics of dynamic graphs which evolve according to rules depending on the local geometry.  Until recently the study of statistics of the random geometric graph has been confined to the static case. 
Under appropriate conditions, our general set-up should provide the  limit theory for statistics of the RGG  on dynamic point sets, including those given by interacting particle system models.  Alternatively, one could potentially use Theorem \ref{t:cltmarkedpp} to establish the  limit theory for statistics of the RGG when the input $\P$  undergoes  dynamic re-positioning. The points in  $\P$  are allowed to undergo perturbations which may depend on nearby point configurations. This gives rise to non-static RGG which  play a role in modeling.  They  are used in the study of mobile networks, a topic considered by \citet{diaz2009large},  as well as the case of mobile random geometric graphs as considered by \citet{peres2013mobile}.
More recently,  they also  feature in modeling SIS epidemics, particularly  when studying the contact process on dynamic point sets.  See the paper of  \citet{baccelli2020computational}.  It is also of interest to consider dynamic repositioning via a dependent spatial birth-death process as
in \citet{qi2008functional,onaran2022functional,onaran2023functional}, which prove functional central limit theorems in the time domain.
  
\item  {\em Dynamic Voronoi tessellations.}  The study of statistics of the Voronoi tessellation  has been mainly confined to the static case,  namely for random input which does not evolve with time.   In many situations, including the study of algorithms in nearest neighbor search queries \citet{KMRSS},  sites are dynamically inserted and deleted, giving rise to Voronoi diagrams which evolve in time. The sites may be viewed as the realization of an interacting particle system.  Subject to appropriate conditions on the dynamics, one may potentially  use Theorem \dy{\ref{t:cltmarkedpp}} to establish the limit theory for statistics of the  Voronoi tessellation on dynamically changing point sets.  

\item {\em Shot-noise and Markov random fields.} Apart from the empirical random field in Section \ref{ss:Geostat},  two other random field models  constructed via marked point processes include  shot-noise random fields and Markov random fields (MRF).  Normal approximation of Poissonian shot-noise random fields has been investigated in \citet{lachieze2019shotnoise,lachieze2021asymptotics,last2016normal} and we expect the  results in Section \ref{s:clta} to be useful in studying general point process counterparts. To date we are unaware of a systematic investigation of the limit theory for stabilizing statistics of MRF in the continuum.  
When the underlying graph is a stabilizing interaction graph as in Section \ref{mainspinresults}, 
 we expect that the approach  described here could lead to the limit theory for general statistics of MRF, including those arising in the detection of  Gauss MRF as in \cite{CBLV}, as well as those going beyond functionals  occuring in specialized settings,
as  in \cite{AYTS}.

\item {\em Models of population genetics.} The birth-death-migration example described in Section \ref{s:exIPS} can be considered as a simple evolutionary model in population genetics. 
Using our general framework, especially Theorem  \ref{thmIPS}, we may potentially establish the  limit theory for statistics of more advanced models on general point processes, such as the Kimura's stepping stone model (see \citet{cox2002stepping}) and the  Fleming-Viot model (see \citet{dawson2014spatial}). 

\end{enumerate}

\bookmarksetup{startatroot}
 
\section*{Acknowledgements} 
\addcontentsline{toc}{section}{Acknowledgements}
The authors are grateful to Jesper M\"oller and G\"unter Last for their comments on  preliminary drafts of this work. 
B.~B{\l}aszczyszyn and D.~Yogeshwaran
 were partially supported by the IFCAM project ``Geometric statistics of stationary point processes''. B.~B{\l}aszczyszyn research was partially supported by  the European Research Council (ERC-NEMO-788851).  D.~Yogeshwaran's research was partially supported by SERB-MATRICS Grant and CPDA from the Indian Statistical Institute. Part of the work was done when DY was hosted by CNRS and Laboratoire MAP5 of Universit\'e Paris Cit\'e. He is thankful to Rajat Subhra Hazra for pointing out some references on sandpile models, to Manjunath Krishnapur for discussions on variance asymptotics, to Mathew Joseph and Kavita Ramanan for explaining various calculations regarding interacting diffusions and to Sunder Sethuraman for comments on CLT for particle systems, to Piyush Srivastava for numerous elaborate discussions on spin systems, particularly regarding Section \ref{r.Hard-core-Sinclair}. The research of J. Yukich was partially supported by a Simons collaboration grant.  He is grateful for support from the CNRS and the Institut Henri Poincar\'e, where part of this research was conducted. 
 He thanks Matthias Schulte for explaining the truncation approach in the
 proof of Theorem \ref{t:multcltmarkedpp}.
 He also thanks Anima Anandkumar for discussions and who pointed out long ago relationships between  weak spatial mixing and  decay of correlations.

\begin{appendices}

\section{Mixing point processes and interaction graphs in our framework}
\label{s.admppinteraction}

We provide examples of point processes satisfying specific assumptions 
underlying our main theoretical results in Part~\ref{part:theory-foundations} and also the study of some applications 
in Part~\ref{part:applications}. 
These include the assumptions of summable exponential mixing correlations as in 
Definition~\ref{def.A2}, which is stronger than 
the merely fast mixing correlations of 
Definition~\ref{d.omegahmixing}, as well as  bounds on the Palm correlation functions ~\eqref{bdedreduced} and~\eqref{e:bd_palmcorr}, 
which are stronger than the `usual' bounds on  correlation functions required in Item~(i) of Assumption~\ref{Ass1} and needed throughout the general theory. 
Stabilizing interaction graphs, introduced in Section~\ref{sss.graph-stabilization}, 
are employed to demonstrate the localization or stabilization of marking functions 
in the models considered in Sections~\ref{s:gibbsmarking}, 
\ref{NEW-s:id_sprg}, and~\ref{s:applnsips}.

\subsection{Mixing point processes}
As indicated after Definition \ref{def.A2}, marked point processes $\tP$ with summable exponential mixing correlations can be constructed by combining a ground point process $\P$ satisfying such an assumption together  with appropriately mixing marks, such as i.i.d. marks. Hence we now provide examples of point processes $\P$ with summable exponential mixing correlations. {\em The stationary Poisson point process is a trivial example.}
\vskip.3cm
\begin{exe}[Point processes
with summable exponential mixing correlations and bounded Palm intensity functions]

 Point processes $\P$ with summable exponential mixing correlations \eqref{eqn:sum}--\eqref{phibd} include stationary $\alpha$-determinantal point processes ($-1/\alpha \in \N$) with exponentially decreasing kernels,  
and certain Cox point processes; see \cite[Section 2.2.2]{BYY19}. Moreover, mixing of correlations for \emph{subcritical} Gibbs point processes 
is established in~\cite[Section~3.3.5]{betsch2022point}. 
Under suitable assumptions---in particular, boundedness of the Papangelou intensity 
and exponential decay of cluster-dependent interactions---using Theorem~3.34 therein, it  can be shown that subcritical Gibbs point processes satisfy~\eqref{eqn:sum}--\eqref{phibd}. 

For these Gibbs point processes, the bound on Palm correlation functions 
in~\eqref{e:bd_palmcorr} still remains to be verified. 
In what follows, we establish this property for stationary 
$\alpha$-determinantal point processes.
For  stationary determinantal point process (i.e.,  $\alpha = -1$) as follows: Using the determinantal structure, Fischer's inequality for determinants \cite[Theorem 7.8.5]{horn2012matrix} twice and stationarity, we have that
$$ \rho^{(p+q)}(\bk x1p, \bk y1q) \leq \rho^{(p)}(\bk x1p)\rho^{(q)}(\bk y1q) \leq \kappa_1^{p+q},$$
and now plugging this into \eqref{e:palmcorrprod} gives that  $\rho^{(q)}_{\bk x1p}(\bk y1q) \leq \kappa_1^q$.
Now using \cite[(1.11)]{BYY19s},  we can show  that \eqref{e:bd_palmcorr} holds for  $\alpha$-determinantal point processes with $-1/\alpha \in \N$, as the latter is an independent superposition of $-1/\alpha$ many determinantal point processes. 
\end{exe}

\subsection{Stabilizing interaction graphs}
Recall the definition of stabilizing interaction graphs given in Definition~\ref{d:stabilizing-graph},
where the interaction range $S_n(x,\mu)$ is defined as the radius of stabilization of the neighborhood marking function $\bar N$; see~\eqref{e:S_n=R_Wn} and Definition~\ref{def.Stab.Radius}.

Graphs with finite-range interactions (i.e., $x \sim y$ implies $|x-y| \leq r$ for some $r>0$) trivially admit the bound
$S_n(x,\mu) \leq r+1$ (recall that the radius of stabilization is defined to take integer values) and hence satisfy~\eqref{sdecayspin} with $\varphi'_p(s)=0$ for $s>r+1$.

We next mention some non-trivial examples.
Though these graphs are standard in the stochastic geometry literature, we shall briefly describe their constructions and indicate the relevant known results that help us to verify the condition \eqref{sdecayspin}. In the three graph examples below, we use a classical argument involving conical regions. Though similar bounds are available in the literature, we give details in planar case (i.e., $d = 2$) here tailored to our stabilization definition and for self-containment. The examples explain how the requisite stabilization bounds can be derived for stationary $\alpha$-determinantal point processes with $-1/\alpha \in \N$ and with exponentially decreasing kernels. The key to this is availability of Palm void probability bounds or concentration bounds for number counts (for example, see \eqref{e.Proba-sector}) of this point process. Though we expect such bounds to be true for more point processes, we do not pursue this here. For example, Palm void probability bounds of first two orders is derived for permanental and Gibbs point processes in \cite[Section 6.3]{klatt2025invariant}.  

Recall $W_n = [-\frac{1}{2} n^{1/d}, \frac{1}{2} n^{1/d}]^d$ and $W_{\infty} = \R^d$. Below we consider $W_n, n \in \N \cup \{\infty\}$.

\begin{exe}[Undirected $k$-nearest neighbors graph]  \label{exe:KNN}
Let $\G(\mu)$ be the  
undirected $k$-nearest neighbors graph on $\mu \in \cN_{\R^d}$,  obtained by including $x\sim y$ as an edge whenever
$y$ is one of the  $k$ nearest neighbors of $x$ and/or $x$ is one of the  $k$ nearest neighbors of $y$.
We first construct an upper bound for the interaction range $S_n(x,\mu)$
and then show that it satisfies \eqref{sdecayspin} for stationary $\alpha$-determinantal point processes,
$-1/\alpha \in \N$. Notice that it is not enough to let $S_n(x,\mu)$ be the distance to the $k\,$th nearest neighbor of $x$ in the point set $\mu$, as \eqref{neighborulespin} 
would fail in general.

Instead, in dimension $d=2$ we proceed as follows. For each $t>0$ and
$x \in \R^2$, let $\cT^{(j)}(x,t)$, $1 \leq j \leq 6$, denote six disjoint
sectors with apex at $x$, radius $t$, and opening angle $\pi/3$
(their common rotation will be fixed later).
Let
\[
\cT^{(j)}(x) := \bigcup_{t>0} \cT^{(j)}(x,t), \quad   1 \leq j \leq 6,
\]
denote the corresponding infinite cones with apex at $x$ and opening angle
$\pi/3$.
Furthermore, for $n \in \N$, define $t_n^{(j)}(x)$ as the smallest radius $t$
for which the truncated sector $\cT^{(j)}(x,t)$ already covers the entire
intersection of the infinite cone with the window $W_n$, that is,
\[
t_n^{(j)}(x)
:= \inf\Bigl\{ t \ge 0 : \bigl(\cT^{(j)}(x) \setminus \cT^{(j)}(x,t)\bigr)
\cap W_n = \emptyset \Bigr\}.
\]
We say that the sector $\cT^{(j)}(x,t)$ is \emph{saturated on $W_n$} if $t \ge t_n^{(j)}(x)$.

For $n\in\N$, $x\in W_n$ and all $1 \leq j \leq 6$,  we put 
\begin{equation}
T^{(j)}_n(x):=
\begin{cases}
\inf \{t \in [0,\infty):  |\mu \cap \cT^{(j)}(x,t) \cap W_n| \geq k + 1\}& \text{if $ |\mu \cap \cT^{(j)}(x)\cap W_n )|\ge k+1$}\\
 t^{(j)}_n(x)&\text{otherwise.}
\end{cases}\label{e.T-KNN}
\end{equation}
We claim that the interaction range of the $k$-nearest neighbors graph $\G(\mu)$ satisfies for $x\in W_n$
\begin{equation}\label{e.S<T-KNN}
S_n(x,\mu)\;\le\; 2\max_{j=1,\ldots,6} T_n^{(j)}(x,\mu)+1.
\end{equation}
Indeed, let
\(
T := \max_{j=1,\ldots,6} T_n^{(j)}(x,\mu)
\),
and in the case that 
\(
T = \max_{j=1,\ldots,6} t_n^{(j)}(x)
\),
then $2T +1 \geq \lceil \text{diam}(W_n) \rceil$ and hence trivially $S_n \leq 2T+1$.

In the opposite case, when
\(
T < \max_{j=1,\ldots,6} t_n^{(j)}(x)
\),
there exist sectors $\cT^{(j)}(x,T) $ which are not saturated on $W_n$. We restrict attention to modifications of the configuration $\mu$ in the regions of $W_n$ not covered by these sectors, and examine how such modifications may affect the edges incident to $x$ in the graph $\G(\mu)$, assuming $\mu\subset W_n$.

Observe first that all \emph{primary-type neighbors} of $x$ in $\G(\mu)$---those which are among $k$-nearest neighbors of $x$---are, clearly, in the ball $B_x(T)\cap W_n$. Moreover, this ball contains also all points of $\mu$ for which $x$ is among their $k$-nearest neighbors; we call such points \emph{secondary-type neighbors} of $x$ in $\G(\mu)$. 
Indeed, consider a point
\[
y \in W_n\cap (\mathcal{T}^{(j)}(x) \setminus B_x(T))
\]
for some sector $j$ for which the intersection  $W_n\cap (\mathcal{T}^{(j)}(x) \setminus B_x(T))$ is not empty. 
Consider the open ball $B_y^\circ(|x-y|)$ centered at $y$
with radius $|x-y|$. This ball completely covers the  sector
$\mathcal{T}^{(j)}(x,|x-y|)$ except its center $x$, and consequently
\[B_y^\circ(|x-y|)\supseteq
\mathcal{T}^{(j)}(x,|x-y|)\setminus\{x\}
\;\supseteq\;
\mathcal{T}^{(j)}(x,T)\setminus\{x\}
\;\supseteq\;
\mathcal{T}^{(j)}\bigl(x, T_n^{(j)}(x,\mu)\bigr)\setminus\{x\}.
\]
By the definition of $T_n^{(j)}(x,\mu)$, the smallest set in the above chain of inclusions 
contains (exactly)~$k$ points of~$\mu$. Consequently, $B_y^\circ(|x-y|)\setminus\{y\}$ contains at least $k$ points of $\mu$
closer to $y$ than $x$, so $x$ cannot be among the $k$ nearest neighbors of $y$.
Hence, $y$ cannot be a secondary-type neighbor of $x$ in $\G(\mu)$.

It follows that adding points to $\mu$ in $W_n$ outside $B_x(T)$ cannot create
new edges incident to $x$ in $\G(\mu)$. However, adding points outside $B_x(T)$
may destroy existing edges incident to $x$ by `attracting' points that were
previously secondary-type neighbors of $x$, and removing points there may
create new edges incident to $x$ by the inverse mechanism. Nevertheless, neither
of these effects is possible for modifications of $\mu$ outside $B_x(2T)$. 
Indeed, any point that could become, or cease to be, a secondary-type neighbor
of $x$ in $\mu$ must lie in $B_x(T)$, and is therefore closer to $x$ than any
point $y \in B_x^c(2T)$. Consequently, changes to the configuration outside
$B_x(2T)$ cannot influence the adjacency relations of $x$ in $\G(\mu)$.
Finally, adding~$1$ in~\eqref{e.S<T-KNN} ensures that $2T+1$ bounds the interaction
range $S_n(x,\mu)$ of the graph $\G(\mu)$, where $S_n(x,\mu)$ is the integer
radius of stabilization as in~\eqref{e:S_n=R_Wn} and
Definition~\ref{def.Stab.Radius}.

We now establish the probabilistic bound~\eqref{sdecayspin} under the assumption that $\P$ is a stationary $\alpha$-determinantal point process with an exponentially decreasing kernel. The first steps of the argument are general, whereas the later, model-specific steps can be adapted to other classes of input point processes along similar lines.

From the previous considerations, for $t \ge 0$ we have the bound,
$$
\mP_{\bk x1p} \bigl(S_n(x_1, \P_n) > t \bigr) 
\leq 
\mP_{\bk x1p} \bigl(\max_{j=1,\ldots,6} T_n^{(j)}(x_1, \P_n) > (t-1)/2 \bigr) 
\le
\sum_{j=1}^6 \mP_{\bk x1p} \bigl(   T_n^{(j)}(x_1, \P_n) > (t-1)/2 \bigr).
$$
Moreover,   
\begin{align}
\mP_{\bk x1p} \bigl(   T_n^{(j)}(x_1, \P_n) > t \bigr)
&=
\begin{cases}
\mP_{\bk x1p}\left( \P(\cT^{(j)}(x_1,t) \cap W_n) \le k  \right) & \text{for  $t\le t^{(j)}_n(x_1)$}  \\
0 & \text{for $t> t^{(j)}_n(x_1)$.}
\end{cases}
\label{e.Proba-sector}
\end{align}

Now we fix the orientation of the six-sector partition so that the sector boundaries are inclined at an angle of $\pi/12$ relative to the coordinate axes. With this choice, there exists a constant $c_0 = \pi/24$ such that, uniformly in $n \in \N$, $x \in W_n$, and $j = 1,\ldots,6$,
\[
\bigl|\cT^{(j)}(x,t) \cap W_n\bigr| \ge c_0\, t^2,
\qquad \text{for all } t < t_n^{(j)}(x).
\]
Then the probability in~\eqref{e.Proba-sector} admits a fast decreasing  upper bound
as a function of~$t$, which follows from standard concentration inequalities
(such as Chernoff bounds). This is because the number $k$ is fixed, whereas the
volume of the sets $\cT^{(j)}(x,t)\cap W_n$ grows with~$t$.

A further technical point is to ensure that these bounds hold uniformly with
respect to the Palm conditioning at $\bk{x}{1}{p}$. Such uniformity is established,
for instance, for stationary $\alpha$-determinantal point processes with
exponentially decaying kernels; see~\cite[Theorem~2.5]{BYY19} and
\cite[Corollary~1.10]{BYY19s}. This completes the proof of stabilization of the
$k$-nearest-neighbor graph in the plane with respect to the aforementioned input processes; i.e., it satisfies conditions of Definition~\ref{d:stabilizing-graph}.

A similar geometric construction applies in dimensions $d \geq 3$, with planar
sectors replaced by cones. Hence, $\G(\cdot)$ is a stabilizing interaction graph
on finite windows in higher dimensions as well, provided that the above
concentration bounds for the Palm distributions of the input process~$\P$ hold.
\end{exe}

\vskip.3cm

\begin{exe}[Delaunay graph] 
The Voronoi cell associated with $x \in \mu \in \cN_{\R^d}$ is the set $C(x, \mu):= \{y \in \R^d: |x - y| \leq |z - y|, \text{for all } z \in \mu, z \neq x\}$. The point $x$  is called the generator of the cell $C(x, \mu)$. The Voronoi tessellation of $\R^d$ induced by $\mu$ is the collection of cells $\{C(x,\mu)\}_{x \in \mu}$.
Given $\mu \in \cN_{\R^d}$ with $\mu$ in general position, let $\G(\mu)$ be the Delaunay graph on $\mu \in \cN_{\R^d}$, i.e. $x \sim y$ iff $C(x,\mu) \cap C(y,\mu) \neq \emptyset$.

It is well known (indeed, this is the original definition of~\citet{delaunay1934}, based on the empty circumsphere property) that 
$x \sim y$ in $\G(\mu)$ if and only if there exist points $z_1,\ldots,z_{d-1} \in \mu$, all distinct from $x$ and $y$ and pairwise distinct, such that the open ball
$B^\circ(x,y,z_1,\ldots,z_{d-1})$ circumscribed about the points $(x,y,z_1,\ldots,z_{d-1})$ contains no further points of $\mu$, that is,
\[
\mu\bigl(B^\circ(x,y,z_1,\ldots,z_{d-1})\bigr) = 0 .
\]

Using this characterization, one can bound the interaction range $S_n(x,\mu)$. We illustrate this argument in dimension $d=2$ as before.  Taking $T_n^{(j)}(x,\mu)$ as given in~\eqref{e.T-KNN} with $k=1$, we can show that
\begin{equation}
 \label{e:Sndelgraph}   
S_n(x,\mu)  \leq 2T := 2 \max_{j=1,\ldots,6} T_n^{(j)}(x,\mu).
\end{equation}
Indeed, let $y \in \R^2$ satisfy
\[
|x-y| > 2T,
\]
with $T := \max_{j=1,\ldots,6} T_n^{(j)}(x,\mu)$. Then, for any ball having both $x$ and $y$ on its boundary, the corresponding radius $r$ must satisfy $r > T$, and hence the center $w$ must satisfy $|w-x| > T$.
Using the same geometric argument as in Example~\ref{exe:KNN}, the open ball $B^\circ_w(r)$ contains the sector
\[
\cT^{(j)}\bigl(x, T^{(j)}_n(x)\bigr)\setminus \{x\}\subset B^\circ_w(r),
\]
for the index $j$ such that $w \in \cT^{(j)}(x)$.
By the definition of $T^{(j)}_n(x)$, we have
\[
\mu\Bigl( \cT^{(j)}\bigl(x, T^{(j)}_n(x)\bigr) \setminus \{x\} \Bigr) = 1,
\]
which implies that the open ball $B_w^\circ(r)$ is nonempty with respect to the points of $\mu$.
This rules out $y$ as a neighbor of $x$ in the Delaunay graph on $\mu$ and thereby verifying \eqref{e:Sndelgraph}. This argument can be extended to $d \geq 3$. 

Following verbatim the discussion in Example~\ref{exe:KNN}, it follows that $S_n(x,\P)$ satisfies~\eqref{sdecayspin} if $\P$ is a stationary $\alpha$-determinantal point process, $-1/\alpha \in \N$, with an exponentially decaying kernel. Consequently, $\G(\P)$ is a stabilizing graph on finite windows with respect to~$\P$ for such processes.
\end{exe}

\begin{exe}[Sphere of influence graph]  Given a point set $\mu \in \cN_{\R^d}$, the sphere of influence graph
 $\G(\mu)$   is constructed as follows.  For each $x \in \mu$, let $B_\textnormal{NN}(x)$ be the closed ball centered at $x$ with radius equal to the distance between $x$ and its nearest neighbor in $\mu$. The sphere of influence graph puts an edge between $x$ and $y$ iff  $B_\textnormal{NN}(x) \cap B_\textnormal{NN}(y) \neq \emptyset$.  To demonstrate the stabilization of this graph on some interesting  input processes,  we need an upper bound for $S_n$,  for which it again suffices to modify  Example~\ref{exe:KNN}

More precisely, in dimension $d=2$, for $x \in \R^2$, we consider the same six disjoint sectors
$\cT^{(j)}(x,t)$, $j=1,\ldots,6$ with apex at $x$, radius $t$, and opening angle $\pi/3$
(their common rotation will be fixed later). As before, let their saturation radii on the windows $W_n$ be $t^{(j)}_n(x)$. For $n\in\N$, $x\in W_n$, and $1 \leq j \leq 6$, we define
\begin{equation}
T^{(j)}_n(x):=
\begin{cases}
\displaystyle
\min \Bigl\{ k\in\N : 1 \leq k \leq t^{(j)}_n(x),\ 
\bigl|\mu \cap \cT^{(j)}(x,k/3) \cap W_n\bigr| \geq 1,\\[0.4em]
\hspace{3em}\text{and}\ 
\bigl|\mu \cap \bigl(\cT^{(j)}(x,k)\setminus \cT^{(j)}(x,2k/3)\bigr) \cap W_n\bigr| \geq 1
\Bigr\},
& \text{if this set is nonempty},\\[0.6em]
t^{(j)}_n(x), & \text{otherwise.}
\end{cases}
\label{e.T-NN}
\end{equation}

In other words, $T_n^{(j)}(x)=k$ is the smallest integer such that, within the sector
$\cT^{(j)}(x,k)$, there is at least one point of $\mu$ in both the inner part
$\cT^{(j)}(x,k/3)$ and the outer annular part
$\cT^{(j)}(x,k)\setminus \cT^{(j)}(x,2k/3)$  in $W_n$. If no such $k$ exists, we set
$T_n^{(j)}(x)$ equal to the saturation radius $t^{(j)}_n(x)$ of this sector on $W_n$.

We claim that the interaction range of the sphere of influence  graph $\G(\mu)$ satisfies for $x\in W_n$
\begin{equation}\label{e.S<T-NN}
S_n(x,\mu)\;\le\; \max_{j=1,\ldots,6} T_n^{(j)}(x,\mu)+1.
\end{equation}
Indeed, the presence of a point in the inner part of each sector $\cT^{(j)}(x,k)$, $j=1,\ldots,6$, when $\max_{j=1,\ldots,6}T_n^{(j)}(x)=k$, yields that $B_{NN}(x) \subset B_{k/3}(x)$. On the other hand, the presence of points in the outer parts of the sectors
$\cT^{(j)}(x,T_n^{(j)}(x))$ ensures that, for any point $y$ with $|y-x| > k$, we have that $B_{NN}(y) \subset B_{k/3}^c(x)$ i.e., the corresponding ball $B_{\textnormal{NN}}(y)$ does not intersect $B_{\textnormal{NN}}(x)$.

To establish the probabilistic bound~\eqref{sdecayspin}, we observe again for $t\ge0$ 
$$
\mP_{\bk x1p} \bigl(S_n(x_1, \P_n) > t \bigr) \le  
\sum_{j=1}^6 \mP_{\bk x1p} \bigl(   T_n^{(j)}(x_1, \P_n) > t-1 \bigr)
$$
and  for $t\ge1$
\begin{align}\non
&\mP_{\bk x1p} \bigl(   T_n^{(j)}(x_1, \P_n) > t \bigr)\\
&\le 
\begin{cases}
\mP_{\bk x1p}\left(\bigcap_{k\in \N\cap [1,t]} \, \Big\{ \P(\cT^{(j)}(x_1,k/3)\cap W_n)=0  \right.\\
\left.  \hspace{4.5em}\text{\ or\  }
\P((\cT^{(j)}(x_1,k)\setminus \cT^{(j)}(x_1,2k/3))\cap W_n)=0 \Big\}  \right)& \text{for  $t\le t^{(j)}_n(x_1)$}\\
0 & \text{for $t>t^{(j)}_n(x_1)$.}
\end{cases}
\label{e.Proba-sector-NN}
\end{align}
Finally 
\begin{align*}
&\mP_{\bk x1p}\left(\bigcap_{k\in \N\cap [1,t]} \, \Big\{ \P(\cT^{(j)}(x_1,k/3)\cap W_n)=0
\text{\ or\  }
\P((\cT^{(j)}(x_1,k)\setminus \cT^{(j)}(x_1,2k/3))\cap W_n)=0 \Big\}  \right)\\
&\le \mP_{\bk x1p}\left(\P(\cT^{(j)}(x_1,\lfloor t\rfloor /3)\cap W_n)=0\right)
+\mP_{\bk x1p}\left(
\P((\cT^{(j)}(x_1,\lfloor t\rfloor)\setminus \cT^{(j)}(x_1,2\lfloor t\rfloor/3))\cap W_n)=0   \right).
\end{align*}

Following the discussion in Example~\ref{exe:KNN}, the probability in~\eqref{e.Proba-sector-NN} admits a fast decreasing upper bound as a function of~$t$ for input processes $\P$ whose Palm distributions satisfy standard concentration inequalities. This includes, in particular, stationary $\alpha$-determinantal point processes, $-1/\alpha \in \N$, with exponentially decaying kernels. 

Consequently, the sphere of influence graph $\G(\cdot)$ on the plane is a stabilizing graph on finite windows with respect to such processes~$\P$, for instance for stationary $\alpha$-determinantal point processes with exponentially decaying kernels. A similar geometric construction applies in dimensions $d \geq 3$.

\end{exe}

\section{Some consequences of stabilization via stopping sets}
\label{s:fast_Stab_finite_windows}

In this appendix, we develop some consequences of stopping set stabilization, as defined in Section~\ref{s:strongstab}. Though our results mainly use the weaker notion of \BL-localization, some of our results may be strengthened under the stronger assumption of stopping set stabilization.  stopping set stabilization is satisfied by many statistics of spatial random models (for example, see Sections \ref{s:applnsips}, \ref{ss:geoBM}). As it is central to many applications we therefore include formal statements and proofs of these properties.

Specifically, we work under the assumption of a real-valued marking (score) function defined on finite windows, $\xi: W_n \times \K \times \hat{\cN}_{W_n \times \K} \longrightarrow \R$, and consider its radii of stabilization $R^\xi_{W_n}$ as in Definition~\ref{def.Stab.Radius}.
Although this definition introduces radii of stabilization for general marking functions taking values in a Polish space, in this section we restrict attention to real-valued scores. However, all results, except those concerning moments of $\xi$, remain valid in this more general setting.

For $\tmu \in \cN_{\R^d \times K}$ and $x \in \tmu$, we put
\[  R^\xi_{\infty}(\tx,\tmu ):=\limsup_{n \to \infty}
  R^{\xi}_{W_n}(\tx,\tmu_n) \]
and when $R^{\xi}_{\infty}(\tx,\tmu) < \infty$, we define the {\em infinite input score}
\be 
\label{scoreconvergence}
\xi_{\infty}(\tx, \tmu) := \xi(\tx,\tmu \cap B_{R^{\xi}_{\infty}(\tx,\tmu)}(x)).
\ee
The following  result justifies defining the input score via \eqref{scoreconvergence} and provides additional properties of $\xi_{\infty}$ and $R^{\xi}_{\infty}$.

\begin{lemma}[Properties of the infinite input score $\xi_{\infty}(\tx, \tP)$]
\label{l.xi-infnity} 
Let $\xi :\R^d \times \K \times \hat\cN_{\R^d \times \K} \to \R$ be a stabilizing marking function satisfying  $R^\xi_\infty(\tx,\tP)<\infty$ a.s. for all $\tx\in\tP$ and consider $\xi_\infty$ as defined in \eqref{scoreconvergence}. The following statements hold.
\begin{enumerate}[wide,label=(\roman*),  labelindent=0pt]
\item  \label{i.SS-limit} There exist (finite) random variables $n_0:= n_0(\tx, \tP)$  such that for all $k \geq n_0$  we have  $R^{\xi}_{W_k}(\tx,\tP_k)=R^\xi_{\infty}(\tx,\tP )$  and $ \xi(\tx,\tP_k)=\xi_{\infty}(\tx, \tP)$. In particular, a.s. for all $\tx\in\tP$, we have that $\xi_{\infty}(\tx, \tP) = \lim_{n \to \infty} \xi(\tx,\tP_n) \ $ i.e., the equivalence \eqref{e:xi-to-infty} holds.
\item  \label{i.SS-limit-radius} $R^\xi_\infty$ is a.s. an upper bound for the radius of stabilization for $\xi_\infty$ as in Definition~\ref{def.Stab.Radius}, i.e., $R^{\xi_{\infty}}(\tx,\tP) \leq R^\xi_\infty(\tx,\tP)$ a.s..
\item \label{i.SS-limit-moment} If $\xi$  satisfies the $p$-moment condition~\eqref{e:xinpmom} on finite windows then ~\eqref{e:xinpmom} also holds
  for $\xi_\infty$ (i.e., replacing $\xi_{i,n}$ 
  by  
$\xi_\infty(\tilde x_i,\tP))$  with the same moment bound $M_p^\xi$.
For this implication to hold, it is enough that $\xi_{\infty}$ is defined as a limit as in \eqref{e:xi-to-infty}, without requiring that $R^\xi_\infty(\tx,\tP)<\infty$.  
\item \label{i.SS-limit-stat} If $\xi$ is translation invariant  then is so $R^\xi_\infty$ and hence $\xi_\infty$.
\end{enumerate}
Suppose now that $\xi$ is defined on  $\R^d \times \K \times
  \cN_{\R^d \times \K}$ and satisfies stopping set stabilization on $\tP$ ~\eqref{stab-infinite}  with $p=1$. Then 
\begin{enumerate}[wide,label=(\roman*),  labelindent=0pt]
\setcounter{enumi}{4}
\item \label{i.SS-2limits} $R^\xi_\infty(\tx,\tP)<\infty$ for all $\tx\in\tP$, a.s.. Consequently all properties (i)--(iv) hold  and moreover $R^\xi_\infty=R^\xi$ and $\xi_\infty=\xi$ for all $\tx\in\tP$, a.s. 
\end{enumerate}
\end{lemma}

  \begin{proof}
\noindent \ref{i.SS-limit}  By assumption, a.s. for all $\tx\in\tP$,
$\limsup_{n\to\infty}R^{\xi}_{W_n}(\tx,\tP_n)= R^\xi_\infty(\tx,\tP)<\infty$
and hence  there exists a finite random variable $n_0:= n_0(\tx, \tP)$  such that $R^\xi_{W_{n_0}}(\tx,\tP_{n_0}) =R^\xi_\infty(\tx,\tP)$ (remember, the stabilizing radii are defined as integers), $B_{R^\xi_{W_{n_0}}(\tx,\tP_{n_0})}(x)\subset W_{n_0}$ and
for all $k \geq n_0$  we have $R^{\xi}_{W_k}(\tx,\tP_k)\le
R^\xi_{W_{n_0}}(\tx,\tP_{n_0})$. Consequently, for $k\ge n_0$, by the stopping property of $R^{\xi}_{W_{n_0}}$,
$R^{\xi}_{W_k}(\tx,\tP_k)= R^\xi_{W_{n_0}}(\tx,\tP_{n_0})=R^\xi_\infty(\tx,\tP)$ and consequently 
$\xi_\infty(\tx,\tP) =\xi(\tx,\tP_k) =\xi\bigl(\tx,\tP\cap
B_{R_\infty^\xi(\tx,\tP)}(x)\bigr)$ guaranteeing also the validity of \eqref{e:xi-to-infty}.

\medskip
\noindent \ref{i.SS-limit-radius} 
We begin by noting the following general consequence of the stopping property~\eqref{e:stopsetstabrad}, which holds for all $k \in \mathbb{N}$, $\tilde\mu\in\hat\cN_{W_k\times\K}$, and $\tx\in\tilde\mu$:
\begin{equation}\label{e.R-stopping-new}
R_k(\tx,\tmu)=R_k\Big( \tx,\bigl(\tmu\cap B_{R_k(\tx,\tmu)}(x)\bigr)\cup  \bigl(B^c_{R_k(\tx,\tmu)}(x)\cap\tilde\nu\big) \Big)
\end{equation}
for all $\tilde\nu\in\hat\cN_{{W_k}\times\K}$,
where for $k \in \N \cup \infty$ we make the notational simplification $R_k:=R^\xi_{W_k}$. Hence, the stabilization radius $R_k$ stabilizes itself, a direct consequence of defining $R$ as the {\em smallest} $r \in \mathbb{N}$ satisfying~\eqref{eq:stopping set}. 

To prove that $R_\infty := R^\xi_\infty$ is an upper bound for the radius of stabilization of $\xi_\infty$, we need to justify  
the equality \begin{align}\label{e.SS-limit-radius}
&\xi_\infty\Bigl(\tx,\bigl(\tP\cap B_{R_\infty(\tx,\tP)}(x)\bigr)\Bigr)=\xi_\infty\Bigl(\tx,\bigl(\tP\cap
B_{R_\infty(\tx,\tP)}(x)\bigr)\cup\bigl(\tilde\nu\cap
B^c_{R_\infty(\tx,\tP)}(x)\bigr)\Bigr)
\end{align}
for almost all  input  $\tP$ and  
 all (possibly infinite)  $\tilde\nu\in\cN_{\R^d \times \K}$.
In this regard, denote the point process arguments on both sides of~\eqref{e.SS-limit-radius} respectively 
by $\Q_1:=\tP\cap B_{R_\infty(\tx,\tP)}(x)$
and $\Q_2:=\bigl(\tP\cap
B_{R_\infty(\tx,\tP)}(x)\bigr)\cup\bigl(\tilde\nu\cap
B^c_{R_\infty(\tx,\tP)}(x)\bigr).$
With the definition of $\xi_\infty$ at ~\eqref{scoreconvergence} and this notation,  equality \eqref{e.SS-limit-radius}  is  equivalent to
\begin{equation}\label{e.SS-limit-radius1} 
 \xi(\tx, \Q_1\cap B_{R_\infty(\tx,\Q_1)})
=\xi(\tx, \Q_2\cap B_{R_\infty(\tx,\Q_2)}),
\end{equation}
which is itself implied by
\begin{equation}\label{e.SS-limit-radius2} 
\Q_1\cap B_{R_\infty(\tx,\Q_1)} = \Q_2\cap B_{R_\infty(\tx,\Q_2)} =\tP\cap B_{R^\xi_\infty(\tx,\tP)}.
\end{equation}
To prove this last equality we study first  the function  $R_\infty(\tx,\cdot) :=\limsup_{n \to \infty}R_n(\tx,\cdot)$ applied to $\Q_1$ and $\Q_2$. 
As in point (i) above, we take $n_0:= n_0(\tx, \tP)$  such that $R_{n_0}(\tx,\tP_{n_0}) =R_\infty(\tx,\tP)$,  $B_{R_{n_0}(\tx,\tP_{n_0})}(x)\subset W_{n_0}$, and
for all $k \geq n_0$  we have $R_k(\tx,\tP_k)\le
R_{n_0}(\tx,\tP_{n_0})$.
Then, for $k>n_0$,  by $\Q_1\in \hat\cN_{W_k\times\K}$, we have 
$$R_k(\tx,\Q_1)=R_{k}(\tx,\tP_k)\le R_{n_0}(\tx,\tP_{n_0})=R_\infty(\tx,\tP),$$
where we used~\eqref{e.R-stopping-new} for the first equality.
Since there are infinitely many values $k\in\N$ such that $R_{k}(\tx,\tP_k)= R_{n_0}(\tx,\tP_{n_0})=R_\infty(\tx,\tP)$ we have  
 $R_\infty(\tx,\Q_1)=  R_\infty(\tx,\tP)$.

Similarly, $R_\infty(\tx, \Q_2)= R_\infty(\tx, \tP)$ given that $\tilde\nu$ is finite. Indeed, in this case, it suffices to increase the value of $k$ used for $\Q_1$, if necessary,  choosing it such that $\tilde\nu\in \hat{\cN}_{W_k\times\K}$. 
With the definition of $\Q_1$ and $\Q_2$,  this 
 justifies~\eqref{e.SS-limit-radius2} for finite $\tilde\nu$.

Suppose now $\tilde\nu$ is infinite. Then recall 
$R_\infty(\tx,\Q_2) = \limsup_{n \to \infty}R_{n}(\tx,\Q_2\cap W_n)$.
For all large enough  $n\in\N$, by the previous argument, we have 
$R_\infty(\tx,\Q_2\cap W_n)=R_\infty(\tx,\tP)$, which establishes \eqref{e.SS-limit-radius2} for $\tilde\nu$ infinite and thus concludes the proof 
that $R_\infty$ is an upper bound for the radius of stabilization of $\xi_\infty$.

\medskip
\noindent \ref{i.SS-limit-moment} 
For $1\le q\le p$, consider a bounded, non-negative function $h:\R^q\to[0,\infty)$  with bounded support. Using the  Campbell-Little-Mecke formula for
  $\xi_\infty$ and $\xi$, we have 
\begin{align*}
&\hspace{-5em}\int_{\R^q}h(\bk x1l)\E_{\bk x1l}[\max(1,|\xi_{\infty}(\tilde x_1,
     \tP)|^q)]\, \rho^{(q)}(\bk x1l) \, \md \bk x1l\\
&=
  \E\Bigl[\sum_{\text{distinct } \bk x1l\subset \P}h(\bk x1l)
     \max(1,|\xi_{\infty}(\tilde X_1, \tP)|^q)\Bigr] \\
  &=\E\Bigl[\sum_{\text{distinct } \bk x1l\subset \P}h(\bk x1l)
    \lim_{n\to\infty} \max(1,|\xi(\tilde X_1, \tP_n)|^q)\Bigr]  \\
 &\le\liminf_{n\to\infty}\E\Bigl[\sum_{\text{distinct }  \bk x1l\subset \P}h(\bk x1l)
     \max(1,|\xi(\tilde X_1, \tP_n)|^q)\Bigr]  \\
   &=\liminf_{n\to\infty} \int_{\R^q}h(\bk x1l)\E_{\bk x1l}[\max(1,|\xi(\tilde x_1,
     \tP_n)|^q)] \, \rho^{(q)}(\bk x1l) \, \md \bk x1l\\
   &\le M^\xi_p \int_{\R^q}h(\bk x1l)\, \rho^{(q)}(\bk x1l) \,\md \bk x1l
\end{align*}
where the first inequality holds by Fatou's lemma and the last one by
the $p$-moment assumption~\eqref{e:xinpmom} for $\xi$ on finite  windows.  
This shows that $\xi_\infty$ satisfies  the $p$-moment condition, provided it is a.s.   well defined as the limit  (or even~$\liminf$) in  \eqref{e:xi-to-infty}.

\medskip
\noindent \ref{i.SS-limit-stat} Assume $\xi$ is translation invariant and
consider  $a\in\R^d$. Note, in general, for some $n$ we might have 
$R^{\xi}_{W_n}(\tx,\tP\cap W_n)\not= R^{\xi}_{W_n}(\tx+a,(\tP+a)\cap W_n)$. However, when both radii  $R^\xi_\infty(\tx,\tP)$ and
$R^\xi_\infty(\tx+a,\tP+a)$ are finite (which holds
a.s.) then for all  $k$ large enough 
$R^{\xi}_{W_k}(\tx,\tP_k)=
R^\xi_{W_{n_0}}(\tx,\tP_{n_0})=R^\xi_\infty(\tx,\tP)$ 
and $R^{\xi}_{W_k}(\tx+a,(\tP+a)\cap W_k)=
R^\xi_{W_{n^a_0}}(\tx+a,(\tP+a)\cap W_{n^a_0})=R^\xi_\infty(\tx+a,\tP+a)$,
with $n_0= n_0(\tx, \tP)$ and   
$n^a_0:= n_0(\tx+a, \tP+a)$ as in point (i).
 Moreover,  $B_{R^\xi_{W_{n^a_0}}(\tx+a,\tP\cap W_{n^a_0})}(x)\subset W_k$
and $B_{R^\xi_{W_{n_0}}(\tx,\tP_{n_0})}(x+a)\subset W_k$.
Then, by the stopping property of  $R^{\xi}_{W_k}(\tx,\tP_k)$ and $R^{\xi}_{W_k}(\tx+a,\tP+a)$    and the translation invariance of $\xi$, we have 
$$
R^{\xi}_{W_k}(\tx,\tP_k)=
R^\xi_{W_k}(\tx+a,(\tP+a)\cap W_k)
$$ and, consequently $\xi(\tx,\tP_k)=\xi(\tx+a,(\tP+a)\cap W_k)=\xi_\infty(\tx,\tP_k)=\xi_\infty(\tx+a,\tP+a)$.

\medskip
\noindent \ref{i.SS-2limits} 
By~\eqref{stab-infinite} with $p=1$ and  the Campbell-Little-Mecke formula, we have $R^\xi(\tx,\tP)<\infty$ a.s. for all $\tx\in\tP$ under $\sP$.  Consequently for some finite (random variable) $n_0:= n_0(\tx, \tP)$, we have $B_{R^\xi(\tx,\tP)}(x)\subset W_{n_0}$ and hence, for $k\ge n_0$, a.s. 
we have $R^{\xi}_{W_k}(\tx,\tP_k) \le 
R^\xi(\tx,\tP)<\infty$ , and consequently $\xi(\tx,\tP_k)=\xi(\tx,\tP)$. Thus parts (i)-(iv) follow as well as the identities $R^{\xi}_{\infty}(\tx,\tP) = R^{\xi}(\tx,\tP), \xi_{\infty}(\tx,\tP) = \xi(\tx,\tP)$ for all $\tx \in \tP$, a.s..
\end{proof}  

Finally, we summarize the relationship between  stopping set  stabilization and BL-localization. The main message is that the former implies the latter.

\begin{lemma}[Stopping set stabilization implies \BL-localization]
\label{l:xi-infty-BL} 
Let  $\tP$  be a marked point process on $\R^d\times\K$ with non-null, finite intensity $\rho$.  Let $\xi : \R^d \times \K \times \hat\cN_{\R^d \times \K} \to \R$ be a marking function. 
\begin{enumerate}[wide,label=(\roman*),labelindent=0pt]
\item \label{i.SS-BL-general}
If $\xi$ is stabilizing (resp. fast stabilizing) on finite windows as in  Definition~\ref{def.stabilizing_marking}\ref{i.stabilizing-windows} (resp. ~\ref{i.stabilizing-fast}) then $\xi$ is \BL-localizing (resp. fast \BL-localizing) on finite windows as in ~Definition~\ref{def.Lp-stabilizing_marking}\ref{i.BL-localizing-windows} (resp.~\ref{i.BL-localizing-fast}). The same implications hold if $\xi$ is defined on the whole space $\R^d \times \K \times \cN_{\R^d \times \K}$, involving the stopping set Definition~\ref{def.stabilizing_marking}\ref{i.stabilizing} and the \BL \,stabilizing Definition~\ref{def.Lp-stabilizing_marking} \ref{i.BL-localizing}.

\item \label{i.SS-BL-infity}
If  $\xi$ is stabilizing (resp. fast stabilizing) on finite windows and if moreover 
$$
\tilde{\xi}_{\infty}(\tx_i,\tP) := \lim_{n \to \infty} \xi(\tx_i,\tP_n)
$$
exists,  $\sP_{x_1,\ldots,x_p}$-a.s. for all $p \in \N$ and all $x_1,\ldots,x_p \in \mR^d$ then $\tilde{\xi}_{\infty}$ is \BL-localizing  (resp. fast \BL-localizing) as in \eqref{Lp-stab-infinite}.  
\end{enumerate}
\end{lemma}

The main thrust of part (ii) above is that \BL-localization of the limiting score can be deduced without imposing finiteness of the radius of stabilization $R_{\infty}^{\xi}$. This offers the advantage that we can derive limit theorems for fast (stopping set) stabilizing score functions on finite windows without requiring uniform bounds on the tail of the stabilization radii.

\begin{proof}
\noindent \ref{i.SS-BL-general}
The arguments are given in  Remark~\ref{rem:comparison_stabilization}(ii). \\

\noindent \ref{i.SS-BL-infity} Fix $p \in \N$, $x_1,\ldots,x_p \in \mR^d$ and $f \in \BL(\R^p)$. Without loss of generality, let $r \geq 1$ and choose $n > (2(r +  \max_{i=1,\ldots,p}|x_i|))^d$. Thus we have that $\P_n \cap B_r(\bk x1p) = \P \cap B_r(\bk x1p)$. Using the fast \BL-localization of $\xi$ on finite windows (which follows from part (i) and the equivalence $\tilde{\xi}_{\infty} = \xi$ on finite point-sets, we obtain
\begin{align*}
& \bigl| \mE_{\bk x1p}
\bigl[ f \bigl( \bk {\txi_{\infty}}1p (\tx,\tP)\bigr) \bigr] - \mE_{\bk x1p} \bigl[ f \bigl( \bk {\txi_{\infty}}1p (\tx,\tP\cap B_r(\bk x1p) ) \bigr)
\bigr] \bigr|  \\
& \leq \bigl|\mE_{\bk x1p}
\bigl[f\bigl(\bk {\txi_{\infty}}1p( \tx,\tP)\bigr) \bigr] - \mE_{\bk x1p}
\bigl[f\bigl(\bk\xi1p( \tx,\tP_n )\bigr)
\bigr] \bigr| \\
& \quad +  \bigl|\mE_{\bk x1p}
\bigl[f\bigl(\bk {\txi_{\infty}}1p( \tx,\tP_n)\bigr) \bigr] - \mE_{\bk x1p}
\bigl[f\bigl(\bk\xi1p( \tx,\tP_n \cap B_r(\bk x1p) )\bigr)
\bigr] \bigr|  \\
& \leq \bigl|\mE_{\bk x1p}
\bigl[f\bigl(\bk {\txi_{\infty}}1p( \tx,\tP)\bigr) \bigr] - \mE_{\bk x1p}
\bigl[f\bigl(\bk\xi1p( \tx,\tP_n )\bigr)
\bigr] \bigr| + 2 \varphi_p(r).
\end{align*}
Letting $n \to \infty$ in the first term and using the  assumed a.s. convergence of $\xi( \tx_i,\tP_n)$ to $\txi_{\infty}( \tx_i,\tP)$ for $i=1,\ldots,p$, we obtain fast \BL-localization of $\txi_{\infty}$.
\end{proof}
As advertised after Theorem \ref{t:multcltmarkedpp}, we now show that the localization assumptions of all Palm orders therein can be substituted by stabilization assumptions of Palm orders $p \in \{1,2\}$. 
\begin{corollary}
\label{c:multicltstab12}
Consider a stationary marked point process  $\tP$ on $\R^d\times\K$ with non-null, finite intensity $\rho$, having summable exponential $\B$-mixing correlations as in Definition~\ref{def.A2}.  Let $\xi_i : \R^d \times \K \times \hat\cN_{\R^d \times \K} \to \R$,  $i =1,\ldots,{m}$ be translation invariant  score functions such that they are fast-stabilizing on finite windows of $\tP$ as in \eqref{stab} for $p\in \{1,2\}$. Assume that for {\em some} $\varepsilon > 0$, the $(2 + \varepsilon)$-moment condition \eqref{e:xinpmomcopied} holds uniformly for  $r > 0$ large enough, for all $\xitr_i$, $i = 1,\ldots, m$.
Then the conclusions of Theorem \ref{t:multcltmarkedpp} hold.
\end{corollary}
\begin{proof}
\smallskip
We consider the linear combination   \(\mu_n^{\xi}(f) = \sum_{i=1}^{m} t_i \mu_n^{\xi_i}(f)\) of the score functions \(\xi_i\) with given fixed coefficients \(t_i \in \mathbb{R}\), \(i=1,\ldots,{m}\) and 
for some given function $f \in \B(W_1)$. Now, instead of truncating the value of  \(\mu_n^{\xi}(f)\) we  restrict individually the score functions 
$\xi_i^{(r)}(\tx, \tmu):= \xi_i( \tx,\tmu \cap B_r(x))$, \(i=1,\ldots,{m}\), for $r>0$, and consider  their differences  $\Delta\xitr_i:=\xi_i-\xitr_i$. Furthermore, we consider the respective linear combinations 
 $\xitr:=\sum_{i=1}^{m} t_i \xitr_i$,  $\Delta\xitr:=\sum_{i=1}^{m} t_i \Delta\xitr_i$, and their integrals 
\(\mu_n^{\xitr}(f) = \sum_{i=1}^{m} t_i \mu_n^{\xitr_i}(f)\) and 
\(\mu_n^{\Delta \xitr}(f)=
\sum_{i=1}^{m} t_i   \mu_n^{\Delta\xitr_i}(f)\).  
The limit~\eqref{e.var-xi} for  the variance of $\mu_n^{\xi}(f)$ exists as before (under fast stabilization---hence fast \BL-localization---under Palm distributions of order $p =\{1,2\}$) and gives 
$$
\lim_{n \to \infty} n^{-1} \Var \,\mu_n^{\xi}(f) =: \sigma^2. 
$$
We 
consider two other limits:
\begin{align*}
\lim_{n \to \infty} n^{-1} \Var \,\mu_n^{\xitr}(f) &=: \sigma^2_r,\\
\lim_{n \to \infty} n^{-1} \Var \,\mu_n^{\Delta\xitr}(f) &=: v^2_r. 
\end{align*}
These limits exist because the fast stabilization of $\xi_i$ of order  $p \in \{1,2\}$  implies stabilization of the same order of $\xitr_i$ and $\Delta\xitr_i$.
(This statement for $\Delta\xitr_i$ does not necessarily hold under \BL-localization of $\xi_i$.)
Moreover, 
since $\xitr_i$ is stabilizing  at fixed radius $r$ with respect to all Palm orders (and hence fast $\BL$-localizing) and since, by assumption, it satisfies a $(2+\epsilon)$- 
moment condition, it follows by Theorem~\ref{t:multcltmarkedpp} that 
\(\mu_n^{\xitr}(f)\) satisfies the CLT as \(n \to \infty\), provided $\sigma^2_r>0$.   
Now, using the same arguments developed in the proof of the central limit theorem for  $(\mu_n^{\xi}(f))_{n \in \N}$ in Theorem \ref{t:multcltmarkedpp} via three variance asymptotics,  
it remains  to establish the limits~\eqref{V1} and~\eqref{V2} for $\sigma^2_r$ and~$v^2_r$ as $r\to\infty$. 

In this regard, observe that the vector $\bxi^{(r)}:=\bigl(\xi_i,\xitr_i,\Delta\xitr_i\bigr)_{i\in\{1,\ldots,m\}}$ is fast stabilizing, hence (jointly) fast  \BL-localizing of Palm order $p \in \{1,2\}$ (this restriction of $p$ is due to the assumption on $\xi$).
Denote by  $\bxil^{(r)}_{\0}$ and  $\bxil^{(r)}_{(\0,z)}$ the \BL-limits of Palm distributions of order 1 and 2, of
$\bxi^{(r)}$:
\begin{align*}
\bxi^{(r)}(\tilde\0, \tP_n) &\xrightarrow[n \to \infty]{\BL, \0} \bxil^{(r)}_{\0},\\ 
\Bigl(\bxi^{(r)}(\tilde\0, \tP_n), \bxi^{(r)}(\tz, \tP_n)\Bigr) &\xrightarrow[n \to \infty]{\BL, (\0,z)} \bxil^{(r)}_{(\0,z)} \stackrel{\rm{law}}{=} \Bigl(\bxil^{(r)}_{(\0, z)}[1],(\bxil^{(r)}_{(\0, z)}[2]\Bigr).
\end{align*}
By the $(2+\epsilon)$-moment conditions the above limits hold also for the expectations
and all cross-products. Moreover, by the triangle inequality, the moments of limits involving $\Delta\xitr_i$ are zero. 
Indeed, observe for $i\in\{1,\ldots,m\}$,  \BL-distance between $\xi_i^{(r)}(\tx_n,\tP_n)$ and $\xi_i(\tx_n,\tP_n)$  goes to~0 for $r\to\infty$ uniformly in $n\in\N$  (this is \BL-localization, here considered for two Palm orders only).

In order to lift the above observations to the linear combinations of $\xi$ and $\Delta\xitr$,  consider the function  $h$ mapping $(x^1_i,x^2_i,x^3_i)_{i\in\{1,\ldots,m\}}\in (\R^{3})^m$ to $\sum_{i=1}^m t_i (x^1_i-x^2_i)\in\R$.
This function is Lipschitz. Consequently,  by the $(2+\epsilon)$-moment condition which is satisfied by  all  co-ordinates of $\bxi^{(r)}$, the \BL-limits
of random variables $\Delta\xitr$ satisfy  
\begin{align*}
\Delta\xitr(\tilde\0, \tP_n) &\xrightarrow[n \to \infty]{\BL, \0} h^{(r)}(\xil^{(r)}_{\0}),\\ 
\Bigl(\Delta\xi^{(r)}(\tilde\0, \tP_n), \Delta\xi^{(r)}(\tz, \tP_n)\Bigr) &\xrightarrow[n \to \infty]{\BL, (\0,z)}  \Bigl(h^{(r)}(\xil^{(r)}_{(\0, z)}[1])),(h^{(r)}(\xil^{(r)}_{(\0, z)}[2]))\Bigr).
\end{align*}
Again, by $(2+\epsilon)$-moment conditions the above limits hold also for the expectations
and all cross-products and  the cross-moments of limits involving $\Delta\xitr$ are zero, by the previous observation for the individual score functions $\xi_i$.

To establish analogs of the limits~\eqref{V1} and~\eqref{V2} for $\sigma^2_r$ and $v^2_r$, it remains  to justify the passage of \(\lim_{r \to \infty}\) under the integral \(\int_{\R^d} [\ldots] \, \md z\) in~\eqref{sigdef-M} where respective variables  $\xil$ represent now the corresponding \BL-limits of $\xitr$ and $\Delta\xitr$.
This is justified again by the dominated convergence theorem, evoking, for both \(\xitr\) and \(\Delta\xitr\), \BL-mixing correlations of order 2, uniformly in \(r>0\); that is, these properties involve the common functions \(\omega_k\) in~\eqref{e.mppmixing}.   
\end{proof}

\section{Some consequences of \texorpdfstring{$L^1$}{L1}-stabilization}
\label{s:stab_lemma_Lp}

We now consider a notion of stabilization of real-valued scores $\xi$, referred to as  
{\em $L^1$-stabilization},  which is stronger than \BL-localization but weaker than stopping set stabilization. It enables us to show that the point processes $\tP^\xi_n := \sum_{x \in \P_n} \delta_{(x, \xi(\tx, \tP_n))}$ converge as $n \to \infty$ to the ground process equipped with new marks $\tP^\xil := \sum_{x \in \P} \delta_{(x, \xil(x))}$ over all of $\R^d$. In this approach, the existence of \(\tP^\xil \cap K\) arises as a consequence of 
$(\tP^\xi_n \cap K)_{n \in \N}$ forming a Cauchy sequence in  the bounded Lipschitz distance for finite point processes introduced in \citet{schuhmacher2008new}  (see also \citet{barbour1992stein}). This Cauchy property is proved using \(L^1\)-stabilization. Once this is done, we demonstrate the consistency of the limits across different subsets \(K \subset \R^d \) by leveraging the common ground processes, thus guaranteeing the new marking \(\tP^\xil\) over the whole space. Finally, we show that the limits \(\xil_{\bk x1p}\) of the Palm distributions developed in Lemma~\ref{l:BL-limits} indeed correspond to the Palm mark distributions \(\cM_{\bk x1p}^{\xil}\) of the limiting point process  \(\tP^\xil\).  The approach involving the Cauchy property is similar to the convergence of the Palm distributions of $\xi$  stated in Lemma~\ref{l:BL-limits} under \BL-localization. However this bounded \BL \,convergence for point processes holds under the stronger \(L^1\)-stabilization of the score function $\xi$ on finite windows of $\tP$  and is clearly implied by the stopping set stabilization (cf. Remark~\ref{rem:comparison_stabilization}).

A general   definition of this new stabilization is a  modification of the \BL-localization formulation in Definition~\ref{def.Lp-stabilizing_marking}~\ref{i.BL-localizing-windows}. Specifically, for $q\ge1$,
we say that   $\xi$  is a {\em  $L^q$-stabilizing marking function on finite windows of $\tP$},  if 
for all $p \in \N$ there are decreasing
 functions
 $\varphi_{p}=\varphi_{p,q}: (0, \infty) \to [0,\infty)$ such that 
 \begin{align} \label{Lq-stab-windows}
  \sup_{1 \leq n < \infty} \sup_{x_1,\ldots,x_p \in W_n} \E_{\bk x1p}[|\xi(\tx_1,\tP_n)-\xi(\tx_1,\tP_n \cap B_r(\bk x1p))|^q]^{1/q}  \leq \varphi_{p,q}(r), \ r > 0.
\end{align}
As in the case of stopping set stabilization (cf.~\eqref{stab}),
$L^q$-stabilization is formulated only for marginal Palm distributions,
provided the input process satisfies suitable moment conditions.  Also, recall that \eqref{Lq-stab-windows} implies \BL-localization as per~\eqref{Lp-stab}; see \eqref{e.1-product-norm-used}.

\begin{lemma}[Convergence of point processes with $L^1$-stabilizing scores]
\label{l:BL-limits-pp}
Let  $\tP$  be a  marked point process on $\R^d\times\K$ with uniformly bounded  first correlation function  $\rho^{(1)}\le \kappa_1<\infty$.   Let $\xi : \R^d \times \K \times \hat\cN_{\R^d \times \K} \to \R$ be a real-valued marking function (score) which is $L^1$-stabilizing on finite windows of $\tP$ as in~\eqref{Lq-stab-windows} with $p=1$ and  satisfies the first  moment condition on finite windows as in ~\eqref{e:xinpmom} with $p=1$. 
\begin{enumerate}[wide,label=(\roman*),labelindent=0pt]
\item \label{i.BL-limit-pp} Then the ground point process marked by the score function 
$\xi$ on finite windows converges weakly in the vague topology to a stationary marking  $\tP^\xil$ of this ground point process on the entire space, i.e.,
\begin{equation}\label{e.BL-limit-pp}
\tP^\xi_n=\sum_{x \in \P_n} \delta_{(x, \xi(\tx, \tP_n))}
\xrightarrow[n \to \infty]{d}
\tP^\xil=\sum_{x \in \P} \delta_{(x, \xil(x))}.
\end{equation}

\item \label{i.BL-limit-pp-moments} If moreover  for $p\in\N$ the input process $\P$ has finite moments of order $p+\epsilon$ for some $\epsilon>0$; i.e., $\E[\P(B)^{p+\epsilon}]<\infty$ for all $B\in\mathcal{B}_b$, with a uniformly bounded correlation function $\rho^{(p)}\le\kappa_p<\infty$, and $\xi$ satisfies~\eqref{Lq-stab-windows} then the  Palm distributions $\cM_{\bk x1p}^{\xil}$ of $\tP^\xil$ of order $p$ satisfy 
$$\cM_{\bk x1p}^{\xil}(\cdot) = \sP\{\xil_{\bk x1p} \in \cdot\}\quad \text{for $\rho^{(p)}(\bk x1p)\,\md \bk x1p$-almost all $\bk x1p\in\R^{dp}$,}$$
where $\xil_{\bk x1p}$  are limits of  the Palm distributions of the score function $\bk{\xi}1p(\tx,\tP_n) \xrightarrow[n \to \infty]{\BL,\bk x1p} \xil_{\bk x1p}$ stated in~\eqref{e:BL-n-limit}.
\end{enumerate}
\end{lemma}

A key element of the proof of the convergence \eqref{e.BL-limit-pp} is demonstrating the Cauchy property in the bounded Lipschitz distance $d_{2,BL}$ (defined in~\eqref{e.d2-Wasserstein} below)
for each sequence $\tP^\xi_{n,K} := \tP^\xi_n \cap K$ on  compact windows $K \subset \R^d$ as $n \to \infty$. This metric $d_{2,BL}$ completely metrizes the weak convergence of finite point processes. Thus using  the Cauchy property, we will construct the limit $\tP^\xil$ on compact windows $K$ and, by Kolmogorov consistency, extend the limit to the entire space $\R^d$. We now define the required metrics and discuss some of their properties before proving the above lemma.

\vskip.3cm

{\em Compactification of points and marks:}  
To use the framework of the upcoming \BL~metric of $d_{2,BL}$
for finite point processes in a compact space, we consider the points of the ground process within a compact window \(K \subset \R^d\). Moreover, to handle the real marks \(\xi(\tx, \tP_n) \in \R\) of these points, we use the one-point (Alexandroff) compactification \(\bar{\R} := \R \cup \{\infty\}\) of \(\R\).
On the compact space \(\bar{\R}\), we consider a complete metric \(\bar{d}\) induced and bounded by the Euclidean distance on \(\R\): for \(x, y \in \R\),
\begin{equation}\label{e:d-on-bar-R} 
\bar{d}(x, y) := \min(|x - y|,  \hbar(x), \hbar(y)) 
\end{equation}
with \(\hbar(z) = 1/(1 + |z|)\) for \(z \in \R\), and \(\bar{d}(\infty, x) := \hbar(x)\); cf~\cite{mandelkern1989metrization}.
With a compact subset \(K \subset \R^d\) representing the locations of points and the compactification \(\bar\R\) of the real space of their marks, we next consider a metric \(d_0\) on \(\frakX := K \times \bar\R\) defined as:
\begin{equation}\label{e:d0-on-X} 
 d_0(\tx, \tx') := D_{K}(|x - x'| + \bar{d}(u, u')), \quad \tx = (x, u), \, \, \tx' = (x', u') \in \frakX,
 \end{equation}
where \(D_K\) is a constant ensuring \(d_0 \leq 1\). The compact metric space \((\frakX,d_0)\) of `marked points' is both separable and complete.
\medskip

{\em Metric on the counting measures:} Following~\cite{schuhmacher2008new},  we define a metric \(d_1\) on the and unify  space of finite counting measures $\hat{\cN}_\frakX$ on \(\frakX\). Specifically, we define the symmetric map $\hat{\cN}^2_\frakX\to [0,\infty)$ for \(\tmu = \sum_{\tx \in \tmu} \delta_{\hx}\) and \(\tmu' = \sum_{\tx' \in \tmu'} \delta_{\tx'}\)  with  $\tmu'(\frakX)\ge\max(\tmu(\frakX),1)$ by 
\begin{equation} \label{e:d1-on-N}
 d_1(\tmu, \tmu') := \frac{1}{n} \left( \min_{\pi \in \Pi_n} \sum_{i = 1}^{m} d_0(\tx_i, \tx'_{\pi(i)}) + (n - m) \right)
\end{equation}
  where \(n :=  \tmu'(\frakX)\) , \(m := \tmu(\frakX)\), and  \(\Pi_n\) denotes the set of permutations of \(\{1, 2, \ldots, n\}\), and we put $d_1(\tmu, \tmu')=0$ if both counting measures are null.  
  
In what follows,  we will consider marked points  sharing the same locations on \(K\), hence we will always have \(\tmu(\frakX) = \tmu'(\frakX)\), as in \eqref{e.BL-limit-X-d2} below. In this case, the metric \(d_1\) simplifies to the one considered in~\citet{barbour1992stein}.
The metric \(d_1\) metrizes the vague topology (the same as the weak topology) of the finite counting measures $\hat{\cN}_\frakX$,  making $(\hat{\cN}_\frakX, d_1)$ locally compact, complete and separable~\cite[Proposition 2.B]{schuhmacher2008new}. 
\medskip

{\em Bounded Lipschiz distance on marked point processes:}
Finally, consider the set \(\BL(\hat{\cN}_\frakX)\) of bounded Lip(1) functions on \((\hat{\cN}_\frakX, d_1)\), the space of finite counting measures  on \(\frakX\). We define a metric \(d_{2,BL}\)
for the probability distributions on \(\hat{\cN}_\frakX\), capturing the distance between the distributions of finite point processes \(\tP,\tP'\) on \(\frakX\):
\begin{equation} \label{e.d2-Wasserstein}
d_{2,BL} (\tP,\tP') := \sup_{\mathfrak{f} \in \BL(\hat{\cN}_\frakX)} \Bigl|\sE[\frakf(\tP)] - \sE[\frakf(\tP')]\Bigr| .
\end{equation}
This metric is denoted by $d_{2}$  in~\cite{schuhmacher2008new} but we have chosen the above notation to align with our general terminology for \BL \,metrics.
We remark that $d_{2,BL}$ is a distance between the laws of the marked point processes rather than the point processes themselves. Since $(\hat{\cN}_\frakX, d_1)$ is complete and separable, the metric \(d_{2,BL}\) completely metrizes  the weak topology of point process distributions on \( \frakX\) ; \cite[Theorems 8.3.2 and 8.10.43]{bogachev2007measure},  see also \cite[Proposition 2.C]{schuhmacher2008new}. This makes it suitable for constructing limits of the marked point processes \((\tP^\xi_{n,K})_{n \in \N}\)  via the Cauchy property.
 
\medskip
\begin{proof}[Proof of Lemma \ref{l:BL-limits-pp}]
{\em Cauchy property and $d_{2,BL}$ convergence of  $\tP^\xi_{n,K}$ \jy{as $n \to \infty$}:}
For a given compact $K\subset \R^d$ we prove that  the sequence of point processes \((\tP^\xi_{n,K})_{n \in \N}\) on $\frakX=K\times\bar \R$, is a Cauchy sequence in the $d_{2,BL}$ metric, thus justifying the convergence 
 \begin{equation}\label{e.BL-limit-X-d2}
\tP^\xi_{n,K}=\sum_{x \in \P_n\cap K} \delta_{(x, \xi(\tx, \tP_n))}
\xrightarrow[n \to \infty]{d_{2,BL}}
\tP^\xil_K 
\end{equation}
to some random counting  measure $\tP^\xil_K\in\hat{\cN}_\frakX$.
Consider $1\le n\le n'$ and  a function $\frakf\in \BL(\hat{\cN}_\frakX)$.  
Coupling  $\tP^\xi_{n,K}$ and $\tP^\xi_{n',K}$ on  the common process $\tP$,  using in~\eqref{e:d1-on-N} the  identity permutation  $\pi \in
\Pi_{\P(K)}$ with respect to any numbering of the ground process, and  \eqref{e:d-on-bar-R}, \eqref{e:d0-on-X} 
we have 
\begin{align}\label{e.BL-cauchy-d2}
&\hspace{-3em}\Bigl|\sE[\frakf(\tP^\xi_{n,K})] - \sE[\frakf(\tP^\xi_{n',K})]\Bigr|\\
&\le\sE\Bigl[\frac{\1{\P(K)>0}}{\P(K)}\sum_{x \in \P\cap K}D_K \bar d(\xi(\tx, \tP_n),\xi(\tx, \tP_{n'}))\Bigr] \nonumber\\
&\le D_K\sE\Bigl[\sum_{x \in \P\cap K} |\xi(\tx, \tP_n)-\xi(\tx, \tP_{n'})|\Bigr] \nonumber\\
&=\kappa_1 D_K\int_{K}\sE_{x}[|\xi(\tx,\tP_n)-\xi(\tx,\tP_{n'})|] \,\md x. \nonumber
\end{align}
We now use $L^1$-stabilization, with 
\eqref{Lq-stab-windows} holding for $p = 1$ and with a  decreasing function $\varphi_{1}=\varphi_{1,1} \in[0,1]$. For $r=d(x,\partial W_n)$, we obtain the bound:
\begin{align}\label{e.BL-limit-X-d2-2}
&\hspace{-3em}\E_x[|\xi(\tx,\tP_n)-\xi(\tx,\tP_{n'})|] \\\nonumber
&\le \E_x[|\xitr(\tx,\tP_n)-\xi(\tx,\tP_{n})|]+\E_x[|\xitr(\tx,\tP_{n'})-\xi(\tx,\tP_{n'})|] \\\nonumber
&\le 2 \varphi_1(d(x,\partial W_n)).
\end{align}
 Substituting into~\eqref{e.BL-cauchy-d2} and applying the Monotone Convergence Theorem, we have:
\begin{equation}\label{e.BL-limit-X-d2-3}
\lim_{n\to\infty}\sup_{\mathfrak{f} \in \BL(\hat{\cN}_\frakX)} \Bigl|\E[\frakf(\tP^\xi_{n,K})] - \E[\frakf(\tP^\xi_{n',K})]\Bigr|\le
2\kappa_1D_K \int_{K}\lim_{n\to\infty}\varphi_1(d(x,\partial W_n))\,\md x=0.
\end{equation}
This confirms statement~\eqref{e.BL-limit-X-d2}.

\vskip.3cm

{\em Real marked point process  \(\tP^\xil_K\): }  The random counting measure $\tP^{\xil}_{K} \in \hat{\mathcal{N}}_\frakX$ is indeed a simple ground process in $K$ marked by values in $\overline{\mathbb{R}}$. This follows from the fact that the projection of $\tP^\xi_{n,K}$ onto $K$ is equal in distribution to $\P \cap K$ whenever $K \subset W_n$. However, in principle, the marked point process limit $\tP^{\xil}_K$ could have marks at $\infty \in \overline{\mathbb{R}}$. This is ruled out by the  moment condition, which implies the tightness of the point processes $\tP^\xi_{n,K}$ 
uniformly in $n \in \N$.
Indeed, by Campbell's formula and the \(p = 1\) moment condition \eqref{e:xinpmom}, we have for any $t > 0$,
$$
\sP(\exists x \in \P \cap K: \xi(x, \tP_n) > t) \leq 
\kappa_1 \int_{K} \sP_{x}(\xi(x, \tP_n) > t) \, \md x \leq 
\kappa_1 \, \Vol(K) \, M^\xi_1 / t.
$$
This along with \eqref{e.BL-limit-X-d2}, justifies that \(\tP^{\xil}_K\) is a  marked point processes on \(K \times \R\).
\medskip

{\em Consistency of \(\tP^\xil_K\) as \(K \uparrow \R^d\):} 
Note that the distributions of \(\tP^\xil_K\) are consistent across \(K\). Indeed, consider \(K \subset K'\), take \(\epsilon > 0\), and due to \eqref{e.BL-limit-X-d2}, choose \(n\) large enough such that 
\(d_{2,BL}(\tP^\xi_{n,K}, \tP^\xil_K) < \epsilon\) and similarly for \(K'\). By the triangle inequality  we obtain
$$
d_{2,BL}(\tP^{\xil}_K, \tP^\xil_{K'} \cap K) \leq 2\epsilon + d_{2,BL}(\tP^\xi_{n,K}, \tP^{\xi}_{n,K'} \cap K)=2\epsilon,
$$
where the last equality is by consistency of distributions of \(\tP^\xi_{n,K}\) across \(K\). This establishes the existence of the marked  point process \(\tP^\xil\) on $\R^d\times\R$ and concludes the proof of statement~\ref{i.BL-limit-pp}.
\medskip

 We next prove statement~\ref{i.BL-limit-pp-moments}. To this end, observe first
that $L^1$-stabilization as in~\eqref{Lq-stab-windows} implies \BL-localization
as in~\eqref{Lp-stab}. Hence, for every $p \in \N$,
Lemma~\ref{l:BL-limits}~\ref{i.BL-n-limit} yields the existence of limits of
the Palm distributions of the score function,
\[
\bk{\xi}{1}{p}(\tx,\tP_n)
\xrightarrow[n \to \infty]{\BL,\bk x1p}
\xil_{\bk x1p},
\]
as stated in~\eqref{e:BL-n-limit}.

These limits are pointwise limits of probability (Palm) kernels and therefore
define a probability kernel $\xil_{\bk x1p}$ from $(\R^d)^p$ to $\R^p$. It
remains to verify that $\xil_{\bk x1p}$ satisfies the
Campbell--Little--Mecke formula~\eqref{CLM} for the marks of $\tP^{\xil}$, namely the formula
\begin{align*}
\int_{(\R^d)^p} &f(\bk x1p)
\int_{\R^p} h(\bk u1p)\,
\cM_{\bk x1p}^{\xil}(\md \bk u1p)\,
\rho^{(p)}(\bk x1p)\,\md \bk x1p 
=
\int_{(\R^d)^p} f(\bk x1p)\,
\E\bigl[h(\xil_{\bk x1p})\bigr]\,
\rho^{(p)}(\bk x1p)\,\md \bk x1p,
\end{align*}
for all bounded continuous functions $f$ and $h$ with bounded support.

To justify this identity, observe first that the left-hand side gives 
\[
\E\Bigl[
\sum_{\bk{(x,u)}1p \in (\tP^{\xil})^{(p)}}
f(\bk \tx1p)\, h(\bk u1p)
\Bigr]
=
\lim_{n\to\infty}
\E\Bigl[
\sum_{\bk{(x,u)}1p \in (\tP_n^\xi)^{(p)}}
f(\bk \tx1p)\, h(\bk u1p)
\Bigr].
\]
Indeed, by~\eqref{e.BL-limit-X-d2}, the marked point processes
$\tP_n^\xi$ converge vaguely to $\tP^{\xil}$. Consequently, the above sums,
viewed as random variables, converge in distribution. The convergence of their
expectations follows from the assumed  $(p+\epsilon)$-moment of $\P$ together with the boundedness and compact support of $f$ and $h$.

Next, the passage in the right-hand side 
\begin{align*}
\lim_{n\to\infty}
\int_{(\R^d)^p} &f(\bk x1p)\,
\E_{\bk x1p}\bigl[
h(\bk{\xi}{1}{p}(\tx,\tP_n))
\bigr]\,
\rho^{(p)}(\bk x1p)\,\md \bk x1p \\
&=
\int_{(\R^d)^p} f(\bk x1p)\,
\lim_{n\to\infty}
\E_{\bk x1p}\bigl[
h(\bk{\xi}{1}{p}(\tx,\tP_n))
\bigr]\,
\rho^{(p)}(\bk x1p)\,\md \bk x1p
\end{align*}
is justified by the bounded convergence theorem. Finally,
\[
\lim_{n\to\infty}
\E_{\bk x1p}\bigl[
h(\bk{\xi}{1}{p}(\tx,\tP_n))
\bigr]
=
\E\bigl[h(\xil_{\bk x1p})\bigr]
\]
follows from the pointwise convergence of the Palm kernels and another application of the bounded convergence theorem.
This establishes
statement~\ref{i.BL-limit-pp-moments} and completes the proof of
Lemma~\ref{l:BL-limits-pp}.
\end{proof}
\end{appendices}

\phantomsection
 \addcontentsline{toc}{section}{References}
%\pdfbookmark[0]{References}{References}  
\bookmarksetup{startatroot}
\bibliographystyle{plainnat}
\bibliography{CLTMark}

\end{document}